\documentclass[12pt,leqno]{article}
\usepackage{amsmath}
\usepackage{amstext}
\usepackage{amssymb}
\usepackage{amsfonts}
\usepackage{amsthm}
\usepackage{amsopn}
\usepackage{amsbsy}
\usepackage{layout}
\makeindex

\def\newtheorems{\newtheorem{theorem}{Theorem}[section]
                 \newtheorem{cor}[theorem]{Corollary}
                 \newtheorem{prop}[theorem]{Proposition}
                 \newtheorem{lemma}[theorem]{Lemma}
                 \newtheorem{defn}[theorem]{Definition}
                 \newtheorem{definition}[theorem]{Definition}
                 \newtheorem{notation}[theorem]{Notations}

                 \newtheorem{example}[theorem]{Example}
                 \newtheorem{examples}[theorem]{Examples}
                 \newtheorem{remark}[theorem]{Remark}
                 
                 \newtheorem{question}[theorem]{Question}
                 }

\newtheorems

\font\bitlarge=cmbx12 scaled 1100

\font\nmini=cmr10 scaled 550
\font\smaller=cmti10 scaled 1000
\font\tinier=cmsy6 scaled 730

\font\srbf=cmbx8 scaled 920
\font\ssbf=cmbx8 scaled 700
\font\scmu=cmu10 scaled 600
\font\sscmu=cmu8 scaled 600

\font\slcmssifont=lcmssi8 scaled 900
\font\lcmssifont=lcmssi8 scaled 1100
\font\blcmssifont=lcmssi8 scaled 1350

\def\fs#1{\mbox{\it #1\kern 1.3pt}}
\def\rfs#1{\mbox{\rm #1\kern 1.3pt}}
\def\bfs#1{\mbox{\bf #1\kern 1.3pt}}
\def\fss#1{\mbox{\scriptsize\it #1\kern 1.3pt}}
\def\fst#1{\mbox{\tiny\it #1\kern 1.1pt}}
\def\sifs#1{\mbox{\scriptsize\it #1\kern 1.3pt}}
\def\srfs#1{\mbox{\kern0.7pt\scriptsize\rm #1\kern 1.3pt}}
\def\trfs#1{\mbox{\kern0.7pt\tiny\rm #1\kern 1.3pt}}
\def\sbfs#1{\mbox{\kern0.7pt\srbf #1\kern -0.6pt}}
\def\srbfs#1{\mbox{\kern0.7pt\srbf #1\kern -0.6pt}}
\def\spfs#1{\mbox{\kern0.7pt\scmu #1\kern 1.3pt}}
\def\sspfs#1{\mbox{\kern0.5pt\sscmu #1\kern 1.1pt}}
\def\ssbfs#1{\mbox{\kern0.7pt\ssbf #1\kern 1.3pt}}
\def\fsm#1{\mbox{\tiny\it #1\kern 1.0pt}}

\newcommand{\ssharp}{\spfs{\symbol{35}}}
\newcommand{\sssharp}{\sspfs{\symbol{35}}}

\newcommand{\comp}{\hbox{$<\kern -3pt >$}}
\newcommand{\ncomp}
		{\;\hbox{\hbox{/}\kern -9.5pt \hbox{$<\kern -3pt >$}}}

\newcommand{\meet}
	       {\hbox{$\wedge \kern -5.75pt \raise 1.5pt \hbox{$.$}\,$}}
\newcommand{\Meet}
	     {\hbox{$\bigwedge \kern -8pt \raise 0.75pt \hbox{$.$}\:$}}
\newcommand{\ld}
	       {\hbox{$< \kern -6pt \raise 2pt \hbox{$.$}\,$}}
\newcommand{\sss}{\: \hbox{$
\underline{\hbox{$\subset$}}\kern -4pt\raise -2pt \hbox{$\tiny |$}  
$}\: }
\newcommand{\almostcontained}{\: \hbox{$
\raise 1.5pt \hbox{\scriptsize $\subset$}\kern -6.3pt\raise -3.5pt
\hbox{\scriptsize $\sim$}  
$}\: }
\newcommand{\rraro}[2]{\hbox{$\kern 3pt\raise 2pt \hbox{$\raro$}
 \kern -14pt \raise
-3.5pt\hbox{\tiny{$#1\raro #2$}}$}}

\newcommand{\ct}{\centerline}

\newcommand{\frc}{\hbox{$\parallel \kern -5.7pt \hbox{$-$}$}}

\newcommand{\nrestriction}{\kern-2.5pt\upharpoonright\kern-2.5pt}
\newcommand{\nnrestriction}{\kern-4.5pt\upharpoonright\kern-4.5pt}

\newcommand{\rharp}{\upharpoonright}

\newcommand{\nfrc}{\not \kern -5pt \frc}
\newcommand{\rest}{\vbox{\hbox{$\:\kern -2pt\mathbin{\vert\kern-3.1pt\lower-1pt
   \hbox{$\mathsurround=0pt\mathchar"0012$}\kern-4pt}\:$}}}

\newcommand{\drest}{\vbox{\hbox{$\:\kern -2pt\mathbin{\rest\kern-4.6pt\lower1.7pt
   \hbox{$\mathsurround=0pt\mathchar"0012$}\kern-4pt}\:$}}}
\newcommand{\vdrest}{\vbox{\hbox{$\:\kern -2pt\mathbin{\rest\kern-3.95pt\lower1.7pt
   \hbox{$\mathsurround=0pt\mathchar"0012$}\kern-4pt}\:$}}}
   \newcommand{\smdrest}{\vbox{\hbox{$\:\kern -2pt\mathbin{\smrest\kern-2.7pt\lower1.7pt
      \hbox{$\mathsurround=0pt\mathchar"0012$}\kern-4pt}\:$\kern 1pt}}}

\newcommand{\rests}{\vbox{\hbox{\scriptsize$\:\kern
-1.4pt\mathbin{\vert\kern-2.4pt\lower-1pt
\hbox{\scriptsize$\mathsurround=0pt\mathchar"0012$}\kern-
3.0pt}\:$}}}

\newcommand{\empt}{\emptyset}

\newcommand{\sbeq}{\subseteq}

\newcommand{\pcntda}{\lower5pt\hbox{$\stackrel{\subset}{\neq}$}}
\newcommand{\pcntdb}{\lower5pt\hbox{$\stackrel{\supset}{\neq}$}}
\newcommand{\pcntdc}{\lower5pt\hbox{$\stackrel{\subseteq}{\neq}$}}

\newcommand{\m}[1]{\hbox{$ #1 $}}

\newcommand{\itPhi}{{\it \Phi}}
\newcommand{\itPsi}{{\it \Psi}}
\newcommand{\itGamma}{{\it \Gamma}}
\newcommand{\itDelta}{{\it \Delta}}
\newcommand{\itUpsilon}{{\it \Upsilon}}
\newcommand{\itSigma}{{\it \Sigma}}
\newcommand{\itOmega}{{\it \Omega}}
\newcommand{\gG}{{\it \Gamma}}

\newcommand{\gvp}{\varphi}

\newcommand{\gl}{\lambda}

\newcommand{\Kappa}{\hbox{\Large$\kappa$}}

\newcommand{\gammaprime}{\gamma\kern1pt'}

\newcommand{\bfA}{\hbox{\bf A}}

\newcommand{\bfB}{\hbox{\bf B}}

\newcommand{\calA}{{\cal A}}

\newcommand{\calB}{{\cal B}}

\newcommand{\calC}{{\cal C}}

\newcommand{\calD}{{\cal D}}

\newcommand{\calE}{{\cal E}}

\newcommand{\calF}{{\cal F}}

\newcommand{\calI}{{\cal I}}

\newcommand{\calK}{{\cal K}}

\newcommand{\calL}{{\cal L}}

\newcommand{\calM}{{\cal M}}

\newcommand{\calN}{{\cal N}}

\newcommand{\calO}{{\cal O}}

\newcommand{\calP}{{\cal P}}

\newcommand{\calQ}{{\cal Q}}

\newcommand{\calS}{{\cal S}}

\newcommand{\calT}{{\cal T}}

\newcommand{\calU}{{\cal U}}

\newcommand{\calX}{{\cal X}}

\newcommand{\baralpha}{\bar{\alpha}}

\newcommand{\hate}{\hat{e}}

\newcommand{\hatg}{\hat{g}}
\newcommand{\hath}{\hat{h}}
\newcommand{\hatk}{\hat{k}}

\newcommand{\hatu}{\hat{u}}
\newcommand{\hatv}{\hat{v}}
\newcommand{\hatw}{\hat{w}}
\newcommand{\hatx}{\hat{x}}
\newcommand{\haty}{\hat{y}}
\newcommand{\hatz}{\hat{z}}

\newcommand{\hatbeta}{\hat{\beta}}

\newcommand{\whatD}{\widehat{D}}
\newcommand{\whatF}{\widehat{F}}

\newcommand{\whatM}{\widehat{M}}

\newcommand{\whatR}{\widehat{R}}
\newcommand{\whatT}{\widehat{T}}

\newcommand{\whatX}{\widehat{X}}
\newcommand{\swhatX}{\kern6pt\widehat{\rule{0pt}{8pt}}\kern-5.9pt X}
\newcommand{\whatY}{\widehat{Y}}

\newcommand{\hattau}{\hat{\tau}}
\newcommand{\hateta}{\hat{\eta}}

\newcommand{\vecc}{\vec{c}}
\newcommand{\vecC}{\vec{C}}

\newcommand{\vecD}{\vec{D}}

\newcommand{\vecJ}{\vec{J}}

\newcommand{\vecL}{\vec{L}}
\newcommand{\vecm}{\vec{m}}
\newcommand{\vecn}{\vec{n}}
\newcommand{\vecp}{\vec{p}}
\newcommand{\vecq}{\vec{q}}
\newcommand{\vecr}{\vec{r}}
\newcommand{\vecs}{\vec{s}}
\newcommand{\vect}{\vec{t}}
\newcommand{\vecu}{\vec{u}}
\newcommand{\vecv}{\vec{v}}
\newcommand{\vecw}{\vec{w}}
\newcommand{\vecx}{\vec{x}}
\newcommand{\vecy}{\vec{y}}
\newcommand{\vecz}{\vec{z}}

\newcommand{\tildetau}{\tilde{\tau}}
\newcommand{\tildeeta}{\tilde{\eta}}

\newcommand{\tildeg}{\tilde{g}}
\newcommand{\tildeh}{\tilde{h}}

\newcommand{\wtildeF}{\widetilde{F}}

\newcommand{\wtildeX}{\widetilde{X}}
\newcommand{\swtildeX}{\kern6pt\widetilde{\rule{0pt}{8pt}}\kern-5.9pt X}
\newcommand{\wtildeY}{\widetilde{Y}}

\newcommand{\raro}{\rightarrow}

\newcommand{\itmb}[1]{\item[[\kern 1.5pt #1\kern -4pt]]}
\newcommand{\itms}[1]{\item[[#1\kern -5pt]]}

\newcommand{\II}{{\bf I\kern -1pt I}}

\newcommand{\half}{\frac{1}{2}}
\newcommand{\dghalf}{1/2}
\newcommand{\dgeighth}{1/8}
\newcommand{\third}{\frac{1}{3}}
\newcommand{\twothirds}{\frac{2}{3}}
\newcommand{\fourthirds}{\frac{4}{3}}
\newcommand{\threehalves}{\frac{3}{2}}
\newcommand{\dgthreehalves}{\dgfrac{3}{2}}
\newcommand{\quarter}{\frac{1}{4}}
\newcommand{\dgquarter}{1 / 4}

\newcommand{\fifth}{\frac{1}{5}}

\newcommand{\eighth}{\frac{1}{8}}
\newcommand{\ninth}{\frac{1}{9}}
\newcommand{\dgninth}{\dgfrac{1}{9}}

\newcommand{\dgfrac}[2]{#1 / #2}

\newcommand{\eqdf}{\hbox{\bf \,:=\,}}

\newcommand{\surj}{\vbox{\hbox{$\longrightarrow $
                  \kern -22pt \hbox{\lower 2.5pt  \hbox{\tiny onto}}
                  \kern -16pt \hbox{\raise 5pt  \hbox{\tiny 1-1}}
                  \kern 3pt}}}
\newcommand{\uarrow}[2]{\vbox{\hbox{$\longrightarrow $
                  \kern -16pt \hbox{\raise 5pt  \hbox{\tiny $#1$}}
                  \kern 10pt }}}

\newcommand{\spri}{\vbox{\hbox{\raise 2pt \hbox{\tiny $\|$}}}}
\newcommand{\tlr}{\vbox{\hbox{\raise 2pt \hbox{\tiny $\leftarrow$}}}}
\newcommand{\trr}{\vbox{\hbox{\raise 2pt \hbox{\tiny $\rightarrow$}}}}
\newcommand{\spr}{\mathrel{\vbox{\hbox{\tlr \kern -1.7pt \spri \kern -1.7pt \trr}}}}

\newcommand{\mspri}{\vbox{\hbox{\raise 2pt \hbox{$\scriptscriptstyle \|$}}}}

\newcommand{\mtlr}{\vbox{\hbox{\raise 2pt \hbox{$\scriptscriptstyle \leftarrow$}}}}

\newcommand{\mtrr}{\vbox{\hbox{\raise 2pt \hbox{$\scriptscriptstyle \rightarrow$}}}}

\newcommand{\mspr}{\mathrel{\vbox{\hbox{\mtlr \kern -1.7pt \mspri \kern -1.7pt \mtrr}}}}

\newcommand{\Nm}{\hbox{\kern -1.3pt \em I\kern-.2300em N\,}}
\newcommand{\Na}{\hbox{\kern -1.3pt \it I\kern-.2300em N\,}}
\newcommand{\Aa}{\hbox{\it A\kern -6.8pt\lower 1.0pt\hbox{-}\,}}
\newcommand{\Ba}{\hbox{\kern -1.3pt \it I\kern-.2300em B\,}}
\newcommand{\Da}{\hbox{\kern -1.3pt \it I\kern-.2300em D\,}}
\newcommand{\Ka}{\hbox{\kern -1.3pt \it I\kern-.2300em K\,}}
\newcommand{\La}{\hbox{\kern -1.3pt \it I\kern-.2300em L\,}}
\newcommand{\Ta}{\hbox{\kern 0.5pt \it T\kern-.5550em T\,}}
\newcommand{\jTa}{\hbox{\kern 0.5pt \it T\kern-.5550em T\,}}
\newcommand{\nTa}{\hbox{\kern 0.5pt \it T\kern-.6300em T\,}}
\newcommand{\Nmt}{\hbox{\kern -1.3pt I\kern-.2300em N\,}}
\newcommand{\N}{\hbox{I\kern-.1500em \hbox{\sf N}}}
\newcommand{\As}{\hbox{\scriptsize\it A\kern -5.3pt\lower 0.7pt\hbox{\it -}\,}}
\newcommand{\Ns}{\hbox{{\scriptsize\it I\kern-.1500em N}}}
\newcommand{\Nss}{\hbox{{\tiny\it I\kern-.1500em N}}}
\newcommand{\Ls}{\hbox{{\scriptsize\it I\kern-.1700em L}}}
\newcommand{\Rm}{\hbox{\kern -1.3pt \em I\kern-.1950em R\,}}
\newcommand{\Ra}{\hbox{\kern -1.3pt \it I\kern-.1950em R\,}}
\newcommand{\R}{\hbox{I\kern-.1500em \hbox{\sf R}}}
\newcommand{\sR}{\hbox{\tiny \hbox{I\kern-.1500em \hbox{\sf R}}}}
\newcommand{\msR}{\hbox{\tiny \hbox{\it I\kern-.2200em \hbox{\it R}}}}
\newcommand{\Rs}{\hbox{\scriptsize\it \hbox{I\kern-.2100em \hbox{R}}}}

\newcommand{\Q}
   {\hbox{${\rm Q} \kern -7.5pt \raise 2pt \hbox{\tiny$|$}\kern 7.5pt$}}
\newcommand{\C}
   {\hbox{${\rm C} \kern -6.5pt \raise 2pt \hbox{\tiny$|$}\kern 6.5pt$}}
\newcommand{\sC}
   {\hbox{\tiny \hbox{${\rm C} \kern -5.0pt \raise 1pt \hbox{\tiny$|$}\kern 7.5pt$}}}
\newcommand{\Z}{\hbox{\sf Z\kern-0.720em\hbox{ Z}}}
\newcommand{\sZ}{\hbox{\tiny\hbox{ \sf Z\kern-0.720em\hbox{ Z}}}}

\newcommand{\myqed}{\kern 5pt\vrule height7.5pt width6.9pt depth0.2pt
\rule{1.3pt}{0pt}}

\newcommand{\proofend}
{
\hbox{
\rule{0.8pt}{7pt}
\kern-3.9pt
\raise6.3pt
\hbox{\rule{5.0pt}{0.8pt}}
\kern-3.8pt
\rule{0.8pt}{7pt}
\kern-9.8pt
\rule{5.8pt}{0.8pt}
}
\kern-06pt}

\newcommand{\proofends}{\proofend\kern-07pt}

\newcommand{\proofendeol}
{
\hbox{
\rule{0.8pt}{7pt}
\kern-3.9pt
\raise6.3pt
\hbox{\rule{4.1pt}{0.8pt}}
\kern-3.8pt
\rule{0.8pt}{7pt}
\kern-9.0pt
\rule{4.9pt}{0.8pt}
\kern-5pt
}
\kern-8pt
}

\newcommand{\sngltn}[1]{\{ #1 \}}
\newcommand{\dbltn}[2]{\{ #1, #2 \}}

\newcommand{\limti}[1]{\lim_{\kern1pt#1 \rightarrow \infty}}
\newcommand{\limtmi}[1]{\lim_{#1 \rightarrow -\infty}}
\newcommand{\liminfti}[1]{\liminf_{#1 \rightarrow \infty}}
\newcommand{\limtz}[1]{\lim_{\kern1.0pt#1 \rightarrow 0}}

\newcommand{\fomseq}[2]{\setm{#1_{#2}}{#2 \in \Nm}}

\newcommand{\setm}[2]{\{#1\:|\:#2\}}

\newcommand{\fsetn}[2]{\{\,#1,\ldots ,#2\}}

\newcommand{\braces}[1]{\{ #1 \}}

\newcommand{\dline}[1]{\| #1 \|}
\newcommand{\norm}[1]{\| #1 \|}
\newcommand{\boldnorm}[1]
{
\kern1pt\lower3pt\hbox{\rule{1.0pt}{12pt}\kern1.4pt}
#1
\lower3pt\hbox{\kern1.6pt\rule{1.0pt}{12pt}\kern2.0pt}
                                            }

\newcommand{\bboldnorm}[1]
{
\kern1pt\lower3pt\hbox{\rule{1.2pt}{12pt}\kern1.4pt}
#1
\lower3pt\hbox{\kern1.6pt\rule{1.2pt}{12pt}\kern2.0pt}
                                            }

\newcommand{\sboldnorm}[1]
{
\kern1pt\lower1.6pt\hbox{\rule{0.9pt}{07pt}\kern1.1pt}
#1
\lower1.6pt\hbox{\kern1.1pt\rule{0.9pt}{07pt}\kern1.8pt}
                                            }

\newcommand{\sbboldnorm}[1]
{
\kern1pt\lower0.5pt\hbox{\rule{1.0pt}{06pt}\kern1.1pt}
#1
\lower0.5pt\hbox{\kern1.1pt\rule{1.0pt}{06pt}\kern1.8pt}
                                            }
\newcommand{\abs}[1]{| #1 |}

\newcommand{\pair}[2]{\langle#1 ,#2\kern1.4pt\rangle}

\newcommand{\rpair}[2]{(#1 ,#2)}
\newcommand{\scolonpair}[2]{(#1 ;#2)}

\newcommand{\uopair}[2]{\{#1 ,#2\}}
\newcommand{\trpl}[3]{\langle#1 ,#2 ,#3 \rangle}
\newcommand{\qdrpl}[4]{\langle#1 ,#2 ,#3 ,#4 \rangle}
\newcommand{\fortpl}[4]{\langle#1 ,#2 ,#3 ,#4\kern1.4pt \rangle}
\newcommand{\frtpl}[4]{\langle#1 ,#2 ,#3 ,#4 \rangle}
\newcommand{\fourtpl}[4]{\langle #1 ,#2 ,#3 ,#4 \rangle}
\newcommand{\fivtpl}[5]{\langle #1 ,#2 ,#3 ,#4 ,#5 \rangle}
\newcommand{\sixtpl}[6]{\langle #1 ,#2 ,#3 ,#4 ,#5 ,#6 \rangle}
\newcommand{\semisixtpl}[6]{\langle #1 ,#2 ,#3 ;#4 ,#5 ,#6 \rangle}

\newcommand{\change}[1]{\lower 0.7pt \hbox{\mbox{\large $#1$}}}

\newcommand{\fnn}[3]{#1:#2 \raro #3}
\newcommand{\tendin}[3]
{\lim_{#1\kern 1pt \ni\kern 1.5pt #2\kern 1pt \rightarrow\kern 1.5pt #3}}
\newcommand{\tend}[2]
{\lim_{#1\kern 1pt \rightarrow\kern 1.5pt #2}}

\newcommand{\iso}[3]{\m{ #1 : #2 \cong #3 }}
\newcommand{\isobreak}[3]{\m{ #1 : #2 \cong}\break #3 }

\newcommand{\titles}[1]{\bigskip\ct {\Large\bf #1}\bigskip \bigskip}

\newcommand{\authorss}[4]{\bigskip\ct {#1} \ct{#2}
\ct{#3} \ct{#4} }

\newcommand{\minus}{^{-1}}
\newcommand{\inverse}{^{-1}}

\newcommand{\Frechet}{Fr\'{e}chet\ }
\newcommand{\tsubseteq}{\lower2.1pt\hbox{\tiny$\subseteq$}}
\newcommand{\temptyset}{\lower1.6pt\hbox{\tiny$\emptyset$}}
\newcommand{\tSml}{\lower1.6pt\hbox{\tiny\it Sml}}
\newcommand{\tSpprtd}{\lower1.6pt\hbox{\tiny\it Spprtd}}
\newcommand{\tprec}{\lower1.6pt\hbox{\tiny$\prec$}}
\newcommand{\tcong}{\lower1.6pt\hbox{\tiny$\cong$}}
\newcommand{\tspr}{\lower3.5pt\hbox{\tiny$\spr$}}
\newcommand{\indexentry}[2]{#1 \hfill #2\\}

\newcommand{\onetoone}{\hbox{1 \kern-3.2pt -- \kern-3.0pt 1}}

\newcommand{\onetoonen}{\hbox
{1 \kern-3.2pt \rule[3.5pt]{4.2pt}{0.4pt} \kern-2.8pt 1}}
\newcommand{\onetoonet}{\hbox
{1 \kern-2.3pt \rule[3.5pt]{4.2pt}{0.4pt} \kern-3.8pt 1}}

\newcommand{\num}[1]{$(#1)$}

\newcommand{\bard}{\bar{d}}

\newcommand{\barg}{\bar{g}}
\newcommand{\barh}{\bar{h}}

\newcommand{\barn}{\bar{n}}

\newcommand{\bars}{\bar{s}}
\newcommand{\bart}{\bar{t}}
\newcommand{\baru}{\bar{u}}
\newcommand{\barv}{\bar{v}}

\newcommand{\barx}{\bar{x}}
\newcommand{\bary}{\bar{y}}
\newcommand{\barz}{\bar{z}}
\newcommand{\barzero}{\bar{0}}

\newcommand{\barnu}{\bar{\nu}}

\newcommand{\barB}{\bar{B}}

\newcommand{\overbbA}
               {\kern3.1pt
	       \overline{\kern-1.1pt\bbA\kern-0.6pt}\kern0.6pt}
\newcommand{\oversbbA}
               {\kern2.6pt
	       \overline{\kern-2.6ptsboldbbA\kern-0.3pt}\kern0.3pt}
\newcommand{\overB}
               {\kern3.1pt\overline{\kern-3.1ptB\kern-0.6pt}\kern0.6pt}
\newcommand{\overE}
               {\kern3.1pt\overline{\kern-3.1ptE\kern-0.6pt}\kern0.6pt}
\newcommand{\oversE}
               {\kern2.6pt\overline{\kern-2.6ptE\kern-0.3pt}\kern0.3pt}
\newcommand{\overssE}
               {\kern2.1pt
	       \overline{\kern-1.2ptE\rule{0.3pt}{0pt}}
	       \kern-0.3pt}
\newcommand{\overF}
               {\kern3.1pt\overline{\kern-3.1ptF\kern-0.6pt}\kern0.6pt}
\newcommand{\oversF}
               {\kern2.6pt\overline{\kern-2.6ptF\kern-0.3pt}\kern0.3pt}
\newcommand{\overG}
               {\kern3.1pt\overline{\kern-3.1ptG\kern-0.6pt}\kern0.6pt}
\newcommand{\overH}
               {\kern3.1pt\overline{\kern-3.1ptH\kern-0.6pt}\kern0.6pt}
\newcommand{\overV}
               {\kern3.1pt\overline{\kern-0.9ptV\kern-0.3pt}\kern0.6pt}
\newcommand{\overW}
               {\kern3.1pt\overline{\kern-0.9ptW\kern-0.3pt}\kern0.6pt}
\newcommand{\overX}
               {\kern3.1pt\overline{\kern-3.1ptX\kern-0.6pt}\kern0.6pt}
\newcommand{\oversG}
               {\kern2.6pt\overline{\kern-2.6ptG\kern-0.3pt}\kern0.3pt}
\newcommand{\oversX}
               {\kern2.6pt\overline{\kern-2.6ptX\kern-0.3pt}\kern0.3pt}
\newcommand{\overfs}[1]
               {\overline{\rfs{#1}}}
\newcommand{\overfss}[1]
               {\overline{\srfs{#1}}}
\newcommand{\overcalB}
               {\kern1.5pt{\overline{\kern-1.1pt\cal B}\kern-0.6pt}
	       \kern0.6pt}
\newcommand{\overcalO}
               {\kern3.1pt\overline{\kern-3.1pt{\cal O}\kern-0.6pt}
	       \kern0.6pt}

\newcommand{\barvarphi}{\bar{\varphi}}
\newcommand{\bartau}{\bar{\tau}}

\newcommand{\sstar}{\hbox{\scriptsize \kern1pt$\star$\kern0pt}}
\newcommand{\rsstar}{\raise1pt\hbox{\scriptsize \kern1pt$\star$\kern1pt}}
\newcommand{\raisedsstar}{\kern0.5pt\raise1.0pt\hbox{\scriptsize$*$}}
\newcommand{\sperp}{\raise2pt\hbox{\kern0pt\tiny $\perp$}}
\newcommand{\svee}{\raise0pt\hbox{\kern0.6pt\tiny $\vee\kern-4.8pt\vee$}}
\newcommand{\ssperp}{\raise1.5pt\hbox{\kern0pt\tinier \symbol{63}}}
\newcommand{\spreceq}{\hbox{\tiny $\preceq$}}

\newcommand{\nperp}
{
\hbox{\kern1.0pt
\hbox{\raise1.0pt \hbox{\rule{4.6pt}{0.3pt}}}\kern-5.1pt$\sperp$}
}
\newcommand{\snperp}
{
\hbox{\kern1.0pt
\hbox{\raise0.4pt \hbox{\rule{3.3pt}{0.2pt}}}\kern-3.65pt$\ssperp$}
}

\newcommand{\sor}{\raise2pt\hbox{\tiny\rm Or}}

\newcommand{\srs}[1]{\raise1.3pt\hbox{\tiny\rm #1}}
\newcommand{\ssrfs}[1]{\raise1.3pt\hbox{\tiny\rm #1}}
\newcommand{\ssrs}[1]{\raise1.3pt\hbox{\nmini #1}}
\newcommand{\scirc}
{\raise1pt\hbox{\scriptsize\kern1.5pt$\circ$\kern1.5pt}}
\newcommand{\sscirc}
{\hbox{\tiny\kern1.5pt$\circ$\kern0.3pt}}

\newcommand{\bcirc}
{\mathop{\lower1pt\hbox{\large\kern1.5pt$\circ$}}}
\newcommand{\bbcirc}
{\mathop{\lower1.8pt\hbox{\Large\kern1.5pt$\circ$}}}

\newcommand{\raisedstar}
{\raise-3.1pt\hbox{$\rule{0.5pt}{0pt}^*$}}
\newcommand{\rstar}
{\raise-3.1pt\hbox{$\rule{0.5pt}{0pt}^*$}}
\newcommand{\raisedtinystar}
{\kern0.7pt\raise1.5pt\hbox{\tiny $*$}}
\newcommand{\rtstar}
{\kern0.7pt\raise1.0pt\hbox{\tiny $*$}\kern-2.0pt}
\font\sgreek=cmr10 scaled 900
\newcommand{\sLambda}{\hbox{\sgreek \symbol{3}}}
\newcommand{\bfie}{\hbox{\blcmssifont \symbol{101}\kern1pt}}
\newcommand{\sbfie}{\hbox{\lcmssifont \symbol{101}\kern0.7pt}}
\newcommand{\bfif}{\hbox{\blcmssifont \symbol{102}\kern1pt}}
\newcommand{\sbfif}{\hbox{\lcmssifont \symbol{102}\kern0.7pt}}
\newcommand{\bfig}{\hbox{\blcmssifont \symbol{103}}}
\newcommand{\bfih}{\hbox{\blcmssifont \symbol{104}\kern1pt}}
\newcommand{\sbfih}{\hbox{\lcmssifont \symbol{104}\kern0.7pt}}

\newcommand{\bfiv}{\hbox{\blcmssifont \symbol{118}}}

\newcommand{\bfix}{\hbox{\blcmssifont \symbol{120}}}
\newcommand{\bfiy}{\hbox{\blcmssifont \symbol{121}}}

\newcommand{\bfiA}{\hbox{\blcmssifont \symbol{65}}}
\newcommand{\bfiB}{\hbox{\blcmssifont \symbol{66}}}
\newcommand{\bfiL}{\hbox{\blcmssifont \symbol{76}}}
\newcommand{\bfiM}{\hbox{\blcmssifont \symbol{77}}}
\newcommand{\bfiN}{\hbox{\blcmssifont \symbol{78}}}
\newcommand{\sbfix}{\hbox{\lcmssifont \symbol{120}}}
\newcommand{\sbfiy}{\hbox{\lcmssifont \symbol{121}}}

\newcommand{\sbfiN}{\hbox{\lcmssifont \symbol{78}}}

\newcommand{\ssbfie}{\hbox{\slcmssifont \symbol{101}}}

\newcommand{\ssbfih}{\hbox{\slcmssifont \symbol{104}}}

\newcommand{\ssbfiM}{\hbox{\slcmssifont \symbol{77}}}
\newcommand{\ssbfiN}{\hbox{\slcmssifont \symbol{78}}}

\newcommand{\hatbfix}{\hat{\bfix}}
\newcommand{\hatbfiy}{\hat{\bfiy}}

\newcommand{\mcdot}{\mathbin{\kern-1.4pt\cdot\kern-1.0pt}}
\newcommand{\lcdot}{\hbox{$\kern0.1pt\cdot\kern1.0pt$}}
\newcommand{\ncdot}{\hbox{$\kern-0.2pt\cdot\kern0.6pt$}}
\newcommand{\fr}{\kern1.0ptr}
\newcommand{\fzero}{\kern1.3pt0}
\newcommand{\fone}{\kern1.0pt1}
\newcommand{\ftwo}{\kern1.3pt2}
\newcommand{\fprime}{\kern1.3pt'}
\newcommand{\farsu}[1]{^{\kern1.3pt#1}}
\newcommand{\sprime}{\raise1.3pt\hbox{\scriptsize\kern1.4pt$'$}}
\newcommand{\fprimesub}{\kern1.5pt'\kern-3.8pt}
\newcommand{\fprimew}{\kern0.7pt'}
\newcommand{\fprimei}{\kern1.5pt'\kern-3.0pt}
\newcommand{\fnprime}{\kern1.5pt'\kern-3.0pt}

\newcommand{\clubsign}{{\raise1.3pt\hbox{\tiny $\clubsuit$}}}
\newcommand{\sclub}{{\raise1.3pt\hbox{\tiny $\clubsuit$}}}
\newcommand{\srsp}[1]{\hbox{\kern1pt\tiny $(#1)$}}
\newcommand{\sminus}{\rule{5.6pt}{0.3pt}}
\newcommand{\ssim}{\hbox{\scriptsize $\sim$}}
\newcommand{\neweq}{
\mathrel{
\hbox{
\lower-0.0pt\hbox{\sminus}
\kern-9.95pt\ssim\kern-6.2pt\raise3.8pt\hbox{\sminus}
}}}

\newcommand{\dcup}{\mathbin{\hbox{$\cup$
\kern-10.43pt\raise1.4pt \hbox{\tiny $\cup$}}}}

\newcommand{\sprt}[1]
{\lower1.7pt\hbox{\kern1.5pt\rule{0.5pt}{10pt}}
\kern-0.5pt\underline{\kern1pt#1\kern1pt}
\lower1.7pt\hbox{\kern-0.5pt\rule{0.5pt}{10pt}}\kern1pt}

\newcommand{\sprtd}[2]
{#1\lower1.7pt\hbox{\kern1.5pt\rule{0.5pt}{10pt}}
\kern-0.5pt\underline{\kern1pt#2\rule{1.9pt}{0pt}}
\lower1.7pt\hbox{\kern-0.5pt\rule{0.5pt}{10pt}}\kern1pt}

\newcommand{\sprtl}[1]
{\lower4.4pt\hbox{\kern1.5pt\rule{0.5pt}{14pt}}
\kern-0.5pt\underline{\kern1pt#1\kern1pt}
\lower4.4pt\hbox{\kern-0.5pt\rule{0.5pt}{14pt}}\kern1pt}

\newcommand{\sprtll}[1]
{\lower5.8pt\hbox{\kern1.5pt\rule{0.5pt}{15.4pt}}
\kern-0.5pt\underline{\kern1pt#1\kern1pt}
\lower5.8pt\hbox{\kern-0.5pt\rule{0.5pt}{15.4pt}}\kern1pt}

\newcommand{\sprtdl}[2]
{#1\lower4.0pt\hbox{\kern1.5pt\rule{0.5pt}{14pt}}
\kern-0.5pt\underline{\kern1pt#2\kern1pt}
\lower4.1pt\hbox{\kern0.0pt\rule{0.5pt}{14pt}}\kern1pt}

\newcommand{\sprtm}[1]
{\lower3.3pt\hbox{\kern1.5pt\rule{0.5pt}{14pt}}
\kern-0.5pt\underline{\kern1pt#1\kern1pt}
\lower3.4pt\hbox{\kern0.0pt\rule{0.5pt}{14pt}}\kern1pt}

\newcommand{\sprtdm}[2]
{#1\lower3.3pt\hbox{\kern1.5pt\rule{0.5pt}{14pt}}
\kern-0.5pt\underline{\kern1pt#2\kern1pt}
\lower3.4pt\hbox{\kern0.0pt\rule{0.5pt}{14pt}}\kern1pt}

\DeclareSymbolFont{AMSb}{U}{msb}{m}{n}
\DeclareMathSymbol{\bbA}{\mathbin}{AMSb}{"41}
\DeclareMathSymbol{\bbB}{\mathbin}{AMSb}{"42}
\DeclareMathSymbol{\bbC}{\mathbin}{AMSb}{"43}
\DeclareMathSymbol{\bbD}{\mathbin}{AMSb}{"44}
\DeclareMathSymbol{\bbE}{\mathbin}{AMSb}{"45}
\DeclareMathSymbol{\bbF}{\mathbin}{AMSb}{"46}
\DeclareMathSymbol{\bbG}{\mathbin}{AMSb}{"47}
\DeclareMathSymbol{\bbH}{\mathbin}{AMSb}{"48}
\DeclareMathSymbol{\bbI}{\mathbin}{AMSb}{"49}
\DeclareMathSymbol{\bbJ}{\mathbin}{AMSb}{"4A}
\DeclareMathSymbol{\bbK}{\mathbin}{AMSb}{"4B}
\DeclareMathSymbol{\bbL}{\mathbin}{AMSb}{"4C}
\DeclareMathSymbol{\bbM}{\mathbin}{AMSb}{"4D}
\DeclareMathSymbol{\bbN}{\mathbin}{AMSb}{"4E}
\DeclareMathSymbol{\bbO}{\mathbin}{AMSb}{"4F}
\DeclareMathSymbol{\bbP}{\mathbin}{AMSb}{"50}
\DeclareMathSymbol{\bbQ}{\mathbin}{AMSb}{"51}
\DeclareMathSymbol{\bbR}{\mathbin}{AMSb}{"52}
\DeclareMathSymbol{\bbS}{\mathbin}{AMSb}{"53}
\DeclareMathSymbol{\bbT}{\mathbin}{AMSb}{"54}
\DeclareMathSymbol{\bbU}{\mathbin}{AMSb}{"55}
\DeclareMathSymbol{\bbV}{\mathbin}{AMSb}{"56}
\DeclareMathSymbol{\bbW}{\mathbin}{AMSb}{"57}
\DeclareMathSymbol{\bbX}{\mathbin}{AMSb}{"58}
\DeclareMathSymbol{\bbY}{\mathbin}{AMSb}{"59}
\DeclareMathSymbol{\bbZ}{\mathbin}{AMSb}{"5A}
\newcommand{\boldbbA}{\mathbb{A}\kern-9pt\mathbb{A}}
\newcommand{\boldbbL}{\mathbb{L}\kern-8.3pt\mathbb{L}}
\newcommand{\boldbbR}{\mathbb{R}\kern-8.3pt\mathbb{R}}

\newcommand{\sboldbbA}
                {\hbox{\kern0.7pt\scriptsize
		$\mathbb{A}\kern-6pt\mathbb{A}$}}
\newcommand{\sboldbbL}
                {\hbox{\kern0.7pt\scriptsize
		$\mathbb{L}\kern-5.68pt\mathbb{L}$}}
\newcommand{\sboldbbN}
                {\hbox{\kern0.7pt\scriptsize
		$\mathbb{N}\kern-5.60pt\mathbb{N}$}}
\newcommand{\sboldbbR}
                {\hbox{\kern0.7pt\scriptsize
		$\mathbb{R}\kern-5.68pt\mathbb{R}$}}
\newcommand{\sboldbbT}
                {\hbox{\kern0.7pt\scriptsize
		$\mathbb{T}\kern-5.60pt\mathbb{T}$}}

\newcommand{\itbfH}{\boldsymbol{H}}
\newcommand{\itbfK}{\boldsymbol{K}}
\newcommand{\itbfL}{\boldsymbol{L}}

\newcommand{\wpm}{\kern0.7pt\pm}
\newcommand{\tLambda}{\hbox{\tiny $\Lambda$}}
\newcommand{\wcat}{\kern 1.8pt{\mathaccent 94 {\,\ }}\kern 0.7pt}
\newcommand{\concat}{\kern1.83pt\hat{\ }\kern1.4pt}
\newcommand{\sconcat}{\kern0.63pt\hat{\ }\kern0.2pt}
\newcommand{\ontpl}[1]{\langle \kern1pt #1 \kern1pt \rangle}
\newcommand{\onetpl}[1]{\langle \kern1pt #1 \kern1pt \rangle}
\newcommand{\newprec}
{\hbox{$\prec$\kern-2.6pt\raise6.3pt\hbox{\rule{2.5pt}{0.3pt}}
\kern-2.5pt\raise-0.5pt\hbox{\rule{2.5pt}{0.3pt}}\kern3pt}}

\newcommand{\iseg}
{\mathrel{\vbox{\hbox{$<$
\kern-5.1pt
\lower1.4pt\hbox{\rule{0.3pt}{2.5pt}}
\kern -4.1pt
\raise 5.0pt\hbox{\rule{0.3pt}{2.5pt}}\kern 0.5pt}}}}

\newcommand{\risegeq}
{\mathrel{\vbox{\hbox{$\geq$
\kern-12.3pt
\lower0.15pt\hbox{\rule{0.3pt}{2.5pt}}
\kern-4.25pt
\raise 6.3pt\hbox{\rule{0.3pt}{2.5pt}}\kern 0.5pt}}}}

\newcommand{\srisegeq}
{\mathrel{\vbox{\hbox{\rule{2pt}{0pt}\scriptsize$\geq$
\kern-8.7pt
\lower0.20pt\hbox{\rule{0.2pt}{1.7pt}}
\kern-3.00pt
\raise 4.3pt\hbox{\rule{0.2pt}{1.7pt}}\kern 0.5pt}}\kern6.7pt}}

\newcommand{\rme}{{\rm e}}
\newcommand{\Ehat}{E\hat{\rule{6pt}{0pt}}}

\begin{document}
\baselineskip 16pt

\setcounter{page}{1}
\titles{Reconstruction of manifolds and subsets}
\vspace{-1cm}
\titles{of normed spaces from subgroups}
\vspace{-1cm}
\titles{of their homeomorphism groups\kern0.5pt\normalsize
\footnote{
AMS Subject Classification 2000:  Primary 57N20, 46B99, 58B99;
Secondary 54E40.
\newline
Keywords: homeomorphism group, reconstruction,
uniformly continuous, bilipschitz.
}}

    \bigskip

    \authorss
   {\bf Matatyahu Rubin}
{Department of Mathematics}{Ben Gurion University of the
    Negev}{Beer Sheva, Israel\kern1pt\small
\footnote{Part of this work was carried out when the first author
visited the University of Colorado Boulder. The first author thanks the
University of Colorado for making his visit possible.
}}

    \bigskip

    \ct{ and}

    \authorss
   {\bf Yosef Yomdin}
{Department of Mathematics}{Weizmann Institute}{Rechovot, Israel}
    \bigskip

\newpage

\baselineskip 16.5pt
\begin{small}
\noindent
{\bitlarge Abstract}

\kern3mm

\noindent
This work concerns with topological spaces of the following
types: open subsets of normed vector spaces,
manifolds over normed vector spaces, the closures of open subsets
of normed vector spaces and some other types of topological
spaces related to the above.
We show that such spaces are determined by various subgroups of
their auto-homeomorphism groups.
Theorems 1\,-\,3 below are typical examples of the results obtained
in this work.
\smallskip

{\bf Theorem 1 } For a metric space $X$ let $\rfs{UC}(X)$ denote the
group of all auto-homeomorphisms $h$ of $X$
such that $h$ and $h\inverse$ are uniformly continuous.
Let $X$ be an open subset of a Banach space with the following property:
for every $\varepsilon > 0$ there is $\delta > 0$ such that for every
$u,v \in X$: if $\norm{u - v} < \delta$, then there is an arc
$L \subseteq X$ connecting $u$ and $v$ such that
$\rfs{diameter}(L) < \varepsilon$.
Suppose that the same holds for $Y$.
Let $\varphi$ be a group isomorphism between $\rfs{UC}(X)$
and $\rfs{UC}(Y)$. Then there is a homeomorphism
$\tau$ between $X$ and $Y$
such that $\tau$ and $\tau\inverse$ are uniformly continuous
and for every $g \in \rfs{UC}(X)$,
$\varphi(g) = \tau \scirc g \scirc \tau\inverse$.

See Corollaries \ref{metr-bldr-c5.7} and \ref{metr-bldr-c2.26}.
\smallskip

{\bf Theorem 2 } Let $H(X)$ denote the group of auto-homeomorphisms
of a topological space $X$.
Let $X$ be a bounded open subset of a Banach space $E$,
and denote by $\rfs{cl}(X)$ the closure of $X$ in $E$.
Suppose that $X$ has the following properties:
(1) There is $d$ such that for every $u,v \in X$ there is
a rectifiable arc $L \subseteq X$ connecting $u$ and $v$
such that $\rfs{length}(L) < d$;
such that $\rfs{length}(L) < d$;
(2) for every point $w$ in the boundary of $X$ and $\varepsilon > 0$,
there is $\delta > 0$ such that for every $u,v \in X$:
if $\norm{u - w},\norm{v - w} < \delta$,
then there is an arc $L \subseteq X$ connecting $u$ and $v$ such that
$\rfs{diameter}(L) < \varepsilon$.
Suppose that the same holds for $Y$.
Let $\varphi$ be a group isomorphism between $H(\rfs{cl}(X))$
and $H(\rfs{cl}(Y))$. Then there is a homeomorphism $\tau$ between
$\rfs{cl}(X)$ and $\rfs{cl}(Y)$ such that for every
$g \in H(\rfs{cl}(X))$,
$\varphi(g) = \tau \scirc g \scirc \tau\inverse$.

See Theorems \ref{t4.35} and \ref{metr-bldr-t6.5}(b) and
Proposition \ref{metr-bldr-p6.3}(c).
\smallskip

{\bf Theorem 3 } Let $\rfs{LIP}(X)$ denote the group of
bilipschitz auto-homeomorphisms of a metric space $X$.
Suppose that $F,K$ are the closure of bounded open subsets of $\bbR^n$,
and suppose further that $F,K$ are manifolds with boundary with an
atlas consisting of bilipschitz charts.
Let $\varphi$ be a group isomorphism between $\rfs{LIP}(F)$
and $\rfs{LIP}(K)$. Then there is a bilipschitz homeomorphism
$\tau$ between $F$ and $K$ such that
$\varphi(g) = \tau \scirc g \scirc \tau\inverse$
for every $g \in \rfs{LIP}(F)$.

See Corollary \ref{editor-d8.5}.

\end{small}

\newpage

\noindent
{\Large\bf Contents}

\contentsline {section}{\numberline {1}\rule{4pt}{0pt}Introduction}{6}
\contentsline {subsection}{\numberline {}
Plan of the introduction}{6}
\contentsline {subsection}{\numberline {1.1}
\rule{1pt}{-1.5pt}\kern0pt General description}{7}
\contentsline {subsection}{\numberline {1.2}
Faithfulness of classes of space-group pairs}{11}
\contentsline {subsection}{\numberline {1.3}
Moduli of continuity and groups of locally uniformly~continuous\\
\rule{4pt}{0pt}homeomorphisms}{14}
\contentsline {subsection}{\numberline {1.4}
Other groups of uniformly continuous homeomorphisms}{17}
\contentsline {subsection}{\numberline {1.5}
Groups of extendible homeomorphisms and the group of homeomorphisms\\
\rule{1pt}{0pt}\kern4pt of the closure of an open set}{20}
\contentsline {subsection}{\numberline {1.6}
Local uniform continuity at the boundary of an open set}{23}
\contentsline {subsection}{\numberline {1.7}
Further questions and discussion}{25}
\contentsline {subsection}{\numberline {1.8}
Some more facts about reconstruction theorems}{28}
\contentsline {section}{\numberline {2}
Obtaining a homeomorphism from a group isomorphism}{33}
\contentsline {subsection}{\numberline {2.1}
Capturing the action of the group on the regular open sets}{33}
\contentsline {subsection}{\numberline {2.2}
Faithfulness in locally compact spaces}{34}
\contentsline {subsection}{\numberline {2.3}
Faithfulness in normed and Banach spaces}{37}
\contentsline {subsection}{\numberline {2.4}
Faithfulness of normed manifolds}{57}
\contentsline {subsection}{\numberline {2.5}
The faithfulness of some smaller subgroups}{60}
\contentsline {section}{\numberline {3}
The local $\Gamma $-continuity of a conjugating homeomorphism}{62}
\contentsline {subsection}{\numberline {3.1}
General description}{62}
\contentsline {subsection}{\numberline {3.2}
Partial actions and decayability}{63}
\contentsline {subsection}{\numberline {3.3}
Translation-like partial actions}{77}
\contentsline {subsection}{\numberline {3.4}
Affine-like partial actions}{88}
\contentsline {subsection}{\numberline {3.5}
Summary and questions}{102}
\contentsline {subsection}{\numberline {3.6}
Normed manifolds}{105}
\contentsline {section}{\numberline {4}
The local uniform continuity of conjugating homeomorphisms}{108}
\contentsline {section}{\numberline {5}
Other groups defined by properties related to uniform continuity}{122}
\contentsline {subsection}{\numberline {5.1}
General description}{122}
\contentsline {subsection}{\numberline {5.2}
The group of uniformly continuous homeomorphisms}{126}
\contentsline {subsection}{\numberline {5.3}
The group of homeomorphisms which are uniformly continuous on\\
\rule{1pt}{0pt}\kern4pt every bounded set}{134}
\contentsline {subsection}{\numberline {5.4}
Groups of homeomorphisms which are uniformly continuous
on every\\
\rule{1pt}{0pt}\kern4pt bounded positively distanced set}
{139}
\contentsline {section}{\numberline {6}
Groups of extendible homeomorphisms and the reconstruction of\\
\rule{1pt}{0pt}\kern4pt the closure of open sets}{167}
\contentsline {subsection}{\numberline {6.1}
General description}{167}
\contentsline {subsection}{\numberline {6.2}
Groups of extendible homeomorphisms}{168}
\contentsline {subsection}{\numberline {6.3}
Completely locally uniformly continuous homeomorphism groups}{194}
\contentsline {subsection}{\numberline {6.4}
The reconstruction of
\hbox{ cl\kern 1.3pt}$(X)$ from
$H(\hbox {cl\kern 1.3pt}(X))$}{197}
\contentsline {subsection}{\numberline {6.5}
Generalizations to manifolds and to nearly open sets}{202}
\contentsline {section}{\numberline {7}
Groups which are not of the same type are not isomorphic}{206}
\contentsline {subsection}{\numberline {7.1}
The group \hbox{ UC\kern 1.3pt}$(X)$ revisited}{206}
\contentsline {subsection}{\numberline {7.2}
The non-existence of isomorphisms between groups of different types}
{221}
\contentsline {section}{\numberline {8}
The group of locally $\Gamma$-continuous homeomorphisms of the closure\\
\rule{1pt}{0pt}\kern4pt of an open set}{225}
\contentsline {subsection}{\numberline {8.1}
General description}{225}
\contentsline {subsection}{\numberline {8.2}
Statement of the main theorems and the plan of the proof}{226}
\contentsline {section}{\numberline {9}
The Uniform Continuity Constant}{232}
\contentsline {subsection}
{\numberline {9.1}Preliminary lemmas about the existence of certain
constants}{232}
\contentsline {subsection}{\numberline {9.2}The main construction}{242}
\contentsline {section}{\numberline {10}
\rule{1pt}{0pt}\kern0pt${\bf 1}$-dimensional
boundaries}{273}
\contentsline {section}{\numberline {11}
Extending the inducing homeomorphism to the boundary}{281}
\contentsline {section}{\numberline {12}
The complete $\Gamma$-bicontinuity of the inducing homeomorphism}{296}
\contentsline {subsection}{\numberline {12.1}
$\itGamma$-continuity in directions
parallel to the boundary of $X$}{296} 
\contentsline {subsection}{\numberline {12.2}
${\it \Gamma }$-continuity for submerged pairs
and the star operation}{304}
\contentsline {subsection}{\numberline {12.3}Final results}{317}
\newpage
\contentsline {section}{\numberline {}
References}{320}
\contentsline {section}{\numberline {}
Index of symbols}{324}
\contentsline {section}{\numberline {}
Index of notations}{325}
\contentsline {section}{\numberline {}
Index of definitions}{332}

\newpage
\baselineskip 18pt

\newpage
\section{\kern-6.7mm. Introduction}
\label{s1}

\noindent
\subsection*{Plan of the Introduction.}
\setcounter{subsection}{0}
{\bf Section \ref{ss1.1}}
\newline
The section starts with a statement of two theorems which exemplify
the type of results proved in this work.
The notions of a faithful class and of a determining category are then
introduced.
A class of topological spaces is said to be faithful
if its members are reconstructible from their homeomorphism groups.
Example \ref{nnr1.2} contains a short survey of older reconstructibily
theorems,
and Example \ref{nnr1.3} mentions several determining categories.
We then describe the precise forms of the theorems which will be
proved in this work. 
\newline
{\bf Section \ref{ss1.2}}
\newline
This section summarizes Chapter \ref{s2}. The theorems described in
\ref{ss1.2} have the form: If for $i = 1,2$, $G_i \leq H(X_i)$
and $\varphi$ is an isomorphism between $G_1$ and
$G_2$, then there is
a homeomorphism $\tau$ between $X_1$ and $X_2$ such that $\tau$
induces $\varphi$.
\newline
{\bf Section \ref{ss1.3}}
\newline
This section is a summary of Chapters \ref{s3} and \ref{s4}.
It starts with the definition of a modulus of continuity.
A modulus of continuity $\itGamma$ is a set of functions from
$[0,\infty)$ to $[0,\infty)$ which serves as a measure for the
continuity of a uniformly continuous function.
With $\itGamma$ one associates the group $H_{\itGamma}^{\srfs{LC}}(X)$
of locally $\itGamma$-bicontinuous homeomorphisms
of $X$.
The reconstruction result for groups of type
$H_{\itGamma}^{\srfs{LC}}(X)$ says that any isomorphism between
$H_{\itGamma}^{\srfs{LC}}(X)$ and $H_{\itGamma}^{\srfs{LC}}(Y)$
is induced by a locally $\itGamma$-bicontinuous homeomorphism between
$X$ and $Y$.
\newline
{\bf Section \ref{ss1.4}}
\newline
Section \ref{ss1.4} summarizes the reconstruction theorems for the
group $\rfs{UC}(X)$ of uniformly bicontinuous homeomorphisms of $X$.
These theorems appear in Chapter \ref{s5}.
\newline
{\bf Section \ref{ss1.5}}
\newline
The previous sections dealt mainly with spaces which are an open
subset of a normed vector space.
This section describes the reconstruction theorems for spaces which
are the closure of an open subset of a normed vector space.
These theorems appear in Chapter \ref{s6}.
Section \ref{ss1.5} also includes a survey of the results of
Chapter \ref{s7}.
\newline
{\bf Section \ref{ss1.6}}
\newline
Let $X$ be the closure of an open subset of a normed space.
Chapters \ref{s8} - \ref{s12} deal with the group
$H_{\itGamma}^{\srfs{LC}}(X)$ when $X$ is such a space.
Section \ref{ss1.6} describes the results obtained in these chapters.
\newline
{\bf Section \ref{ss1.7}}
\newline
This section contains a discussion and open problems.
\newline
{\bf Section \ref{ss1.8}}
\newline
This section contains a short historical survey.

\subsection{\kern-6mm. General description.}
\label{ss1.1}
\label{ss1.1-general-description}

This work concerns with groups of auto-homeomorphisms of open subsets
of normed vector spaces and of manifolds over normed vector spaces.
Mainly, we consider groups whose definition is based on the metric
of the normed space, for example, the group of all bilipschitz
auto-homeomorphisms of such a space.

Two types of results are proved.
The following statement is an example of the first type.

1.
Suppose that $X_1,X_2$ are open subsets of the Banach spaces spaces
$E_1$ and $E_2$ respectively.
For $i = 1,2$ let $G_i$ be a group of auto-homeomorphisms of $X_i$
such that every bilipschitz homeomorphism of $X_i$ belongs to $G_i$.
Suppose that $\varphi$ is a group isomorphism between $G_1$ and $G_2$.
Then there is a homeomorphism $\tau$ between $X_1$ and $X_2$ such that
for every $g \in G_1$, $\varphi(g) = \tau \scirc g \scirc \tau\inverse$.

\kern1.5mm

An example of the second type of results is as follows.

2.
$\rfs{BL}(E)$ denotes the group of all auto-homeomorphisms $f$
of a Banach space $E$ such that $f$ and $f\inverse$ are Lipschitz
on every bounded set, and $\rfs{BUC}(E)$ denotes the group of all
auto-homeomorphisms $f$ of $E$ such that $f$ and $f\inverse$
are uniformly continuous on every bounded set.
These groups determine the spaces they act upon in the
following sense.\vspace{-3.5mm}

\begin{enumerate}
\item[(a)] Suppose that $E_1$ and $E_2$ are Banach spaces,
and $\varphi$ is a group isomorphism between
$\rfs{BL}(E_1)$ and $\rfs{BL}(E_2)$.
Then
there is a unique
homeomorphism $\tau$ between $E_1$ and $E_2$ such that
for every $f\in \rfs{BL}(E_1)$,
$\varphi(f) = \tau \scirc f \scirc \tau^{-1}$.
Also, $\tau$ and $\tau\inverse$
are Lipschitz on every bounded set ($\tau$ is BL).

\item[(b)] The same holds for groups of the type $\rfs{BUC}(E)$.
That is, the statement obtained from (a) by replacing BL by BUC is true.

\item[(c)] For every $E_1$ and $E_2$, \,$\rfs{BL}(E_1)$
and $\rfs{BUC}(E_2)$
are not isomorphic.
\smallskip
\end{enumerate}

\noindent
{\bf Terminology }
The notation $\iso{f}{X}{Y}$ means that $f$ is a homeomorphism between
the topological spaces $X$ and $Y$.
That is, $f$ is bijective, and $f$ and $f\inverse$ are continuous.
		   \index{N@AAAA@@$\iso{f}{X}{Y}$.
		   This means $f$ is a homeomorphism between
		   $X$ and $Y$}
Let $H(X) = \setm{f}{\iso{f}{X}{X}}$.
		   \index{N@h00@@$H(X)$. The group of all auto-homeomorphisms
		   of $X$}
If $G,H$ are groups, then $\iso{\varphi}{G}{H}$ means that
$\varphi$ is an isomorphism between $G$ and $H$.
The ordered pair with elements $a$ and $b$ is denoted by $\pair{a}{b}$.

   \index{N@AAAA@@$\iso{\varphi}{G}{H}$. This means
   $\varphi$ is an isomorphism between $G$ and $H$}
   \index{N@AAAA@@$\pair{a}{b}$. Notation of an ordered pair}

\begin{defn}\label{d1.1}
\begin{rm}
\num{a}
A pair $\pair{X}{G}$ consisting of a topological space
$X$ and a group $G$ of auto-homeomorphisms of $X$ is called a
{\it space-group pair}.
   \index{space-group pair. $\pair{X}{G}$ is a space-group pair,
	  if $X$ is a topological space and $G \leq H(X)$}
Let $K$ be a class of space-group pairs.
$K$ is {\it faithful}, if for every
   \index{faithful class of space-group pairs}
$\pair{X_1}{G_1}, \pair{X_2}{G_2} \in K$
and $\iso{\varphi}{G_1}{G_2}$ there exists
$\iso{\tau}{X_1}{X_2}$ which induces $\gvp$. That is, for every
$f \in G_1$, $\varphi(f) = \tau \scirc f \scirc \tau\inverse$.

A class $K$ of topological spaces is {\it faithful},
   \index{faithful class of topological spaces}
if
$\setm{\pair{X}{H(X)}}{X \in K}$ is faithful.

   \index{restricted topological category}
\num{b}
A {\it restricted topological category} is a category $\itbfK$
whose objects are topological spaces,
in which every morphism between two objects $X$ and $Y$ of $\itbfK$
is a homeomorphism from $X$ onto $Y$,
and in which for every morphism $g$ of $\itbfK$, $g \inverse$
also belongs to $\itbfK$.
For every $X,Y \in \itbfK$ let $\rfs{Iso}_{\itbfK}(X,Y)$
denote the set of morphisms between $X$ and $Y$
and $\rfs{Aut}_{\itbfK}(X) = \rfs{Iso}_{\itbfK}(X,X)$.

   \index{determining category}
We say that $\itbfK$ is a {\it determining category} if for every
$X,Y \in \itbfK$ and a group isomorphism
$\iso{\varphi}{\rfs{Aut}_{\itbfK}(X)}{\rfs{Aut}_{\itbfK}(Y)}$
there is $\tau \in \rfs{Iso}_{\itbfK}(X,Y)$ such that
$\varphi(g) = \tau \scirc g \scirc \tau\inverse$ for every
$g \in \rfs{Aut}_{\itbfK}(X)$.

   \index{distinguishable categories}
Let $\itbfK,\itbfL$ be restricted topological categories.
$\itbfK,\itbfL$ are said to be {\it distinguishable}\break
if for every $X \in \itbfK$ and $Y \in \itbfL$:
if $\rfs{Aut}_{\itbfK}(X) \cong \rfs{Aut}_{\itbfL}(Y)$,
then
$$
X \in \itbfL \mbox{ and }
\rfs{Aut}_{\itbfK}(X) = \rfs{Aut}_{\itbfL}(X)
\mbox{\ \ or \ \ }
Y \in \itbfK \mbox{ and }
\rfs{Aut}_{\itbfL}(Y) = \rfs{Aut}_{\itbfK}(Y).
$$
\hfill\proofend
\end{rm}
\end{defn}

The above notions provide a convenient way for stating the second type
of results in this work.
However, we shall not use other notions or any techniques
from category theory.

Some faithful classes of topological spaces and some determining
categories are listed in the next two examples.
The lists are not exhaustive.
\begin{examples}\label{nnr1.2}
\begin{rm}
The following classes are faithful.

(a) The class of Euclidean manifolds.
This was proved by J. Whittaker \cite{W} (published 1963).

(b) The class of manifolds over the Hilbert cube.
This was proved by R. McCoy \cite{McC} (published 1972).

(c) The class Euclidean manifolds with boundary.
This was proved by M. Rubin \cite{Ru1} (published 1989).

(d) The class of all spaces $\pair{X}{\tau}$ such that:
\begin{itemize}
\addtolength{\parskip}{-11pt}
\addtolength{\itemsep}{06pt}
\item[(1)] $X$ is a polyhedron,
and $\tau$ is either the metric or the coherent topology of $X$,
\item[(2)] the simplicial complex defining $X$ does not have an infinite
increasing \kern3pt(with respect to inclusion) sequence of simplexes,
\item[(3)] for every $x \in X$, $\setm{h(x)}{h \in H(X)}$
has no isolated points.
\vspace{-05.7pt}
\end{itemize}
This was proved by M. Rubin \cite{Ru1}.

(e) The class of all manifolds over normed vector spaces.
This was proved by M. Rubin \cite{Ru1}.

(f) The class of manifolds over the class of real topological vector
spaces which are locally convex, normal
and have a nonempty open set which intersects every straight line
in a bounded set.
This was proved by A. Leiderman and M. Rubin \cite{LR} (published 1999).
\end{rm}
\end{examples}

\begin{examples}\label{nnr1.3}
\begin{rm}
The following are determining categories.

\num{a} For $n \leq \infty$ let $\itbfK^C_n$ be the category of
$C^k$-smooth manifolds. The morphisms of
$\itbfK^C_n$ are the homeomorphisms $f$ such that $f$ and $f\inverse$
are $k$ times continuously differentiable.
This was proved in \cite{Fi} (R. Filipkiewicz 1982),
but was earlier proved
by W. Ling in \cite{Lg1} and \cite{Lg2} (unpublished preprint, 1980).
See the topic ``Reconstruction questions for related groups''
in Subsection~\ref{ss1.7} of the Introduction.

\num{b} The categories arising from $C^k$-smooth Euclidean manifolds
carrying various types of additional structure. The morphisms being
the $C^k$-diffeomorphisms which preserve that structure.
These are determining categories.
This includes e.g.\ foliated manifolds (Ling \cite{Lg1} and \cite{Lg2})
and symplectic manifolds (Banyaga \cite{Ba1} 1997).
See the topic ``Reconstruction questions for related groups''
in Subsection~\ref{ss1.7} for more details.

\num{c} The category of open subsets of $\bbR^n$ with quasi-conformal
homeomorphisms as morphisms.
This was proved by V. Gol'dshtein and M. Rubin \cite{GR}
(1995).
\hfill\proofend
\end{rm}
\end{examples}

Continuing the investigaton of faithful classes and
determining categories,
we consider topological spaces with extra structure.
The spaces considered in this work are open subsets of a
normed vector space,
and more generally, manifolds over normed vector spaces.
We Also consider sets which are the closures of open subsets of
a normed space.

If $X$ is an open subset of a normed space $E$, the ``extra structure''
attached to $X$ is usually 
the object $\trpl{X}{\rfs{bd}^E(X)}{d}$, where $\rfs{bd}^E(X)$
is the boundary
of $X$ in $E$, and $d$ is the metric on $\rfs{cl}^E(X)$
inherited from $E$, ($\rfs{cl}^E(X)$ denotes the closure of $X$ in $E$).
The methods of this work can be applied
to more general ``extra structures''. See Remarks \ref{metr-bldr-r6.23}
and \ref{r4.6}.

This extra structure is used to define various subgroups of $H(X)$.
The groups $\rfs{BL}(X)$ and $\rfs{BUC}(X)$ defined at the beginning of
Subsection \ref{ss1.1} are examples of such subgroups.
Another typical example is as follows.
Let $X,Y$ be open subsets of the normed spaces
$E$ and $F$ respectively.
A homeomorphism $\iso{h}{X}{Y}$ is said to be {\it extendible}
if there is a continuous function
$\fnn{\barh}{\rfs{cl}(X)}{\rfs{cl}(Y)}$
such that $\barh$ extends $h$.
We consider the group
$\rfs{EXT}(X) \eqdf
\setm{h \in H(X)}{h \mbox{ and } h\inverse \mbox{ are extendible}}$.

   \index{extendible homeomorphism}

A homeomorphism $\iso{h}{X}{Y}$ is said to be
{\it completely locally uniformly continuous (CMP.LUC)}
if $h$ is extendible,
and for every $x \in \rfs{cl}(X)$
there is a neighborhood $U$ of $x$ in $\rfs{cl}(X)$
such that $h \nrestriction (U \cap X)$ is uniformly continuous.
We also consider the group
$$\rfs{CMP.LUC}(X) \eqdf
\setm{h \in H(X)}{h \mbox{ and } h\inverse \mbox{ are CMP.LUC}}.
$$

   \index{N@cmpluc@@$\rfs{CMP.LUC}(X)$. The group of biextendible
          homeomorphisms of $X$ which are bi-uniformly
	  \newline\indent
	  continuous at every $x \in \rfs{cl}(X)$}

The setting is thus as follows.
We shall have a class $\calM$ of topological spaces.
Usually this class consists of spaces $X$ such that
either $X$ is an open subset or the closure of an open subset
of normed vector space,
or even more generally,
$X$ can be the closure of an open subset of a manifold over a
normed vector space.
$\calP$ and $\calQ$ are properties
of maps between $X$ and $Y$ defined for objects of the form
$\trpl{X}{\rfs{bd}(X)}{d}$.
The set $\calP(X)$ of all homeomorphisms $f \in H(X)$
such that $f$ and $f\inverse$ have property $\calP$ is a subgroup of
$H(X)$, and the same holds for $\calQ(Y)$.
The final results have the following form.

If $X,Y \in \calM$ and $\iso{\varphi}{\calP(X)}{\calQ(Y)}$,
then
\begin{enumerate}
\addtolength{\parskip}{-11pt}
\addtolength{\itemsep}{06pt}
\item[(1)]
$\varphi$ is induced by a unique homeomorphism $\iso{\tau}{X}{Y}$,
\item[(2)]
$\calP(X) = \calQ(X)$
and $\tau$ and $\tau\inverse$ have property~$\calQ$,
or $\calP(Y) = \calQ(Y)$
and $\tau$ and $\tau\inverse$ have property~$\calP$.
\vspace{-05.7pt}
\end{enumerate}

   \index{N@k00@@$\itbfK_{\calM,\calP}$.
The category
$\pair{\calM}{\setm{g}{\iso{g}{X}{Y},\ X,Y \in \calM \mbox{ and }
g,g\inverse \mbox{ have property } \calP}}$
}

Let $\itbfK_{\calM,\calP}$ be the following category.

\begin{enumerate}
\addtolength{\parskip}{-11pt}
\addtolength{\itemsep}{06pt}
\item[{\rm (a)}]
The class of objects of $\itbfK_{\calM,\calP}$
is $\calM$.
\item[{\rm (b)}]
The class of morphisms of $\itbfK_{\calM,\calP}$
is
$\setm{\iso{g}{X}{Y}}{X,Y \in \calM \mbox{ and } g \mbox{ and }
g\inverse \mbox{ have}\break \mbox{property } \calP}$.
\vspace{-05.7pt}
\end{enumerate}

Conclusion (1)\,-\,(2) is the same as saying that
$\itbfK_{\calM,\calP}$ and $\itbfK_{\calM,\calQ}$
are determining categories and 
$\itbfK_{\calM,\calP}$ and $\itbfK_{\calM,\calQ}$ are distinguishable.

\kern2mm

This work uses only elementary facts. It is self-contained except for
Theorem \ref{t2.3} which is taken from \cite{Ru5};
it is stated there as Corollary 1.4
on page 122, and it is proved there in Corollary 2.10 on page 131.

Theorem \ref{t2.3} says that given a pair $(X,G)$, where $G$
is a subgroup of $H(X)$ satisfying certain weak transitivity
requirements, it is possible to recover from $G$ the Boolean algebra
$\rfs{Ro}(X)$ of regular open subsets of $X$, together with
the action of $G$ on $\rfs{Ro}(X)$. (A set $U$ is regular open
if $U$ is equal to the interior of its closure).
 
Consider the structures $(G,\, \rfs{Ro}(X);\,\gl^{\srfs{Ro}(X)}_G)$ and
$(G,\, X;\,\gl^X_G)$, where $\gl^{\srfs{Ro}(X)}_G$ and $\gl^X_G$ denote
the action of $G$ on $\rfs{Ro}(X)$ and on $X$ respectively.
The essence of Chapter 2 is showing that
for appropriate classes of $(X,\,G)$'s, \,
$(G,\, X;\,\gl^X_G)$ can be recoverd from
$(G,\, X;\,\gl^{\srfs{Ro}(X)}_G)$.
This kind of argument appears in Theorems
\ref{t2.2}, \ref{metr-bldr-t2.8}, \ref{t2.15} and \ref {t-bddlip-1.8}.

\noindent
\subsection{\kern-6mm. Faithfulness of classes of space-group pairs.}
\label{ss1.2}
\label{ss1.2-space-group-pairs}

Chapter \ref{s2} deals with the faithfulness of classes of
space-group pairs.
We introduce some terminology.

\begin{defn}\label{d1.2}
\begin{rm}
(a) A homeomorphism $h$ between two metric spaces $\rpair{X}{d^X}$
and $\rpair{Y}{d^Y}$ is {\it Lipschitz} if there is $K > 0$ such that
$d^Y(h(u),h(v)) \leq K d^X(u,v)$, for every $u,v \in X$.
We say that $h$ is {\it bilipschitz} if both $h$ and $h\inverse$ are
Lipschitz homeomorphisms. Define
$$\rfs{LIP}(X) \eqdf \setm{h \in H(X)}{h \mbox{ is bilipschitz}}.$$
   \index{lipschitz homeomorphism@@Lipschitz homeomorphism}
   \index{bilipschitz homeomorphism}
   \index{N@lip00@@$\rfs{LIP}(X)$. The group of bilipschitz
   auto-homeomorphisms of a metric space $X$}

(b) Let $X,Y$ be metric spaces. A homeomorphism $h$ between $X$ and $Y$
is {\it locally Lipschitz} if for every $u \in X$
there is a neighborhood $U$ of $u$ such that $h \nrestriction U$
is Lipschitz. $h$ is {\it locally bilipschitz}
if both $h$ and $h\inverse$ are locally Lipschitz.
Define
$$\rfs{LIP}^{\srfs{LC}}(X) \eqdf \setm{h \in H(X)}{h
\mbox{ is locally bilipschitz}}.$$

   \index{locally lipschitz homeomorphism@@locally
          Lipschitz homeomorphism}
   \index{locally bilipschitz homeomorphism} 
   \index{N@liplc00@@$\rfs{LIP}^{\srfs{LC}}(X)$.
          The group of locally bilipschitz auto-homeomorphisms of $X$}

   \index{N@lip01@@$\rfs{LIP}(X,S)$. For $S \subseteq X$,
          $\rfs{LIP}(X,S) =
          \setm{h \in \rfs{LIP}(X)}
	  {h \nrestriction (X - S) = \rfs{Id}}$}

(c) If $S \subseteq X$ is open, then
$$\rfs{LIP}(X,S) \eqdf
\setm{h \in \rfs{LIP}(X)}{h \nrestriction (X - S) = \rfs{Id}}.$$

   \index{N@lip02@@$\rfs{LIP}(X;F)$.
          For a normed space $E$,
	  $X \subseteq E$ and a dense linear subspace
          \newline\indent
          $F$ of $E$, \ 
          $\rfs{LIP}(X;F) =
	  \setm{h \in \rfs{LIP}(X)}{h(X \cap F) = X \cap F}$}

(d) Let $E$ be a normed vector space,
$F$ be dense linear subspace of $E$,
and $X$ be an open subset of $E$. Set
$$\rfs{LIP}(X;F) \eqdf
\setm{h \in \rfs{LIP}(X)}{h(X \cap F) = X \cap F}.$$

(e) For $E$, $F$, $X$, $S$ as above we define
   \index{N@lip03@@$\rfs{LIP}(X;S,F) =
	  \rfs{LIP}(X;F) \cap \rfs{LIP}(X,S)$}
$$\rfs{LIP}(X;S,F) \eqdf \rfs{LIP}(X;F) \cap \rfs{LIP}(X,S).$$

(f) $\rfs{LIP}^{\srfs{LC}}(X,S)$, $\rfs{LIP}^{\srfs{LC}}(X;F)$
and $\rfs{LIP}^{\srfs{LC}}(X;S,F)$ are defined analogously.

   \index{N@liplc01@@$\rfs{LIP}^{\srfs{LC}}(X,S)$. For $S \subseteq X$,
          $\rfs{LIP}(X,S)^{\srfs{LC}} =
          \setm{h \in \rfs{LIP}^{\srfs{LC}}(X)}{h \nrestriction (X - S)
	  = \rfs{Id}}$}
   \index{N@liplc02@@$\rfs{LIP}^{\srfs{LC}}(X;F)$.
          For a normed space $E$, $X \subseteq E$ and a dense linear
	  subspace
          \newline\indent
          $F$ of $E$, \ 
          $\rfs{LIP}^{\srfs{LC}}(X;F) =
          \setm{h \in \rfs{LIP}^{\srfs{LC}}(X)}{h(X \cap F) =
	  X \cap F}$}

   \index{N@liplc03@@$\rfs{LIP}^{\srfs{LC}}(X;S,F) =
          \rfs{LIP}(X;F) \cap \rfs{LIP}(X,S)$}

(g) Let $G \leq H$ mean that $G$ is a subgroup of $H$.
   \index{N@AAAA@@$G \leq H$. $G$ is a subgroup of $H$}

(h) For a normed vector space $E$, $x \in E$ and $r > 0$ let
$$B^E(x,r) = \setm{y \in E}{\norm{y - x} < r}.$$
Note that $\rfs{LIP}(X,S)$ and $\rfs{LIP}(X;F)$ are subgroups of $H(X)$.
\hfill\proofend
\end{rm}
\end{defn}

The main result of \,Chapter \ref{s2} is Part (c) of the next theorem.
It is restated as Theorem~\ref{t2.4}(b).
Parts (a) and (b) of Theorem~\ref{metr-bldr-t1.4}
are special cases of
Part (c).
They are more frequently used, and are more readable.

\begin{theorem}\label{metr-bldr-t1.4}
$(a)$ Let $K$ be the class of all pairs $\pair{X}{G}$ such that
$X$ is an open subset of some Banach space
and $\rfs{LIP}(X) \leq G \leq H(X)$.
Then $K$ is faithful.

$(b)$ Let $K$ be the class of all pairs $\pair{X}{G}$ such that
$X$ is an open subset of some normed vector space
and $\rfs{LIP}^{\srfs{LC}}(X) \leq G \leq H(X)$.
Then $K$ is faithful.

$(c)$ The class $K$ of all pairs $\pair{X}{G}$ which satisfy
\num{1} and \num{2}, \hbox{or \num{3} and \num{4}} below is faithful.
\begin{itemize}
\addtolength{\parskip}{-11pt}
\addtolength{\itemsep}{06pt}
\item[\num{1}] $X$ is an open subset of some Banach space $E$
and $G \leq H(X)$.
\item[\num{2}] For every $x \in X$ there are
an open set $S \subseteq X$ containing $x$
and a dense linear subspace $F \subseteq E$
such that $\rfs{LIP}(X;S,F) \leq G$.
\item[\num{3}] $X$ is an open subset of some normed vector space $E$
and $G \leq H(X)$.
\item[\num{4}] For every $x \in X$ there are
an open set $S \subseteq X$ containing $x$
and a dense linear subspace $F \subseteq E$
such that $\rfs{LIP}^{\srfs{LC}}(X;S,F) \leq G$.
\vspace{-05.7pt}
\end{itemize}
\end{theorem}

Compare Parts (a) and (b) of Theorem \ref{metr-bldr-t1.4}.
Part (a) deals with Banach spaces, and assumes that
$\rfs{LIP}(X) \leq G$.
Part (b) deals with normed spaces, but assumes that
$\rfs{LIP}^{\srfs{LC}}(X) \leq G$.
It is unknown whether in (b), assuming only that
$\rfs{LIP}(X) \leq G$ suffices.
The following theorem contains the strongest known fact regarding
this question.
It is restated as Corollary \ref{nc2.23}.

For a metric space $Z$, $x \in Z$ and $r > 0$
let $B^Z(x,r)$ denote the open ball in $Z$ determined by $x$ and $r$.
Let $X$ be an open subset of a normed space $E$.
Let $\overE$ denote the completion of $E$.
Define
$\overfs{int}(X) =
\bigcup \setm{B^{\overssE}(x,r)}{B^E(x,r) \subseteq X}$
and
$$
\rfs{IXT}(X) = \setm{h \nrestriction X}{h \in H(\overfs{int}(X))
\mbox{ and } h(X) = X}.
$$

\begin{theorem}\label{metr-bldr-t1.5}
Let $K$ be the class of all space-group pairs $\pair{X}{G}$
such that
\begin{itemize}
\addtolength{\parskip}{-11pt}
\addtolength{\itemsep}{06pt}
\item[\num{1}]
$X$ is an open subset of a Banach space,
or $X$ is an open subset
of a normed vector space which is a topological space
of the first category,
\item[\num{2}]
$\rfs{LIP}(X) \leq G \leq \rfs{IXT}(X)$.
\vspace{-05.7pt}
\end{itemize}
Then $K$ is faithful.
\end{theorem}

Theorem \ref{metr-bldr-t1.4} deals with open subsets of normed spaces.
However, the method of proof transfers without substantial change
to the more cumbersome setting of manifolds over
normed vector spaces ({\it normed manifolds}). 
This is dealt with in Theorem \ref{t2.15}.
In fact, Theorem \ref{t2.15} deals even with
normed manifolds with boundary and with spaces which are the closures of
open subsets of normed spaces. For such spaces
Theorem \ref{t2.15} says that the ``extended normed interior'' of the
space can be reconstructed from the group.
See Definition \ref{d2.14}.
An additional step is needed in order to recover the entire space.
This step is carried out under various assumptions in Theorems
\ref{t4.4}, \ref{metr-bldr-t6.20}, \ref{metr-bldr-t6.22},
\ref{metr-bldr-t6.25}(a) and \ref{metr-bldr-t6.28}.

For reasons of exposition and accessibility we include
in Chapter~\ref{s2} a theorem from \cite{Ru1}.
It says that $K_{\srfs{LCM}}$ is faithful, where $K_{\srfs{LCM}}$
is the class of all space-group pairs
$\pair{X}{G}$ which satisfy:\vspace{2pt}
\begin{itemize}
\addtolength{\parskip}{-11pt}
\addtolength{\itemsep}{06pt}
\item[(i)]
$X$ is a locally compact Hausdorff space without isolated points.
\item[(ii)]
$G$ has the property that for every nonempty open subset $U$ of $X$
and $x \in U$ the closure of the set
$\setm{g(x)}{g \in G \mbox{ and } g \nrestriction (X - U) = \rfs{Id}}$
has a nonempty interior.
\vspace{-05.7pt}
\end{itemize}

This result appears here as Theorem~\ref{t2.2}.

\vbox{
\noindent
\subsection{\kern-5.8mm. Moduli of continuity
and groups of locally uniformly continuous homeomorphisms.}
\label{ss1.3}
\label{ss1.3-LUC-homeomorphisms}
}

\noindent
Chapters \ref{s3}, \ref{s4} and \ref{s5} deal with groups consisting of
uniformly continuous homeomorphisms.
The uniform continuity of a function $f$ can be measured by a real
function which determines the bound of $d(f(x),f(y))$
as a function of $d(x,y)$.
Using semigroups of such real functions we obtain a hierarchy of
subgroups
of $H(X)$.

\begin{defn}\label{d1.5}
\begin{rm} $\rfs{MC}$ denotes the set of functions
   \index{N@mc@@$\rfs{MC} =
   \setm{h \in H([0,\infty))}{h \mbox{ is concave}}$}
$\alpha \in H([0,\infty))$
such that
for every\break
$x,y \in [0,\infty)$ and $0 \leq \lambda \leq 1$
$$
\alpha(\lambda x + (1 - \lambda)y) \geq
\lambda \alpha(x) + (1 - \lambda)\alpha(y).
$$
That is, $\rfs{MC}$ is the set of all concave homeomorphisms of
$[0,\infty)$.
\hfill\proofend
\end{rm}
\end{defn}

It is trivial that if $\alpha \in \rfs{MC}$,
then \,$\alpha(cx) \geq c \alpha(x)$ and $\alpha(dx) \leq d \alpha(x)$,
for every $0 \leq c\leq 1$ and $d \geq 1$.

\begin{defn}\label{d1.3}
\begin{rm} 
Let $f$ be a function from a metric space $(X,d^X)$ to a
metric space $(Y,d^Y)$.
Let $\alpha \in \rfs{MC}$. We say that $f$ is {\it $\alpha$-continuous}
   \index{continuous.
          $\alpha$-continuous. $f$ is $\alpha$-continuous,
          if for every $x,y$, \, $d(f(x),f(y)) \leq d(x,y)$}
if $d^Y(f(u),f(v)) \leq \alpha(d^X(u,v))$ 
for every $u,v \in X$. 

If $\fnn{f,g}{A}{\bbR \cup \sngltn{\infty}}$,
then $f \leq g$ means that $f(a)~\leq~g(a)$
for every $a \in A$.
   \index{N@AAAA@@$f \leq g$. This means for every $t$,
          $f(t) \leq g(t)$}

Let $\fnn{\alpha,\beta}{[0,\infty)}{\bbR \cup \sngltn{\infty}}$.
Then $\alpha \preceq \beta$ means that 
   \index{N@AAAA@@$\alpha \preceq \beta\ \ \equiv\ \ $ 
	  for some $t > 0$, \, $\alpha \nrestriction [0,t] \leq
	  \beta \nrestriction [0,t]$}
there is $a > 0$ such that
$\alpha \nrestriction [0,a] \leq \beta \nrestriction [0,a]$.

For $\itGamma \subseteq \rfs{MC}$ we define
$$\rfs{cl}_{\spreceq}(\itGamma) =
\setm{\alpha \in \rfs{MC}}{\mbox{ for some } \gamma \in \itGamma,\ \,
\alpha \preceq \gamma}.$$
\rule{0pt}{1pt}\index{N@cl00@@$\rfs{cl}_{\spreceq}(\itGamma) =
   \setm{\alpha \in \rfs{MC}}{\mbox{ for some } \gamma \in \itGamma,\ \,
   \alpha \preceq \gamma}$.}\kern-77pt
\hfill\proofendeol
\end{rm}
\end{defn}

Note that if $K > 0$, then the function $y = Kx$ belongs to \rfs{MC}. 
Also, if $\alpha,\beta \in \rfs{MC}$,
then $\alpha + \beta,\ \alpha \scirc \beta \in \rfs{MC}$.

\begin{defn}\label{nn1.3}
\begin{rm}
Let $\itGamma$ denote a subset of $\rfs{MC}$ containing 
$\rfs{Id}_{[0,\infty)}$.
We define the following properties of $\itGamma$.
\begin{itemize}
\addtolength{\parskip}{-11pt}
\addtolength{\itemsep}{06pt}
\item[M1]
For every $\alpha \in \itGamma$ and $\beta \in \rfs{MC}$: 
if $\beta \preceq \alpha$, then $\beta \in \itGamma$.
\item[M2]
For every $\alpha \in \itGamma$ and $K > 0$:  
$K \alpha,\ \alpha(Kx) \in \itGamma$.
\item[M3]
For every $\alpha,\beta \in \itGamma$:
$\alpha + \beta \in \itGamma$.
\item[M4] For every $\alpha,\beta \in \itGamma$:
$\alpha \scirc \beta \in \itGamma$.
\item[M5]
$\itGamma$ is {\it countably generated}. This means that there is a
countable set $\itGamma_0 \subseteq \itGamma$ such that 
$\itGamma \subseteq \rfs{cl}_{\spreceq}(\itGamma_0)$.
\item[M6]
Let $\alpha^{\sscirc n}$ denote $\alpha \scirc \ldots \scirc \alpha$,
$n$ times.
We say that $\itGamma$ is {\it principal}
if there is $\alpha \in \itGamma$ such that
$\itGamma \subseteq
\rfs{cl}_{\spreceq}(\setm{\alpha^{\sscirc n}}{n \in \bbN})$.
\vspace{-05.7pt}
\end{itemize}
   \index{countably generated. $\itGamma$ is countably generated,
          if for some
          countable $\itGamma_0 \subseteq \itGamma$,
	  \newline\indent\rule{0pt}{1pt}\kern00mm
	  $\itGamma \subseteq
	  \setm{\alpha \in \rfs{MC}}{(\exists \gamma \in \itGamma_0)
	  (\alpha \preceq \gamma)}$}
   \index{N@AAAA@@$g^{\sscirc n}$. Notation for
          $g \scirc \ldots \scirc g$ \ $n$ times}
   \index{prinicipal. $\itGamma$ is principal, if for some
          $\alpha \in \itGamma$,
          $\itGamma \subseteq
          \rfs{cl}_{\spreceq}(\setm{\alpha^{\sscirc n}}{n \in \bbN})$}
\end{rm}
\end{defn}

\begin{example}\label{e1.4}
\begin{rm}
(a) The set
$\itGamma^{\srfs{LIP}} \eqdf
\setm{\alpha \in \rfs{MC}}{\alpha \preceq Kx \mbox{ for some $K > 0$}}$
satis\-fies M1~-~M6, and it is called the {\it Lipschitz modulus}.
   \index{N@gammalip@@$\itGamma^{\srfs{LIP}} =
   \setm{\alpha \in \rfs{MC}}{\mbox{ for some $K > 0$, }
   \alpha \preceq Kx}$}

(b) For $0<r\leq 1$ the set
$\itGamma^{\srfs{HLD}}_r \eqdf
\setm{\alpha \in \rfs{MC}}
{\alpha \preceq K x^r \mbox{ for some $K > 0$}}$
is called the {\it $r$-H\"{o}lder set},
and it satisfies M1 - M3 and M5.
   \index{N@gammahld00@@$\itGamma^{\srfs{HLD}}_r =
   \setm{\alpha \in \rfs{MC}}{\mbox{ for some $K > 0$, }
   \alpha \preceq K x^r}$}

(c) The set $\itGamma^{\srfs{HLD}} \eqdf
\bigcup \setm{\itGamma^{\srfs{HLD}}_r}{r \in (0,1]}$
is called the {\it H\"{o}lder modulus},
and it satisifies M1 - M6.
   \index{N@gammahld01@@$\itGamma^{\srfs{HLD}} =
          \bigcup \setm{\itGamma^{\srfs{HLD}}_r}{r \in (0,1]}$.
          The H\"{o}lder modulus}
\end{rm}
\end{example}

\begin{prop}\label{p1.5}
\num{a} If $\itGamma \supseteq \itGamma^{\srfs{LIP}}$
and $\itGamma$ satisfies M1 and M4, then it satisfies M3.

\num{b} If $\itGamma$ satisfies M1 and M3, then it satisfies M2.
\end{prop}

\noindent
{\bf Proof }
Left to the reader.
\hfill\myqed

\begin{defn}\label{d1.6}
\begin{rm}
(a) Let $\itGamma\subseteq \rfs{MC}$ and $f$ be a
function from a metric space $X$ to a metric space $Y$.
Then $f$ is {\it locally $\itGamma$-continuous}
   \index{locally $\itGamma$-continuous}
if for every $x \in X$ there is a neighborhood $U$ of $x$
and $\alpha \in \itGamma$ such that $f \nrestriction U$
is $\alpha$-continuous.
$f$ is {\it locally $\itGamma$-bicontinuous},
if $f$ is a homeomorphism between
   \index{locally $\itGamma$-bicontinuous}
$X$ and $\rfs{Rng}(f)$, and both $f$ and $f\inverse$ are locally
$\itGamma$-continuous.

(b) Let $\itGamma \subseteq \rfs{MC}$. Then $\itGamma$ is called
a {\it modulus of continuity}
   \index{modulus of continuity}
if $\rfs{Id}_{[0,\infty)} \in \itGamma$
and $\itGamma$ satisfies\break
M1 - M4. Hence $\itGamma^{\srfs{LIP}} \subseteq \itGamma$.

(c) Let $\itGamma$ be a modulus of continuity,
and $X$ be a metric space.
	  \index{N@hlc00@@$H_{\itGamma}^{\srfs{LC}}(X)$.
	  The group of locally $\itGamma$-bicontinuous
	  auto-homeomorphisms of $X$}
$H_{\itGamma}^{\srfs{LC}}(X)$
denotes the set of locally $\itGamma$-bicontinuous homeomorphisms
from $X$ onto $X$.
\rule{1pt}{0pt}\hfill\proofend
\end{rm}
\end{defn}

Obviously, $\pair{H_{\itGamma}^{\srfs{LC}}(X)}{\scirc}$
is a group.\smallskip

Chapters \ref{s3} and \ref{s4} deal with groups of type
$H_{\itGamma}^{\srfs{LC}}(X)$.
The main result on such groups is stated in
Theorem \ref{metr-bldr-t4.1}(a),
and is proved at the end of Chapter~\ref{s4}.
The part of that theorem which deals with moduli of continuity
different from $\rfs{MC}$ appears in Corollary~\ref{metr-bldr-t3.42}(a).

The following theorem captures much of the contents of
\ref{metr-bldr-t4.1}(a). The full statement of \ref{metr-bldr-t4.1}(a)
requires more terminology.

\begin{theorem}\label{t1.7}\label{metr-bldr-t1.12}
For $\ell = 1,2$
let $\itGamma_{\ell}$ be a modulus of continuity such that
either $\itGamma_{\ell}$ is countably generated
or $\itGamma_{\ell} = \rfs{MC}$; let
$E_{\ell}$ be a normed space and
$X_{\ell}$ be a nonempty open subset of $E_{\ell}$.
Let
$\iso{\varphi}
{H_{\itGamma_{1}}^{\srfs{LC}}(X_1)}{H_{\itGamma_{2}}^{\srfs{LC}}(X_2)}$.
Then $\itGamma_1 = \itGamma_2$,
{\thickmuskip=2mu \medmuskip=1mu \thinmuskip=1mu 
and there is a locally $\itGamma_1$-bicontinuous
homeomor\-phism $\tau$
such that $\tau$ induces $\varphi$.
That is,
$\varphi(f) = \tau \scirc f \scirc \tau\inverse$
for every $f \in H_{\itGamma_{1}}^{\srfs{LC}}(X)$.
}
\end{theorem}

   \index{N@k01@@$\itbfK_{\itGamma}$}
Let $\itbfK_{\itGamma}$ denote the restricted topological category
in which the objects are open subsets of normed vector spaces,
and the morphisms are
locally $\itGamma$-bicontinuous homeomorphisms between such sets.
The above theorem says that for every $\itGamma$ as above
$\itbfK_{\itGamma}$ is a determining category,
that $\itbfK_{\itGamma_1}$ and $\itbfK_{\itGamma_2}$ are
distinguishable,
and that for every nonempty open subset of a normed vector space $X$
and distinct $\itGamma_{1}$ and $\itGamma_{2}$,
$H_{\itGamma_{1}}^{\srfs{LC}}(X) \neq H_{\itGamma_{2}}^{\srfs{LC}}(X)$.

The proof of \ref{t1.7} has two main steps.
In the first step we apply
Theorem \ref{metr-bldr-t1.4}
and deduce that there is $\iso{\tau}{X}{Y}$ such that $\tau$
induces $\varphi$.
This part of the argument 
is used repeatedly for the other groups which are dealt with in this
work.

The following statement constitutes the second step in the proof of
\ref{t1.7}.

\begin{theorem}\label{t1.12} Let $X$ and $\,Y$ be open subsets of the
normed spaces $E$ and $\,F$ respectively and $\iso{\tau}{X}{Y}$.
Let $\itGamma$ be a countably generated modulus of continuity.
If $\rfs{LIP}(X)^{\tau} \subseteq H_{\itGamma}^{\srfs{LC}}(Y)$,
then $\tau$ is locally $\itGamma$-bicontinuous.
\end{theorem}

The above theorem is restated as Theorem \ref{metr-bldr-t3.27}.

\begin{remark}\label{r1.8}
\begin{rm}
(a) Theorem \ref{t1.7} is stated only
for open subsets of normed spaces. But it is also true for normed
manifolds.
See Definitions \ref{d2.14} and \ref{metr-bldr-d3.46}
and Corollary \ref{metr-bldr-c3.48}(a).
In fact, if $\pair{X}{\itPhi}$ is a normed manifold with an atlas
$\itPhi$
such that for every $\varphi_1,\varphi_2 \in \itPhi$,
$\varphi_1 \scirc \varphi_2\inverse$ is locally $\itGamma$-continuous,
then $H_{\itGamma}^{\srfs{LC}}(X)$ can be defined,
and Theorem~\ref{t1.7} remains true. The proof remains 
essentially unchanged.

(b) Theorem \ref{t1.7} has the obvious shortcoming of assuming 
that $\itGamma$ is countably generated.
In fact, the assumption on $\itGamma$ in
Theorem \ref{metr-bldr-t4.1}(a) is weaker.
For example, for open subsets $X,Y \subseteq \ell_{\infty}$
the conclusion of\ Theorem~\ref{t1.7} is true for every
modulus of continuity.
Note though that the two natural moduli which motivated 
\ref{t1.7}, the Lipschitz and the H\"{o}lder moduli
are countably generated, and hence are covered by \ref{t1.7}.
But the question of whether Theorem \ref{t1.7} is true for every
modulus of continuity remains open.
\end{rm}
\end{remark}

\subsection{\kern-5.7mm. Other groups of uniformly continuous
homeomorphisms.}
\label{ss1.4}
\label{ss1.4-other-UC-groups}

A priori it seems natural to deal with the group
$\rfs{UC}(X)$ of all uniformly bicontinuous homeomorphisms of $X$
rather than with $H^{\srfs{LC}}_{\srfs{MC}}(X)$.
(A homeomorphism $h$ is {\it uniformly bicontinuous} if for every
$\varepsilon > 0$ there is $\delta > 0$
such that if $d(x,y) < \delta$, then
$d(h(x),h(y)) < \varepsilon$,
and if $d(h(x),h(y)) < \delta$,
then $d(x,y) < \varepsilon$).

Similarly, the group $H_{\itGamma}(X)$ of all $\itGamma$-bicontinuous
homeomorphisms of $X$ seems to be more natural than
$H^{\srfs{LC}}_{\itGamma}(X)$.
(A homeomorphism $h$ is {\it $\itGamma$-bicontinuous}, if there is
$\gamma \in \itGamma$ such that $h$ and $h\inverse$
are $\gamma$-continuous).
It turns out that $\rfs{UC}(X)$ and $H_{\itGamma}(X)$
pose more problems than their counterparts.
Chapter \ref{s5} addresses these groups and some related groups.

   \index{bicontinuous. $\itGamma$-bicontinuous.
   $h$ is $\itGamma$-bicontinuous,
   if $(\exists \gamma \in \itGamma)(
   h,\ h\inverse \mbox{ are } \gamma\mbox{-continuous})$}
   \index{N@h01@@$H_{\itGamma}(X)$. The group of all
   $\itGamma$-bicontinuous auto-hoeomorphisms of $X$}

Let $\calP$ be a property of maps and $X,Y$ be topological spaces.
Define
$\calP(X,Y) = \setm{h}{\iso{h}{X}{Y}
\mbox{ and $h$ has property }\calP}$.
If $H$ is a set of $\onetoonen$ functions,
then $H\inverse \eqdf \setm{h\inverse}{h \in H}$.
Define $\calP^{\pm}(X,Y) = \calP(X,Y) \cap (\calP(Y,X))\inverse$
and
$\calP(X) = \calP^{\pm}(X,X)$. We consider only $\calP$'s
such that $\calP(X)$ is a group.
The final results of Chapter \ref{s5} have the following form.
\begin{itemize}
\addtolength{\parskip}{-11pt}
\addtolength{\itemsep}{06pt}
\item[$(*)$] Suppose that $\iso{\varphi}{\calP(X)}{\calP(Y)}$.
Then there is $\tau \in \calP^{\pm}(X,Y)$ such that
$\tau$ induces $\varphi$.\kern-5pt
\vspace{-05.7pt}
\end{itemize}
{\thickmuskip=2mu \medmuskip=1mu \thinmuskip=1mu 
A class $\calM$ of topological spaces is called
{\it $\calP$\,-\,determined}, if $(*)$ holds for every $X,Y \in~K$,
}
that is, if the category $\itbfK_{\calM,\calP}$
whose objects are the members of $\calM$
and whose morphisms are the members of $\calP^{\pm}(X,Y)$
for $X,Y \in \calM$ is a determining category.
   \index{N@p00@@$\calP(X,Y) = \setm{h}{\iso{h}{X}{Y}
          \mbox{ and $h$ has property }\calP}$}
   \index{N@AAAA@@$H\inverse = \setm{h\inverse}{h \in H}$}
   \index{N@p01@@$\calP^{\pm}(X,Y) = \calP(X,Y) \cap
          (\calP(Y,X))\inverse$}
   \index{N@p02@@$\calP(X) = \calP^{\pm}(X,X)$}
   \index{determined class. $\calP$\,-\,determined class of
          topological spaces} 

The first result in Chapter \ref{s5} is about groups of type
$\rfs{UC}(X)$.
Denote the {\it diameter} of a subset $A$ of a metric space
by $\rfs{diam}(A)$.
A metric space $\pair{X}{d}$
is {\it uniformly\,-\,in\,-\,diameter arcwise\,-\,connected}
if for every $\varepsilon > 0$ there is $\delta > 0$
such that for every $x,y \in X$:
if $d(x,y) < \delta$, then there is an arc $L \subseteq X$ connecting
$x$ and $y$ such that $\rfs{diam}(L) < \varepsilon$.
The following statement is the main result on groups of type
$\rfs{UC}(X)$. It is restated as Corollary~\ref{metr-bldr-c5.7}.

\begin{theorem}\label{metr-bldr-t1.15}
Let $X$ be an open subset of a Banach space or of a normed vector space
of the first caregory. Suppose that the same holds for $Y$.
Suppose further that $X$ and $Y$ are uniformly\,-\,in\,-\,diameter
arcwise\,-\,connected.
Let $\iso{\varphi}{\rfs{UC}(X)}{\rfs{UC}(Y)}$.
Then there is
$\tau \in \rfs{UC}^{\pm}(X,Y)$ such that $\tau$
induces $\varphi$.
\end{theorem}

The following theorem restated later as \ref{t4.4} is a corollary of
\ref{metr-bldr-t1.15}.

\begin{theorem}\label{t1.12new}\label{metr-bldr-t1.16}
Let $F$ and $K$ be the closures of uniformly\,-in\,-\,diameter
arcwise\,-\,connected open bounded
subsets of $\kern1pt\bbR^m$ and $\kern1pt\bbR^n$ respectively. Let
$\iso{\varphi}{H(F)}{H(K)}$. Then $\varphi$
is induced by a homeomorphism between $F$ and~$K$.
\end{theorem}

Theorem \ref{metr-bldr-t1.16} is considerably stronger than the
analogous statement for Euclidean manifolds with boundary.
This is so, since
uniformly\,-\,in\,-\,diameter arcwise\,-\,connected open subsets
of $\bbR^n$
may have a boundary which is more complicated than the boundary of
a manifold with boundary.

$\rfs{UC}(X)$ is a special case of the groups $H_{\itGamma}(X)$.
But the analogue of Theorem \ref{metr-bldr-t1.15} is not true for
$H_{\itGamma}(X)$. In Example \ref{metr-bldr-e5.10}
it is shown that for every normed space $E$
there is
$\tau \in H(E)$ such that $(\rfs{LIP}(E))^{\tau} = \rfs{LIP}(E)$
but $\tau \not\in \rfs{LIP}(E)$.

Chapter \ref{s5} proves $\calP$\,-\,determined-ness for several other
$\calP$'s. Definition~\ref{metr-bldr-d5.5} lists eight types of groups
for which $\calP$\,-\,determined-ness can be proved.
But we have chosen to deal only with properties $\calP$
which occur in other mathematical contexts.

\begin{defn}\label{metr-bldr-d1.17}
\begin{rm}
(a) Let $\rfs{BUC}(X,Y)$ denote the set of homeomorphisms
$\iso{g}{X}{Y}$ such that $g$ takes bounded sets to bounded sets
and for every bounded $B \subseteq X$, $g \nrestriction B$ is
uniformly continuous.
   \index{N@buc00@@$\rfs{BUC}(X,Y) = \setm{g \in H(X,Y)}
   {g \mbox{ is boundedness preserving and } g \nrestriction A
   \mbox{ is UC for every}
   \newline\phantom{\hbox{$\rfs{BUC}(X,Y) =\ $}}\mbox{bounded set }
   A \subseteq X}$}

(b) Let $X$ be a metric space.
$X$ is
{\it boundedly uniformly\,-\,in\,-\,diameter arcwise\,-\,connected}
if for every bounded set $B \subseteq X$ and $\varepsilon > 0$
there is $\delta >0$ such that for every $x,y \in B$:
if $d(x,y) < \delta$, then there is an arc $L \subseteq X$ connecting
$x$ and $y$ such that $\rfs{diam}(L) < \varepsilon$.

(c) If $\fnn{h}{[0,1] \times X}{X}$ and $t_0 \in [0,1]$,
then the function $f$ from $X$ to $X$ defined by 
$f(x) = h(t_0,x)$ is denoted by $h_{t_0}$.
$X$ has {\it Property MV1} if for every bounded $B \subseteq X$
there are \hbox{$r = r_B > 0$}
and $\alpha = \alpha_B \in \rfs{MC}$
such that for every $x \in B$ and $0 < s \leq r$,
there is an $\alpha$-continuous function $\fnn{h}{[0,1] \times X}{X}$
such that:
for every $t \in [0,1]$, $h_t(x) \in H(X)$ and
$h_t\inverse$ is $\alpha$-continuous;
$h_0 = \rfs{Id}$ and $d(x,h_1(x)) = s$;
and
$h_t \nrestriction (X - B(x,2s)) = \rfs{Id}$
for every $t \in [0,1]$.
\hfill\proofend
\end{rm}
\end{defn}

The following $\calP$\,-\,determined-ness theorem is restated as
Theorem~\ref{metr-bldr-t5.19}.

\begin{theorem}\label{metr-bldr-t1.18}
Let $K$ be the class of all $X$ such that $X$ is an open subset of
a Banach space or $X$ is an open subset of a normed space of the first
category,
$X$ is boundedly uniformly\,-\,in\,-\,diameter arcwise\,-\,connected,
and $X$ has Property MV1.
\underline{Then} $K$ is BUC\,-\,determined.
\end{theorem}

There is of course the $\itGamma$ variant of $\rfs{BUC}(X)$.
For a modulus of continuity $\itGamma$ define
\vspace{1.5mm}
\newline
\renewcommand{\arraystretch}{1.5}
\addtolength{\arraycolsep}{-3pt}
\centerline{
$
\begin{array}{lll}
H^{\srfs{BD}}_{\itGamma}(X)
&
=
&
\setm{h \in H(X)}
{\mbox{for every bounded } A \subseteq X
\mbox{ there is }\gamma \in \itGamma \mbox{ such that}
\\
&&
h \nrestriction A \mbox{ is }
\gamma\mbox{-bicontinuous}}
\end{array}
$
}
\renewcommand{\arraystretch}{1.0}
\addtolength{\arraycolsep}{3pt}

\noindent
When $X$ is a subset of a finite-dimensional normed space
and $\itGamma$ is principal, then Theorem~\ref{ams-bddly-lip-bldr-t1.4}
provides a faithfulness result for this type of groups.

We do not know a more general theorem in this direction.

\kern2.0pt

The last type of groups considered in Chapter \ref{s5}
are groups of homeomorphisms $g$ such that $g \nrestriction B$
is uniformly continuous for every $B \subseteq X$
such that $B$ is bounded, and the distance of $B$ from the boundary
of $X$ is positive.
The $\calP$\,-\,determined-ness in this situation is proved
in Theorems \ref{metr-bldr-t5.30} and \ref{metr-bldr-t5.35}.

These theorems are not quoted here because their statement requires
terminology that has not yet been introduced.

Throughout Chapter \ref{s5} one encounters two types of intermediate
results.
\begin{itemize}
\addtolength{\parskip}{-11pt}
\addtolength{\itemsep}{06pt}
\item[(1)] 
Let $\iso{\tau}{X}{Y}$ be such that
$(\calP(X))^{\tau} = \calP(Y)$. Then $\tau \in \calP^{\pm}(X,Y)$.
\item[(2)] 
Let $\iso{\tau}{X}{Y}$ be such that
$(\calP(X))^{\tau} \subseteq \calP(Y)$.
Then $\tau \in \calP^{\pm}(X,Y)$.
\vspace{-05.7pt}
\end{itemize}

Results of type (2) are stronger, but they are not true for all
$\calP$'s which we consider.~Results~of type (2)
are needed in order to show that $\calP(X)$ cannot be
isomorphic to $\calQ(Y)$ when $\calP$ is different from
$\calQ$.

\noindent
\subsection{\kern-5.7mm. Groups of extendible homeomorphisms
and the group of\\
homeomorphisms of the closure of an open set.}
\label{ss1.5}
\label{ss1.5-extendible-homeomorphisms}

\kern1.8mm

Chapter \ref{s6} is concerned with the faithfulness of groups
of the form $H(\rfs{cl}(X))$ and with groups of the form $\rfs{EXT}(X)$,
where $X$ is an open subset of a normed vector space.
The group $\rfs{EXT}(X)$ is defined below.

Let $X,Y$ be open subsets of the normed spaces $E$ and $F$.
A continuous function\break
{\thickmuskip=2mu \medmuskip=1mu \thinmuskip=1mu 
$\fnn{g}{X}{Y}$ is called an {\it extendible function},
if there is a continuous function
$\fnn{\hatg}{\rfs{cl}(X)}{\rfs{cl}(Y)}$
}
such that $\hatg$ extends $g$.
The set of extendible homeomorphisms between $X$ and $Y$
is denoted by $\rfs{EXT}(X,Y)$.
Accordingly,
$\rfs{EXT}(X) = \setm{g \in H(X)}{g \mbox { and } g\inverse
\mbox{ are extendible}}$.
Note that if $X$ is a regular open subset of $\bbR^n$,
then $\rfs{EXT}(X) = H(\rfs{cl}(X))$.
Recall that a set is called {\it regular open}
if it is equal to the interior of
its closure.

   \index{N@ext00@@$\rfs{EXT}(X) =
          \setm{g \in H(X)}
	  {g \mbox { and } g\inverse \mbox{ are extendible}}$}

The goal is to find large classes $K$ of open subsets of a normed
space containing\break
the commonly encountered open sets
and containing also exotic open sets for which\break
$\setm{\rfs{cl}(X)}{X \in K}$ is faithful.
It is not true, though, that for every open subsets of
$X,Y \subseteq \bbR^n$, if
$\iso{\varphi}{H(\rfs{cl}(X))}{H(\rfs{cl}(Y))}$, then there is
$\iso{\tau}{\rfs{cl}(X)}{\rfs{cl}(Y)}$ such that $\tau$ induces
$\varphi$.
Example \ref{metr-bldr-e5.9} demonstrates this phenomenon
in two different ways.

The following theorem gives the flavor of the type of results proved in
Chapter \ref{s6}.

\begin{theorem}\label{metr-bldr-t1.19}
Let $X,Y$ be open bounded subsets of the Banach spaces $E$ and $F$.
Assume that:
\begin{itemize}
\addtolength{\parskip}{-11pt}
\addtolength{\itemsep}{06pt}
\item[\num{1}] There is $d$ such that for every $u,v \in X$ there is
a rectifiable arc $L \subseteq X$ connecting $u$ and $v$ such that
$\rfs{length}(L) \leq d$.
\item[\num{2}] For every point $w$ in the boundary of $X$ and for every
$\varepsilon > 0$
there is $\delta > 0$ such that for every $u,v \in X$:
if $\norm{u - w},\norm{v - w} < \delta$,
then there is an arc $L \subseteq X$ connecting $u$ and $v$
such that $\rfs{diam}(L) < \varepsilon$.
\item[\num{3}] Conditions \num{1} and \num{2} hold for $Y$.
\vspace{-05.7pt}
\end{itemize}
\underline{Then}
\begin{itemize}
\addtolength{\parskip}{-11pt}
\addtolength{\itemsep}{06pt}
\item[\num{a}] If $\iso{\varphi}{H(\rfs{cl}(X))}{H(\rfs{cl}(Y))}$,
then there is $\iso{\tau}{\rfs{cl}(X)}{\rfs{cl}(Y)}$
such that $\tau$ induces $\varphi$.
\item[\num{b}] If $\iso{\varphi}{\rfs{EXT}(X)}{\rfs{EXT}(Y)}$,
then there is $\tau \in \rfs{EXT}^{\pm}(X,Y)$
such that $\tau$ induces $\varphi$.
\vspace{-05.7pt}
\end{itemize}
\end{theorem}

Part (a) of the above theorem is an excerpt from
Theorem \ref{metr-bldr-t6.20},
and (b) is an excerpt from Theorem \ref{metr-bldr-t6.5}(a).

The class of spaces defined in Theorem \ref{metr-bldr-t1.19}
contains some spaces whose boundary is quite complicated.
Also, such spaces may have boundary points which are fixed under
$H(\rfs{cl}(X))$. Here is an example of a possibly not well-behaved set
which is covered by Theorem \ref{metr-bldr-t6.20}.

\begin{example}\label{bddly-lip-bldr-e1.19}\label{metr-bldr-e1.20}
\begin{rm}
Let $B$ and $S$ be the open unit ball and the unit sphere
in a Banach space $E$, 
and $\setm{\overB_i}{i \in I}$ be a family
of pairwise disjoint closed balls such that $\overB_i \subseteq B$
for every $i \in I$.
Suppose that for every $x \in E$: if every neighborhood
of $x$ intersects infinitely many $B_i$'s, then $x \in S$.
Then the set
$X \eqdf B - \bigcup_{i \in I} \overB_i$,
satisfies Clauses~(1)~and~(2) of Theorem~\ref{metr-bldr-t1.19}.
Note that even in the case of $E = \bbR^n$, the boundary of
$X$ can be complicated.
\hfill\proofend
\end{rm}
\end{example}

\kern1mm

Clause (2) in Theorem \ref{metr-bldr-t1.19} implies that
$\rfs{cl}(X)$ is arcwise connected.
Consider the open set $X$ described in the following example.
Its closure is not locally arcwise connected.

\begin{example}\label{metr-bldr-e1.21}
\begin{rm}
Let
$X = \setm{(r,\theta)}{\theta \in (\pi,\infty) \mbox{ and }
1 - \frac{1}{\theta - \pi/2} < r <
1 - \frac{1}{\theta + \pi/2}}$.
($X$ is described in polar coordinates).
Note that $X$ is an open spiral strip converging to the circle
$S(0,1)$.
\rule{0pt}{0pt}\hfill\proofend
\end{rm}
\end{example}

Example \ref{metr-bldr-e1.21} is not covered by
Theorem \ref{metr-bldr-t1.19} but it is included in the class
considered in the following theorem.

\begin{theorem}\label{metr-bldr-t1.22}
Let $X,Y$ be open bounded subsets of the normed spaces $E$ and $F$.
Assume that:
\begin{itemize}
\addtolength{\parskip}{-11pt}
\addtolength{\itemsep}{06pt}
\item[\num{1}]
For every sequence $\vecx = \setm{x_n}{n \in \bbN} \subseteq X$
there are a subsequence $\vecy$ of $\vecx$,
a sequence $\vecz$ such that $\vecz$ is convergent in $E$
and a sequence of rectifiable arcs $L_n \subseteq X$,\ %
$n \in \bbN$, such that $\sup_{n \in \bbN} \rfs{length}(L_n) < \infty$
and $L_n$ connects $y_n$ and $z_n$.
\item[\num{2}] For every $x \in \rfs{bd}(X)$ and $r > 0$
there is a continuous function
$\fnn{h_t(x)}{[0,1] \times \rfs{cl}(X)}{\rfs{cl}(X)}$
such that $h_0 = \rfs{Id}$, $h_1(x) \neq x$,
and for every $t \in [0,1]$,
$h_t \nrestriction X \in \rfs{EXT}(X)$ and
$h_t \nrestriction (\rfs{cl}(X) - B(x,r)) = \rfs{Id}$.
\item[\num{3}] Conditions \num{1} and \num{2} hold for $Y$.
\vspace{-05.7pt}
\end{itemize}
\underline{Then}
\begin{itemize}
\addtolength{\parskip}{-11pt}
\addtolength{\itemsep}{06pt}
\item[\num{a}] If $\iso{\varphi}{H(\rfs{cl}(X))}{H(\rfs{cl}(Y))}$,
then there is $\iso{\tau}{\rfs{cl}(X)}{\rfs{cl}(Y)}$
such that $\tau$ induces $\varphi$.
\item[\num{b}] If $\iso{\varphi}{\rfs{EXT}(X)}{\rfs{EXT}(Y)}$,
then there is $\tau \in \rfs{EXT}^{\pm}(X,Y)$
such that $\tau$ induces $\varphi$.
\vspace{-05.7pt}
\end{itemize}
\end{theorem}

Theorem \ref{metr-bldr-t1.22}(a) is an excerpt from
Theorem \ref{metr-bldr-t6.22},
and \ref{metr-bldr-t1.22}(b) is an excerpt from \ref{metr-bldr-t6.15}.
Example \ref{metr-bldr-e1.21} is restated as \ref{metr-bldr-e6.12}(a).
Other examples which are covered by Theorems \ref{metr-bldr-t6.22}
and \ref{metr-bldr-t6.15},
but have a non-locally arcwise connected closure appear in
\ref{metr-bldr-e6.7} and \ref{metr-bldr-e6.12}(b).
Another EXT\,-\,determined class is described in
Theorem \ref{metr-bldr-t6.10}.

Chapter \ref{s6} also deals with groups of type $\rfs{CMP.LUC}(X)$
defined in Subsection \ref{ss1.1-general-description}.
CMP.LUC\,-\,determined-ness is proved in
Theorem \ref{metr-bldr-t6.17}(a).
It completes the picture given in \hbox{Chapters \ref{s8}\,-\ref{s12}.}
The following is a special case of \ref{metr-bldr-t6.17}(a).

\begin{theorem}\label{metr-bldr-t1.23}
Let $X,Y$ be open bounded subsets of the normed spaces $E$ and $F$.
Assume that:
\begin{itemize}
\addtolength{\parskip}{-11pt}
\addtolength{\itemsep}{06pt}
\item[\num{1}]
For every sequence $\vecx = \setm{x_n}{n \in \bbN} \subseteq X$
there are a subsequence $\vecy$ of $\vecx$,
a sequence $\vecz$ such that $\vecz$ is convergent in $E$
and a sequence of rectifiable arcs $L_n \subseteq X$,\ %
$n \in \bbN$, such that $\sup_{n \in \bbN} \rfs{length}(L_n) < \infty$
and $L_n$ connects $y_n$ and $z_n$.
\item[\num{2}]
For every $x \in \rfs{bd}(X)$ there is $r > 0$
such that for every $\varepsilon > 0$ there is $\delta > 0$ such that
for every $u,v \in B^E(x,r) \cap X$: if $d(u,v) < \delta$,
then there is an arc $L \subseteq X$ connecting $u$ and $v$
such that $\rfs{diam}(L) < \varepsilon$.
\item[\num{3}] Conditions \num{1} and \num{2} hold for $Y$.
\vspace{-05.7pt}
\end{itemize}
\underline{Then} if $\iso{\varphi}{\rfs{CMP.LUC}(X)}{\rfs{CMP.LUC}(Y)}$,
then there is $\tau \in \rfs{CMP.LUC}^{\pm}(X,Y)$ such that
$\tau$ induces $\varphi$.
\end{theorem}

Two extensions of the results of Chapter \ref{s6} are presented at the
end of that chapter. These extensions cover some natural spaces
which are not covered by the original classes.
Also, the faithful class dealt with in Extension 2
contains $2^{2^{\aleph_0}}$ subsets of $\bbR^3$.

(1) The original classes considered in Chapter \ref{s6}
consist of open subsets of normed vector spaces, and the closures of
such sets. However, all the results obtained for these classes
translate to the class of open subsets of manifolds over
normed vector spaces and the closures of such sets.
See Example \ref{metr-bldr-e6.26} and Theorem \ref{metr-bldr-t6.28}.

(2) The results obtained for the class of closures of open subsets of
a normed vector space extend to the class of all subsets $Z$ of a normed
vector space which satisfy\break
$Z \subseteq \rfs{cl}(\rfs{int}(Z))$.
See Example \ref{metr-bldr-e6.24} and Theorem \ref{metr-bldr-t6.25}.

\kern1mm

{\thickmuskip=2mu \medmuskip=1mu \thinmuskip=1mu 
\hbox{Chapter \ref{s7} contains theorems of the following type.
Suppose that
$\iso{\varphi}{\calP(X)}{\calQ(Y)}$.
Then
}
}
\noindent
(i) There is $\iso{\tau}{X}{Y}$ such that $\tau$ induces $\varphi$.
\newline
(ii)
$\calP(X) = \calQ(X)$ and
$\tau \in \calQ^{\pm}(X,Y)$,
or $\calP(Y) = \calQ(Y)$ and
$\tau \in \calP^{\pm}(X,Y)$.

These results appear in Corollary \ref{metr-bldr-c7.11}.
As an example of such results we quote \ref{metr-bldr-c7.11}(e).

\begin{theorem}\label{metr-bldr-t1.24}
If $X$ and $Y$ are nonempty open subsets of an infinite-dimensional
Banach space,
then $\rfs{UC}(X) \not\cong \rfs{EXT}(Y)$.
\end{theorem}

\subsection{\kern-5.7mm. Local uniform continuity at the
boundary of an open set.}
\label{ss1.6}
\label{ss1.6-LUC-at-the-boundary}

Let $X \subseteq \bbR^n$ and $Y \subseteq \bbR^m$ be open sets
and suppose that
$\iso{\varphi}{\rfs{LIP}(\rfs{cl}(X))}{\rfs{LIP}(\rfs{cl}(Y))}$.
Can we conclude that there is
$\iso{\tau}{\rfs{cl}(X)}{\rfs{cl}(Y)}$ such that $\tau$ is bilipschitz
and $\tau$ induces $\varphi$?
This question motivates the work presented in
Chapters \ref{s8}\,-\,\ref{s12}.
Indeed, if the boundaries of $X$ and $Y$ are well-behaved,
then the answer to the above question is positive.

Let $X,Y$ be open subsets of the normed spaces $E$ and $F$,
and $\itGamma$ be a modulus of continuity.
For $g \in \rfs{EXT}(X,Y)$ let
$g^{\srfs{cl}}$ denote the continuous extension
of $g$ to $\rfs{cl}(X)$.
Define
$$
H^{\srfs{CMP.LC}}_{\itGamma}(X,Y) =
\setm{g \in \rfs{EXT}(X,Y)}{g^{\srfs{cl}}
\mbox{ is locally $\itGamma$-continuous}}
$$
and $H^{\srfs{CMP.LC}}_{\itGamma}(X) =
(H^{\srfs{CMP.LC}}_{\itGamma})^{\pm}(X,X)$.

   \index{N@hcmplc@@$H^{\srfs{CMP.LC}}_{\itGamma}(X,Y) =
          \setm{g \in \rfs{EXT}(X,Y)}{g^{\srfs{cl}}
          \mbox{ is locally $\itGamma$-continuous}}$}

Note that the group $\rfs{CMP.LUC}(X)$ discussed
in Subsection \ref{ss1.5} is a special case of groups of the form
$H^{\srfs{CMP.LC}}_{\itGamma}(X)$.
Indeed, $\rfs{CMP.LUC}(X) = H^{\srfs{CMP.LC}}_{\srfs{MC}}(X)$.
In the special case that
$X \subseteq \bbR^n$ is a regular open bounded set we have
$\rfs{LIP}(\rfs{cl}(X)) =
H^{\srfs{CMP.LC}}_{\itGamma^{\srs{LIP}}}(X)$.
More generally,
$H_{\itGamma}(\rfs{cl}(X)) =
H^{\srfs{CMP.LC}}_{\itGamma}(X)$.
So a determining-ness result for the property
$\calP = \hbox{\rm CMP.LC$_{\itGamma^{\srfs{LIP}}}$}$
implies such a result for the class $\itbfK_{\calM,\calP}$,
where $\calP = {\rm LIP}$
and $\calM$ is the class of bounded regular open subsets of
finite-dimensional spaces.
\smallskip

Chapters \ref{s8}\,-\,\ref{s12} are devoted to the
proof of the following statement about 
$H^{\srfs{CMP.LC}}_{\itGamma}(X)$.
\begin{description}
\item[$(*)$]
If $\iso{\varphi}{H^{\srfs{CMP.LC}}_{\itGamma}(X)}
{H^{\srfs{CMP.LC}}_{\itDelta}(Y)}$,
then $\itGamma = \itDelta$, and there is
$\tau \in (H^{\srfs{CMP.LC}}_{\itGamma})^{\pm}(X,Y)$\break
\rule{0pt}{0pt}\kern-2.8mm such that $\tau$ induces $\varphi$.
\end{description}

Statement $(*)$ is proved for $X,Y$, $\itGamma$ and $\itDelta$
which satisfy the following assumptions.
\begin{itemize}
\addtolength{\parskip}{-11pt}
\addtolength{\itemsep}{06pt}
\item[(1)] 
$\itGamma$ is principal, (see M6 in Definition \ref{nn1.3}).
\item[(2)] 
$X$ is locally $\itGamma$-LIN-bordered,
and $Y$ is locally $\itDelta$-LIN-bordered,
(see Definition \ref{editor-d8.1}(b)).
\vspace{-05.7pt}
\end{itemize}
The exact definition of local LIN-borderedness is a bit long,
but a main special case
is the class open sets whose closure is a manifold with boundary with
a $\itGamma$-bicontinuous atlas.

Statement $(*)$ is restated in Theorem
\ref{ams-bddly-lip-bldr-t1.4}(a).
The proof of \ref{ams-bddly-lip-bldr-t1.4}(a) has four steps.
The two major steps are Steps 3 and 4,
which are stated as Theorems \ref{t-bddlip-1.8}
and \ref{ams-bddly-lip-bldr-t5.22}.
The following theorem is the conclusion of the first
three steps combined together.
The prinicipality of $\itGamma$ is not needed here.
It is needed only at Step 4.

\begin{theorem}\label{metr-bldr-t1.25}
Let $\itGamma,\itDelta$ be countably generated moduli of continuity,
$E$ and $F$ be normed spaces and
$X \subseteq E$, $Y \subseteq F$ be open.
Suppose that $X$ is locally $\itGamma$-LIN-bordered,
and $Y$ is locally $\itDelta$-LIN-bordered.
Let
{\thickmuskip=2mu \medmuskip=1mu \thinmuskip=1mu 
$\iso{\varphi}{H^{\srfs{CMP.LC}}_{\itGamma}(X)}
{H^{\srfs{CMP.LC}}_{\itDelta}(Y)}$.
Then there is
$\tau \in \rfs{EXT}^{\pm}(X,Y)$
}
such that $\tau$ induces $\varphi$.
\end{theorem}

The proof of Theorem~\ref{metr-bldr-t1.25} requires much technical
work. This work is carried out in
Chapters \ref{s9} and \ref{s10}.
The proof of \ref{metr-bldr-t1.25}
appears at the end of Chapter \ref{s11}.

Step 4 of the proof of Theorem~\ref{ams-bddly-lip-bldr-t1.4}(a)
says that if in Theorem~\ref{metr-bldr-t1.25},
$\itGamma$ is principal,
then the homeomorphism $\tau$ obtained in \ref{metr-bldr-t1.25}
belongs to
$(H^{\srfs{CMP.LC}}_{\itGamma})^{\pm}(X,Y)$.

It should be pointed out that the results mentioned above are true
for open subsets of normed manifolds.
The final result for manifolds is stated in
Theorem~\ref{ams-bddly-lip-bldr-t1.4}(b).

\kern1.5mm

As a byproduct of the proof of the main theorem of
Chapters \ref{s8}\,-\,\ref{s12}, we also obtain a determining-ness
result for the group defined below.
Let $X$ be an open subset of a normed space $E$. Define
\vspace{1mm}
\newline
\indent
$H^{\srfs{BDR.LC}}_{\itGamma}(X) =
\setm{g \in \rfs{EXT}(X)}
{\mbox{every } x \in \rfs{cl}(X) - X
\mbox{ has a neighborhood } U \mbox{ in } \rfs{cl}(X)\\
\rule{35mm}{0pt}
\mbox{ such that }
g^{\srfs{cl}} \nrestriction U
\mbox{ is } \itGamma\mbox{-bicontinuous}}$.
\newline
Theorem \ref{ams-bddly-lip-bldr-t6.1}(b) contains
a determining-ness result for the property
$\calP = {\rm BDR.LC}_{\itGamma}$.

\subsection{\kern-5.7mm. Further questions and discussion.}
\label{ss1.7}
\label{ss1.7-questions-discussion}

This work leaves many unsolved questions,
which we mention at the point where they naturally arise.
In what follows we highlight the questions we regard to be more
central.

\kern2.0mm

\noindent
{\bf The countable generatedness of $\itGamma$.}

\begin{question}\label{metr-bldr-q1.26}
\begin{rm}
Can Theorem \ref{metr-bldr-t1.12} be proved for every pair of
moduli of continuity, regardless of whether they are 
countably generated or not?
That is, we ask if the following statement true?

For $\ell = 1,2$
let $\itGamma_{\ell}$ be a modulus of continuity.
Let
$E_{\ell}$ be a normed space and
$X_{\ell}$ be an open subset of $E_{\ell}$.
Let
$\iso{\varphi}
{H_{\itGamma_{1}}^{\srfs{LC}}(X_1)}{H_{\itGamma_{2}}^{\srfs{LC}}(X_2)}$.
\hbox{Then $\itGamma_1 = \itGamma_2$,}
and there is a locally $\itGamma_1$-bicontinuous homeomor\-phism $\tau$
such that $\tau$ induces~$\varphi$.
\hfill\proofend
\end{rm}
\end{question}

Note that the assumption in Theorem \ref{metr-bldr-t4.1}
is in fact somewhat weaker than countable generatedness.
We ask Question \ref{metr-bldr-q1.26} also for the other
theorems in which $\itGamma$ is required to be countably generated.
See e.g.\ Parts (a) and (b) of Theorem \ref{metr-bldr-t5.23}.

\kern2.0mm

\noindent
{\bf The principality of $\itGamma$ in the theorem about
$H^{\srfs{CMP.LC}}_{\itGamma}(X)$.}

\begin{question}\label{metr-bldr-q1.27}
\begin{rm}
Is Theorem \ref{ams-bddly-lip-bldr-t6.1}(a) true without the assumption
that $\itGamma$ is principal?
That is, we ask if the following statement true?
\newline
Let $X,Y$ be open subsets of a normed space,
and $\itGamma,\itDelta$ be moduli of continuity.\break
Assume that
$X$ is locally $\itGamma$-LIN-bordered,
and $Y$ is locally $\itDelta$-LIN-bordered.
If\break
$\iso{\varphi}{H^{\srfs{CMP.LC}}_{\itGamma}(X)}
{H^{\srfs{CMP.LC}}_{\itDelta}(Y)}$,
then
$\itGamma = \itDelta$, and there is
$\tau \in (H^{\srfs{CMP.LC}}_{\itGamma})^{\pm}(X,Y)$ such that $\tau$
induces $\varphi$.
\hfill\proofend
\end{rm}
\end{question}

Obviously, the case that $\itGamma$ and $\itDelta$ are
countably generated is also
unknown.

\kern2.0mm

\noindent
{\bf A possible stronger way of distinguishing between the
$H^{\srfs{LC}}_{\itGamma}(X)$'s.}

The fact that
$H^{\srfs{LC}}_{\itGamma}(X) \not\cong H^{\srfs{LC}}_{\itDelta}(Y)$
for $\itGamma \neq \itDelta$ may have a stronger reason.
That is, maybe there is a locally $\itDelta$-bicontinuous
homeomorphism which is not conjugate to any
locally $\itGamma$-bicontinuous homeomorphism.
So a positive answer to the following question together with
the faithfulness result of Theorem \ref{metr-bldr-t1.4}(a) will imply
the distinguishability of the $\itbfK_{\itGamma}$'s.

\begin{question}\label{metr-bldr-q1.29.1}
\begin{rm}
Let $\itGamma,\itDelta$ be moduli of continuity such that
$\itDelta \not\subseteq \itGamma$ and let $X$ be a nonempty open subset
of a normed space of dimension $> 1$.
Is there a locally $\itDelta$-bicontinuous homeomorphism $g$ of $X$
such that $g$ is not conjugate to any
$\itGamma$-bicontinuous homeomorphism?
\hfill\proofend
\end{rm}
\end{question}

In the space $\bbR$,
every homeomorphism is conjugate to a Lipschitz homemorphism.

\kern2.0mm

\noindent
{\bf Relaxing the assumption on the boundary
in the theorem about $H^{\srfs{CMP.LC}}_{\itGamma}(X)$.}\break
\indent
Let $X_0 =
\setm{(x,y) \in \bbR^2}{x > 0, -x^2 < y < x^2}$.
The set $X_0$ is not $\itGamma^{\srfs{LIP}}$-LIN-bordered.
Our general question is whether Theorem \ref{ams-bddly-lip-bldr-t6.1}(a)
can be strengthened to classes which include sets similar to $X_0$.
We may ask the following concrete question.

\begin{question}\label{metr-bldr-q1.28}
\begin{rm}
Let
$\varphi \in \rfs{Aut}(H^{\srfs{CMP.LC}}_{\itGamma^{\srs{LIP}}}(X_0))$.
Is $\varphi$ an inner automorphism?
\rule{1pt}{0pt}\hfill\proofend
\end{rm}
\end{question}

Question~\ref{editor-q8.10}
introduces the notion of a locally $\itGamma$-almost-linearly-bor\-dered
set (locally $\itGamma$-ALIN-bordered set).
It seems that Theorem \ref{ams-bddly-lip-bldr-t6.1}(a)
can be extended to the class of locally $\itGamma$-ALIN-bordered sets.
This requires a more detailed technical analysis similar to
the work carried out in Chap\-ters~\ref{s9}\,-\,\ref{s11}.

However, we do not know how to handle the type of singularity
at the boundary point $(0,0)$ of $X_0$ above.

\kern2.0mm

\noindent
{\bf A variant of the group $H^{\srfs{CMP.LC}}_{\itGamma}(X)$.}

Let $X,Y$ be open subsets of the normed spaces $E$ and $F$,
$\fnn{f}{X}{Y}$
and $\itGamma$ be a modulus of continuity.
$f$ is {\it completely weakly $\itGamma$-continuous
(CMP.WK $\itGamma$-continuous)}, if $f$ is extendible,
and there is $\gamma \in \itGamma$
such that for every $x \in \rfs{cl}(X)$ there is a neighborhood $U$
of $x$ such that
$f^{\srfs{cl}} \nrestriction U$ is $\gamma$-continuous.
As usual,
\newline
{\thickmuskip=09mu \medmuskip=08mu \thinmuskip=7mu 
$H^{\srfs{CMP.WK}}_{\itGamma}(X,Y) \ \eqdf \ %
\setm{f}{f \mbox{ is a homeomorphism between\kern5pt$X$\kern5pt
and\kern5pt $Y$\kern5pt}
\mbox{and $f$ is CMP.WK}\break
\rule{42mm}{0pt}\itGamma\mbox{-continuous}}$.
}

\begin{question}\label{metr-bldr-q1.29}
\begin{rm}
Prove the analogue of Theorem \ref{ams-bddly-lip-bldr-t6.1}(a)
for the groups of type
\newline
$H^{\srfs{CMP.WK}}_{\itGamma}(X)$.

Naturally, the definition of local $\itGamma$-LIN-borderedness has to
be replaced by the analogous notion of weak $\itGamma$-LIN-borderedness.
\vspace{-2.0mm}
\end{rm}
\end{question}

It seems that the main difficulty in proving
CMP.WK$_{\itGamma}$\,-\,determined-ness is the counterpart of
Theorem \ref{metr-bldr-t1.25}.

\kern2.0mm

\noindent
{\bf Groups which fit into the framework
but have not been investigated.}

\begin{defn}\label{d1.13}\label{metr-bldr-d1.30}
\begin{rm}
Let $\itGamma$ be a modulus of continuity and
$\fnn{f}{X}{Y}$.

(a) $f$ is {\it regionally $\itGamma$-continuous} if for every 
nonempty open $U \subseteq X$ there is a nonempty $V \subseteq U$ and
$\alpha \in \itGamma$ such that
$f \nrestriction V$ is $\alpha$-continuous.

   \index{regionally $\itGamma$-continuous}

(b) $f$ is {\it pointwise $\itGamma$-continuous}
if for every $x \in X$ there is a neighborhood $V$ of $x$ and
$\alpha \in \itGamma$ such that
$d(f(y),f(x)) \leq \alpha(d(y,x))$
for every $y \in V$.

Note that ``pointwise \rfs{MC}-continuous'' is just ``continuous''.

   \index{pointwise $\itGamma$-continuous}

(c) $f$ is {\it boundedly $\itGamma$-continuous}
if for every bounded set
$V \subseteq X$ there is  $\alpha \in \itGamma$
such that $f \nrestriction V$ is $\alpha$-continuous.

   \index{boundedly $\itGamma$-continuous}

Let $H_{\itGamma}^{\srfs{RG}}(X)$, $H^{\srfs{PW}}_{\itGamma}(X)$
and $H_{\itGamma}^{\srfs{BD}}(X)$ denote the groups of homeomorphisms
corresponding to the notions introduced in (a)\,-\,(c).

   \index{N@hrg@@$H_{\itGamma}^{\srfs{RG}}(X)$}
   \index{N@hpw@@$H_{\itGamma}^{\srfs{PW}}(X)$}
   \index{N@hbd@@@$H_{\itGamma}^{\srfs{BD}}(X)$}

\end{rm}
\end{defn}

\begin{prop}\label{p1.14}\label{metr-bldr-1.31}
\num{a} Let $X$ be a metric space and $\itGamma$ be a modulus
of continuity.
Then
\newline
\num{i}
$H^{\srfs{BD}}_{\itGamma}(X) \subseteq
H^{\srfs{LC}}_{\itGamma}(X) \subseteq
H^{\srfs{PW}}_{\itGamma}(X)$;
\newline
\num{ii}
$H^{\srfs{LC}}_{\itGamma}(X) \subseteq H^{\srfs{RG}}_{\itGamma}(X)$.

\num{b} Let $X$ be an open subset of a Banach space and
$\itGamma$ be a countably 
generated modulus of continuity. Then
$H^{\srfs{PW}}_{\itGamma}(X)\subseteq H_{\itGamma}^{\srfs{RG}}(X)$. 
\end{prop}

\noindent
{\bf Proof }
(a) Part (a) follows from the definitions.

\num{b} Suppose that $\fnn{f}{X}{Y}$ is not regionally
$\itGamma$-continuous. 
Let
$\setm{\alpha_i}{i \in \bbN}$ generate $\itGamma$. 
Let $U \subseteq X$ be an open ball which shows that $f$ is not
regionally $\itGamma$-continuous.
We define by induction $x_i,y_i \in U$.
Let $x_0,y_0$ be such that $d(f(x_0),f(y_0)) > 2 \alpha_0(d(x_0,y_0))$.
Suppose that $x_i,y_i$ have been defined.
Let
$x_{i + 1},y_{i + 1}\in B\left((x_i + y_i)/2,\
d(x_i,y_i)/2^i\right)$ be such that
$d(f(x_{i + 1}),f(y_{i + 1})) >
2 \alpha_{i + 1}(d(x_{i + 1},y_{i + 1}))$.
Since $\setm{x_i}{i \in \bbN}$ is a Cauchy sequence it converges,
say to $z$. Hence $\lim_i\,y_i = z$. 
We may assume that
$d(f(z),f(x_i)) \geq
\dgfrac{d(f(x_i),f(y_i))}{2}$ for every  $i \in \bbN$.
So for $i \in \bbN$,

$$d(f(z),f(x_i)) \geq \hbox{$\half$} d(f(x_i),f(y_i)) > 
\hbox{$\frac{1}{2}$} \cdot 2 \alpha_i(d(x_i,y_i)) > \alpha_i(d(z,x_i)).
$$

\kern1mm

\noindent
Hence $z$ shows that $f$ is not pointwise $\itGamma$-continuous.
\hfill\myqed

\kern2mm

Let
$$
K = \setm{X}{X \mbox{ is an open subset of a separable normed space
of the second category}}.
$$

\noindent
Using an argument similar to the one used in
Theorem \ref{metr-bldr-t3.41},
one can prove the 
analogues of \ref{t1.7} and \ref{t1.12} for the class
$$
\setm{H^{\srfs{RG}}_{\itGamma}(X)}{X \in K \mbox{ and $\itGamma$ is
a countably generated modulus of continuity}}.
$$
It was not checked whether other arguments used for
$H^{\srfs{LC}}_{\itGamma}(X)$ can be applied to
$H^{\srfs{RG}}_{\itGamma}(X)$.

\begin{question}\label{metr-bldr-q1.32}
\begin{rm}
Prove the analogues of \ref{t1.7} and \ref{t1.12} 
for the class
$\setm{H^{\srfs{RG}}_{\itGamma}(X)}{X
\mbox{ is an open}\break
\mbox{subset of a normed space, and }
\itGamma \mbox{ is a countably generated modulus of continuity}}$.
\hfill\proofend
\end{rm}
\end{question}

It is easy to see that a reconstruction theorem for the class of
$H^{\srfs{RG}}_{\itGamma}(X)$'s
implies reconstruction theorems for the classes of
$H^{\srfs{WK}}_{\itGamma}(X)$'s and $H^{\srfs{BD}}_{\itGamma}(X)$'s.

\kern2.0mm

\subsection{\kern-5.7mm.
Some more facts about reconstruction theorems.}
\label{ss1.8}
\label{ss1.8-additional-information}

\noindent
{\bf Reconstruction questions for related groups.}

Much work has been done on the analogous problems for
diffeomorphism groups. It seems that the first work in this direction
was carried out by F. Takens \cite{Ta}.

Soon afterwards there was an unpublished extensive work
by W. Ling \cite{Lg1} and \cite{Lg2}.
Ling proved that many types of structures on a Euclidean manifold give
rise to a determining category (or to an appropriate variant of
this notion).
Some of these categories are:
\newline
(1) The category of $k$-smooth Euclidean manifolds with $k$-smooth
diffeomorphisms.
\newline
(2) The category of $k$-smooth Euclidean manifolds with a $k$-smooth
volume form with diffeomorphisms preserving the form.
\newline
(3) The category of $k$-smooth foliated Euclidean manifolds with
the foliation preserving diffeomorphisms.
\newline
(4) Differentiable manifolds with a contact form.
\newline
(5) Manifolds with a piecewise linear structure, and homeomorphisms
preserving this structure.

The authors in \cite{RY} (unpublished) reproved Result (1) from Ling's
work, and proved some additional facts. For example, they showed
that the category of Euclidean differentiable manifolds with
diffeomorphisms that have a locally $\itGamma$-continuous
$k$'th derivative is a determining category, for every countably
generated modulus of continuity $\itGamma$.

The next work was by R. Filipkiewicz \cite{Fi}.
He proved that the category of
$k$-smooth manifolds with $k$-smooth diffeomorphisms is a determining
category.

Further work on this subject has been done more recently by a number
of authors.

A. Banyaga \cite{Ba1}, \cite{Ba2} proved the determining-ness for
the categories arising from
differentiable structures,
unimodular structures,
symplectic structures,
and contact structures.
Also, he established an analogous result for
measure preserving homeomorphisms.

T. Rybicki \cite{Ryb} presented an axiomatic approach to groups of
$C^{\infty}$ diffeomorphisms which determine a $C^{\infty}$ manifold.

Recent progress on reconstruction problems was obtained by 
J. Borzellino and V. Brunsden \cite{BB}.
They proved faithfulness for the class of spaces which are locally
compact orbifolds.

Results on differentiabilty obtained by the authors of this work
which refine older results and which also deal with Fr\'echet
differentiabilty in infinite-dimensional spaces, will appear in
a subsequent work.

V. Gol'dshtein and M. Rubin obtained analogous results for
quasi-conformal homeomorphism groups.
Part of these results appeared in \cite{GR}.
The results for quasi-conformal homeomorphism groups apply to
finite and infinite-dimensional spaces.
The full work on this subject will be presented in a separate article.

Another interesting theorem on a determining category appears
in the works of M. G. Brin and of Brin and F. Guzm\'an
on the Thompson group.
Let $G \leq H([0,1])$ be the group of all homeomorphisms $h$
such that:
(1) $h$ is piecewise linear;
(2) every slope of $h$ is an integral power of $2$;
(3) every breakpoint of $h$ is a diadic number.
It is clear that
$G \in K_{\srfs{LCM}}$, (see \ref{d2.4} and \ref{t2.2}).
Hence $\sngltn{\pair{[0,1]}{G}}$ is faithful.
Interestingly, $G$ is a finitely presented group.

One of Brin's results from
\cite{Br1} is as follows.
\newline
(1) Every automorphism of $G$ is
induced by a homeomorphism $f \in H([0,1])$ such that
for every $a < b$ in $[0,1]$,
$f \nrestriction [a,b]$ satisfies (1)-(3) above.
\newline
(2) Every such a homeomorphism induces an automorphism of $G$.

Denote by $G^+$ the group of all $f \in H([0,1])$ such that
conjugation by $f$ is an automorphism of $G$.
Brin also proves that $\sngltn{\pair{[0,1]}{G^+}}$
is a determining category.
See also Brin \cite{Br2} and Brin and F. Guzm\'an \cite{BG}.

\kern2.0mm

\noindent
{\bf Reconstruction theorems in other areas.}

The theme of reconstructing a structure from its automorphism group
was investigated in several other areas.

The recovery of a vector space from its group of linear isomorphisms
has a long history. Mackey \cite{Mac} proved in 1942 that a normed
vector space $X$ can be reconstructed from its group $L(X)$
of isomorphisms, (that is, bijective bounded linear transformations
from the space to itself). More precisely, Mackey showed that
if $X$ is finite-dimensional and $L(X) \cong L(Y)$, then
$\rfs{dim}(X) = \rfs{dim}(Y)$.
In the case that $X$ is infinite-dimensional an isomorphism between
$L(X)$ and $L(Y)$ is induced by an isomorphism between $X$ and $Y$.
In the case that $X$ is reflexive
an isomorphism between $L(X)$ and $L(Y)$ can also be induced by an
isomorphism between $X^*$ and $Y$.

Let $F_1,F_2$ be division rings and $n_1,n_2 > 2$ be integers.
If the linear groups $\rfs{GL}(n_1,F_1)$ and $\rfs{GL}(n_2,F_2)$
are isomorphic, then $n_1 = n_2$
and either $F_1 \cong F_2$ or $F_1 \cong F_2^{\rm op}$,
where $F^{\rm op}$ is the division ring obtained from $F$
by reversing the multiplication. That is,
$a \kern2pt\cdot^{F^{\rm op}}\kern2pt b = b \cdot^F a$.
This fact is due to
J. Dieudonn\'e \cite{Di1} (1947) and \cite{Di2} (1951).

For infinite-dimensional vector spaces, $V_1$ over $F_1$ and
$V_2$ over $F_2$, every isomorphism between $\rfs{Aut}(V_1)$ and
$\rfs{Aut}(V_2)$ is induced by isomorphisms between $F_1$ and $F_2$
and between $V_1$ and $V_2$. 
A strong theorem concerning this, but not exactly this fact,
was proved by C. E. Rickart in \cite{Ri1} - \cite{Ri3} (1950 - 1951).
The theorem of Dieudonn\'e for finite dimension is a special case
of Rickart's Theorem.
O. O'Meara \cite{Om} (1977)
proved the reconstruction theorem for infinite dimension.
Another proof was found by V. Tolstykh \cite{To1} (2000).

Free groups are also reconstructible from their automorphism groups.
That $\rfs{Aut}(F_n) \not\cong \rfs{Aut}(F_m)$
for $n \neq m$, can be deduced from
the work of J. Dyer and, G. P. Scott \cite{DS} (1975).
$F_n$ denotes the free group with $n$ generators
(in the variety of all groups).
E.~Formanek in \cite{Fo} (1990) proved that
$\rfs{Inn}(F_n)$ is the only normal free subgroup of rank $n$ of
$\rfs{Aut}(F_n)$. This implies immediately the reconstruction result
for finitely generated free groups.
V. Tolstykh in \cite{To2} (2000) proved that if $\lambda$ is an infinite
cardinal then $\rfs{Inn}(F_{\lambda})$
is definable in $\rfs{Aut}(F_{\lambda})$.
This implies the reconstruction result for free groups with infinite
rank.

Another body of reconstruction results for groups of
linear transformations is due to
M. Droste and M. G\"obel \cite{DG1} (1995) and \cite{DG2} (1996).
Given a ring $R$ with unity and a poset $P$
one can define the generalized McLain group $G(R,P)$ of $R$ and $P$.
Droste and G\"obel reconstruct $R$ and $P$ from $G(R,P)$.

The symmetric group is another important instance.
It is the automorphism group of a structure with no relations and no
operations. $\rfs{Sym}(6)$ is the only symmetric group which has outer
automorphisms.
A proof that $A$ is recoverable from $\rfs{Sym}(A)$ appears in
McKenzie \cite{McK} (1971). This had been known before.
See Scott \cite{Sc} p.311.

Automorphism groups of various types of ordered structures
were also extensively investigated.
We mention some of the more recent references.
Reconstruction theorems for trees appear in Rubin \cite{Ru3} (1993).
Linear orders and related structures are considered in
Rubin \cite{Ru5} (1996) and in \cite{MR}.
And Boolean algebras are reconstructed in Rubin \cite{Ru2} (1989).

The reconstruction of measure algebras is dealt with in \cite{Ru2}.
The  group of measure preserving transformations of $[0,1]$
is considered by S. Eigen in \cite{Ei} (1982).

Rubin \cite{Ru4} (1994) deals with the reconstruction of
$\aleph_0$-categorical structures.

\kern2.0mm

\noindent
{\bf Acknowledgements.}

We would like to thank Vladimir Fonf, Wieslaw Kubis,
Arkady Leiderman and Michael Levin for many very helpful discussions
and for informing us about various facts from topology and functional
anaysis which were relevant to this work. Their involvement and interest
was very valuable.

The fact that the principal modulus of continuity generated by
$\alpha$ is $\alpha$-star-closed was proved by
Kubis. See Definition \ref{ams-bddly-lip-bldr-d5.14}(d)
and Proposition \ref{ams-bddly-lip-bldr-p5.15}(a).

We also thank Vladimir Tolstykh for his help in surveying
reconstruction theorems in algebra.
His thorough survey was a great help.
Also thanked for helpful discussions are Yoav Benyamini
and Edmund Ben Ami.

\newpage

\section{Obtaining a homeomorphism from a group isomorphism}
\label{s2}

\subsection{\kern-5.7mm. Capturing the action of the group on
the regular open sets.}
\label{ss2.1}
\label{ss2.1-regular-open-sets}

Let $G \leq H(X)$. In order to prove that $X$ is
reconstructible from $G$, we shall first show that the action of $G$ on
the set of regular open subsets of $X$ is
reconstructible from $G$.

We next introduce some notations, recall some basic definitions, and
present some notions specific to this work.

\begin{defn}\label{d2.1}
\begin{rm}
Let $X$ be a topological space $U \subseteq X$ and $G \leq H(X)$.

(a) Let $\rfs{int}^X(U)$, $\rfs{cl}^X(U)$, $\rfs{bd}^X(U)$
and $\rfs{acc}^X(U)$
denote respectively the interior,
closure, boundary and the set of accumulation points of $U$ in $X$.
The boundary, $\rfs{bd}^X(U)$ is defined by
$\rfs{bd}^X(U) \eqdf \rfs{cl}^X(U) \cap \rfs{cl}^X(X - U)$.
The superscript $X$ is omitted when $X$ is understood from the context.
   \index{N@int00@@$\rfs{int}^X(U)$ interior of $U$ in $X$}
   \index{N@cl01@@$\rfs{cl}^X(U)$ closure of $U$ in $X$}
   \index{N@bd@@$\rfs{bd}^X(U)$ boundary of $U$ in $X$}
   \index{N@acc@@$\rfs{acc}^X(U)$ the set of accumulation points
          of $U$ in $X$}

(b) $U$ is {\it regular open}
   \index{regular open. A set is regular open, if it is equal to
          the interior of its closure}
if $\,U = \rfs{int}(\rfs{cl}(U))$.
$\rfs{Ro}(X)$ denotes the set of regular open subsets of $X$.
   \index{N@ro@@$\rfs{Ro}(X)$. The set of regular open subsets of $X$}
We equip $\rfs{Ro}(X)$ with the operations:
$U + V \eqdf \rfs{int}(\rfs{cl}(U \cup V))$,
$U \cdot V \eqdf U \cap V$
and $-U \eqdf \rfs{int}(X - U)$.
Then $\qdrpl{\rfs{Ro}(X)}{+}{\cdot}{-}$
is a complete Boolean algebra. Obviously,\ 
$0^{\srfs{Ro}(X)} = \emptyset$,\ 
$1^{\srfs{Ro}(X)} = X$, and the induced partial
ordering of $\rfs{Ro}(X)$ is \ $\leq^{\srfs{Ro}(X)} \ =\  \subseteq$.
We regard $\rfs{Ro}(X)$ both as a set and as a Boolean algebra.

(c) If $\iso{g}{X}{Y}$ then $g$ induces an isomorphism
$g^{\srfs{Ro}}$ between $\rfs{Ro}(X)$ and $\rfs{Ro}(Y)$:
$g^{\srfs{Ro}}(U) = \setm{g(x)}{x \in U}$.
For $G \leq H(X)$ let
$G^{\srfs{Ro}} \eqdf \setm{g^{\srfs{Ro}}}{g \in G}$.
Then $G^{\srfs{Ro}} \leq \rfs{Aut}(\rfs{Ro}(X))$
and if $X$ is Hausdorff,
then $g \mapsto g^{\srfs{Ro}}$
is an embedding of $G$ into $\rfs{Aut}(\rfs{Ro}(X))$.
We assume that $X$ is Hausdorff and identify $G$ with $G^{\srfs{Ro}}$.
So $H(X)$ is regarded as a subgroup of $\rfs{Aut}(\rfs{Ro}(X))$.

(d) $G$ is a {\it locally moving} subgroup of $H(X)$
   \index{locally moving subgroup of $H(X)$}
if for every nonempty open $V \subseteq X$ there is
$g \in G - \sngltn{Id}$ such that $g \nrestriction (X - V) = \rfs{Id}$.
In that case $\pair{X}{G}$ is called a
{\it topological local movement system}.
   \index{topological local movement system}

(e) Let $\fnn{\rfs{Ap}}{G \times \rfs{Ro}(X)}{X}$
be the {\it application function}. That is,
$\rfs{Ap}(g,V) = g(V)$.
The structure $\rfs{MR}(X,G)$ is defined as follows.
$$
\rfs{MR}(X,G) = \sixtpl{\rfs{Ro}(X)}{G}{+}{\cdot}{-}{\rfs{Ap}}.
$$

   \index{N@mr@@$\rfs{MR}(X,G) =
          \sixtpl{\rfs{Ro}(X)}{G}{+}{\cdot}{-}{\rfs{Ap}}$}

(f) $\iso{\eta}{\rfs{MR}(X,G)}{\rfs{MR}(Y,H)}$
means that $\eta$ is an isomorphism
   \index{N@AAAA@@$\iso{\eta}{\rfs{MR}(X,G)}{\rfs{MR}(Y,H)}$}
between $\rfs{MR}(X,G)$ and $\rfs{MR}(Y,H)$.
That is, $\eta$ is a bijection between $\rfs{Ro}(X) \cup G$
and $\rfs{Ro}(Y) \cup H$, $\eta(G) = H$,
and $\eta$ preserves $+,\cdot,-$ and $\rfs{Ap}$.

(g) If $\fnn{\eta}{A}{B}$ is a bijection and $\fnn{g}{A}{A}$,
then the {\it conjugation} of $g$ by $\eta$ is defined as
$g^{\eta} \eqdf \eta \scirc g \scirc \eta\inverse$.
\end{rm}
\end{defn}

\begin{prop}\label{p2.2}
Let $X,Y$ be Hausdorff spaces, $G \leq H(X)$ and $H \leq H(Y)$.
Suppose that $\iso{\varphi}{G}{H}$ and
$\iso{\eta}{\rfs{Ro}(X)}{\rfs{Ro}(Y)}$. Then
$\iso{(\varphi \cup \eta)}{\rfs{MR}(X,G)}{\rfs{MR}(Y,H)}$
\ iff\ \ $\varphi(g) = g^{\eta}$ for every $g \in G$.
\end{prop}

The next theorem says that for topological local movement systems
the action of $G$ on $\rfs{Ro}(X)$ can be reconstructed from $G$.
This theorem is proved in \cite{Ru5}.

\begin{theorem}\label{t2.3}
{\rm The reconstruction theorem for topological local
movement systems. }
Let $\pair{X}{G}$ and $\pair{Y}{H}$ be topological local movement
systems and $\iso{\varphi}{G}{H}$. Then there is a unique
$\iso{\eta}{\rfs{Ro}(X)}{\rfs{Ro}(Y)}$ such that
$\iso{(\varphi \cup \eta)}{\rfs{MR}(X,G)}{\rfs{MR}(Y,H)}$.
\newline
That is, there is a unique $\iso{\eta}{\rfs{Ro}(X)}{\rfs{Ro}(Y)}$
such that $\varphi(g) = g^{\eta}$ for every $g \in G$.
\end{theorem}

\noindent
{\bf Proof }
See \cite{Ru5} Definition 1.2, Corollary 1.4 or Corollary 2.10
and Proposition 1.8.
\hfill\myqed

\subsection{\kern-5.7mm. Faithfulness in locally compact spaces.}
\label{ss2.2}
\label{ss2.2-locallt-compact-spaces}

The first faithfulness theorem to be presented is about locally
compact spaces. It is taken from \cite{Ru1} and brought here for the
sake of completeness.
It is the conjunction of parts (a), (b) and (c) of Theorem 3.5 there.

\begin{defn}\label{d2.4}
\begin{rm}
(a) For $G \leq H(X)$, $g \in H(X)$ and $x \in X$,
let
$G(x) \eqdf \setm{g(x)}{g\in G}$.
   \index{N@AAAA@@$G(x) = \setm{g(x)}{g\in G}$}
A set $A \subseteq X$ is {\it somewhere dense}
if $int(\rfs{cl}(A)) \neq \emptyset$.
   \index{somewhere dense set. A set whose closure contains a
          nonempty open set}
$X$ is a {\it perfect space} if there is no $x \in X$ such that
$\sngltn{x}$ is open.
Suppose that $G$ is a set of permutations of a set
\hbox{$A$ and $B \subseteq A$. \kern2pt Define}
$\sprtd{G}{B} \eqdf \setm{g \in G}{g \nrestriction (A - B) = \rfs{Id}}$.

   \index{N@AAAA@@$\sprtd{G}{B}$.
          If $G \subseteq \setm{g}{\fnn{g}{A}{A}}$,
          then $\sprtd{G}{B} \eqdf
          \setm{g \in G}{g \nrestriction (A - B) = \rfs{Id}}$}

$(b)$
Let
\addtolength{\arraycolsep}{0pt}
\newline
$
\rule{0pt}{20pt}\begin{array}{lll}
\rule{0pt}{0pt}\kern-4pt
K_{\srfs{LCM}}
\kern-7pt& \eqdf \kern-7pt& \setm{\pair{X}{G}}
{\kern2pt X \mbox{ is\kern5.3pt a\kern5.3pt perfect\kern5.3pt
locally\kern5.3pt compact\kern5.3pt Hausdorff\kern5.3pt
space,\kern5.3pt
and}
\\
&&
\rule{0pt}{14pt}\kern1pt\mbox{for every open }
V \subseteq X \mbox{ and }
x \in V,\  G\sprt{V}(x) \mbox{ is somewhere dense}}.
\end{array}
$
\addtolength{\arraycolsep}{0pt}

   \index{N@klcm@@$K_{\srfs{LCM}} =
          \setm{\pair{X}{G}} { \mbox{ is Haussdorf, perfect
          locally compact and for every}\newline\indent
          \rule{25pt}{0pt}\mbox{ open } V \subseteq X \mbox{ and }
          x \in V,\  G\sprt{V}(x) \mbox{ is somewhere dense}}$}
\end{rm}
\end{defn}

\begin{theorem}\label{t2.2} {\rm (Rubin \cite{Ru1} 1989)}
$K_{\srfs{LCM}}$ is faithful.
\vspace{-3.1mm}
\end{theorem}

\noindent
{\bf Proof }
It follows easily from the definitions that
for every $\pair{X}{G} \in K_{\srfs{LCM}}$,
$\pair{X}{G}$ is a topological local movement system.

A subset $p$ of a Boolean algebra $B$ is called an {\it ultrafilter} if:
(i) $0 \not\in p$;
(ii) if $ a_1,\ldots,a_n \in p$, then
$\prod_{i = 1}^n a_i \in p$;
(iii) if $a \in p$ and $b \geq a$, then $b \in p$;
(iv) for every $a \in B$ either $a \in p$ or $-a \in p$.

By Zorn's lemma, every subset of $B$ satisfying (i)-(ii) is 
contained in an ultrafilter.
For an ultrafilter $p$ in $\rfs{Ro}(X)$, let
$A_p \eqdf \bigcap\setm{\rfs{cl}(V)}{V \in p}$.
Let $\pair{X}{G} \in K_{\srfs{LCM}}$.
We say that an ultrafilter $p$ in $\rfs{Ro}(X)$
is {\it good} if $A_p$ is a singleton.
If $p$ is good and $A_p = \sngltn{x}$, then we write $x = x_p$.
The following facts can be easily checked.

(a) $A_p = \sngltn{x}$
iff $p$ contains all regular open neighborhoods of $x$.

(b)
$p$ is good iff
there is $W \in \rfs{Ro}(X) - \sngltn{\emptyset}$
such that for every
$V \in \rfs{Ro}(X) - \sngltn{\emptyset}$:\break
if $V\subseteq W$, then there is $g \in G$ such that $g(V) \in p$.

(c) Let $p$ and $q$ be good ultrafilters.
Then $x_p \neq x_q$ iff
\newline\centerline{
\rule{0pt}{36pt}
$
\begin{array}{c}
(\exists U \in p)(\exists V \in q)
\left(\rule{0pt}{17pt}
(U \cap V = \emptyset) \wedge
(\forall U_1 \subseteq U)
\left(\rule{0pt}{15pt}
U_1 \neq\emptyset \raro
\right.\right.
\\
\rule{0pt}{11pt}\kern-12pt\left.\left.
(\exists f \in G)
\left(\rule{0pt}{11pt}
V \in f(q)\wedge U_1 \in f(p)
\right)\rule{0pt}{15pt}
\right)\rule{0pt}{17pt}\right).
\end{array}\newline
$
}

(d)\rule{0pt}{17pt}
Let $p$ be a good ultrafilter, and $U \in \rfs{Ro}(X)$. Then
$x_p \in U$ iff for every good ultrafilter $q$\,:
\,if $x_q = x_p$, \,then $U \in q$.

(e) Let $p,q$ be good ultrafilters, and $g \in G$.
Then $g(x_p) = x_q$ iff $x_{g(p)} = x_q$.

(f) If $p$ is a good ultrafilter and $g \in G$, then $g(p)$ is
a good ultrafilter.

(g) If $x \in X$, then there is a good ultrafilter $p$ such that
$x_p = x$.

Clearly, the fact that $p$ is an ultrafilter is expressible in terms of
the operations of 
$\fourtpl{\rfs{Ro}(X)}{+}{\cdot}{-}$.

(1) By (b), the fact that $p$ is a good ultrafilter is expressible
in terms of the operations of 
$\rfs{MR}(X,G)$.

(2) By (c), for good ultrafilters $p$ and $q$, the fact that $x_p = x_q$
is expressible in terms of the operations of $\rfs{MR}(X,G)$.

(3) By (d), for a good ultrafilter $p$ and $U \in \rfs{Ro}(X)$,
the fact that $x_p \in U$
is expressible in terms of the operations of $\rfs{MR}(X,G)$.

(4) By (e), for good ultrafilters $p$ and $q$ and $g \in G$,
the fact that $g(x_p) = x_q$
is expressible in terms of the operations of $\rfs{MR}(X,G)$.

Let $\pair{X}{G}, \pair{Y}{H} \in K_{\srfs{LCM}}$,
and let $\iso{\varphi}{G}{H}$. By Theorem \ref{t2.3},
there is\break
$\iso{\eta}{\rfs{Ro}(X)}{\rfs{Ro}(Y)}$ such that
$\iso{(\varphi \cup \eta)}{\rfs{MR}(X,G)}{\rfs{MR}(Y,H)}$.
Let $\psi = \varphi \cup \eta$.
We define $\fnn{\tau}{X}{Y}$.
Let $x \in X$. By (g), there is an ultrafilter $p$ such that
$x_p = x$. By (1), $\psi(p)$ is a good ultrafilter.

We define $\tau(x) = x_{\psi(p)}$.
If $q$ is a good ultrafilter such that also $x_q = x$, then by (2),
$x_{\psi(q)} = x_{\psi(p)}$. So the definition of $\tau$ is valid.

We check that $\tau$ is a bijection between $X$ and $Y$.
Suppose that $x_p \neq x_q$.
By (2), $\tau(x_p) = x_{\psi(p)} \neq x_{\psi(q)} = \tau(x_q)$.
So $\tau$ is injective.

Let $y \in Y$. By (g), there is an ultrafilter $q$ such that
$x_q = y$. By~(1), $p \eqdf \psi\inverse(q)$ is a good ultrafilter.
So $\tau(x_p) = x_{\psi(p)} = x_q = y$.
So $\tau$ is surjective.

Let $\tau(A)$ denote $\setm{\tau(a)}{a \in A}$. In order to show that
$\tau$ is a homeomorphisms, it suffices to show that for some
open base $\calB$ of $X$, $\setm{\tau(U)}{U \in \calB}$
is an open base for $Y$.
Since $X$ and $Y$ are locally compact, they are regular spaces. So
$\rfs{Ro}(X)$ and $\rfs{Ro}(Y)$ are open bases of $X$ and $Y$
repectively.
So it suffices to show that
$\setm{\tau(U)}{U \in \rfs{Ro}(X)} = \rfs{Ro}(Y)$.
Let $x \in X$ and $U \in \rfs{Ro}(X)$.
Let $p$ be an ultrafilter such that $x_p = x$.
By (3), $x_p \in U$ iff $x_{\psi(p)} \in \psi(U)$.
That is, $x \in U$ iff $\tau(x) \in \psi(U)$.
So $\tau(U) = \psi(U)$ for every $U \in \rfs{Ro}(X)$.
Hence
$\setm{\tau(U)}{U \in \rfs{Ro}(X)} =
\setm{\psi(U)}{U \in \rfs{Ro}(X)} = \rfs{Ro}(Y)$.
So $\tau$ is a homeomorphism.

It remains to show that $\tau$ induces $\varphi$.
Let $g \in G$ and $y \in Y$.
Let $q$ be an ultrafilter in $\rfs{Ro}(Y)$ such that $x_q = y$.
Then $g^{\tau}(y) = \tau \scirc g \scirc \tau\inverse(x_q) =
\tau \scirc g(x_{\psi\inverse(q)}) =
\tau(x_{g(\psi\inverse(q))}) =
x_{\psi(g(\psi\inverse(q)))} =
x_{\eta(g(\eta\inverse(q)))} =
x_{g^{\eta}(q)}$.
But by Proposition \ref{p2.2}, $g^{\eta} = \varphi(g)$.
So $x_{g^{\eta}(q)} = x_{\varphi(g)(q)}$.
However, if $x_q = y$, then trivially $x_{h(q)} = h(y)$
for every $h \in H$.
In particular, $x_{\varphi(g)(q)} = \varphi(g)(y)$.

We have shown that $g^{\tau}(y) = \varphi(g)(y)$ for every $y \in Y$.
So $g^{\tau} = \varphi(g)$.
\vspace{1mm}\hfill\myqed

{\bf Remark } In the above proof the existence of the inducing
homeomorphism $\tau$ was deduced from facts (b) - (e) which showed
that point representation, equality, belonging and application
were expressible in $\rfs{MR}(X,G)$.
The toil of deducing the existence of $\tau$ from (b) - (e)
could have been spared by using certain general machinery called the
method of interpretation.
The notion of interpretation is not introduced here,
since it is used only twice.
Interpretations are described e.g. in \cite{Ru2} Section~2 or in
\cite{MR} Section~6.

Theorem \ref{t2.2} has many applications in the Euclidean case.
For example, it applies to $m$ times continuously differentiable
Euclidean manifolds.

\begin{cor}\label{c2.3} {\rm \cite{Ru1}}
Let $K_D = \{\pair{X}{G}|$
for some $0 \leq m \leq \infty$, $X$ is a Euclidean $C^m$-manifold
and $G$ contains all homeomorphisms $f$ such that both $f$ and
$f^{-1}$ are $C^m$ homeomorphisms$\}$. Then $K_D$ is faithful.
\vspace{-2mm}
\end{cor}

\noindent
{\bf Proof }
$K_D\sbeq K_{\srfs{LCM}}$.\vspace{1.5mm} \hfill\myqed

Theorem \ref{t2.2} also applies to Hilbert cube manfolds, and in fact
to manifolds over $[0,1]^{\lambda}$ for any cardinal $\lambda$.

The class of Menger manifolds is also a subclass of $K_{\srfs{LCM}}$,
and hence it is faithful. See Kawamura \cite{K}.

\noindent

The finitely presented subgroups of $H(\bbR)$ defined by
R. Thompson (see \cite{Br1}, \cite{Br2} and \cite{BG}) also belong
to $K_{\srfs{LCM}}$.

\subsection{\kern-5.7mm. Faithfulness in normed and Banach spaces.}
\label{ss2.3}
\label{ss2.3-banach-and-normed-faithfulness}

We now turn to the context of normed vector spaces and Banach spaces.

To avoid notational complications, we shall mainly deal
with open subsets of normed and Banach spaces and not with manifolds
over such spaces.
Nevertheless, all theorems and proofs transfer
(with a correct translation)
to manifolds.
In this section, Definition \ref{d2.14} and Theorem \ref{t2.15}
deal with the setting of manifolds (and indeed with a somewhat
more general setting).

Manifolds are considered again at the end of Chapter \ref{s3}
starting from Definition \ref{metr-bldr-d3.46}.
\vspace{1.5mm}

Recall that for a metric space $X$,
$\rfs{LIP}(X) = \setm{h \in H(X)}{h \mbox{ is bilipschitz}}$
and\break
$\rfs{LIP}^{\srfs{LC}}(X) = \setm{h \in H(X)}{h
\mbox{ is locally bilipschitz}}$.

For a normed space $E$, an open set $S \subseteq E$ and a dense
linear subspace
$F \subseteq E$, we shall use the notations
$\rfs{LIP}(X;S,F)$, $\rfs{LIP}^{\srfs{LC}}(X;S,F)$,
$\rfs{LIP}(X;F)$ and $\rfs{LIP}^{\srfs{LC}}(X,F)$
introduced in Definition \ref{d1.2}.

We shall prove the faithfulness of the classes $K_{\srfs{B}}$
and $K_{\srfs{N}}$
defined below. However, these faithfulness results do not suffice
for some of the continuations. To this end we define the bigger
class $K_{\srfs{BNO}}$ and prove its faithfulness.

\begin{defn}\label{d2.7}
\begin{rm}
Let $E$ be a normed space, $X \subseteq E$ be open,
$\calS$ be a set of open subsets of $X$ and
$\calF = \setm{F_S}{S \in \calS}$ be a family of dense linear
subspaces of $E$ indexed by $\calS$.
Then $\calF$ is called a {\it subspace choice} for $\calS$.
If $\calS$ is a cover of $X$,
then $\frtpl{E}{X}{\calS}{\calF}$ is called
a {\it subspace choice system}.
   \index{subspace choice}
   \index{subspace choice system}

   \index{N@lip04@@$\rfs{LIP}(X;\calS,\calF)$.
	  Subgroup of $H(X)$ generated by
          $\bigcup \setm{\rfs{LIP}(X;S,F_S)}{S \in \calS}$}
   \index{N@liplc04@@
          {\thickmuskip=2mu \medmuskip=1mu \thinmuskip=1mu 
	  $\rfs{LIP}^{\srfs{LC}}(X;\calS,\calF)$.
	  Subgroup of $H(X)$ generated by
          $\bigcup \setm{\rfs{LIP}^{\srfs{LC}}(X;S,F_S)}{S \in \calS}$}}

(a)
$\rfs{LIP}(X;\calS,\calF)$ denotes the subgroup of $H(X)$ generated by
$\bigcup \setm{\rfs{LIP}(X;S,F_S)}{S \in \calS}$.
$\rfs{LIP}^{\srfs{LC}}(X;\calS,\calF)$
denotes the subgroup of $H(X)$ generated by
$\bigcup \setm{\rfs{LIP}^{\srfs{LC}}(X;S,F_S)}{S \in \calS}$.
Also, $\rfs{LIP}(X,\calS)$ denotes the subgroup of $H(X)$ generated by
$\bigcup \setm{\rfs{LIP}(X,S)}{S \in \calS}$.\break
The group $\rfs{LIP}^{\srfs{LC}}(X,\calS)$ is defined analogously.
   \index{N@lip05@@$\rfs{LIP}(X,\calS)$.
The subgroup of $H(X)$ generated by
$\bigcup \setm{\rfs{LIP}(X,S)}{S \in \calS}$}
   \index{N@liplc05@@$\rfs{LIP}^{\srfs{LC}}(X,\calS)$.
The subgroup of $H(X)$ generated by
$\bigcup \setm{\rfs{LIP}^{\srfs{LC}}(X,S)}{S \in \calS}$}

(b)
Let $K_{\srfs{B}}$ be the class of all $\pair{X}{G}$'s
such that $X$ is an open subset of some Banach space,
and $\rfs{LIP}(X) \leq G \leq H(X)$.
   \index{N@kb@@$K_{\srfs{B}} =
          \setm{\pair{X}{G}}{X \mbox{ is an open subset of a
          Banach space and } \rfs{LIP}(X) \leq G \leq H(X)}$}

Let $K_{\srfs{N}}$ be the class of all $\pair{X}{G}$'s
such that $X$ is an open subset of some normed space,
and $\rfs{LIP}^{\srfs{LC}}(X) \leq G \leq H(X)$.
   \index{N@kn@@$K_{\srfs{N}} =
          \setm{\pair{X}{G}}{X \mbox{ is an open subset of a
          normed space and }
          \rfs{LIP}^{\srfs{LC}}(X) \leq G \leq H(X)}$}

Let $K_{\srfs{BO}}$ be the class of all $\pair{X}{G}$'s
such that:
\begin{itemize}
\addtolength{\parskip}{-11pt}
\addtolength{\itemsep}{06pt}
\item[(1)] $X$ is an open subset of some Banach space $E$,
\item[(2)] there are an open cover $\calS$ of $X$
and a subspace choice $\calF$ for $\calS$ such that\break
$\rfs{LIP}(X;\calS,\calF) \leq G \leq H(X)$.
\vspace{-05.7pt}
\end{itemize}
   \index{N@kbo@@$K_{\srfs{BO}}$}

Let $K_{\srfs{NO}}$ be the class of all $\pair{X}{G}$'s
such that:
\begin{itemize}
\addtolength{\parskip}{-11pt}
\addtolength{\itemsep}{06pt}
\item[(1)] $X$ is an open subset of some normed space $E$,
\item[(2)] there are an open cover $\calS$ of $X$
and a subspace choice $\calF$ for $\calS$ such that\break
$\rfs{LIP}^{\srfs{LC}}(X;\calS,\calF) \leq G \leq H(X)$.
\vspace{-05.7pt}
\end{itemize}
   \index{N@kno@@$K_{\srfs{NO}}$}

Let $K_{\srfs{BNO}} = K_{\srfs{BO}} \cup K_{\srfs{NO}}$.
   \index{N@kbno@@$K_{\srfs{BNO}}$}
If $\pair{X}{G} \in K_{\srfs{BNO}}$ and $E,\calS,\calF$ are as above,
then\break
$\fivtpl{E}{X}{\calS}{\calF}{G}$ is called a
{\it BNO-system}.
   \index{bnosyst@@BNO-system}
\end{rm}
\end{defn}

\begin{theorem}\label{t2.4}\label{metr-bldr-t2.8}
\num{a} $K_{\srfs{B}} \cup K_{\srfs{N}}$ is faithful.

\num{b} $K_{\srfs{BNO}}$ is faithful.
\end{theorem}

Note that $K_{\srfs{B}} \cup K_{\srfs{N}} \subseteq K_{\srfs{BNO}}$.
So only (b) has to be proved.

\begin{remark}\label{metr-bldr-r2.9}
\begin{rm}
(a) Dealing with the larger but less natural classes of groups
$\rfs{LIP}(X;\calS,\calF)$ and $\rfs{LIP}^{\srfs{LC}}(X;\calS,\calF)$
needs justification.
Certainly the groups $\rfs{LIP}(X)$ and $\rfs{LIP}^{\srfs{LC}}(X)$
are those that come to mind first.
There are two classes of groups which merit attention
for which Theorem \ref{t2.4}(a) does not suffice, 
but Theorem \ref{t2.4}(b) does.

Let $E$ be a normed vector space and $X \subseteq E$ be open.
The group of extendible homeomorphisms of $X$ is defined as follows:
$$
\rfs{EXT}^E(X) =
\setm{h \nrestriction X}{h \in H(\rfs{cl}^E(X)) \mbox{ and }
h \nrestriction X \in H(X)}.
$$
If $E$ is a Banach space, then $\rfs{LIP}(X) \subseteq \rfs{EXT}^E(X)$.
However, if $E$ is not complete, then
$\rfs{LIP}(X) \not\subseteq \rfs{EXT}^E(X)$.

For $h \in \rfs{EXT}(X)$ let $h^{\srfs{cl}}$ denote the extension of $h$
to $\rfs{cl}^E(X)$.
Let $\itGamma$ be a modulus of continuity.
Define
{\thickmuskip=2.9mu \medmuskip=2mu \thinmuskip=1mu 
$$
\rule{0pt}{1pt}
\kern2pt H_{\itGamma}^{\srfs{CMP.LC}}(X) =
\setm{h \in \rfs{EXT}(X)}{\mbox{for some } \alpha \in \itGamma,\ 
h^{\srfs{cl}}
\mbox{ is locally $\alpha$-bicontinuous}}.
$$
}
Then Theorem \ref{t2.4}(a) does not apply to
$H_{\itGamma}^{\srfs{CMP.LC}}(X)$, but Theorem \ref{t2.4}(b) does.

Another such example is the following group.
Let $E$ be a finite-dimensional normed space,
$X \subseteq E$ be open and
$$
H = \setm{h \in H(X)}{\rfs{cl}(\setm{x \in X}{h(x) \neq x})
\mbox{ is compact}}.
$$
Then $G \not\supseteq \rfs{LIP}(X)$, but nevertheless $X$ is
reconstructible from $G$.

The reason for introducing the group $\rfs{LIP}(X;F)$ is as follows.
For an incomplete normed space $X$, we give a proof that $X$ is
reconstructible from $G$'s which contain $\rfs{LIP}^{\srfs{LC}}(X)$.
but we do not know whether $X$ is reconstructible from
$\rfs{LIP}(X)$.
In fact, every member of $\rfs{LIP}(X)$ can be uniquely
extended to a homeomorphism of $\overX$, the completion of $X$.
So $\rfs{LIP}(X)$ can be regarded as a subgroup of $H(\overX)$.
By considering $\rfs{LIP}(\overX;X)$ we prove that $\overX$
is reconstructible from $\rfs{LIP}(X)$.
It remains open
(except for spaces of the first caregory)
whether $X$ is reconstructible from $\rfs{LIP}(X)$.

(b) The groups $\rfs{LIP}(X;\calS,\calF)$ and $\rfs{LIP}(X)$
in Theorem \ref{t2.4} can be replaced by the following smaller groups.
Suppose that a normed or a Banach space $E$ has an equivalent norm
which is
$C^m$, $m \leq \infty$, that is, a norm which is $m$ times
continuously \Frechet differentiable at every $x \neq 0$.
We define $\rfs{Diff}^m(X)$ to be the group of all homeomorphisms
$g$ of $X$ such that $g,g\inverse$ are $C^m$,
and whose first derivative is bounded.
The group $\rfs{Diff}^m(X;\calF,\calS)$ is defined in analogy to
$\rfs{LIP}(X;\calS,\calF)$,
and the classes $K_{\srfs{BD}^m}$, $K_{\srfs{ND}^m}$ and
$K_{\srfs{BNOD}^m}$ are defined in analogy
to $K_{\srfs{B}}$, $K_{\srfs{N}}$ and $K_{\srfs{BNO}}$.
Then Theorem \ref{t2.4} remains true.
The proof remains the same.
The only difference is that the homeomorphisms which are constructed
in the proof of Theorem \ref{t2.4} have to be in this case $C^m$
and not just bilipschitz.

This variant of Theorem \ref{t2.4} will be needed in a subsequent work
where groups of Fr\'echet differentiable homeomomorphisms will be
considered.
\end{rm}
\end{remark}

{\bf An explanation of the method of proof of Theorem \ref{t2.4}.}

We show that there is a property 
$P(x,y)$ of pairs $\pair{x}{y}$ which is expressible in terms
of the operations of $\rfs{MR}(X,G)$
such that
for every $\pair{X}{G} \in K_{\srfs{BNO}}$ and $U,\,V \in \rfs{Ro}(X)$:
\newline\centerline{
$P(U,V)$ holds in $\rfs{MR}(X,G)$ 
\ iff \ $\rfs{cl}(U) \cap \rfs{cl}(V)$ is a singleton.
}
A pair $\pair{U}{V}$ satisfying $P$ is called a
{\it point representing pair}.

We shall then prove two similar facts.
\begin{itemize}
\addtolength{\parskip}{-11pt}
\addtolength{\itemsep}{06pt}
\item[(1)] 
There is a property $Q(x_1,y_1,x_2,y_2)$ expressible in
terms of the operations of\break
$\rfs{MR}(X,G)$
such that for every $\pair{X}{G} \in K_{\srfs{BNO}}$
and point representing pairs
$\pair{U_1}{V_1},\break\pair{U_2}{V_2} \in (\rfs{Ro}(X))^2$:
\newline\centerline{\rule{0pt}{16pt}
$Q(U_1,V_1,U_2,V_2)$ holds in $\rfs{MR}(X,G)$ \ iff \ 
$\rfs{cl}(U_1) \cap \rfs{cl}(V_1) = \rfs{cl}(U_2) \cap \rfs{cl}(V_2)$.
}
\item[(2)]
There is a property 
$S(x,y,z)$ expressible in terms of the operations of $\rfs{MR}(X,G)$
such that for every
$\pair{X}{G} \in K_{\srfs{BNO}}$,
a point representing pair
$\pair{U}{V} \in (\rfs{Ro}(X))^2$ and $W \in \rfs{Ro}(X)$:
\newline\rule{12mm}{0pt}\rule{0pt}{16pt}
$S(U,V,W)$ holds in $\rfs{MR}(X,G)$ \ iff \ 
$\rfs{cl}(U) \cap \rfs{cl}(V) \subseteq W$.
\hfill
\vspace{-05.7pt}
\end{itemize}
As in the proof of \ref{t2.2}, the existence of properties
$P$, $Q$ and $S$ implies that every isomorphism between 
$\rfs{MR}(X,G)$ and $\rfs{MR}(Y,H)$ is induced by a 
homeomorphism between $X$~and~$Y$.

The following conventions are kept through Lemma \ref{l2.13}
and the proof of Theorem~\ref{t2.4}.

\kern3pt

(a) In what follows,
$\fivtpl{E}{X}{\calS}{\calF}{G}$ denotes a BNO-system.
That is, $E$ denotes a normed vector space,
$X$ is an open subset of $E$,
$\calS$ is a cover of $X$,
$\calF$ is a subspace choice for $\calS$ and $G \leq H(X)$.
If $E$ is a Banach space, then $\rfs{LIP}(X;\calS,\calF) \leq G$,
and if $E$ is incomplete,
then $\rfs{LIP}^{\srfs{LC}}(X;\calS,\calF) \leq G$.

If $X$ is an open subset of $E$
and $\pair{X}{G} \in K_{\srfs{B}} \cup K_{\srfs{N}}$,
then $\pair{X}{G}$ is regarded as a BNO-system with
$\calS = \sngltn{X}$ and $F_X = E$.

(b) Also, $U,V,W$ denote members of $\rfs{Ro}(X)$.
If $A \subseteq X$, then $\rfs{cl}^X(A)$ and $\rfs{int}^X(A)$
are abbreviated by $\rfs{cl}(A)$ and $\rfs{int}(A)$ respectively.

\begin{defn}\label{d2.5}
\begin{rm}
(a) For $U,\,V \in \rfs{Ro}(X)$ let $U \cong V$ denote that
$(\exists g \in G)(g(U) = V)$.

   \index{N@AAAA@@$U \cong V$. This means $(\exists g \in G)(g(U) = V)$}

(b) $U$ is a {\it small set},
   \index{small set}
if there is $W \neq \emptyset$
such that for every $\emptyset \neq W' \subseteq W$
there is $U' \cong U$ such that $U' \subseteq W'$.

(c) $U$ is {\it strongly small} in $V$\ $(U \prec V)$,
   \index{strongly small set}
   \index{N@AAAA@@$U \prec V$. This means $U$ is strongly small in $V$}
if there is $\emptyset \neq W \subseteq V$
such that for every $\emptyset \neq W_1 \subseteq W$ there is
$g \in \sprtd{G}{V}$
such that $g(U) \subseteq W_1$.

(d) $U$ is {\it strongly separated} from $W$ ($U \spr W$),
   \index{strongly separated. $U$ is strongly separated from $V$,
   if $\exists W(U \prec W \mbox{ and } W \cap V = \emptyset$)}
   \index{N@AAAA@@$U \spr V$.
          This means $U$ is strongly separated from $V$}
if there is $V \in \rfs{Ro}(X)$ such that $U \prec V$
and $V \cap W = \emptyset$.
\end{rm}
\end{defn}

\begin{remark}\label{nr2.11}
\begin{rm}
Properties ``$U \cong V$'', ``$U$ is a small set'', 
``$U \prec V$'' and ``$U \spr W$''
are expressible in terms of the operations of $\rfs{MR}(X,G)$.
Formally this means the following statements.
\newline
(1)
Let $\chi_{\tcong}(x,y) \equiv (\exists z \in G)(\rfs{Ap}(z,x) = y)$.
Then $U,V$ satisfy $\chi_{\tcong}$ in $\rfs{MR}(X,G)$ iff \ 
$U \cong V$.
\newline
(2) Let $\chi_{\tsubseteq}(x,y) \ \ \equiv \ \ x \cdot y = x$.
Then $U,V$ satisfy $\chi_{\tsubseteq}$ in $\rfs{MR}(X,G)$ iff
$U \subseteq V$.
\newline
(3) Let  $\chi_{\temptyset}(x) \ \equiv \ 
(\forall y \in \rfs{Ro}(X))(x \cdot y = x)$.
Then $U$ satisfies $\chi_{\temptyset}$ in $\rfs{MR}(X,G)$ iff
$U = \emptyset$.
\newline
(4) Let
\newline
$\chi_{\tSml}(x) \ \equiv \ 
(\exists y \in \rfs{Ro}(X))
\left(\rule{0pt}{14pt}\neg \chi_{\temptyset}(y) \wedge
(\forall y' \in \rfs{Ro}(X))
\left(\rule{0pt}{14pt}\kern-3pt
\left(\rule{0pt}{10pt}
\chi_{\tsubseteq}(y',y) \wedge \neg \chi_{\temptyset}(y')
\right)
\rightarrow\right.\right.$
\newline\rule{17.5mm}{0pt}
$\left.\left.
\rule{2pt}{0pt}(\exists x' \in \rfs{Ro}(X))
\left(\rule{0pt}{10pt}
\chi_{\tcong}(x',x) \wedge \chi_{\tsubseteq}(x',y')
\right)
\kern-3pt\rule{0pt}{14pt}\right)
\kern-3pt\right)$.
\newline
\rule{0pt}{14pt}\kern-0.0pt
Then $U$ satisfies $\chi_{\tSml}$ in $\rfs{MR}(X,G)$ iff
$U$ is small.
\newline
(5) Let 
$\chi_{\tSpprtd}(x,y) \equiv
(\forall z \in \rfs{Ro}(X))(\chi_{\temptyset}(z \cdot y) \rightarrow
(\rfs{Ap}(x,z) = z))$.
Then $g,V$ satisfy $\chi_{\tSpprtd}$ in $\rfs{MR}(X,G)$ iff
$g \in G\sprt{V}$.

Similar formulas $\chi_{\tprec}$ and $\chi_{\tspr}$
can be written for $U \prec V$ and for $U \spr V$.
The above formulas use only the operations
$+$, $\cdot$, $-$ and $\rfs{Ap}$.
So if $\chi$ is any of the above formulas,
$\iso{\psi}{\rfs{MR}(X,G)}{\rfs{MR}(Y,H)}$ and $U,V \in \rfs{Ro}(X)$,
then $U,V$ satisfy $\chi$ in $\rfs{MR}(X,G)$ iff
$\psi(U),\psi(V)$ satisfy $\chi$ in $\rfs{MR}(Y,H)$.
So smallness, $\prec$, $\spr$ etc.\ are preserved under isomorphisms.
\end{rm}
\end{remark}

\begin{defn} \label{nd2.11}
\begin{rm}
(a) For a metric space $\rpair{Z}{d}$, $x \in Z$ and $r > 0$
we define
   \index{N@b00@@$B^Z(x,r) = \setm{y \in Z}{\,d(x,y) < r}$}
   \index{N@s00@@$S^Z(x,r) = \setm{y \in Z}{\,d(x,y) = r}$}
   \index{N@b01@@$\rule{0pt}{0pt}\overB^E(x,r) =
          \setm{y \in E}{\,d(x,y) \leq r}$}
$B^Z(x,r) \eqdf \setm{y \in Z}{\,d(x,y) < r}$,
$S^Z(x,r) \eqdf \setm{y \in Z}{\,d(x,y) = r}$
and
$\overB^Z(x,r) \eqdf \setm{y \in Z}{\,d(x,y) \leq r}$.
If $A \subseteq Z$, then
$B^Z(A,r) \eqdf \bigcup_{x \in A} B^Z(x,r)$.

   \index{N@b02@@$B^Z(A,r) = \bigcup_{x \in A} B^Z(x,r)$}

In the context of this section there are two metric spaces involved:
a normed space $E$ and an open subset $X \subseteq E$.
We use $B(x,r)$, $S(x,r)$ and $\overB(x,r)$ as
abbreviations of
$B^X(x,r)$, $S^X(x,r)$ and $\overB^X(x,r)$.

   \index{N@b03@@\rule{0pt}{0pt}$B(x,r)$.
          An abbrviation of $B^X(x,r)$}

For $x,y \in E$, $[x,y]$ denotes the line segment connecting $x$
and $y$.
   \index{N@AAAA@@$[x,y]$. The line segment with endpoints $x$ and $y$}
For $v \in E$ let $\fnn{\rfs{tr}_v^E}{E}{E}$
be the translation by $v$, that is, $\rfs{tr}_v^E(x) = v + x$.
Whenever $E$ can be understood from the context,
$\rfs{tr}_v^E$ is abbreviated by $\rfs{tr}_v$.
   \index{N@tr00@@$\rfs{tr}_v^E$. Translation by $v$. For $v,x \in E$,
          $\rfs{tr}_v^E(x) = v + x$}
   \index{N@tr01@@$\rfs{tr}_v$. Abbreviation of $\rfs{tr}_v^E$}

(b) Let $\calN = \fivtpl{E}{X}{\calS}{\calF}{G}$ be a BNO-system
and $B = B^E(x,r)$ be a ball of $E$.
$B$ is a {\it manageable ball} of $X$ (with respect to $\calN$)
   \index{manageable ball (with respect to a BNO-system)}
if there are $S \in \calS$ and $\varepsilon > 0$
such that $x \in S \cap F_S$ and $B^E(x,r + \varepsilon) \subseteq S$.
In such a case we say that {\it $B$ is based on $S$}.
   \index{manageable ball $B$ based on $S$}
Note that if $B = B^E(x,r)$ is a manageable ball,
then $B^E(x,r) = B^X(x,r)$ and $\rfs{cl}^E(B) = \rfs{cl}^X(B)$.

(c) For a topological space $Y$ and $h \in H(Y)$, the support of $h$
its defined as
   \index{N@supp@@$\rfs{supp}(h) = \setm{y \in Y}{h(y) \neq y}$}
$$
\rfs{supp}(h) = \setm{y \in Y}{h(y) \neq y}.
$$
\end{rm}
\end{defn}

\begin{prop}\label{np2.12}
\num{a} Suppose that $Y$ is any topological space,
and let $H \leq H(Y)$.\break
For $k \in H$ let $\fnn{\psi_k}{\rfs{MR}(Y,H)}{\rfs{MR}(Y,H)}$
be defined as follows. For every $h \in H$,\break
{\thickmuskip=3.5mu \medmuskip=2.5mu \thinmuskip=1.5mu 
$\psi_k(h) = h^k$,
and for every $U \in \rfs{Ro}(Y)$,
$\psi_k(U) = \setm{h(x)}{x \in U\kern-1pt}$.
Then $\psi_k \in \rfs{Aut}(\rfs{MR}(Y,H))$.}

\num{b} Let $Y$ be any topological space.
\begin{itemize}
\addtolength{\parskip}{-11pt}
\addtolength{\itemsep}{06pt}
\item[\num{i}] If $F \subseteq Y$ is closed,
then $\rfs{int}(F) \in \rfs{Ro}(Y)$.
\item[\num{ii}] $\rfs{int}(\rfs{cl}(A)) \in \rfs{Ro}(Y)$
for every $A \subseteq Y$.
\item[\num{iii}] $\rfs{int}(\rfs{cl}(A))$
is the minimal regular open set containing $A$.
\item[\num{iv}] If $T,S \subseteq Y$ are disjoint open sets,
then $\rfs{int}(\rfs{cl}(T)) \cap S = \emptyset$.
\vspace{-05.7pt}
\end{itemize}
\end{prop}

{\bf Proof } Trivial. \hfill\myqed

We shall next construct certain homeomorphisms in
$\rfs{LIP}(X;\calS,\calF)$.
Geometrically, their existence is quite obvious.
However, the formal proof requires some computation.

All balls mentioned in the next lemma are manageable.
For such balls we write $B^E(x,r) = B(x,r)$.
Part (b)(ii) of the lemma will be used in Chapter \ref{s3}.
See Proposition \ref{metr-bldr-p3.4}.

\begin{lemma}\label{l2.6}
\num{a} Suppose that $B = B(x_0,r_0)$ is a manageable ball based on~$S$,
$x_0 \in F_S$ and
$0 < s_0 < s_1 < r_0$.
Then there is $h \in \sprtd{\rfs{LIP}(X;\calS,\calF)}{B}$ such that
$h(B(x_0,s_1)) = B(x_0,s_0)$.

\num{b} Suppose that $B = B(x_0,r_0)$ is a manageable ball based on $S$,
$x_0,v \in B \cap F_S$, $0 < r < r_0$
and $0 < s < r_0 - \norm{v - x_0}$.
Then
\begin{itemize}
\addtolength{\parskip}{-11pt}
\addtolength{\itemsep}{06pt}
\item[\num{i}] There is $h \in \sprtd{\rfs{LIP}(X;\calS,\calF)}{B}$
such that $h(B(x_0,r)) = B(v,s)$.
\item[\num{ii}] If also $r = s$,
then $h$ is $(1 + \frac{\norm{v}}{r_0 - r - \norm{v}})$-bilipschitz
and $h \nrestriction B(x_0,r) = \rfs{tr}_v \nrestriction B(x_0,r)$.
\vspace{-05.7pt}
\end{itemize}

\num{c} Let $B$ be a manageable ball based on $S$,
$x,y \in B \cap F_S$ and $r > 0$. Assume that $B([x,y],r) \subseteq B$.
Then there is
$h \in \sprtdl{\rfs{LIP}(X;\calS,\calF)}{B([x,y],r)}$
such
that $h \nrestriction B(x,\dgfrac{2r}{3}) =
\rfs{tr}_{y - x} \nrestriction B(x,\dgfrac{2r}{3})$.
Moreover, there is a function $K_{\srfs{seg}}(\ell,t)$
increasing in $\ell$ and
decreasing in $t$ such that the above $h$ is
$K_{\srfs{seg}}(\norm{x - y},r)$-bilpschitz.

\num{d} Let $U \subseteq X$ be open,
$\fnn{\gamma}{[0,1]}{U}$ be continuous and $\onetoonen$
and $s \in (0,1]$.
Then there is $h \in \rfs{LIP}(X)$ such that $h(\gamma(0)) = \gamma(0)$,
$h(\gamma(s)) = \gamma(1)$ and $\rfs{supp}(h) \subseteq U$.
\end{lemma}

\noindent
{\bf Proof }
(a) Assume for simplicity that $x_0 = 0$.
Let $g \in H([0,\infty))$ be the piecewise linear function with
breakpoints at $s_0$ and $r_0$ such that
$g(s_0) = s_1$
and $g(t) = t$ for every
$t \geq r_0$. Then $g$ is $K$-bilipschitz with
$K = \max(\frac{s_1}{s_0},\frac{r_0 - s_0}{r_0 - s_1})$.

We define $\fnn{h}{E}{E}$.
$$
h(x) = g(\norm{x}) \frac{x}{\norm{x}} \mbox{ if } x \neq 0
\mbox{ \ and \ } h(0) = 0.
$$
Let $x,y \in E$.
We may assume that $0 \neq \norm{y} \leq \norm{x}$.
Let $z = \norm{y} \frac{x}{\norm{x}}$.
Then
$\norm{x - z} = \norm{x} - \norm{y} \leq \norm{x - y}$
and
$\norm{z - y} \leq \norm{z - x} + \norm{x - y} \leq 2 \norm{x - y}$.
So
\vspace{1.5mm}
\newline
\rule{7pt}{0pt}
\renewcommand{\arraystretch}{1.5}
\addtolength{\arraycolsep}{-6pt}
$
\begin{array}{ll}
&
\norm{h(x) - h(y)} \leq
\norm{h(x) - h(z)} + \norm{h(z) - h(y)} \leq
K \norm{x - z} + \frac{g(\norm{y})}{\norm{y}} \norm{z - y}
\\
\leq
\rule{5pt}{0pt}
&
K \norm{x - y} + K \cdot 2 \norm{x - y} = 3K \norm{x - y}.
\vspace{1.7mm}
\end{array}
$
\renewcommand{\arraystretch}{1.0}
\addtolength{\arraycolsep}{6pt}
\newline
An identical argument shows that $h\inverse$ is $3K$-Lipschitz.

It is obvious that $h(F) = F$ and that $h(B(0,s_0)) = B(0,s_1)$.
So $h\inverse \nrestriction X$ is as required.

(b) Assume for simplicity that $x_0 = 0$.
By (a), we may assume that $r = s$.
Define $\fnn{g}{[0,\infty)}{[0,1]}$ as follows:
$$
g(t) =
\left\{
\begin{array}{ll}
1                                       & 0 \leq t \leq r, \medskip\\

\frac{r_0 - t}{r_0 - r}\hspace{8mm}   & r \leq t \leq r_0, \medskip\\

0                                       & r_0 \leq t.
\end{array}
\right.
$$
Suppose that $a > r_0$ and $\overB(0,a) \subseteq X$.
We define $\fnn{h}{\overB(0,a)}{E}$ by
$h(x) = x + g(\norm{x}) \cdot v$.
Obviously, $h(B(0,r)) = B(v,r)$.

We show that $h$ is Lipschitz.
At first we see that
$h \nrestriction (\overB(0,r_0) - B(0,r))$ is Lipschitz.
Let $x,y \in \overB(0,r_0) - B(0,r)$. Then
$h(x) - h(y) = x - y + \frac{\norm{y} - \norm{x}}{r_0 - r} \cdot v$.
It follows that
$$\hbox{$\norm{h(x) - h(y)} \leq
\norm{x - y} + 
\frac{\abs{\kern1.5pt\norm{y} - \norm{x}\kern1.5pt}}{r_0 - r} \cdot
\norm{v} \leq
\norm{x - y} +
\frac{\norm{x - y}}{r_0 - r} \cdot \norm{v} =
(1 + \frac{\norm{v}}{r_0 - r}) \cdot \norm{x - y}$.}
$$
Let $x,y \in \overB(0,a)$. If $x,y \in B(0,r)$ or
$x,y \in \overB(0,r_0) - B(0,r)$ or
$x,y \in \overB(0,a) - B(0,r_0)$, then
$\norm{h(x) - h(y)} \leq
(1 + \frac{\norm{v}}{r_0 - r}) \cdot \norm{x - y}$.

If $x \in B(0,r)$ and $y \in \overB(0,r_0) - B(0,r)$, let
$z \in [x,y] \cap S(0,r)$.
Then
\vspace{1.5mm}
\newline
\rule{7pt}{0pt}
\renewcommand{\arraystretch}{1.5}
\addtolength{\arraycolsep}{-6pt}
$
\begin{array}{ll}
&
\norm{h(x) - h(y)} \leq
\norm{h(x) - h(z)} + \norm{h(z) - h(y)} \leq
\norm{x - z} + (1 + \frac{\norm{v}}{r_0 - r}) \cdot \norm{z - y}
\\
\leq
\rule{5pt}{0pt}
&
(1 + \frac{\norm{v}}{r_0 - r}) \cdot (\norm{x - z} + \norm{z - y}) =
(1 + \frac{\norm{v}}{r_0 - r}) \cdot \norm{x - y}.
\vspace{1.7mm}
\end{array}
$
\renewcommand{\arraystretch}{1.0}
\addtolength{\arraycolsep}{6pt}
\newline
The other cases are dealt with similarly. So $h$ is
$(1 + \frac{\norm{v}}{r_0 - r})$-Lipschitz.

In order to show that $h$ is $\onetoonen$ and that $h\inverse$ is
Lipschitz,
we first check that there is $K$ such that
$\norm{x - y} \leq K \ncdot \norm{h(x) - h(y)}$
for every $x,y \in \overB(0,r_0) - B(0,r)$. Indeed,
\vspace{1.5mm}
\newline
\rule{7pt}{0pt}
\renewcommand{\arraystretch}{1.5}
\addtolength{\arraycolsep}{-6pt}
$
\begin{array}{ll}
&
\norm{h(x) - h(y)} \geq \norm{x - y} -
\frac{\abs{\kern1.5pt\norm{y} - \norm{x}\kern1.5pt}}{r_0 - r} \cdot
\norm{v} \geq
\norm{x - y} - \frac{\norm{y - x}}{r_0 - r} \cdot \norm{v}
\\
=
\rule{5pt}{0pt}
&
(1 - \frac{\norm{v}}{r_0 - r}) \cdot \norm{x - y} =
\frac{r_0 - r - \norm{v}}{r_0 - r} \cdot \norm{x - y}.
\vspace{1.7mm}
\end{array}
$
\renewcommand{\arraystretch}{1.0}
\addtolength{\arraycolsep}{6pt}
\newline
Clearly, $\frac{r_0 - r - \norm{v}}{r_0 - r} > 0$.
Let $K = \frac{r_0 - r}{r_0 - r - \norm{v}}$.
Then $\norm{x - y} \leq K \cdot \norm{h(x) - h(y)}$.
This implies that $h \nrestriction (\overB(0,r_0) - B(0,r))$ is
$\onetoonen$.

We next check that
$h(\overB(0,r_0) - B(0,r)) = \overB(0,r_0) - B(v,r)$.
Let
$x \in \overB(0,r_0) - B(0,r)$. There are
$x_1,x_2 \in \rfs{bd}(\overB(0,r_0) - B(0,r))$ such that
$x \in [x_1,x_2] \subseteq \overB(0,r_0) - B(0,r)$,
and $x_2 = x_1 + \lambda v$ for some $\lambda \geq 0$.
Suppose first that $x_1,x_2 \in S(0,r_0)$. Clearly, $h([x_1,x_2])$
is a line segment. Since $h \nrestriction [x_1,x_2]$ is $\onetoonen$
and $h(x_i) = x_i, \ i = 1,2$, we have $h([x_1,x_2]) = [x_1,x_2]$.

A similar argument shows that if $x_1 \in S(0,r_0)$ and
$x_2 \in S(0,r)$,
then
$h([x_1,x_2]) = [x_1,x_2 + v] \subseteq \overB(0,r_0) - B(0,r)$.
Also if $x_1 \in S(0,r)$ and $x_2 \in S(0,r_0)$,
then
$h([x_1,x_2]) = [x_1 + v,x_2] \subseteq \overB(0,r_0) - B(0,r)$.

It follows that
$h(\overB(0,r_0) - B(0,r)) \subseteq \overB(0,r_0) - B(v,r)$.
A similar consideration shows that
$\overB(0,r_0) - B(v,r) \subseteq h(\overB(0,r_0) - B(0,r))$.
Also, $h(B(0,r)) = B(0,v)$,
$h(\overB(0,a) - \overB(0,r_0)) = \overB(0,a) - \overB(0,r_0)$ and
$h \nrestriction ((\overB(0,a) - \overB(0,r_0) \cup B(0,r))$ is
$\onetoonen$. So $h$ is a bijection and
\hbox{$\rfs{Rng}(h) = \overB(0,a)$.}
We have proved that
$h\inverse \nrestriction (\overB(0,r_0) - B(0,r))$ is
$\frac{r_0 - r}{r_0 - r - \norm{v}}$\,-\,Lipschitz.
The argument that $h\inverse$ is
$\frac{r_0 - r}{r_0 - r - \norm{v}}$\,-\,Lipschitz,
is the same one used in
showing that $h$ is Lipschitz.

Clearly,
$\frac{r_0 - r}{r_0 - r - \norm{v}} =
1 + \frac{\norm{v}}{r_0 - r - \norm{v}}$
and
$1 + \frac{\norm{v}}{r_0 - r} \leq
1 + \frac{\norm{v}}{r_0 - r - \norm{v}}$.
So $h$ is $1 + \frac{\norm{v}}{r_0 - r - \norm{v}}$-bilipschitz.
As in the preceding arguments, this implies that
$h \cup \rfs{Id} \nrestriction (X - \overB(0,a))$ is
$1 + \frac{\norm{v}}{r_0 - r - \norm{v}}$-bilipschitz.

For every $x \in \overB(0,a)$,
$h(x) - x \in \rfs{span}(\sngltn{v}) \subseteq F$.
So $x \in F$ iff $h(x) \in F$.
Hence $h \cup \rfs{Id} \nrestriction (X - \overB(0,a)) \in
\rfs{LIP}(X;\calF,\calS)$.
Note also that
$h \nrestriction B(0,r) = \rfs{tr}_v \nrestriction B(0,r)$.
So $h \cup \rfs{Id} \nrestriction (X - \overB(0,a))$
fulfills the requirements of (i) and (ii).

(c) Let $x_0, \ldots, x_n \in [x,y]$ be such that $x_0 = x$,
$x_1 = y$ and
$\norm{x_i - x_{i + 1}} < \dgfrac{r}{4}$ for every $i < n$.
By (b), for every $i < n$ there is
$h_i \in \sprtdl{\rfs{LIP}(X;\calS,\calF)}{B(x_i,r)}$ such that
$h_i \nrestriction B(x_i,\dgfrac{2r}{3})) =
\rfs{tr}_{x_{i + 1} - x_i} \nrestriction B(x_i,\dgfrac{2r}{3})$.
Let $h = h_0 \scirc \ldots \scirc h_{n - 1}$. Then $h$ is as required.
Note that $n$ can be chosen to be $[4 \dgfrac{\norm{x - y}}{r]} + 1$.
By (b) each $h_i$ is
$(1 + \frac{\dgfrac{r}{4}}{r - \dgfrac{2r}{3}
- \dgfrac{r}{4}})$-bilipschitz.
That is, $h_i$ is $4$-bilipschitz. Hence
$K_{\srfs{seg}}(\ell,t) = 4^{[\frac{4 \ell}{t}] + 1}$.

(d) Let $x = \gamma(s)$, $y = \gamma(1)$, $L = \gamma([s,1])$
and $r = d(L,(X - U) \cup \sngltn{s(0)})$.
There\break
is a sequence of balls $B(x_1,r),\ldots,B(x_n,r)$ such that
$x_1,\ldots,x_n \in L$ and
$\bigcup_{i = 1}^n B(x_i,r) \supseteq L$.
We may assume that
{\thickmuskip=2mu \medmuskip=1mu \thinmuskip=1mu 
$x \in B(x_1,r)$,
$y \in B(x_n,r)$,
and $B(x_i,r) \cap B(x_{i + 1},r) \neq \emptyset$ for every $i < n$.
}
For every $i < n$ let $y_i \in B(x_i,r) \cap B(x_{i + 1},r)$.
Set $y_0 = x$ and $y_n = y$.
By (b), for every $i = 1,\ldots,n$
there is $h_i \in \rfs{LIP}(X)$ such that $h_i(y_{i - 1}) = y_i$
and $\rfs{supp}(h_i) \subseteq B(x_i,r)$.
Clearly, $h_n \scirc \ldots \scirc h_1$ is as required.
\smallskip\hfill\myqed

The following observation will be used in many arguments.
Its proof is left to the reader.

\begin{prop}\label{p2.6.0}\label{p4.7}
\num{a} Let $X$ be a metric space, and $\vecx$ be a sequence in $X$.
Then either $\vecx$ has a Cauchy subsequence,
or there are $r > 0$ and a subsequence $\setm{y_n}{n \in \bbN}$ of
$\vecx$ such that for distinct $i,j \in \bbN$,
$d(y_i,y_j) \geq r$.

\num{b}
Let $X$ be a metric space
and $\setm{x_i}{i \in \bbN} \subseteq X$ be a bounded sequence.
Then either $\setm{x_i}{i \in \bbN}$ has a Cauchy subsequence,
or there is a subsequence
$\setm{y_i}{i \in \bbN}$ of $\setm{x_i}{i \in \bbN}$ and $r > 0$
such that for every $\varepsilon > 0$ there is $N \in \bbN$
such that $\abs{d(y_i,y_j) - r} < \varepsilon$
for distinct $n,m > N$.
\end{prop}

\begin{prop}\label{p2.6.1}
\num{a} If \,$U_1 \subseteq U \prec V \subseteq V_1$,
then $U_1 \prec V_1$.

\num{b} If \,$U \prec V$ for some $V$, then $U$ is small.

\num{c} Let $B(x,r)$ and $B(y,s)$
be manageable balls based on the same $S$.
If
$\rfs{cl}(B(x,r)) \subseteq B(y,s)$, then $B(x,r) \prec B(y,s)$.

\num{d} If $U \in \rfs{Ro}(X)$ is a subset of a manageable ball,
then $U$ is small.

\num{e} If $U \prec V$, then $\rfs{cl}(U) \subseteq V$.

\num{f} If $B$ is a manageable ball of $X$, then $B \in \rfs{Ro}(X)$
and $B$ is small.
\end{prop}

\noindent
{\bf Proof }
Parts (a) and (b) follow trivially from the definitions.

(c) Note that if $\rfs{cl}(B(x,r)) \subseteq B(y,s)$,
then $\norm{x - y} + r < s$.
So (c) follows from Lemma \ref{l2.6}(b).

(d) Suppose that $U \subseteq B$, and $B$ is a manageable ball.
There is a mangable ball $B'$ with the same center as $B$ such that
$\rfs{cl}(B) \subseteq B'$.
Obviously, $B$ and $B'$ are based on the same
$S$. So by (c), $B \prec B'$. By (a), $U \prec B'$.
By (b), $U$ is small.

(e) Suppose that $x \in \rfs{cl}(U) - V$.
Let $\emptyset \neq W \subseteq V$.
Then there is $\emptyset \neq W' \subseteq W$
such that $\rfs{cl}(W') \subseteq W$.
Let $g \in \sprtd{G}{V}$. Then $g(x) = x$.
Suppose by contradiction that $g(U) \subseteq W'$.
Then
$g(x) \in g(\rfs{cl}(U)) \subseteq \rfs{cl}(W') \subseteq W \not\ni x$.
A contradiction.

(f) $B \in \rfs{Ro}(E)$ and
$\rfs{int}(\rfs{cl}(B)) = \rfs{int}^E(\rfs{cl}^E(B))$.
So $B \in \rfs{Ro}(X)$.
\vspace{1.5mm}\hfill\myqed

Let $\calU \subseteq \rfs{Ro}(X)$.
We use $\sum\,\, \calU$ to denote the supremum of $\calU$
in the complete Boolean algebra $\rfs{Ro}(X)$.
It is easy to check that
$\sum\,\, \calU = \rfs{int}(\rfs{cl}(\bigcup\,\, \calU))$.

\begin{defn}\label{d2.16}
\begin{rm}
(a) Let $U \subseteq V$ and $\calU \subseteq \rfs{Ro}(X)$.
$\calU$ is called a {\it $V$-small semicover} of $U$,
if $\,\sum\,\, \calU = U$
   \index{small semicover. $V$-small semicover}
and \,$U' \prec V$ for every $U'\in\calU$.

(b) Let $\calU$ be a $V$-small semicover of $U$, and let
$\setm{U_i}{i \in I}$ be a
$\onetoonen$ enumeration of~$\calU$. We say that $\calU$ is
a {\it $V$-good semicover} of $U$, if the following holds.
   \index{good semicover. $V$-good semicover}
For every $J \subseteq I$
and $\setm{W_j}{j \in J} \subseteq \rfs{Ro}(X)$:
if $J$ is infinite and $\emptyset \neq W_j \subseteq U_j$
for every $j \in J$,
then there are pairwise disjoint infinite
$J_1,J_2 \subseteq J$ and
$\setm{W'_j}{j \in J_1 \cup J_2} \subseteq \rfs{Ro}(X)$
such that $\emptyset \neq W'_j \subseteq W_j$
for every $j \in J_1 \cup J_2$
and 
$\sum_{j \in J_1} W'_j \spr \sum_{j \in J_2} W'_j$.

(c) For a normed vector space $E$ let $\overE$ denote the completion
of $E$. So $\overE$ is a Banach space.
   \index{N@AAAA@@$\rule{1.0pt}{0pt}\overE$.
          The completion of a normed vector space $E$}

(d) Let $Z$ be a topological space.
Suppose that $F \subseteq H(Z)$ and
$\rfs{supp}(f) \cap \rfs{supp}(g) = \emptyset$
for distinct $f,g \in F$.
We define
$$
\mbox{\small
$\bcirc \kern-1.7pt F \eqdf
\bigcup \setm{f \nrestriction \rfs{supp}(f)}{f \in F} \cup
\rfs{Id} \nrestriction (Z - \bigcup \setm{\rfs{supp}(f)}{f \in F}).$
}
$$
Let $F = \setm{f_n}{n \in \bbN} \subseteq H(Z)$ be such that for every
distinct \hbox{$m,n \in \bbN$,}
Then $\bcirc_{n \in \sboldbbN}f_n \eqdf \bcirc\kern-1.7pt F$.

   \index{N@AAAA@@$\bcirc \kern-1ptF =
          \bigcup \setm{f \nrestriction \rfs{supp}(f)}{f \in F} \cup
          \rfs{Id} \nrestriction
	  (Z - \bigcup \setm{\rfs{supp}(f)}{f \in F})$}
   \index{N@AAAA@@$\bcirc_{n \in \sboldbbN} h_n$}
\end{rm}
\end{defn}

\begin{lemma}\label{l2.8} Let $V$ be a small set. Then for every
$U \in \rfs{Ro}(X)$:\ $\rfs{cl}(U) \subseteq V$
iff \,$U$ has a $V$-good semicover.
\end{lemma}

\noindent
{\bf Proof }
Suppose that $\rfs{cl}(U) \not\subseteq V$.
Let $\calU$ be a $V$-small semicover of $U$;
we show that $\calU$ is not $V$-good.
The fact that $V$ is small is not used in the proof of this direction.
Let $x \in \rfs{cl}(U) - V$. If $U' \in \calU$,
then by \ref{p2.6.1}(e), $\rfs{cl}(U') \subseteq V$.
By induction on $i \in \bbN$ we define
$U_i \in \calU$ and \hbox{$W_i \subseteq U_i$.}
Let $U_0$ be any member of $\calU$ and $W_0 = U_0$.
Suppose $U_0,\ldots,U_{i - 1}$ and $W_0,\ldots,W_{i - 1}$
have been defined. Let $B_i$ be a ball with center at $x$ and radius
$< \dgfrac{1}{i}$ such that $B_i \cap \bigcup_{j < i} U_j = \emptyset$.
Let $U_i \in \calU$ be such that $B_i \cap U_i \neq \emptyset$,
and let $W_i = U_i \cap \rfs{int}(\rfs{cl}(B_i))$.
So $W_i \in \rfs{Ro}(X)$.
For every infinite $J'\sbeq \bbN$ and
$\setm{W'_j}{j \in J'} \subseteq \rfs{Ro}(X)$:
if $\emptyset \neq W_j' \subseteq W_j$
for every $j \in J'$,
then $x \in \rfs{cl} \left(\sum_{j \in J'} W'_j \right)$.
Suppose by contradiction that $\calU$ is $V$-good.
The family $\setm{U_i}{i \in \bbN}$ is an infinite subset of $\,\calU$,
and $W_i \subseteq U_i$ for every $i \in \bbN$.
So let $J_1,J_2$ and $\setm{W'_j}{j \in J_1 \cup J_2}$
be as required in the definition of $V$-goodness
for $\setm{U_i}{i \in \bbN}$ and $\setm{W_i}{i \in \bbN}$,
and let $W$ strongly separate $\sum_{j \in J_1} W'_j$ from
$\sum_{j \in J_2} W'_j$.
Since $x \in \rfs{cl}\left(\sum_{j \in J_2} W'_j \right)$ and
$W \cap \sum_{j \in J_2} W'_j = \emptyset$,
it follows that $x \not\in W$.
But $x \in \rfs{cl}\left(\sum_{j \in J_1} W'_j \right)$.
So by \ref{p2.6.1}(e), $\sum_{j \in J_1} W'_j \not\prec W$.
A contradiction.

Assume next that $V$ is small and that $\rfs{cl}(U)\sbeq V$;
we will construct a $V$-good semicover $\calU$ of $U$.
Since $V$ is small,
there is $g \in G$ such that $g(V)$ is contained in a manageable ball.
Obviously $\rfs{cl}(g(U)) \subseteq g(V)$.
Clearly, $g(U)$ has a $g(V)$-good semicover iff 
$U$ has a $V$-good semicover.
In fact, this follows from Proposition \ref{np2.12}(a).
We may thus assume that $V$ is contained in a manageable ball.
This means that $\rfs{cl}(U) = \rfs{cl}^E(U)$.

We may further assume that there is a manageable ball
$B^{\rstar} = B^E(x^*,r^*)$ such that
$V \subseteq B^E(x^*,\dgfrac{r^*}{16})$.
Suppose that $B^{\rstar}$ is based on $S^{\rstar}$,
and denote $F_{S^{\raisedsstar}}$ by $F^{\rstar}$.
We may assume that $x^* \in F^{\rstar}$.
For every $x \in \rfs{cl}(U)$ let $W_x \in \rfs{Ro}(X)$ be such that
$x \in W_x \prec V$.
The existence of $W_x$ follows from Proposition \ref{p2.6.1} (c), (a)
and (f).
Since $\rfs{cl}(U)$ is paracompact,
there is an open locally finite refinement
$\calT$ of
$\setm{W_x}{x \in \rfs{cl}(U)}$ such that
$\rfs{cl}(U) \subseteq \bigcup \calT$.
Let $\calU = \setm{\rfs{int}(\rfs{cl}(T)) \cap U}{T \in \calT}$.
By Proposition \ref{np2.12}(b)(ii),
$\calU \subseteq \rfs{Ro}(X)$.
Clearly, $\bigcup \calU = U$. So $\sum \calU = U$.

We show that for every $x \in \rfs{cl}(U)$ there is a neighborhood $S_x$
such that\break
$\setm{U' \in \calU}{U' \cap S_x \neq \emptyset}$ is finite.
For $x \in \rfs{cl}(U)$ let $S_x$ be an open neighborhood of $x$
such that $\setm{T \in \calT}{T \cap S_x \neq \emptyset}$ is finite.
By Proposition \ref{np2.12}(b)(iv),
$\setm{T \in \calT}{\rfs{int}(\rfs{cl}(T)) \cap S_x \neq \emptyset}$
is finite.
So
$\setm{T \in \calT}
{(\rfs{int}(\rfs{cl}(T)) \cap U) \cap S_x \neq \emptyset}$
is finite.
That is, $\setm{U' \in \calU}{U' \cap S_x \neq \emptyset}$ is finite.

We show that $\calU$ is $V$-good.
Let $\setm{U_i}{i\in \bbN} \subseteq \calU$ be such that
$U_i\neq U_j$ for every $i \neq j$;
and let $\empt\neq W_i\sbeq U_i$.
We shall find $J_1,\,J_2$ and $\setm{W'_j}{j \in J_1 \cup J_2}$
as required in the definition of $V$-goodness.
For every $i \in \bbN$ let $x_i \in W_i \cap F^{\rstar}$.

{\bf Claim 1.} $\setm{x_i}{i \in \bbN}$ does not contain
a convergent subsequence.

{\bf Proof } If $x$ is a limit of such a subsequence,
then $x \in \rfs{cl}(U)$, but then $S_x$ intersects only finitely many
$U_i$'s. So $\setm{i \in \bbN}{x_i \in S_x}$ is finite.
A contradiction, so the claim is proved.
\smallskip

By Claim 1 and Proposition \ref{p2.6.0}(b),
either (i) or (ii) below happen:
\indent
\begin{itemize}
\item[(i)] $E$ is incomplete,
there is an infinite $J \subseteq \bbN$ such that
$\setm{x_i}{i \in J}$ is a Cauchy sequence,
and $\setm{x_i}{i \in J}$ is not convergent in $\rfs{cl}^E(X)$.
\item[(ii)] There is infinite $J \subseteq \bbN$ and an $r > 0$
such that for every  distinct $i,j \in J$,\break
$r < \norm{x_i - x_j} < \dgfrac{9r}{8}$.
\end{itemize}

{\bf Case (i) } Assume that (i) happens.
Let $\barx = \lim^{\oversE}_{i \in J} x_i$.
Hence $\barx \in \rfs{cl}^{\oversE}(V) - X$.
Since $V \subseteq B^E(x^*,\dgfrac{r^*}{16})$,
there is $r > 0$ such that
$B^{\oversE}(\barx,r) \cap E \subseteq B^E(x^*,\dgfrac{r^*}{8})$.
So $\barx \not\in E$.
We may assume that $x_i \in B^{\oversE}(\barx,\dgfrac{r}{8})$
for every $i \in J$.
Let $v \in F^{\rstar}$ and $\norm{v} = \dgfrac{r}{2}$.
Let $L_i = [x_i,x_i + v]$ and $L = [\barx,\barx + v]$.
So $L_i \subseteq B^{F^{\raisedsstar\kern-2pt}}(x^*,\dgfrac{r^*}{8})$
for every $i \in J$.
Also, $L \subseteq \overE - E$.
{\thickmuskip=2mu \medmuskip=1mu \thinmuskip=1mu 
One can choose an infinite subset $J_0 \subseteq J$
and a sequence $\setm{r_i}{i \in J_0} \subseteq (0,\dgfrac{r}{8})$
}
such that \hbox{$B^E(x_i,r_i) \subseteq W_i$}
for every $i \in J_0$,
and $\rfs{cl}^E(B(L_i,r_i)) \cap \rfs{cl}^E(B(L_j,r_j)) = \emptyset$
for distinct $i,j \in J_0$.

For every $i \in J_0$ let $W'_i = B(x_i,\dgfrac{r_i}{3})$.
Let $J_1 \subseteq J_0$ be such that $J_1$ and $J_0 - J_1$ are infinite,
and let $J_2 = J_0 - J_1$.
For $\ell = 1,2$ let $W^{\ell} = \sum_{i \in J_{\ell}} W'_i$.
We shall show that $W^1 \spr W^2$.

For every $i \in J_1$,
$\norm{x_i - \barx} < \dgfrac{r}{8}$ and $r_i < \dgfrac{r}{8}$,
and for every $u \in L_i$,
$\norm{u - x_i} \leq \norm{(x_i + v) - x_i} = \dgfrac{r}{2}$.
It follows that for every $u \in B(L_i,r_i)$,
$\norm{u - \barx} < \frac{r}{8} + \frac{r}{2} + \frac{r}{8} =
\frac{3r}{4}$.
So
$B(L_i,r_i) \subseteq
B(\barx,r) \subseteq B(x^*,r^*) \subseteq S^{\rstar}$.

By Lemma \ref{l2.6}(c), for every $i \in J_1$ there is
$h_i \in \sprtdl{\rfs{LIP}(X;\calS,\calF)}{B(L_i,r_i)}$
such that\break
$h_i(B(x_i,\dgfrac{r_i}{3})) = B(x_i + v,\dgfrac{r_i}{3})$.
Let $h = \bcirc_{i \in J_1\kern-2.5pt} h_i$.
We show that
$h \in \rfs{LIP}^{\srfs{LC}}(X;\calS,\calF)$.
Clearly,
$\rfs{supp}(h) =
\bigcup_{i \in J_1} \rfs{supp}(h_i) \subseteq S^{\rstar}$.
We show that for every $u \in E$, there is a neighborhood $V_u$ of $u$
such that
$\abs{\setm{i \in J_1}{B(L_i,r_i) \cap V_u \neq \emptyset}} \leq 1$.
Suppose that $u$ is a counter-example.
Since $\setm{x_i}{i \in \bbN}$ is a Cauchy sequence and the
$B(x_i,r_i)$'s are pairwise disjoint,
$\lim_i r_i = 0$.
Since for $i \neq j$,
$\rfs{cl}^E(B(L_i,r_i)) \cap \rfs{cl}^E(B(L_j,r_j)) = \emptyset$,
there is at most one $i$ such that
$u \in \rfs{cl}^E(B(L_i,r_i))$.
Hence there is an infinite set $J_3 \subseteq J_1$
and a sequence $\setm{u_i}{i \in J_3}$ such that $u_i \in B(L_i,r_i)$
for every $i \in J_3$,
and $\lim_{i \in J_3} u_i = u$.
There are $y_i \in L_i$ such that $\norm{y_i - u_i} < r_i$.
Hence $\lim_{i \in J_3} y_i = u$.
Let $y_i = x_i + t_i v$.
Since $\sngltn{x_i}$ and $\sngltn{y_i}$ converge in $\overE$,
$\lim_{i \in J_3} t_i$ exists.
Also, $\lim_{i \in J_3} t_i \in [0,1]$.
So $u \in [\barx,\barx + v]$. Hence $u \not\in E$, a contradiction.

Let $u \in X$. Then  there is $i \in J_1$ such that
$h \nrestriction V_u = h_i \nrestriction V_u$.
So $h \nrestriction V_u$ is bilipschitz.
This means that $h \in \rfs{LIP}^{\srfs{LC}}(X;\calS,\calF)$.
Since $E$ is incomplete,
$\rfs{LIP}^{\srfs{LC}}(X;\calS,\calF) \subseteq G$. So $h \in G$.

We shall prove that $h(W^1) \spr W^2$.
Let us first see that
$h(W^1) \subseteq B^{\oversE}(\barx + v,\dgfrac{r}{6})$.
We have
$h(W^1) = \bigcup_{i \in J_1} h_i(W'_i) =
\bigcup_{i \in J_1} h_i(B(x_i,\dgfrac{r_i}{3})) =
\bigcup_{i \in J_1} B(x_i + v,\dgfrac{r_i}{3})$.
Also,\break
$\norm{(x_i + v) - (\barx + v)} = \norm{x_i - \barx} < \dgfrac{r}{8}$.
Since $\norm{x_i - \barx} < \dgfrac{r}{8}$
and $\barx \not\in B^{\oversE}(x_i,r_i)$,
It follows that $r_i < \dgfrac{r}{8}$.
So
$B(x_i + v,\dgfrac{r_i}{3}) \subseteq
B^{\oversE}(\barx + v,\dgfrac{r}{6})$.
That is, $h(W'_i) \subseteq B^{\oversE}(\barx + v,\dgfrac{r}{6})$.
Hence $h(W^1) \subseteq B^{\oversE}(\barx + v,\dgfrac{r}{6})$.

Similarly, $W^2 \subseteq B^{\oversE}(\barx,\dgfrac{r}{6})$.
Since $\norm{(\barx + v) - \barx} = \dgfrac{r}{2} > \dgfrac{r}{3}$,
there are $\hatx \in E$ and $0 < s_0 < s_1$ such that
$B^{\oversE}(\barx + v,\dgfrac{r}{6}) \subseteq B^{\oversE}(\hatx,s_0)$
and $B^{\oversE}(\hatx,s_1) \cap B^{\oversE}(\barx,\dgfrac{r}{6}) =
\emptyset$.
So
$h(W^1) \subseteq B^E(\hatx,s_0)$.
By Propo\-sitions \ref{p2.6.1}(c) and \ref{p2.6.1}(a),
$h(W^1) \prec B^E(\hatx,s_1)$.
Since $B^E(\hatx,s_1)  \cap W^2 = \emptyset$,
it follows that $h(W^1) \spr W^2$.

Note that $h(W^2) = W^2$.
By Proposition \ref{np2.12}(a),
$h\inverse(h(W^1)) \spr h\inverse(W^2)$.
But
$h^{-1}(h(W^1)) = W^1$
and $W^2 = h\inverse(h(W^2))$.
So $W^1 \spr W^2$.

\medskip

{\bf Case (ii) } Assume that (ii) happens.
Since the $x_i$'s belong to $B^E(x^*,\dgfrac{r^*}{16})$
and $r < \norm{x_i - x_j}$, it follows that $r < \dgfrac{r^*}{8}$.
Let $i_0 \in J$
and $J_1$ and $J_2$ be disjoint infinite subsets of $J$
not containing~$i_0$. For every $i \in J_1 \cup J_2$
let $B_i = B^E(x_i,\dgfrac{r}{8})$ and $W'_i = B_i \cap W_i$.
Clearly, $B_i \subseteq B^E(x^*,\dgfrac{3 r^*}{16})$.
So $B_i \subseteq X$, and hence $W'_i \in \rfs{Ro}(X)$.
For $\ell = 1,2$ let $W{^\ell} = \sum_{i \in J_{\ell}}W'_i$,
and let $W = B(x_{i_0},2r)$.

We shall show that:
($*$) There is $h \in \rfs{LIP}(E;B^E(x_{i_0},3r),F^{\rstar})$
such that
$h \nrestriction W^1 = \rfs{Id}$ and $h(W^2) \cap W = \emptyset$.

But at first we prove that ($*$) implies that $W^1 \spr W^2$.
If $x \in B^E(x_{i_0},3r)$,
then
$\norm{x - x^*} \leq \norm{x - x^*} + 3r <
\dgfrac{r^*}{16} + \dgfrac{3 r^*}{8} = \dgfrac{7 r^*}{16}$.
So $B^E(x_{i_0},3r) \subseteq B^{\rstar} \subseteq S^{\rstar}$.
Hence
$h \nrestriction X \in \rfs{LIP}(X;S^{\rstar},F^{\rstar}) \subseteq
\rfs{LIP}(X;\calS,\calF)$.
Now, $W^1 \subseteq B(0,\dgfrac{5r}{4})$,
so by \ref{p2.6.1}(c) and \ref{p2.6.1}(a),
$W^1 \prec W$.
Also
$h(W^2) \cap W = \emptyset$.
Hence $W$ strongly separates $W^1$ from $h(W^2)$. 
That is, $W^1 \spr h(W^2)$.
By Proposition \ref{np2.12}(a),
$h\inverse(W^1) \spr h\inverse(h(W^2))$.
But
$h^{-1}(W^1) = W^1$
and $W^2 = h\inverse(h(W^2))$.
So $W^1 \spr W^2$.

To complete the proof, it remains to show that ($*$) holds.
For simplicity let us assume that $x_{i_{0}} = 0$
and that $r = 1$.
We define a function $\fnn{g}{[0,3] \times [0,\infty)}{\bbR}$
as follows.
For every $s_0 \in [0,3]$, $g(s_0,t)$ will be a piecewise linear
homeomorphism of $[0,\infty]$.
Let $a(s)$ be the linear function such that
$a(\dgfrac{3}{8}) = \dgfrac{3}{4}$ and $a(\dgfrac{5}{8}) = 2$.

If $s_0 \leq \dgfrac{3}{8}$, then $g(s_0,t) = t$.
If $\dgfrac{3}{8} \leq s_0 \leq \dgfrac{5}{8}$, then
$$
g(s_0,t) =
\left\{
\begin{array}{ll}
t                     & t \leq \half,                        \medskip\\
\change{ \frac{a(s_0)-\half}{\frac{3}{4}-\half}}(t-\half)+\half,
\hspace{5mm}
                       & \half \leq t \leq \frac{3}{4}       \medskip\\
\change{ \frac{3-a(s_0)}{3-\frac{3}{4}}}(t-3)+3,
		       & \frac{3}{4} \leq t \leq 3           \medskip\\
t                     & 3 \leq t.                            \medskip\\
\end{array}
\right.
$$
If $\dgfrac{5}{8} \leq s_0 \leq 3$,
then $g(s_0,t) = g(\dgfrac{5}{8},t)$.

Let $F = \setm{x_i}{i\in J_1}$ and
$$
h(x) =
\left\{
\begin{array}{ll}
g\left( d(\frac{x}{\dline{x}},F),\dline{x}\right)
\cdot \frac{x}{\dline{x}} \hspace{8mm}
					      & x \neq 0,   \medskip\\
0                                             & x   =  0.
\end{array}
\right.
$$
We leave it to the reader to check that
$h \in \rfs{LIP}(E;B^E(0,3),F_{S^*})$.

If $i \in J_1 \cup J_2$ and $x \in B_i$, then
$\norm{\frac{x}{\norm{x}} - x_i} \leq
\norm{\frac{x}{\norm{x}} - x} + \norm{x - x_i} <
\dgfrac{1}{4} + \dgfrac{1}{8} =  \dgfrac{3}{8}$.
Let $x \in W^1$. 
There is $i \in J_1$ such that $x \in B_i$.
Hence
$d \left(\frac{x}{\norm{x}},F \right) \leq
\norm{\frac{x}{\norm{x}} - x_i} < \dgfrac{3}{8}$.
So $g \left(d(\frac{x}{\norm{x}},F),\norm{x} \right) = \norm{x}$,
and hence $h(x) = x$.
Let $x \in W^2$. 
There is $i \in J_2$ such that $x \in B_i$.
So
$d \left(\frac{x}{\norm{x}},F \right) \geq
d(x_i, F) - \norm{\frac{x}{\norm{x}} - x_i} >
1 - \dgfrac{3}{8} = \dgfrac{5}{8}$.
Also, $\dline{x} > \dgfrac{7}{8}$.
Hence
$$\hbox{
$\norm{h(x)} =
\left\| g \left(d(\frac{x}{\norm{x}},F),\norm{x} \right) \cdot
\frac{x}{\norm{x}} \right\| =
g \left(d(\frac{x}{\norm{x}},F),\norm{x} \right) =
g(\dgfrac{5}{8},\norm{x}) > g(\dgfrac{5}{8},\dgfrac{3}{4}) = 2$.
}$$
We have proved ($*$),
so the proof of the lemma is complete. \hfill\myqed
\vspace{-0.0mm}

\begin{lemma} \label{l2.9}
Let $V$ be a small set. then for every $U$:
$\rfs{cl}(U) \cap \rfs{cl}(V) \neq \emptyset$ iff for every small $V_1$:
if $\rfs{cl}(V) \subseteq V_1$, then $V_1 \cap U \neq \emptyset$.
\end{lemma}

\noindent
{\bf Proof }  If $\rfs{cl}(U) \cap \rfs{cl}(V) \neq \emptyset$,
then clearly $V_1 \cap U \neq \emptyset$
for every $V_1 \supseteq \rfs{cl}(V)$.
Conversely, suppose that $V$ is small
and $\rfs{cl}(V) \cap \rfs{cl}(U) = \emptyset$.
Let $V'$ be a small set such that $\rfs{cl}(V) \subseteq V'$,
and let $V_1 = V' \cap \rfs{int}(X - U)$.
Since $\rfs{int}(X - U) \supseteq \rfs{cl}(V),\ V_1 \supseteq
\rfs{cl}(V)$,
hence $V_1$ is as required.
\rule{0pt}{0pt}\hfill\myqed

\begin{lemma}\label{l2.10}
Let $U$ and $V$ be small sets. Then
$|\rfs{cl}(U) \cap \rfs{cl}(V)| = 1$ iff the following holds.
\begin{itemize}
\addtolength{\parskip}{-11pt}
\addtolength{\itemsep}{06pt}
\item[\num{i}]
$\rfs{cl}(U) \cap \rfs{cl}(V) \neq \emptyset$,
\item[\num{ii}]
for every small $W_1$ and $W_2$:
if $\rfs{cl}(U \cap W_1) \cap \rfs{cl}(V \cap W_1) \neq \emptyset$ and
$\rfs{cl}(U \cap W_2) \cap\break
\rfs{cl}(V \cap W_2) \neq \emptyset$, then
$\rfs{cl}(W_1) \cap \rfs{cl}(W_2) \neq \emptyset$.
\vspace{-05.7pt}
\end{itemize}
\end{lemma}

\noindent
{\bf Proof }
Suppose that $x_1,x_2 \in \rfs{cl}(U) \cap \rfs{cl}(V)$
and $x_1 \neq x_2$.
For $i = 1,2$ let $W_i \in \rfs{Ro}(X)$ be a neighborhood of $x_i$
such that
$W_i$ is small
and $W_i \subseteq B^X(x_i,\frac{1}{3}\norm{x_2 - x_1})$.
Then
$\rfs{cl}(U \cap W_i) \cap \rfs{cl}(V \cap W_i) \neq \emptyset$
for $i = 1,2$,
but $\rfs{cl}(W_1) \cap \rfs{cl}(W_2) = \emptyset$.

Suppose that $\rfs{cl}(U) \cap \rfs{cl}(V) = \sngltn{x}$
and let $W_i,\ i = 1,2$, be such that
$\rfs{cl}(U \cap W_i) \cap\break \rfs{cl}(V \cap W_i) \neq\emptyset$.
Hence $x \in \rfs{cl}(W_1) \cap \rfs{cl}(W_2)$.\hfill\myqed

\begin{lemma}\label{l2.11}
For $i = 1,2$ let $U_i,V_i$ be small sets such that
$|\rfs{cl}(U_i) \cap \rfs{cl}(V_i)| = 1$.
Then
$\rfs{cl}(U_1) \cap \rfs{cl}(V_1) = \rfs{cl}(U_2) \cap \rfs{cl}(V_2)$
iff $(*)$ 
for any small $W_1,W_2$:
if $\rfs{cl}(U_i \cap W_i) \cap \rfs{cl}(V_i \cap W_i) \neq \emptyset$,
$i = 1,2$,
then $\rfs{cl}(W_1) \cap \rfs{cl}(W_2) \neq \emptyset$.
\end{lemma}

\noindent
{\bf Proof } Similar to \ref{l2.10}.\hfill\myqed

\begin{lemma}\label{l2.12}
Let $U,\,V$ be small sets such that
$\rfs{cl}(U) \cap \rfs{cl}(V) = \sngltn{x}$
and $W \in \rfs{Ro}(X)$. Then $x \in W$ iff
$(*)$ for any small $U',\,V'$:
if $\rfs{cl}(U') \cap \rfs{cl}(V') = \rfs{cl}(U) \cap \rfs{cl}(V)$,
then $U' \cap W \neq \emptyset$.
\end{lemma}

\noindent
{\bf Proof }
It is trivial that if $x \in W$, then ($*$) holds.

Suppose that $x \not\in W$. Since $W$ is regular open,
$x \in \rfs{cl}(X - \rfs{cl}(W))$.
Let $B$ be a manageable ball containing $x$.
So let $\setm{x_i}{i \in \bbN} \subseteq B$
be a $\onetoonen$ sequence converging to $x$
and disjoint from $\rfs{cl}(W)$.
Let
$r_i =
\third
\min(\dgfrac{1}{i},d(x_i,\setm{x_j}{j \neq i} \cup W \cup (X - B)))$.
Let $U' = \bigcup \setm{B^E(x_i,r_i)}{i \mbox{ is odd}}$
and $V' = \bigcup \setm{B^E(x_i,r_i)}{i \mbox{ is even}}$.
Then $U',\,V' \subseteq B \subseteq X$. It is easy to see that
$U',\,V' \in \rfs{Ro}(X)$.
Also, since $U',\,V' \subseteq B$, they are small. We have\break
$\rfs{cl}(U') \cap \rfs{cl}(V') = \sngltn{x} =
\rfs{cl}(U) \cap \rfs{cl}(V)$,
and $U' \cap W = \emptyset$. So ($*$) does not hold.\hfill\myqed

\begin{lemma}\label{l2.13}
For every $x \in X$ there are small $U,V$ such that
$\rfs{cl}(U) \cap \rfs{cl}(V) = \sngltn{x}$.
\end{lemma}

\noindent
{\bf Proof } Use the construction of \ref{l2.12}.
\vspace{1mm}\hfill\myqed

\medskip
\noindent{\bf Proof of Theorem \ref{t2.4} }
Recall that \ref{t2.4}(a) is a special case of \ref{t2.4}(b).
We prove~(b).
Let $\pair{X_1}{G_1},\
\pair{X_2}{G_2} \in K_{\srfs{BNO}}$ and $\iso{\varphi}{G_1}{G_2}$.
It is trivial that $\pair{\rfs{Ro}(X_i)}{G_i}$ are topological
local movement systems.
Indeed, this follows from Lemma \ref{l2.6}(a).
Hence by Theorem \ref{t2.3}, there is
{\thickmuskip=3.9mu \medmuskip=2.9mu \thinmuskip=1.9mu 
$\iso{\eta}{\rfs{Ro}(X_1)}{\rfs{Ro}(X_2)}$ such that
$\iso{(\varphi \cup \eta)}{\rfs{MR}(X_1,G_1)}{\rfs{MR}(X_2,G_2)}$.
}
Let $\psi = \varphi \cup \eta$.

As in Remark \ref{nr2.11}
the property of $\calU$ being a $V$-small
semicover of $U$ is expressed in terms of
the operations of $\rfs{MR}(X,G)$.
That is, there is a formula
$\varphi_{\srfs{sm-sc}}(\calX,x,y)$
expressed in terms of the operations of $\rfs{MR}(X,G)$
such that for every
$\pair{X}{G} \in K_{\srfs{BNO}}$, $\calU \subseteq \rfs{Ro}(X)$
and $U,V \in \rfs{Ro}(X)$, $\trpl{\,\calU}{U}{V}$ satisfies
$\varphi_{\srfs{sm-sc}}(\calX,x,y)$ in $\rfs{MR}(X,G)$
iff $\calU$ is a $V$-small semicover of $U$.
Hence, if $\calU$ is a $V$-small semicover of $U$
in $\rfs{MR}(X_1,G_1)$,
then \hbox{$\psi(\calU) \eqdf \setm{\psi(U')}{U' \in \calU}$}
is a $\psi(V)$-small semicover of $\psi(U)$ in $\rfs{MR}(X_2,G_2)$.

The same fact is true for the property of being a $V$-good semicover.

Lemmas \ref{l2.8} - \ref{l2.12}, and the existence of the formulas
$\chi_{\tSml}$ etc.\ of Remark~\ref{nr2.11}
imply that the following properties are 
expressible in terms of the operations of $\rfs{MR}(X,G)$.
\begin{itemize}
\addtolength{\parskip}{-11pt}
\addtolength{\itemsep}{06pt}
\item[(1)] $U$ and $V$ are small,
and $\rfs{cl}(U) \cap \rfs{cl}(V)$ is a singleton.
\item[(2)] $U_1,\,V_1,\,U_2,\,V_2$ are small,
$\rfs{cl}(U_1) \cap \rfs{cl}(V_1)$ is a singleton,
and
$\rfs{cl}(U_1) \cap \rfs{cl}(V_1) = \rfs{cl}(U_2) \cap \rfs{cl}(V_2)$.
\item[(3)] $U$ and $V$ are small,
$\rfs{cl}(U) \cap \rfs{cl}(V)$ is a singleton,
and $\rfs{cl}(U) \cap \rfs{cl}(V) \subseteq W$.
\vspace{-05.7pt}
\end{itemize}
A word of caution. In (1) - (3) smallness cannot be omitted.
This is so, since in Lemmas \ref{l2.8} - \ref{l2.12} the equivalence
of (1) - (3) to the expressible properties mentioned there was proved
only under the assumption that the sets in question are small.

We are ready to define $\fnn{\tau}{X_1}{X_2}$. Let $x \in X_1$.
By Lemma \ref{l2.13}, there are small $U$ and $V$
such that $\sngltn{x} = \rfs{cl}(U) \cap \rfs{cl}(V)$.
Since $\psi$ is an isomorphism between
$\rfs{MR}(X_1,G_1)$ and $\rfs{MR}(X_2,G_2)$,
and by the expressibility of (1) above,
$\rfs{cl}(\psi(U)) \cap \rfs{cl}(\psi(V))$ is a singleton.
Denote it by $\sngltn{y}$ and define $\tau(x) = y$.

By the expressibility of (2) above:
if $U',\,V'$ are small and
$\sngltn{x} = \rfs{cl}(U') \cap \rfs{cl}(V')$,
then $\rfs{cl}(\psi(U') \cap \rfs{cl}(\psi(V')) = \sngltn{y}$.
So the definition of $\tau$ is valid.
As in the proof of Theorem~\ref{t2.2},
Lemma~\ref{l2.13} and the expressibility of (1) and (2)
imply that $\tau$ is $\onetoonen$ and onto.
As in the proof of Theorem~\ref{t2.2},
the expressibility of (3) implies that $\tau$ is a
homeomorphism and that $\tau$ induces $\varphi$.
This completes the proof of Theorem~\ref{t2.4}.\medskip\hfill\myqed

Consider the class
   \index{N@knl@@$K_{\srfs{NL}}$}
$$
K_{\srfs{NL}} = \setm{\pair{X}{G}}{X
\mbox{ is an open subset of a normed space and }
\rfs{LIP}(X) \,\leq\, G \,\leq\, H(X)}.
$$
It is not known whether $K_{\srfs{NL}}$ is faithful.
But we can show the faithfulness of the subclass of $K_{\srfs{NL}}$
consisting of those $\pair{X}{G}$'s in which $X$ is a first category
topological space and $G$ is internally extendible. (See below).
To this end we have strengthened the original statement of
Theorem \ref{t2.4}, and included $G$'s which are required to contain
$\rfs{LIP}(X;F)$ rather than $\rfs{LIP}(X)$.
Since $\rfs{LIP}(X;F) \subseteq \rfs{LIP}(X)$,
this is a stronger result.

\begin{defn}\label{metr-bldr-d2.23}
\begin{rm}
Suppose that $E$ is a normed vector space,
and that $X \subseteq E$ is open.

(a) The {\it complete interior} of $X$ in $E$ is defined by
\newline
\centerline{
$\overline{\rfs{int}}^E(X) =
\bigcup \setm{B^{\oversE}(x,r)}{x \in E \mbox{ and }
B^E(x,r) \subseteq X}$.
}
   \index{N@int01@@
          $\overline{\rfs{int}}^E(X) =
          \bigcup \setm{B^{\oversE}(x,r)}{x \in X \mbox{ and }
          B^E(x,r) \subseteq X}$}
Note that $\overline{\rfs{int}}^E(X)$ is open in $\overE$.

(b) Let $h \in H(X)$. We say that $h$ is
{\it internally extendible in $E$},
if there is $\barh \in H(\overfs{int}^E(X))$
such that $\barh$ extends $h$.
Let $\rfs{IXT}^E(X)$ denote the group of internally extendible
homeomorphisms of $X$.
   \index{internally extendible in $E$.
          A homeomrphism of $X \subseteq E$ which extends to a
          \newline\indent\rule{0pt}{1pt}\kern00mm
          continuous function on $\overline{\rfs{int}}^E(X)$}

   \index{N@ixt@@$\rfs{IXT}^E(X)$. The group of bi-externally-extendible
          auto-homeomorphisms of $X$}

(c) Let $X$ be an open subset of a normed space $E$,
and $\calU$ be a set
of open subsets of $X$.
Then $\calU$ is a {\it complete cover} of $X$
if $\bigcup \setm{\overfs{int}(U)}{U \in \calU} = \overfs{int}(X)$.

   \index{complete cover. $\calU$ is a complete cover of $X$,
          if $\bigcup \setm{\overfs{int}(U)}{U \in \calU} =
	  \overfs{int}(X)$}

(d) For a subset $A$ of a metric space denote the
{\it diameter of $A$} by
$\rfs{diam}(A)$. That is,
$\rfs{diam}(A) = \sup_{x,y \in A} d(x,y)$.
So $\rfs{diam}(A) \in \bbR \cup \sngltn{\infty}$.
   \index{N@diam@@$\rfs{diam}(A) = \sup_{x,y \in A} d(x,y)$}
\hfill\proofend
\end{rm}
\end{defn}

The following proposition is known. See \cite{BP},
the chapter on incomplete norms.
We present a proof here.

\begin{prop}\label{metr-bldr-p2.25}
\num{a} Let $E$ be a normed space and
$x,y \in B^{\oversE}(0,a) - E$. Then there is
$h \in \rfs{LIP}(\overE;E)\sprtl{B^{\oversE}(0,a)}$ such that
$h(x) = y$.

\num{b} Let $E$ be a normed space, $x \in B^E(0,a)$ and
$y \in B^{\oversE}(0,a) - E$. Then there is
$h \in \rfs{LIP}(\overE)\sprtl{B^{\oversE}(0,a)}$ such that
$h(E - \sngltn{x}) = E$ and $h(x) = y$.
\end{prop}

\noindent
{\bf Proof }
(a) We leave the straight-forward proof of the following claim
to the reader.

{\bf Claim 1.} Let $E$ be a normed space.
Let $\setm{K_n}{n \in \bbN} \subseteq (1,\infty)$ be such that
$\prod_{n \in \bbN} K_n < \infty$
and $\setm{g_n}{n \in \bbN} \subseteq \rfs{LIP}(\overE;E)$
be such that:
\begin{itemize}
\addtolength{\parskip}{-11pt}
\addtolength{\itemsep}{06pt}
\item[(1)]
$g_n$ is $K_n$-bilipschitz;
\item[(2)]
$\sum_{n \in \bbN} \rfs{diam}(\rfs{supp}(g_n)) < \infty$;
\item[(3)]
there is $x_0 \in \overE - E$ and a sequence
$\setm{r_n}{n \in \bbN} \subseteq (0,\infty)$ converging to~$0$
such that for every $n \in \bbN$,
$\rfs{supp}(g_n) \subseteq
g_{n - 1} \scirc \ldots \scirc g_0(B^{\overE}(x_0,r_n))$.
\vspace{-05.7pt}
\end{itemize}
Let $h_n = g_{n - 1} \scirc \ldots \scirc g_0$.
\underline{Then} for every $x \in \overE$,
$\lim_{n \rightarrow \infty} h_n(x)$ exists.
Define
$h(x) = \lim_{n \rightarrow \infty} h_n(x)$.
Then $h \in \rfs{LIP}(\overE;E)$.
\smallskip

We construct $g_n$'s which satisfy the assumptions of Claim 1.
Let
$\setm{M_n}{n \in \bbN} \subseteq (3,\infty)$ be such that
$\prod_{n \in \bbN} (1 + \dgfrac{1}{(M_n - 3})) < \infty$.
We may assume that $\norm{x - y} \mcdot M_0 < a$.
Set $x = x_0$ and $\norm{x - y} = d_0$.
Define $d_n$ by induction as follows:
$d_{n + 1} = \dgfrac{d_n}{M_{n + 1}}$.

We shall apply Proposition \ref{l2.6}(b)(ii). The normed space $E$ of
\ref{l2.6} is taken to be $\overE$, $\calS = \sngltn{\overE}$,
$F_{\oversE} = \overE$ and $a$ of \ref{l2.6}(b) is $a$ here.
The homeomorphism $h$ constructed in Proposition \ref{l2.6}(b)
depended on the vectors $x_0$ and $v$ and on the radii $r_0$ and $r$.
Denote that $h$ by $h_{x_0,v,r_0,r}$.

We define $g_n$ and $x_{n + 1}$ by induction.
Suppose that $x_n$ has been defined.
Let
$$
u_n = d_{n + 1} \mcdot \frac{y - x_n}{\norm{y - x_n}}
\mbox{ \,and\, }
f_n = h_{x_n,u_n,M_n d_n,2 d_n}.
$$
So
$\rfs{supp}(f_n) \subseteq B(x_n,M_n d_n)$.
Note that $f_n$ is
$(1 + \frac{d_n}{M_n d_n - 2 d_n - d_{n + 1}})$-bilipschitz.
Since $d_{n + 1} < d_n$, we have
$\frac{d_n}{M_n d_n - 2 d_n - d_{n + 1}} > \frac{1}{M_n - 3}$.
So
\begin{list}{}
{\setlength{\leftmargin}{39pt}
\setlength{\labelsep}{08pt}
\setlength{\labelwidth}{20pt}
\setlength{\itemindent}{-00pt}
\addtolength{\topsep}{-04pt}
\addtolength{\parskip}{-02pt}
\addtolength{\itemsep}{-05pt}
}
\item[(1.1)] 
$\norm{y - f_n(x_n)} = d_{n + 1} < 2 d_{n + 1}$,
\item[(1.2)] 
$f_n \nrestriction B(x_n,2 d_n) =
\rfs{tr}_{u_n} \nrestriction B(x_n,2 d_n)$,
\item[(1.3)] 
for some $\varepsilon > 0$,
$f_n$ is $(1 + \frac{1}{M_n - 3} + \varepsilon)$-bilipschitz,
\item[(1.4)] 
if $n > 0$, then for some $\varepsilon > 0$,
$\rfs{supp}(f_n) \subseteq B(x_n,d_{n - 1} - \varepsilon)$.
\vspace{-02.0pt}
\end{list}
Choose $y_n,v_n \in E$ close enough to $x_n$ and $u_n$ respectively
so that for $g_n$ defined by $g_n = h_{y_n,v_n,M_n d_n,2 d_n}$
the following holds:
\begin{list}{}
{\setlength{\leftmargin}{39pt}
\setlength{\labelsep}{08pt}
\setlength{\labelwidth}{20pt}
\setlength{\itemindent}{-00pt}
\addtolength{\topsep}{-04pt}
\addtolength{\parskip}{-02pt}
\addtolength{\itemsep}{-05pt}
}
\item[(2.1)]
(2.1) $\norm{y - g_n(x_n)} < 2 d_{n + 1}$,
\item[(2.2)]
$g_n \nrestriction B(x_n,d_n) =
\rfs{tr}_{v_n} \nrestriction B(x_n,d_n)$,
\item[(2.3)]
$g_n$ is $(1 + \frac{1}{M_n - 3})$-bilipschitz,
\item[(2.4)]
if $n > 0$, then $\rfs{supp}(g_n) \subseteq B(x_n,d_{n - 1})$.
\vspace{-02.0pt}
\end{list}
Let $x_{n + 1} = g_n(x_n)$.
So $x_{n + 1} = x_n + v_n$.
Also,
$g_n \in \rfs{LIP}(\overE;E)$

We check that (1)-(3) of Claim 1 are fulfilled.
Clearly,
$K_n = 1 + \frac{1}{M_n - 3}, \ n \in \bbN$ fulfill
Clause (1).
Since $d_{n + 1} < \dgfrac{d_n}{3}$,
we have $\sum_{n \in \bbN} d_n < \infty$. So
$\sum_{n \in \bbN} \rfs{diam}(\rfs{supp}(g_n)) <
\sum_{n \in \bbN} 2 d_n < \infty$, proving (2).

Let $h_n = g_n \scirc \ldots \scirc g_0$
and $w_n = \sum_{i \leq n} v_i$. We show by induction that
\begin{list}{}
{\setlength{\leftmargin}{39pt}
\setlength{\labelsep}{08pt}
\setlength{\labelwidth}{20pt}
\setlength{\itemindent}{-00pt}
\addtolength{\topsep}{-04pt}
\addtolength{\parskip}{-02pt}
\addtolength{\itemsep}{-05pt}
}
\item[(2.5)]
$h_n \nrestriction B(x_0,d_n) = \rfs{tr}_{w_n} \nrestriction B(x_0,d_n)$
for every $n \in \bbN$.
\vspace{-02.0pt}
\end{list}
By (2.2), this is true for $n = 0$.
Assume it is true for $n$.
Hence $x_{n + 1} = h_n(x_0) = x_0 + w_n$.
For $n + 1$ we have
\vspace{1.5mm}
\newline
\rule{0pt}{0pt}
\renewcommand{\arraystretch}{1.5}
\addtolength{\arraycolsep}{-6pt}
$
\begin{array}{ll}
&
h_{n + 1} \nrestriction B(x_0,d_{n + 1}) =
(g_{n + 1} \scirc h_n) \nrestriction B(x_0,d_{n + 1}) =
g_{n + 1} \nrestriction h_n(B(x_0,d_{n + 1}))  \scirc
\rfs{tr}_{w_n} \nrestriction B(x_0,d_{n + 1})
\\
=
\rule{5pt}{0pt}
&
g_{n + 1} \nrestriction B(x_0 + w_n,d_{n + 1}) \scirc
\rfs{tr}_{w_n} \nrestriction B(x_0,d_{n + 1}) =
g_{n + 1} \nrestriction B(x_{n + 1},d_{n + 1}) \scirc
\rfs{tr}_{w_n} \nrestriction B(x_0,d_{n + 1})
\\
=
\rule{5pt}{0pt}
&
\rfs{tr}_{v_{n + 1}} \nrestriction B(x_{n + 1},d_{n + 1}) \scirc
\rfs{tr}_{w_n} \nrestriction B(x_0,d_{n + 1}) =
\rfs{tr}_{w_{n + 1}} \nrestriction B(x_0,d_{n + 1}).
\end{array}
$
\renewcommand{\arraystretch}{1.0}
\addtolength{\arraycolsep}{6pt}
\newline
It follows from (2.4) and (2.5) that
$\rfs{supp}(g_{n + 1}) \subseteq B(x_{n + 1},d_n) =
h_n(B(x_0,d_n))$.
Since $\lim_{n \rightarrow \infty} d_n = 0$,
Clause (3) of Claim 1 holds.
Let $h$ be as constructed in Claim 1.
So $h \in \rfs{LIP}(\overE;E)$.

Since $\norm{y - x_n} = d_n$ and 
$\lim_{n \rightarrow \infty} d_n = 0$, we have $h(x) = y$.
We show that $\rfs{supp}(g_n) \subseteq B(x,a)$
for every $n \in \bbN$.
For $n = 0$,
$\rfs{supp}(g_0) \subseteq B(x,M_0 d_0) \subseteq B(x,a)$.
Suppose that $n > 0$. Then
$\rfs{supp}(g_n) \subseteq B(x_n,M_n d_n) \subseteq
B(x,M_n d_n + \norm{x_n - x})$.
Since
$$M_n d_n + \norm{x_n - x} \leq
M_n d_n + \norm{x_n - y} + \norm{y - x} <
d_{n - 1} + 2 d_n + d_0 < 3 d_0 < M_0 d_0 < a,
$$
we have $\rfs{supp}(g_n) \subseteq B(x,a)$.
It follows that $\rfs{supp}(h) \subseteq B(x,a)$.
So $h$ is as required.

(b) The proof very similar to the proof of (a).
\rule{0pt}{1pt}\hfill\myqed

\begin{cor}\label{nc2.23}\label{metr-bldr-c2.26}
Let $K_{\srfs{NFCB}}$ be the class of all space-group pairs
$\pair{X}{G}$ for which there is a normed space $E$ such that
$X$ is an open subset of $E$ and
\begin{enumerate}
\item[\num{1}] $E$ is of the first category, or $E$ is a Banach space;
\item[\num{2}] There is a complete cover $\calU$ of $X$ such that
$\rfs{LIP}(X,\calU) \leq G \leq \rfs{IXT}(X)$.
\end{enumerate}
Then $K_{\srfs{NFCB}}$ is faithful.

   \index{N@knfcb@@$K_{\srfs{NFCB}}$}

\end{cor}

\noindent
{\bf Proof }
Let $\pair {X}{G} \in K_{\srfs{NFCB}}$. For $g \in G$ let 
$\barg$ be the extension of $g$ to $\overfs{int}(X)$ and
$\rule{1pt}{0pt}\overG = \setm{\barg}{g \in G}$.
Then $\pair{\overfs{int}^E(X)}{\overG} \in K_{\srfs{BO}}$.

Let $\overcalO(X,G)$ be the set of orbits of $\overG$.
That is,
$\overcalO(X,G) = \setm{\overG(x)}{x \in \overfs{int}(X)}$.
It follows from Proposition \ref{metr-bldr-p2.25}(a) that
if $X$ is an open subset of an incomplete normed space,
then for every $O \in \overcalO(X,G)$
there is a set $\calC$ of connected components of
$\overfs{int}(X)$ such that
$O = E \cap \bigcup \calC$ or
$O = (\overE - E) \cap \bigcup \calC$.
Clearly, if $X$ is an open subset of a Banach space,
then for every $O \in \overcalO(X,G)$ there is a set of connected
components of $X$ such that $O = \bigcup \calC$.
Let
$\rfs{FC}(X,G) = \bigcup \setm{O \in \overcalO(X,G)}
{O \mbox{ is a first category set}}$.
If $X$ is of the first category, then $X = \rfs{FC}(X,G)$.

For $i = 1,2$ let $\pair{E_i}{G_i} \in K_{\srfs{NFCB}}$,
and let $\iso{\varphi}{G_1}{G_2}$.
Let $\fnn{\barvarphi}{\overG_1}{\overG_2}$ be defined by
$\barvarphi(\barg) = \overline{\varphi(g)}$.
Then $\iso{\barvarphi}{\overG_1}{\overG_2}$.
By Theorem \ref{t2.4}(b), there is $\iso{\bartau}{\overE_1}{\overE_2}$
which induces $\barvarphi$.
Obviously,
$\bartau$ takes orbits of $\overG_1$ to orbits of $\overG_2$.
So $\overcalO(X,G_1)$ contains members of the first category
iff
$\overcalO(X,G_2)$ contains members of the first category.

It is obvious that $\bartau$ takes every first category orbit of
$\overG_1$ to a first category orbit of $\overG_2$.
So if $X_1$ is of the first category,
then
$\bartau(X_1) = \bartau(\rfs{FC}(X_1,G_1)) = \rfs{FC}(X_2,G_2) = X_2$,
and hence $\iso{\tau}{X_1}{X_2}$.
If $X_1$ is an open subset of a Banach space,
then $\bartau = \tau$ and hence $\iso{\tau}{X_1}{X_2}$.
\vspace{-2mm}
\hfill\myqed

\begin{remark}\label{nnr2.24}
\begin{rm}
If $E$ has a countable Hamel basis, then it is of the first category.
The space $\ell_1$ is a linear subspace of $\ell_2$,
and it is of the first category in $\ell_2$.

This is a special case of the following fact.
If $\fnn{T}{F}{E}$ is a bounded linear operator
from a Banach space $F$ to a Banach space $E$, and
$\rfs{Rng}(T)$ is a proper dense subset of $E$,
then $\rfs{Rng}(T)$ is of the first category in $E$.
This follows from the proof of the Open Mapping Theorem.
If $\rfs{Rng}(T)$ is of the second category, then for some ball
$B = B^F(0,n)$, $T(B)$ is somewhere dense.
Hence $T(B)$ is dense in some ball of the form $B^E(0,r)$.
It can then be proved that $T(B) \supseteq B^E(0,r)$.
This implies that $\rfs{Rng}(T) = E$.
\hfill\proofend
\end{rm}
\end{remark}

In Corollary \ref{nc2.23} the assumptions that $E$ is of the
first category, and that $G$ is completely extendible are undesirable.
We do not know whether they can be dispensed with.

The final reconstruction results of Chapter \ref{s5}
are proved for open subsets
of first category normed vector spaces and for open subsets of
Banach spaces.
The proofs of all intermediate theorems are valid for open subsets
of any normed space.
If Parts (c) or (d) of the following question have a negative answer,
then the final results of Chapter \ref{s5} will be true for open
subsets of any normed vector space.

On the other hand, examples answering (c) or (d) below 
in the affirmative,
imply that certain results in Chapter \ref{s5} are not true for
arbitrary normed spaces.

\begin{question}\label{nq2.25}
\begin{rm}
(a) Is $K_{\srfs{NL}}$ faithful?

(b) Let $K_{\srfs{NLIX}}$ be the subclass of
$K_{\srfs{NL}}$ consisting of all
$\pair{E}{G}$'s in which $G$ is internally extendible.
Is $K_{\srfs{NLIX}}$ faithful?

(c)
Are there normed spaces $E$ and $F$ and a homeomorphism
$\iso{\tau}{\overE}{\overF}$ such that $\tau(E) = \overF - F$?

\smallskip

Note that the answer to (b) is positive iff
the answer to (c) is negative.

\smallskip

(d)
Are there normed spaces $E$ and $F$ and a uniformly
bicontinuous homeomorphism
$\iso{\tau}{\overE}{\overF}$ such that $\tau(E) = \overF - F$?
\vspace{-2mm}
\end{rm}
\end{question}

\subsection{Faithfulness of normed manifolds.}
\label{ss2.4}
\label{ss2.4-faithfulness-of-normed-manifolds}

As has been mentioned, the proof of \kern1pt Theorem \ref{t2.4}
extends without change to manifolds over normed vector spaces.
This class contains some new instances.
The unit sphere of a normed space is one,
and spaces which are a finite product of manifolds is another.

We extend the results a bit further, in order to allow
the inclusion of manifolds with boundary over a normed vector space.
To this end we introduce the notion of a
``regionally normed manifold''.
By combining Remark \ref{r2.15} with the various results on extendible
homeomorphism groups appearing in Chapter \ref{s5},
one obtains reconstruction theorems for manifolds with boundary.

It should be pointed out that no new arguments are needed for this
new framework,

\begin{defn}\label{d2.14}
\begin{rm}
(a) Let $X$ be a topological space. A family of mappings
$\itPhi$ is called a {\it regional normed atlas} for $X$
if the following holds.
\smallskip
   \index{regional normed atlas for $X$}

(1) $\bigcup \setm{\rfs{Rng}(\varphi)}{\varphi \in \itPhi}$
is a dense subset of $X$.

{\thickmuskip=2.7mu \medmuskip=1mu \thinmuskip=1mu 
(2) For every $\varphi \in \itPhi$
there is a normed space
$E = E_{\varphi}$, $x = x_{\varphi} \in E$
and $r = r_{\varphi} > 0$
such that:
}
\begin{itemize}
\addtolength{\parskip}{-11pt}
\addtolength{\itemsep}{06pt}
\item[(i)] $\fnn{\varphi}{\overB^E(x,r)}{X}$,
\item[(ii)] $\varphi$ is a homeomorphism between $\rfs{Dom}(\varphi)$
and $\rfs{Rng}(\varphi)$,
\item[(iii)] $\rfs{Rng}(\varphi)$ is closed in $X$, and 
$\varphi(B^E(x,r))$ is open in $X$.
\vspace{-05.7pt}
\end{itemize}
If $\itPhi$ is a regional normed atlas for $X$, then
$\pair{X}{\itPhi}$ is called a {\it regionally normed manifold (RNM)}.
   \index{regionally normed manifold (RNM)}
   \index{rnm   @@RNM. A regionally normed manifold}
If
$X = \bigcup \setm{\varphi(B^{E_{\varphi}}(x_{\varphi},r_{\varphi}))}
{\varphi \in \itPhi}$,
then $\pair{X}{\itPhi}$ is called a {\it normed manifold}.
   \index{normed manifold}
Let $\pair{X}{\itPhi}$ be an RNM.
If for every $\varphi \in \itPhi$, $E_{\varphi}$ is a Banach space,
then $\pair{X}{\itPhi}$ is said to be a
{\it regional Banach manifold (RBM)}.
   \index{regional Banach manifold (RBM)}
   \index{rbm   @@RBM. A regional Banach manifold}
   \index{banach manifold@@Banach manifold}
A normed manifold which is an RBM is called a {\it Banach manifold}.

(b) Recall that for a metric space $(Y,d)$, $x \in Y$ and $r > 0$,
$S^Y(x,r)$ denotes
\newline
$\setm{y \in Y}{d(x,y) = r}$. 
For a normed space $E$, $x \in E$ and $r > 0$ let
\addtolength{\arraycolsep}{-3pt}
\newline
\indent
$
\begin{array}{lll}
\rule{0pt}{16pt}
L_1(E,x,r)& \eqdf &\setm{h \in H(\overB^E(x,r))}
{h \mbox{ is bilipschitz, and }
h \nrestriction S(x,r) = \rfs{Id}},
\\
\rule{0pt}{16pt}
L_1^{\srfs{LC}}(E,x,r)& \eqdf& \setm{h \in H(\overB^E(x,r))}
{h \mbox{ is locally bilipschitz, and }
h \nrestriction S(x,r) = \rfs{Id}}.
\medskip
\end{array}
$
\newline
Let $F$ be a dense linear subspace of $E$. Define
\newline
\indent
$
\begin{array}{lll}
\rule{0pt}{16pt}
L_1(E,x,r;F)& \eqdf& \setm{h \in L_1(E,x,r)}
{h(\overB^E(x,r) \cap F) = \overB^E(x,r) \cap F},
\\
\rule{0pt}{16pt}
L_1^{\srfs{LC}}(E,x,r;F)& \eqdf& \setm{h \in L_1^{\srfs{LC}}(E,x,r)}
{h(\overB^E(x,r) \cap F) = \overB^E(x,r) \cap F}.
\medskip
\end{array}
$
\newline
\addtolength{\arraycolsep}{3pt}
If $\pair{X}{\itPhi}$ is an RNM, $\varphi \in \itPhi$ and
$h \in L_1^{\srfs{LC}}(E_{\varphi},x_{\varphi},r_{\varphi})$, then
\newline
$h^{[\varphi]} \eqdf 
h^{\varphi} \cup \rfs{Id} \nrestriction
(X - \rfs{Rng}(\varphi)) \in H(X)$.
Suppose that $\calF \eqdf \setm{F_{\varphi}}{\varphi \in \itPhi}$
is a family of linear spaces such that for every $\varphi \in \itPhi$,
$F_{\varphi}$ is a dense subspace of $E_{\varphi}$.
Then $\calF$ is called a {\it subspace choice for $\pair{X}{\itPhi}$}.
   \index{subspace choice for $\pair{X}{\itPhi}$}

Let
$\rfs{LIP}(X;\itPhi,\calF)$ denote the subgroup of $H(X)$ generated by
$\setm{h^{[\varphi]}}{\varphi \in \itPhi,\ %
h \in L_1(E_{\varphi},x_{\varphi},r_{\varphi};F_{\varphi})}$.
   \index{N@lip06@@$\rfs{LIP}(X;\itPhi,\calF)$}
Let $\rfs{LIP}^{\srfs{LC}}(X;\itPhi,\calF)$
denote the subgroup of $H(X)$ generated by
$\setm{h^{[\varphi]}}{\varphi \in \itPhi,\ %
h \in
L_1^{\srfs{LC}}(E_{\varphi},x_{\varphi},r_{\varphi};F_{\varphi})}$.
   \index{N@liplc06@@$\rfs{LIP}^{\srfs{LC}}(X;\itPhi,\calF)$}
If $F_{\varphi} = E_{\varphi}$ for every $\varphi \in \itPhi$,
then $\rfs{LIP}(X;\itPhi,\calF)$ and
$\rfs{LIP}^{\srfs{LC}}(X;\itPhi,\calF)$
are denoted by
$\rfs{LIP}(X;\itPhi)$ and $\rfs{LIP}^{\srfs{LC}}(X;\itPhi)$
respectively.

   \index{N@lip07@@$\rfs{LIP}(X;\itPhi)$}
   \index{N@liplc07@@$\rfs{LIP}^{\srfs{LC}}(X;\itPhi)$}

Remark: Even though the groups considered below contain
$\rfs{LIP}(X;\itPhi,\calF)$,
we do not have to require at this point that the transition maps 
in the atlas be Lipschitz.
That is, we do not require that
$\varphi\inverse \scirc \psi$ is bilipschitz
for every $\varphi,\psi \in \itPhi$.

(c) Let $K_{\srfs{BM}}$ be the class of all $\pair{X}{G}$'s
which satisfy the following:
There are $\itPhi$ and $\calF$ such that
\begin{itemize}
\addtolength{\parskip}{-11pt}
\addtolength{\itemsep}{06pt}
\item[(1)] $\pair{X}{\itPhi}$ is a Banach manifold and $\calF$
is a subspace choice for $\itPhi$,
\item[(2)] $\rfs{LIP}(X;\itPhi,\calF) \leq G \leq H(X)$.\smallskip
\vspace{-05.7pt}
\end{itemize}
   \index{N@kbm@@$K_{\srfs{BM}}$}
Let $K_{\srfs{NM}}$ be the class of all $\pair{X}{G}$'s
which satisfy the following:
There are $\itPhi$ and $\calF$ such that
\begin{itemize}
\addtolength{\parskip}{-11pt}
\addtolength{\itemsep}{06pt}
\item[(1)] $\pair{X}{\itPhi}$ is a normed manifold and $\calF$
is a subspace choice for $\itPhi$,
\item[(2)] $\rfs{LIP}^{\srfs{LC}}(X;\itPhi,\calF) \leq G \leq H(X)$.
\vspace{-05.7pt}
\end{itemize}
   \index{N@knm@@$K_{\srfs{NM}}$}
Let $K_{\srfs{BNM}} = K_{\srfs{BM}} \cup K_{\srfs{NM}}$.
	             \index{N@kbnm@@$K_{\srfs{BNM}}$}

(d) Let $\pair{X}{\itPhi}$ be an RNM.
{\thickmuskip=2mu \medmuskip=1mu \thinmuskip=1mu 
The set \hbox{$\rfs{NI}(X,\itPhi) \eqdf
\bigcup \setm{\varphi(B^{E_{\varphi}}(x_{\varphi},r_{\varphi}))}
{\varphi \in \itPhi}$}
}
is called the {\it normed interior} of $\pair{X}{\itPhi}$.
     \index{N@ni@@$\rfs{NI}(X,\itPhi) =
     \bigcup \setm{\varphi(B^{E_{\varphi}}(x_{\varphi},r_{\varphi}))}
     {\varphi \in \itPhi}$. The normed interior of
     $\pair{X}{\itPhi}$}

Let $G \leq H(X)$.
The
{\it extended normed interior} of $\trpl{X}{\itPhi}{G}$ is defined as
$$\rfs{ENI}(X,\itPhi,G) \eqdf
\setm{g(x)}{x \in \rfs{NI}(X,\itPhi) \mbox{ and } g \in G}.$$
Also, $\rfs{ENI}(X,\itPhi,H(X))$ is denoted by $\rfs{ENI}(X,\itPhi)$.

If $X$ is a subset of a normed space $E$ and
$\rfs{int}^E(X)$ is dense in $X$,
then $X$ is a regional normed manifold.
As a regional normed atlas for $X$ we take the set
$\itPhi =
\setm{\rfs{Id} \nrestriction
\rule{-1pt}{0pt}\overB^E(x,r)}{\overB^E(x,r) \subseteq X}$.
We denote $\rfs{ENI}(X,\itPhi)$ by $\rfs{ENI}(X)$.
Hence
$\rfs{ENI}(X) =
\setm{h(x)}{x \in \rfs{int}^E(X),\ h \in H(X)}$.

   \index{N@eni00@@$\rfs{ENI}(X,\itPhi,G) =
   \setm{g(x)}{x \in \rfs{NI}(X,\itPhi) \mbox{ and } g \in G}$.
   Extended normed interior of $\trpl{X}{\itPhi}{G}$}
   \index{N@eni01@@$\rfs{ENI}(X,\itPhi) = \rfs{ENI}(X,\itPhi,H(X))$}
   \index{N@eni02@@$\rfs{ENI}(X) =
   \setm{h(x)}{x \in \rfs{int}^E(X),\ h \in H(X)}$}

\end{rm}
\end{defn}

\begin{theorem}\label{t2.15}
\num{a} $K_{\srfs{BNM}}$ is faithful.

\num{b}
For $i = 0,1$ let $\pair{X_i}{\itPhi_i}$ be an RNM and
$\calF_i$ be a subspace choice for $\pair{X_i}{\itPhi_i}$.
Let $G_i \leq H(X_i)$.
Suppose that for $i = 0,1$:
\begin{itemize}
\addtolength{\parskip}{-11pt}
\addtolength{\itemsep}{06pt}
\item[\num{1}] if $\pair{X_i}{\itPhi_i}$ is an RBM,
then $\rfs{LIP}(X_i,\itPhi_i;\calF_i) \leq G_i$,
\item[\num{2}] if $\pair{X_i}{\itPhi_i}$ is not an RBM,
then $\rfs{LIP}^{\srfs{LC}}(X_i,\itPhi_i;\calF_i) \leq G_i$.
\vspace{-05.7pt}
\end{itemize}
Let $\iso{\varphi}{G_1}{G_2}$. Then there is
$\iso{\tau}{\rfs{ENI}(X_1,\itPhi_1,G_1)}{\rfs{ENI}(X_2,\itPhi_2,G_2)}$
such that $\tau$ induces~$\varphi$.
That is,
$\varphi(g) \nrestriction \rfs{ENI}(X_2,\itPhi_2,G_2) =
(g \nrestriction \rfs{ENI}(X_1,\itPhi_1,G_1))^{\tau}$
for every $g \in G_1$.

\num{c} Let $X$ be a subset of a normed space $E$ and $Y$ be a subset of
a normed space $F$ such that $\rfs{int}^E(X)$ is dense in $X$ and
$\rfs{int}^F(Y)$ is dense in $Y$.
Suppose that $\iso{\varphi}{H(X)}{H(Y)}$.
Then there is $\iso{\tau}{\rfs{ENI}(X)}{\rfs{ENI}(Y)}$
such that $\tau$ induces~$\varphi$.
That is, for every $g \in H(X)$,
$\varphi(g) \nrestriction \rfs{ENI}(Y) =
(g \nrestriction \rfs{ENI}(X))^{\tau}$.
\vspace{-2mm}
\end{theorem}

\noindent
{\bf Proof }
(a) If $\pair{X}{G} \in K_{\srfs{BNM}}$
and $\itPhi$ is a normed regional atlas for $X$ which demonstrates that
$X$ is a normed manifold,
then
$\rfs{NI}(X,\itPhi) = X$.
So $\rfs{ENI}(X,\itPhi,G) = X$.
Hence (b) implies (a). 

(b) The proof of Theorem \ref{t2.4} applies without change.

(c) This is a special case of (b).
\hfill\myqed

\begin{remark}\label{r2.15}
\begin{rm}
The proof  of the above theorem applies to RNM's too.
The statement that is proved for  RNM's is as follows.
If $\iso{\varphi}{G_1}{G_2}$, then there is
$\tau\,:\,\rfs{ENI}(X_1,\itPhi_1,G_1) \cong\break
\rfs{ENI}(X_2,\itPhi_2,G_2)$
such that $\tau$ induces $\varphi$.
That is, for every $g \in G_1$, \ 
$\varphi(g) \nrestriction \rfs{ENI}(X_2,\itPhi_2,G_2) =
(g \nrestriction \rfs{ENI}(X_1,\itPhi_1,G_1))^{\tau}$.
\hfill\proofend
\end{rm}
\end{remark}


Manifolds with boundary, closures of open subsets of a normed space
and closures of open subsets of a normed manifold
are obviously RNM's. Note that in the above theorem, the groups $G_i$
are not assumed to preserve the boundary of $X_i$.
Indeed, when the $X_i$'s are infinite-dimensional, it may happen
that their boundary is not preserved.

\subsection{The faithfulness of some smaller subgroups.}
\label{ss2.5}
\label{ss2.5-faithfulness-of-smaller-subgroups}

The homeomorphisms constructed in Lemma \ref{l2.6}(b) suggest
some new types of subgroups of $H(X)$ which may be
interesting in the context of reconstruction and in other
contexts involving homeomorphisms of infinite-dimensional spaces.

\begin{defn}\label{metr-bldr-d2.31}
\begin{rm}
Let $X$ be an open subset of a normed vector space $E$ and $g \in H(X)$.
\smallskip

(a) We call $g$ a
{\it ``finite-dimensional difference'' homeomorphism},
if there is a finite-dimensional subspace $F$ of $E$ such that
$g(x) - x \in F$ for every $x \in X$.

   \index{finite-dimensional difference homeomorphism}

Let $\rfs{FD}(X)$ denote the set of
``finite-dimensional difference'' homeomorphisms of $X$
and $\rfs{FD.LIP}(X) \eqdf \rfs{FD}(X) \cap \rfs{LIP}(X)$.
   \index{N@fd@@$\rfs{FD}(X)$}
   \index{N@fdlip@@$\rfs{FD.LIP}(X)$}

(b) We call $g$ a
{\it weakly ``finite-dimensional difference'' homeomorphism},
if there is a finite-dimensional subspace $F$ of $E$ such that
for every $x \in X$ there is $a \in \bbR - \sngltn{0}$
such that $g(x) - ax \in F$.
   \index{weakly ``finite-dimensional difference'' homeomorphism}

Let $\rfs{WFD}(X)$ denote the set of
weakly ``finite-dimensional difference'' homeomorphisms of $X$
and $\rfs{WFD.LIP}(X) \eqdf \rfs{WFD}(X) \cap \rfs{LIP}(X)$.
   \index{N@wfd@@$\rfs{WFD}(X)$}
   \index{N@wfdlip01@@$\rfs{WFD.LIP}(X)$}
For a subspace choice system $\frtpl{E}{X}{\calS}{\calF}$
define $\rfs{WFD.LIP}(X;\calS,\calF)$
and
$\rfs{WFD.LIP}^{\srfs{LC}}(X;\calS,\calF)$
in analogy to the definition
of\break
$\rfs{LIP}(X;\calS,\calF)$. See Definition \ref{d2.7}(a).
Also, define $K_{\srfs{WFD.BNO}}$ in analogy to the definition of
$K_{\srfs{BNO}}$.
   \index{N@wfdlip01@@$\rfs{WFD.LIP}(X;\calS,\calF)$}
   \index{N@wfdlclip@@$\rfs{WFD.LIP}^{\srfs{LC}}(X;\calS,\calF)$}
   \index{N@kwfdbno@@$K_{\srfs{WFD.BNO}}$}
\rule{1pt}{0pt}\hfill\proofend
\end{rm}
\end{defn}

It is easy to check that $\rfs{FD}(X)$ and $\rfs{WFD}(X)$ are groups.
The following is a corollary of the proof of Theorem \ref{t2.4}.

\begin{cor}\label{metr-bldr-c2.32}
$K_{\srfs{WFD.BNO}}$ is faithful.
\vspace{-1.0mm}
\end{cor}

\noindent
{\bf Proof } The proof of Theorem \ref{t2.4} applies,
since it uses only homeomorphisms belonging to $\rfs{WFD}(X)$.
\hfill\myqed

By Lemma \ref{l2.6}(b), $\rfs{FD.LIP}(X)$ is locally moving.
In fact, the construction of \ref{l2.6}(b)
can be used to show that $\rfs{FD.LIP}(X)$ is transitive in
the following sense.
There is an open base $\calB$ of $X$
such that for every $B \in \calB$ and for every
finite injective function $\rho$ whose domain and range are subsets of
$B$ there is $g \in \sprtd{G}{B}$ such that $g$ extends $\rho$.
In fact,
$\calB$ can be taken to be
\hbox{$\setm{B^E(x,r)}{B^E(x,r) \subseteq X}$.}

\begin{question}\label{nq2.28}
\begin{rm}
Are any of the classes related to $\rfs{FD}(X)$ faithful?
For example, is the class
$K_{\srfs{BFD}} \eqdf
\setm{\pair{E}{G}}{E \mbox{ is a Banach space,}
\mbox{ and } \rfs{FD}(E) \leq G \leq H(E)}$
faithful?
\end{rm}
\end{question}
\newpage


\section{The local $\Gamma$-continuity of a
conjugating homeomorphism}
\label{s3}

\subsection{General description.}
\label{ss3.1}
\label{ss3.1-general-description}

The Main Result of this section is the statement that
if $X_1,X_2$ are open subsets of normed spaces $E_1$ and $E_2$
respectively,
$\itGamma_1$ and $\itGamma_2$ are
countably generated moduli of continuity, and $\iso{\tau}{X_1}{X_2}$
is such that
$(H_{\itGamma_{1}}^{\srfs{LC}}(X_1))^{\tau} =
H_{\itGamma_{2}}^{\srfs{LC}}(X_2)$,
then $\itGamma_1 = \itGamma_2$ and $\tau$ is
locally $\itGamma_1$-bicontinuous.
This is proved in Theorem \ref{metr-bldr-t3.19}(a).
Equally central are the four results stated in
Corollary \ref{metr-bldr-c3.43}.

The conjunction of the final results of
Chapters \ref{s2} and \ref{s3} is stated in
Theorem \ref{metr-bldr-t3.42}. It says that
the existence of an isomorphism $\varphi$ between the groups 
$H_{\itGamma_{1}}^{\srfs{LC}}(X_1)$ and
$H_{\itGamma_{2}}^{\srfs{LC}}(X_2)$ implies that
$\itGamma_1 = \itGamma_2$, and that $\varphi$ is induced by a locally
$\itGamma_1$-bicontinuous homeomorphism $\tau$ between $X_1$ and $X_2$.

As in Chapter \ref{s2}, the results quoted above are in fact
special cases of a more general setting. The groups which are
actually being considered are of the type
$H_{\itGamma}^{\srfs{LC}}(X;\calS,\calF)$.
See Definition \ref{metr-bldr-d3.17}.

There are two methods of proving the Main Result.
The central intermediate lemma in Method I roughly says
that if $X_1$, $X_2$ are normed vector spaces,
$\iso{\tau}{X_1}{X_2}$,
and for every translation $\rfs{tr}_v$ of $X_1$,
$(\rfs{tr}_v)^{\tau} \in H_{\itGamma_2}^{\srfs{LC}}(X_2)$,
then $\tau\inverse$
is locally $\itGamma_2$-continuous.
This is in fact the hidden contents of Theorem \ref{metr-bldr-t3.15}.
A variant of this statement which works only for second category
spaces, but yields a slightly stronger result
is proved in Theorem \ref{metr-bldr-t3.26}.

The main lemma in Method II says roughly
that if $X_1$, $X_2$ are normed vector spaces,
$\iso{\tau}{X_1}{X_2}$,
and for every bounded affine isomorphism $T$ of $X_1$,
$T^{\tau} \in H_{\itGamma_2}^{\srfs{LC}}(X_2)$,
then $\tau$ is locally $\itGamma_2$-bicontinuous.

Going back to the Main Result,
we in fact prove a stronger statement.
Suppose that $\fourtpl{E}{X}{\calS}{\calE}$
and $\fourtpl{F}{Y}{\calT}{\calF}$ are subspace choice systems,
\hbox{$\itGamma,\itDelta$ are}
countably generated moduli of continuity,
$\iso{\tau}{X}{Y}$, \hbox{and the following holds:}
$$
(H_{\itGamma}(X;\calS,\calF))^{\tau} \subseteq
H_{\itDelta}^{\srfs{LC}}(Y)
\mbox{ \ and \ }
(H_{\itDelta}(Y;\calT,\calF))^{\tau\inverse} \subseteq
H_{\itGamma}^{\srfs{LC}}(X).
$$
Then $\itGamma = \itDelta$
and $\tau$ is locally $\itGamma$-bicontinuous.
This is proved in Theorem \ref{metr-bldr-t3.19}(b).
See Definitions \ref{d2.7} and \ref{metr-bldr-d3.17}(a).

Part of this strengthening is needed in the proof that if
$\iso{\tau}{\rfs{cl}(X)}{\rfs{cl}(Y)}$ and
\newline
$(H_{\itGamma}^{\srfs{LC}}(\rfs{cl}(X)))^{\tau} =
H_{\itDelta}^{\srfs{LC}}(Y)$,
then $\itGamma = \itDelta$ and $\tau$ is
locally $\itGamma$-bicontinuous.

There are two situations in which we use Method I and we cannot use
Method II.
The first one appears in Chapter \ref{s11}, where 
the reconstruction of the closure of an open set is considered.
Method I is used again in the proof that
the derivative of a conjugating homeomorphism is $\itGamma$-continuous.
Such results will appear in a subsequent work.

\subsection{Partial actions and decayability.}
\label{ss3.2}
\label{ss3.2-partial-actions-decayability}
If $X$ is a proper open subset of a normed space $E$, then $X$ is not
closed under the group of translations $\bbT(E)$ of $E$.
So there is no natural action of $\bbT(E)$ on $X$.
But for every $x \in X$ there are neighborhoods $B_x$ and $V_x$ of
$x$ in $X$ and $\rfs{Id}^E$ in $\bbT(E)$ respectively
such that the action of every $\rfs{tr}_v \in V_x$
on $B_x$ is defined.
Moreover, $H(X)$ contains a homeomorphism which coincides
with $\rfs{tr}_v$ on $B_x$ and which is the identity outside
some bigger neighborhood of $x$. Indeed, even $\rfs{LIP}(X)$ contains
such a homeomorphism.
We shall use such homeomorphisms. To this end we introduce two notions:
the notion of a partial action of a topological group on a topological
space, and the notion of decayability of partial actions.

\begin{defn}\label{metr-bldr-d3.1}\label{was-d3.16}
\begin{rm}
(a) Let $X$ be a topological space and $x \in X$.
Set
$\rfs{Nbr}^X(x) \eqdf \setm{U}{x \in U \subseteq X
\mbox{ and $U$ is open}}$
   \index{N@nbr@@$\rfs{Nbr}^X(x) = \setm{U}{x \in U \subseteq X
          \mbox{ and $U$ is open}}$}
and
$\rfs{MBC} =
\setm{\alpha \in \rfs{MC}}{\rfs{Id}_{[0,\infty)} \leq \alpha}$.
   \index{N@mbc@@$\rfs{MBC} =
          \setm{\alpha \in \rfs{MC}}{\rfs{Id}_{[0,\infty)} \leq
	  \alpha}$}
Let $\alpha \in \rfs{MBC}$, $X,Y$ be metric spaces
and $\iso{\tau}{X}{Y}$.
We say that $\tau$ is
{\it $\alpha$-bicontinuous},
if $\tau$ and $\tau\inverse$ are $\alpha$-continuous.
Let $x \in X$. We say that
$\tau$ is {\it $\alpha$-continuous at $x$}, if for
some $U \in \rfs{Nbr}(x)$,
$\tau \nrestriction U$ is $\alpha$-continuous.
Also, $\tau$ is said to be {\it $\alpha$-bicontinuous at $x$}, if for
some $U \in \rfs{Nbr}(x)$, $\tau \nrestriction U$
is $\alpha$-bicontinuous.
Let $\itGamma \subseteq \rfs{MC}$.
We say that $\tau$ is {\it $\itGamma$-continuous
(resp.\ $\itGamma$-bicontinuous) at $x$} if for some
$\alpha \in \itGamma$, $\tau$ is $\alpha$-continuous
(resp.\ $\alpha$-bicontinuous) at $x$.
   \index{bicontinuous. $\alpha$-bicontinuous homeomorphism}
   \index{continuous. $\alpha$-continuous at $x$.
There is $U \in \rfs{Nbr}(x)$ such that $f \nrestriction U$
is $\alpha$-continuous}
   \index{bicontinuous. $\alpha$-bicontinuous at $x$}
   \index{continuous. $\itGamma$-continuous at $x$. There is
          $\alpha \in \itGamma$ such that
	  $f$ is $\alpha$-continuous at $x$}
   \index{bicontinuous. $\itGamma$-bicontinuous at $x$}

If $H$ is a group, then $e_H$ denotes the unit of $H$.
   \index{N@e@@$e_H$. The unit of a group $H$}

(b) Let $H$ be a topological group, $X$ be a toplogical space and
$\lambda$ be a function such that
$\rfs{Dom}(\lambda) \subseteq H \times X$ and
$\rfs{Rng}(\lambda) \subseteq X$. We say that $\lambda$
is a {\it partial action} of $H$ on $X$,
if the following conditions hold.
\begin{itemize}
\item[(1)] $\lambda$ is continuous.
\item[(2)] $\rfs{Dom}(\lambda)$ is open in $H \times X$.
\item[(3)] For $g \in H$ let $g_{\lambda}$ be the function defined by
$g_{\lambda}(x) = \lambda(g,x)$. Then $g_{\lambda}$ is
a homeomorphism between $\rfs{Dom}(g_{\lambda})$ and
$\rfs{Rng}(g_{\lambda})$.
\item[(4)] 
$(e_H)_{\lambda} = \rfs{Id}_{\srfs{Dom}((e_H)_{\lambda})}$.
\item[(5)] For every $g \in H$,
$(g\inverse)_{\lambda} = (g_{\lambda})\inverse$.
\item[(6)] For every $g,h \in H$ and $x \in X$:
if $g_{\lambda}(x)$ and $h_{\lambda}(g_{\lambda}(x))$ are defined,
then $(hg)_{\lambda}(x)$
is defined and 
$(hg)_{\lambda}(x) = h_{\lambda}(g_{\lambda}(x))$.
\end{itemize}
   \index{partial action of a topological group on a topological space}

Define $\rfs{Fld}(\lambda) \eqdf \rfs{Dom}((e_H)_{\lambda})$.
Note that by (5) and (6),
$\rfs{Dom}(g_{\lambda}) \subseteq \rfs{Fld}(\lambda)$
for every $g \in H$.
   \index{N@fld@@$\rfs{Fld}(\lambda) = \rfs{Dom}(e_{\lambda})$}

(c) Let $\alpha \in \rfs{MBC}$, $a \in (0,1)$,
$H$ be a topological group,
$\lambda$ be a partial action of $H$ on a metric space $X$,
$G \leq H(X)$
and $x \in \rfs{Fld}(\lambda)$.
{\thickmuskip=4mu \medmuskip=2mu \thinmuskip=1mu 
Then $\lambda$ is called an
{\it $(a,\alpha,G)$-decayable action at~$x$},
if there is $r_x > 0$ such that for every
$r \in (0,r_x)$ there is
$V = V_{x,r} \in \rfs{Nbr}(e_H)$ such that:
}
\begin{itemize}
\item[(i)] $V \times B(x,ar) \subseteq \rfs{Dom}(\lambda)$;
\item[(ii)] for every $h \in V$ there is $g \in G$ such that:
$g$ is $\alpha$-bicontinuous,
$g \nrestriction B(x,ar) =\break h_{\lambda} \nrestriction B(x,ar)$
and $\rfs{supp}(g) \subseteq B(x,r)$.
\end{itemize}
   \index{decayable action.
          $\lambda$ is an $(a,\alpha,G)$-decayable action at $x$}

Let $A \subseteq \rfs{Fld}(\lambda)$.
We say that $\lambda$ is an
{\it $(a,\alpha,G)$-decayable action in $A$,}
if it is $(a,\alpha,G)$-decayable
at every $x \in A$;
   \index{decayable action. $\lambda$ is an
          $(a,\alpha,G)$-decayable action in $A$}
$\lambda$ is {\it $(a,\alpha,G)$-decayable}
if it is $(a,\alpha,G)$-decayable in $\rfs{Fld}(\lambda)$.
Suppose that $\itGamma$ is a modulus of continuity.
Then $\lambda$ is called {\it $(a,\itGamma,G)$-decayable}
if $\lambda$ is $(a,\alpha,G)$-decayable
for some $\alpha \in \itGamma$.
   \index{decayable action. $\lambda$ is an
          $(a,\alpha,G)$-decayable action}
   \index{decayable action. $\lambda$ is
          an $(a,\itGamma,G)$-decayable action}

If in the above $a = \dghalf$, then we omit its mention.
So ``$\lambda$ is $(\alpha,G)$-decayable at $x$'' means
``$\lambda$ is $(\dghalf,\alpha,G)$-decayable at $x$'' etc.
If $a = \dghalf$ and $G = H(X)$, then we omit the mention
of $a$ and $G$.
So ``$\lambda$ is $\alpha$-decayable at $x$'' means
``$\lambda$ is $(\dghalf,\alpha,H(X))$-decayable at $x$'',
``$\lambda$ is $\alpha$-decayable in $A$'' means
``$\lambda$ is $(\dghalf,\alpha,H(X))$-decayable in $A$'' etc.

   \index{decayable action. $(\alpha,G)$-decayable at $x$. This means
$(\dghalf,\alpha,G)$-decayable at $x$}
   \index{decayable action. $\alpha$-decayable at $x$. This means
$(\dghalf,\alpha,H(X))$-decayable at $x$}

(d) Let $\lambda$ be a partial action of a topological group $H$
on a topological space $X$, $A \subseteq H$ and $x \in X$.
We write $A_{\lambda}(x) = \setm{h_{\lambda}(x)}{h \in A}$.
We say that $x$ is a {\it $\lambda$-limit-point},
if $x \in \rfs{acc}(V_{\lambda}(x))$
for every $V \in \rfs{Nbr}(e_{H})$.
   \index{N@AAAA@@$A_{\lambda}(x) = \setm{h_{\lambda}(x)}{h \in A}$}
   \index{limit-point. $\lambda$-limit-point}
\hfill\proofend
\end{rm}
\end{defn}

Note that if $\lambda$ is $(a,\alpha,G)$ decayable partial action
of $H$ at $x$, then there are $V \in \rfs{Nbr}(e^H)$ and
$U \in \rfs{Nbr}(x)$ such that $h_{\lambda} \nrestriction U$
is $\alpha$-bicontinuous for every $h \in V$.

The partial actions appearing in this section are obtained by
restricting a full group action on a space $E$ to an open subset
of $E$.
This is described (a) below.

\begin{prop}\label{metr-bldr-p3.2}
\num{a} Suppose that $\lambda$ is a partial action of a topological
group $H$ on a topological space $E$.
Let $X \subseteq \rfs{Fld}(\lambda)$ be open,
and define $\lambda \drest X$ by setting
$\rfs{Dom}(\lambda \drest X) =
\setm{\pair{h}{x}}{h \in H \mbox { and } x,h_{\lambda}(x) \in X}$
and $(\lambda \drest X)(h,x) = \lambda(h,x)$.
Then $\lambda \drest X$ is a partial action of $H$ on $X$.
   \index{N@AAAA@@$\lambda \vdrest X$.
The restriction of a partial group action $\lambda$ to an open set $X$}

\num{b} Let $\lambda$ be a partial action of $H$ on $X$,
$G \leq H(X)$,
$D \subseteq C \subseteq \rfs{Fld}(\lambda)$,
$a \in (0,1)$, $\alpha \in \rfs{MBC}$, $r_0 > 0$
and let $V_r \in \rfs{Nbr}(e_H)$ for every $r \in (0,r_0)$.
Assume that: \num{i} $D$ is a dense subset of $C$,
\num{ii} $\lambda$ is $(a,\alpha,G)$-decayable in~$D$,
\num{iii} $r_x \geq r_0$ for every $x \in D$,
\num{iv} $V_{x,r} \supseteq V_r$
for every $x \in D$ and $r \in (0,r_0)$.
\underline{Then} $\lambda$ is $(a,\alpha,G)$-decayable in $C$,
$r_x \geq r_0$ for every $x \in C$,
and $V_{x,r} \supseteq V_r$ for every $x \in C$ and $r \in (0,r_0)$.
\end{prop}

\noindent
{\bf Proof } The proof of both parts is trivial.
\smallskip\hfill\myqed

Suppose that $X$ is an open subset of a normed space $E$.
We shall be interested in two partial actions on $X$:
the partial action of the group $\bbT(E)$ of translations of $E$,
and the partial action of the group $\bbA(E)$ of affine
transformations of $E$.
We need to know that these partial actions are decayable.
In fact, we shall show that $\bbA(E)$ is $(\alpha,G)$-decayable,
where $\alpha(t) = 15t$, and $G$ is any group containing
$\rfs{LIP}(X)$.

Obviously, the decayability of $\bbA(E)$ implies the decayability of
both $\bbT(E)$ and the group of bounded linear automorphisms of $E$.
Because we deal with groups containing $\rfs{LIP}(X;F)$,
we shall really need to show that
$\setm{T \in \bbA(E)}{T(F) = F}$ is decayable with respect to any
group $G$ containing $\rfs{LIP}(X;F)$.

\begin{defn}\label{metr-bldr-d3.3}\label{bddly-lip-bldr-d5.4}
\begin{rm}
(a)
Let $E$ be a normed space and $v \in E$.
Define $\rfs{tr}_v^E(x) \eqdf v + x$ and
\hbox{$\bbT(E) = \setm{\rfs{tr}_v^E}{v \in E}$.}
Whenever $E$ can be understood from the context,
we abbreviate $\rfs{tr}_v^E$ by $\rfs{tr}_v$.
   \index{N@t00@@$\bbT(E) = \setm{\rfs{tr}_v^E}{v \in E}$.
          The group of translations of $E$}
We define $d(\rfs{tr}_u,\rfs{tr}_v) = \norm{u - v}$.

(b) Let $E$ be a normed space and $x \in X$.
Denote the group of bounded linear automorphisms of $E$ by $\bbL(E)$
and set $\bbL(E,x) = (\bbL(E))^{\srfs{tr}_x^E}$.
For $S,T \in \bbL(E)$ define
$d(S,T) = \norm{S - T} + \norm{S\inverse - T\inverse}$.
Let
$
\bbA(E) \eqdf \setm{\rfs{tr}^E_{v} \scirc T}
{v \in E, \ T \in \bbL(E)}.
$
   \index{N@l00@@$\bbL(E)$. The group of bounded linear automorphisms
          of $E$}
   \index{N@l01@@$\bbL(E,x) =
          (\bbL(E))^{\srfs{tr}_x^E}$}
   \index{N@a00@@$\bbA(E)$. The group of affine automorphisms of $E$}
That is, $\bbA(E)$ is the group of bounded affine transformations
of $E$.
Suppose that $A = \rfs{tr}^E_{v} \scirc T \in \bbA(E)$.
Then $v$ and $T$ are uniquely determined by $A$.
We denote $v = v_A$ and $T = T_A$.
We may thus define
$$
d(A_1,A_2) = \norm{v_{A_1} - v_{A_2}} +
\norm{T_{A_1} - T_{A_2}} + \norm{T_{A_1}\inverse - T_{A_2}\inverse}.
$$
Then $d$ is a metric on $\bbA(E)$,
$\pair{\bbA(E)}{d}$ is a topological group, and the action of $\bbA(E)$
on $E$ is continuous.
Note that $\bbL(E,x) \leq \bbA(E)$ and the function
$T \mapsto T^{\srfs{tr}_x}, \ T \in \bbL(E)$, is a topological
isomorphism between $\bbL(E)$ and $\bbL(E,x)$.

Let $\lambda_{\bbT}^E$, $\lambda_{\bbL}^E$, $\lambda_{\bbL}^{E,x}$,
$\lambda_{\bbA}^{E,x}$
denote respectively the natural actions of $\bbT(E)$, $\bbL(E)$,
$\bbL(E,x)$ and $\bbA(E)$ on $E$.
   \index{N@lt00@@$\lambda_{\bbT}^E$,
          $\lambda_{\bbL}^E$, $\lambda_{\bbL}^{E,x}$,
	  $\lambda_{\bbA}^E$. Actions of $\bbT(E)$,
	  $\bbL(E)$, $\lambda_{\bbL}^{E,x}$ and $\bbA(E)$ on $E$}

(c) Suppose that $E$ is a normed space, $F$ is a linear subspace
of $E$ and $x \in F$. Define
\smallskip

\indent\indent\indent
$\bbT(E;F) = \setm{\rfs{tr}^E_v}{v \in F}$, \ \ %
$\bbL(E;F) = \setm{T \in \bbL(E)}{T(F) = F}$,
\vspace{1mm}

\indent\indent\indent
$\bbA(E;F) = \setm{A \in \bbA(E)}{A(F) = F}$, \ %
$\bbL(E,x;F) = (\bbL(E;F))^{\srfs{tr}^E_x}$.
\vspace{0.5mm}
\smallskip
\newline
The groups $\bbT(E;F)$, $\bbL(E;F)$, $\bbL(E,x;F)$
and $\bbA(E;F)$ equipped with the metric
they inherit from $\bbT(E)$, $\bbL(E)$, $\bbL(E,x)$ and $\bbA(E)$
respectively are metric topological groups.
\kern-19pt
   \index{N@t01@@$\bbT(E;F) = \setm{\rfs{tr}^E_v}{v \in F}$}
   \index{N@l02@@$\bbL(E;F) = \setm{T \in \bbL(E)}{T(F) = F}$}
   \index{N@a01@@$\bbA(E;F) = \setm{A \in \bbA(E)}{A(F) = F}$}

If $\lambda$ is a partial action of $H$ on $X$ and $H_1 \leq H$, let
$\lambda \drest H_1$ denote the restriction of $\lambda$ to~$H_1$.
   \index{N@AAAA@@$\lambda \vdrest H_1$.
          The restriction of a partial action $\lambda$ of $H$
          to a subgroup
$H_1$ of $H$}
Let $\lambda_{\bbT}^{E;F} = \lambda_{\bbT}^E \drest \bbT(E;F)$.
$\lambda_{\bbL}^{E;F}$, $\lambda_{\bbL}^{E,x;F}$
and $\lambda_{\bbA}^{E;F}$
are defined in a similar way.
   \index{N@lt01@@$\lambda_{\bbT}^{E;F} =
          \lambda_{\bbT}^E \vdrest \bbT(E;F)$}
   \index{N@ll01@@$\lambda_{\bbL}^{E;F} =
          \lambda_{\bbL}^E \vdrest \bbL(E;F)$}
   \index{N@ll02@@$\lambda_{\bbL}^{E,x;F} =
          \lambda_{\bbL}^{E,x} \vdrest \bbL(E,x;F)$}
   \index{N@la00@@$\lambda_{\bbA}^{E;F} =
          \lambda_{\bbA}^E \vdrest \bbA(E;F)$}

(d) Suppose that $X$ is a topological space and $F$ is a set. Define
\smallskip

\indent\indent\indent
$H(X;F) \eqdf \setm{h \in H(X)}{h(X \cap F) = X \cap F}$.
   \index{N@h02@@$H(X;F) = \setm{h \in H(X)}{h(X \cap F) = X \cap F}$}
\end{rm}
\end{defn}

\begin{prop}\label{metr-bldr-p3.4}
Let $E$ be a normed space, $X \subseteq E$ be open, $\calS$
be an open cover of $X$,
$\calF$ be a subspace choice for $\calS$, $S \in \calS$,
$G = \rfs{LIP}(X;S,F_S)$
and $\alpha(t) = 3t$.
\underline{Then} $\lambda_{\bbT}^{E;F_S} \drest S$
is $(\dgfrac{5}{8},\alpha,G)$-decayable.
In particular,
$\lambda_{\bbT}^{E;F_S} \drest S$
is
$(\alpha,G)$-decayable.
\end{prop}

\noindent
{\bf Proof }
We show that if $x \in S \cap F_S$,
then  $\lambda_{\bbT}^{E;F_S} \drest S$
is  $(\dgfrac{5}{8},\alpha,G)$-decayable
at~$x$, $r_x = d(x,E - S)$, and for every $r \in (0,r_x)$,
$V_{x,r} = B^{\bbT(E;F_S)}(\rfs{Id}_E,\dgfrac{r}{4})$.

Let $r < r_x$.
Let $\rfs{tr}^E_v \in V_{x,r}$.
So $v \in F_S$ and $\norm{v} < \dgfrac{r}{4}$.
We apply Lemma \ref{l2.6}(b). Choose $r_0$ of \ref{l2.6}(b) to be
$r$, choose $r$ and $s$ of \ref{l2.6}(b) to be $\dgfrac{5r}{8}$ and
$v$ of \ref{l2.6}(b) to be $v$.
Let $h$ be as assured by \ref{l2.6}(b).
By \ref{l2.6}(b)(ii), $h$ is
$(1 + \frac{\norm{v}}{r - \dgfrac{5r}{8} - \norm{v}})$-bilipschitz.
$(1 + \frac{\norm{v}}{r - \dgfrac{5r}{8} - \norm{v}}) < 3$.
Hence $h$ is $3$-bilipschitz.
It follows from \ref{l2.6}(b)(ii) that $h$ is as required.
By Proposition \ref{metr-bldr-p3.2}(b),
$\lambda_{\bbT}^{E;F_S} \drest S$ is $(\alpha,G)$-decayable.
\hfill\myqed

\begin{prop}\label{metr-bldr-p3.5}
Let $\fnn{\eta}{[0,\infty)}{[0,1]}$.
Suppose that $\eta$ is $K$-Lipschitz and that $\eta(t) = a$
for every $t \geq a$.
Let $E$ be a normed space. Define $\fnn{g}{E}{E}$ by
$g(x) = \eta(\norm{x}) \ncdot x$.
Then $g$ is $(1 + Ka)$-Lipschitz.
\end{prop}

\noindent
{\bf Proof }
Let $x,y \in E$. If $\norm{x},\norm{y} \geq a$,
then $g(x) = x$ and $g(y) = y$, and hence
$\norm{g(x) - g(y)} = \norm{x - y}$.
Assume that $\norm{x} \leq a$ or $\norm{y} \leq a$.
Without loss of generality $\norm{y} \leq a$. Hence
\vspace{1.5mm}
\newline
\rule{7pt}{0pt}
\renewcommand{\arraystretch}{1.5}
\addtolength{\arraycolsep}{-6pt}
$
{\thickmuskip=3.8mu \medmuskip=1.1mu \thinmuskip=1mu 
\begin{array}{ll}
&
\norm{g(x) - g(y)} = 
\norm{\kern3pt\eta(\norm{x}) \ncdot x - \eta(\norm{y}) \ncdot y\kern3pt}
\\
\leq
\rule{5pt}{0pt}
&
\norm{\kern3pt\eta(\norm{x}) \ncdot x - \eta(\norm{x}) \ncdot y\kern3pt}
+
\norm{\kern3pt\eta(\norm{x}) \ncdot y - \eta(\norm{y}) \ncdot y\kern3pt}
\\
=
\rule{5pt}{0pt}
&
\eta(\norm{x}) \ncdot \norm{x - y} +
\abs{\kern3pt\eta(\norm{x}) - \eta(\norm{y})\kern3pt} \ncdot \norm{y}
\leq
\norm{x - y} + K \ncdot \norm{x - y} \ncdot \norm{y} \leq
(1 + Ka) \ncdot \norm{x - y}.
\vspace{1.7mm}
\end{array}
}
$
\renewcommand{\arraystretch}{1.0}
\addtolength{\arraycolsep}{6pt}
\rule{7pt}{0pt}\hfill\myqed

\begin{prop}\label{metr-bldr-p3.6}
Let $E$ be a normed space, $T \in \bbL(E)$,
$\fnn{\eta}{[0,\infty)}{[0,1]}$ and $a~>~0$.
Set $\rfs{Id}_E = I$.
Suppose that $\eta$ is $K$-Lipschitz,
$\eta(t) = t$ for every $t \geq a$
and $\norm{I - T} (1 + Ka) < 1$.
Define $\fnn{h}{E}{E}$ by
$$h(x) = (1 - \eta(\norm{x})) \ncdot T(x) + \eta(\norm{x}) \ncdot x.
$$
Then
\begin{itemize}
\addtolength{\parskip}{-11pt}
\addtolength{\itemsep}{06pt}
\item[\num{i}]
$h \in H(E)$,
$h$ is
$(\norm{T} + \norm{I - T} \ncdot (1 + Ka))$-Lipschitz,
and $h\inverse$ is
$\max(\frac{\norm{T\inverse}}{1 - \norm{I - T} \ncdot (1 + Ka)},
1)$-Lipschitz.
\item[\num{ii}]
If $F$ is a linear subspace of $E$,
and $T \in \bbL(E;F)$, then $h \in H(E;F)$.
\vspace{-05.7pt}
\end{itemize}
\end{prop}

\noindent
{\bf Proof }
(i) We prove that $h$ is Lipschitz.
Let $x,y \in E$. Then
\vspace{1.5mm}
\newline
\rule{7pt}{0pt}
\renewcommand{\arraystretch}{1.5}
\addtolength{\arraycolsep}{-4pt}
$
\begin{array}{ll}
&
h(x) - h(y) =
(1 - \eta(\norm{x})) \ncdot T(x) + \eta(\norm{x}) \ncdot x -
((1 - \eta(\norm{y})) \ncdot T(y) + \eta(\norm{y}) \ncdot y)
\\
=
\rule{01pt}{0pt}
&
T(x - y) + (I - T)(\eta(\norm{x}) \ncdot x - \eta(\norm{y}) \ncdot y).
\vspace{1.7mm}
\end{array}
$
\renewcommand{\arraystretch}{1.0}
\addtolength{\arraycolsep}{4pt}
\newline
By Proposition \ref{metr-bldr-p3.5},
{\thickmuskip=2.8mu \medmuskip=1.2mu \thinmuskip=1mu 
$$
\norm{h(x) - h(y)} \leq
\norm{T} \ncdot \norm{x - y} +
\norm{I - T} \ncdot (1 + Ka) \ncdot \norm{x - y} \leq
(\norm{T} + \norm{I - T} \ncdot (1 + Ka)) \ncdot \norm{x - y}.
$$
}
Hence $h$ is $(\norm{T} + \norm{I - T} \ncdot (1 + Ka))$-Lipschitz.

We prove that $h\inverse$ is Lipschitz.
Let $x,y \in E$.
By the above,
\vspace{1.5mm}
\newline
\rule{06pt}{0pt}
\renewcommand{\arraystretch}{1.5}
\addtolength{\arraycolsep}{-5pt}
$
\begin{array}{ll}
&
T\inverse(h(x) - h(y)) =
(x - y) +
T\inverse(I - T)(\eta(\norm{x}) \ncdot x - \eta(\norm{y}) \ncdot y)
\\
=
\rule{2pt}{0pt}
&
(x - y) + (T - I)(\eta(\norm{x}) \ncdot x - \eta(\norm{y}) \ncdot y).
\vspace{1.7mm}
\end{array}
$
\renewcommand{\arraystretch}{1.0}
\addtolength{\arraycolsep}{5pt}
\newline
So
{\thickmuskip=2mu \medmuskip=1mu \thinmuskip=1mu 
\vspace{1.5mm}
\newline
\rule{7pt}{0pt}
\renewcommand{\arraystretch}{1.5}
\addtolength{\arraycolsep}{-2pt}
$
\begin{array}{ll}
&
\norm{T\inverse} \ncdot \norm{h(x) - h(y)} \geq
\norm{T\inverse(h(x) - h(y))} \geq
\norm{x - y} -
\norm{(T - I)(\eta(\norm{x}) \ncdot x - \eta(\norm{y}) \ncdot y)}
\\
\geq
&
\norm{x - y} - \norm{(T - I)} \ncdot (1 + Ka) \ncdot \norm{x - y} =
(1 - \norm{T - I} \ncdot (1 + Ka)) \ncdot \norm{x - y}.
\vspace{0.7mm}
\end{array}
$
}
\renewcommand{\arraystretch}{1.0}
\addtolength{\arraycolsep}{2pt}
\newline
That is,
$\norm{x - y} \leq
\frac{\norm{T\inverse}}
{1 - \norm{T - I} \ncdot (1 + Ka)} \mcdot \norm{h(x) - h(y)}$.
\smallskip

(ii) Let $x \in F$.
Set $T_x = (1 - \eta(\norm{x}) T + \eta(\norm{x}) I$.
Then $h(x) = T_x(x)$ and $T_x(F) = F$.
\rule{8pt}{0pt}\hfill\myqed

\begin{lemma}\label{metr-bldr-l3.7}
Let $E$ be a normed space, $X \subseteq E$ be open, $\calS$ be an open
cover of $X$,
$\calF$ be a subspace choice for $\calS$, $S \in \calS$,
$x \in S \cap F_S$,
$G = \rfs{LIP}(X;S,F_S)$
and $\alpha(t) = 5t$.
Then $\lambda_{\bbL}^{E,x;F_S} \drest S$ is $(\alpha,G)$-decayable
at $x$, $r_x = d(x,E - S)$,
and $V_{x,r} = (B^{\bbL(E;F)}(\rfs{Id}_E,\dgquarter)^{\srfs{tr}^E_x}$
for every $r \in (0,r_x)$.
\end{lemma}

\noindent
{\bf Proof }
We may assume that $0^E \in S$ and $x = 0^E$.
Set $I = \rfs{Id}_E$.
Let $r_0 = d(0^E,E - S)$ and $V = B^{\bbL(E;F_S)}(I,\dgquarter)$.
Let $r < r_0$ and $T \in V$.
We show that $T$ is ``decayable''.
Define
$\fnn{\eta(t)}{[0,\infty)}{[0,1]}$ to be the following piecewise
linear function. The breakpoints of $\eta$ are $\dgfrac{r}{2}$ and $r$;
$\eta(t) = 0$ for every $t \in [0,\dgfrac{r}{2}]$
and $\eta(t) = 1$ for every $t \geq r$.
Clearly, $\eta$ is $\dgfrac{2}{r}$-Lipschitz.

Define $\fnn{h}{E}{E}$ by
$h(y) = (1 - \eta(\norm{y})) \ncdot T(y) +
\eta(\norm{y}) \ncdot y$.
We check that Proposition~\ref{metr-bldr-p3.6} applies to $h$.
Set $K = \dgfrac{2}{r}$. So $\eta$ is $K$-Lipschitz.
Since $\norm{I - T} < \dgquarter$ and $Ka = \frac{2}{r} \ncdot r = 2$,
it follows that
$\norm{I - T} \ncdot (1 + Ka) < \quarter \ncdot (1 + 2) =
\dgfrac{3}{4} < 1$.
It thus follows from \ref{metr-bldr-p3.6}(i) that
$h \in H(E)$
and $h$ is $\norm{T} + \norm{I - T} \ncdot (1 + Ka)$-Lipschitz.
By the above,
$\norm{T} + \norm{I - T} \ncdot (1 + Ka) <
\dgfrac{5}{4} + \dgfrac{3}{4} = 2$.
So $h$ is $2$-Lipschitz.
\rule{0pt}{12pt}\kern-0.5pt
Since $\norm{T\inverse} < \dgfrac{5}{4}$, it follows that
$\frac{\norm{T\inverse}}{1 - \norm{I - T} \ncdot (1 + Ka)} <
\frac{\dgfrac{5}{4}}{1 - \dgfrac{3}{4}} = 5$.
By \ref{metr-bldr-p3.6}(i), $h\inverse$ is $5$-Lipschitz.
So $h$ is $5$-bilipschitz.

\rule{0pt}{12pt}
Clearly, $\rfs{supp}(h) \subseteq B(0^E,r) \subseteq X$.
So $h \nrestriction X \in H(X)$.
Also,
$h \nrestriction B(0^E,\dgfrac{r}{2})
= T \nrestriction B(0^E,\dgfrac{r}{2})$.
By \ref{metr-bldr-p3.6}(ii), $h(E \cap F_S) = F_S$.
Hence $h \nrestriction X$ is as required.
\hfill\myqed

\begin{lemma}\label{metr-bldr-l3.8}
Let $E$ be a normed space, $X \subseteq E$ be open, $\calS$ be an open
cover of $X$,
$\calF$ be a subspace choice for $\calS$,
$\rfs{LIP}(X;\calS,\calF) \leq G \leq H(X)$
and $\alpha(t) = 15 t$.
Let $S \in \calS$.
Then $\lambda_{\bbA}^{E;F_S} \drest S$
is $(\alpha,G)$-decayable.
\end{lemma}

\noindent
{\bf Proof } 
Set $I = \rfs{Id}_E$.
Let $x \in S \cap F_S$, $r_x = d(x,E - S)$ and $r \in (0,r_x)$.
If $x \neq 0^E$ let
$a_r = \min(\dgquarter,\dgfrac{r}{8},\frac{r}{8 \norm{x}})$ and
if $x = 0^E$ let $a_r = \min(\dgquarter,\dgfrac{r}{8})$.
Let
$V_{x,r} = B^{\bbA(E;F)}(I,a_r)$.
We show that
\begin{equation}
\tag{$*$}
V_{x,r} \subseteq B^{\bbT(E;F_S)}(I,\dgfrac{r}{4}) \scirc
(B^{\bbL(E;F_S)}(I,\dgquarter))^{\srfs{tr}^E_x}.
\end{equation}
If $A \in \bbA(E;F)$, then $A$ can be uniquely represented in the form
$A = \rfs{tr}_{u_{A,x}} \scirc (T_{A,x})^{\srfs{tr}_x}$,
where $T_{A,x} \in \bbL(E;F)$.
Let $A = \rfs{tr}_{v_A} \scirc T_A$, where $T_A \in \bbL(E;F)$.
Then $T_{A,x} = T_A$ and
$u_{A,x} = v_A + (T_A - I)(x)$.
Set $T = T_A$, $v = v_A$ and $u = u_{A,x}$.
Suppose that $A \in V_{x,r}$.
Then
$d(T,I) < a_r < \dgquarter$. So $T \in B^{\bbL(E;F_S)}(I,\dgquarter)$.
Hence
$T^{\srfs{tr}_x} \in (B^{\bbL(E;F_S)}(I,\dgquarter))^{\srfs{tr}_x}$.
Suppose that $x \neq 0$. Then
$
\norm{u} \leq \norm{v} + \norm{T - I} \ncdot \norm{x} \leq
\dgfrac{r}{8} + \frac{r}{8 \norm{x}} \ncdot \norm{x} = \dgfrac{r}{4}.
$
If $x = 0$, then $u = v$.
So $\norm{u} < \dgfrac{r}{4}$.
In both cases $u \in B^{\bbT(E;F_S)}(I,\dgfrac{r}{4})$.
This proves ($*$).

Let $A \in V_{x,r}$. Let $T$ and $u$ be as above.
By Lemma \ref{metr-bldr-l3.7},
there is
$h_1 \in H(X;F_S)\sprtl{B(x,r)}$ such that
$h_1 \nrestriction B(x,\dgfrac{r}{2}) =
T^{\srfs{tr}_x} \nrestriction B(x,\dgfrac{r}{2})$
and $h_2$ is $5$-bilipschitz.
By Proposition \ref{metr-bldr-p3.4},
there is $h_2 \in H(X;F_S)\sprtl{B(x,r)}$ such that
$h_2 \nrestriction B(x,\dgfrac{5r}{8}) =
\rfs{tr}_u \nrestriction B(x,\dgfrac{5r}{8})$
and $h_1$ is $3$-bilipschitz.
Let $h = h_2 \scirc h_1$. So $h \in H(X;F_S)$,
$\rfs{supp}(h) \subseteq B(x,r)$ and $h$ is $15$-bilipschitz.
It remains to show that
$h \nrestriction B(x,\dgfrac{r}{2}) =
A \nrestriction B(x,\dgfrac{r}{2})$.
Let
$y \in B(x,\dgfrac{r}{2})$. Then $h_1(y) = T^{\srfs{tr}_x}(y)$.
Since $T \in B^{\bbL(E;F_S)}(I,\dgquarter)$,
$\norm{T} \leq \dgfrac{5}{4}$.
So
$\norm{T(y - x)} \leq \frac{5}{4} \norm{y - x}$.
That is,
$d(T(y - x),0) \leq \frac{5}{4} \norm{y - x}$.
Since $\rfs{tr}_x$ is an isometry,
$d(T^{\srfs{tr}_x}(\rfs{tr}_x(y - x)),\rfs{tr}_x(0)) \leq
\frac{5}{4} \norm{y - x}$.
That is,
$\norm{T^{\srfs{tr}_x}(y) - x} \leq \frac{5}{4} \norm{y - x}$.
Since $y \in B(x,\dgfrac{r}{2})$,
$\norm{T^{\srfs{tr}_x}(y) - x} \leq \dgfrac{5r}{8}$.
Hence
$h_2(T^{\srfs{tr}_x}(y)) = \rfs{tr}_u(T^{\srfs{tr}_x}(y))$.
So $h(y) = h_2(h_1(y)) = A(y)$.
We have shown that if $x \in S \cap F_S$, then
$\lambda_{\bbA}^{E;F_S} \drest S$ is $(\alpha,G)$-decayable at $x$.

Let $x \in S - F_S$. Then $x \in \rfs{acc}(S \cap F_S)$.
Define $r_x = \half d(x,E - S)$. For $r \in (0,r_x)$ let\break
$a_r = \half \min(\dgquarter,\dgfrac{r}{8},\frac{r}{8 \norm{x}})$
and $V_{x,r} = B^{\bbA(E;F)}(x,a_r)$.
Let $D = B(x,\dgfrac{r}{3}) \cap F_S$.
By the above argument, for every $y \in D$:
$\lambda_{\bbA}^{E;F_S} \drest S$
is $(\alpha,G)$-decayable at $y$,
$r_y \geq r_x$,
and $V_{y,r} \supseteq V_{x,r}$ for every $r \in (0,r_x)$.
By Proposition \ref{metr-bldr-p3.2}(b),
$\lambda_{\bbA}^{E;F_S} \drest S$ is $(\alpha,G)$-decayable at $x$.
\smallskip\hfill\myqed

Recall that in this section we prove that
if $(H_{\itGamma}^{\srfs{LC}}(E))^{\tau} = H_{\itGamma}^{\srfs{LC}}(F)$,
then $\tau$ is locally $\itGamma$-bicontinuous.
If $\itGamma$ is countably generated or if $\itGamma = \rfs{MC}$,
then the above is true for any $E$ and $F$.
For $\itGamma$'s which are not countably generated,
we have only a partial answer.
We know how to prove that $\tau$ is locally $\itGamma$-bicontinuous
only for $\itGamma$'s which are $\kappa(E)$-generated.
See the definition below.

\begin{defn}\label{metr-bldr-d3.9}
\begin{rm}
(a) Let $X$ be a metric space and $r > 0$.
A family ${\cal A}$ of subsets of $X$ is {\it $r$-spaced},
if $d(A,B) \geq r$ for any distinct $A,B \in {\cal A}$.
A subset $C \subseteq X$ is {\it $r$-spaced}, if
$\setm{\sngltn{x}}{x \in C}$ is $r$-spaced.
A set $C$ is {\it spaced}, if $C$ is $r$-spaced for some $r > 0$.
   \index{spaced set of sets. $r$-spaced set of sets}
   \index{spaced subset of $X$. $A \subseteq X$ is spaced if
          $(\exists r > 0)(\forall x,y \in A)((x \neq y)
          \rightarrow (d(x,y) \geq r))$}

(b) Let $X$ be a metric space $x \in X$ and $A \subseteq X$.
We define the set of cardinals $\Kappa^X(x,A)$ as follows:
$\kappa \in \Kappa^X(x,A)$ iff for every $U \in \rfs{Nbr}(x)$ there is
$B \subseteq A \cap U$ such that
$\abs{B} = \kappa$ and $B$ is spaced.
Let
$$
\kappa^X(x,A) = \sup(\Kappa^X(x,A)), \ \ \ %
\kappa(X) = \min\limits_{x \in X} \kappa^X(x,X).
$$

   \index{N@k02@@$\Kappa^X(x,A) =
          \setm{\kappa}{(\forall U \in \rfs{Nbr}(x))
          (\exists B \subseteq A \cap U)
          (\abs{B} = \kappa \mbox{ and } B \mbox{ is spaced}}$}
   \index{N@k03@@$\kappa^X(x,A) = \sup(\Kappa^X(x,A))$}
   \index{N@k04@@$\kappa(X) = \min_{x \in X} \kappa^X(x,X)$}

(c) Let $\itGamma$ be a modulus of continuity.
We say that $\itGamma_0$ {\it generates} $\itGamma$ if
$\itGamma = \rfs{cl}_{\preceq}(\itGamma_0)$.
We say that $\itGamma$ is
{\it $(\leq\kern-3pt \kappa)$-generated} if there is
$\itGamma_0$ such that $\abs{\itGamma_0} \leq \kappa$
and $\itGamma = \rfs{cl}_{\preceq}(\itGamma_0)$.
   \index{generates. $\itGamma_0$ generates $\itGamma$.
          This means
          $\itGamma = \rfs{cl}_{\preceq}(\itGamma_0)$}
   \index{generated. $\itGamma$ is $(\leq\kern-3pt \kappa)$-generated.
          This means $\exists \itGamma_0(\abs{\itGamma_0} \leq \kappa
          \mbox{ and } \itGamma = \rfs{cl}_{\preceq}(\itGamma_0))$}

(d) Let $\gamma \in \rfs{MC}$ and $a,b \in [0,\infty)$. 
Then $a \approx^{\gamma} b$ means that $a \leq \gamma(b)$ and
$b \leq \gamma(a)$.
   \index{N@AAAA@@$a \approx^{\gamma} b$. This means that
          $a \leq \gamma(b)$ and $b \leq \gamma(a)$}

(e) Let $X$ be a metric space, $x \in X$, $G \leq H(X)$ and
$\alpha \in \rfs{MBC}$.
We say that $G$ is {\it $\alpha$-infinitely-closed at $x$}
if there is $U \in \rfs{Nbr}(x)$ such that if
$F \subseteq G$ and $F$ satisfies:
\begin{itemize}
\addtolength{\parskip}{-11pt}
\addtolength{\itemsep}{06pt}
\item[(1)] for every $f \in F$, $f$ is $\alpha$-bicontinuous,
\item[(2)] for every $f \in F$, $\rfs{supp}(f) \subseteq U$
and $x \not\in \rfs{cl}(\rfs{supp}(f))$,
\item[(3)] for any distinct $f,g \in G$,
$\rfs{cl}(\rfs{supp}(f)) \cap \rfs{cl}(\rfs{supp}(g)) = \emptyset$,
\item[(4)] $\rfs{cl}(\bigcup_{f \in F} \rfs{supp}(f)) =
\sngltn{x} \cup \bigcup_{f \in F} \rfs{cl}(\rfs{supp}(f))$,
\vspace{-05.7pt}
\end{itemize}
\underline{then} \,$\bcirc\kern-1pt F \in G$.
\smallskip
   \index{infinitely-closed. $\alpha$-infinitely-closed at $x$''}

Note that if $F$ is as above, then $\bcirc\kern-1pt F \in H(X)$.
So $H(X)$ is $\alpha$-infinitely-closed at $x$ for every
$\alpha \in \rfs{MBC}$.

(f) When dealing with partial actions, we often wish to perform
a composition $g \scirc f$,
where $\rfs{Rng}(f) \not\subseteq \rfs{Dom}(g)$.
Such a composition is considered to be legal.
The domain of the resulting function is
$f\inverse(\rfs{Rng}(f) \cap \rfs{Dom}(g))$.

If $f,g$ are functions and $\rho$ is a $\onetoonen$ function,
then $f \sim^{\rho} g$ means that
$$
\rfs{Dom}(f) \cup \rfs{Rng}(f) \subseteq \rfs{Dom}(\rho), \ %
g = \rho \scirc f \scirc \rho\inverse.
$$
   \index{N@AAAA@@$f \sim^{\rho} g$. This means that
          $\rfs{Dom}(f) \cup \rfs{Rng}(f) \subseteq \rfs{Dom}(\rho)$
          and $g = \rho \scirc f \scirc \rho\inverse$}
\end{rm}
\end{defn}

\begin{prop}\label{metr-bldr-p3.10}
\num{a} If $X$ is a metric space, $A \subseteq X$ and
$x \in \rfs{acc}(A)$, then $\kappa(x,A) \geq \aleph_0$.

\num{b} If $E$ is a normed space, then
$\kappa(x,E) =
\min(\setm{\abs{D}}{D \mbox{ is a dense subset of } E})$
for every $x \in E$.

\num{c} If $E = \ell_{\infty}$,
then $\kappa(E) = 2^{\aleph_0}$.

\num{d} If $E$ is a Hilbert space with an orthonormal base of
cardinality $\nu$, then $\kappa(E) = \nu$.
\end{prop}

\noindent
{\bf Proof } The proof is trivial.
\smallskip\hfill\myqed

The next lemma says roughly that if for every $h \in H$,
$(h_{\lambda})^{\tau}$ is
$\itGamma$-bicontinuous at~$x$,
then there are $\gamma \in \itGamma$ and neighborhoods $T,V$
of $x$ and $e_H$ respectively
such that
$(h_{\lambda})^{\tau} \nrestriction T$ is $\gamma$-bicontinuous
for every $h \in V$.
This is proved under the assumption that $H$ is $G$-decayable,
where $G$ is an infinitely-closed subgroup of $H(X)$.

For countably generated $\itGamma$'s the conclusion of the lemma is
true for every metric space $X$. If however, $\itGamma$ is not
countably generated, then we need to assume that $\itGamma$ has
a generating set of size $\leq \kappa(X)$.
The lemma will be applied to $\bbT(E;F)$ and $\bbA(E;F)$.

\begin{lemma}\label{metr-bldr-l3.11}\label{bddly-lip-bldr-p5.4}
Suppose that:
\begin{itemize}
\item[\num{i}] $X$ is a metric space, $G \leq H(X)$,
$H$ is a topological group,
$\lambda$ is a partial action\break
of $H$ on $X$,
$x \in \rfs{Fld}(\lambda)$,
$x$ is a $\lambda$-limit-point,
$\alpha \in \rfs{MBC}$,
$G$ is $\alpha$-infinitely-closed at $x$,\break
and for some
$N \in \rfs{Nbr}(x)$,
$\lambda$ is $(\alpha,G)$-decayable at every point
$y \in H_{\lambda}(x) \cap N$.
Set
$$
\kappa = \min(\setm{\kappa(x,V_{\lambda}(x))}{V \in \rfs{Nbr}(e_H)}).
$$
\item[\num{ii}] $Y$ is metric space and $\iso{\tau}{X}{Y}$.
\item[\num{iii}] $\itGamma$ is a modulus of continuity,
and $\itGamma$ is $(\leq\kern-3pt\kappa)$-generated.
\item[\num{iv}] There is $U \in \rfs{Nbr}(x)$ such that
for every $g \in G\sprt{\kern1ptU}$:
if $g$ is $\alpha \scirc \alpha$-bicon\-tinuous,
then $g^{\tau}$ is $\itGamma$-bicontinuous at $\tau(x)$.
\end{itemize}
\underline{Then} $P(x)$ holds, where
\smallskip
\newline
\rule{5pt}{0pt}$P(x)$: \,%
There are $T \in \rfs{Nbr}(x)$, $V \in \rfs{Nbr}(e_H)$
and $\gamma \in \itGamma$ such that
for every $h \in V$,\break
\phantom{\rule{5pt}{0pt}$P(x)$: \,}
$T \subseteq \rfs{Dom}(h_{\lambda})$
and $(h_{\lambda})^{\tau} \nrestriction \tau(T)$
is $\gamma$-bicontinuous.
\end{lemma}

\noindent
{\bf Proof }
Let $U_1 \in \rfs{Nbr}(x)$ be as assured by the 
$\alpha$-infinite-closedness of $G$ at $x$.
Let $r_x$ be as assured by the decayability of $H$ at $x$.
Let $r \in (0,r_x)$ be such that $B(x,r) \subseteq U_1 \cap U \cap N$,
and $W = V_{x,r}$ be as assured by the decayability of $H$ at $x$.
So $W \in \rfs{Nbr}(e_H)$,
$W \times B(x,r) \subseteq \rfs{Dom}(\lambda)$
and $W_{\lambda}(x) \subseteq B(x,r)$.
First we prove the following claim.

{\bf Claim 1.} There is
$y \in B(x,\dgfrac{r}{2}) \cap W_{\lambda}(x)$ such that $P(y)$ holds.

{\bf Proof } Suppose by contradiction that there is no such $y$.
Let $\itGamma_0$ be as assured by Clause~(iii).
We distinguish between two cases.

{\bf Case 1 } $\abs{\itGamma_0} = \aleph_0$.
Let $\vecx = \setm{x_i}{i \in \bbN}$ be a $\onetoonen$ sequence
tending to $x$ and contained in
$B(x,\dgfrac{r}{2}) \cap W_{\lambda}(x) - \sngltn{x}$.
Let $\setm{\gamma_i}{i \in \bbN}$ be an enumeration of $\itGamma_0$
such that $\setm{j}{\gamma_j = \gamma_i}$ is infinite
for every~$i$.
Let $r_{x_i} > 0$ be as assured by the decayability of $\lambda$
at $x_i$.
Let $\setm{r_i}{i \in \bbN}$ be a sequence
such that for any distinct $i,j \in \bbN$
we have $0 < r_i < r_{x_i}$,
$B(x_i,r_i) \subseteq B(x,r)$, $d(x_i,x) > r_i$
and $\rfs{cl}(B(x_i,r_i)) \cap \rfs{cl}(B(x_j,r_j)) = \emptyset$.

Let $W_i = V_{x_i,r_i}$ be as assured by the
decayability of $\lambda$ at $x_i$.
That is, $W_i \in \rfs{Nbr}(e_H)$ and
$\rfs{Dom}(h_{\lambda}) \supseteq B(x_i,\dgfrac{r_i}{2})$
for every $h \in W_i$,
and there is $g \in G$ such that $g$ is
$\alpha$-bicontinuous,
$g \nrestriction B(x_i,\dgfrac{r_i}{2}) =
h_{\lambda} \nrestriction B(x_i,\dgfrac{r_i}{2})$
and $\rfs{supp}(g) \subseteq B(x_i,r_i)$.

Let $V_i = B(x_i,\dgfrac{r_i}{2})$.
Then $\rfs{Dom}(h_{\lambda}) \supseteq V_i$ for every $h \in W_i$.
Since $\neg P(x_i)$ holds,\break
there is $h_i \in W_i$ such that
$((h_i)_{\lambda})^{\tau} \nrestriction \tau(V_i)$
is not $\gamma_i$-bicontinuous.
Let $g_i \in G$ be such that $g_i$ is $\alpha$-bicontinuous,
$g_i \nrestriction B(x_i,\dgfrac{r_i}{2}) =
(h_i)_{\lambda} \nrestriction B(x_i,\dgfrac{r_i}{2})$
and $\rfs{supp}(g_i) \subseteq B(x_i,r_i)$.
Clearly, $F \eqdf \setm{g_i}{i \in \bbN}$ satisfies Clauses (1)\,-\,(4)
in the definition of $\alpha$-infinite-closedness,
so $g \eqdf \bcirc_{i \in \sboldbbN}\, g_i \in G$.
For every $u,v \in X$ there are $i,j \in \bbN$ such that
$g(u) = g_i \scirc g_j(u)$ and $g(v) = g_i \scirc g_j(v)$.
So $g$ is $\alpha \scirc \alpha$-continuous. Similarly,
$g\inverse$ is $\alpha \scirc \alpha$-continuous.
Since $\rfs{supp}(g) \subseteq U$, by Clause~(iv),
$g^{\tau}$ is $\itGamma$-bicontinuous at $\tau(x)$.
That is, there are $\gamma \in \itGamma$ and $T \in \rfs{Nbr}(\tau(x))$
such that
\begin{equation}
\tag{1.1}
g^{\tau} \nrestriction T \mbox{ is } \gamma\mbox{-bicontinuous.}
\end{equation}
Let $i$ be such that $\gamma \preceq \gamma_i$,
and let $t > 0$ be such that
$\gamma \nrestriction [0,t] \leq \gamma_i \nrestriction [0,t]$.
There is $j$ such that
$\gamma_j = \gamma_i$,
$\tau(B(x_j,r_j)) \subseteq T$
and ($\dagger$) $\rfs{diam}(\tau(B(x_j,r_j))) < t$.

Set $k = (h_j)_{\lambda}$.
Now,
$g \nrestriction V_j = g_j \nrestriction V_j = k \nrestriction V_j$.
So
\begin{equation}
\tag{1.2}
g^{\tau} \nrestriction \tau(V_j) =
(g_j)^{\tau} \nrestriction \tau(V_j) =
k^{\tau} \nrestriction \tau(V_j).
\end{equation}
Recall that $k^{\tau} \nrestriction \tau(V_j)$
is not $\gamma_j$-bicontinuous.
So there are $u,v \in \tau(V_j)$ such that
\newline
$d^Y(k^{\tau}(u),k^{\tau}(v)) \not\approx^{\gamma_j} d^Y(u,v)$.
By (1.2),
\begin{equation}
\tag{1.3}
d^Y(g^{\tau}(u),g^{\tau}(v)) \not\approx^{\gamma_j} d^Y(u,v).
\end{equation}
Let $u_1 = \tau\inverse(u)$ and $v_1 = \tau\inverse(v)$.
So $u_1,v_1 \in B(x_j,\dgfrac{r_j}{2})$.
Since
$k \nrestriction B(x_j,\dgfrac{r_j}{2}) =
g_j \nrestriction B(x_j,\dgfrac{r_j}{2})$
and $\rfs{supp}(g_j) \subseteq B(x_j,r_j)$,
we have $k(u_1),k(v_1) \in B(x_j,r_j)$.
By ($\dagger$), \hbox{$d^Y(\tau(k(u_1)),\tau(k(v_1))) < t$.}
Also, $\tau(k(u_1)) = k^{\tau}(u)$, and the same holds for
$v$ and $v_1$.
So $d^Y(u,v) < t$ and $d^Y(k^{\tau}(u),k^{\tau}(v)) < t$.
By (1.2),
\begin{equation}
\tag{1.4}
d^Y(u,v) < t, \mbox{ \ } d^Y(g^{\tau}(u),g^{\tau}(v)) < t.
\end{equation}
Recall that
$\gamma \nrestriction [0,t] \leq \gamma_j \nrestriction [0,t]$.
Hence by (1.3) and (1.4),
\begin{equation}
\tag{1.5}
d^Y(g^{\tau}(u),g^{\tau}(v)) \not\approx^{\gamma} d^Y(u,v).
\end{equation}
Recall that $u,v \in \tau(V_j) \subseteq T$.
Hence (1.1) and (1.5) are contradictory.
So there is
$y \in B(x,\dgfrac{r}{2}) \cap W_{\lambda}(x)$ such that $P(y)$ holds.
\vspace{1.5mm}

{\bf Case 2 } $\abs{\itGamma_0} > \aleph_0$.
Let $L = W_{\lambda}(x)$ and $\kappa = \kappa^X(x,L)$.
We prove that there are sequences
$\setm{r_i}{i \in \bbN} \subseteq (0,\infty)$ and
$\setm{L_i}{i \in \bbN}$ such that:
\begin{itemize}
\item[(i)] $r_0 = \dgfrac{r}{2}$
and $\setm{r_i}{i \in \bbN}$ is a strictly decreasing sequence
converging to $0$;
\item[(ii)] for every $i \in \bbN$,
$L_i \subseteq L \cap (B(x,r_i) - B(x,r_{i + 1}))$ and $L_i$ is spaced;
\item[(iii)] $\abs{\bigcup \setm{L_i}{i \in \bbN}} = \kappa$.
\end{itemize}

Suppose first that $\rfs{cf}(\kappa) = \aleph_0$.
(That is, there is a countable set of cardinals $\Kappa$ such that
for every $\kappa' \in \Kappa$, $\kappa' < \kappa$ and
$\sum \Kappa = \kappa$).
Let $\Kappa = \setm{\kappa_i}{i \in \bbN}$ and $r_0 = \dgfrac{r}{2}$.
We may assume that each $\kappa_i$ is infinite.
We define $L_i$ and $r_{i + 1}$ by induction on $i$.
Suppose that $r_i$ has been defined.
Since $\kappa_i < \kappa^X(x,L)$
there is $L_i \subseteq L \cap B(x,r_i)$ such that $L_i$ is spaced
and $\abs{L_i} = \kappa_i$.
Suppose that $L_i$ is $s_i$-spaced.
There is at most one member $y \in L_i$ such that
$d(x,y) < \dgfrac{s_i}{2}$. So by removing this member we may assume
that $d(L_i,x) \geq \dgfrac{s_i}{2}$.
Let $r_{i + 1} = \min(\frac{s_i}{2},\frac{1}{i + 1})$.
Evidently, $\setm{r_i}{i \in \bbN}$, $\setm{L_i}{i \in \bbN}$
fulfill (i)\,-\,(iii).
\vspace{1.5mm}

Suppose that $\rfs{cf}(\kappa) > \aleph_0$.
First we show that
\begin{itemize}
\item[$(*)$] For every $s > 0$ there is
$M \subseteq L \cap B(x,s)$ such that $\abs{M} = \kappa$
and $M$ is spaced.
\end{itemize}
Suppose not, and let $s$ be a counter-example.
For every $n > 0$ let $\Kappa_n$ be the set of all $\kappa'$ such that
there is $M \subseteq L \cap B(x,s)$ such that $\abs{M} = \kappa'$
and $M$ is $\dgfrac{1}{n}$-spaced.
Then there is $n$ such that $\Kappa_n$ is unbounded in~$\kappa$. 
Let $N$ be a maximal $\frac{1}{2n}$-spaced subset of
$L \cap B(x,s)$.
Then $\abs{N} < \kappa$. So there is $\kappa' \in \Kappa_n$
such that $\abs{N} < \kappa'$.
Let $M$ be a $\dgfrac{1}{n}$-spaced subset of $L \cap B(x,s)$
of cardinality $\kappa'$.
Then there are $y \in N$ and $z_1,z_2 \in M$ such that
$z_1,z_2 \in B(y,\frac{1}{2n})$. A contradiction, so $(*)$ holds.

As in the case that $\rfs{cf}(\kappa) = \aleph_0$ we define
a sequence $\setm{\kappa_i}{i \in \bbN}$.
Indeed, we set $\kappa_i = \kappa$ for every $i \in \bbN$.
The $L_i$'s and $r_i$'s are now constructed as in the case
$\rfs{cf}(\kappa) = \aleph_0$,
and they obviously fulfill Clauses (i)\,-\,(iii).

We really need sequences
$\setm{r_i}{i \in \bbN} \subseteq (0,\infty)$ and
$\setm{L_i}{i \in \bbN}$ which fulfill the following conditions:
\begin{itemize}
\item[(i)] $r_0 = \dgfrac{r}{2}$
and $\setm{r_i}{i \in \bbN}$ is a strictly decreasing sequence
converging to $0$;
\item[(ii)] for every $i \in \bbN$,
$L_i \subseteq L \cap (B(x,\dgfrac{r_i}{2}) - B(x,2 r_{i + 1}))$
and $L_i$ is spaced,
and $\abs{L_i} \leq \abs{L_j}$
for every $i < j$;
\item[(iii)] $\abs{\bigcup \setm{L_i}{i \in \bbN}} = \abs{\itGamma_0}$.
\end{itemize}
Such sequences can be obtained from the original
$\setm{r_i}{i \in \bbN}$ and $\setm{L_i}{i \in \bbN}$
by taking an appropriate subsequence of $\setm{r_i}{i \in \bbN}$
and by replacing $L_i$ by a subset of $L_i$ if necessary.

Let $s_i > 0$ be such that $L_i$ is $s_i$-spaced.
Set $M = \bigcup \setm{L_i}{i \in \bbN}$, and let
$\fnn{\iota}{M}{\itGamma_0}$ be a function such that for every
$\gamma \in \itGamma_0$ there is $n \in \bbN$ such that
$\gamma \in \iota(L_m)$ for every $m \geq n$.
Define $\gamma_y = \iota(y)$.
Let $r_y$ be as assured by the decayability of $H$ at $y$.
For every $y \in M$ we define $s_y > 0$.
If $y \in L_i$,
choose $s_y < \min(r_y,r_{i + 1},\dgfrac{s_i}{3})$.
Note that for distinct
$y,z \in L_i$, $B(y,s_y) \subseteq B(x,r_i) -  B(x,r_{i + 1})$ and
$\rfs{cl}(B(y,s_y)) \cap \rfs{cl}(B(z,s_z)) = \emptyset$.
So for distinct
$y,z \in M$, $\rfs{cl}(B(y,s_y)) \cap \rfs{cl}(B(z,s_z)) = \emptyset$.

For every $y \in M$ let $W_y = V_{y,s_y}$ be as assured by the
decayability of $\lambda$ at~$y$.
That is, $W_y \in \rfs{Nbr}(e_H)$,
$\rfs{Dom}(h_{\lambda}) \supseteq B(y,\dgfrac{s_y}{2})$
for every $h \in W_y$,
and there is $g \in G$ such that $g$ is
$\alpha$-bicontinuous,
$g \nrestriction B(y,\dgfrac{s_y}{2}) =
h_{\lambda} \nrestriction B(y,\dgfrac{s_y}{2})$
and $\rfs{supp}(g) \subseteq B(y,s_y)$.

Let $V_y = B(y,\dgfrac{s_y}{2})$.
So $\rfs{Dom}(h_{\lambda}) \supseteq V_y$ for every $h \in W_y$.
Since $\neg P(y)$ holds,
there is $h_y \in W_y$ such that
$((h_y)_{\lambda})^{\tau} \nrestriction \tau(V_y)$
is not $\gamma_y$-bicontinuous.
Let $g_y \in G$ be such that $g_y$ is $\alpha$-bicontinuous,
$g_y \nrestriction B(y,\dgfrac{s_y}{2}) =
(h_y)_{\lambda} \nrestriction B(y,\dgfrac{s_y}{2})$
and $\rfs{supp}(g_y) \subseteq B(y,s_y)$.
For any distinct $y,z \in M$,
$\rfs{supp}(g_y) \cap \rfs{supp}(g_z) =~\emptyset$.
Clearly, $F \eqdf \setm{g_y}{y \in M}$ satisfies Clauses~(1)\,-\,(4)
in the definition of $\alpha$-infinite-closedness,
so $g = \bcirc_{y \in M}\, g_y \in G$.
The rest of the argument is identical to the one given in Case 1.
We have proved Claim 1.
\smallskip

Let $y$ be as assured by Claim 1.
Since $y \in W_{\lambda}(x)$, there is $\hath \in W$ such that
$y = \hath_{\lambda}(x)$.
Since $W = V_{x,r}$,
there is $g \in G$ such that $g$ is $\alpha$-bicontinuous,
$g \nrestriction B(x,\dgfrac{r}{2}) =
\hath_{\lambda} \nrestriction B(x,\dgfrac{r}{2})$
and $\rfs{supp}(g) \subseteq B(x,r)$. So $g(x) = y$.
Since $\alpha \in \rfs{MBC}$, we have
$\alpha \leq \alpha \scirc \alpha$,
and hence $g$ is $\alpha \scirc \alpha$-bicontinuous.
The bicontinuity of $g$
and the fact $\rfs{supp}(g) \subseteq B(x,r) \subseteq U$
imply that $g^{\tau}$ is $\itGamma$-bicontinuous at $\tau(x)$.
Let $R \in \rfs{Nbr}(\tau(x))$ and $\beta \in \itGamma$ be such that
$g^{\tau} \nrestriction R$ is $\beta$-bicontinuous.
We may assume that
\begin{equation}
\tag{2.1}
\tau\inverse(R) \subseteq B(x,\dgfrac{r}{2}).
\end{equation}
Hence
$g^{\tau} \nrestriction R = (\hath_{\lambda})^{\tau} \nrestriction R$.
So
\begin{equation}
\tag{2.2}
(\hath_{\lambda})^{\tau} \nrestriction R \mbox{ is }
\beta\mbox{-bicontinuous.}
\end{equation}
Note that if $\,T',V' ,\gamma'$ fulfill the requirements
of $P(y)$ and $T' \supseteq T'' \in \rfs{Nbr}(y)$,
then $T'',V',\gamma'$ fulfill the requirements of $P(y)$.
Since $P(y)$ holds, there are $S_1 \in \rfs{Nbr}(y)$,
$V_1 \in \rfs{Nbr}(e_H)$
and $\gamma_1 \in \itGamma$ such that for every
$h \in V_1$,
\begin{equation}
\tag{2.3}
S_1 \subseteq \rfs{Dom}(h_{\lambda}), \ \ %
(h_{\lambda})^{\tau} \nrestriction \tau(S_1)
\mbox{ is } \gamma_1\mbox{-bicontinuous.}
\end{equation}
Since $\hath_{\lambda}(x) = y$ and $\tau\inverse(R) \in \rfs{Nbr}(x)$,
we may assume that
\begin{equation}
\tag{2.4}
S_1 \subseteq \hath_{\lambda}(\tau\inverse(R)).
\end{equation}
So $S_1 \subseteq \hath_{\lambda}(B(x,\dgfrac{r}{2}))$.
Let $S_2 \in \rfs{Nbr}(y)$ and
$V_2 \in \rfs{Nbr}(e_H)$ be such that
\begin{equation}
\tag{2.5}
S_2 \subseteq S_1, \ \ V_2 \subseteq V_1, \ \ %
\lambda(V_2 \times S_2) \subseteq S_1.
\end{equation}
Note that $S_2 \subseteq \rfs{Rng}(\hath_{\lambda})$.
We define
$T = (\hath_{\lambda})\inverse(S_2)$,
$V = \hath\inverse \cdot V_2 \cdot \hath$
and $\gamma = \beta \scirc \gamma_1 \scirc \beta$
and show that $T,V,\gamma$ satisfy the requirements of $P(x)$.
Since $\beta,\gamma_1 \in \itGamma$, we have
\begin{equation}
\tag{2.6}
\gamma \in \itGamma.
\end{equation}
We verify that if $h \in V$, then
\begin{equation}
\tag{2.7}
T \subseteq \rfs{Dom}(h_{\lambda}) \mbox{ \ and\ \ }
(h^{\hath})_{\lambda} \nrestriction S_2
\ \sim^{\rho\inverse}\ %
h_{\lambda} \nrestriction T,
\mbox{ where } \rho = \hath_{\lambda} \nrestriction \tau\inverse(R).
\end{equation}
Let $\barh = h^{\hath}$. Then $\barh \in V_2$
and $h = \hath\inverse \cdot \barh \cdot \hath$.
We show that
$\hath_{\lambda}(z)$, $\barh_{\lambda}(\hath_{\lambda}(z))$
and $(\hath\inverse)_{\lambda}(\barh_{\lambda}(\hath_{\lambda}(z)))$
are defined for every $z \in T$.
Clearly, $T \subseteq \rfs{Dom}(\hath_{\lambda})$
and $\hath_{\lambda}(T) = S_2$.
So by (2.5),
\begin{itemize}
\item[(i)] for every $z \in T$,
$\barh_{\lambda}(\hath_{\lambda}(z))$ is defined and
$\barh_{\lambda}(\hath_{\lambda}(z)) \in S_1$.
\end{itemize}
By (2.4), $S_1  \subseteq \rfs{Rng}(\hath_{\lambda})$.
So
$(\hath_{\lambda})\inverse(\barh_{\lambda}(\hath_{\lambda}(z)))$
is defined. Since $h = \barh^{\hath}$
and by the definition of a partial action,
it follows that
\begin{itemize}
\item[(ii)] for every $z \in T$, $h_{\lambda}(z)$ is defined and
$h_{\lambda}(z) =
(\hath_{\lambda})\inverse \scirc \barh_{\lambda} \scirc
\hath_{\lambda}(z)$.
\end{itemize}
By (ii), $T \subseteq \rfs{Dom}(h_{\lambda})$,
and by (2.1), $\tau\inverse(R) \subseteq \rfs{Dom}(\hath_{\lambda})$. So
$\rfs{Dom}(\rho\inverse) = \rfs{Rng}(\rho) =
\hath_{\lambda}(\tau\inverse(R))$.
Since $\barh \in V_2$,
we have that $S_2 \subseteq \rfs{Dom}(\barh_{\lambda})$,
hence $\rfs{Dom}(\barh_{\lambda} \nrestriction S_2) = S_2$.
By (2.4) and (2.5), $S_2 \subseteq \hath_{\lambda}(\tau\inverse(R))$.
So
$\rfs{Dom}(\barh_{\lambda} \nrestriction S_2) \subseteq
\rfs{Dom}(\rho\inverse)$.
We have that
$\rfs{Rng}(\barh_{\lambda} \nrestriction S_2) = \barh_{\lambda}(S_2)$,
and from (2.5) and the fact that $\barh \in V_2$,
it follows that $\barh_{\lambda}(S_2) \subseteq S_1$.
By (2.4), $S_1 \subseteq \hath_{\lambda}(\tau\inverse(R))$,
so
$\rfs{Rng}(\barh_{\lambda} \nrestriction S_2) \subseteq
\rfs{Dom}(\rho\inverse)$.
Note that $T \subseteq \tau\inverse(R)$;
indeed, this follows from the definition of $T$,
(2.4) and (2.5). So
\begin{itemize}
\item[(iii)] for every $z \in T$,
$\hath_{\lambda}(z) =
(\hath_{\lambda} \nrestriction \tau\inverse(R))(z) = \rho(z)$.
\end{itemize}
Also,
\begin{itemize}
\item[(iv)] for every $z \in T$,
$\barh_{\lambda}(\hath_{\lambda}(z)) =
(\barh_{\lambda} \nrestriction S_2)(\hath_{\lambda}(z))$.
\end{itemize}
Let $z \in T$ and denote $u = \barh_{\lambda}(\hath_{\lambda}(z))$.
By (i) and (2.4),
$u \in S_1 \subseteq \hath_{\lambda}(\tau\inverse(R))  =
\rfs{Dom}(\rho\inverse)$.
Hence $(\hath_{\lambda})\inverse(u) = \rho\inverse(u)$.
We conclude that
\begin{itemize}
\item[(v)] for every $z \in T$,
$(\hath_{\lambda})\inverse(\barh_{\lambda}(\hath_{\lambda}(z))) =
\rho\inverse(\barh_{\lambda}(\hath_{\lambda}(z)))$.
\end{itemize}
It follows from (ii)\,-\,(v) that
$h_{\lambda} \nrestriction T =
\rho\inverse \scirc (\barh_{\lambda} \nrestriction S_2) \scirc \rho$.
We have verified (2.7).
Next conjugate (2.7) by $\tau$. We obtain that
\begin{equation}
\tag{2.8}
((h^{\hath})_{\lambda} \nrestriction S_2){^\tau}
\ \sim^{(\rho\inverse)^{\tau}}\ %
(h_{\lambda} \nrestriction T)^{\tau}.
\end{equation}
Clearly,
$((h^{\hath})_{\lambda} \nrestriction S_2)^{\tau} =
((h^{\hath})_{\lambda})^{\tau} \nrestriction \tau(S_2)$.
Since $h \in V$, we have $h^{\hath} \in V^{\hath} = V_2$.
So by (2.3),
\begin{equation}
\tag{2.9}
((h^{\hath})_{\lambda} \nrestriction S_2){^\tau}
\mbox{ is } \gamma_1\mbox{-bicontinuous.}
\end{equation}
Fact (2.8) has the form
$f \sim^{\sigma\inverse} k$,
where $f = ((h^{\hath})_{\lambda} \nrestriction S_2){^\tau}$,
$k = (h_{\lambda} \nrestriction T)^{\tau}$
and
$\sigma  = \rho^{\tau} = ({\hath_{\lambda}})^{\tau} \nrestriction R$.
By (2.9), $f$ is $\gamma_1$-bicontinuous,
and by (2.2) $\sigma$ is $\beta$-bicontinuous.
Since $k = \sigma\inverse \scirc f \scirc \sigma$,
it follows that $k$ is
$\beta \scirc \gamma_1 \scirc \beta$-bicontinuous.
Recall that $\gamma = \beta \scirc \gamma_1 \scirc \beta$
and
$k = (h_{\lambda} \nrestriction T)^{\tau} =
(h_{\lambda})^{\tau} \nrestriction \tau(T)$.
Hence $(h_{\lambda})^{\tau} \nrestriction \tau(T)$
is $\gamma$-bicontinuous.

We have shown that for every $h \in V$,
$\rfs{Dom}(h_{\lambda}) \supseteq T$ and
$(h_{\lambda})^{\tau} \nrestriction \tau(T)$ is $\gamma$-bicontinuous.
So $T,V,\gamma$ satisfy the requirements of the lemma.
\hfill\myqed

\subsection{Translation-like partial actions.}
\label{ss3.3}
\label{ss3.3-translation-like-actions}

We have isolated the properties of \kern0.6mm$\bbT(E)$ and $\bbA(E)$
which are used in the proof that $\tau$ is $\itGamma$-continuous.
The following definition deals with the properties of $\bbT(E)$.
Partial actions having these properties are called
translation-like partial actions.
In fact, the definition captures the properties of
\kern0.6mm$\bbT(E;F)$,
where $F$ is any dense linear subspace of $E$.
The properties of $\bbA(E)$ to be used, appear in
Definition \ref{metr-bldr-d3.28}(b).

\begin{defn}\label{metr-bldr-d3.12}
\begin{rm}
Suppose that $X$ is a metric space, $H$ is a topological group,
and $\lambda$ is a partial action of $H$ on $X$.
Let $x \in \rfs{Fld}(\lambda)$. We say that $\lambda$ is
{\it a translation-like partial action at $x$},
if for every $V \in \rfs{Nbr}(e_H)$ there are:
\vspace{3pt}
\newline\rule{01mm}{0pt}
\renewcommand{\arraystretch}{1.2}
\begin{tabular}{rl}
(i)& $U = U_{x,V} \in \rfs{Nbr}(x)$,
and a dense subset of $U$, $D = D_{x,V}$;\\
(ii)& a radius $r = r_{x,V} > 0$ and a constant $K = K_{x,V} > 0$;
\vspace{5pt}\\
\end{tabular}
\renewcommand{\arraystretch}{1.0}
\newline
such that the following holds.\newline
For any distinct $\barx_0,\barx_1 \in D$ there are
$n \leq K \ncdot \frac{r}{d(\barx_0,\barx_1)}$\,,
a sequence $\barx_0 = x_0, x_1, \ldots, x_n \in X$ and
$h_1, \ldots, h_n \in V$ such that
$x_n \not \in B(x,r)$, and for every $i = 1, \ldots, n$,\ \ %
$\barx_0,\barx_1 \in \rfs{Dom}((h_i)_{\lambda})$,\ \,%
$(h_i)_{\lambda}(\barx_0) = x_{i - 1}$
and $(h_i)_{\lambda}(\barx_1) = x_i$.

A partial action $\lambda$ is {\it translation-like},
if for every $x \in \rfs{Fld}(\lambda)$,
$\lambda$ is translation-like at $x$.
   \index{translation-like partial action at $x$}
   \index{translation-like partial action}
\end{rm}
\end{defn}

\begin{prop}\label{metr-bldr-p3.13}\label{was-p3.2}
Let $E$ be a normed space, $F$ be a dense linear subspace of $E$ and
$X \subseteq E$ be open.
Then $\lambda_{\bbT}^{E;F} \drest X$ is a translation-like
partial action.
\end{prop}
 
\noindent
{\bf Proof }
For $x \in X$
and $V \in \rfs{Nbr}^{\bbT(E;F)}(\rfs{Id})$ we define $U = U_{x,V}$,
$D = D_{x,V}$ etc.\ as follows.
Let $r_0 > 0$ be such that $B^E(x,r_0) \subseteq X$
and $\setm{\rfs{tr}_v}{v \in B^F(0,r_0)} \subseteq V$.
Now define
$U = B(x,\dgfrac{r_0}{4})$, $D = F \cap U$, $r = \dgfrac{r_0}{2}$
and $K = 2$.

For distinct $\barx_0, \barx_1 \in D$ we define $n$,
$x_0, \ldots, x_n$ and $h_1, \ldots, h_n$
as required in Definition~\ref{metr-bldr-d3.12}.
Let $n$ be the least integer such that
$n \ncdot \norm{\barx_1 - \barx_0} \geq r$.
For $i = 0, \ldots, n$ let $x_i = \barx_0 + i(\barx_1 - \barx_0)$
and for $i = 1, \ldots, n$ let
$h_i = \rfs{tr}_{(i - 1)(\barx_1 - \barx_0)}$.
It is easily checked that $n$, the $x_i$'s and the $h_i$'s are
as required.\medskip\hfill\myqed

We let $X$ and $Y$ denote metric spaces. Their metrics are denoted
by $d^X$ and $d^Y$. However, in most cases we write $d(x,y)$
as an abreviation of both $d^X(x,y)$ and $d^Y(x,y)$.

\begin{lemma}\label{metr-bldr-l3.14}
Let $X$ be a metric space and $\lambda$ be a partial action of
$H$ on~$X$. Suppose that $x \in \rfs{Fld}(\lambda)$ and $\lambda$ is
traslation-like at $x$.
Let $Y$ be a metric space and $\iso{\tau}{X}{Y}$.
Let $\itGamma \subseteq \rfs{MC}$, and suppose that for every
$\gamma \in \itGamma$ and $K > 0$, $K \mcdot \gamma \in \itGamma$.
Suppose that
$P(x)$ of Lemma \ref{metr-bldr-l3.11} holds.
That is,
there are $T \in \rfs{Nbr}(x)$, $V \in \rfs{Nbr}(e_H)$
and $\gamma \in \itGamma$ such that
for every $h \in V$,
$T \subseteq \rfs{Dom}(h_{\lambda})$ and
$(h_{\lambda})^{\tau} \nrestriction \tau(T)$
is $\gamma$-bicontinuous.
\underline{Then} $\tau\inverse$ is $\itGamma$-continuous at $\tau(x)$.
\end{lemma}

\noindent
{\bf Proof }
Let $U = U_{x,V}$, $D = D_{x,V}$, $r = r_{x,V}$ and $K = K_{x,V}$
{\thickmuskip=2mu \medmuskip=1mu \thinmuskip=1mu 
be as assured by the translation-likeness of $H$ at $x$.
Set $y = \tau(x)$, $B = B(x,r)$ and $C = \tau(B)$.
}
Since $C \in \rfs{Nbr}(y)$,
we have $e \eqdf d(y,Y - C) > 0$.
Let $R = \tau(T \cap U) \cap B(y,\dgfrac{e}{2})$
and $M = \dgfrac{2Kr}{e}$.
Since $\gamma \in \itGamma$, we have $M \ncdot \gamma \in \itGamma$.

We show that $\tau\inverse \nrestriction R$ is
$M \ncdot \gamma$-continuous. Suppose by way of contradiction
that this is not true.
\hbox{Hence there are  $\bary_0,\bary_1 \in R$} such that
$d(\tau\inverse(\bary_0),\tau\inverse(\bary_1)) >
M \ncdot \gamma(d(\bary_0,\bary_1))$.
Since $D$ is dense in $U$ and $\bary_0,\bary_1 \in \tau(U)$,
we may assume that $\bary_0,\bary_1 \in \tau(D)$.
For every $h \in H$ let $\hath$ denote $h_{\lambda}$,
and for $\ell = 0,1$ let
\hbox{$\barx_{\ell} = \tau\inverse(\bary_{\ell})$.}
Hence \hbox{$\barx_0,\barx_1 \in D$.}
So there are $n \leq \dgfrac{Kr}{d(\barx_0,\barx_1)}$,
$\barx_0 = x_0, x_1, \dots, x_n$
and
$h_1, \ldots, h_n \in V$ such that $x_n \not\in B$, and for every
$i = 1, \ldots, n$,
$\barx_0,\barx_1 \in \rfs{Dom}(h_i)$,
$\hath_i(\barx_0) = x_{i - 1}$
and $\hath_i(\barx_1) = x_i$.
For $i = 1, \ldots, n$
let $y_i = \tau(x_i)$.

In the space $Y$ we thus have the following situation:
\vspace{3pt}
\newline\rule{01mm}{0pt}
\renewcommand{\arraystretch}{1.2}
\begin{tabular}{rl}
(i)& $d(y,y_0) < \dgfrac{e}{2}$;\\
(ii)& for every $i = 1, \ldots, n$,\,
$\hath_i^{\tau}(\bary_0) = y_{i - 1}$
and $\hath_i^{\tau}(\bary_1) = y_i$;\\
(iii)& $y_n \not\in C$.
\vspace{5pt}\\
\end{tabular}
\renewcommand{\arraystretch}{1.0}
\newline
Every $h_i$ belongs to $V$,
hence $\hath_i^{\tau} \restriction \tau(T)$ is
$\gamma$-bicontinuous.
Also, $\bary_0,\bary_1 \in \tau(T)$,~so
\vspace{3pt}
\newline\rule{01mm}{0pt}
\renewcommand{\arraystretch}{1.2}
\begin{tabular}{rl}
(iv)& $d(y_{i - 1},y_i) \leq \gamma(d(\bary_0,\bary_1))$.
\vspace{5pt}\\
\end{tabular}
\renewcommand{\arraystretch}{1.0}
\newline
Hence
\vspace{1.5mm}
\newline
\rule{7pt}{0pt}
\renewcommand{\arraystretch}{1.5}
\addtolength{\arraycolsep}{-3.3pt}
$
\begin{array}{lll}
e
&
=
&
d(y,Y - C) \leq d(y,y_n) \leq d(y,y_0) +
\sum^n_{i = 1}d(y_{i - 1},y_i) <
\dgfrac{e}{2} + n \cdot \gamma(d(\bary_0,\bary_1))
\\
&
\leq
&
\dgfrac{e}{2} + \frac{Kr}{d(\barx_0,\barx_1)} \cdot
\gamma(d(\bary_0,\bary_1))
<
\dgfrac{e}{2} +
\frac{Kr}{M \cdot \gamma(d(\bary_0,\bary_1))} \cdot
\gamma(d(\bary_0,\bary_1))
=
\dgfrac{e}{2} + \frac{Kr}{\dgfrac{2Kr}{e}} = e.
\vspace{1.7mm}
\end{array}
$
\renewcommand{\arraystretch}{1.0}
\addtolength{\arraycolsep}{3.3pt}
\newline
\medskip\smallskip
A contradiction, so the lemma is proved.
\hfill\myqed

The following theorem is the conjunction of 
Lemmas \ref{metr-bldr-l3.11} and \ref{metr-bldr-l3.14}.
It will be used in Theorem \ref{metr-bldr-t3.16}.
The statement of Theorem \ref{metr-bldr-t3.15} is rather technical.
So it seems worthwhile to explain its main application.
Let $X$ be an open subset of a normed space $E$ and
$G \leq H(X)$.
Suppose that for every $x \in X$ and $r > 0$
there are $s \in (0,r)$ and $K > 0$ such that
for every $v \in B_E(0,s)$ there is $g \in G$ such that
$g \nrestriction B(x,s) = \rfs{tr}_v \nrestriction B(x,s)$,
$g$~is $K$-bilipschitz and $\rfs{supp}(g) \subseteq B(x,r)$.
Assume further that $G$ is $\alpha$-infinitely closed for every
$\alpha$ of the form $y = Mt$.
Then if $\tau$ is a homeomorphism between $X$ and a metric space $Y$,
$\itGamma$ is a countably generated modulus of continuity and
$G^{\tau} \subseteq \rfs{LIP}^{\srfs{LC}}_{\itGamma}(Y)$,
then $\tau\inverse$ is locally $\itGamma$-continuous.

\begin{theorem}\label{metr-bldr-t3.15}
Suppose that:
\begin{list}{}
{\setlength{\leftmargin}{39pt}
\setlength{\labelsep}{06pt}
\setlength{\labelwidth}{30pt}
\setlength{\itemindent}{0pt}
\addtolength{\topsep}{-04pt}
\addtolength{\parskip}{-02pt}
\addtolength{\itemsep}{-05pt}
}
\item[\num{i}] $X$ is a metric space, $G \leq H(X)$,
$H$ is a topological group, $\lambda$ is a partial action
of $H$ on~$X$, $x \in \rfs{Fld}(\lambda)$
and $\alpha \in \rfs{MBC}$;
\item[\num{ii}] $G$ is $\alpha$-infinitely-closed at $x$;
\item[\num{iii}] $x$ is a $\lambda$-limit-point;
\item[\num{iv}] for some $N \in \rfs{Nbr}(x)$,
$\lambda$ is $(\alpha,G)$-decayable in $H_{\lambda}(x) \cap N$;
\item[\num{v}] $\lambda$ is traslation-like at $x$;
\item[\num{vi}] $\itGamma$ is a modulus of continuity
and $\itGamma$ is $(\leq\kern-3pt \kappa)$-generated,
where
$\kappa = \min(\{\kappa(x,V_{\lambda}(x))\,|\break
V \in \rfs{Nbr}(e_H))\})$;
\item[\num{vii}] $Y$ is a metric space and $\iso{\tau}{X}{Y}$;
\item[\num{viii}]
there is $U \in \rfs{Nbr}(x)$ such that
for every $g \in G\sprt{\kern1ptU}$:
if $g$ is $\alpha \scirc \alpha$-bicontinuous,
then $g^{\tau}$ is $\itGamma$-bicontinuous at $\tau(x)$.
\vspace{2.0pt}
\end{list}
\underline{Then} $\tau\inverse$ is $\itGamma$-continuous at $\tau(x)$.
\end{theorem}

\noindent
{\bf Proof } Combine Lemmas \ref{bddly-lip-bldr-p5.4}
and \ref{metr-bldr-l3.14}.
\smallskip\hfill\myqed

The above lemma will be used in the proof
that the derivative of a diffeomorphism $\tau$ is locally
$\itGamma$-continuous.
For groups of type $H_{\itGamma}^{\srfs{LC}}(X)$,
Theorem \ref{metr-bldr-t3.15} yields a result which is
slightly weaker than the result obtained in
Theorem \ref{metr-bldr-t3.27}, where the action is assumed to be
``affine-like'' rather than just ``translation-like''.

\begin{theorem}\label{metr-bldr-t3.16}
Let $\fourtpl{E}{X}{\calS}{\calF}$ be a subspace choice system,
$\itGamma$ be a \,\hbox{$(\leq\kern-0pt \kappa(E))$}-generated
modulus of continuity,
$Y$ be a metric space and $\iso{\tau}{X}{Y}$.
Suppose that
$(\rfs{LIP}(X;\calS,\calF))^{\tau} \subseteq
H_{\itGamma}^{\srfs{LC}}(Y)$.
Then $\tau\inverse$ is locally $\itGamma$-continuous.
\end{theorem}

\noindent
{\bf Proof }
Let $x \in X$ and $S \in \calS$ be such that $x \in S$.
Write $H = \bbT(E;F_S)$,
$\lambda = \lambda^{E;F_S}_{\bbT} \drest S$,
$G = \rfs{LIP}(X;S,F_S)$ and $\alpha(t) = 3t$.
We shall apply Theorem~\ref{metr-bldr-t3.15}.

By Lemma \ref{metr-bldr-p3.4}, $\lambda$ is
$(\alpha,G)$-decayable. So \ref{metr-bldr-t3.15}(iv) holds.
Let $V \in \rfs{Nbr}(e_H)$.
Then there is $r > 0$ such that $V_{\lambda}(x) \supseteq B^{F_S}(x,r)$.
Since $F_S$ is dense in $E$, $\kappa(F_S) = \kappa(E)$.
So $\kappa(x,V_{\lambda}(x)) = \kappa(F_S) = \kappa(E)$.
It follows that
$\min(\setm{\kappa(x,V_{\lambda}(x))}{V \in \rfs{Nbr}(e_H)}) =
\kappa(E)$.
Since $\itGamma$ is \hbox{$(\leq\kern-0pt \kappa(E))$}-generated,
\ref{metr-bldr-t3.15}(vi) holds.

Take $U$ in the definition of $\alpha$-infinite-closedness to be $S$.
Let $L$ be a subset of $G$ which satisfies Clauses (1)\,-\,(4) in the
definition of $\alpha$-infinite-closedness.
(See Definition \ref{metr-bldr-d3.9}(e)).
Then $\bcirc \kern-1.5ptL$ is $\alpha \scirc \alpha$-bicontinuous,
which implies that $\bcirc \kern-1.5ptL \in G$.
So $G$ is $\alpha$-infinitely-closed at $x$.
That is, \ref{metr-bldr-t3.15}(ii) holds.

Since for every $V \in \rfs{Nbr}(e_H)$
there is $r > 0$ such that $V_{\lambda}(x) \supseteq B^{F_S}(x,r)$,
$x$ is a $\lambda$-limit-point.
That is, \ref{metr-bldr-t3.15}(iii) holds.
By Proposition \ref{metr-bldr-p3.13}, $\lambda$ is translation-like
at $x$.  That is, \ref{metr-bldr-t3.15}(v) holds.
By the assumptions of this theorem,
\ref{metr-bldr-t3.15}(vii) and (viii) hold.

We have seen that all the assumptions of Theorem \ref{metr-bldr-t3.15}
are fulfilled, so $\tau\inverse$ is $\itGamma$-continuous
at $\tau(x)$.
\hfill\myqed

\begin{defn}\label{metr-bldr-d3.17}
\begin{rm}
{\thickmuskip=2mu \medmuskip=1mu \thinmuskip=1mu 
(a)
Let $E$ be a normed space, $S \subseteq X \subseteq E$ be open subsets
and $F$ be a dense linear subspace of $E$.
}
Let $\itGamma$ be a modulus of continuity.
We define
$$
H_{\itGamma}(X) =
\setm{h \in H(X)}{\mbox{there is } \gamma \in \itGamma
\mbox{ such that } h \mbox{ is } \gamma\mbox{-bicontinuous}},$$
$$H_{\itGamma}(X,S) =  H_{\itGamma}(X)\sprt{S},$$
$$
H_{\itGamma}(X;F) =
\setm{h \in  H_{\itGamma}(X)}{h(X \cap F) = X \cap F}
$$
and 
$$
H_{\itGamma}(X;S,F) = H_{\itGamma}(X,S) \cap H_{\itGamma}(X;F).
$$
Similarly, let
$H_{\itGamma}^{\srfs{LC}}(X,S) =  H_{\itGamma}^{\srfs{LC}}(X)\sprt{S}$,
$H_{\itGamma}^{\srfs{LC}}(X;F) =
\setm{h \in  H_{\itGamma}^{\srfs{LC}}(X)}{h(X \cap F) = X \cap F}$
and 
$H_{\itGamma}^{\srfs{LC}}(X;S,F) =
H_{\itGamma}^{\srfs{LC}}(X,S) \cap H_{\itGamma}^{\srfs{LC}}(X;F)$.

Let $\fourtpl{E}{X}{\calS}{\calF}$ be a subspace choice system.
We define
$H_{\itGamma}(X;\calS,\calF)$
to be the subgroup of $H(X)$ generated by
$\bigcup \setm{H_{\itGamma}(X;S,F_S)}{S \in \calS}$.
Analogously, the group $H_{\itGamma}^{\srfs{LC}}(X;\calS,\calF)$
is defined to be the subgroup of $H(X)$ generated by
$\bigcup \setm{H_{\itGamma}^{\srfs{LC}}(X;S,F_S)}{S \in \calS}$.

   \index{N@h03@@$H_{\itGamma}(X) =
          \setm{h \in H(X)}{\mbox{there is } \gamma \in \itGamma
          \mbox{ such that } h \mbox{ is } \gamma\mbox{-bicontinuous}}$}
   \index{N@h4@@$H_{\itGamma}(X,S) = H_{\itGamma}(X)\sprt{S}$}
   \index{N@h05@@$H_{\itGamma}(X;F) =
          \setm{h \in H_{\itGamma}(X)}{h(F \cap X) = F \cap X}$}
   \index{N@h06@@$H_{\itGamma}(X;S,F) =
          H_{\itGamma}(X,S) \cap H_{\itGamma}(X;F)$}
   \index{N@h07@@$H_{\itGamma}(X;\calS,\calF)$.
          The subgroup of $H(X)$ generated by
          {\thickmuskip=2mu \medmuskip=1mu \thinmuskip=1mu 
          $\bigcup \setm{H_{\itGamma}(X;S,F_S)}{S \in \calS}$}}
   \index{N@hlc01@@$H_{\itGamma}^{\srfs{LC}}(X)$.
          The group of locally $\itGamma$-bicontinuous homeomorphisms
	  of $X$}
   \index{N@hlc02@@$H_{\itGamma}^{\srfs{LC}}(X,S) =
          H_{\itGamma}^{\srfs{LC}}(X)\sprt{S}$}
   \index{N@hlc03@@$H_{\itGamma}^{\srfs{LC}}(X;F) =
          \setm{h \in H_{\itGamma}^{\srfs{LC}}(X)}{h(F \cap X) =
          F \cap X}$}
   \index{N@hlc04@@$H_{\itGamma}^{\srfs{LC}}(X;S,F) =
          H_{\itGamma}^{\srfs{LC}}(X,S) \cap
	  H_{\itGamma}^{\srfs{LC}}(X;F)$}
   \index{N@hlc05@@ {\thickmuskip=2mu \medmuskip=1mu \thinmuskip=1mu
          $H_{\itGamma}^{\srfs{LC}}(X;\calS,\calF)$.}
The subgroup $H(X)$ generated by
{\thickmuskip=2mu \medmuskip=1mu \thinmuskip=1mu 
$\bigcup \setm{H_{\itGamma}^{\srfs{LC}}(X;S,F_S)}{S \in \calS}$}}

(b) Let $E$ be a normed space, $z \in E$ and $\eta \in H([0,\infty))$.
Define
$h = \rfs{Rad}^E_{\eta,z}$ as follows.
$$
h(x) = z + \eta(\norm{x - z}) \frac{x - z}{\norm{x - z}}, \ \,%
x \neq z
$$
and $h(z) = z$.
Clearly, $h \in H(E)$. We call $h$ the
{\it radial homeomorphism based on $\eta,z$}.
   \index{radial homeomorphism. $\rfs{Rad}^E_{\eta,z}$.
          The radial homeomorphism based on $\eta,z$}
   \index{N@rad00@@$\rfs{Rad}^E_{\eta,z} =
          z + \eta(\norm{x - z}) \frac{x - z}{\norm{x - z}}$.
          The radial homeomorphism based on $\eta,z$}
Also, denote $\rfs{Rad}^E_{\eta,0^E}$ by $\rfs{Rad}^E_{\eta}$,
and call it the {\it radial homeomorphism based on $\eta$}.
          \index{N@rad01@@$\rfs{Rad}^E_{\eta} = \rfs{Rad}^E_{\eta,0^E}$}
   \index{radial homeomorphism based on $\eta$. $\rfs{Rad}^E_{\eta}$.}
\end{rm}
\end{defn}

\noindent
{\bf Remark } Note the following facts.
\begin{itemize}
\item[(1)] $H_{\itGamma}(X)$ is a special case of
$H_{\itGamma}(X;\calS,\calF)$,
where $\calS = \sngltn{X}$ and $F_X = E$. The same holds for
$H_{\itGamma}^{\srfs{LC}}(X)$.
\item[(2)]
$H_{\itGamma}(X,S), \,%
H_{\itGamma}(X;F), \,%
H_{\itGamma}(X;S,F), \,%
H_{\itGamma}(X;\calS,\calF) \subseteq
H_{\itGamma}(X)$.
\item[(3)]
$H_{\itGamma}^{\srfs{LC}}(X,S), \,%
H_{\itGamma}^{\srfs{LC}}(X;F), \,%
H_{\itGamma}^{\srfs{LC}}(X;S,F), \,%
H_{\itGamma}^{\srfs{LC}}(X;\calS,\calF) \subseteq
H_{\itGamma}^{\srfs{LC}}(X)$.
\end{itemize}

\begin{prop}\label{metr-bldr-p3.18}
Let $E$ be a normed space, $z \in E$
and $\eta \in H([0,\infty))$.
Suppose that $\eta$ is $\alpha$-bicontinuous.
Then $h_{\eta,z}$ is $3 \mcdot \alpha$-bicontinuous.
\end{prop}

\noindent
{\bf Proof }
Set $h = \rfs{Rad}_{\eta,z}$. We may assume that $z = 0$.
Note that $\eta(t) \leq \alpha(t)$ for every $t \geq 0$.
Since $\alpha$ is concave, it follows that
$\frac{\alpha(t)}{t} \mcdot s \leq \alpha(s)$
for every $0 < s \leq t$.

Let $u,v \in E - \sngltn{0}$. Assume that $\norm{u} \leq \norm{v}$
and set $w = \frac{\norm{u}}{\norm{v}} v$.
Then
$\norm{w - u} \leq \norm{u} + \norm{w} = 2 \norm{u}$.
So $\norm{\frac{w - u}{2}} \leq \norm{u}$.
Also,
$\norm{v - w} = \norm{v} - \norm{u} \leq \norm{v - u}$.
So
$\norm{w - u} \leq \norm{v - u} + \norm{v - w} \leq 2 \norm{v - u}$.
Hence
\vspace{1.5mm}
\newline
\rule{7pt}{0pt}
\renewcommand{\arraystretch}{1.5}
\addtolength{\arraycolsep}{-6pt}
$
\begin{array}{ll}
&
\norm{h(v) - h(u)} \leq \norm{h(v) - h(w)} + \norm{h(w) - h(u)}
\\
=
\rule{1pt}{0pt}
&
(\eta(\norm{v}) - \eta(\norm{w})) +
\frac{\eta(\norm{u})}{\norm{u}} \norm{w - u} =
(\eta(\norm{v}) - \eta(\norm{u})) +
2 \mcdot \frac{\eta(\norm{u})}{\norm{u}} \norm{\frac{w - u}{2}}
\\
\leq
\rule{1pt}{0pt}
&
\alpha(\norm{v} - \norm{u}) + 2 \alpha(\frac{\norm{w - u}}{2}) \leq
\alpha(\norm{v - u}) + 2 \alpha(\norm{v - u}) =
3 \alpha(\norm{v - u}).
\vspace{1.7mm}
\end{array}
$
\newline
\renewcommand{\arraystretch}{1.0}
\addtolength{\arraycolsep}{6pt}
So $h$ is $3 \alpha$-continuous.
Since $h\inverse = \rfs{Rad}_{\eta\inverse,z}$,
it follows that $h\inverse$ is $3 \alpha$-continuous.
\smallskip\hfill\myqed

The main result of the next theorem is Part (a).
It is a more readable special case of~(b).
Part (b) is a trivial corollary of (c).
The proof of (c) is more than just collecting some of the
previous lemmas together. It requires an additional argument.

\begin{theorem}\label{metr-bldr-t3.19}
\num{a} Let $X,Y$ be open subsets of the normed spaces $E$ and
$F$ respectively.
Write $\kappa = \kappa(E)$ and let $\itGamma,\itDelta$ be
\hbox{$(\leq \kappa)$}-generated moduli of continuity.
Let $\iso{\tau}{X}{Y}$, and suppose that
$(H_{\itGamma}^{\srfs{LC}}(X))^{\tau} = H_{\itDelta}^{\srfs{LC}}(Y)$.
Then $\itGamma = \itDelta$ and $\tau$ is locally
$\itGamma$-bicontinuous.
\smallskip
\smallskip

\num{b} Let $\fourtpl{E}{X}{\calS}{\calE}$
and $\fourtpl{F}{Y}{\calT}{\calF}$ be subspace choice systems.
Write $\kappa = \kappa(E)$ and let $\itGamma,\itDelta$ be
\hbox{$(\leq \kappa)$}-generated moduli of continuity.
Let $\iso{\tau}{X}{Y}$, and suppose that:
\newline
\num{i} $(H_{\itGamma}(X;\calS,\calF))^{\tau} \subseteq
H_{\itDelta}^{\srfs{LC}}(Y)$,
\num{ii}
$(H_{\itDelta}(Y;\calT,\calF))^{\tau\inverse} \subseteq
H_{\itGamma}^{\srfs{LC}}(X)$.
\newline
\underline{Then} $\itGamma = \itDelta$ and $\tau$ is
locally $\itGamma$-bicontinuous.
\smallskip

\num{c} Let $\fourtpl{E}{X}{\calS}{\calE}$
and $\fourtpl{F}{Y}{\calT}{\calF}$ be subspace choice systems.
Write $\kappa = \kappa(E)$ and let $\itGamma,\itDelta$ be
\hbox{$(\leq \kappa)$}-generated moduli of continuity.
Let $\iso{\tau}{X}{Y}$, and suppose that:
\begin{itemize}
\item[\num{i}]
$(\rfs{LIP}(X;\calS,\calF))^{\tau} \subseteq
H_{\itDelta}^{\srfs{LC}}(Y)$,
\item[\num{ii}]
$(H_{\itDelta}(Y;\calT,\calF))^{\tau\inverse} \subseteq
H_{\itGamma}^{\srfs{LC}}(X)$.
\end{itemize}
Then $\itDelta \subseteq \itGamma$ and $\tau$ is
locally $\itGamma$-bicontinuous.
\end{theorem}

\noindent
{\bf Proof }
Part (a) is a special case of (b),
and (b) is concluded by applying (c) twice:
once to $X,Y$ and once to $Y,X$. 
So it suffices prove (c).

(c)
Since $X$ and $Y$ are homeomorphic, $\kappa(F) = \kappa(E) = \kappa$.
Suppose by way of contradiction that $\itDelta \not\subseteq \itGamma$.
Pick any $T \in \calT$ and $y \in T \cap F_T$, and set
$x = \tau\inverse(y)$.
(Recall that $F_T$ denotes the dense subspace of $F$ assigned to $T$
by the subspace choice system).
Let $x \in S \in \calS$.
By Theorem~\ref{metr-bldr-t3.16} and Clause (c)(i),
for some $\delta \in \itDelta$,
$\tau\inverse$ is $\delta$-continuous at $\tau(x)$.
There is $\alpha \in (\itDelta - \itGamma) \cap \rfs{MBC}$ such that
$\delta \preceq \alpha$.
So $\tau\inverse$ is $\alpha$-continuous at $\tau(x)$.
Choose $r > 0$ be such that
$\tau\inverse \nrestriction B^F(y,r)$ is $\alpha$-continuous and
$B^F(y,r) \subseteq \tau(S) \cap T$,
and let $e$ be such that $\alpha \scirc \alpha(e) = \dgfrac{r}{2}$.
We define $\fnn{\eta}{[0,\infty)}{[0,\infty)}$ as follows.
For $t \in [0,e]$, $\eta(t) = \alpha \scirc \alpha(t)$,
for $t \in [r,\infty)$, $\eta(t) = t$,
$\eta \nrestriction [e,r]$ is a linear function,
and $\eta$ is continuous.
Clearly, $\eta \in H([0,\infty))$,
and it is easily seen that
$\eta$ is $4 \mcdot \alpha \scirc \alpha$-continuous
and that $\eta\inverse$ is $2$-Lipschitz.
So $\eta$ is $4 \mcdot \alpha \scirc \alpha$-bicontinuous.
Let $h = \rfs{Rad}_{\eta,y} \nrestriction Y$.
By Proposition~\ref{metr-bldr-p3.18},
$h$ is $12 \mcdot \alpha \scirc \alpha$-bicontinuous,
hence $h \in H_{\itDelta}(Y)$.
Since $y \in F_T$, we have $h(Y \cap F_T) = Y \cap F_T$,
and so $h \in H_{\itDelta}(Y;\calT,\calF)$.
By Clause~(c)(ii),
$g \eqdf h^{\tau\inverse}$ is locally $\itGamma$-bicontinuous,
and by Theorem \ref{metr-bldr-t3.16} and Clause (c)(ii),
$\tau$ is locally $\itGamma$-continuous.
This implies that $\tau \scirc g$ is locally $\itGamma$-continuous.
Since $h \scirc \tau = \tau \scirc g$,
we conclude that $h \scirc \tau$ is locally $\itGamma$-continuous.
Let $\gamma \in \itGamma$ be such that
$h \scirc \tau$ is $\gamma$-continuous at~$x$,
and choose $s$ such that $h \scirc \tau \nrestriction B^E(x,s)$
is $\gamma$-continuous.
We may assume that $\tau(B^E(x,s)) \subseteq B^F(y,\dgfrac{r}{2})$.

Since $\alpha \not\in \itGamma$,
there is $t < s$ such that $\alpha(t) > \gamma(t)$.
Choose $w$ such that $\norm{w - x} = t$
and set $z = \tau(w)$. Then $z \in B^F(y,\dgfrac{r}{2})$
and hence $\norm{h(z) - h(y)} = \alpha \scirc \alpha(\norm{z - y})$.
Now,\break
$\norm{w - x} = \norm{\tau\inverse(z) - \tau\inverse(y)} \leq 
\alpha(\norm{z - y})$.
So
$\alpha\inverse(\norm{w - x}) \leq \norm{z - y}$
and hence\break

\kern-29pt
$$
\norm{h(z) - h(y)} = \alpha \scirc \alpha(\norm{z - y}) \geq
\alpha \scirc \alpha(\alpha\inverse(\norm{w - x})) =
\alpha(\norm{w - x}).
$$
That is,
$\norm{h \scirc \tau(w) - h \scirc \tau(x)} \geq
\alpha(\norm{w - x}) > \gamma(\norm{w - x})$.
This contradicts the fact that $h \scirc \tau \nrestriction B^E(x,s)$
is $\gamma$-continuous. So $\itDelta \subseteq \itGamma$.

Since $\tau\inverse$ is locally $\itDelta$-continuous,
$\tau\inverse$ is locally $\itGamma$-continuous.
Recall also that $\tau$ is locally $\itGamma$-continuous.
So $\tau$ is locally $\itGamma$-bicontinuous.
\smallskip\hfill\myqed

\noindent

\begin{remark}\label{metr-bldr-r3.20}
\begin{rm}
The assumptions of Theorem \ref{metr-bldr-t3.19}(c)
probably imply that $\tau$ is locally $\itDelta$-bicontinuous.
We do not know to prove this fact.
However, the final result is not affected.\break
We also do not know to prove Theorem \ref{metr-bldr-t3.19}(a)
without the assumption that
$\itGamma,\itDelta$ are \hbox{$(\leq \kappa(E))$}-generated.
\hfill\proofend
\end{rm}
\end{remark}

There is a variant of translation-likeness which we shall use
in the context of diffeomorphisms.
Suppose that $f,g \in \rfs{Diff}([0,1])$.
If the derivative $f'$ of $f$ is $\alpha$-continuous
and $g'$ is $\beta$-continuous,
then (i) for some $K,L > 0$,
$(f \scirc g)'$ is $(K \mcdot \alpha + L \mcdot \beta)$-continuous.
Also, (ii) for some $M > 0$,
$(f\inverse)'$ is $M \mcdot \alpha$-continuous.
(iii) A similar fact holds for higher derivatives.

Let $\itGamma \subseteq \rfs{MC}$, and assume that
$K \mcdot \alpha + L \mcdot \beta \in \itGamma$
for every $\alpha,\beta \in \itGamma$ and $K,L > 0$.
Consider the set
$G_{\itGamma} = \setm{f \in \rfs{Diff}([0,1])}{\mbox{for some }
\alpha \in \itGamma, \ f' \mbox{ is } \alpha\mbox{-continuous}}$.
By (i)-(ii), $G_{\itGamma}$ is a group,
and by (iii), the analogous fact for
$\rfs{Diff}^{\kern1.5pt n}([0,1])$ is also true.
So $\itGamma$ need not be a modulus of continuity in order for
$G_{\itGamma}$ to be a group.
Let us call such a $\itGamma$ a modulus of differentiability.

We do not deal in this work with differentiability,
but we shall show here that if $\itGamma$ is a modulus of
differentiability and
$(\rfs{LIP}(X))^{\tau} \subseteq H_{\itGamma}^{\srfs{LC}}(Y)$,
then $\tau\inverse$ is locally $\itGamma$-continuous.
This is the analogue of Theorem \ref{metr-bldr-t3.16},
and Theorem \ref{metr-bldr-t3.15} has an analogue too.
The proofs use the additional assumptions that $X$ is of the second
category, and that $\itGamma$ is countably generated.
On the other hand, the infinite-closedness of $G$ is not needed,
and the assumption of decayability is replaced by a much weaker
property.

\begin{defn}\label{metr-bldr-d3.21}
\begin{rm}
Let $X$ be a topological space, $\lambda$ be a partial action of
a topological group $H$ on $X$ and $G \leq H(X)$.
Let $x \in X$. We say that $\lambda$ is
{\it compatible with $G$ at $x$},
if there is $W \in \rfs{Nbr}(e_H)$ such that for every $h \in W$
there are $U \in \rfs{Nbr}(x)$ and $g \in G$ such that
$U \subseteq \rfs{Dom}(h_{\lambda})$ and
$h_{\lambda} \nrestriction U = g \nrestriction U$.
   \index{compatible. $\lambda$ is compatible with $G$ at $x$}

We say that $\lambda$ is {\it compatible with $G$},
if $\lambda$ is compatible with $G$ at every
$x \in \rfs{Fld}(\lambda)$.
\hfill\proofend
   \index{compatible. $\lambda$ is compatible with $G$}
\end{rm}
\end{defn}

The following lemma replaces Lemma \ref{metr-bldr-l3.11}.

\begin{lemma}\label{metr-bldr-l3.22}
Suppose that:
\begin{itemize}
\addtolength{\parskip}{-09pt}
\addtolength{\itemsep}{04pt}
\item[\num{i}]
$X$ is a metric space, $G \leq H(X)$,
$H$ is a topological group and $H$ is of the second category,
$\lambda$ is a partial action of $H$ on $X$,
$x \in \rfs{Fld}(\lambda)$,
and $\lambda$ is  compatible with $G$ at~$x$.
\item[\num{ii}]
$Y$ is metric space and $\iso{\tau}{X}{Y}$.
\item[\num{iii}]
$\itGamma$ is a countably generated subset of $ \rfs{MC}$,
$\rfs{cl}_{\preceq}(\sngltn{\gamma}) \subseteq \itGamma$
and $K \mcdot \gamma \in \itGamma$
for every $\gamma \in \itGamma$ and $K > 0$.
\item[\num{iv}]
For every $g \in G$,
$g^{\tau}$ is $\itGamma$-bicontinuous at $\tau(x)$.
\end{itemize}
\underline{Then} $Q(x)$ holds, where
\smallskip
\newline
$Q(x)$: \,%
For every $W \in \rfs{Nbr}(e_H)$ there are
$T \in \rfs{Nbr}(x)$, a nonempty open subset $V \subseteq W$
and $\gamma \in \itGamma$ such that for every $h \in V$:
$T \subseteq \rfs{Dom}(h_{\lambda})$
and $(h_{\lambda})^{\tau} \nrestriction \tau(T)$
is $\gamma$-bicontinuous.
\end{lemma}

\noindent
{\bf Proof }
For every $h \in H$ denote $h_{\lambda}$ by $\hath$.
Let $W \in \rfs{Nbr}(e_H)$.
We may assume that for every $h \in W$
there are $U_h \in \rfs{Nbr}(x)$ and $g_h \in G$ such that
$U_h \subseteq \rfs{Dom}(\hath)$ and
$h_{\lambda} \nrestriction U = g_h \nrestriction U$.

We verify that $(*)$ for every $h \in W$ there are $r_h > 0$
and $\gamma_h \in \itGamma$ such that
$B(x,r_h) \subseteq \rfs{Dom}(\hath)$ and
$\hath^{\tau} \nrestriction \tau(B(x,r_h))$ is $\gamma_h$-bicontinuous.
Let $U_h,g_h$ be as above. Then $(g_h)^{\tau}$
is $\itGamma$-bicontinuous at $\tau(x)$.
Let $\gamma_h \in \itGamma$ and $T \in \rfs{Nbr}(\tau(x))$
be such that $(g_h)^{\tau} \nrestriction T$ is $\gamma_h$-bicontinuous,
and let $r_h > 0$ be such that
$B(x,r_h) \subseteq U_h$ and $\tau(B(x,r_h)) \subseteq T$.
Obviously,
$$
\hath^{\tau} \nrestriction \tau(B(x,r_h)) =
(\hath \nrestriction B(x,r_h))^{\tau} =
(g_h \nrestriction B(x,r_h))^{\tau} =
(g_h)^{\tau} \nrestriction \tau(B(x,r_h)).
$$
So $\hath^{\tau} \nrestriction \tau(B(x,r_h)$ is $\gamma$-bicontinuous.
That is, $(*)$ holds.

Let $ \itGamma_0 = \setm{\gamma_i}{i \in \bbN}$ be such that
$\itGamma = \rfs{cl}_{\preceq}(\itGamma_0)$, and assume that
$\setm{j}{\gamma_j = \gamma_i}$ is infinite
for every $i \in \bbN$.
Set
$$
K_i = \setm{h \in W}
{B(x,\frac{1}{i + 1}) \subseteq \rfs{Dom}(\hath) \mbox{ and }
\hath^{\tau} \nrestriction \tau(B(x,\frac{1}{i + 1}))
\mbox{ is } \gamma_i\mbox{-bicontinuous}}.
$$
By $(*)$, \,$\bigcup_{i \in \sboldbbN} K_i = W$.
We show that for every $i$, $K_i$ is closed in $W$.
Set $B_i = B(x,\frac{1}{i + 1}))$.
Let $h \in W - K_i$.
So there are $y_1,y_2 \in \tau(B_i)$ such that
(i) $d(\hath^{\tau}(y_1),\hath^{\tau}(y_2)) > \gamma_i(d(y_1,y_2))$
or
(ii)
$d(\hath^{\tau}(y_1),\hath^{\tau}(y_2)) <
\gamma_i\inverse(d(y_1,y_2))$.
We may assume that (i) happens.
For $\ell = 1,2$ let $T_{\ell}$ be an open neighbourhood of
$\hath^{\tau}(y_\ell)$ such that
$d(T_1,T_2) > \gamma_i(d(y_1,y_2))$.
Set $S_{\ell} = \tau\inverse(T_{\ell})$
and $x_{\ell} = \tau\inverse(y_{\ell})$.
Let
$V_0 =
\setm{k \in W}
{x_1,x_2 \in \rfs{Dom}(\hatk)\,,\ %
\hatk(x_1) \in S_1 \mbox{ and } \hatk(x_2) \in S_2}$.
So $V_0$ is open.
We show that $V_0$
contains $h$ and is disjoint from $K_i$.
Clearly, $\hath(x_{\ell}) =
\tau\inverse(\hath^{\tau}(y_{\ell})) \in \tau\inverse(T_{\ell}) =
S_{\ell}$,
hence $h \in V_0$.
If $k \in V_0$, then $\hatk(x_{\ell}) \in S_{\ell}$ and so
$\hatk^{\tau}(y_l) \in \tau(S_{\ell}) = T_{\ell}$.
Hence $\hatk^{\tau} \nrestriction \tau(B_i)$ is not
$\gamma_i$-continuous, namely, $k \not\in K_i$.
Since $W$ is of the second category and every $K_n$ is closed,
there is $n$ such that $\rfs{int}(K_n) \neq \emptyset$.
Define $V = \rfs{int}(K_n)$, $T = B(x,\frac{1}{n + 1})$
and $\gamma = \gamma_n$.
Then $V$, $T$ and $\gamma$ are as required in the lemma.
\hfill\myqed

\begin{defn}\label{d3.1}\label{metr-bldr-d.23}
\begin{rm}
Let $X$ be a metric space, $H$ be a topological group $\lambda$ be a
partial action of $H$ on $X$ and $x \in \rfs{Fld}(\lambda)$.
The action $\lambda$ is said to be
{\it regionally translation-like at $x$},
   \index{regionally translation-like at $x$}
if there is $W_x \in \rfs{Nbr}(e_H)$ such that
for every nonempty open $V \subseteq W_x$ such that
$V \times \sngltn{x} \subseteq \rfs{Dom}(\lambda)$ there are:
\begin{itemize}
\addtolength{\parskip}{-09pt}
\addtolength{\itemsep}{04pt}
\item[(i)] $U = U_{x,V} \in \rfs{Nbr}(x)$
and a dense subset of $U$, $D = D_{x,V}$;
\item[(ii)] a point $z = z_{x,V}$, a radius $r = r_{x,V} > 0$,
and a constant $K =  K_{x,V} > 0$;
\end{itemize}
such that for every distinct $\barx_0,\barx_1 \in U \cap D$
there are $n \leq K \ncdot \frac{r}{d(\barx_0,\barx_1)}$\,,
a sequence $z = z_0, z_1, \ldots, z_n \in X$ and
$h_1, \ldots, h_n \in V$
such that
$z_n \not \in B(z,r)$, and for every $i = 1, \ldots, n$,\ \ %
$\barx_0,\barx_1 \in \rfs{Dom}((h_i)_{\lambda})$, \ %
$(h_i)_{\lambda}(\barx_0) = z_{i - 1}$
and $(h_i)_{\lambda}(\barx_1) = z_i$.

If $\lambda$ is regionally translation-like at every
$x \in \rfs{Fld}(\lambda)$, then $\lambda$ is said to be
a {\it regionally translation-like action}.
   \index{regionally translation-like action}
\end{rm}
\end{defn}

The next proposition is a counterpart of
Proposition~\ref{metr-bldr-p3.13}.

\begin{prop}\label{metr-bldr-p3.24}
Let $E$ be a normed vector space, $F$ be a dense linear subspace
of $E$ and $X$ be an open subset of $E$.
Then $\lambda_{\bbT}^{E;F} \drest X$ is regionally translation-like.
\end{prop}

\noindent
{\bf Proof }
Write $\lambda = \lambda_{\bbT}^{E;F} \drest X$
and define $W_x = \bbT(E;F)$. Let $V \subseteq W_x$ be open and
nonempty,
and suppose that $V \times \sngltn{x} \subseteq \rfs{Dom}(\lambda)$.
Choose $v \in F$ and $s > 0$ such that
$V_1 \eqdf \setm{\rfs{tr}^E_u}{u \in B^F(v,s)} \subseteq V$
and $V_1 \times B^E(x,s) \subseteq \rfs{Dom}(\lambda)$.
Define $z_{x,V} = v + x$, $r = r_{x,V} = \dgfrac{s}{2}$,\break
$U_{x,V} = B(x,\dgfrac{s}{4})$, $D_{x,V} = U_{x,V} \cap (x + F)$
and $K_{x,V} = 2$. It is left to the reader to verify that the above
satisfy the requirements of regional translation-likeness
of $\lambda$ at $x$.
\smallskip\hfill\myqed

The following lemma is a counterpart of Lemma \ref{metr-bldr-l3.14}.

\begin{lemma}\label{metr-bldr-l3.25}
Let $X$ be a metric space, and $\lambda$ be a partial action of
$H$ on~$X$.
Suppose that $x \in \rfs{Fld}(\lambda)$, and $\lambda$ is regionally
traslation-like at $x$.
Let $Y$ be a metric space and $\iso{\tau}{X}{Y}$.
Let $\itGamma \subseteq \rfs{MC}$, and suppose that for every
$\gamma \in \itGamma$ and $K > 0$, $K \mcdot \gamma \in \itGamma$.
Also assume that
$Q(x)$ of Lemma \ref{metr-bldr-l3.22} holds.
That is,
for every $W \in \rfs{Nbr}(e_H)$ there are
$U \in \rfs{Nbr}(x)$, a nonempty open subset $V \subseteq W$
and $\gamma \in \itGamma$ such that
$U \subseteq \rfs{Dom}(h_{\lambda})$ and
$(h_{\lambda})^{\tau} \nrestriction \tau(U)$ is
$\gamma$-bicon\-tinuous for every $h \in V$.
\underline{Then} $\tau\inverse$ is $\itGamma$-continuous at $\tau(x)$.
\end{lemma}

\noindent
{\bf Proof }
Let $W_x$ be as assured by the regional translation-likeness of
$\lambda$ at $x$.
By $Q(x)$, there are $U \in \rfs{Nbr}(x)$,
a nonempty open $V \subseteq W_x$
and $\gamma \in \itGamma$ such that
for every $h \in V$: $U \subseteq \rfs{Dom}(h_{\lambda})$
and $(h_{\lambda})^{\tau} \nrestriction \tau(U)$
is $\gamma$-bicon\-tinuous.
So $V \subseteq W_x$
and $V \times \sngltn{x} \subseteq \rfs{Dom}(\lambda)$.
We apply the definition of regional translation-likeness to $V$.
Write $S = U_{x,V}$, $D = D_{x,V}$, $z = z_{x,V}$, $r = r_{x,V}$
and $K =  K_{x,V}$.

Let
$w = \tau(z)$, $B = B(z,r)$ and $C = \tau(B)$.
Since $C \in \rfs{Nbr}(w)$, we conclude that $e \eqdf d(w,Y - C) > 0$.
Let $R = \tau(U \cap S)$
and $M = \dgfrac{Kr}{e}$.
Since $\gamma \in \itGamma$, we have $M \cdot \gamma \in \itGamma$.

We show that $\tau\inverse \nrestriction R$ is
$M \ncdot \gamma$-continuous. Suppose by contradiction
that this is not true.
For $h \in H$ denote $h_{\lambda}$ by $\hath$.
\hbox{Hence there are  $\bary_0,\bary_1 \in R$} such that
$d(\tau\inverse(\bary_0),\tau\inverse(\bary_1)) >
M \ncdot \gamma(d(\bary_0,\bary_1))$.
Since $D$ is dense in $S$ and
$\bary_0,\bary_1 \in \tau(S)$,
we may assume that
$\bary_0,\bary_1 \in \tau(D)$.
For $\ell = 0,1$ let
\hbox{$\barx_{\ell} = \tau\inverse(\bary_{\ell})$.}
Hence \hbox{$\barx_0,\barx_1 \in D$.}
So there are $n \leq \frac{Kr}{d(\barx_0,\barx_1)}$,
$z = z_0, z_1, \dots, z_n$
and
$h_1, \ldots, h_n \in V$ such that $z_n \not\in B$,
and for every $i = 1, \ldots, n$,
$\barx_0,\barx_1 \in \rfs{Dom}(\hath_i)$,
$\hath_i(\barx_0) = z_{i - 1}$
and $\hath_i(\barx_1) = z_i$.
For $i = 1, \ldots, n$\ \,%
let $w_i = \tau(z_i)$.

In the space $Y$ we have:
$w_0 = w$;
for every $i = 1, \ldots, n$,
$\hath_i^{\tau}(\bary_0) = w_{i - 1}$ and
$\hath_i^{\tau}(\bary_1) = w_i$;
and $w_n \not\in C$.
Every $h_i$ belongs to $V$, hence $\hath_i^{\tau} \,\rest\, \tau(U)$
is $\gamma$-bicontinuous.
Also, $\bary_0,\bary_1 \in \tau(U)$,
so $d(w_{i - 1},w_i) \leq \gamma(d(\bary_0,\bary_1))$.
Hence \medskip\\
\centerline{
\renewcommand{\arraystretch}{1.5}
\addtolength{\arraycolsep}{-3pt}
$
\begin{array}{rcl}
e &=&d(w,Y - C) \leq d(w,w_n)  = d(w_0,w_n) \leq
\sum^n_{i = 1} d(w_{i - 1},w_i)\\
&\leq&
n \cdot \gamma(d(\bary_0,\bary_1)) \leq
\frac{Kr}{d(\barx_0,\barx_1)} \cdot \gamma(d(\bary_0,\bary_1)) <
\frac{Kr}{M \cdot \gamma(d(\bary_0,\bary_1))} \cdot
\gamma(d(\bary_0,\bary_1))
=
\frac{Kr}{\dgfrac{Kr}{e}} = e.
\medskip
\end{array}
$
}
\renewcommand{\arraystretch}{1.0}
\addtolength{\arraycolsep}{3pt}
A contradiction, so the lemma is proved.
\hfill\myqed

\begin{theorem}\label{metr-bldr-t3.26}
Assume the following facts.
\begin{itemize}
\addtolength{\parskip}{-09pt}
\addtolength{\itemsep}{04pt}
\item[\num{i}]
$X$ is a metric space, $G \leq H(X)$, $H$ is a topological group
and $H$ is of the second category,
$\lambda$ is a partial action of $H$ on $X$
and $x \in \rfs{Fld}(\lambda)$.
\item[\num{ii}]
$\lambda$ is compatible with $G$ at $x$.
\item[\num{iii}]
$\lambda$ is regionally traslation-like at $x$.
\item[\num{iv}]
$\itGamma$ is a countably generated subset of $\rfs{MC}$,
$\rfs{cl}_{\preceq}(\sngltn{\gamma}) \subseteq \itGamma$,
and $K \mcdot \gamma \in \itGamma$
for every $\gamma \in \itGamma$ and $K > 0$.
\item[\num{v}]
$Y$ is metric space and $\iso{\tau}{X}{Y}$.
\item[\num{vi}]
For every $g \in G$, $g^{\tau}$ is $\itGamma$-bicontinuous at $\tau(x)$.
\end{itemize}
\underline{Then} $\tau\inverse$ is $\itGamma$-continuous at $\tau(x)$.
\end{theorem}

\noindent
{\bf Proof } Combine Lemmas \ref{metr-bldr-l3.22} and
\ref{metr-bldr-l3.25}.
\hfill\myqed

\subsection{Affine-like partial actions.}
\label{ss3.4}
\label{ss3.4-affine-like-actions}

The goal of this part of the chapter is the following final theorem.

\begin{theorem}\label{metr-bldr-t3.27}
{\thickmuskip=2mu \medmuskip=1mu \thinmuskip=1mu 
Let $\fourtpl{E}{X}{\calS}{\calE}$
be a subspace choice system with $\rfs{dim}(E) > 1$,
}
$Y$ be an open subset of a normed space $F$,
$\itGamma$ be a $(\leq \kappa(E))$-generated modulus of continuity
and $\iso{\tau}{X}{Y}$.
Suppose that
$(\rfs{LIP}(X,\calS,\calE))^{\tau} \subseteq
H_{\itGamma}^{\srfs{LC}}(Y)$.
Then $\tau$ is locally $\itGamma$-bicontinuous.
\end{theorem}

This parallels Theorem \ref{metr-bldr-t3.16},
but has a stronger conclusion.
Whereas in \ref{metr-bldr-t3.16} the conclusion is that $\tau\inverse$
is locally $\itGamma$-continuous, \ref{metr-bldr-t3.27} says that
$\tau$ is locally $\itGamma$-bicontinuous.

\begin{defn}\label{d3.4}\label{metr-bldr-d3.28}
\begin{rm}
(a) A subset $D$ of a metric space $X$ is called a
{\it metrically dense subset} of $X$,
if for every $x,y \in X$ and $\varepsilon > 0$ the are
$x_1 \in B(x,\varepsilon) \cap D$ and $y_1 \in B(y,\varepsilon) \cap D$
such that $d(x_1,y_1) = d(x,y)$. Note that metric density implies
density.
   \index{metrically dense subset}

(b) Let $X$ be a metric space,
$H$ be a topological group and $\lambda$ be a partial action of
$H$ on $X$. For $h \in H$ denote $h_{\lambda}$ by $\hath$.
Let $x \in X$.
We say that $\lambda$ is an {\it affine-like partial action at $x$},
if the following holds.
For every $V \in \rfs{Nbr}(e_H)$ and $U \in \rfs{Nbr}(x)$
there are
$n = n(x,V,U) \in \bbN$,
$U_0 = U_0(x,V,U) \in \rfs{Nbr}(x)$
and $D = D(x,V,U) \subseteq U_0$ such that
$U_0 \subseteq U$, $D$ is metrically dense in $U_0$,
and for every
$x_1,y_1,x_2,y_2 \in D$: if
$d(x_1,y_1) = d(x_2,y_2)$,
then there are $h_1, \ldots, h_n \in V$ such that
$\hath_1 \scirc \ldots \scirc \hath_n(x_1) = x_2$, \,%
$\hath_1 \scirc \ldots \scirc \hath_n(y_1) = y_2$
and $\hath_i \scirc \hath_{i + 1} \scirc \ldots \scirc
\hath_n(\braces{x_1,y_1}) \subseteq U$
for every $1 \leq i \leq n$.
   \index{affine-like partial action at $x$}

If $\lambda$ is affine-like at every $x \in \rfs{Fld}(\lambda)$,
then $\lambda$ is said to be an {\it affine-like partial action}.
   \index{affine-like partial action}

(c)
If $H$ is a group, $A \subseteq H$ and $n \in \bbN$,
then $A^n = \setm{a_1 \cdot \ldots \cdot a_n}{a_1, \ldots, a_n \in A}$.
Let $\lambda$ be a partial action of a topological group $H$ on a
topological space $X$. If $h \in H$ then $h_{\lambda}$ is denoted
by $\hath$.
For $U \subseteq H$ and $W_1,W_2 \subseteq X$ define
\medskip
\newline
\renewcommand{\arraystretch}{1.5}
\addtolength{\arraycolsep}{-6pt}
$
\begin{array}{rcl}
U^{[n;W_1,W_2]}&
=
&
\setm{h_1 \mcdot \ldots \mcdot h_n}{h_1, \ldots, h_n \in U,\ \,%
W_1 \subseteq \rfs{Dom}(\hath_i \scirc \ldots \scirc \hath_n)
\mbox{ and }
\\
&&
\hath_i \scirc \ldots \scirc \hath_n(W_1) \subseteq W_2
\mbox{ for every } i = 1, \ldots, n}.
\end{array}
$
\renewcommand{\arraystretch}{1.0}
\addtolength{\arraycolsep}{6pt}
\newline
\rule{0pt}{0pt}\hfill\proofend
   \index{N@AAAA@@$U^{[n;W_1,W_2]}$}
\end{rm}
\end{defn}

We shall prove two intermediate main facts. They roughly say the
following.

\begin{list}{}
{\addtolength{\parskip}{-00pt}
\addtolength{\topsep}{-07pt}
\addtolength{\itemsep}{-3pt}}
\item[(a)] If $X$ is an open subset of a normed space $E$,
and $F$ is a dense linear subspace of $E$,
then $\lambda_{\sboldbbA}^{E;F} \drest X$ is affine-like.
\item[(b)] Suppose that $\lambda$ is a decayable affine-like
partial action of $H$ on $X$,
$\iso{\tau}{X}{Y}$, $\itGamma$ is a countably generated modulus of
continuity,
and $(h_{\lambda})^{\tau}$ is locally $\itGamma$-bicontinuous
for every $h \in H$.
Then $\tau$ is locally $\itGamma$-bicontinuous.
\end{list}

We start with the proof of (a).
When proving the affine-likeness of $\lambda_{\sboldbbA}^{E;F} \drest X$
at $x$, it is easier to deal first with $x$'s which belong to
$F \cap X$.
To conclude that $\lambda_{\sboldbbA}^{E;F} \drest X$
is affine-like at every $x \in X$, we use an observation
which says that if $\lambda$ is affine-like at every $x \in C$,
and $U_0(x,V,U)$ and $n(x,V,U)$ depend on $x \in C$ and $V$
in some uniform way,
then $\lambda$ is affine-like at every $x \in \rfs{cl}(C)$.

\begin{prop}\label{metr-bldr-p3.29}
Assume the following facts.
\begin{itemize}
\addtolength{\parskip}{-09pt}
\addtolength{\itemsep}{04pt}
\item[\num{i}] 
$X$ is a metric space,
$\lambda$ is a partial action of $H$ on $X$,
$C \subseteq \rfs{Fld}(\lambda)$,
$r_0 > 0$,\break
$\fnn{\iota}{\rfs{Nbr}(e_H) \times C}{\rfs{Nbr}(e_H)}$,
$\fnn{\barn}{\rfs{Nbr}(e_H) \times (0,r_0)}{\bbN}$
and
$\fnn{\bars}{\rfs{Nbr}(e_H) \times (0,r_0)}{(0,\infty)}$.
Denote $\iota(V,y)$ by $V_y$.
\item[\num{ii}] 
For every $y \in C$, $\lambda$ is affine-like at $y$,
and for every $V \in \rfs{Nbr}(e_H)$ and $r \in (0,r_0)$,
$n(y,V_y,B(y,r)) \leq \barn(V,r)$ and
$U_0(y,V_y,B(y,r)) \supseteq B(y,\bars(V,r))$.
\item[\num{iii}] 
For every $x \in \rfs{cl}(C)$ and $W \in \rfs{Nbr}(e_H)$
there are  $U_1 \in \rfs{Nbr}(x)$ and $V \in \rfs{Nbr}(e_H)$
such that for every $y \in C \cap U_1$, $V_y \subseteq W$.
\end{itemize}
\underline{Then} for every $x \in \rfs{cl}(C)$,
$\lambda$ is affine-like at $x$.
Also, if $r < r_0$, then $n(x,V,B(x,r))$ and $U_0(x,V,B(x,r))$
can be taken to be $\barn(V,\dgfrac{r}{2})$ and
$B(x, \half \bars(V,\dgfrac{r}{2}))$ respectively.
\end{prop}

\noindent
{\bf Proof }
Let $x \in \rfs{cl}(C)$, $W \in \rfs{Nbr}(e_H)$,
$r \in (0,r_0)$
and $U = B(x,r)$.
There is $V \in \rfs{Nbr}(e_H)$ and $U_1 \in \rfs{Nbr}(x)$
such that for every $y \in U_1 \cap C$, $V_y \subseteq W$.
Define
$U_0 =  U_0(x,W,U)$ to be $B(x,\half \bars(V,\dgfrac{r}{2}))$.
Let
$y \in C \cap U_0 \cap U_1$.
Then
$U_0 \subseteq B(y,\bars(V,\dgfrac{r}{2})) \subseteq
U_0(y,V_y,B(y,\dgfrac{r}{2}))$.
Hence $D(y,V_y,B(y,\dgfrac{r}{2})) \cap U_0$
is metrically dense in $U_0$.
Let $D = D(x,W,V) = D(y,V_y,B(y,\dgfrac{r}{2})) \cap U_0$
and $n = n(x,W,U) = \barn(V,\dgfrac{r}{2})$.
We show that $U_0,D$ and $n$ fullfill the requirements
of affine-likeness.

{\thickmuskip=3mu \medmuskip=2mu \thinmuskip=1mu 
Let $x_1,x_2,y_1,y_2 \in D$ be such that
$d(x_1,y_1) = d(x_2,y_2)$.
Let $h_1,\ldots,h_n\in V_y$
}
be as assured by the affine-likeness of
$\lambda$ at~$y$.
So for every $i = 1,\ldots,n$,
$\hath_i \scirc \ldots \scirc \hath_n(\dbltn{x_1}{y_1}) \subseteq
B(y,\dgfrac{r}{2})$.
Clearly, $\bars(V,\dgfrac{r}{2}) \leq \dgfrac{r}{2}$
and $d(x,y) < \bars(V,\dgfrac{r}{2})$.
So $B(y,\dgfrac{r}{2}) \subseteq B(x,r) = U$.
Since $y \in U_1$, $V_y \subseteq W$.
So $h_1,\ldots,h_n$ fulfill the requirements needed in demonstrating
that $\lambda$ is affine-like at $x$.
\hfill\myqed

If $X$ is an open subset of $\bbR$, then $\bbA(\bbR) \drest X$ is not 
affine-like. So in what follows we assume that $\rfs{dim}(E) > 1$.

The group $\bbL(E)$ has a property similar to affine-likeness.
But the ``affine-likeness'' of $\bbL(E)$ applies only to pairs of pairs
$x_1,y_1,x_2,y_2$ in which $x_1 = x_2 = 0^E$.

\kern1mm

\begin{lemma}\label{metr-bldr-l3.30}
Let $E$ be a normed space with
dimension $> 1$, $E_1$ \hbox{be a dense}
linear subspace of $E$
and $V \in \rfs{Nbr}^{\sboldbbL(E;E_1)}(\rfs{Id})$.
\kern-3ptThen there is $n = n(V) \in \bbN$
{\thickmuskip=2.5mu \medmuskip=2mu \thinmuskip=1mu 
with the following property:\kern2.2pt
$(*)$\kern0.9pt For every \,$W_1 \in \rfs{Nbr}^E(0)$
there is $W_2 \in \rfs{Nbr}^E(0)$
}
such that $W_2 \subseteq W_1$
and for every $x_1,x_2 \in W_2 \cap E_1$:
if $\norm{x_1} = \norm{x_2}$,
then there is $S \in V^{[n;W_2,W_1]}$ such that $S(x_1) = x_2$.

Moreover, if in the above $V = B^{\sboldbbL(E;E_1)}(\rfs{Id},r)$ and
$W_1 = B^E(0,s)$, then $W_2$ can be taken to
be $B^E\kern-3.0pt\left(0,\frac{s}{(1 + r)^{n(V)}}\right)$.
\end{lemma}

\noindent
{\bf Proof } The proof of the lemma relies on three easy claims.

{\bf Claim 1.} Let $\itbfH^2$ be the $2$-dimensional
Hilbert space.
For every $K \geq 1$ and\break
$V \in \rfs{Nbr}^{\sboldbbL(\itbfH^2)}(\rfs{Id})$
there is $n = n(V,K) \in \bbN$
such that for every $x_1,x_2 \in \itbfH^2$:
if\break
$\frac{1}{K} \leq \frac{\raise 2.5pt \hbox{\scriptsize$\dline{x_1}$}}
{\lower 2.7pt \hbox{\scriptsize$\dline{x_2}$}} \leq K$,
then there is $T \in V^n$ such that $T(x_1) = x_2$.

{\bf Proof }
We may assume that $V = V\inverse$. For some angle $\gamma_0 > 0$,
$U$ contains all rotations $Rot_{\gamma}$, $\gamma \in [0,\gamma_0]$.
For some $\varepsilon_0 > 0$, $U$ contains all isomorphisms
$T_{\varepsilon}(x) = (1 + \varepsilon)x$
where $\varepsilon \in [0,\varepsilon_0]$.
It is left to the reader to verify that
$n(U,K) =
[\dgfrac{\pi}{\gamma_0}] +
\dgfrac{\log K}{\log(\varepsilon_0 + 1)} + 2$
is as required.\smallskip

We do not prove Claim 2 which is well-known and easy. In fact,
the best possible constant in Claim 2 is $\sqrt{2}$.

{\bf Claim 2.} For every $2$-dimensional normed space $E$ there is
an isomorphism $T$ between $E$ and the $2$-dimensional Hilbert space
$\itbfH^2$ such that
$\norm{T} \leq 1$ and $\norm{T\inverse} \leq 3\sqrt{2}$.

{\bf Claim 3.} Let $E$ be a normed space, $E_1$ be a dense linear
subspace of $E$, $F$ be a $2$-dimensional linear subspace of $E_1$
and $T \in \bbL(F)$, then there is $T_1 \in \bbL(E;E_1)$
extending $T$ such that
$d(T_1,\rfs{Id}_E) \leq 3 d(T,\rfs{Id}_F)$.

{\bf Proof } Let $x_1,x_2$ be a basis for $F$ such that
$\norm{x_1} = d(x_1,\rfs{span}(\sngltn{x_2}))$.
For $i =1,2$ let $\varphi_1,\varphi_2 \in F^{\raisedstar}$ be
such that $\varphi_i(x_j) = \delta_{i,j} \mcdot \norm{x_j}$,
and let
$\psi_i \in E^{\raisedstar}$ be such that
$\psi_i$ extends $\varphi_i$ and $\norm{\psi_i} = \norm{\varphi_i}$.
Set $F_1 =
\bigcap_{\kern0.8pt i = 1}^{\kern0.8pt 2} \rfs{ker}(\varphi_i)$,
hence $F \oplus F_1 = E$.
For $x \in E$ let $\hatx \in F$ and $\barx \in F_1$ denote the
components of $x$ in $F$ and $F_1$ respectively.
If $\hatx = ax_1 + bx_2$, denote $ax_1$ and $bx_2$ by $\hatx^1$ and
$\hatx^2$ respectively.
Let $x \in F$. Then
$\abs{\varphi_1(\hatx)} = \norm{\hatx^1} =
d(\hatx,\rfs{span}(\sngltn{x_2})) \leq \norm{\hatx}$.\break
So $\norm{\varphi_1} \leq 1$. Hence $\norm{\psi_1} \leq 1$.
It follows that
$\norm{\hatx^1} = \abs{\psi_1(x)} \leq \norm{x}$.
Also,\break
$\norm{\hatx^2} \leq \norm{\hatx} + \norm{\hatx^1} \leq 2 \norm{\hatx}$.
Hence $\abs{\varphi_2(\hatx)} = \norm{\hatx^2} \leq 2 \norm{\hatx}$.
So
$\norm{\psi_2} = \norm{\varphi_2} \leq 2$.
Hence $\norm{\hatx^2} = \abs{\psi_2(x)} \leq 2 \norm{x}$.
So $\norm{\hatx} \leq \norm{\hatx^1} + \norm{\hatx^2} \leq 3 \norm{x}$.

Let $T_1$ be defined by $T_1(x) = T(\hatx) + \barx$.
Hence $T_1\inverse(x) = T\inverse(\hatx) + \barx$.
Then for every $x \in E$,
$\norm{(T_1 - \rfs{Id}_E)(x)} =  \norm{(T - \rfs{Id}_F)(\hatx)} \leq
\norm{T - \rfs{Id}_F} \mcdot \norm{\hatx} \leq
3 \norm{T - \rfs{Id}_F} \mcdot \norm{x}$.
That is, $\norm{T_1 - \rfs{Id}_E} \leq 3 \norm{T - \rfs{Id}_F}$.
A similar computation shows that
$\norm{T_1\inverse - \rfs{Id}_E} \leq 3 \norm{T\inverse - \rfs{Id}_F}$.
So $d(T_1,\rfs{Id}_E) \leq 3 d(T,\rfs{Id}_F)$.

Also for every $x \in E$, $T_1(x) - x \in F \subseteq E_1$.
So $T_1(E_1) = E_1$, that is, $T_1 \in \bbL(E;E_1)$.
This proves Claim 3.
\smallskip

We return to the proof of the lemma.
Let $V \in \rfs{Nbr}^{\sboldbbL(E;E_1)}(\rfs{Id})$.
We may assume that $V = B^{\sboldbbL(E;E_1)}(\rfs{Id}_E,r)$.
Let
$n = n\kern-1.5pt
\left(B^{\sboldbbL(\itbfH^2)}(\rfs{Id}_{\itbfH^2},\frac{r}{9 \sqrt{2}}),
3 \sqrt{2}\right)$ be as assured by Claim 1.

Let $x_1,x_2 \in E_1$ be such that $\norm{x_1} = \norm{x_2} \neq 0$.
We show that there is $S \in V^n$ such that $S(x_1) = x_2$.
Let $F$ be a $2$-dimensional subspace of $E_1$ containing $x_1$
and $x_2$,\break
and $\fnn{T}{F}{\itbfH^2}$ be as assured by Claim~2.
Since $\norm{T} \leq 1$ and $\norm{T\inverse} \leq 3 \sqrt{2}$,\break
it follows that
$\frac{1}{3 \sqrt{2}} \leq
\frac{\raise 2.0pt\hbox{\scriptsize$\norm{T(x_1)}$}}
{\lower 1.1pt\hbox{\scriptsize$\dline{T(x_2)}$}} \leq 3 \sqrt{2}$.
Hence there is
$S_0 \in \left(B^{\sboldbbL
(\itbfH^2)}(\rfs{Id}_{\itbfH^2},\frac{r}{9 \sqrt{2}})\right)^n$
such that $S_0(T(x_1)) = T(x_2)$.
Let
$S_0 = S_{0,1} \scirc \ldots \scirc S_{0,n}$, where
$S_{0,i} \in
B^{\sboldbbL(\itbfH^2)}(\rfs{Id}_{\itbfH^2},\frac{r}{9 \sqrt{2}})$,
and define $S_1 = T\inverse S_0 T$ and $S_{1,i} = T\inverse S_{0,i} T$.
Then $S_1(x_1) = x_2$ and $S_1 = S_{1,1} \scirc \ldots \scirc S_{1,n}$.
Clearly,
$S_{1,i} - \rfs{Id}_F = T\inverse (S_{0,i} - \rfs{Id}_{\itbfH^2}) T$,
and hence
$$
\norm{S_{1,i} - \rfs{Id}_F} \leq
\norm{T\inverse} \mcdot \norm{(S_{0,i} - \rfs{Id}_{\itbfH^2})} \mcdot
\norm{T} \leq 3 \sqrt{2} \mcdot \norm{(S_{0,i} - \rfs{Id}_{\itbfH^2})}.
$$
The same inequality holds for $(S_{1,i})\inverse$. So
\vspace{0.9mm}
\newline\rule{4pt}{0pt}
{\thickmuskip=3mu \medmuskip=1mu \thinmuskip=1mu 
\renewcommand{\arraystretch}{1.5}
\addtolength{\arraycolsep}{-6pt}
$
\begin{array}{lll}
d(\rfs{Id}_F,S_{1,i})
&
=
&
\norm{S_{1,i} - \rfs{Id}_F} + \norm{(S_{1,i})\inverse - \rfs{Id}_F}
\\
&
\leq
&
3 \sqrt{2} \cdot \norm{(S_{0,i} - \rfs{Id}_{\itbfH^2})} +
3 \sqrt{2} \cdot \norm{((S_{0,i})\inverse - \rfs{Id}_{\itbfH^2})} =
3 \sqrt{2} \cdot d(S_{0,i},\rfs{Id}_{\itbfH^2}) < \dgfrac{r}{3}.
\end{array}
$
}
\renewcommand{\arraystretch}{1.0}
\addtolength{\arraycolsep}{6pt}

\noindent
By Claim 3, there are $S_{2,i} \in \bbL(E;E_1)$ extending
$S_{1,i}$ such that
$d(\rfs{Id}_E,S_{2,i}) \leq 3 \mcdot d(\rfs{Id}_F,S_{1,i})$.
Hence $S_{2,i} \in B^{\sboldbbL(E;E_1)}(\rfs{Id}_E,r)$, and so
$S \eqdf S_{2,1} \scirc \ldots \scirc S_{2,n}
\in \left(B^{\sboldbbL(E;E_1)}(\rfs{Id}_E,r)\right)^n = V^n$.

Let $W_1 \in \rfs{Nbr}^E(0)$,
and suppose that $W_1 \supseteq B^E(0,s)$.
Set $W_2 = B^E(0,\frac{s}{(1 + r)^n})$.
For any $L \in V$,
$\norm{L} < 1 + r$, hence for every $i \leq n$ and $L' \in V^i$,
$\norm{L'} < (1 + r)^i$.
So $L'(W_2) \subseteq W_1$ for every $i \leq n$ and $L' \in V^i$.
This proves that $n$ fulfills the requirements of the lemma.
\smallskip\rule{0pt}{0pt}\hfill\myqed

The following lemma is analogous to Proposition \ref{metr-bldr-p3.13}.

\begin{lemma}\label{metr-bldr-l3.31}
Let $E$ be a normed space with dimension $> 1$,
$F$ be a dense linear subspace of $E$ and $X \subseteq E$ be open.
Then $\lambda_{\bbA}^{E;F} \drest X$ is an affine-like
partial action.
\end{lemma}

\noindent
{\bf Proof }
At first we show that for every $x \in X \cap F$,
$\lambda_{\bbA}^{E;F} \drest X$ is affine-like at $x$.
Let $Y = X - x$.
The function $\chi$ from $\bbA(E;F) \cup X$ to $\bbA(E;F) \cup Y$
defined by: $\chi(u)~=~u~-~x$, \  $x \in X$;
and
$\chi(h) = h^{\srfs{tr}_{\kern-1.3pt-x}}$, \ $h \in \bbA(E;F)$,
is an isomorphism between the partial actions
$\lambda_{\bbA}^{E;F} \drest X$ and $\lambda_{\bbA}^{E;F} \drest Y$.
Also, $\chi \nrestriction X$ is an isometry.
So it suffices to prove that
$\lambda_{\bbA}^{E;F} \drest Y$ is
affine-like at~$0^E$. We rename $Y$ and call it $X$.

Denote $\bbA(E;F)$ by $\bbA$, $\bbT(E;F)$ by $\bbT$ and $\bbL(E;F)$
by $\bbL$.
Let $r,s  > 0$, $V = B^{\sboldbbA}(\rfs{Id},r)$,
$U = B^E(0,s)$, and assume that $U \subseteq X$.
We shall find $n = n(0^E,V,U)$, $U_0 = U_0(0^E,V,U)$ and
$D = D(0^E,V,U)$ which demonstrate that $\bbA$ is affine-like
at $0^E$.
Let $m = n(B^{\sboldbbL}(\rfs{Id},r))$ be as assured by
Lemma \ref{metr-bldr-l3.30}.
Set $t = \dgfrac{\min(r,s)}{2}$, $W_1 = B^E(0,t)$
and $W_2 = B^E\kern-3.0pt\left(0,\frac{t}{(1 + r)^m}\right)$,
and define $n = m + 2$, $U_0 = \half W_2$ and $D = U_0 \cap F$.

It is obvious that $D$ is metrically dense in $U_0$.
Let $x_1,y_1,x_2,y_2 \in D$ be such that
$\norm{x_1 - y_1} = \norm{x_2 - y_2}$.
For $\ell = 1,2$
let $g_{\ell} = \rfs{tr}^E_{\kern-1pt-x_{\ell}}$.
Since $\norm{x_1},\norm{x_2} < \frac{r}{2}$,
{\thickmuskip=2mu \medmuskip=1mu \thinmuskip=1mu 
it follows that $g_1,g_2 \in B^{\sboldbbT}(\rfs{Id},r)$.
Clearly, $g_{\ell}(x_{\ell}) = 0$,
and since $x_{\ell},y_{\ell} \in U_0 = \half W_2$, it follows that
$g_{\ell}(y_{\ell}) \in W_2$.
}
By Lemma \ref{metr-bldr-l3.30},
there are $h_1,\ldots,h_m \in B^{\sboldbbL}(\rfs{Id},r)$,
such that $h_1 \scirc \ldots \scirc h_m(g_1(y_1)) = g_2(y_2)$
and for every $i = 1,\ldots,m$,
$h_i \scirc \ldots \scirc h_m(g_1(y_1)) \in W_1$.
It follows that $g_2\inverse,h_1,\ldots,h_m,g_1$ are as required in
the definition of affine-likeness.

To show that $\bbA$ is affine-like at points that do not belong to $F$
we shall apply Proposition~\ref{metr-bldr-p3.29}.
Let $x \in X$.
Choose $r_0 > 0$ such that $B(x,2r_0) \subseteq X$
and set $C = B(x,r_0) \cap F$.
By the preceding argument, $\bbA$ is affine-like at every $y \in C$.
For $y \in C$ and $V \in \rfs{Nbr}^{\sboldbbA}(\rfs{Id})$
we define $V_y = V^{\srfs{tr}_y}$.

We next define functions
$\fnn{\barn}{\rfs{Nbr}^{\sboldbbA}(\rfs{Id}) \times (0,r_0)}{\bbN}$
and
$\fnn{\bars}{\rfs{Nbr}^{\sboldbbA}(\rfs{Id}) \times (0,r_0)}
{(0,\infty)}$
as needed in \ref{metr-bldr-p3.29}.
Let $V = B^{\sboldbbA}(\rfs{Id},r)$ and $s \in (0,r_0)$.
Set $m = n(B^{\sboldbbL}(\rfs{Id},r))$,
where $n(B^{\sboldbbL}(\rfs{Id},r))$ is as assured by
Lemma \ref{metr-bldr-l3.30}. Define $\barn(V,s) = m + 2$,
set $t = \dgfrac{\min(r,s)}{2}$
and define $\bars(V,s) = \frac{t}{2 (1 + r)^m}$.
It was proved in the preceding argument that
$$
n(0^E,V,B(0^E,s)) = \barn(V,s) \ \ \mbox{and} \ \ %
U_0(0^E,V,B(0^E,s)) = B(0^E,\bars(V,s)).
$$
Since $\rfs{tr}_y$ defines an isomorphism of partial actions,
which is an isometry on $X$, and  since $\rfs{tr}_y$ takes $0^E$ to $y$,
it can be concluded that
$$
n(y,V^{\srfs{tr}_y},B^E(y,s)) = \barn(V,s) \ \ \mbox{and} \ \ %
U_0(y,V^{\srfs{tr}_y},B^E(y,s)) = B(y,\bars(V,s)).
$$
We have shown that Clauses (i) and (ii) of
Proposition \ref{metr-bldr-p3.29} hold.

Recall that $x \in X$, $B(x,2r_0) \subseteq X$ and
$C = B(x,r_0) \cap F$.
Let $r > 0$ and $W = B^{\sboldbbA}(\rfs{Id},r)$.
We shall find $U_1$ and $V$ as required in Clause (iii) of
\ref{metr-bldr-p3.29}.
Let $\overbbA = \bbT(E) \mcdot \bbL(E;F)$.
Clearly, $\overbbA \leq \bbA(E)$. Also, $\bbA$ is dense in $\overbbA$.
Let $\overW = B^{\overbbA}(\rfs{Id},r)$, $g = \rfs{tr}_x$
and $\overV_1 = \overW^{g\inverse}$.
Note that $W = \overW \cap \bbA$.
Let $t > 0$ be such that
$(B^{\overbbA}(\rfs{Id},t))^3 \subseteq \overV_1$
and denote $\overV = B^{\overbbA}(\rfs{Id},t)$.
Define
$$
V = B^{\sboldbbA}(\rfs{Id},t) \ \ \mbox{and} \ \ U_1 = x + B^E(0,t).
$$
Let $y \in U_1$.
Then $\rfs{tr}_y \in g \mcdot \overV$ and so
$$
(\overV)^{\srfs{tr}_y} \subseteq g \cdot \overV \mcdot \overV
\cdot (\overV)\inverse \cdot g\inverse =
g \cdot (\overV)^3 \cdot g\inverse \subseteq (\overV_1)^g =
\overW.
$$
That is, $(\overV)^{\srfs{tr}_y} \subseteq \overW$.
If $y \in F$, then $V^{\srfs{tr}_y} \subseteq \bbA$.
In particular, if $y \in U_1 \cap F$,
then $V^{\srfs{tr}_y} \subseteq \overW \cap \bbA = W$.
This implies that Clause (iii) of Proposition \ref{metr-bldr-p3.29}
holds.
By Proposition \ref{metr-bldr-p3.29}, $\bbA$ is affine-like at $x$.
\rule{0pt}{1pt}\hfill\myqed

\begin{defn}\label{temp-d3.5}\label{metr-bldr-d3.32}
\begin{rm}

(a) Let $X$ be a metric space and $x \in X$.
We say that
$X$ has the {\it discrete path property at $x$ ($X$ is DPT at $x$)},
if the following holds.
There is $U \in \rfs{Nbr}(x)$ and $K \geq 1$ such that
($*$) for every $y,z \in U$ and $d \in (0, d(y,z))$
there are $n \in \bbN$ and $u_0,\ldots,u_n \in X$
such that $n \leq K \cdot \frac{d(y,z)}{d}$,
$d(y,u_0),\,d(u_n,z) < d$
and \,$d(u_{i - 1},u_i) = d$ for every $i = 1, \ldots, n$.

If $X$ is DPT at every $x \in X$, then $X$ is called a {\it DPT space}.
   \index{dpt@@DPT. A metric space $X$ is DPT at $x \in X$}
   \index{dpt@@DPT. A metric space is DPT}

(b) Let $X$ be a metric space and $x \in X$.
$X$ has {\it connectivity property 1 at $x$, ($X$ is CP1 at $x$)},
if for every $r > 0$ there is $r^* \in (0,r)$
such that for every $x' \in X$ and $r' > 0$:
if $B(x',r') \subseteq B(x,r^*)$ and $C$ is a connected component
of $B(x,r) - B(x',r')$,
then $C \cap (B(x,r) - B(x,r^*)) \neq \emptyset$.

If $X$ is CP1 at every $x \in X$, then $X$ is called a {\it CP1 space}.
   \index{cp1@@CP1. $X$ is CP1 at $x$}
   \index{cp1 space@@CP1 space}
\end{rm}
\end{defn}

\begin{prop}\label{metr-bldr-p3.33}
Let $X$ be an open subset of a normed space $E$.
Then $X$ is DPT and CP1.
\end{prop}

\noindent
{\bf Proof } Let $x \in X$ and $s > 0$
be such that $B^E(x,s) \subseteq X$.
At first we show that $X$ is DPT at $x$.
Let $y,z \in B^E(x,s)$ and $d \in (0,\norm{z - y})$.
The points
$u_i = y + i \mcdot \frac{d(z - y)}{\norm{z - y}}$,\ \ %
$i = 0,\ldots, \dgfrac{[\kern1pt\norm{z - y}\kern1pt}{d}]$
demonstrate the DPT-ness at $x$. So $K = 1$.

Let $r > 0$. Take $r^*$ to be any member of $(0,\min(r,s))$.
Let $x\fprime$ and $r\fprime < r^*$ be such that
$B^E(x\fprime,r\fprime) \subseteq B^E(x,r^*)$.
It is trivial that $B^E(x,s) - B^E(x\fprime,r\fprime)$ is connected.
So there is only one component $C$ of $B(x,r) - B^E(x\fprime,r\fprime)$
which intersects $B^E(x,s)$, and $C$ contains
$B^E(x,s) - B^E(x\fprime,r\fprime)$.
So $C$ intersects $B(x,r) - B^E(x^*,r)$.
Trivially, any connected component of
$B(x,r) - B^E(x\fprime,r\fprime)$ which is disjoint from $B^E(x,s)$
intersects $B(x,r) - B^E(x^*,r)$.
\smallskip\rule{10pt}{0pt}\hfill\myqed

Suppose that $X$ is an open subset of a normed space $E$, $G \leq H(X)$,
$\iso{\tau}{X}{Y}$ and\break
$G^{\tau} \subseteq H_{\itGamma}^{\srfs{LC}}(Y)$.
Loosely speaking we shall prove that if
\,($\dagger$)\, $\bbA(E) \nrestriction X \subseteq G$,
then $\tau$ is locally $\itGamma$-bicontinuous.
Obviously, ($\dagger$) is flawed because
$\bbA(E) \nrestriction X$ is not a set of homeomorphisms of $X$,
and hence not a subset of $G$.
The correct statement which replaces ($\dagger$) has the assumption that
$\lambda^{E}_{\sboldbbA}$ is compatible with $G$.
We do not know that this assumption suffices unless $E$ is a normed
space of the second category, or in particular, a Banach space.
Instead we assume that $\lambda^{E}_{\sboldbbA}$ is $G$-decayable,
and that $G$ is infinitely closed. These assumptions work for every
normed space $E$.

The following remains open.

\begin{question}\label{bddly-lip-bldr-q3.34}
\begin{rm}
Let $E,F$ be normed spaces of the first category,
$\iso{\tau}{E}{F}$ and $\itGamma$ be a countably generated
modulus of continuity.
Suppose that
$\bbA(E)^{\tau} \subseteq H_{\itGamma}^{\srfs{LC}}(F)$.
Are $\tau$ or $\tau\inverse$ or both locally $\itGamma$-continuous?
\end{rm}
\end{question}

The core fact that leads to the final result of \ref{metr-bldr-t3.27}
is stated in the following theorem.

\begin{theorem}\label{temp-t3.6}\label{metr-bldr-t3.35}
Assume the following facts.
\begin{itemize}
\addtolength{\parskip}{-11pt}
\addtolength{\itemsep}{06pt}
\item[\num{i}] 
$X$ and $Y$ are metric spaces, $x \in X$ and $\iso{\tau}{X}{Y}$.
Also, $X$ is DPT at $x$,
and $Y$ is DPT and CP1 at $\tau(x)$.
\item[\num{ii}] 
$G \leq H(X)$, $\lambda$ is a partial action of a topological
group $H$ on $X$, $\alpha \in \rfs{MBC}$,
$x \in \rfs{Fld}(\lambda)$, $x$ is a $\lambda$-limit-point,
$G$ is $\alpha$-infinitely-closed at $x$
and for some $N \in \rfs{Nbr}(x)$,
$\lambda$ is $(\alpha,G)$-decayable in $H_{\lambda}(x) \cap N$.
\item[\num{iii}] 
$\itGamma$ is a modulus of continuity,
and $\itGamma$ is $(\leq\kern-3pt\kappa)$-generated, where
$\kappa = \min(\setm{\kappa(x,V_{\lambda}(x))}{\break
V \in \rfs{Nbr}(e_H)})$.
\item[\num{iv}] 
There is $U \in \rfs{Nbr}(x)$ such that
for every $g \in G\sprt{\kern1ptU}$:
if $g$ is $\alpha \scirc \alpha$-bicon\-tinuous,
then $g^{\tau}$ is $\itGamma$-bicontinuous at $\tau(x)$.
\end{itemize}
\underline{Then} $\tau$ is $\itGamma$-bicontinuous at $x$.
\end{theorem}

We next introduce the notion of almost $\itGamma$-continuity.
The proof of \kern1ptTheorem \ref{metr-bldr-t3.35}
is broken into two claims.
The first one, Lemma \ref{metr-bldr-l3.37}(b), says that if
$G$ fulfills assumptions (i)\,-\,(iv) of \ref{metr-bldr-t3.35}
and $G^{\tau} \subseteq H_{\itGamma}^{\srfs{LC}}(Y)$,
then $\tau$ is locally almost $\itGamma$-continuous.
This part of the proof does not use the DPT-ness or the CP1-ness
of $X$ or $Y$.
The second claim is stated in Theorem \ref{metr-bldr-t3.40}.
It says that if $X$ and $Y$ are DPT and PC1 metric spaces,
and $\iso{\tau}{X}{Y}$ is
locally almost $\itGamma$-continuous,
then $\tau$ is locally $\itGamma$-bicontinuous.  

\begin{defn}\label{temp-d3.7}\label{metr-bldr-d3.36}
\begin{rm}
(a) Let $X,\,Y$ be metric spaces, $\alpha \in \rfs{MC}$,
$\itGamma$ be a modulus of continuity and
$\fnn{f}{X}{Y}$.
We say that $f$ is {\it almost $\alpha$-continuous},
if $f$ is continuous,
and for every \hbox{$x_1,y_1, x_2,y_2 \in X$:} if
$d(x_1,y_1) = d(x_2,y_2)$, then
$d(f(x_2),f(y_2)) \leq \alpha(d(f(x_1),f(y_1)))$.
The notion $f$ is {\it almost $\alpha$-continuous at $x$}
means that there is
$U \in \rfs{Nbr}(x)$ such that $f \nrestriction U$ is
almost $\alpha$-continuous.
We say that $f$ is {\it almost $\itGamma$-continuous at $x$},
if for some $\gamma \in \itGamma$,\break
$f$ is almost $\gamma$-continuous at~$x$,
and $f$ is said to be {\it locally almost $\itGamma$-continuous},
if for every $x \in X$, $f$ is almost $\itGamma$-continuous at $x$.

   \index{almost $\alpha$-continuous}
   \index{almost $\alpha$-continuous at $x$}
   \index{almost $\itGamma$-continuous at $x$}
   \index{locally almost $\itGamma$-continuous}

(b) If $\fnn{g}{A}{A}$, then $g^{\sscirc n}$
denotes $\,\overbrace{g \scirc \ldots \scirc g}^n$.
   \index{N@AAAA@@$g^{\sscirc n}$. This means $n$ times composition
          of $g$}
\hfill\proofend
\end{rm}
\end{defn}

The following lemma has also a variant in which $H$ is assumed to be
of the second category, but decayability is replaced by compatibility,
and infinite-closedness is dropped.

\begin{lemma}\label{temp-l3.8}\label{metr-bldr-l3.37}
\num{a}
Suppose that the following facts hold.

\kern5.0pt

\begin{list}{}
{\setlength{\leftmargin}{33pt}
\setlength{\labelsep}{08pt}
\setlength{\labelwidth}{20pt}
\setlength{\itemindent}{-00pt}
\addtolength{\topsep}{-11pt}
\addtolength{\parskip}{-12pt}
\addtolength{\itemsep}{-05pt}
}
\item[\num{i}] 
$X$ and $Y$ are metric spaces, $x \in X$ and $\iso{\tau}{X}{Y}$.
\item[\num{ii}] 
$\lambda$ is a partial action of a topological group $H$ on $X$,
$x \in \rfs{Fld}(\lambda)$ and $\lambda$ is affine-like at~$x$.
\item[\num{iii}] 
$\itGamma$ is a modulus of continuity and
$\gamma \in \itGamma$.
\item[\num{iv}] 
$T \in \rfs{Nbr}(x)$, $V \in \rfs{Nbr}(e_H)$,
$V \times T \subseteq \rfs{Dom}(\lambda)$
and for every $h \in V$, $(h_{\lambda})^{\tau} \nrestriction \tau(T)$
is $\gamma$-bicontinuous.
\end{list}
\underline{Then} $\tau$ is almost $\itGamma$-continuous at $x$.

\kern3pt

\num{b}
Suppose that the following facts hold.
\begin{itemize}
\addtolength{\parskip}{-11pt}
\addtolength{\itemsep}{06pt}
\item[\num{i}] 
$X$ and $Y$ are metric spaces, $x \in X$ and $\iso{\tau}{X}{Y}$.
\item[\num{ii}] 
$G \leq H(X)$, $\lambda$ is a partial action of a topological
group $H$ on $X$ and $\alpha \in \rfs{MBC}$.
Also, $x \in \rfs{Fld}(\lambda)$, $x$ is a $\lambda$-limit-point,
$G$ is $\alpha$-infinitely-closed at $x$
and for some $N \in \rfs{Nbr}(x)$,
$\lambda$ is $(\alpha,G)$-decayable in $H_{\lambda}(x) \cap N$.
\item[\num{iii}] 
$\itGamma$ is a $(\leq\kern-3pt\kappa)$-generated modulus of continuity,
where
$\kappa = \min(\setm{\kappa(x,V_{\lambda}(x))}{V \in \rfs{Nbr}(e_H)})$.
\item[\num{iv}] 
There is $U \in \rfs{Nbr}(x)$ such that
for every $g \in G\sprt{\kern1ptU}$:
if $g$ is $\alpha \scirc \alpha$-bicon\-tinuous,
then $g^{\tau}$ is $\itGamma$-bicontinuous at $\tau(x)$.
\end{itemize}
\underline{Then} $\tau$ is almost $\itGamma$-continuous at $x$.
\end{lemma}

\noindent
{\bf Proof }
(a) Let $n  = n(x,V,T)$, $U_0 = U_0(x,V,T)$ and $D = D(x,V,T)$
be as assured by the definition of affine-likeness,
(Definition \ref{metr-bldr-d3.28}(a)).
For $h \in H$ denote $h_{\lambda}$ by $\hath$.
Set $\beta = \gamma^{\sscirc n}$, so $\beta \in \itGamma$.
Suppose that $x_1,y_1,x_2,y_2 \in D$ and
$d(x_1,y_1) = d(x_2,y_2)$.
Choose $h_1,\ldots,h_n \in V$ as assured by the
definition affine-likeness,
and define $h = \bcirc_{i = 1}^n h_i$.
So $\hath(x_1) = x_2$, $\hath(y_1) = y_2$ and
$\hath_i \scirc \ldots \scirc \hath_n(\dbltn{x_1}{x_2}) \subseteq T$
for every $i = 1,\ldots,n$.
Also, for every $i = 1,\ldots,n$,
$(\hath_i)^{\tau} \nrestriction \tau(T)$ is $\gamma$-continuous.
Hence
$d(\tau(x_2),\tau(y_2)) =
d((\tau(x_1))^{\hath},(\tau(y_2))^{\hath}) \leq
\beta(d(\tau(x_2),\tau(y_2)))$.
We have shown that $\tau \nrestriction D$ is almost $\beta$-continuous.
Relying on the fact that $D$ is metrically dense in $U_0$ we conclude
that $\tau \nrestriction U_0$ is almost $\beta$-continuous.
\hbox{So $\tau$ is almost $\itGamma$-continuous at $x$.}

(b) By Lemma \ref{metr-bldr-l3.11},
there are $T \in \rfs{Nbr}(x)$, $V \in \rfs{Nbr}(e_H)$
and $\gamma \in \itGamma$ such that
for every $h \in V$:
$T \subseteq \rfs{Dom}(h_{\lambda})$ and
$(h_{\lambda})^{\tau} \nrestriction \tau(T)$ is $\gamma$-bicontinuous.
By Part (a), $\tau$ is almost $\itGamma$-continuous at $x$.
\smallskip\hfill\myqed

The next two propositions are ingredients in the proof
of Theorem \ref{metr-bldr-t3.40}.

\begin{prop}\label{temp-p3.9.1}\label{metr-bldr-p3.38}
Let $x$ belong to a metric space $X$,
and suppose that $X$ is DPT at $x$,
that $K$ and $U$ satisfy condition
$(*)$ of Definition \ref{metr-bldr-d3.32}(a) and that
$W \in \rfs{Nbr}(x)$.
Then there is $T \in \rfs{Nbr}(x)$
such that:
$(**)$ $T \subseteq W$, and
for every $y,z \in T$ and $d \in (0, d(y,z))$
there are $n \in \bbN$ and $u_0,\ldots, u_n \in W$
such that
$n \leq K \cdot \frac{d(x,y)}{d}$,
$d(x,u_0),\,d(u_n,y) < d$, and
$d(u_i,u_{i + 1}) = d$
for every $i = 0, \ldots, n - 1$.
\end{prop}

\noindent
{\bf Proof }
Let $s > 0$ be such that
$B(x,(2K + 3)s) \subseteq U \cap W$. We show that $T \eqdf B(x,s)$ is
as required. Let $y,z \in T$ and $d \in (0,d(y,z))$. Let $n \in \bbN$
and $u_0, \ldots, u_n$
be as assured in $(*)$ of \ref{metr-bldr-d3.32}(a).
Then for every $i = 1, \ldots, n$,
\kern0.2mm
\newline\indent
\renewcommand{\arraystretch}{1.5}
\addtolength{\arraycolsep}{-12pt}
$
\begin{array}{lll}
d(u_i,x)
&
\leq
&
d(u_i,u_0) + d(u_0,y) + d(y,x) < id + d + s \leq nd + d + s
\\&
\leq
&
K d(x,y) + 2s + s < K \mcdot 2s + 2s + s < (2K + 3) s.
\end{array}
$
\renewcommand{\arraystretch}{1.0}
\addtolength{\arraycolsep}{12pt}

\kern0.2mm

\noindent
So $u_i \in W$.\hfill\myqed

\begin{prop}\label{temp-p3.9.2}\label{metr-bldr-p3.39}
Let $X,Y$ be metric spaces and $\iso{\tau}{X}{Y}$.
Suppose that $x \in X$, $\tau$ is almost $\alpha$-continuous at $x$,
and $Y$ is CP1 at $\tau(x)$.
Then there is $U \in \rfs{Nbr}(x)$ such that
every $x_1,y_1,x_2,y_2 \in U$:
if $d(x_2,y_2) \leq d(x_1,y_1)$,
then $d(\tau(x_2),\tau(y_2)) \leq \alpha(d(\tau(x_1),\tau(y_1)))$.
\end{prop}

\noindent
{\bf Proof }
Let $T \in \rfs{Nbr}(x)$ be such that $\tau \nrestriction T$ is almost
$\alpha$-continuous,
and $s > 0$ be such that
$B(\tau(x),s) \subseteq \tau(T)$.
Choose $s^* \in (0,s)$ such that for every $y \in Y$ and $t > 0$:
if $B(y,t) \subseteq B(\tau(x),s^*)$, then every connected component of
$B(\tau(x),s) - B(y,t)$ intersects $B(\tau(x),s) - B(\tau(x),s^*)$.
Let $r^* > 0$ be such that
\begin{list}{}
{\setlength{\leftmargin}{33pt}
\setlength{\labelsep}{08pt}
\setlength{\labelwidth}{20pt}
\setlength{\itemindent}{-00pt}
\addtolength{\topsep}{-10pt}
\addtolength{\parskip}{-04pt}
}
\item[(i)] 
$\tau(B(x,r^*)) \subseteq B(\tau(x),s^*)$,
\end{list}
\kern2pt
and let $r \in (0,\dgfrac{r^*}{3})$ be such that $U \eqdf B(x,r)$
satisfies the following condition:
\begin{list}{}
{\setlength{\leftmargin}{36pt}
\setlength{\labelsep}{08pt}
\setlength{\labelwidth}{27pt}
\setlength{\itemindent}{-00pt}
\addtolength{\topsep}{-10pt}
\addtolength{\parskip}{-04pt}
}
\item[(ii)] 
$\rfs{diam}(\tau(U)) + \alpha(\rfs{diam}(\tau(U))) < s^*$.
\end{list}
\kern2pt
We show that $U$ is as required.
Let $x_1,y_1,x_2,y_2 \in U$ and  $d(x_2,y_2) \leq d(x_1,y_1)$.
If $d(x_2,y_2) = d(x_1,y_1)$, then by the choice of $T$,
and since $U \subseteq T$,
$d(\tau(x_2),\tau(y_2)) \leq \alpha(d(\tau(x_1),\tau(y_1)))$.
Suppose next that $d(x_2,y_2) < d(x_1,y_1)$.
Let $r_1 = d(x_1,y_1)$,
and set $s_1 = \alpha(d(\tau(x_1),\tau(y_1)))$.
By the almost $\alpha$-continuity of $\tau \nrestriction T$,
\begin{list}{}
{\setlength{\leftmargin}{36pt}
\setlength{\labelsep}{08pt}
\setlength{\labelwidth}{27pt}
\setlength{\itemindent}{-00pt}
\addtolength{\topsep}{-10pt}
\addtolength{\parskip}{-04pt}
}
\item[(iii)] 
$\tau(S(x_2,r_1)) \subseteq B(\tau(x_2),s_1 + \varepsilon)$
for every $\varepsilon > 0$.
\end{list}
\kern2pt
Since $r < \dgfrac{r^*}{3}$, \ $d(x_2,x) < r$
and $r_1 = d(x_1,y_1) < 2r$, we have
\begin{list}{}
{\setlength{\leftmargin}{36pt}
\setlength{\labelsep}{08pt}
\setlength{\labelwidth}{27pt}
\setlength{\itemindent}{-00pt}
\addtolength{\topsep}{-10pt}
\addtolength{\parskip}{-04pt}
}
\item[(iv)] 
$B(x_2,r_1) \subseteq B(x,r^*)$.
\end{list}
\kern2pt
The following three facts:
$d(\tau(x),\tau(x_2)) \leq \rfs{diam}(\tau(U))$,
$s_1 \leq \alpha(\rfs{diam}(\tau(U)))$ and\break
$\rfs{diam}(\tau(U)) + \alpha(\rfs{diam}(\tau(U))) < s^*$,
imply that
\begin{list}{}
{\setlength{\leftmargin}{36pt}
\setlength{\labelsep}{08pt}
\setlength{\labelwidth}{27pt}
\setlength{\itemindent}{-00pt}
\addtolength{\topsep}{-10pt}
\addtolength{\parskip}{-04pt}
}
\item[(v)] 
for all sufficiently small $\varepsilon$'s,
$B(\tau(x_2),s_1 + \varepsilon) \subseteq B(\tau(x),s^*)$.
\end{list}
\kern2pt
Let $z \in Y - B(\tau(x_2),s_1 + \varepsilon)$. We show that
$\tau\inverse(z) \not\in B(x_2,r_1)$.
If $z \not\in B(\tau(x),s)$,
then $\tau\inverse(z) \not\in B(x,r^*) \supseteq B(x_2,r_1)$.
Suppose that $z \in B(\tau(x),s)$,
and let $C$ be the connected component of $z$ in
\hbox{$B(\tau(x),s) - B(\tau(x_2),s_1 + \varepsilon)$. Hence}
\begin{list}{}
{\setlength{\leftmargin}{36pt}
\setlength{\labelsep}{08pt}
\setlength{\labelwidth}{27pt}
\setlength{\itemindent}{-00pt}
\addtolength{\topsep}{-10pt}
\addtolength{\parskip}{-04pt}
}
\item[(vi)] 
$C \cap (B(\tau(x),s) - B(\tau(x),s^*)) \neq \emptyset$.
\end{list}
\kern2pt
Since $\tau(B(x_2,r_1)) \subseteq \tau(B(x,r^*)) \subseteq
B(\tau(x),s^*)$,
it follows that
\begin{list}{}
{\setlength{\leftmargin}{36pt}
\setlength{\labelsep}{08pt}
\setlength{\labelwidth}{27pt}
\setlength{\itemindent}{-00pt}
\addtolength{\topsep}{-10pt}
\addtolength{\parskip}{-04pt}
}
\item[(vii)] 
$\tau^{-1}(C) \cap (X - B(x_2,r_1)) \neq \emptyset$.
\end{list}
\kern2pt
From the facts:
$\tau(S(x_2,r_1)) \subseteq B(\tau(x_2),s_1 + \varepsilon)$
and $C \cap B(\tau(x_2),s_1 + \varepsilon) = \emptyset$,
\hbox{we conclude that}

\begin{list}{}
{\setlength{\leftmargin}{36pt}
\setlength{\labelsep}{08pt}
\setlength{\labelwidth}{27pt}
\setlength{\itemindent}{-00pt}
\addtolength{\topsep}{-10pt}
\addtolength{\parskip}{-04pt}
}
\item[(viii)] 
$\tau^{-1}(C) \cap S(x_2,r_1) = \emptyset$.
\end{list}
\kern2pt
The connectedness of $C$ and hence of $\tau\inverse(C)$
and facts (vii) and (viii) imply that
\begin{list}{}
{\setlength{\leftmargin}{36pt}
\setlength{\labelsep}{08pt}
\setlength{\labelwidth}{27pt}
\setlength{\itemindent}{-00pt}
\addtolength{\topsep}{-10pt}
\addtolength{\parskip}{-04pt}
}
\item[(ix)] 
$\tau^{-1}(C) \cap B(x_2,r) = \emptyset$.
\end{list}
\kern2pt
This implies that $\tau^{-1}(z) \not\in B(x_2,r_1)$.
Since the above argument holds for all sufficiently small
$\varepsilon$'s, it follows that
for every $z \in Y$: if $z \not\in \overB(\tau(x_2),s_1)$,
then $\tau\inverse(z) \not \in B(x_2,r_1)$.
But $y_2 \in B(x_2,r_1)$, so $\tau(y_2) \in \overB(\tau(x_2),s_1)$.
That is, $d(\tau(x_2),\tau(y_2)) \leq s_1 = \alpha(d(x_1,y_1))$.
\rule{0pt}{0pt}\hfill\myqed

\begin{theorem}\label{temp-l3.9}\label{metr-bldr-t3.40}
Let $X$ and $Y$ be metric spaces, $x_0 \in X$,
$\iso{\tau}{X}{Y}$ and $\alpha \in \rfs{MBC}$.
Suppose that $X$ is DPT at $x_0$, $Y$ is DPT and CP1 at $\tau(x_0)$,
and $\tau$ is almost $\alpha$-continuous at $x_0$.
Then there is $M > 0$ such that $\tau$ is $M \mcdot \alpha$-bicontinuous
at $x_0$.
\end{theorem}

\noindent
{\bf Proof }
{\thickmuskip=2mu \medmuskip=1mu \thinmuskip=1mu 
We first show that there is some $M > 0$ such that $\tau^{-1}$ is
$M \cdot\alpha$-con\-tinuous at $\tau(x_0)$.}
By Proposition \ref{metr-bldr-p3.39}, by the fact that $Y$ is CP1,
and since $\tau$ is almost $\alpha$-continuous at $x_0$,
there is $U \in \rfs{Nbr}(x_0)$,
such that for every $x_1,y_1,x_2,y_2 \in U$:
if $d(x_2,y_2) \leq d(x_1,y_1)$,
then $d(\tau(x_2),\tau(y_2)) \leq \alpha(d(\tau(x_1),\tau(y_1)))$.
It is given that $X$ is DPT at $x_0$,
so there are $W \in \rfs{Nbr}(x_0)$ and $K \geq 1$ such that
$W \subseteq U$, and $W,\,K$
satisfy condition ($*$) of Definition \ref{metr-bldr-d3.32}(a).
Let $V \subseteq W$ be an open neighborhood of $x_0$
satisfying condition ($**$) of Proposition~\ref{metr-bldr-p3.38}.
Fix any distinct $x_1,y_1 \in V$ and set $d_1 = d(x_1,y_1)$,
$e_1 = d(\tau(x_1),\tau(y_1))$,
$V_1 = B(x_0,\dgfrac{d_1}{2}) \cap V$ and $V_2 = \tau(V_1)$.

We show that $\tau^{\inverse} \nrestriction V_2$ is
$\frac{d_1}{e_1} \mcdot (K + 2) \cdot \alpha$-continuous.
Let $u,v \in V_2$ be distinct
and set
$d = d(\tau^{-1}(u),\,\tau^{-1}(v))$.
Since $\tau^{-1}(u),\,\tau^{-1}(v) \in V_1$, \,
$d < d_1 = d(x_1,y_1)$.
So there are $n \leq K \mcdot \frac{d(x_1,y_1)}{d}$
and $z_0, \ldots, z_n \in U$ such that
$d(x_1,z_0),\ d(z_n,y_1) < d$ and $d(z_i,z_{i + 1}) = d$\break
for all $i = 0, \ldots, n - 1$.
By the choice of $U$,
$$
d(\tau(x_1),\tau(z_0)),\,d(\tau(z_n),\tau(y_1)),\,
d(\tau(z_i),\tau(z_{i + 1})) \leq
\alpha(d(\tau\tau\inverse(u),\tau\tau\inverse(v))) = \alpha(d(u,v)).
$$
Hence  
\vspace{1.5mm}
\newline
\rule{40pt}{0pt}
\renewcommand{\arraystretch}{1.5}
\addtolength{\arraycolsep}{-18pt}
$
\begin{array}{ll}
&
d(\tau(x_1),\tau(y_1)) \leq
d(\tau(x_1),\tau(z_0)) +
\sum\limits_{i = 0}^{n - 1} d(\tau(z_i),\tau(z_{i + 1})) +
d(\tau(z_n),\tau(y_1))
\\
\leq
\rule{5pt}{0pt}
&
(n + 2) \alpha(d(u,v)) \leq
\left(K\frac{d(x_1,y_1)}{d(\tau^{-1}(u),\tau^{-1}(v))} + 2\right)
\alpha(d(u,v)).
\\
\end{array}
$
\renewcommand{\arraystretch}{1.0}
\addtolength{\arraycolsep}{18pt}
\vspace{2.1mm}\newline
It follows from the above inequality that
\vspace{1.5mm}
\newline
\rule{40pt}{0pt}
\renewcommand{\arraystretch}{1.6}
\addtolength{\arraycolsep}{-18pt}
$
\begin{array}{ll}
&
d(\tau\minus(u),\tau\minus(v)) \leq \frac{K d(x_1,y_1) +
2d(\tau\inverse(u),\tau\inverse(v))}{d(\tau(x_1),\tau(y_1))}
\alpha(d(u,v))
\\
\leq
\rule{5pt}{0pt}
&
\frac{Kd_1 + 2d_1}{e_1} \cdot \alpha(d(u,v))  =
\frac{d_1}{e_1} \mcdot (K + 2) \alpha(d(u,v)).
\vspace{1.7mm}
\end{array}
$
\renewcommand{\arraystretch}{1.0}
\addtolength{\arraycolsep}{18pt}
\newline
So $\tau\inverse \nrestriction V_2$ is
$\frac{d_1}{e_1} \mcdot (K + 2) \mcdot \alpha$-continuous,
and hence $\tau\inverse$ is locally $\itGamma$-continuous.

Note that in the above proof we only used the facts that $X$ is
DPT at $x_0$, and that $Y$ is CP1 at $\tau(x_0)$.

We now turn to the proof that there is $M > 0$
such that $\tau$ is $M \mcdot \alpha$-continuous at $x_0$.
In this part we use the facts that $Y$ is DPT and CP1 at $\tau(x_0)$.
Let $U_1 \in \rfs{Nbr}(x_0)$ and $K \geq 1$ be such that
$\tau(U_1)$ and $K$ satisfy condition~$(*)$ of \ref{metr-bldr-d3.32}(a)
applied to $\tau(x_0)$.
By Proposition~\ref{metr-bldr-p3.39}, there is $U_0 \in \rfs{Nbr}(x_0)$
such that $U_0 \subseteq U_1$, and
\begin{itemize}
\addtolength{\parskip}{-11pt}
\addtolength{\itemsep}{06pt}
\item[(1)] 
for every $x_1,y_1,x_2,y_2 \in U_0$: if $d(x_2,y_2) \leq d(x_1,y_1)$,
then
$d(\tau(x_2),\tau(y_2)) \leq\break
\alpha(d(\tau(x_1),\tau(y_1)))$.
\vspace{-05.7pt}
\end{itemize}
\kern2pt
We apply Proposition \ref{metr-bldr-p3.38}
to $\tau(x_0)$ and $\tau(U_0)$,
and obtain $T \subseteq Y$ satisfying condition ($**$)
of Propostion \ref{metr-bldr-p3.38}.
Let $U = \tau^{-1}(T)$.
We may assume that
\begin{itemize}
\addtolength{\parskip}{-11pt}
\addtolength{\itemsep}{06pt}
\item[(2)] 
$K \geq 2$.
\vspace{-05.7pt}
\end{itemize}
Let $x,y \in U$ be distinct.
Set $N = \frac{4K d(\tau(x),\tau(y))}{d(x,y)}$
and $M = \max(1,N)$.
\rule{0pt}{13pt}\kern-0pt
We show that if $x',y' \in U$ and $d(x',y') < \frac{d(x,y)}{4K}$,
then
$d(\tau(x'),\tau(y')) \leq M \mcdot \alpha(d(x',y'))$.
Obviously, this implies that
$\tau \nrestriction (B(x_0,\frac{d(x,y)}{8K}) \cap U)$
is $M \mcdot \alpha$-continuous.

Let $x',y' \in U$ be such that $d(x',y') < \frac{d(x,y)}{4K}$
and $n = \left[\frac{d(x,y)}{K d(x',y')}\right] - 2$.
Hence $n \geq 2$.\break
Let $d = \dgfrac{d(\tau(x),\tau(y))}{n}$.
So there are $m \leq Kn$ and $z_0, \ldots, z_m \in \tau(U_0)$
such that\break
{\thickmuskip=4mu \medmuskip=2.5mu \thinmuskip=1.5mu 
$d(\tau(x),z_0),\ d(z_m,\tau(y)) < d$
and $d(z_{i - 1},z_i) = d$, $i = 1,\ldots,m$.
Let $x_i = \tau\inverse(z_i)$.
}
Denote $x$ by $x_{-1}$ and $y$ by $x_{m + 1}$.
For $\ell \in \dbltn{-1}{m + 1}$ let $z_{\ell} = \tau(x_{\ell})$.
The number of $x_j$'s is $m + 3$. So the number of distances
between consecutive $x_j$'s is $m + 2$.
Hence for some $i \in \fsetn{0}{m + 1}$,\break

\kern-15.7pt
\begin{itemize}
\addtolength{\parskip}{-11pt}
\addtolength{\itemsep}{06pt}
\item[(3)] 
$d(x_{i - 1},x_i) \geq \frac{d(x,y)}{m + 2}$.
\vspace{-05.7pt}
\end{itemize}
It follows from (3) and (2) that
$$
\hbox{
$d(x_{i - 1},x_i) \geq \frac{d(x,y)}{m + 2} \geq 
\frac{d(x,y)}{K (\left[\frac{d(x,y)}{K d(x',y')}\right] - 2) + 2} \geq
\frac{d(x,y)}{K (\frac{d(x,y)}{K d(x',y')} + 1 - 2) + K} \geq
\frac{d(x,y)}{K \cdot \frac{d(x,y)}{K d(x',y')}} = d(x',y').$}
$$
That is,
\begin{itemize}
\addtolength{\parskip}{-11pt}
\addtolength{\itemsep}{06pt}
\item[(4)] 
$d(x',y') \leq d(x_{i - 1},x_i)$.
\vspace{-05.7pt}
\end{itemize}

Since the $z_i$'s belong to $\tau(U_0)$, the $x_i$'s belong to
$U_0$. This is also true for $x_{-1} = x$ and $x_{m + 1} = y$ because
they belong to $U \subseteq U_0$.
By (1) and (4),
\begin{itemize}
\addtolength{\parskip}{-11pt}
\addtolength{\itemsep}{06pt}
\item[(5)] 
$d(\tau(x'),\tau(y')) \leq \alpha(d(z_{i - 1},z_i)) = \alpha(d)$.
\vspace{-05.7pt}
\end{itemize}
Also,

\kern1mm
\renewcommand{\arraystretch}{1.8}
\addtolength{\arraycolsep}{-12pt}
$
\begin{array}{lll}
d
&
=
&
\frac{1}{n} d(\tau(x),\tau(y)) =
\frac{1}{\left[\frac{d(x,y)}{K d(x',y')}\right] - 2} d(\tau(x),\tau(y))
\leq
\frac{1}{\frac{d(x,y)}{K d(x',y')} - 1 - 2} d(\tau(x),\tau(y)) 
\\
&
=
&
\frac{K d(x',y')}{d(x,y) -3K d(x',y')} d(\tau(x),\tau(y)) \leq
\frac{K d(x',y')}{d(x,y) -3K \frac{d(x,y)}{4K}} d(\tau(x),\tau(y))
\\
&=
&
\frac{4K d(\tau(x),\tau(y))}{d(x,y)} d(x',y') = N d(x',y').
\end{array}
\renewcommand{\arraystretch}{1.0}
\addtolength{\arraycolsep}{12pt}
$\vspace{4pt}
\newline
By (5), by the fact $M \geq 1,N$ and by the concavity of $\alpha$,
$$
d(\tau(x'),\tau(y')) \leq \alpha(d) \leq \alpha(N d(x',y')) \leq
\alpha(M \cdot d(x',y')) \leq M \cdot \alpha(d(x',y')).
$$
We have thus shown that
$\tau \nrestriction (B(x_0,\frac{d(x,y)}{8K}) \cap U)$
is $M \mcdot \alpha$-continuous.
\medskip\hfill\myqed

\noindent
{\bf Proof of Theorem \ref{metr-bldr-t3.35} }
Let $X$, $x$, $Y$, $\tau$, $\itGamma$ etc.\ fulfill the premises of
\ref{metr-bldr-t3.35}. Then the assumptions of
Lemma \ref{metr-bldr-l3.37}(b) are satisfied.
So $\tau$ is almost $\itGamma$-continuous at $x$.
By Theorem~\ref{metr-bldr-t3.40},
$\tau$ is $\itGamma$-bicontinuous at $x$.
\smallskip\hfill\myqed

\noindent
{\bf Proof of Theorem \ref{metr-bldr-t3.27} }
Let $\fourtpl{E}{X}{\calS}{\calE}$ be a subspace choice system,
$Y$ be an open subset of a normed space $F$,
$\itGamma$ be a $(\leq \kappa(E))$-generated modulus of continuity
and $\iso{\tau}{X}{Y}$.
Suppose that
$(\rfs{LIP}(X,\calS,\calE))^{\tau} \subseteq
H_{\itGamma}^{\srfs{LC}}(Y)$,
and we prove that $\tau$ is locally $\itGamma$-bicontinuous.

For $x \in X$ choose $S \in \calS$ such that $x \in S$ and denote
$E_S$ by $D$.
We wish to apply Theorem~\ref{metr-bldr-t3.35}
to $G = \rfs{LIP}(X,\calS,\calE)$, $H = \bbA(E;D)$,
$\alpha(t) = 15t$ and
$\lambda = \lambda^{E;D}_{\sboldbbA} \drest S$,
so we check that Clauses (i)-(iv) of Theorem \ref{metr-bldr-t3.35} hold.

In Clause (i) we have to check that $X$ is DPT at $x$ and that
$Y$ is DPT and CP1 at $\tau(x)$, and this was proved in
Proposition \ref{metr-bldr-p3.33}.
In Clause (ii) we have to check: (1) $x$ is a $\lambda$-limit-point;
(2) $G$ is $\alpha$-infinitely-closed at~$x$;
(3)  for some $N \in \rfs{Nbr}(x)$, $\lambda$ is
$(\alpha,G)$-decayable in $N \cap H_{\lambda}(x)$.

(1) Obviously, for every $V \in \rfs{Nbr}^H(\rfs{Id})$,
$V_{\lambda}(x)$ contains a ball with center at $x$. So
$x$ is a $\lambda$-limit-point.

(2) Suppose that $\beta \in \rfs{MC}$,
$K \subseteq H_{\sngltn{\beta}}(Z)$ and for every
distinct $k_1,k_2 \in K$,
$\rfs{supp}(k_1) \cap \rfs{supp}(k_2) = \emptyset$.
Then $k \eqdf \bcirc K \in H(Z)$,
and $k$ is $\beta \scirc \beta$-bicontinuous.
Also, if $M \subseteq Z$, and $k'(M) = M$ for every $k' \in K$,
then $k(M) = M$.
These observations imply that $G$ is $\alpha$-infinitely-closed.

(3) The $(\alpha,G)$-decayabilty of $\lambda$ at every point of $S$
was proved in Lemma~\ref{metr-bldr-l3.8}.

Clause (iii) is given,
and Clause (iv) holds, since it is given that
$G^{\tau} \subseteq H_{\itGamma}^{\srfs{LC}}(Y)$.

By Theorem \ref{metr-bldr-t3.35}, $\tau$ is $\itGamma$-bicontinuous
at $x$.
We have shown that $\tau$ is locally $\itGamma$-bicon\-tinuous.
\hfill\myqed

In Theorem \ref{metr-bldr-t3.26} we have presented an alternative
argument for showing the local\break
$\itGamma$-continuity of $\tau\inverse$.
This method used Baire Caregory Theorem, but did not require the
assumptions of decayability of $\lambda$ and the infinite-closedness
of $G$.
The same alternative argument can be employed in the context of
affine-like partial actions.
It is presented in the following theorem.

\begin{theorem}\label{metr-bldr-t3.41}
Assume that the following facts hold.
\begin{list}{}
{\setlength{\leftmargin}{39pt}
\setlength{\labelsep}{08pt}
\setlength{\labelwidth}{29pt}
\setlength{\itemindent}{-00pt}
\addtolength{\topsep}{-08pt}
\addtolength{\parskip}{-12pt}
\addtolength{\itemsep}{-05pt}
}
\item[\num{i}] 
$X$ is a metric space, $G \leq H(X)$,
$H$ is a topological group and $H$ is of the second category,
$\lambda$ is a partial action of $H$ on $X$
and $x \in \rfs{Fld}(\lambda)$.
\item[\num{ii}] 
$\lambda$ is compatible with $G$ at $x$.
\item[\num{iii}] 
$\lambda$ is affine-like at $x$.
\item[\num{iv}] 
$\itGamma$ is a countably generated modulus of continuity.
\item[\num{v}] 
$Y$ is metric space and $\iso{\tau}{X}{Y}$.
\item[\num{vi}] 
For every $g \in G$,
$g^{\tau}$ is $\itGamma$-bicontinuous at $\tau(x)$.
\item[\num{vii}] 
$X$ is DPT at $x$ and $Y$ is DPT and CP1 at $\tau(x)$.
\end{list}
\underline{Then} $\tau$ is $\itGamma$-bicontinuous at $x$.
\end{theorem}

\noindent
{\bf Proof }
For $h \in H$ write $h_{\lambda} = \hath$.
The assumptions of Lemma \ref{metr-bldr-l3.22} hold, so
there are $T \in \rfs{Nbr}(x)$, a nonempty open subset $V \subseteq H$
and $\gamma \in \itGamma$ such that
$V \times T \subseteq \rfs{Dom}(\lambda)$ and
$\hath^{\tau} \nrestriction \tau(T)$ is $\gamma$-bicon\-tinuous
for every $h \in V$.
Note that $(\hath\inverse)^{\tau} \nrestriction \tau(\hath(T))$
is $\gamma$-bicon\-tinuous for every $h \in V$.

Let $h_0 \in V$.
There are $S \in \rfs{Nbr}(x)$ and $V_1 \in \rfs{Nbr}(h_0)$ such that
$V_1 \subseteq V$, $S \subseteq T$
and $\lambda(V_1 \times S) \subseteq \hath_0(T)$.
Set $W = h_0\inverse \mcdot V_1$.
Clearly, $W \in \rfs{Nbr}(e_H)$
and $W \times S \subseteq \rfs{Dom}(\lambda)$.
Let $h \in W$. So for some $h_1 \in V_1$ we have
$h = h_0\inverse \mcdot h_1$.
From the facts $h_1 \in V_1 \subseteq V$ and $S \subseteq T$,
it follows that
(1) ($\hath_1)^{\tau}\nrestriction \tau(S)$ is $\gamma$-bicontinuous,
and since $\hath_1(S) \subseteq \hath_0(T)$
and $h_1\inverse \in V\inverse$, we conclude that
(2) $(\hath_0\inverse)^{\tau} \nrestriction \tau(\hath_1(S))$
is $\gamma$-bicontinuous.
(1) and (2) imply that $\hath^{\tau} \nrestriction \tau(S)$ is
$\gamma \scirc \gamma$-bicontinuous.

{\thickmuskip=2mu \medmuskip=1mu \thinmuskip=1mu 
We have shown that there are $W \in \rfs{Nbr}(e_H)$
and $S \in \rfs{Nbr}(x)$ such that
$W \times S \subseteq \rfs{Dom}(\lambda)$,
}
and for every $h \in W$, $\hath^{\tau} \nrestriction \tau(W)$ is
$\gamma \scirc \gamma$-bicontinuous.
By Lemma \ref{metr-bldr-l3.37}(a), $\tau$ is almost
$\itGamma$-continuous at $x$,
and by Theorem \ref{metr-bldr-t3.40},
$\tau$ is $\itGamma$-bicontinuous at $x$.
\hfill\myqed

\subsection{Summary and questions.}
\label{ss3.5}
\label{ss3.5-summary-and-questions}

The following final theorem combines the results of
Chapters \ref{s2} and \ref{s3}.
Note that Part (a) of \ref{metr-bldr-t3.42} is not a special case of
Part (b).

\begin{theorem}\label{metr-bldr-t3.42}
\num{a} Let $X,Y$ be open subsets of the normed spaces
$E$ and $F$ respectively,
$\itGamma,\itDelta$ be moduli of continuity
and
$\iso{\varphi}{H_{\itGamma}^{\srfs{LC}}(X)}
{H_{\itDelta}^{\srfs{LC}}(Y)}$.
Suppose that $\itGamma$ is $(\leq\kern-1pt \kappa(E))$-generated.
\underline{Then} $\itGamma = \itDelta$,
there is $\iso{\tau}{X}{Y}$ such that
$\varphi(h) = h^{\tau}$ for every
$h \in H_{\itGamma}^{\srfs{LC}}(X)$,
and $\tau$ is locally $\itGamma$-bicontinuous.

\num{b}
Let $\fourtpl{E}{X}{\calS}{\calE}$ and $\fourtpl{F}{Y}{\calT}{\calF}$
be subspace choice systems, $\itGamma,\itDelta$ be moduli of continuity
and
$\iso{\varphi}
{H_{\itGamma}^{\srfs{LC}}(X;\calS,\calE)}
{H_{\itDelta}^{\srfs{LC}}(Y;\calT,\calF)}$.
Suppose that $\itGamma$ and $\itDelta$ are
$(\leq\kern-1pt \kappa(E))$-generated.
{\thickmuskip=4mu \medmuskip=3mu \thinmuskip=2mu 
\underline{Then} $\itGamma = \itDelta$,
there is $\iso{\tau}{X}{Y}$ such that
$\varphi(h) = h^{\tau}$ for every
$h \in H_{\itGamma}^{\srfs{LC}}(X;\calS,\calE)$,
and $\tau$ is locally $\itGamma$-bicontinuous.
}
\end{theorem}

\noindent
{\bf Proof } (a) 
$\rfs{LIP}^{\srfs{LC}}(X) \subseteq
H_{\itGamma}^{\srfs{LC}}(X) \subseteq H(X)$
and the same holds for $Y$.
So by Theorem \ref{t2.4}(a) there is $\iso{\tau}{X}{Y}$
such that $\tau$ induces $\varphi$.
Hence
$(H_{\itDelta}^{\srfs{LC}}
(Y))^{\tau\inverse} =
H_{\itGamma}^{\srfs{LC}}(X)$.
In particular,
$(\rfs{LIP}(Y))^{\tau\inverse} \subseteq
H_{\itGamma}^{\srfs{LC}}(X)$.
Since $X \cong Y$, $\kappa(F) = \kappa(E)$.
So $\itGamma$ is
$(\leq\kern-1pt \kappa(F))$-generated.
By Theorem \ref{metr-bldr-t3.27}, $\tau\inverse$ is locally
$\itGamma$-bicontinuous.
That is, $\tau$ is locally $\itGamma$-bicontinuous.
Hence
$H_{\itDelta}^{\srfs{LC}}(Y) =
(H_{\itGamma}^{\srfs{LC}}(X))^{\tau} \subseteq
H_{\itGamma}^{\srfs{LC}}(Y)$\kern-0.8pt.
It is easy to see that if
$\alpha \in \itDelta - \itGamma$, then there
is $h \in H(Y)$
such that $h$ is $\alpha$-bicontinuous and $h$ is not locally
$\itGamma$-continuous. This implies that
$\itDelta \subseteq \itGamma$.

Suppose by contradiction that $\itGamma - \itDelta \neq \emptyset$.
It is easy to see that there is
$h \in H_{\itGamma}^{\srfs{LC}}(Y) - H_{\itDelta}^{\srfs{LC}}(Y)$.
So $g \eqdf h^{\tau\inverse} \in  H_{\itGamma}^{\srfs{LC}}(X)$.
However, $g^{\tau} = h \not\in H_{\itDelta}^{\srfs{LC}}(Y)$.
A contradiction. So $\itGamma = \itDelta$.

(b) $\rfs{LIP}^{\srfs{LC}}(X;\calS,\calE) \subseteq
H_{\itGamma}^{\srfs{LC}}(X;\calS,\calE) \subseteq H(X)$
and the same holds for $Y$.
So by Theorem \ref{t2.4}(b) there is $\iso{\tau}{X}{Y}$
such that $\tau$ induces $\varphi$.
Hence
$(H_{\itGamma}^{\srfs{LC}}
(X;\calS,\calE))^{\tau} =
H_{\itGamma}^{\srfs{LC}}(Y;\calT,\calF)$.
In particular,
$(\rfs{LIP}(X;\calS,\calE))^{\tau} \subseteq
H_{\itDelta}^{\srfs{LC}}(Y)$ and
$(\rfs{LIP}(Y;\calT,\calF))^{\tau\inverse} \subseteq
H_{\itGamma}^{\srfs{LC}}(X)$.
By Theorem~\ref{metr-bldr-t3.19}(b),
$\itGamma = \itDelta$ and $\tau$ is locally $\itGamma$-bicontinuous.
\smallskip\hfill\myqed

The technical and abstract formulation of
Theorems
\ref{metr-bldr-t3.15},
\ref{metr-bldr-t3.26},
\ref{metr-bldr-t3.35} and
\ref{metr-bldr-t3.41}
hinders the understanding of their essence.
The above theorems are better understood through their application
to normed spaces, as stated in the following corollary.

\begin{cor}\label{metr-bldr-c3.43}
Suppose that
\begin{list}{}
{\setlength{\leftmargin}{27pt}
\setlength{\labelsep}{-00pt}
\setlength{\labelwidth}{00pt}
\setlength{\itemindent}{-00pt}
\addtolength{\topsep}{-08pt}
\addtolength{\parskip}{-12pt}
\addtolength{\itemsep}{-05pt}
}
\item{\num{1}}
$\fourtpl{E}{X}{\calS}{\calE}$ is a subspace choice system
and $G \leq H(X)$,
\item{\num{2}}
$\alpha \in \rfs{MBC}$ and $\itGamma \subseteq \rfs{MC}$,
\item{\num{3}}
$F$ is a normed space, $Y \subseteq F$ is open
and $\iso{\tau}{X}{Y}$,
\item{\num{4}}
for every $g \in G$, $g^{\tau}$ is locally $\itGamma$-bicontinuous.
\end{list}

\num{a} Assume that
in addition to \num{1}\,-\,\num{4}
the following conditions are fulfilled.
\begin{list}{}
{\setlength{\leftmargin}{27pt}
\setlength{\labelsep}{-00pt}
\setlength{\labelwidth}{00pt}
\setlength{\itemindent}{-00pt}
\addtolength{\topsep}{-08pt}
\addtolength{\parskip}{-12pt}
\addtolength{\itemsep}{-05pt}
}
\item{\num{a1}} For every $x \in X$, if $x \in S \in \calS$,
then $\lambda_{\sboldbbT}^{E;E_S} \drest S$
is $(\alpha,G)$-decayable at $x$.
\item{\num{a2}} For every $x \in X$,
$G$ is $\alpha$-infinitely closed at $x$.
\item{\num{a3}} $\itGamma$ is a modulus of continuity.
\item{\num{a4}} $\itGamma$ is $(\leq\kern-1pt \kappa(E))$-generated.
\end{list}
\underline{Then} $\tau\inverse$ is locally $\itGamma$-continuous.

\num{b} Assume that
in addition to \num{1}\,-\,\num{4}
the following conditions are fulfilled.
\begin{list}{}
{\setlength{\leftmargin}{27pt}
\setlength{\labelsep}{-00pt}
\setlength{\labelwidth}{00pt}
\setlength{\itemindent}{-00pt}
\addtolength{\topsep}{-08pt}
\addtolength{\parskip}{-12pt}
\addtolength{\itemsep}{-05pt}
}
\item{\num{b1}} For every $x \in X$, if $x \in S \in \calS$,
then $\lambda_{\sboldbbT}^{E;E_S} \drest S$
is compatible with $G$ at $x$.
\item{\num{b2}} For every $S \in \calS$,
$E_S$ is of the second category.
\item{\num{b3}} For every $\gamma \in \itGamma$ and $K > 0$,
$K \mcdot \gamma \in \itGamma$.
\item{\num{b4}} $\itGamma$ is countably generated.
\end{list}
\underline{Then} $\tau\inverse$ is locally $\itGamma$-continuous.

\num{c} Assume that
in addition to \num{1}\,-\,\num{4}
the following conditions are fulfilled.
\begin{list}{}
{\setlength{\leftmargin}{27pt}
\setlength{\labelsep}{-00pt}
\setlength{\labelwidth}{00pt}
\setlength{\itemindent}{-00pt}
\addtolength{\topsep}{-08pt}
\addtolength{\parskip}{-12pt}
\addtolength{\itemsep}{-05pt}
}
\item{\num{c1}} For every $x \in X$, if $x \in S \in \calS$,
then $\lambda_{\sboldbbA}^{E;E_S} \drest S$
is $(\alpha,G)$-decayable at $x$.
\item{\num{c2}} For every $x \in X$,
$G$ is $\alpha$-infinitely closed at $x$.
\item{\num{c3}} $\itGamma$ is a modulus of continuity.
\item{\num{c4}} $\itGamma$ is $(\leq\kern-1pt \kappa(E))$-generated.
\end{list}
\underline{Then} $\tau$ is locally $\itGamma$-bicontinuous.

\num{d} Assume that
in addition to \num{1}\,-\,\num{4}
the following conditions are fulfilled.
\begin{list}{}
{\setlength{\leftmargin}{27pt}
\setlength{\labelsep}{-00pt}
\setlength{\labelwidth}{00pt}
\setlength{\itemindent}{-00pt}
\addtolength{\topsep}{-08pt}
\addtolength{\parskip}{-12pt}
\addtolength{\itemsep}{-05pt}
}
\item{\num{d1}} For every $x \in X$, if $x \in S \in \calS$,
then $\lambda_{\sboldbbA}^{E;E_S} \drest S$
is compatible with $G$ at $x$.
\item{\num{d2}} For every $S \in \calS$,
$E_S$ is of the second category.
\item{\num{d3}} $\itGamma$ is a modulus of continuity.
\item{\num{d4}} $\itGamma$ is countably generated.
\end{list}
\underline{Then} $\tau$ is locally $\itGamma$-bicontinuous.
\end{cor}

\noindent
{\bf Proof } Parts (a), (b), (c) and (d) follow respectively
from theorems \ref{metr-bldr-t3.15}, \ref{metr-bldr-t3.26},
\ref{metr-bldr-t3.35} and \ref{metr-bldr-t3.41}.
\smallskip\hfill\myqed

There are cases in which the action is translation-like but not
affine-like. In such situations Parts (a) or (b) are applicable
but Parts (c) and (d) are not,
and hence we can only prove the $\itGamma$-continuity of
$\tau\inverse$.

For spaces of the first category only Parts (a) and (c) are applicable.
Part (c) has a conclusion stronger than that of Part (a).
However, the final theorem about groups of the form
$H_{\itGamma}^{\srfs{LC}}(X)$, (Theorem \ref{metr-bldr-t3.19}),
can be inferred from either Part (a) or Part (c).

The conclusion of Part (c) is stronger than that of Part (d).
But the assumptions of Part (c) are stronger in some respect than
those of Part~(d).
Nevertheless, we do not know to construct a group $G$
to which the reconstruction methods of Chapter \ref{s2} apply,
and for which Part~(d) can be applied but Part (c) cannot.

There are two outstanding open questions. The first is whether the
assumption that $\itGamma$ is 
$(\leq\kern-1pt \kappa(E))$-generated, is needed.
The second is whether translation-likeness of the action implies the
$\itGamma$-continuity of $\tau$.

\begin{question}\label{bddly-lip-bldr-q3.44}
\begin{rm}
Let $X,Y$ be open subsets of the normed spaces $E,F$,
and $\itGamma$ be a modulus of continuity.
Suppose that $\iso{\tau}{X}{Y}$ and that
$(H_{\itGamma}^{\srfs{LC}}(X))^{\tau} = H_{\itGamma}^{\srfs{LC}}(Y)$
is $\tau$ locally $\itGamma$-bicontinuous?
\end{rm}
\end{question}

\begin{question}\label{bddly-lip-bldr-q3.45}
\begin{rm}
Let $E$ and $F$ be normed space, $\iso{\tau}{X}{Y}$ and
$\itGamma$ be a countably generated modulus of continuity.
Suppose that $(\bbT(E))^{\tau} \subseteq H_{\itGamma}^{\srfs{LC}}(Y)$.
Is $\tau$ locally $\itGamma$-continuous?
Is the above true when $E,F$ are Banach spaces?
\end{rm}
\end{question}

%
%

\subsection{Normed manifolds.}
\label{ss3.6}
\label{ss3.6-normed-manifolds}

As in Chapter \ref{s2}, the results of this section extend to
normed manifolds. Also, the proofs presented to this point
transfer without change to this new context.
We now state some of these results explicitly.

\begin{defn}\label{metr-bldr-d3.46}
\begin{rm}
(a) Let $\pair{X}{\itPhi}$ be a normed manifold.
We say
that $\pair{X}{\itPhi}$ is a {\it locally Lipschitz normed manifold},
if for every $\varphi,\psi \in \itPhi$, $\varphi\inverse \scirc \psi$
is a bilipschitz function.

   \index{locally lipschitz normed manifold@@locally
          Lipschitz normed manifold}

(b) Let $\pair{X}{\itPhi}$ and $\pair{Y}{\itPsi}$ be locally
Lipschitz normed manifods and $\iso{\tau}{X}{Y}$.
We say that $\tau$ is
{\it Lipschitz with respect to $\itPhi$ and $\itPsi$},
if there is $K$ such that for every
$x \in X$ there are $\varphi \in \itPhi$ and $\psi \in \itPsi$
such that $x \in \rfs{int}(\rfs{Rng}(\varphi))$,
$\tau(x) \in \rfs{int}(\rfs{Rng}(\psi))$
and $\psi\inverse \scirc \tau \scirc \varphi$ is $K$-Lipschitz.
We say that
$\tau$ is {\it bilipschitz with respect to $\itPhi$ and $\itPsi$},
if both $\tau$ and $\tau\inverse$ are Lipschitz.
Define
$$
\rfs{LIP}(X,\itPhi) = \setm{h \in H(X)}
{h \mbox{ is bilipschitz with respect to } \itPhi}.
$$

   \index{lipschitz function between@@Lipschitz function between
   locally Lipschitz normed manifolds}

   \index{bilipschitz homeomorphism between
   locally Lipschitz normed manifolds}

   \index{N@lip08@@$\rfs{LIP}(X,\itPhi) = \setm{h \in H(X)}
   {\exists K \forall x(\exists \varphi,\psi \in \itPhi)
   (x \in \rfs{int}(\rfs{Rng}(\varphi)),\ 
   h(x) \in \rfs{int}(\rfs{Rng}(\psi))
   \mbox{ and }
   \newline\rule{24mm}{0pt}\psi\inverse \scirc h \scirc \varphi
   \mbox{ is $K$-bilipschitz}}$}

(c) Let $\pair{X}{\itPhi}$ and $\pair{Y}{\itPsi}$
be locally Lipschitz normed manifolds
and $\itGamma$ be a modulus of continuity.
A homeomorphism $\iso{\tau}{X}{Y}$
is {\it locally $\itGamma$-continuous with respect to
$\itPhi$ and $\itPsi$},
if for every
$x \in X$ there are $\varphi \in \itPhi$, $\psi \in \itPsi$,
$U \in \rfs{Nbr}(\varphi\inverse(x))$ and $\gamma \in \itGamma$
such that
$x \in \rfs{int}(\rfs{Rng}(\varphi))$,
$\tau(x) \in \rfs{int}(\rfs{Rng}(\psi))$,
$U \subseteq \rfs{Dom}(\varphi)$ and
$(\psi\inverse \scirc \tau \scirc \varphi) \nrestriction U$ is
$\gamma$-continuous.
   \index{locally $\itGamma$-continuous
   with respect to $\itPhi$ and $\itPsi$}
We say that
$\tau$ is locally $\itGamma$-bicontinuous, if $\tau$ and $\tau\inverse$
are locally $\itGamma$-continuous. Define
$$
H_{\itGamma}^{\srfs{LC}}(X,\itPhi) =
\setm{h \in H(X)}
{h \mbox{ is locally } \itGamma\mbox{-bicontinuous with respect to }
\itPhi}.
$$

   \index{locally $\itGamma$-bicontinuous
   with respect to $\itPhi$ and $\itPsi$}
   \index{N@hlc05@@$H_{\itGamma}^{\srfs{LC}}(X,\itPhi) =
   \setm{h \in H(X)}
   {\forall x(\exists \varphi,\psi \in \itPhi)
   (x \in \rfs{int}(\rfs{Rng}(\varphi)),\ 
   h(x) \in \rfs{int}(\rfs{Rng}(\psi))
   \mbox{ and }
   \newline\rule{25mm}{0pt}\psi\inverse \scirc h \scirc \varphi
   \mbox{ is $\itGamma$-bicontinuous at $\varphi\inverse(x)$}}$}

(d) Let $\pair{X}{\itPhi}$ be a locally Lipschitz normed manifold,
$\calS$ be an open cover of $X$
and $\itGamma$ be a modulus of continuity.
Define $\rfs{LIP}(X,\itPhi,\calS)$ to be the group generated by
$\bigcup \setm{\rfs{LIP}(X,\itPhi)\sprt{S}}{S \in \calS}$
and $H^{\srfs{LC}}_{\itGamma}(X,\itPhi,\calS)$ to be
the group generated by
$\bigcup
\setm{H^{\srfs{LC}}_{\itGamma}(X,\itPhi)\sprt{S}}{\break S \in \calS}$.
   \index{N@lip09@@$\rfs{LIP}(X,\itPhi,\calS)$}
   \index{N@hlc06@@$H^{\srfs{LC}}_{\itGamma}(X,\itPhi,\calS)$}
\end{rm}
\end{defn}

\begin{theorem}\label{t3.24}\label{metr-bldr-t3.47}
Let $\pair{X}{\itPhi}$ and $\pair{Y}{\itPsi}$
be normed manifolds with locally Lipschitz atlasses
and $\iso{\tau}{X}{Y}$.
Let $\itGamma$ be a countably generated modulus of continuity.

\num{a}
Suppose that
$(\rfs{LIP}(X,\itPhi))^{\tau} \subseteq
H_{\itGamma}^{\srfs{LC}}(Y,\itPsi)$.
Then $\tau$ is locally
$\itGamma$-bicon\-tinuous
with respect to $\itPhi$ and $\itPsi$.

\num{b} Let $\calS$ be an open cover of $X$,
and suppose that
$(\rfs{LIP}(X,\itPhi,\calS))^{\tau} \subseteq
H_{\itGamma}^{\srfs{LC}}(Y,\itPsi)$.
Then $\tau$ is locally $\itGamma$-bicontinuous
with respect to $\itPhi$ and $\itPsi$.
\end{theorem}

Note that Part (a) is a special case of Part (b).

We simplify the notation below by omitting the mention of $\itPhi$
and $\itPsi$.

\begin{cor}\label{c3.25}\label{metr-bldr-c3.48}
Let $\pair{X}{\itPhi}$ and $\pair{Y}{\itPsi}$
be normed manifolds with locally Lipschitz atlasses.

\num{a}
{\thickmuskip=2mu \medmuskip=1mu \thinmuskip=1mu 
Let $\itGamma$ and $\itDelta$ be countably generated
moduli of continuity,
and
$\iso{\varphi}{H_{\itGamma}^{\srfs{LC}}(X)}
{H_{\itDelta}^{\srfs{LC}}(Y)}$.
Then $\itGamma = \itDelta$
and there is $\iso{\tau}{X}{Y}$ such that $\tau$ induces $\varphi$,
and $\tau$ is locally $\itGamma$-bicontinuous.
}

\num{b} Let $\itGamma$ be a countably generated modulus of continuity,
$\calS$ an open cover of $X$, and $G \leq H(X)$.
Assume that if $\pair{X}{\itPhi}$ is a Banach manifold,
then $\rfs{LIP}(X,\calS) \leq G$,
and if $\pair{X}{\itPhi}$ is not a Banach manifold,
then $\rfs{LIP}^{\srfs{LC}}(X,\calS) \leq G$.
Suppose that $\iso{\varphi}{G}{H_{\itGamma}^{\srfs{LC}}(Y)}$.
\newline
Then $G = H_{\itGamma}^{\srfs{LC}}(X)$
and there is $\iso{\tau}{X}{Y}$ such that
$\tau$ induces $\varphi$, and $\tau$ is locally
$\itGamma$-bicon\-tinuous.
\end{cor}

\noindent
{\bf Proof }
(a) Note that if
$H_{\itGamma}^{\srfs{LC}}(X) = H_{\itDelta}^{\srfs{LC}}(X)$,
then $\itGamma = \itDelta$.
Hence Part (a) can be concluded from Part (b).

(b) We shall apply Theorem \ref{t2.15}(a).
Clearly,
$\rfs{LIP}^{\srfs{LC}}(Y;\itPsi) \leq H_{\itGamma}^{\srfs{LC}}(Y)$.
(See Definition~\ref{d2.14}(b)).
There is an atlas $\itPhi'$ for $X$ such that
if $\pair{X}{\itPhi}$
is a Banach manifold,
then $\rfs{LIP}(X,\itPhi') \leq G$,
and if $\pair{X}{\itPhi}$ is not a Banach manifold,
then $\rfs{LIP}^{\srfs{LC}}(X,\itPhi') \leq G$. Indeed,
{\thickmuskip=4mu \medmuskip=3mu \thinmuskip=2mu 
$\itPhi' =
\setm{\psi \nrestriction \overB(x,r)}
{\psi \in~\kern-5pt\itPhi,\ \overB(x,r) \subseteq \rfs{Dom}(\psi)
\mbox{ and there is } U \in \calS
\mbox{ such that }
\psi(\overB(x,r))~\subseteq~U}$.
}
By Theorem \ref{t2.15}(a), there is $\iso{\tau}{X}{Y}$
such that $\tau$ induces $\varphi$.
So $G^{\tau} = H_{\itGamma}^{\srfs{LC}}(Y)$.
In particular,
$(\rfs{LIP}(X,\calS))^{\tau} \subseteq  H_{\itGamma}^{\srfs{LC}}(Y)$.
By Theorem \ref{metr-bldr-t3.47}(b), $\tau$ is
locally $\itGamma$-bicontinuous.
So $G = H_{\itGamma}^{\srfs{LC}}(X)$.
\hfill\myqed

\noindent
{\bf Question } In the above theorem does it suffice to assume that
$\rfs{LIP}(X,\calS) \leq G$, regardless of whether
$\pair{X}{\itPhi}$ is a Banach Manifold?

\newpage


\section{The local uniform continuity of conjugating
homeomorphisms}
\label{s4}

To complete the picture on the local $\itGamma$-bicontinuity of
conjugating homeomorphisms, we now deal with the group
$H_{\srfs{MC}}^{\srfs{LC}}(X)$ of locally bi-uniformly-continuous
homeomorphisms.
(Note that $\rfs{MC}$ is a modulus of continuity, so the notation
$H_{\srfs{MC}}^{\srfs{LC}}(X)$ is
a special case of definition \ref{d1.6}(c)).
The methods employed in dealing with $H_{\srfs{MC}}^{\srfs{LC}}(X)$
are quite different from those used in the previous section.

We shall prove the following extension of
Theorem \ref{metr-bldr-t3.42}

\begin{theorem}\label{metr-bldr-t4.1}
\num{a} Let $X,Y$ be open subsets of the normed spaces
$E$ and $F$ respectively,
$\itGamma,\itDelta$ be moduli of continuity
and
$\iso{\varphi}{H_{\itGamma}^{\srfs{LC}}(X)}
{H_{\itDelta}^{\srfs{LC}}(Y)}$.
Suppose that $\itGamma$ is $(\leq\kern-1pt \kappa(E))$-generated
or $\itGamma = \rfs{MC}$.
\underline{Then} $\itGamma = \itDelta$,
there is $\iso{\tau}{X}{Y}$ such that
$\varphi(h) = h^{\tau}$ for every
$h \in H_{\itGamma}^{\srfs{LC}}(X)$,
and $\tau$ is locally $\itGamma$-bicontinuous.
 
\num{b}
Let $\fourtpl{E}{X}{\calS}{\calE}$ and $\fourtpl{F}{Y}{\calT}{\calF}$
be subspace choice systems, $\itGamma,\itDelta$ be moduli of continuity
and
$\iso{\varphi}
{H_{\itGamma}^{\srfs{LC}}(X;\calS,\calE)}
{H_{\itDelta}^{\srfs{LC}}(Y;\calT,\calF)}$.
Suppose that $\itGamma$
is $(\leq\kern-1pt \kappa(E))$-generated or $\itGamma = \rfs{MC}$,
and the same holds for $\itDelta$.
\underline{Then} $\itGamma = \itDelta$,
there is $\iso{\tau}{X}{Y}$ such that
$\varphi(h) = h^{\tau}$ for every
$h \in H_{\itGamma}^{\srfs{LC}}(X;\calS,\calE)$,
and $\tau$ is locally $\itGamma$-bicontinuous.
\end{theorem}

Note that Part (a) is not a special case of Part (b),
since in (b) $\itDelta$ is assumed to be
$(\leq\kern-1pt \kappa(E))$-generated or equal to $\rfs{MC}$,
and this is not assumed in (a).
The key intermediate step in the proof of Theorem \ref{metr-bldr-t4.1}
is Theorem \ref{metr-bldr-t4.8}.

There are several ways of defining uniform continuity.
We sort this matter out in the next definition and proposition.

\begin{defn}\label{metr-bldr-d4.2}
\begin{rm}
{\thickmuskip=2mu \medmuskip=1mu \thinmuskip=1mu 
(a) Let $\pair{X}{d^X}$ and $\pair{Y}{d^Y}$ be metric spaces,
and $\fnn{f}{X}{Y}$.
}
We say that $f$ is {\it uniformly continuous} ($f$ is UC),
if for every $\varepsilon > 0$ there is $\delta > 0$
such that for every $x,y \in X$:
if $d^X(x,y) < \delta$, then $d^Y(f(x),f(y)) < \varepsilon$.
   \index{uniformly continuous}
   \index{uc@@UC. Abbreviation of uniformly continuous}
If $\iso{f}{X}{f(X)}$ and both $f$ and $f\inverse$ are
uniformly continuous, then $f$ is said to be
{\it bi-uniformly-continuous (bi-UC)}.
   \index{bi-uniformly-continuous}
   \index{bi-UC. Abbreviation of bi-uniformly-continuous}

(b) Let $\alpha \in \rfs{MC}$ and $r > 0$. We say that
$\fnn{f}{X}{Y}$ is {\it $(r,\alpha)$-continuous},
if for every $x,y \in X$:
if $d^X(x,y) < r$, then $d^Y(f(x),f(y)) \leq \alpha(d^X(x,y))$.
   \index{continuous. $(r,\alpha)$-continuous}

(c) We say that $\fnn{f}{X}{Y}$ is
{\it uniformly continuous for all distances},
if there is $\alpha \in \rfs{MC}$ such that $f$ is $\alpha$-continuous.
   \index{uniformly continuous for all distances}

(d) Let $\fnn{f}{X}{Y}$ and $x \in X$.
Say that {\it $f$ is uniformly continuous at $x$ ($f$ is UC at $x$)},
if there is $U \in \rfs{Nbr}(x)$ such that $f \nrestriction U$ is UC,
and {\it $f$ is bi-uniformly-continuous at $x$ (bi-UC at $x$)},
if there is $U \in \rfs{Nbr}(x)$ such that $f \nrestriction U$ is bi-UC.

   \index{uniformly continuous at $x$}
   \index{uc at x@@UC at $x$.
          Abbreviation of uniformly continuous at $x$}
   \index{bi-uniformly-continuous at $x$}
   \index{bi-UC at $x$. Abbreviation of bi-uniformly-continuous at $x$}

(e) Let $\fnn{f}{X}{Y}$.
Say that $f$ is {\it locally uniformly continuous (locally UC)},
if $f$ is UC at every $x \in X$,
and $f$ is {\it locally bi-uniformly-continuous (locally bi-UC)},
if $f$ is bi-UC at every $x \in X$.

   \index{locally uniformly continuous}
   \index{locally uc@@locally UC.
          Abbreviation of locally uniformly continuous}
   \index{locally bi-uniformly-continuous}
   \index{locally bi-uc@@locally bi-UC.
          Abbreviation of locally bi-uniformly-continuous}

(f) Let $\pair{X}{d}$ be a metric space.
{\it The discrete path property for large distances} is the following
property of $X$.
There are $a,b > 0$
such that for every $x,y \in X$ and $r > 0$
there are $n \in \bbN$ and $x = x_0,\,x_1,\dots,\,x_n = y$ in $X$
such that for every $i < n$, $d(x_i,x_{i + 1}) < r$ and
$\sum_{i < n}d(x_i,x_{i + 1}) \leq ad(x,y) + b$.

   \index{discrete path property for large distances}

\end{rm}
\end{defn}

\begin{prop}\label{metr-bldr-p4.3}
\num{a} Let $\fnn{f}{X}{Y}$. Then $f$ is UC
iff for some $\alpha \in \rfs{MC}$ and $r > 0$,\break
$f$ is $(r,\alpha)$-continuous.

\num{b} Suppose that $X$ has the discrete path property for large
distances. Let $\fnn{f}{X}{Y}$.
Then $f$ is UC iff $f$ is uniformly continuous for all distances.

\num{c} Suppose that $\fnn{f}{X}{Y}$, $f$ is UC and
$\rfs{Rng}(f)$ is bounded.
Then $f$ is uniformly continuous for all distances.

\num{d}
{\thickmuskip=2mu \medmuskip=1mu \thinmuskip=1mu 
Let $\fnn{f}{X}{Y}$ and $x \in X$.
Then $f$ is UC at $x$, iff for some $\alpha \in \rfs{MC}$,
$f$ is $\alpha$-continuous at $x$.
}
\end{prop}

\noindent
{\bf Proof } All parts are trivial.
However, the proof of implication $\Rightarrow$ in (a)
requires the following fact.
If $\fnn{\eta}{(0,a]}{[0,\infty)}$, and $\limtz{t} \eta(t) = 0$,
then there is $\alpha \in \rfs{MC}$ such that
$\eta \leq \alpha \nrestriction (0,a]$.
The verification of this fact is left to the reader.
\rule{0pt}{0pt}\hfill\myqed

\begin{defn}\label{metr-bldr-d4.4}
\begin{rm}
(a) Suppose that $X$, $Y$ are topological spaces $D \subseteq X$. Define
$H(X,Y) = \setm{h}{\iso{h}{X}{Y}}$ and
$H(X;D) = \setm{h \in H(X)}{h(D) = D}$.

   \index{N@h08@@$H(X,Y) = \setm{h}{\iso{h}{X}{Y}}$}
   \index{N@h09@@$H(X;D) = \setm{h \in H(X)}{h(D) = D}$}

(b) For metric spaces $X,Y$ define
$\rfs{UC}(X,Y) = \setm{h \in H(X,Y)}{h \mbox{ is UC}}$,
$\rfs{UC}^{\wpm}(X,Y) = \setm{h \in H(X,Y)}{h \mbox{ is bi-UC}}$
and $\rfs{UC}(X) = \rfs{UC}^{\wpm}(X,X)$.
For $x \in X$ let
$\rfs{PNT.UC}(X,x) = \setm{h \in H(X)}
{h(x) = x \mbox{ and } h \mbox{ is bi-UC at } x}$.

   \index{N@uc00@@$\rfs{UC}(X,Y) = \setm{h \in H(X)}{h \mbox{ is UC}}$}
   \index{N@uc01@@$\rfs{UC}^{\wpm}(X,Y) =
          \setm{h \in H(X)}{h \mbox{ is bi-UC}}$}
   \index{N@uc02@@$\rfs{UC}(X) = \rfs{UC}^{\wpm}(X,X)$}
   \index{N@pntuc@@$\rule{1pt}{0pt}\rfs{PNT.UC}(X,x) = \setm{h \in H(X)}
          {h(x) = x \mbox{ and } h \mbox{ is bi-UC at } x}$}

(c) Let $X$ be an open subset of a normed space $E$, $S \subseteq X$
be open,
and $F$ be a dense linear subspace of $E$. Define
$\rfs{UC}(X;F) = \setm{h \in \rfs{UC}(X)}{h(X \cap F) = X \cap F}$
and
$\rfs{UC}(X;S,F) = \rfs{UC}(X)\sprt{S} \cap \rfs{UC}(X;F)$.
For $x \in S$ \,let
$\rfs{UC}(X;S,F,x) = \setm{h \in \rfs{UC}(X;S,F)}{h(x) = x}$.

   \index{N@uc03@@$\rfs{UC}(X;F) =
          \setm{h \in \rfs{UC}(X)}{h(X \cap F) = X \cap F}$}
   \index{N@uc04@@$\rfs{UC}(X;S,F) =
          \rfs{UC}(X)\sprt{S} \cap \rfs{UC}(X;F)$}
   \index{N@uc05@@$\rfs{UC}(X;S,F,x) =
          \setm{h \in \rfs{UC}(X;S,F)}{h(x) = x}$}

(d) Let $\frtpl{E}{X}{\calS}{\calF}$ be a subspace choice system.
Then $\rfs{UC}(X,\calS)$
denotes the subgroup of $H(X)$ generated by
$\bigcup \setm{\rfs{UC}(X)\sprt{S}}{S \in \calS}$,
and $\rfs{UC}(X;\calS,\calF)$
denotes the subgroup of $H(X)$ generated by
$\bigcup \setm{\rfs{UC}(X;S,F_S)}{S \in \calS}$.

   \index{N@uc06@@$\rfs{UC}(X,\calS)$.
   The subgroup of $H(X)$ generated by
   $\bigcup \setm{\rfs{UC}(X)\sprt{S}}{S \in \calS}$}
   \index{N@uc06@@$\rfs{UC}(X;\calS,\calF)$.
   {\thickmuskip=2mu \medmuskip=1mu \thinmuskip=1mu 
   The subgroup of $H(X)$ generated by
   $\bigcup \setm{\rfs{UC}(X;S,F_S)}{S \in \calS}$}}

(e) For metric spaces $X,Y$ let
$\rfs{LUC}(X,Y) = \setm{h \in H(X,Y)}{h \mbox{ is locally UC}}$.
   \index{N@luc00@@$\rfs{LUC}(X,Y) =
   \setm{h \in H(X,Y)}{h \mbox{ is locally UC}}$}
   \index{N@luc01@@$\rfs{LUC}^{\wpm}(X,Y) =
   \setm{h \in H(X,Y)}{h \mbox{ is locally bi-UC}}$}
   \index{N@luc02@@$\rfs{LUC}(X) = \rfs{LUC}^{\wpm}(X,X)$}
As usual we define
$\rfs{LUC}^{\wpm}(X,Y) =
\setm{h \in H(X,Y)}{h \mbox{ is locally bi-UC}}$
and $\rfs{LUC}(X) = \rfs{LUC}^{\wpm}(X,X)$.
\end{rm}
\end{defn}

\noindent
{\bf Remark } Note that $H_{\srfs{MC}}(X) \leq \rfs{UC}(X)$ but equality
need not hold. See Proposition \ref{metr-bldr-p4.3}.
It is the group $H_{\srfs{MC}}(X)$ that fits into the framework better,
but the group which has been traditionally considered is $\rfs{UC}(X)$.
We based the above definitions on $\rfs{UC}(X)$ rather than
on $H_{\srfs{MC}}(X)$.
As for local uniform continuity, the two ways of defining
this notion are equivalent.
Hence $\rfs{LUC}(X) = H_{\srfs{MC}}^{\srfs{LC}}(X)$
for every metric space $X$. This fact is a triviality.

The following easy proposition will be used extensively.

\begin{prop}\label{p4.9}\label{metr-bldr-p4.5}
Let $X$ be a metric space and $\setm{U_n}{n \in \bbN}$ be a sequence
of open sets in $X$ such that
$\lim_{n \rightarrow \infty} \rfs{diam}(U_n) = 0$, and
for every distinct $m,n \in \bbN$, $d(U_m,U_n) > 0$.
For every $n \in \bbN$ let $h_n \in \rfs{UC}(X)$ be such that
$\rfs{supp}(h_n) \subseteq U_n$.
Then $\bcirc_{n \in \sboldbbN} h_n \in \rfs{UC}(X)$.
\end{prop}

\noindent
{\bf Proof }
Let $h = \bcirc_{n \in \sboldbbN} h_n$.
Let $\varepsilon > 0$. Let $N \in \bbN$ be such that for every
$m \geq N$,
$\rfs{diam}(U_m) < \dgfrac{\varepsilon}{3}$.
Let $\delta_1 > 0$ be such that
for every $i < N$ and $x,y \in X$:
if $d(x,y) < \delta_1$,
then $d(h_i(x), h_i(y)) < \dgfrac{\varepsilon}{3}$.
Let $\delta_2 = \rfs{min}(\setm{d(U_i,U_j)}{ i < j < N})$,
and let
$\delta = \rfs{min}(\delta_1, \delta_2, \dgfrac{\varepsilon}{3})$.

Suppose that $d(x,y) < \delta$,
and we show that $d(h(x),h(y)) < \varepsilon$.
Since for every distinct $i,j < N$, $d(x,y) < d(U_i,U_j)$, there are
no distinct $i,j < N$ such that $x \in U_i$ and $y \in U_j$. So
we may assume that one of the following occurs:
(i) for some $i < N$, $x \in U_i$ and
$y \not\in \bigcup \setm{U_j}{j \neq i}$;
(ii) for some $i < N$ and $j \geq N$, $x \in U_i$ and
$y \in U_j$;
(iii) for some $i \geq N$, $x \in U_i$ and
$y \not\in \bigcup \setm{U_j}{j \neq i}$;
(iv) for some distinct $i,j \geq N$, $x \in U_i$ and $y \in U_j$;\break
(v) $x,y \not\in \bigcup \setm{U_i}{i \in \bbN}$.

In Case (i), $h(x) = h_i(x)$ and $h(y) = h_i(y)$,
so since $d(x,y) < \delta_1$, it follows that
$d(h(x),h(y)) < \varepsilon$.
In Case (ii), 
$$
d(h(x),h(y)) \leq 
d(h(x),y) + d(y,h(y)) =
d(h_i(x),h_i(y)) + d(y,h_j(y)) < \dgfrac{\varepsilon}{3} +
\dgfrac{\varepsilon}{3} < \varepsilon.
$$
In Case (iii),
$$
d(h(x),h(y)) =
d(h_i(x),h_i(y)) \leq 
d(h_i(x),x) + d(x,y) + d(y,h_i(y)) <
\dgfrac{\varepsilon}{3} + \dgfrac{\varepsilon}{3} +
\dgfrac{\varepsilon}{3} = \varepsilon.
$$
Case (iv) is similar to Case (iii), and Case (v) is trivial.
\hfill\myqed

\begin{defn}\label{metr-bldr-d4.6}
\begin{rm}
Let $M$ be a topological space and $N$ be a Hausdorff space.

(a) Let $A \subseteq M$ and $\fnn{g}{A}{N}$ be continuous.
For every $x \in \rfs{cl}^M(A)$ there is at most one $y \in N$
such that $g \cup \sngltn{\pair{x}{y}}$ is a continuous function.
Let
$$
g^{\srfs{cl}}_{M,N} = \setm{\pair{x}{y}}
{x \in \rfs{cl}^M(A), \, y \in N \mbox{ and }
g \cup \sngltn{\pair{x}{y}} \mbox{ is a continuous function}}.
$$
Obviously, $g^{\srfs{cl}}_{M,N}$ extends $g$,
and $\rfs{Rng}(g^{\srfs{cl}}_{M,N}) \subseteq
\rfs{cl}^N(\rfs{Rng}(g))$.
When possible, $g^{\srfs{cl}}_{M,N}$ is abbreviated by $g^{\srfs{cl}}$,
and if $M = N$, then $g^{\srfs{cl}}_{M,N}$ is denoted by
$g^{\srfs{cl}}_M$.
If $H$ is a set of continuous functions from $A$ to $B$,
then $H^{\srfs{cl}}$ denotes $\setm{h^{\srfs{cl}}}{h \in H}$. 

(b) Let $X \subseteq M$ and $Y \subseteq N$.
We define
$$
\rfs{EXT}^{M,N}(X,Y) =
\setm{h \in H(X,Y)}{\rfs{Dom}(h^{\srfs{cl}}_{M,N}) = \rfs{cl}^M(X)}.
$$
When possible, we abbreviate
$\rfs{EXT}^{M,N}(X,Y)$ by
$\rfs{EXT}(X,Y)$.
The notation $\rfs{EXT}^M(X)$
stands for $(\rfs{EXT}^{M,M})^{\pm}(X,X)$.

   \index{N@AAAA@@
{\thickmuskip=3mu \medmuskip=2mu \thinmuskip=1mu 
$g^{\srfs{cl}}_{M,N} = \setm{\pair{x}{y}}
{x \in \rfs{cl}^M(A), \, y \in N \mbox{ and }
g \cup \sngltn{\pair{x}{y}} \mbox{ is a continuous function}}$}}
   \index{N@AAAA@@$g^{\srfs{cl}}$.
          Abbreviation of $g^{\srfs{cl}}_{M,N}$}
   \index{N@AAAA@@$g^{\srfs{cl}}_M$.
          Abbreviation of $g^{\srfs{cl}}_{M,M}$}
   \index{N@AAAA@@$H^{\srfs{cl}}  = \setm{h^{\srfs{cl}}}{h \in H}$} 
   \index{N@ext01@@$\rfs{EXT}^{M,N}(X,Y) =
   \setm{h \in H(X,Y)}{\rfs{Dom}(h^{\srfs{cl}}_{M,N}) = \rfs{cl}^M(X)}$}
   \index{N@ext02@@$\rfs{EXT}(X,Y)$.
          Abbreviation of $\rfs{EXT}^{M,N}(X,Y)$}
   \index{N@ext03@@$\rfs{EXT}^M(X)$.
          Abbreviation of $(\rfs{EXT}^{M,M})^{\pm}(X,X)$}
\end{rm}
\end{defn}

\begin{prop}\label{metr-bldr-p4.7}\label{p4.28}
\num{a} \num{i}
Let $X$ be a topological space, $D \subseteq X$ be dense,
$Y$  be a regular topological space and $\fnn{h}{D}{Y}$ be continuous.
Suppose that for every $x \in X$
there is a continuous function $\fnn{h_x}{D \cup \sngltn{x}}{Y}$
extending $h$. Then $\bigcup \setm{h_x}{x \in X}$ is continuous.

\num{ii} Let $M$ be a topological space, $N$ be a regular space
$A \subseteq M$ and $\fnn{g}{A}{N}$ be continuous.
Then $g^{\srfs{cl}}_{M,N}$ is continuous.

\num{b} Let $X$ be a metric space, $Y$ be a complete metric space,
$A \subseteq X$,
and $\fnn{g}{A}{Y}$ be a uniformly continuous function.
Then $\rfs{Dom}(g^{\srfs{cl}}) = \rfs{cl}(A)$.

\num{c}
Let $E$ be a normed space, $D$ be a dense linear subspace of $E$,
$X \subseteq E$ be open, $u \in D$, $B^E(u,p) \subseteq X$,
$x,y \in D \cap B^E(u,p)$, $z \in B^{\oversE}(u,p)$,
$\varepsilon > 0$,
$0 < s < \min(\norm{x - z},\norm{y - z})$
and $\max(\norm{x - z},\norm{y - z}) < t < \norm{z - u} + p$.
Then there is $h \in \rfs{UC}(X;D)$
such that:\break
\num{i}
$\rfs{supp}(h^{\srfs{cl}}_{\oversE}) \subseteq
B^{\oversE}(z,t) - B^{\oversE}(z,s)$,
\num{ii} $h(x) = x$ and \num{iii} $h(y) \in B(x,\varepsilon)$.
\end{prop}

\noindent
{\bf Proof }
The proofs of Parts (a) and (b) are trivial.

(c) Write $r' = \norm{z - u} + t$. For every $a \in (0,1)$ there is
$h \in \rfs{LIP}(X;D)\sprtl{B(u,p)}$ such that
$h \nrestriction B^E(u,r')$ is the multiplication by the scalar
$\dgfrac{a}{r'}$,
that is, $h(w) = \frac{a}{r'} w$ for every $w \in B^E(u,r')$.
So we may assume that
$B^{\oversE}(z,t) \subseteq B^{\oversE}(u,ap)$.
Let
$s < \bars < \min(\norm{x - z},\norm{y - z})$,
$t > \bart > \max(\norm{x - z},\norm{y - z})$ and $\barz \in D$
be such that $\norm{\barz - z} < \min(t - \bart,\bars - s)$.
Since $\rfs{tr}_{\barz}^E$ is an isometry belonging to $H(E;D)$,
we may shift $\barz$ to the origin.
That is, we may assume that $\barz = 0$.
$\norm{x} \geq \norm{x - z} - \norm{z} > s - (s - \bars) = \bars$.
The same computation applies to $y$. We conclude that
$\norm{x},\norm{y} > \bars$.
Another similar computation shows that $\norm{x},\norm{y} < \bart$.
It is also obvious that
$B^{\oversE}(z,s) \subseteq B^{\oversE}(0,\bars)$ and that
$B^{\oversE}(z,t) \supseteq B^{\oversE}(0,\bart)$.
It thus remains to show that there is $h \in \rfs{UC}(X;D)$
such that
$\rfs{supp}(h) \subseteq B(0,\bart) - B(0,\bars)$,
and $h$ fulfills Clauses (ii) and (iii).
The construction of such a homeomorphism is routine but long,
so we skip some details.

In the inclusion
$B^{\oversE}(z,t) \subseteq B^{\oversE}(u,ap)$,
choose $a$ small enough so that we have\break
$B^E(0,6 \max(\norm{x},\norm{y})) \subseteq X$.
By an argument similar to the choice of $a$ above,
we may also assume that
(1) \,$\bart > 5 \max(\norm{x},\norm{y})$
and $\bars < \frac{1}{5} \min(\norm{x},\norm{y})$.
Let $F = \rfs{span}(\dbltn{x}{y})$.
As in the proof of Claim 3 in Lemma \ref{metr-bldr-l3.30},
there is $E_1$ such that $F \oplus E_1 = E$,
and $\norm{v_0} + \norm{v_1} \leq 3 \norm{v_0 + v_1}$
for every $v_0 \in F$ and $v_1 \in E_1$.
Let $\norm{\ }^{\srbfs{H}}$ be a Hilbert norm on $F$ such that
(2) $\norm{v} \leq \norm{v}^{\srbfs{H}} \leq 3 \sqrt{2} \norm{v}$
for every $v \in F$.

For $v \in E$ let $v_F$ and $v_{E_1}$ be such that
$v = v_F + v_{E_1}$ and define
$\bboldnorm{v} = \norm{v_F}^{\srbfs{H}} + \norm{v_{E_1}}$.
We may assume that $\norm{y}^{\srbfs{H}} \neq \norm{x}^{\srbfs{H}}$.
Let
$S = \setm{v \in F}{\norm{v}^{\srbfs{H}} = \norm{y}^{\srbfs{H}}}$.
By (1) and (2),
$S \subseteq B^E(0,\bart) - \overB^E(0,\bars)$.
So there is $b > 0$ such that 
$x \not\in B^{\pair{E}{\kern0.5pt\sbboldnorm{\kern1.7pt}}}(S,b)
\subseteq
B^E(0,\bart) - \overB^E(0,\bars)$.

Suppose that the angle between $x$ and $y$ in
$\pair{F}{\norm{\ }^{\srbfs{H}}}$ is $\theta$.
Let $\fnn{\eta}{[0,\infty)}{[0,\infty)}$ be the piecewise linear
function with a unique
breakpoint at $b$ such that $\eta(0) = \theta$ and 
$\eta(b) = 0$.
For $v \in X$ define
$h_1(v) = \rfs{Rot}_{\eta(\sbboldnorm{v})}(v_F) + v_{E_1}$,
where $\rfs{Rot}_{\phi}$ is rotation by angle $\phi$ in $F$.
Obviously,
$h_1 \in \rfs{LIP}(E;D)$,
$\rfs{supp}(h_1) \subseteq B^E(0,\bart) - \overB^E(0,\bars)$,
$h_1(x) = x$ and for some $c > 0$, $h_1(y) = cx$.
It is easy to construct radial homeomorphism $h_2 \in \rfs{LIP}(E;D)$,
such that
$\rfs{supp}(h_2) \subseteq B^E(0,\bart) - \overB^E(0,\bars)$,
$h_2(x) = x$ and $h_2(cy) \in B(x,\varepsilon)$.
So $h = h_2 \scirc h_1$ is as required.
\smallskip\hfill\myqed

Theorem \ref{metr-bldr-t4.8} is phrased in a way that Part (a)
is easiest to read, Part (b) is the main statement of the theorem,
and Part (c) is the ``pointwise'' version of Part (b).
So (c) $\Rightarrow$ (b) $\Rightarrow$ (a),
and we actually prove (c).

Note that Theorem \ref{metr-bldr-t4.8}(b) is analogous
to Theorem \ref{metr-bldr-t3.27},
but the assumption here is that
$(\rfs{UC}(X))^{\tau} \subseteq \rfs{LUC}(Y)$, whereas
in \ref{metr-bldr-t3.27} the weaker assumption that
$(\rfs{LIP}(X))^{\tau} \subseteq H_{\itGamma}^{\srfs{LC}}(Y)$
did suffice.

\begin{theorem}\label{metr-bldr-t4.8}
\num{a}
Let $X,Y$ be open subsets of the normed spaces $E$ and $F$,
and $\tau \in H(X,Y)$ be such that
\hbox{$(\rfs{UC}(X))^{\tau} \subseteq \rfs{LUC}(Y)$.
Then $\tau \in \rfs{LUC}^{\wpm}(X,Y)$.}

\num{b}
Let $\frtpl{E}{X}{\calS}{\calD}$ be a subspace choice system,
$Y$ be open subset of a normed space $F$
and $\tau \in H(X,Y)$. Suppose that
$(\rfs{UC}(X;\calS,\calD))^{\tau} \subseteq \rfs{LUC}(Y)$.
Then $\tau \in \rfs{LUC}^{\wpm}(X,Y)$.

\num{c}
Let $X,Y$ be open subsets of the normed spaces $E$ and $F$,
\hbox{$S \subseteq X$} be open,
$D$ be a dense linear subspace of $E$,
$x^* \in S$ and $\tau \in H(X,Y)$.
Suppose that
$(\rfs{UC}(X;S,D,x^*))^{\tau} \subseteq \rfs{PNT.UC}(Y,\tau(x^*))$.
Then $\tau$ is bi-UC at $x^*$.
\end{theorem}

\noindent
{\bf Proof }
(c) Let $X$, $Y$ etc.\ be as in Part (c).

{\bf Part 1 } $\tau$ is UC at $x^*$.
\newline
Suppose by contradiction that for every $U \in \rfs{Nbr}^X(x^*)$,
$\tau \nrestriction U$ is not UC.
The trivial proof of the following claim
is left to the reader.\smallskip

{\bf Claim 1.} For every $r > 0$ there are sequences
$\vecx,\vecy$ and $d,e > 0$ such that:
\begin{itemize}
\addtolength{\parskip}{-09pt}
\addtolength{\itemsep}{04pt}
\item[(1)]
$\rfs{Rng}(\vecx) \cup \rfs{Rng}(\vecy) \subseteq
B^X(x^*,\dgfrac{r}{2}) \cap D$;
\item[(2)]
$\limti{n} \norm{x_n - y_n} = 0$;
\item[(3)]
either (i) for any distinct $m,n \in \bbN$,
$d(\uopair{x_m}{y_m},\uopair{x_n}{y_n}) \geq e$,
or (ii) $\vecx$ is a Cauchy sequence;
\item[(4)]
$d(\rfs{Rng}(\vecx) \cup \rfs{Rng}(\vecy),x^*) > e$;
\item[(5)]
for every $n \in \bbN$,
$\norm{\tau(x_n) - \tau(y_n)} \geq d$.
\end{itemize}

Let $e_{-1} > 0$ be such that $B^E(x^*,e_{-1}) \subseteq S$.
It is easy to define by induction on $i \in \bbN$
a radius $r_i$, sequences $\vecx^i = \setm{x^i_n}{n \in \bbN}$,
$\vecy^i = \setm{y^i_n}{n \in \bbN}$ and
$d_i,e_i > 0$ such that $r_i = \dgfrac{e_{i - 1}}{8}$
and such that $\vecx^i,\vecy^i,d_i,e_i$
satisfy (1)\,-\,(5) of Claim~1 for $r = r_i$.
By deleting, if necessary, initial segments from each of the
$\vecx\farsu{i}$'s and $\vecy\farsu{i}$'s,
we may further assume that for every $i,n \in \bbN$,
$\norm{x^i_n - y^i_n} < \dgfrac{e_i}{4}$.
We may further assume that either for every $i \in \bbN$,
Clause (3)(i) of Claim 1 holds,
or for every $i \in \bbN$,
Clause (3)(ii) of Claim 1 holds.

{\bf Case 1 } Clause (3)(i) of Claim 1 holds.
Let $\setm{\pair{i(k)}{n(k)}}{k \in \bbN} \subseteq \bbN^2$
be a $\onetoonen$ sequence of pairs such that
$\limti{k} \norm{x^{i(k)}_{n(k)} - y^{i(k)}_{n(k)}} = 0$,
and for every $i \in \bbN$, $\setm{k}{i(k) = i}$ is infinite.
For every $k \in \bbN$ set $u_k = x^{i(k)}_{n(k)}$,
$v_k = y^{i(k)}_{n(k)}$,
$s_k = 2\norm{u_k - v_k}$ and $B_k = B(u_k,s_k)$.
Then it can be easily checked that for every distinct $k,l \in \bbN$,
$B_k \subseteq S$ and $d(B_k,B_l) > \dgfrac{e_{i(k)}}{4}$.
Also, $\limti{k} \rfs{diam}(B_k) = 0$.
Let $w_k \in [u_k,v_k] - \sngltn{u_k}$
be such that $\norm{\tau(w_k) - \tau(u_k)} < \dgfrac{1}{(k + 1)}$.
So $w_k \in B_k \cap D$.
By Lemma \ref{l2.6}(c),
there is $h_k \in \rfs{LIP}(X;S,D)$ such that
$\rfs{supp}(h_k) \subseteq B_k$, $h_k(u_k) = u_k$ and $h_k(w_k) = v_k$.

By Propostion \ref{metr-bldr-p4.5},
$h \eqdf \bcirc_{k \in \sboldbbN} h_k \in \rfs{UC}(X)$, and indeed
$h \in \rfs{UC}(X;S,D,x^*)$. However, we shall now see
that for every $V \in \rfs{Nbr}^Y(\tau(x^*))$,
$h^{\tau} \nrestriction V$ is not
uniformly continuous and hence
$h^{\tau} \not\in \rfs{PNT.UC}(Y,\tau(x^*))$ which is a contradiction.

Write $h^{\tau} = \hath$, $h(u_k) = \hatu_k$,
$h(v_k) = \hatv_k$ and $h(w_k) = \hatw_k$.
Then $\hath(\hatu_k) = \hatu_k$ and $\hatw_k = \hatv_k$.
There is $i$ such that for every
$n$, $\tau([x^i_n,y^i_n]) \subseteq V$.
Define $\sigma = \setm{k \in \bbN}{i(k) = i}$.
Then $\hatu_k, \hatv_k, \hatw_k \in V$ for every $k \in \sigma$.
Now, $\lim_{k \in \sigma} \norm{\hatu_k - \hatw_k} = 0$,
but
$\norm{\hath(\hatu_k) - \hath(\hatw_k)} =
\norm{\hatu_k - \hatv_k} \geq d_i$
for every $k \in \sigma$.
So $\hath \nrestriction V$ is not uniformly continuous.

{\bf Case 2 } Clause (3)(ii) of Claim 1 holds.
Let $\barz_i = \lim \vecx^i$. Note that
{\thickmuskip=3mu \medmuskip=2mu \thinmuskip=1mu 
$\barz_i \in \overE - E$.
{\thickmuskip=2mu \medmuskip=1mu \thinmuskip=1mu 
Clearly, $\barz_i \in B^{\oversE}(x^*,r_i) - B^{\oversE}(x^*,e_i)$.
Fix $i \in \bbN$ and for $j \in \bbN$ set
$t_{i,j} = \max(\norm{x^i_j - \barz_i},\norm{y^i_j - \barz_i})$
}
and $s_{i,j} = \min(\norm{x^i_j - \barz_i},\norm{y^i_j - \barz_i})$.
By taking a subsequence of $\setm{\pair{x^i_j}{y^i_j}}{j \in \bbN}$,
we may assume that for every $j$, $t_{i,j + 1} < s_{i,j}$.
}
Let $\varepsilon_{i,j} > 0$ be such that
for every $u \in B(x^i_j,\varepsilon_{i,j})$,
$\norm{\tau(u) - \tau(x^i_j)} < \frac{1}{j + 1}$.
Choose $\bars_{i,j},\bart_{i,j}$ such that for every $j$,
$s_{i,j} > \bars_{i,j} > \bart_{i,j + 1} > t_{i,j + 1}$.
We may also assume that for every distinct $i$ and $i'$,
$d(B^{\oversE}(\barz_i,\bart_{i,0}),
B^{\oversE}(\barz_{i'},\bart_{i',0})) > 0$
and that
$B^{\oversE}(\barz_0,\bart_{0,0}) \subseteq \rfs{cl}^{\oversE}(S)$.

By Proposition \ref{metr-bldr-p4.7}(c), for every $i,j$
there is $h_{i,j} \in \rfs{UC}(X;D)$
such that
$\rfs{supp}(h_{i,j}) \subseteq
B^{\oversE}(\barz_i,\bart_i) - \overB^{\oversE}(\barz_i,\bars_i)$,
$h_{i,j}(x^i_j) = x^i_j$ and
$h_{i,j}(y^i_j) \in B(x^i_j,\varepsilon_{i,j})$.
Let $h_i = \bcirc_{j \in \sboldbbN} h_{i,j}$.
By Proposition \ref{metr-bldr-p4.5}, $h_i \in \rfs{UC}(X)$.
So $h_i \in \rfs{UC}(X;D)$.
Also, $\rfs{supp}(h_i) \subseteq S$.
Let $h = \bcirc_{i \in \sboldbbN} h_i$.
Applying again Proposition \ref{metr-bldr-p4.5}, we conclude that
$h \in \rfs{UC}(X;S,D,x^*)$.

We check that $h^{\tau}$ is not bi-UC at $\tau(x^*)$.
Let $V \in \rfs{Nbr}^Y(\tau(x^*))$.
For some~$i$,
$\rfs{supp}((h_i)^{\tau}) \subseteq V$.
Denote $u^i_j = \tau(x^i_j)$ and $v^i_j = \tau(y^i_j)$.
So
\begin{itemize}
\addtolength{\parskip}{-11pt}
\addtolength{\itemsep}{06pt}
\item[(1)] 
for every $j$, $\norm{u^i_j - v^i_j} > d_i$.
\vspace{-05.7pt}
\end{itemize}
Since $h_i(y^i_j) \in B(x^i_j,\varepsilon_{i,j})$,
it follows that
$\lim_{j \rightarrow \infty} \norm{\tau(x^i_j) - \tau(h_i(y^i_j)} = 0$.
That is,\break
$\lim_{j \rightarrow \infty}
\norm{(h_i)^{\tau}(u^i_j) - (h_i)^{\tau}(v^i_j)} = 0$.
Hence
\begin{itemize}
\addtolength{\parskip}{-11pt}
\addtolength{\itemsep}{06pt}
\item[(2)] 
$\lim_{j \rightarrow \infty}
\norm{h^{\tau}(u^i_j) - h^{\tau}(v^i_j)} = 0$.
\vspace{-05.7pt}
\end{itemize}
(1) and (2) imply that $h^{\tau} \nrestriction V$ is not bi-UC.
That is, $h^{\tau} \not\in \rfs{PNT.UC}(Y,\tau(x^*))$.
A contradiction.
We have reached a contradiction in both Case 1 and Case 2.
So $\tau$ is UC at~$x^*$.

{\bf Part 2 } $\tau\inverse$ is UC at $\tau(x^*)$.
\newline
Suppose by contradiction that
this is not true.
So for every $V \in \rfs{Nbr}^Y(\tau(x^*))$,
$\tau\inverse \nrestriction V$ is not UC.

{\bf Claim 2.}
For every $k \in \bbN$ there are positive numbers
$r_1^k,\ldots,r_5^k$ and sequences $\vec{x}^k$ and $\vec{y}^k$
which fulfill the following requirements.
\begin{itemize}
\addtolength{\parskip}{-11pt}
\addtolength{\itemsep}{06pt}
\item[(1)] 
$r_1^k > r_2^k \geq r_3^k > r_4^k > r_5^k = 2  r_1^{k + 1}$.
\item[(2)] 
$\limti{i} \norm{x_i^k - x^*} = r_2^k$
and $\limti{i} \norm{y_i^k - x^*} = r_3^k$.
\item[(3)] 
There is $e_k > 0$ such that $\norm{x^k_i - y^k_i} > e_k$
for every $i \in \bbN$.
\item[(4)] 
$\rfs{Rng}(\vecx^k) \cup \rfs{Rng}(\vecy^k) \subseteq D$.
\item[(5)] 
Define
$s_k = \sup(\setm{\norm{\tau(x) - \tau(x^*)}}{x \in B(x^*,r_4^k)})$
and
$t_k = \norm{\tau(x^*) - \tau(\vecx^k)}$.
Then $s_k < t_k$.
\item[(6)] 
$\limti{i} \norm{\tau(x_i^k) - \tau(y_i^k)} = 0$.
\item[(7)] 
Either $\vecx^k$ is a Cauchy sequence or $\vecx^k$ is spaced,
and either $\vecy^k$ is a Cauchy sequence or $\vecy^k$ is spaced.
\vspace{-05.7pt}
\end{itemize}

{\bf Proof }
Let $r_1^0 > 0$ be such that $\overB(x^*,r_1^0) \subseteq S$.
Suppose that $r_1^k$ has been defined,
and we define $r^k_2,\ldots,r^k_5$ and $r^{k + 1}_1$.
Let $r = \dgfrac{r_1^k}{2}$.
Since $\tau\inverse \nrestriction \tau(B(x^*,r))$
is not uniformly continuous, there are $e_k > 0$
and sequences $\vecx,\vecy \subseteq B(x^*,r)$
such that for every $i \in \bbN$,
$\norm{x_i - y_i} > e_k$ and
$\limti{i} \norm{\tau(x_i) - \tau(y_i)} = 0$.
Since $D \cap S$ is dense in $S$, we may assume that
$\rfs{Rng}(\vecx) \cup \rfs{Rng}(\vecy) \subseteq D$.
We may also assume that
$x^* \not\in \rfs{Rng}(\vecx) \cup \rfs{Rng}(\vecy)$.

By interchanging some $x_i$'s with their corresponding $y_i$'s,
we may assume that\break
$\norm{x_i - x^*} \geq \norm{y_i - x^*}$.
Taking subsequences we may assume that
$r_2^k \eqdf \limti{i}\norm{x_i - x^*}$ and
$r_3^k \eqdf \limti{i}\norm{y_i - x^*}$ exist.
Hence $r_3^k \leq r_2^k$.
Taking subsequences again, we may assume that 
either $\vecx$ is a Cauchy sequence or $\vecx$ is spaced, and that
either $\vecy$ is a Cauchy sequence or $\vecy$ is spaced.

Note that $\vecx$ does not contain a convergent subsequence, since if
$x'$ is a limit of a subsequence of $\vecx$,
then $\tau\inverse$ is not continuous at $\tau(x')$.
Also recall that $x^* \not\in \rfs{Rng}(\vecx)$. It thus follows that
$t_k \eqdf \norm{\tau(x^*),\tau(\vecx^k)} > 0$.
Next define $\vecx^k = \vecx$ and $\vecy^k = \vecy$.
Let $r_4^k < r_3^k$ be such that
$s_k \eqdf
\sup(\setm{\norm{\tau(x) - \tau(x^*)}}{x \in B(x^*,r_4^k)}) < t_k$.
Finally, let $r_5^k = \dgfrac{r_4^k}{2}$
and $r^{k + 1}_1 = \dgfrac{r^k_5}{2}$.
This concludes the construction which proves Claim 2.
\smallskip

Since $\lim_{i \rightarrow \infty} \norm{x^k_i} = r_2^k$
and $\lim_{i \rightarrow \infty} \norm{y^k_i} = r_3^k$,
we may assume that
\begin{itemize}
\addtolength{\parskip}{-11pt}
\addtolength{\itemsep}{06pt}
\item[(8)] 
for every \hbox{$i \in \bbN$,}
$r_4^k < \norm{x^k_i - x^*} < \dgfrac{(r_2^k + r_1^k)}{2}$ and
$r_4^k < \norm{y^k_i - x^*} < \dgfrac{(r_2^k + r_1^k)}{2}$.
\vspace{-05.7pt}
\end{itemize}
We may also assume that either for every $k \in \bbN$, $\vecy^k$ is
spaced, or for every $k \in \bbN$, $\vecy^k$ is a Cauchy sequence.

{\bf Case 1 } For every $k \in \bbN$, $\vecy^k$ is spaced.
Fix $k \in \bbN$ and denote $r_i^k$,
$\vecx^k$, $\vecy^k$ and $e_k$ by $r_i$, $\vecx$, $\vecy$ and $e$
respectively.

{\bf Claim 3.} There are subsequences $\setm{x_{i_n}}{n \in \bbN}$
$\setm{y_{i_n}}{n \in \bbN}$
of $\vecx$ and $\vecy$ respectively, such that
$d(\setm{x_{i_n}}{n \in \bbN},\setm{y_{i_n}}{n \in \bbN}) > 0$.

{\bf Proof } 
The claim is trivial if $\vecx$ is a Cauchy sequence.
So suppose $\vecx$ is spaced.
We show that there is a sequence $\setm{i_n}{n \in \bbN}$ such that
(i) $\lim_{n > m \raro \infty} \dline{x_{i_m} - y_{i_n}}$ exists, and 
(ii) $\lim_{n > m \raro \infty} \dline{y_{i_m} - x_{i_n}}$ exists.
By repeatedly applying Ramsey Theorem, we obtain a
decreasing sequence $A_0 \supseteq A_1 \supseteq A_2 \ldots$ of 
infinite subsets of $\bbN$ such that for every $\ell \in \bbN$ and
$m,n,m',n' \in A_{\ell}$: if $m < n$ and $m' < n'$,
then $|\norm{x_m - y_n} - \norm{x_{m'} - y_{n'}}| < 2^{-\ell}$.
Let $\setm{i_n}{n \in \bbN}$ be $\onetoonen$ sequence such that
for every $n \in \bbN$, $i_n \in A_n$.
Then (i) holds for $\setm{i_n}{n \in \bbN}$. 
The same argument is applied to show that (ii) holds.

Let
$\bars_1 = \lim_{n > m \raro \infty} \norm{x_{i_m} - y_{i_n}}$ and 
$\bars_2 = \lim_{n > m \raro \infty} \norm{y_{i_m} - x_{i_n}}$.
It is easy to see that if $\bars_1 = 0$ or $\bars_2 = 0$,
then $\vecx$ is a Cauchy sequence. So $\bars_1,\bars_2 > 0$.
By removing an initial segment from the sequences
$\sngltn{x_{i_n}}_{n \in \bbN}$ and $\sngltn{y_{i_n}}_{n \in \bbN}$
we may assume that for every $n > m$,
$\norm{x_{i_m} - y_{i_n}} > \dgfrac{\bars_1}{2}$ and
$\norm{x_{i_n} - y_{i_m}} > \dgfrac{\bars_2}{2}$.
Recall also that $\norm{x_i - y_i} > e$ for every $i \in \bbN$.
So 
$d(\setm{x_{i_n}}{n \in \bbN},\setm{y_{i_n}}{n \in \bbN}) \geq
\min(\dgfrac{\bars_1}{2},\dgfrac{\bars_2}{2},e)$.
So Claim 3 is proved.\smallskip

We may thus assume that
$d_k  \eqdf d(\rfs{Rng}(\vecx^k),\rfs{Rng}(\vecy^k)) > 0$.

{\bf Claim 4.}
For every $k \in \bbN$ there is $h_k \in \rfs{LIP}(X;D)$
with the following properties:\break
(i) $\rfs{supp}(h_k) \subseteq B(x^*,r_1^k) - B(x^*,r_5^k)$; and
(ii) there is $n_k \in \bbN$ such that for every $i \geq n_k$,
$h_k(x_i^k) = x_i^k$ and $h_k(y_i^k) \in B(x^*,r_4^k)$.

{\bf Proof }
Fix $k$, for $j = 1,\ldots,5$ set
$r^k_j = r_j$, write $\vecx^k = \vecx$, $\vecy^k = \vecy$,
$x^k_i = x_i$, $y^k_i = y_i$
and define $w_i = x_i - x^*$, $z_i = y_i - x^*$
and $u_i = \dgfrac{z_i}{\norm{z_i}}$.
Note that $\lim_{i \in \bbN} \norm{(x^* + r_3 u_i) - y_i} = 0$,
and recall that $d(\rfs{Rng}(\vecx),\rfs{Rng}(\vecy)) > 0$.
From these facts it follows that by removing an initial segment of
$\vecx$ and of $\vecy$,
we may assume that there is $a > 0$
such that
$\norm{x_i - (x^* + r_3 u_j)} \geq a$ for every $i,j \in \bbN$.
Similarly, since $\vecy$ is spaced,
we may assume that $\sngltn{x^* + r_3 u_i}_{i \in \bbN}$ is spaced too.
Certainly we may choose $a$ to be smaller than
$r_3 - r_4$ and $r_1 - r_3$,
and we may assume that for every $i$,
$\norm{w_i} \geq r_3 - \dgfrac{a}{8}$
and $r_3 - \dgfrac{a}{4} < \norm{z_i} < r_3 + \dgfrac{a}{4}$.
Let $L_i = [x^* + r_4 u_i,x^* + (r_3 + \dgfrac{a}{4}) u_i]$.
Note that $y_i \in L_i$.
We show that for every $i,j$, \ $d(x_i,L_j) > \dgfrac{a}{4}$.
Let $y \in L_j$.
If
$y \in [x^* + (r_3 - \dgfrac{a}{2})u_j,x^* + (r_3 + \dgfrac{a}{4})u_j]$,
then
$$
\norm{x_i - y} \geq
\norm{x_i - (x^* + r_3 u_j)} - \norm{(x^* + r_3 u_j) - y} \geq
a - \dgfrac{a}{2} = \dgfrac{a}{2},
$$
and if $y \in [x^*,x^* + (r_3 - \dgfrac{a}{2})u_i]$,
then 
$$
\norm{x_i - y} \geq \norm{x_i - x^*} - \norm{y - x^*} \geq
r_3 - \dgfrac{a}{8} - (r_3 - \dgfrac{a}{2}) = \dgfrac{3a}{8}.
$$
It follows that $d(x_i,L_j) > \dgfrac{a}{4}$.

Let $v_i = x^* + r_4 u_i$,
and let $b > 0$ be such that for every $i \neq j$,
$\norm{v_i - v_j} > b$.
We show that if $i \neq j$, then $d(L_i,L_j) \geq \dgfrac{b}{2}$.
It is easy to see that
$d(L_i,L_j) = d(v_i,L_j)$.
Let $x^* + tu_j \in L_j$. If $t \in [r_4,r_4 + \dgfrac{b}{2}]$,
then
$$
\norm{v_i - (x^* + tu_j)} \geq
\norm{v_i - v_j} - \norm{x^* + tu_j - v_j } > b - \dgfrac{b}{2} =
\dgfrac{b}{2}.
$$
If $t > r_4 + \dgfrac{b}{2}$, then
$$
\norm{v_i - (x^* + tu_j)} \geq \norm{tu_j } - \norm{v_i - x^*} >
r_4 + \dgfrac{b}{2} - r_4 = \dgfrac{b}{2}.
$$

It follows that there is $d > 0$ such that:
\begin{itemize}
\addtolength{\parskip}{-11pt}
\addtolength{\itemsep}{06pt}
\item[(1)] 
for every $i \neq j$, \ $2d < d(L_i,L_j)$;
\item[(2)] 
for every $i \neq j$, \ $d < d(x_i,L_j)$;
\item[(3)] 
$r_3 + \dgfrac{a}{4} + d < r_1$;
\item[(4)] 
$r_4 - d > r_5$.
\vspace{-05.7pt}
\end{itemize}

Let $L_i^1 = [v_i,y_i]$. So $L_i^1 \subseteq L_i$.
Hence
\begin{list}{}
{\setlength{\leftmargin}{39pt}
\setlength{\labelsep}{08pt}
\setlength{\labelwidth}{20pt}
\setlength{\itemindent}{-00pt}
\addtolength{\topsep}{-04pt}
\addtolength{\parskip}{-02pt}
\addtolength{\itemsep}{-05pt}
}
\item[(1.1)]
for every $i \neq j$, \ $2d < d(L_i^1,L_j^1)$;
\item[(1.2)]
for every $i \neq j$, \ $d < d(x_i,L_j^1)$;
\item[(1.3)]
$\norm{y_i - v_i} < r_3 - r_4 + \dgfrac{a}{4}$.
\vspace{-02.0pt}
\end{list}
By (3),
$d(B(L_i^1,d),X - B(x^*,r_1)) > r_1 - (r_3 + \dgfrac{a}{4} + d) > 0$
and by (4), $d(B(L_i^1,d),B(x^*,r_5)) > r_4 - r_5 - d > 0$. So
\begin{list}{}
{\setlength{\leftmargin}{39pt}
\setlength{\labelsep}{08pt}
\setlength{\labelwidth}{20pt}
\setlength{\itemindent}{-00pt}
\addtolength{\topsep}{-04pt}
\addtolength{\parskip}{-02pt}
\addtolength{\itemsep}{-05pt}
}
\item[(1.4)]
$d(B(L_i^1,d),X - (B(x^*,r_1) - B(x^*,r_5))) > 0$
for every $i \in \bbN$.
\vspace{-02.0pt}
\end{list}
Recall that $y_i \in D$, but $v_i$ need not be in $D$.
For every $i$,
choose $v'_i \in D$ sufficiently close to $v_i$ and define
$L'_i = [v'_i,y_i]$.
This can be done in such a way that $L'_i$ satisfy (1.1)\,-\,(1.4).
So indeed choose $v'_i \in D \cap B(x^*,r_4)$
in such a way that
the $L'_i$'s fulfill (1.1)\,-\,(1.4).
Write $v_{k,i} = v'_i$.

Let $K = K_{\srfs{seg}}(r_3 - r_4 + \dgfrac{a}{4},d)$ be as in
\ref{l2.6}(c) and $i \in \bbN$.
By \ref{l2.6}(c), there is
$h'_i \in \rfs{LIP}(X;D)$ such that:
$\rfs{supp}(h'_i) \subseteq B(L'_i,d)$,
$h'_i$ is $K$-bilipschitz,
and $h'_i(y_i) = v'_i$.
Since the $L'_i$'s satisfy (1.1), it follows that
for every $i \neq j$,
$d(\rfs{supp}(h'_i),\rfs{supp}(h'_j)) > 0$.
So $h_k \eqdf \bcirc_{i \in \sboldbbN} h'_i$ is well defined.
Also, $h_k$ is $2K$-bilipschitz.

For every $i$,
$h_k(y_i) = h'_i(y_i) = v'_i \in B(x^*,r_4)$.
By (1.2) applied to the $L'_j$'s,
$x_i \not\in \rfs{supp}(h_k)$. So $h_k(x_i) = x_i$.
By (1.4) applied to $L'_i$, for every $i$,
$\rfs{supp}(h'_j) \subseteq B(x^*,r_1) - B(x^*,r_5)$.
So
$\rfs{supp}(h_k) \subseteq B(x^*,r_1) - B(x^*,r_5)$.
Recall that for every $i$, $h'_i \in H(X;D)$. So $h_k \in H(X;D)$.
We have shown that $h_k$ fulfills the requirements of Claim 4.
\smallskip

Let $h = \bcirc_{k \in \sboldbbN} h_k$.
By Proposition \ref{metr-bldr-p4.5},
$h \in \rfs{UC}(X)$. Since $B(x^*,r^0_1) \subseteq S$,
we obtain that $\rfs{supp}(h) \subseteq S$, and since for every $k$,
$h_k \in H(X;D)$, we conclude that $h \in H(X;D)$.
Also for every $k$, $x^* \not\in \rfs{supp}(h_k)$.
So $h(x^*) = x^*$, that is, $h \in \rfs{UC}(X;S,D,x^*)$.

We shall reach a contradiction by showing that
$h^{\tau} \not\in \rfs{PNT.UC}(Y \tau(x^*))$.
Let $V \in \rfs{Nbr}^Y(\tau(x^*))$.
\smallskip Let $k$ be such that $\tau(B(x^*,r_1^k)) \subseteq V$.
Hence
\begin{itemize}
\addtolength{\parskip}{-11pt}
\addtolength{\itemsep}{06pt}
\item[(i)]
for every~$i \in \bbN$,
$\tau(x_i^k),\tau(y_i^k) \in V$, and
$\limti{i} \norm{\tau(x_i^k) - \tau(y_i^k)} = 0$.
\vspace{-05.7pt}
\end{itemize}
$h^\tau(\tau(x_i^k)) = \tau(x_i^k)$ and
$h^\tau(\tau(y_i^k)) = \tau(h(y_i^k)) \in \tau(B(x^*,r_4^k))$.
So for every $i \in \bbN$,
\smallskip
\newline
$(\dagger)$
\kern-1mm
\renewcommand{\arraystretch}{1.5}
\addtolength{\arraycolsep}{-0pt}
$
\begin{array}[t]{ll}
&\norm{(h^\tau(\tau(x_i^k)) - \tau(x^*))
- (h^\tau(\tau(y_i^k)) - \tau(x^*))}
\\
 = &\norm{(\tau(x_i^k) - \tau(x^*)) - (\tau(h(y_i^k)) - \tau(x^*))} \geq
\norm{\tau(x_i^k) - \tau(x^*)} -
\norm{\tau(h(y_i^k)) - \tau(x^*)}.
\end{array}
$
\renewcommand{\arraystretch}{1.0}
\addtolength{\arraycolsep}{0pt}
\smallskip\newline
Recall that $h(y_i^k) = v_{k,i} \in B(x^*,r^k_4)$.
Let $s_k,t_k$ be as in Clause (5) of Claim 2.
Then $\norm{\tau(h(y_i^k)) - \tau(x^*)} \leq s_k$
and $\norm{\tau(x_i^k) - \tau(x^*)} \geq t_k$.
Denote the right handside of $(\dagger)$ by $A$.
So $A \geq t_k - s_k$.
By Clause (5) in Claim 2, $t_k - s_k > 0$.
We have proved that
\begin{itemize}
\addtolength{\parskip}{-11pt}
\addtolength{\itemsep}{06pt}
\item[(ii)] for every $i \in \bbN$,
$\norm{h^\tau(\tau(x_i^k)) - h^\tau(\tau(y_i^k))} \geq t_k - s_k > 0$.
\vspace{-05.7pt}
\end{itemize}
(i) and (ii) demonstrate that $h^{\tau} \nrestriction V$ is not bi-UC.
We have shown that for every\break
$V \in \rfs{Nbr}(\tau(x^*))$,
$h^{\tau} \nrestriction V$ is not UC. That is,
$h^{\tau} \not\in \rfs{PNT.UC}(Y,\tau(x^*))$.
A contradiction.

{\bf Case 2 } For every $k \in \bbN$, $\vecy^k$ is a Cauchy sequence.

{\bf Claim 5.}
For every $k \in \bbN$
there is $h_k \in \rfs{LIP}(X;D)$ with the
following properties:
(i) $\rfs{supp}(h_k) \subseteq B(x^*,r_1^k) - B(x^*,r_5^k)$; and
(ii) there is $n_k \in \bbN$ such that for every $i \geq n_k$,
$h_k(x_i^k) = x_i^k$ and $h_k(y_i^k) \in B(x^*,r_4^k)$.

{\bf Proof }
Fix $k$, and denote $\vecx^k$, $\vecx^k$, $r^k_j$ etc.\ by
$\vecx$, $\vecy$, $r_j$ etc..
Let $\bary = \lim^{\oversE} \vecy$. Since $\tau\inverse$ is continuous,
$\bary \in \rfs{cl}^{\oversE}(S) - S$.
Also, $\norm{\bary - x^*} = r_3$.
Since $\bary \not\in E$ and $\rfs{Rng}(\vecx) \subseteq E$,
$\rfs{Rng}(\vecx) \cap [x^*,\bary]$ contains at most one element.
By removing this element
we may assume that $\hate \eqdf d(\rfs{Rng}(\vecx),[x^*,\bary]) > 0$.
Let $b = \dgfrac{(r_4 + r_5)}{2}$,
$a = \dgfrac{(r_4 - r_5)}{2}$ and $c = \min(a,\hate,r_1 - r_3)$.
Let $w \in [x^*,\bary]$ be such that $\norm{w - x^*} = b$.
Let $u,v \in D$ be such that
$\norm{u - \bary},\norm{v - w} < \dgfrac{c}{12}$.
By Lemma \ref{l2.6}(c), there is $h \in \rfs{LIP}(X;D)$ such that
$\rfs{supp}(h) \subseteq B([u,v],\dgfrac{c}{4})$,
$h(u) = v$ and $h(B(u,\dgfrac{c}{12})) = B(v,\dgfrac{c}{12})$.
Since
$h$ is bilipschitz, $\rfs{Dom}(h^{\srfs{cl}}) = \rfs{cl}^{\oversE}(X)$.
Denote $\hath = h^{\srfs{cl}}$.
We show that $\hath(\bary) \in B^{\oversE}(x^*,r_4)$.
Since $\bary \in B^{\oversE}(u,\dgfrac{c}{12})$,
$\hath(\bary) \in B^{\oversE}(v,\dgfrac{c}{12})$.
So
$$
\norm{\hath(\bary) - x^*} \leq
\norm{\hath(\bary) - v} + \norm{v - w} + \norm{w - x^*} <
\dgfrac{c}{12} + \dgfrac{c}{12} + b \leq b + \dgfrac{a}{6} <
b + a = r_4.
$$
It follows that
\begin{itemize}
\addtolength{\parskip}{-11pt}
\addtolength{\itemsep}{06pt}
\item[(1)] 
for all but finitely many $i$'s, $h(y_i) \in B(x^*,r_4)$.
\vspace{-05.7pt}
\end{itemize}
For every $i$,
$d(x_i,[u,v]) \geq
d(x_i,[\bary,w]) - (\dgfrac{c}{12} + \dgfrac{c}{12}) \geq
\hate - \dgfrac{c}{6} \geq \dgfrac{c}{4}$.
So $x_i \not\in \rfs{supp}(h)$ and hence
\begin{itemize}
\addtolength{\parskip}{-11pt}
\addtolength{\itemsep}{06pt}
\item[(2)] 
$h(x_i) = x_i$ for all $i \in \bbN$.
\vspace{-05.7pt}
\end{itemize}
$\norm{u - x^*} \leq \dgfrac{c}{12} + r_3 < r_1 - \dgfrac{c}{4}$.
It easily follows that
$B([u,v],\dgfrac{c}{4}) \subseteq B(x^*,r_1)$.\break
$\norm{v - x^*} \geq b - \dgfrac{c}{12} > r_5 + \dgfrac{a}{4}$.
Next we have that
$$
d(B([u,v],\dgfrac{c}{4}),x^*) \geq
d(B([\bary,w],\dgfrac{c}{4}),x^*) - \dgfrac{c}{6} - \dgfrac{c}{4} =
b - \dgfrac{5c}{12} > r_5.
$$
So $B([u,v],\dgfrac{c}{4}) \cap B(x^*,r_5) = \emptyset$.
Similarly, for every $y \in B([u,v],\dgfrac{c}{4})$,
$$
\norm{y} \leq \max(\norm{u},\norm{v}) + \dgfrac{c}{4} \leq
\max(\norm{\bary},\norm{w}) + \dgfrac{c}{12} + \dgfrac{c}{4} =
r_3 + \dgfrac{5c}{12} < r_1.
$$
That is, $\rfs{supp}(h) \subseteq B(x^*,r_1)$.
So
\begin{itemize}
\addtolength{\parskip}{-11pt}
\addtolength{\itemsep}{06pt}
\item[(3)] 
$\rfs{supp}(h) \subseteq B(x^*,r_1) - B(x^*,r_5)$.
\vspace{-05.7pt}
\end{itemize}
It follows that $h_k \eqdf h$ fulfills the requirements
of Claim 5. So Claim 5 is proved.\smallskip

The remaining steps in the proof are identical to those in Case 1. 
So both Case 1 and Case 2 lead to a contradiction.
This means that $\tau\inverse$ is UC at $\tau(x^*)$.
\hfill\myqed

\begin{question}\label{metr-bldr-q4.9}
\begin{rm}
Let $X,Y$ be open subsets of the normed spaces $E$ and $F$
and $\tau \in H(X,Y)$ be such that
$(\rfs{LIP}(X))^{\tau} \subseteq \rfs{LUC}(Y)$.
Is $\tau$ locally UC? Is $\tau\inverse$ locally UC?

Note that by Theorem \ref{metr-bldr-t3.27}, the answer to both
parts of the question is positive for $E$'s such that
$\kappa(E) \geq 2^{\aleph_0}$.
Hence the answer is positive for open subsets of $\ell_{\infty}$.
\end{rm}
\end{question}

\noindent
{\bf Proof of Theorem \ref{metr-bldr-t4.1}}

(a) Let $X,Y$, $\itGamma,\itDelta$ and $\varphi$ be as in Part (a).
Suppose that $\itGamma$ is $(\leq\kern-1pt \kappa(E))$-generated.
Then by Theorem \ref{metr-bldr-t3.42}, $\itGamma = \itDelta$
and there is $\tau \in H(X,Y)$ as required.

Note that for every metric space $X$,
$\rfs{LUC}(X) = H_{\srfs{MC}}^{\srfs{LC}}(X)$.

Suppose that $\itGamma = \rfs{MC}$. By Theorem \ref{t2.4}(a), there is
$\tau \in H(X,Y)$ such that $\tau$ induces~$\varphi$.
$(\rfs{UC}(X))^{\tau} \subseteq \rfs{LUC}(Y)$.
So by Theorem \ref{metr-bldr-t4.8}(a),
$\tau$ is locally bi-UC.
So $(\rfs{LUC}(X))^{\tau} = \rfs{LUC}(Y)$.
Hence $H_{\srfs{MC}}^{\srfs{LC}}(X) = H_{\itDelta}^{\srfs{LC}}(Y)$.
We have seen that the above equality implies that
$\rfs{MC} = \itDelta$. So Part (a) is proved.

(b)
Let $\fourtpl{E}{X}{\calS}{\calE}$, $\fourtpl{F}{Y}{\calT}{\calF}$,
$\itGamma,\itDelta$ be and $\varphi$ be as in Part (b).
If both $\itGamma$ and $\itDelta$ are
$(\leq\kern-1pt \kappa(E))$-generated,
then by Theorem \ref{metr-bldr-t3.42}, $\itGamma = \itDelta$,
and there is $\tau$ which induces $\varphi$.

Suppose that $\itDelta$ or $\itGamma$ are not
$(\leq\kern-1pt \kappa(E))$-generated.
By Theorem \ref{t2.4}(a), there is\break
$\tau \in H(X,Y)$
such that $\tau$ induces $\varphi$.

Suppose by contradiction that
$\itGamma = \rfs{MC}$ and $\itDelta \neq \rfs{MC}$.
Hence $\itDelta$ is $(\leq\kern-1pt \kappa(E))$-generated.
Clearly,
$(\rfs{LIP}(X;\calS,\calE))^{\tau} \subseteq
H_{\itDelta}^{\srfs{LC}}(Y)$.
By Theorem \ref{metr-bldr-t3.27}, $\tau$ is
locally $\itDelta$-bicontinuous.
Hence 
$(H_{\itDelta}^{\srfs{LC}}(Y;\calT,\calF))^{\tau\inverse} \subseteq
H_{\itDelta}^{\srfs{LC}}(X)$.
However,
$(H_{\itDelta}^{\srfs{LC}}(Y;\calT,\calF))^{\tau\inverse} =
H_{\srfs{MC}}^{\srfs{LC}}(X;\calS,\calE)$.
Hence
$H_{\srfs{MC}}^{\srfs{LC}}(X;\calS,\calE) \subseteq
H_{\itDelta}^{\srfs{LC}}(X)$.
A contradiction.
It follows that $\itGamma = \itDelta = \rfs{MC}$.
\smallskip\hfill\myqed

As in Chapter \ref{s3}, the analogous statement for manifolds is
also true.

\begin{cor}\label{metr-bldr-c4.10}
Let $\pair{X}{\itPhi}$ and $\pair{Y}{\itPsi}$
be normed manifolds with locally Lipschitz atlasses.
Let $\itGamma$ and $\itDelta$ be moduli of continuity,
Suppose that $\itGamma$ is countably generated or $\itGamma = \rfs{MC}$,
and the same holds for $\itDelta$.

\num{a}
If
$\iso{\varphi}{H_{\itGamma}^{\srfs{LC}}(X,\itPhi)}
{H_{\itDelta}^{\srfs{LC}}(Y)}$.
Then $\itGamma = \itDelta$
and there is $\iso{\tau}{X}{Y}$ such that $\tau$ induces $\varphi$,
and $\tau$ is locally $\itGamma$-bicontinuous.

\num{b} Let $\calS$ be an open cover of $X$,
$\calT$ be an open cover of\kern1pt\ $Y$
and
$\varphi\,:\,H_{\itGamma}^{\srfs{LC}}(X,\itPhi,\calS) \cong
H_{\itDelta}^{\srfs{LC}}(Y,\itPsi,\calT)$.
Then $\itGamma = \itDelta$,
there is $\iso{\tau}{X}{Y}$ such that
$\varphi(h) = h^{\tau}$ for every
$h \in H_{\itGamma}^{\srfs{LC}}(X;\calS,\calE)$,
and $\tau$ is locally $\itGamma$-bicontinuous.
\end{cor}

\newpage

\section{Other groups defined by properties related to
uniform continuity}\label{s5}

\subsection{General description.}
\label{ss5.1}

The results we have obtained on groups of type
$H_{\itGamma}^{\srfs{LC}}(X)$
are more comprehensive than those obtained for other types
of groups.
We have presented the results on $H_{\itGamma}^{\srfs{LC}}(X)$ in the
quite general framework of ``subspace choice systems''.
We now abandon this framework, and restrict the discussion
to the class of open subsets of normed spaces.

Recall the following notations which were introduced in the
introduction.

\begin{defn}\label{d4.1}\label{metr-bldr-d5.1}
\begin{rm}
(a) For a set $F$ of $\onetoonen$ functions let
$F\inverse = \setm{f\inverse}{f \in F}$.
Suppose that $\calP$ is used as an abbreviation for some property
of maps,
and let $X$ and $Y$ be topological spaces.
We shall use the notation
$\calP(X,Y)$ to denote the set of all homeomorphisms between
$X$ and $Y$ which have property~$\calP$. We denote
$$
\calP^{\pm}(X,Y) \eqdf \calP(X,Y) \cap (\calP(Y,X))\inverse
\mbox{\ \,and\ \ \ } \calP(X) \eqdf \calP^{\pm}(X,X).
$$
Usually but not always this convention will be used for
$\calP$'s which are ``closed under composition''.
($\calP$ is closed under composition, if for every
$\fnn{f}{X}{Y}$ and $\fnn{g}{Y}{Z}$: if $f$ and $g$ fulfill $\calP$,
then $g \scirc f$ fulfills $\calP$).
In such cases $\calP(X)$ is a group.

(b) Let $\pair{X}{d}$ be a metric space.
$X$ is
{\it uniformly\,-\,in\,-\,diameter arcwise\,-\,connected (UD.AC)},
if for every $\varepsilon > 0$ there is $\delta > 0$ such that
for every $ x,y \in X$:
if $d(x,y) < \delta$, then there is an arc $L \subseteq X$ connecting
$x$ and $y$ such that $\rfs{diam}(L) < \varepsilon$.

   \index{uniformly\,-\,in\,-\,diameter arcwise\,-\,connected}
   \index{udac@@UD.AC. Abbreviation of uniformly\,-\,in\,-\,diameter
          arcwise\,-\,connected}

(c)
Let $K^{\calO}_{\srfs{NRM}}$ be the class of all spaces $X$ such that
$X$ is an open subset of a normed space.
Let $K^{\calO}_{\srfs{BNC}}$ be the class of all spaces $X$ such that
$X$ is an open subset of a Banach space.
Let $K^{\calO}_{\srfs{NFCB}}$ be the class of all spaces $X$ such that
$X$ is an open subset of a normed space of the first category,
or $X$ is an open subset of a Banach space.
   \index{N@konrm@@$K^{\calO}_{\srfs{NRM}}$.
          Class of all spaces which are an open subset of
	  a normed space}
   \index{N@kobnc@@$K^{\calO}_{\srfs{BNC}}$.
          Class of all spaces which are an open subset of a
	  Banach space}
   \index{N@konfcb@@$K^{\calO}_{\srfs{NFCB}}$.
          Class of open subsets of first category or complete
	  normed spaces}
\smallskip\hfill\proofend
\end{rm}
\end{defn}

Note that a disconnected space may be UD.AC.
The space $[0,1] \cup [2,3]$ is such an example.

The following statement is a typical example of
some of the final results obtained in this chapter.
It is restated in Corollary \ref{metr-bldr-c5.7}.

\kern1.0mm

\begin{it}
\noindent
{\bf Theorem A.\kern2pt}
Let $X,Y \in K_{\srfs{NFCB}}^{\calO}$.
Suppose that $X$ and $Y$ are UD.AC spaces.
Let
\newline
$\iso{\varphi}{\rfs{UC}(X)}{\rfs{UC}(Y)}$.
Then there is
$\tau \in \rfs{UC}^{\pm}(X,Y)$ such that $\tau$ induces~$\varphi$.
\vspace{0.8mm}
\end{it}

The reason that Theorem A can be proved only for members
of $K_{\srfs{NFCB}}^{\calO}$ and not for all members of
$K_{\srfs{NRM}}^{\calO}$
is that Theorem \ref{metr-bldr-t2.8} cannot be used.
This is so, since in Theorem \ref{metr-bldr-t2.8} we need to know
that $\rfs{LIP}^{\srfs{LC}}(X) \leq G$.
However, $\rfs{LIP}^{\srfs{LC}}(X) \not\leq \rfs{UC}(X)$.

Theorem A assumes that the open sets $X$ and $Y$
are UD.AC.
Different extra assumptions on the open sets in question are often used
in proving other reconstruction results.
We make sure though, that these extra assumptions do not exclude
the known well-behaved open subsets of a normed space.
For example, convex bounded open sets are always included.
Usually the classes for which reconstruction is proved do contain
\hbox{some complicated open sets.}

Theorem A has the following corollary.

\begin{theorem}\label{t4.4}\label{metr-bldr-t5.2}
Let $F$ and $K$ be the closures of UD.AC bounded open subsets
of $\bbR^m$ and $\bbR^n$ respectively.
Let $\iso{\varphi}{H(F)}{H(K)}$. Then $\varphi$ is induced by a
homeomorphism between $F$ and $K$.
\end{theorem}


The proof of Theorem \ref{t4.4} appears after Example \ref{e4.11}.
The boundedness of $F$ and $K$ above is necessary,
see Example \ref{e4.13}.
The analogue of Theorem \ref{metr-bldr-t5.2}
for open subsets of infinite-dimensional normed spaces is proved
in \ref{t4.35}.
The boundedness of $F$ and $K$ is not required in the
infinite-dimensional case.

Let us point out that the closure of a UD.AC open subset of
$\bbR^n$ does not have to be a Euclidean manifod with boundary,
neither does it have to be a polyhedron.
The reconstruction theorems for polyhedra and for Euclidean manifolds
with boundary were proved in \cite{Ru1} 3.34 and 3.43.
Theorem \ref{metr-bldr-t5.2} is not a special case of these
theorems.

\begin{defn}\label{metr-bldr-d5.3}
\begin{rm}
(a) Throughout this section, if not otherwise stated,
$X$ and $Y$ denote nonempty open subsets of normed spaces $E$ and $F$
respectively.
The metrics $d^E$ and $d^F$ are both abbreviated by $d$. 
For $A \subseteq X$,
$\rfs{cl}(A)$, $\rfs{bd}(A)$, $\rfs{acc}(A)$,
$B(A,r)$ etc.\ are abbreviations
for $\rfs{cl}^E(A)$, $\rfs{bd}^E(A)$, $\rfs{acc}^E(A)$, $B^E(A,r)$ etc..
Let $\vecx$, $\vecy$, $\vecx^{\kern0.7pt0}$ etc.\ denote the infinite
sequences $\fomseq{x}{n}$, $\fomseq{y}{n}$, $\fomseq{x^0}{n}$ etc..
So $\vecx \subseteq X$ means that $\fomseq{x}{n} \subseteq X$.

(b) For $A \subseteq X$ define $\delta^X(A) \eqdf d(A, E-X)$.
The notation $\delta^X(x)$ abbreviates $\delta^X(\sngltn{x})$
and $\delta^X(A)$ and $\delta^X(x)$ are abbreviated by $\delta(A)$
and $\delta(x)$.

   \index{N@d00@@$\delta^X(A) = d(A,E - X)$}
   \index{N@d01@@$\delta^X(x) = \delta^X(\sngltn{x})$}
   \index{N@d02@@$\delta(A)$. Abbreviation of $\delta^X(A)$}
   \index{N@d03@@$\delta(x)$. Abbreviation of $\delta^X(x)$}

(c) If $L$ is a rectifiable arc, then $\rfs{lngth}(L)$
denotes the length of $L$.

   \index{N@lngth@@$\rfs{lngth}(L)$. Length of an arc}

(d) Let $A \subseteq X$.
We say that $A$ is a {\it positively distanced set (PD set)},
if $\delta(A) > 0$.
A bounded PD set is called a {\it BPD set}.
A squence $\vecx$ is a {\it BPD sequence} if $\rfs{Rng}(\vecx)$
is a BPD set.

   \index{positively distanced set. A subset of $X$ whose distance from
   the boundary of $X$ is $> 0$}
   \index{pd set@@PD set. A subset of $X$ whose distance from
   the boundary of $X$ is $> 0$}
   \index{bpd set@@BPD set. A bounded subset of $X$ whose distance
   from the boundary of $X$ is positive}
   \index{bpd sequence@@BPD sequence. A sequence $\vecx$ such that
   $\rfs{Rng}(\vecx)$ is a BPD set}

(e) Let $\setm{A_i}{i \in \bbN}$ be a sequence of sets.
We define $\limti{i} A_i = x$, if for every $U \in \rfs{Nbr}(x)$
there is $i_0$ such that for every $i > i_0$, $A_i \subseteq U$.
   \index{N@lim@@$\limti{i} A_i = x$. The limit of a sequence of sets}

(f) Let $\fnn{f}{X}{Y}$.
We say that $f$ is {\it positive distance preserving ($f$ is PD.P)},
if for every PD set $A \subseteq X$, $f(A)$ is a PD subset of $Y$.
The function $f$ is {\it boundedness preserving ($f$ is BDD.P)},
if for every bounded $A \subseteq X$, $f(A)$ is a bounded set,
and $f$ is {\it boundedness positive distance preserving
($f$ is BPD.P)},
if for every bounded PD set $A \subseteq X$,
$f(A)$ is a bounded PD subset of $Y$.

   \index{pdp@@PD.P function. A function which takes PD sets
	  to PD sets}
   \index{boundedness preserving function.
	  A function which takes bounded sets to bounded sets}
   \index{bddp@@BDD.P function. A function which takes bounded
	  sets to bounded sets}
   \index{bpdp@@BPD.P function. A function which takes BPD sets to
	  BPD sets}

(h) Let $\rfs{UC}_0(X) \eqdf
\setm{f \in \rfs{UC}(X)}{\rfs{Dom}(f^{\srfs{cl}}) = \rfs{cl}(X)
\mbox{ and } f^{\srfs{cl}} \nrestriction \rfs{bd}(X) =
\rfs{Id}}$.\kern-10pt
   \index{N@uc07@@$\rfs{UC}_0(X) =
          \setm{f \in \rfs{UC}(X) \cap \rfs{EXT}(X)}{
f^{\srfs{cl}} \nrestriction \rfs{bd}(X) = \rfs{Id}}$}
\hfill\proofend
\end{rm}
\end{defn}

The following definition lists some subgroups of $H(X)$ for which
reconstruction can be proved.

\begin{defn} \label{d4.6}\label{metr-bldr-d5.5}
\begin{rm}
Let $\fnn{f}{X}{Y}$.

(a) $f$ is {\it boundedly UC ($f$ is BUC)},
if $f$ is boundedness preserving,
and for every bounded set $B \subseteq X$,
$f \nrestriction B$ is UC.
According to Definition~\ref{metr-bldr-d5.1}(a),
$\rfs{BUC}(X,Y) = \setm{f \in H(X,Y)}{f \mbox{ is BUC}}$.

   \index{boundedly UC function. A function which is uniformly
	  continuous on every bounded set}
   \index{buc function@@BUC function. A function which is uniformly
	  continuous on every bounded set}

(b) $f$ is {\it extendible},
if $\rfs{Dom}(f^{\srfs{cl}}) = \rfs{cl}(X)$.
According to Definition \ref{metr-bldr-d4.6}(b),
\newline
$\rfs{EXT}(X,Y) \eqdf \setm{f \in H(X,Y)}{f \mbox{ is extendible}}$.

   \index{extendible function. A function from $X$ to $Y$
	  that can be extended to a continuous\\
	  \indent function from
	  $\rfs{cl}(X)$ to $\rfs{cl}(Y)$}

(c) {\it $f$ is bounded positive distance UC ($f$ is BPD.UC)},
if $f$ is BPD.P, and for every BPD set $A \subseteq X$,
$f \nrestriction A$ is UC.

   \index{bounded positive distance UC function}
   \index{bpduc function@@BPD.UC function.
	  A function which is uniformly continuous on every BPD set}

(d) {\it $f$ is positive distance UC ($f$ is PD.UC)},
if $f$ is PD.P,
and for every PD set $A \subseteq X$, $f \nrestriction A$ is UC.

   \index{positive distance UC function}
   \index{pduc function@@PD.UC function. A function which is uniformly
	  continuous on every PD set}

(e) {\it  $f$ is LUC on $\rfs{bd}(X)$
($f$ is BR.LUC)}, if $f$ is extendible,
and for every $x \in \rfs{bd}(X)$
there is $U \in \rfs{Nbr}^{\srfs{cl}(X)}(x)$
such that $f^{cl} \nrestriction U $ is UC.

   \index{luc on@@LUC on $\rfs{bd}(X)$ function}
   \index{brluc function@@BR.LUC function}

(f) {\it  $f$ is completely LUC,
($f$ is CMP.LUC)},
if $f$ is extendible,
and $f^{\srfs{cl}}$ is UC at every $x \in \rfs{cl}(X)$.
That is, for every $x \in \rfs{cl}(X)$
there is $U \in \rfs{Nbr}^{\srfs{cl}(X)}(x)$
such that $f^{cl} \nrestriction U$ is UC.

   \index{completely LUC function}
   \index{cmpluc function@@CMP.LUC function.
	  An extendible function which is
	  UC at every $x \in \rfs{cl}(X)$}

(g) {\it  $f$ is UC around $\rfs{bd}(X)$
($f$ is BDR.UC)}, if $f$ is extendible, and for some $d > 0$,
$f^{cl} \nrestriction \setm{x \in \rfs{cl}(X)}{ \delta(x) < d}$ is UC.

   \index{uc around@@UC around $\rfs{bd}(X)$}
   \index{bdruc function@@BDR.UC function}

(h) Let $A,B \subseteq X$. We say that $f$ is $(A,B)$-UC,
if for every $\varepsilon > 0$ there is $\delta >0$
such that for every
$x \in A$ and $y \in B$: if $d(x,y) < \delta$, then
$d(f(x),f(y)) < \varepsilon$.
The function
{\it $f$ is BI.UC}, if $f$ is extendible, 
and $f^{cl}$ is $(\rfs{bd}(X),X)$-UC.
Note that $f$ is BI.UC iff
for every $\varepsilon > 0$ there is $\delta > 0$
such that for every $x,y \in X$: if $\delta(x), d(x,y) < \delta$,
then $d(f(x),f(y))) < \varepsilon$.\break
\rule{7pt}{0pt}\hfill\proofend
   \index{biuc function@@BI.UC function}
\end{rm}
\end{defn}

Note that if $\calP$ is one of the properties defined in
(a)\,-\,(h), that is, if
\begin{equation*}
\calP =
\rfs{BUC},\kern2.9pt \rfs{EXT},\kern2.9pt \rfs{BPD.UC},\kern2.9pt 
\rfs{PD.UC},\kern2.9pt  \rfs{BR.LUC},\kern2.9pt \rfs{CMP.LUC},\kern2.9pt
\rfs{BDR.UC}, \rfs{BI.UC},
\end{equation*}
then $\calP(X)$ is a group.

For each $\calP$ appearing above we can prove the
following statement.
If $\iso{\varphi}{\calP(X)}{\calP(Y)}$, then there is
$\tau \in \calP^{\pm}(X,Y)$ such that $\tau$ induces $\varphi$.
More precisely, the above statement can be proved, provided that
some additional restrictions are imposed on $X$ and $Y$.

We shall prove the above statement only for $\rfs{UC}(X)$
and the groups $\rfs{BUC}(X)$, $\rfs{EXT}(X)$, $\rfs{BPD.UC}(X)$ and 
$\rfs{CMP.LUC}(X)$
defined in \ref{metr-bldr-d5.5}(a), (b), (c) and (f).
Recall that the group $\rfs{LUC}(X)$ has been already dealt with in
Chapter \ref{s4}.
We omit the proof for the remaining groups,
since the arguments used are similar to those employed
in the proofs that we do present fully.
Also, the groups that we do deal with are defined by properties
that seem to have played a role in other contexts in analysis
and topology.

\smallskip

The group $\rfs{UC}(X)$
and each of the groups in Definition \ref{metr-bldr-d5.5}
except for $\rfs{EXT}(X)$ has a generalization
in which ``uniform continuity''
is replaced by ``$\itGamma$-continuity''.
This type generalization is demonstrated by the following
three examples.

{\bf Example 1 }
The generalization of $\rfs{UC}(X)$ is defined as follows.
Let $\itGamma$ be a modulus of continuity.
We say that $\fnn{f}{X}{Y}$ is nearly $\itGamma$-continuous,
if there are \hbox{$\alpha \in \itGamma$ and $r > 0$}
such that $f$ is $(r,\alpha)$-continuous.
Let $H_{\itGamma}^{\srfs{NR}}(X,Y)$ be the set of $f \in H(X,Y)$
such that $f$ is nearly $\itGamma$-continuous.
In view of Proposition~\ref{metr-bldr-p4.3}(a),
$\rfs{UC}(X) = H_{\srfs{MC}}^{\srfs{NR}}(X)$.

{\bf Example 2 }
The generalization of $\rfs{CMP.LUC}(X)$ is defined as follows.
For a modulus of continuity $\itGamma$ let
$H_{\itGamma}^{\srfs{CMP.LC}}(X) =
\setm{h \in \rfs{EXT}(X)}{\mbox{for every } x \in \rfs{cl}(X),\kern2pt
h^{\srfs{cl}} \mbox{ is }
\itGamma\mbox{-bicontinuous}\break
\mbox{at } x}$.

{\bf Example 3 }
The generalization of $\rfs{BPD.UC}(X)$ is the following group.
For a modulus of continuity $\itGamma$ let
\newline
\rule{0pt}{16pt}\kern-5pt
\kern2pt
\addtolength{\arraycolsep}{-2pt}
$
\begin{array}{lll}
\rule{0pt}{18pt}
H_{\itGamma}^{\srfs{NBPD}}(X) &\kern3.0pt=\kern3.0pt&
\setm{h \in H(X)}{h \mbox{ and } h\inverse \mbox{ are BPD.P,
and for every BPD set }
A \subseteq X,
h \nrestriction A \mbox{ is}
\\
&&
\rule{0pt}{14pt}
\kern2pt\mbox{nearly }
\itGamma\mbox{-bicontinuous}}.
\end{array}
$
\addtolength{\arraycolsep}{2pt}

\kern2mm

The reconstruction problem for these generalizations
has not been investigated thoroughly.
However, an answer for the groups in Example 3
is given in Theorem \ref{metr-bldr-t5.30}.
Example 2 is considerably more difficult to sort out.
It is dealt with in Chapters \ref{s8}\kern2.5pt-\kern1pt\ref{s12}.
The generalization in Example 1 is not true.
A counter-example is presented in Example \ref{metr-bldr-e5.10}.

So far, the reconstruction question arising from Example 2
has only a partial answer. It is proved only for
principal moduli of continuity, (see M6 in Definition \ref{nn1.3}),
and only for $X$'s with a ``well-behaved'' boundary.
This is proved in Theorem \ref{ams-bddly-lip-bldr-t6.1}.
 
\subsection{The group of uniformly continuous homeomorphisms.}
\label{ss5.2}

The first group to be considered is $\rfs{UC}(X)$.
The final reconstruction theorem for such groups is stated in
Corollary \ref{metr-bldr-c5.7}.
The following is the main intermediate theorem.

\begin{theorem}\label{t4.10}\label{metr-bldr-t5.6}
Let $X,Y \in K_{\srfs{NRM}}^{\calO}$.
Suppose that $X$ is UD.AC.
Let
$\tau \in H(X,Y)$ be such that
$(\rfs{UC}_0(X))^{\tau} \subseteq \rfs{UC}(Y)$.
Then\ %
$\tau \in \rfs{UC}(X,Y)$.
\end{theorem}

\noindent
{\bf Proof }
Variants of the argument used in this proof will be applied in
several other proofs.

Suppose by contradiction that $\tau \not\in \rfs{UC}(X,Y)$.
Let $d > 0$ and $\vecx, \vecy \subseteq X$
be such that $\lim_{n \raro \infty} d(x_n,y_n) = 0$, and for every
$n \in \bbN$, $d(\tau(x_n),\tau(y_n)) \geq d$.
Since $\tau$ is continuous, there is no $z \in X$ such that
$\setm{n}{x_n = z}$ is infinite.
So we may assume that $\vecx$ is $\onetoonen$.
We may further assume that for every distinct $m,n \in \bbN$,
$\dbltn{x_m}{y_m} \cap \dbltn{x_n}{y_n} = \emptyset$.
By \ref{p4.7}(a),
we may assume that either (i) $\vecx$ is Cauchy sequence,
or (ii) there is $e > 0$ such that $\vecx$ is $e$-spaced.

{\bf Case 1 } (i) holds. Let $x^* = \lim^{\oversE} \vecx$.
So $x^* \in \overE - X$.
Note that either $x^* \in \overline{\rfs{int}}^E(X)$
or $x^* \in \rfs{cl}^{\oversE}(\rfs{bd}(X))$.
By the UD.AC-ness of $X$
and since $\lim_{n \raro \infty} d(x_n,y_n) = 0$,
we may assume that for every $n \in \bbN$ there is an arc
$L_n \subseteq X$ connecting $x_n$ and $y_n$ such that
$\limti{n} \rfs{diam}(L_n) = 0$.
We define by induction on $k$, $n_k \in \bbN$ and $r_k > 0$.
Let $n_0 = 0$. Suppose that $n_k$ has been defined.
Let $r_k = \quarter d^{\oversE}(L_{n_k},\sngltn{x^*} \cup (E - X))$
and $n_{k + 1}$ be such that
$L_{n_{k + 1}} \subseteq B^{\oversE}(x^*,r_k)$.
We denote $x_{n_k}$, $y_{n_k}$ and $L_{n_k}$ by $u_k$,
$v_k$ and $J_k$ respectively.

Let $U_k = B^X(J_k,r_k)$.
Clearly, $\limti{k} \rfs{diam}(U_k) = 0$,
and for every 
$k \in \bbN$,
$\delta(U_k) > r_k$ and
$d(U_k,\bigcup \setm{U_m}{m \neq k}) > r_k$.
Let $w_k \in J_k - \sngltn{u_k}$ be such that
$d(\tau(u_k),\tau(w_k)) < \dgfrac{1}{(k + 1)}$.
By Lemma \ref{l2.6}(d), there is $h_k \in \rfs{LIP}(X)$ such that
$\rfs{supp}(h_k) \subseteq U_k$, 
$h_k(u_k) = u_k$ and $h_k(w_k) = v_k$. 

Let $h = \bcirc_{k \in \sboldbbN} h_k$. By Proposition \ref{p4.9},
$h \in \rfs{UC}(X)$.
\hbox{Since $\delta(\rfs{supp}(h_k)) > 0$,}
$h \in \rfs{UC}_0(X)$.
We check that $h^{\tau} \not\in \rfs{UC}(Y)$.
Clearly, $h^{\tau}(\tau(u_k)) = \tau(u_k)$ and
$h^{\tau}(\tau(w_k)) = \tau(v_k)$.
However,
$\limti{k} d(\tau(u_k),\tau(w_k)) = 0$,
whereas for every $k \in \bbN$,
$d(\tau(u_k),\tau(v_k)) \geq d$. So $h^{\tau} \not\in \rfs{UC}(Y)$.

{\bf Case 2 } (ii) holds.
By the UD.AC-ness of $X$, and since
$\limti{n} d(x_n,y_n) = 0$,
there is $N \in \bbN$ such that for every $n \geq N$ there is an arc
$L_n \subseteq X$ connecting $x_n$ and $y_n$
such that $\rfs{diam}(L_n) < \dgfrac{e}{6}$
and $\limti{n} \rfs{diam}(L_n) = 0$. We may assume that $N = 0$.
Let
$r_n = \min(\rfs{diam}(L_n), \dgfrac{\delta(L_n)}{2})$ and
$U_n = B(L_n,r_n)$.
So $\delta(U_n) > 0$,
$\limti{n} \rfs{diam}(U_n) = 0$,
and for every distinct $m,n \in \bbN$, $d(U_m,U_n) \geq \dgfrac{e}{3}$. 
The proof now proceeds as in Case 1.\medskip
\hfill\myqed\smallskip

The final result for groups of type $\rfs{UC}(X)$ is at this stage as
follows.

\begin{cor}\label{metr-bldr-c5.7}
Let $X,Y \in K_{\srfs{NFCB}}^{\calO}$.
Suppose that $X$ and $Y$ are UD.AC spaces.
Let $\iso{\varphi}{\rfs{UC}(X)}{\rfs{UC}(Y)}$.
Then there is
$\tau \in \rfs{UC}^{\pm}(X,Y)$ such that $\tau$ induces $\varphi$.
\end{cor}

\noindent
{\bf Proof }
Combine Corollary \ref{metr-bldr-c2.26}\ 
and Theorem \ref{metr-bldr-t5.6}.
\smallskip\hfill\myqed

In the case of local uniform continuity, we deduced from the fact that
$(\rfs{UC}(X))^{\tau} \subseteq \rfs{LUC}(Y)$, that both $\tau$ and 
$\tau\inverse$ are LUC. The analogue of this fact for uniform
continuity is not true.

\begin{example}\label{e4.11}\label{metr-bldr-e5.8}
\num{a} Let $X = Y = (1,\infty)$, and $\fnn{\tau}{X}{Y}$ be defined by
$\tau(x) = \sqrt{x}$.
Then $(\rfs{UC}(X))^{\tau} \subseteq \rfs{UC}(Y)$, but $\tau\inverse$ is
not UC.

\num{b} There are bounded open subsets $X$ and $Y$ of the
Hilbert space $\ell_2$ and $\tau \in H(X,Y)$ such that
$(\rfs{UC}(X))^{\tau} \subseteq \rfs{UC}(Y)$, but $\tau^{-1}$ is not
uniformly continuous. The boundary of both $X$ and $Y$ is the union
of a spaced family of spheres.
\end{example}

\noindent
{\bf Proof }
(a) Clearly $\tau\inverse \not\in \rfs{UC}(X)$.
Let $f \in \rfs{UC}(X)$.
By Proposition~\ref{metr-bldr-p4.3}(b), $f$ is
$\alpha$-continuous for some $\alpha \in \rfs{MC}$.
By the uniform continuity of $f\inverse$, there is $C$ such
that for every $y \in X$, $f\inverse(y + 1) -  f\inverse(y) \leq C$.
Set $K = C + 1$. We check that
$f(x) \geq \dgfrac{x}{K}$ for every $x \in X$. Let $y \in X$. Then
$f\inverse(y) - 1 \leq f\inverse([y] + 1) - f\inverse(1) \leq
[y] \mcdot C \leq y \mcdot C$.
Hence $f\inverse(y) \leq Cy + 1 \leq (C + 1)y$.
That is, $y \leq f((C + 1) y)$.
Write $x = (C + 1) y$. We conclude that
if $x \geq C + 1$, then $\dgfrac{x}{K} \leq f(x)$.
The above inequality holds automatically for $x \leq C + 1$
since $f(x) \geq 1$.

We show that $f^{\tau}$ is $(1,2 \sqrt{K} \alpha)$-continuous.
This trivially implies that $f^{\tau}$ is UC.
Let $y > x \geq 1$ be such that $y - x \leq 1$.
We have $\tau\inverse(y) - \tau\inverse(x) = y^2 - x^2 \leq 2y(y - x)$. 
So $f(\tau\inverse(y)) - f(\tau\inverse(x)) \leq \alpha(2y(y - x)) \leq
2y\alpha(y - x)$.
The last inequality follows from the fact that $2y \geq 1$.
Now, $\tau f \tau\inverse(y) - \tau f \tau\inverse(x) =
\sqrt{f(y^2)} - \sqrt{f(x^2)}$.
There is $c \in (f(x^2),f(y^2))$ such that
$\sqrt{f(y^2)} - \sqrt{f(x^2)} =
\frac{1}{2\sqrt{c}}(f(y^2) - f(x^2))$.
Recall that $f(x^2) \geq \dgfrac{x^2}{K}$.
So\vspace{2pt}
\newline\rule{8pt}{0pt}
\renewcommand{\arraystretch}{1.6}
\addtolength{\arraycolsep}{0pt}
$
\begin{array}{ll}
&
f^{\tau}(y) - f^{\tau}(x) =
\tau f \tau\inverse(y) - \tau f \tau\inverse(x) = 
\frac{1}{2\sqrt{c}}(f(y^2) - f(x^2)) \leq
\frac{1}{2\sqrt{f(x^2)}} \cdot 2y \alpha(y - x)
\\
\leq
&
\frac{1}{\sqrt{\dgfrac{x^2}{K}}} \cdot y \alpha(y - x) \leq
\frac{1}{\sqrt{\dgfrac{x^2}{K}}} \cdot 2x \alpha(y - x) =
2 \sqrt{K} \alpha(y - x).
\end{array}
$
\renewcommand{\arraystretch}{1.0}
\addtolength{\arraycolsep}{0pt}
\smallskip

(b) In $\ell_2$ let
$e_i = (0,\ldots,0,\raise0.7pt\hbox{$\stackrel{i}{1}$},0,\ldots)$
and $a_i = 3\sqrt{2}e_i$.
Let
$X = B(0,6) - \bigcup_{n > 0} B(a_i,1)$ and
$Y = B(0,6) - \bigcup_{n > 0} B(a_i,\dgfrac{1}{n})$.
For every $n > 0$ let $\fnn{h_n}{[0,\infty)}{[0,\infty)}$
be the piecewise linear function with two breakpoints which takes
$0$ to $0$, $1$ to $\dgfrac{1}{n}$,
and such that $h_n(t) = t$ for every $t \geq 2$.
Let $\fnn{\tau_n}{X}{Y}$ be defined by
$\tau_n(x) =
a_n + h_n(\norm{x - a_n})\frac{x - a_n}{\norm{x - a_n}}$,
and $\tau = \bcirc_{n > 0} \tau_n$. It is left to the reader to check
that $\tau$ is as required.
\bigskip\hfill\myqed

We shall later see a finite-dimensional example in which
$(\rfs{UC}(X))^{\tau} \subseteq \rfs{UC}(Y)$,
but $\tau^{-1}$ is not uniformly continuous.
In Example \ref{e4.34}(a) we construct two bounded domains
$X,Y \subseteq \bbR^{2}$ and $\tau \in H(X,Y)$ with these properties.

However, for some sets $X$, which are very well behaved, the fact that
$(\rfs{UC}(X))^{\tau} \subseteq \rfs{UC}(Y)$
does imply that $\tau^{-1}$ is uniformly continuous.
Theorems~\ref{t4.36} and \ref{t4.41}(a)
and Remark~\ref{metr-bldr-r7.8}(b) and (c) prove the above
implication in some special cases involving subsets of a Banach space
or a Banach manifold.
For example, the above implication holds when $X$ and $Y$ are spheres
of a Banach space.

\medskip

\noindent
{\bf Proof of Theorem \ref{t4.4} }
Let $X'$ and $Y'$ be UD.AC open subsets of $\bbR^m$ and $\bbR^n$
respectively, $F = \rfs{cl}(X')$, $K = \rfs{cl}(Y')$
and $\iso{\varphi}{H(F)}{H(K)}$. Let $X = \rfs{int}(F)$
and  $Y = \rfs{int}(K)$.
Clearly, $X$ and $Y$ are regular open sets,
$F = \rfs{cl}(X)$ and $K = \rfs{cl}(Y)$.
It is trivial to check that $X$ and $Y$ are UD.AC.
It is also trivial to check that if $Z$ is a bounded regular open
subset of $\bbR^k$,
then $H(\rfs{cl}(Z)) = \setm{f^{\srfs{cl}}}{f \in \rfs{UC}(Z)}$.
Let $\fnn{\psi}{H(X)}{H(Y)}$ be defined by
$\psi(f) = \varphi(f^{\srfs{cl}}) \nrestriction Y$.
So $\iso{\psi}{\rfs{UC}(X)}{\rfs{UC}(Y)}$.

By Theorem \ref{t2.4}, there is
$\tau \in H(X,Y)$ such that for every $h \in \rfs{UC}(X)$,
$\psi(h) = h^{\tau}$.
Obviously, $(\rfs{UC}(X))^{\tau} = \rfs{UC}(Y)$.
Applying Theorem \ref{t4.10} to $\tau$ and $\tau\inverse$
one concludes that $\tau$ and $\tau\inverse$ are uniformly continuous.
It follows that $\iso{\tau^{\srfs{cl}}}{F}{K}$.
It is trivial that for every $h \in H(F)$,
$\varphi(h) = h^{\tau^{\srs{cl}}}$.
\smallskip\hfill\myqed

Part (a) of the next example shows that in
Theorem~\ref{metr-bldr-t5.2},
the requirement that $F$ and $K$ are bounded cannot be dropped,
and Part (b) shows that in Theorem~\ref{metr-bldr-t5.2},
the requirement that $F$ and $K$ are closures of UD.AC open sets
cannot be dropped.

\begin{example}\label{e4.13}\label{metr-bldr-e5.9}
\num{a} There are regular open connected subsets $X,Y \subseteq \bbR^2$
such that $X,Y$ are UD.AC,
$X$ is bounded,
$\rfs{cl}(X) \not\cong \rfs{cl}(Y)$
but $H(\rfs{cl}(X)) \cong H(\rfs{cl}(Y))$.

\num{b} There are regular open connected subsets $X,Y \subseteq \bbR^2$
such that $X$ is UD.AC,
$X$ and $Y$ are bounded,
$\rfs{cl}(X) \not\cong \rfs{cl}(Y)$
but $H(\rfs{cl}(X)) \cong H(\rfs{cl}(Y))$.
\end{example}

\noindent
{\bf Proof }
(a) Let $x \in S(0,1)$ and
$B_i = B(\dgfrac{x}{2^{2i + 2}}\,,\,\dgfrac{1}{2^{2i + 3}})$.
So
$\bigcup_{i \in \sboldbbN} B_i \subseteq B(0,\dgfrac{1}{2})$,
for every $i \neq j$,
$\rfs{cl}(B_i) \cap \rfs{cl}(B_j) = \emptyset$
and \hbox{$\limti{i} B_i = 0$.}

Let $F = \rfs{cl}(B(0,1)) - \bigcup_{i \in \sboldbbN} B_i$.
Let $\tau(x) \eqdf \dgfrac{x}{\norm{x}^2}$ be the inversion map
in $\bbR^2$ and $K = \tau(F - \sngltn{0})$.
Let $X = \rfs{int}(F)$ and $Y  = \rfs{int}(K)$.
Then $F = \rfs{cl}(X)$ and $K = \rfs{cl}(Y)$.
Clearly, $X,Y$ are UD.AC.
It is easy to see that
$H(K) = \setm{(h \nrestriction (F - \sngltn{0}))^{\tau}}{h \in H(F)}$.
So $H(F) \cong H(K)$.
It is obvious $F \not\cong K$.

(b) Let
$$\hbox{$
X_0 =
\setm{(\theta - \pi,t)}
{\theta \in (0,2\pi),\ %
t \in (1 - \quarter \mcdot \abs{\sin\frac{\theta}{2}},
1 + \quarter \mcdot \abs{\sin\frac{\theta}{2}}
}
$}$$
and
$$\hbox{$
Y_0 = \setm{t \mcdot (\cos\theta,\sin\theta)}
{\theta \in (0,2\pi),\ %
t \in (1 - \quarter \mcdot \abs{\sin\frac{\theta}{2}},
1 + \quarter \mcdot \abs{\sin\frac{\theta}{2}}}.
$}$$
Note that $X_0$ is a strip surrounding the line segment
$((-\pi,0),(\pi,0))$
with width tending to $0$ as $(\theta,0)$ approaches
$(-\pi,0)$ and $(\pi,0)$,
and
$Y_0$ is a strip surrounding the circular arc
$\setm{(\cos\theta,\sin\theta)}{\theta \in (0,2\pi)}$
with width tending to $0$ as $\theta$ approaches $0$ and $2\pi$.
Let $\fnn{\tau}{X_0}{Y_0}$ be defined by
$\tau((\theta - \pi,t)) = t \mcdot (\cos\theta,\sin\theta)$.
Then $\tau \in H(X_0,Y_0)$.

For every $n \in \bbZ$
let $x_n = (\frac{n}{\abs{n} + 1} \ncdot \pi,0)$,
$r_n = \third
\min(\delta^{X_0}(x_n), d(x_n,\setm{x_i}{i \in \bbZ - \sngltn{n}})$
and
$\overB_n = \overB(x_n,r_n)$.
So $\overB_n \subseteq X_0$, for $n \neq m$,
$\overB_n \cap \overB_m = \emptyset$,
$\limti{n} \overB_n = (\pi,0)$ and
$\limtmi{n} \overB_n = (-\pi,0)$.
Let $X = X_0 - \bigcup_{n \in \bbZ} \overB_n$
and $Y = \tau(X)$.
Clearly, $X$ and $Y$ are bounded, connected and regular open.
Hence $H(\rfs{cl}(X)) = (H(X))^{\srfs{cl}}$,
and the same holds for $Y$.
It is also obvious that $\rfs{cl}(X) \not\cong \rfs{cl}(Y)$.
Note that for every $h \in H(\rfs{cl}(X))$,
$h((\pi,0)) \in \dbltn{(\pi,0)}{(-\pi,0)}$ and the same holds for
$(-\pi,0)$.
Also, for every $h \in H(\rfs{cl}(Y))$,
$h((1,0)) = (1,0)$.
It follows that
$h^{\srfs{cl}} \mapsto (h^{\tau})^{\srfs{cl}}, \ \, h \in H(X)$,
is an isomorphism between $H(\rfs{cl}(X))$ and $H(\rfs{cl}(Y))$.
\smallskip\hfill\myqed

Example \ref{metr-bldr-e5.9}(b) calls for the following questions.

\begin{question}\label{metr-bldr-q5.10}
\begin{rm}
A topological space $Z$ has the Perfect Orbit Property,
if for every $z \in Z$, $z \in \rfs{acc}(\setm{h(z)}{h \in H(Z)})$.
Is it true that for every open $X \subseteq \bbR^m$
and $Y \subseteq \bbR^n$:
if $\rfs{cl}(X)$ and $\rfs{cl}(Y)$ have the Perfect Orbit Property
and $\iso{\varphi}{H(\rfs{cl}(X))}{H(\rfs{cl}(Y))}$,
then there is $\tau \in H(\rfs{cl}(X),\rfs{cl}(Y))$
such that $\tau$ induces $\varphi$?

If the above is not true, is the conclusion in the above question true
for open subsets of $\bbR^n$ that have the following stronger property:
For every $x \in \rfs{bd}(X)$
the orbit of $x$ under $H(\rfs{cl}(X))$ is locally arcwise connected.

Is the same true for open subsets of infinite-dimensional normed
spaces?
\hfill\proofend
\end{rm}
\end{question}

The generalization of Corollary \ref{metr-bldr-c5.7} is not true for all
moduli of continuity.
As shown in the next example,
$\itGamma^{\srfs{LIP}}$ is a counter-example.
The question whether Theorem \ref{metr-bldr-c5.7} is true for any
countably generated $\itGamma$ is open.

\begin{question}\label{metr-bldr-q5.11}
\begin{rm}
Is there a countably generated modulus of continuity $\itGamma$
such that for every normed space $E$ and $\tau \in H(E)$:
if $(H_{\itGamma}(E))^{\tau} = H_{\itGamma}(E)$,
then $\tau \in H_{\itGamma}(E)$?
\end{rm}
\end{question}

\begin{example}\label{metr-bldr-e5.10}
\begin{rm}
Let $E$ be a normed space
and $\tau \in H(E)$ be defined by:
$\tau(x) = x$ if $\norm{x} \leq 1$
and $\tau(x) = \norm{x} \mcdot x$ if $\norm{x} > 1$.
Then $(\rfs{LIP}(E))^{\tau} = \rfs{LIP}(E)$
and $\tau \not\in \rfs{LIP}(E,E)$.
\end{rm}
\end{example}

\noindent
{\bf Proof }
Let $g \in \rfs{LIP}(X,X)$.
We show that $g^{\tau}$ is Lipschitz.
Let $r$ be such that
$r \geq 1,\norm{g(0)}$ and $g(B(0,r)) \supseteq B(0,1)$.
We show that $g^{\tau} \nrestriction (E - B(0,r^2))$ is Lipschitz.
Suppose that $g$ is $K$-Lipschitz.
Let $u \in E - B(0,r)$.
Then
$$
\norm{g(u)} \leq \norm{g(u) - g(0)} + \norm{g(0)} \leq
K \norm{u} + \norm{g(0)} \leq K \norm{u} + \norm{u} = (K + 1) \norm{u}.
$$
That is,
\begin{itemize}
\addtolength{\parskip}{-11pt}
\addtolength{\itemsep}{06pt}
\item[(i)] 
$\norm{g(u)} \leq (K + 1) \norm{u}$.
\vspace{-05.7pt}
\end{itemize}
For $u,v \in E - \sngltn{0}$
write $w(u,v) = \frac{\norm{v}}{\norm{u}} u$,
and for $u,v \neq g\inverse(0)$
set
$w_g(u,v) = w(g(u),g(v))$.
Clearly,
\begin{itemize}
\addtolength{\parskip}{-11pt}
\addtolength{\itemsep}{06pt}
\item[(ii)]
$\norm{u - w(u,v)} = \abs{\norm{u} - \norm{v}} \leq \norm{u - v}$,
\item[(iii)]
$\norm{w(u,v) - v} \leq \norm{w(u,v) - u} +
\norm{u - v} \leq 2 \norm{u - v}$,
\vspace{-05.7pt}
\end{itemize}
and it follows that
\begin{itemize}
\addtolength{\parskip}{-11pt}
\addtolength{\itemsep}{06pt}
\item[(iv)]
$\norm{g(u) - w_g(u,v)} \leq K \norm{u - v}$,
\item[(v)] 
$\norm{w_g(u,v) - g(v)} \leq 2K \norm{u - v}$.
\vspace{-05.7pt}
\end{itemize}

{\bf Claim 1.} There is $M$ such that for every $x,y \in E - B(0,r^2)$:
if $y = \lambda x$ for some $\lambda > 1$, then
$\norm{g^{\tau}(y) - g^{\tau}(x)} \leq M \norm{y - x}$.

{\bf Proof }
Let $x = az$ and $y = (a + e)z$, where $\norm{z} = 1$ and $a > 0$.
Clearly, $e > 0$ and hence $\norm{y - x} = e$.
Also, $a \geq r^2$.
Then $\norm{\tau\inverse((a + e)z) - \tau\inverse(az)} =
\sqrt{a + e} - \sqrt{a} \leq \dgfrac{e}{\sqrt{a + e}}$.
Set $u = \tau\inverse((a + e)z)$ and $v = \tau\inverse(az)$.
So $\norm{u - v} \leq \dgfrac{e}{\sqrt{a + e}}$.
The next inequality uses the definitions of $\tau$ and $w_g$,
the $K$-Lipschitz-ness of $g$ and Fact (i).
\newline\rule{8pt}{0pt}
\renewcommand{\arraystretch}{1.6}
\addtolength{\arraycolsep}{-1.0pt}
$
\begin{array}{ll}
&
\norm{\tau(g(u)) - \tau(w_g(u,v))} =
\abs{\norm{g(u)}^2 - \norm{w_g(u,v)}^2} =
\abs{\norm{g(u)}^2 - \norm{g(v)}^2}
\\
=
&
(\norm{g(u)} + \norm{g(v)}) \mcdot \abs{\norm{g(u)} - \norm{g(v)}} \leq
(\norm{g(u)} + \norm{g(v)}) \mcdot \norm{g(u) - g(v)}
\\
\leq
&
(\norm{g(u)} + \norm{g(v)}) \mcdot K \norm{u - v} \leq
(\norm{g(u)} + \norm{g(v)}) \mcdot \frac{Ke}{\sqrt{a + e}}
\\
\leq
&
(K + 1) (\norm{u} + \norm{v}) \mcdot \frac{Ke}{\sqrt{a + e}} =
(K + 1) (\sqrt{a + e} + \sqrt{a}) \mcdot \frac{Ke}{\sqrt{a + e}}
\\
\leq
&
2 (K + 1)^2 \sqrt{a + e} \mcdot \frac{e}{\sqrt{a + e}} =
2 (K + 1)^2 e =
2 (K + 1)^2 \norm{y - x}.
\end{array}
$
\renewcommand{\arraystretch}{1.0}
\addtolength{\arraycolsep}{1.0pt}
\vspace{2mm}\newline
We next find a bound for $\norm{\tau(w_g(u,v)) - \tau(g(v))}$.
Since $g$ is $K$-Lipschitz and by (v),
\newline\rule{8pt}{0pt}
\renewcommand{\arraystretch}{1.6}
\addtolength{\arraycolsep}{-1.5pt}
$
\begin{array}{ll}
&
\norm{\tau(w_g(u,v)) - \tau(g(v))} =
\norm{g(v)} \mcdot \norm{w_g(u,v) - g(v)} \leq
(K + 1) \mcdot \norm{v} \mcdot 2K \mcdot \norm{u - v}
\\
\leq
&
(K + 1) \mcdot \sqrt{a} \mcdot 2K \mcdot \frac{e}{\sqrt{a + e}} \leq
2 (K + 1)^2 \mcdot \norm{y - x}.
\end{array}
$
\renewcommand{\arraystretch}{1.0}
\addtolength{\arraycolsep}{1.5pt}
\vspace{2mm}\newline
Note that $g^{\tau}(y) = \tau(g(u))$ and 
$g^{\tau}(x) = \tau(g(v))$.
It follows that
{\thickmuskip=4mu \medmuskip=3mu \thinmuskip=2mu 
$$\hbox{$
\norm{g^{\tau}(y) - g^{\tau}(x)} \leq
\norm{\tau(g(u)) - \tau(w_g(u,v))} + \norm{\tau(w_g(u,v)) - \tau(g(v))}
\leq
4 (K + 1)^2 \mcdot \norm{y - x}.
$}$$
}
So Claim 1~is proved.
\smallskip

{\bf Claim 2.} There is $M$ such that for every
$x,y \in E - B(0,r^2)$:
if
$\norm{x} = \norm{y}$,
then $\norm{g^{\tau}(x) - g^{\tau}(y)} \leq M \norm{x - y}$.

{\bf Proof }
Let $\norm{x} = \norm{y} = a \geq r^2$.
Set $u = \tau\inverse(x)$ and $v = \tau\inverse(y)$.
Then by (iv),
$\norm{g(u) - w_g(u,v)} \leq K \norm{u - v}$.
So
\newline\rule{2pt}{0pt}
\renewcommand{\arraystretch}{1.6}
\addtolength{\arraycolsep}{-0.0pt}
$
\begin{array}{ll}
&
\norm{\tau(g(u)) - \tau(w_g(u,v))} =
\abs{\norm{g(u)}^2 - \norm{w_g(u,v)}^2} =
\abs{\norm{g(u)}^2 - \norm{g(v)}^2}
\\
=
&
(\norm{g(u)} + \norm{g(v)}) \mcdot \abs{\norm{g(u)} - \norm{g(v)}}
\leq
(K + 1) (\norm{u} + \norm{v}) \mcdot \norm{g(u) - g(v)}
\\
\leq
&
2(K + 1) \sqrt{a} \mcdot K \norm{u - v} = 2 (K + 1)K \sqrt{a} \mcdot
\frac{\norm{x - y}}{\sqrt{a}} \leq 2 (K + 1)^2 \norm{x - y}.
\end{array}
$
\renewcommand{\arraystretch}{1.0}
\addtolength{\arraycolsep}{0.0pt}
\vspace{2mm}\newline
We next find a bound for
$\norm{\tau(w_g(u,v)) - \tau(g(v))}$.
By (iv) we have $\norm{w_g(u,v) - g(v)} \leq 2 K \norm{u - v}$. So
\newline\rule{8pt}{0pt}
\renewcommand{\arraystretch}{1.6}
{\thickmuskip=5.0mu \medmuskip=4mu \thinmuskip=2mu 
$
\begin{array}{ll}
&
\norm{\tau(w_g(u,v)) - \tau(g(v))} =
\norm{g(v)} \mcdot \norm{w_g(u,v) - g(v)}
\leq
(K + 1) \sqrt{a} \mcdot \norm{w_g(u,v) - g(v)}
\\\leq
&
(K + 1) \sqrt{a} \mcdot 2K \norm{u - v} =
(K + 1) \sqrt{a} \mcdot 2K \mcdot \frac{\norm{x - y}}{\sqrt{a}} \leq
2 (K + 1)^2 \norm{x - y}.
\end{array}
$
}
\renewcommand{\arraystretch}{1.0}
\vspace{2mm}\newline
It follows that
$\norm{g^{\tau}(x) - g^{\tau}(y)} \leq 4(K + 1)^2 \norm{x - y}$.
We have proved Claim 2.
\smallskip

Let $x,y \in E - B(0,r^2)$.
By Claims 1 and 2 and by (ii) and (iii),
\newline\rule{5pt}{0pt}
\renewcommand{\arraystretch}{1.6}
\addtolength{\arraycolsep}{0.5pt}
$
\begin{array}{ll}
&
\norm{g^{\tau}(x) - g^{\tau}(y)} \leq
\norm{g^{\tau}(x) - g^{\tau}(w(x,y))} +
\norm{g^{\tau}(w(x,y)) - g^{\tau}(y)}
\\
\leq
&
4(K + 1)^2 \norm{x - w(x,y)} + 4 (K + 1)^2 \norm{w(x,y) - y} \leq
12 (K + 1)^2 \norm{x - y}.
\end{array}
$
\renewcommand{\arraystretch}{1.0}
\addtolength{\arraycolsep}{-0.5pt}
\medskip

We have shown that if $g$ is Lipschitz,
then $g^{\tau} \nrestriction (E - B(0,r^2))$ is Lipschitz.
Since for every bounded set $B$, $\tau \nrestriction B$ is bilipschitz,
it follows that $g^{\tau} \nrestriction \overB(0,r^2)$ is Lipschitz.
It is now esay to conclude that $g^{\tau}$ is Lipschitz.

The proof that $(\rfs{LIP}(E))^{\tau\inverse} \subseteq \rfs{LIP}(E)$
is slightly different. Denote $\tau\inverse$ by $\eta$.
We prove that if $g$ is bilipschitz, then $g^{\eta}$ is Lipschitz.
Let $g \in \rfs{LIP}(X)$,
suppose that $g$ is $K$-bilipschitz
and let $r$ be such that
$r \geq \rfs{max}(1,2K\norm{g(0)})$ and $g(B(0,r)) \supseteq B(0,1)$.
We show that $g^{\eta} \nrestriction (E - B(0,\sqrt{r}))$ is Lipschitz.

We shall use facts (ii)\,-\,(v) from the preceding part of the proof.
In addition, we need the following fact.
Let $u \in E - B(0,r)$. Then
$$
\norm{g(u)} \geq \norm{g(u) - g(0)} - \norm{g(0)} \geq
\dgfrac{\norm{u}}{K} - \norm{g(0)} \geq
\dgfrac{\norm{u}}{K} - \dgfrac{\norm{u}}{2K} = \dgfrac{\norm{u}}{(2K)}.
$$
That is,
\begin{itemize}
\addtolength{\parskip}{-09pt}
\addtolength{\itemsep}{06pt}
\item[(vi)] 
$\norm{g(u)} \geq \dgfrac{\norm{u}}{(2K)}$.
\end{itemize}

{\bf Claim 3.} There is $M$ such that
for every $x,y \in E - B(0,\sqrt{r})$:
if $y = \lambda x$ for some $\lambda > 1$, then
$\norm{g^{\eta}(y) - g^{\eta}(x)} \leq M \norm{y - x}$.

{\bf Proof }
Let $x = az$ and $y = (a + e)z$, where $\norm{z} = 1$ and $a,e > 0$.
Then $\norm{y - x} = e$ and $a \geq \sqrt{r}$.
Set $u = \eta\inverse((a + e)z)$ and $v = \eta\inverse(az)$.
We skip the verification of the following facts:
\begin{equation}\tag{1}
\norm{g^{\eta}(x)) - \eta(w_g(v,u))} \leq 
\sqrt{2}K^{\dgfrac{3}{2}} \norm{x - y},
\end{equation}
\begin{equation}\tag{2}
\norm{\eta(w_g(v,u)) - g^{\eta}(y))} \leq
4 \sqrt{2} K^{\dgfrac{3}{2}} \norm{x - y}.
\end{equation}
From (1) and (2) it follows that
{\thickmuskip=4mu \medmuskip=3mu \thinmuskip=2mu 
$$\hbox{$
\norm{g^{\eta}(x) - g^{\eta}(y)} \leq
\norm{\eta(g(v)) - \eta(w_g(v,u))} +
\norm{\eta(w_g(v,u)) - \eta(g(u))} \leq
5 \sqrt{2} K^{\dgfrac{3}{2}} \norm{x - y}.
$}$$
}
This proves Claim 3.
\smallskip

{\bf Claim 4.} There is $M$ such that for every
$x,y \in E - B(0,\sqrt{r})$:
if
$\norm{x} = \norm{y}$,
then $\norm{g^{\eta}(x) - g^{\eta}(y)} \leq M \norm{x - y}$.

{\bf Proof }
Let $\norm{x} = \norm{y} \geq \sqrt{r}$.
Set $u = \eta\inverse(x)$ and $v = \eta\inverse(y)$.
We skip the verification of the following facts:
\begin{equation}\tag{3}
\norm{\eta(g^{\eta}(y)) - \eta(w_g(v,u))} \leq
(\dgfrac{\sqrt{2}}{2}) K^{\dgfrac{3}{2}} \norm{y - x},
\end{equation}
\begin{equation}\tag{4}
\norm{\eta(w_g(v,u)) - g^{\eta}(x))} \leq
2 \sqrt{2} K^{\dgfrac{3}{2}} \norm{y - x}.
\end{equation}
We conclude that
$$
\norm{g^{\eta}(y) - g^{\eta}(x)} \leq
(\dgfrac{5\sqrt{2}}{2}) K^{\dgfrac{3}{2}} \norm{y - x}.
$$
This proves Claim 4.

\smallskip
The rest of the argument is the same as in the preceding part
of the proof.
\rule{7pt}{0pt}\hfill\myqed

\subsection{The group of homeomorphisms which are uniformly
continuous on every bounded set.}
\label{ss5.3}
\label{ss5.3-BUC-homeomorphisms}

We now turn to the group $\rfs{BUC}(X)$ of all homeomorphisms $f$ of $X$
such that $f$ and $f\inverse$ are boundedness preserving,
and $f$ and $f\inverse$ are uniformly continuous on every bounded subset
of $X$.
The final reconstruction result for such groups is stated in
Theorem \ref{metr-bldr-t5.19}.
The conclusion of \ref{metr-bldr-t5.19} is the statement:
$(*)$ if $\iso{\varphi}{\rfs{BUC}(X)}{\rfs{BUC}(Y)}$,
then there is $\tau \in \rfs{BUC}^{\pm}(X,Y)$
such that $\tau$ induces $\varphi$.
However, $(*)$ is not true for general open subsets of a normed space,
so we shall make some extra assumptions on $X$ and $Y$.
These assumptions are (roughly):
(1) $X$ and $Y$ are
uniformly\,-\,in\,-\,diameter arcwise\,-\,connected;\break
(2) the orbit of every member of $\rfs{bd}(X)$ under the action of
$\rfs{BUC}(X)$ contains an arc, and the same holds for $Y$.

Let
$\rfs{ABUC}(X,Y) = \setm{h \in H(X,Y)}{\mbox{for every bounded set }
A \subseteq X,\ h \nrestriction A \mbox{ is UC}}$.
Recall that $\rfs{ABUC}(X) = \rfs{ABUC}^{\pm}(X,X)$.
Whereas $\rfs{BUC}(X)$ is a group,
it is not always true that $\rfs{ABUC}(X)$ is a group.
It is easy to construct an open set $X$ in a normed space and
$f \in \rfs{ABUC}(X)$ such that $f$ takes a bounded set
to an unbounded set. We can then choose another $g \in \rfs{ABUC}(X)$
such that $g \scirc f \not\in \rfs{ABUC}(X)$.
However, if $X$ has the discrete path property for large
distances, (see \ref{metr-bldr-d4.2}(f)),
then every member of $\rfs{ABUC}(X)$ is boundedness
preserving, and hence $\rfs{ABUC}(X) = \rfs{BUC}(X)$.
So $\rfs{ABUC}(X)$ is a group.

   \index{N@abuc@@$\rfs{ABUC}(X,Y) =
          \setm{h \in H(X,Y)}{\mbox{for every bounded set }
          A \subseteq X,\ h \nrestriction A \mbox{ is UC}}$}

\begin{prop}\label{p4.12}\label{metr-bldr-p5.10}
Let X have the discrete path property for large distances.

\num{a} There are $a_1,b_1 > 0$ such that,
for every $x,y \in X$ and $0 < t < d(x,y)$, there are $n \in \bbN$
and $x = x_0,\,x_1,\,\ldots,\,x_n = y$ such that
$n \leq \dgfrac{(a_1 d(x,y) + b_1)}{t}$, and for every $i < n$,
$d(x_i,x_{i+1}) \leq t$.

\num{b} If $Y$ is a metric space, and $\tau \in \rfs{ABUC}(X,Y)$,
then $\tau$ is boundedness preserving.
(Hence $\tau \in \rfs{BUC}(X,Y)$).

\num{c} $\rfs{BUC}(X) = \rfs{ABUC}(X)$.
\end{prop}

\noindent
{\bf Proof }
(a) Let
$x = z_0,\,z_1,\,\ldots,\,z_m = y$ be such that
$d(z_i,z_{i + 1}) < \dgfrac{t}{2}$ for every $i < m$,
and $\sum_{i < m}d(z_i,z_{i + 1}) \leq ad(x,y) + b$.
There are $n \in \bbN$ and $0 = i_0 < \ldots < i_n \leq m$
such that for every $j < n$,
$\dgfrac{t}{2} \leq d(z_{i_j},z_{i_{j + 1}}) < t$
and $d(z_{i_n},z_m) \leq \dgfrac{t}{2}$.
It follows that
$n \cdot \frac{t}{2} \leq
\sum_{j < i_n}d(z_j,z_{j + 1}) \leq ad(x,y) + b$.
Hence $n \leq \dgfrac{(2ad(x,y) + 2b)}{t}$ and so
$n + 1 \leq \dgfrac{((2a + 1)d(x,y) + 2b)}{t}$.
For $j \leq n$ define $x_j = z_{i_j}$ and define $x_{n + 1} = z_m$.
Then  $n + 1$ and $x_0,\ldots,x_{n + 1}$ are as required. That is,
we may take $a_1$ and $b_1$ to be $2a + 1$ and $2b$.
So (a) is proved.

(b) Let $a_1,b_1$ be the numbers obtained by applying Part (a) to $X$.
Let $C \subseteq X$ be bounded. Define $r = \rfs{diam}(C)$
and $B = B(C,a_1 r + b_1)$. Since $B$ is bounded, there is $\delta > 0$
such that for every $x,y \in B$: if $d(x,y) \leq \delta$,
then $d(\tau(x),\tau(y)) \leq 1$.
Let $x,y \in C$.
If $d(x,y) \leq \delta$, then 
$d(\tau(x),\tau(y)) \leq 1$.
Otherwise, let $n \in \bbN$ and $x = z_0,\ldots,z_n = y$ be such that
$n \leq \dgfrac{(a_1 d(x,y) + b_1)}{\delta}$
and $d(z_i,z_{i + 1})) \leq \delta$ for every $i < n$.
So for every $i \leq n$,
$d(x,z_i) \leq n \delta \leq
\frac{a_1 d(x,y) + b_1}{\delta} \mcdot \delta \leq
a_1 r + b_1$.
So $z_i \in B$ and hence $d(\tau(z_i),\tau(z_{i + 1})) \leq 1$.
Then
$d(\tau(x),\tau(y)) \leq \sum_{i < n} d(\tau(z_i),\tau(z_{i + 1})) \leq
n \leq \dgfrac{(a_1 d(x,y) + b_1)}{\delta} \leq
\dgfrac{(a_1 \cdot \rfs{diam}(C) + b_1)}{\delta}$.
So $\tau(C)$ is bounded.

(c) By Part (b), if $f \in \rfs{ABUC}(X,X)$,
then $f \in \rfs{BUC}(X,X)$.
So $\rfs{ABUC}(X) = \rfs{BUC}(X)$.
\rule{0pt}{0pt}\hfill\myqed\bigskip

\noindent
{\bf Remark } Part (b) of the above proposition
follows trivially from Proposition \ref{metr-bldr-p4.3}(b).
However, the proof of \ref{metr-bldr-p4.3} was left to the reader.
\smallskip

Suppose that $\tau \in H(X,Y)$ and
$(\rfs{UC}(X))^{\tau} \subseteq \rfs{ABUC}(Y)$. Assuming
that $\tau$ is boundedness preserving,
the proof that $\tau \in \rfs{ABUC}(X,Y)$
is just as the proof of \ref{t4.10}.
This is the contents of the next lemma.
The main problem will be to deduce that $\tau$ is boundedness
preserving.

\begin{defn}\label{d4.13}\label{metr-bldr-d5.12}
\begin{rm}
Let $X$ be a metric space.
$X$ is
{\it boundedly uniformly\,-\,in\,-\,diameter arcwise\,-\,connected
($X$ is BUD.AC)},
   \index{boundedly uniformly\,-\,in\,-\,diameter arcwise\,-\,connected}
   \index{budac@@BUD.AC.
   Abbreviation of boundedly uniformly\,-\,in\,-\,diameter
   arcwise\,-\,connected}
if for every bounded set $B \subseteq X$ and $\varepsilon > 0$
there is $\delta >0$ such that for every $x,y \in B$:
if $d(x,y) < \delta$, then there is an arc $L \subseteq X$ connecting
$x$ and $y$ such that $\rfs{diam}(L) < \varepsilon$.
\end{rm}
\end{defn}

\begin{lemma}\label{l4.14}\label{metr-bldr-l5.13}
Let $X$ be BUD.AC,
and $\tau \in H(X,Y)$ be boundedness preserving.
Suppose that $(\rfs{UC}(X))^{\tau} \subseteq \rfs{BUC}(Y)$.
Then $\tau \in \rfs{BUC}(X,Y)$.
\end{lemma}

\noindent
{\bf Proof }
The proof of the lemma is the same as the proof of \ref{t4.10}.
\hfill\myqed\bigskip

The following example is a preparation for
Theorem~\ref{metr-bldr-t5.17}.
It shows that the assumptions of that theorem are ``correct''.

\begin{example}\label{e4.15}\label{metr-bldr-e5.14}
\begin{rm}
(a) Let $X = B^E(0,1) - \sngltn{0}$, $Y = E - \rfs{cl}(B^E(0,1))$, and 
$\tau(x) \eqdf \frac{x}{\norm{x}^2}$ be the inversion map from
$X$ to $Y$.
Then $(\rfs{BUC}(X))^{\tau} = \rfs{BUC}(Y)$, but $\tau$ is not ABUC. 
Note that $0 \in \rfs{bd}(X)$ and for every $h \in \rfs{BUC}(X)$,
$h^{\srfs{cl}}(0) = 0$. In Part (b) we get rid of this pathology.

(b) Let $X$, $Y$ and $\tau$ be as in part (a).
Let $X_1 = X \times \bbR$, $Y_1 = Y \times \bbR$
and $\tau_1(x,y) = (\tau(x),y)$.
Then $(\rfs{BUC}(X_1))^{\tau_1} \subseteq \rfs{BUC}(Y_1)$,
but $\tau_1$ is not ABUC.
In this example, $X$ does not have boundary points fixed
under $\rfs{BUC}(X)$, but we have containment and not equality
between $(\rfs{BUC}(X_1))^{\tau_1}$ and $\rfs{BUC}(Y_1)$.
\end{rm}
\end{example}

We next formulate the movability property of $X$, which will be  used
in the proof that $\tau$ is boundedness preserving.  
It is rather technical but it includes many open sets whose boundary is
not so well-behaved.

\begin{defn}\label{d4.16}\label{metr-bldr-d5.15}
\begin{rm}
For $\fnn{h}{[0,1] \times X}{X}$ and $t \in [0,1]$ we define
$h_t(x) \eqdf h(t,x)$.
We say that $X$ has {\it Property MV1},
if for every bounded $B \subseteq X$
there are $r = r_B > 0$ and $\alpha = \alpha_B \in \rfs{MC}$
such that for every $x \in B$ and $0 < s \leq r$,
there is an $\alpha$-continuous function $\fnn{h}{[0,1] \times X}{X}$
such that: 
(1) for every $t \in [0,1]$, $h_t \in H(X)$
and $h_t\inverse$ is $\alpha$-continuous;
(2) $h_0 = \rfs{Id}$ and $d(x,h_1(x)) = s$; and
(3) For every $t \in [0,1]$, $\rfs{supp}(h_t) \subseteq B(x,2s)$.
\hfill\proofend
   \index{property mv1@@Property MV1}
\end{rm}
\end{defn}

Note that if there is $x \in \rfs{bd}(X)$ such that $f(x) = x$
for every $f \in \rfs{BUC}(X)$, then $X$ does not have Property MV1.
On the other hand, Property MV1 holds for sets whose boundary is,
in a certain sense, well behaved.
Open half spaces, open balls, and complements of closed subspaces
fulfill MV1.

The following family of examples contains open sets $X$ such that
$\rfs{cl}(X)$ is not a manifold with boundary.
Let $U$ be any nonempty open subset of a normed space $E_0$
and $X = U \times \bbR$. Then $X$ has Property MV1.
More generally, $X$ has Property MV1 if the following happens.
Let $E_0$ be a normed space, $E = E_0 \times \bbR$,
$s > 0$ and $\alpha \in \rfs{MBC}$.
Suppose that $X$ is an open subset of $E$ with the following property.
For every $x \in \rfs{bd}(X)$ there are:
an open subset $U \subseteq E_0$,
$x_0 \in \rfs{bd}(U)$
and a homeomorphism $\varphi$
from $\overB^{E_0}(x_0,s) \times [-1,1]$ into $E$,
such that:
\begin{itemize}
\addtolength{\parskip}{-11pt}
\addtolength{\itemsep}{06pt}
\item[(1)] 
$\varphi(x_0,0) = x$,
\item[(2)] 
$\rfs{Rng}(\varphi)$ is closed in $E$,
and $\varphi(B^{E_0}(x_0,s) \times (-1,1))$ is open in $E$,
\item[(3)] 
$X \cap \rfs{Rng}(\varphi) =
\varphi((U \cap \overB^{E_0}(x_0,s)) \times [-1,1])$,
\item[(4)] 
$\varphi$ is $\alpha$-bicontinuous.
\vspace{-05.7pt}
\end{itemize}

\begin{prop}\label{p4.17}\label{metr-bldr-p5.16}
\num{a} Let $X$ be a metric space, $\alpha \in \rfs{MC}$
and $\setm{h_n}{n \in \bbN} \subseteq H(X)$.
Suppose that for every distinct $m,n \in \bbN$,
$h_m$ is $\alpha$-continuous
and $\rfs{supp}(h_m) \cap \rfs{supp}(h_n) = \emptyset$.
Then
$\bcirc_{n \in \sboldbbN} h_n$ is $\alpha \scirc \alpha$-continuous.

\num{b} Let $X$ be a subset of a normed space $E$,
$\alpha \in \rfs{MC}$ and $\setm{h_n}{n \in \bbN} \subseteq H(X)$.
Suppose that for every distinct $m,n \in \bbN$,
$h_m$ is $\alpha$-continuous,
$\rfs{cl}^E(\rfs{supp}(h_n)) \subseteq X$
and $\rfs{supp}(h_m) \cap \rfs{supp}(h_n) = \emptyset$.
Then $\bcirc_{n \in \sboldbbN} h_n$ is $2 \alpha$-continuous.
\end{prop}

\noindent
{\bf Proof }
(a) Denote $h = \bcirc_{n \in \sboldbbN} h_n$.
Let $x,y \in X$. Then there are $m,n \in \bbN$ such that\break
$x,y \in \rfs{supp}(h_m) \cup \rfs{supp}(h_n) \cup
(X - \bigcup_{i \in \sboldbbN} \rfs{supp}(h_i))$.
So $h(x) = h_m \scirc h_n(x)$ and
$h(y) = h_m \scirc h_n(y)$.
Since $h_m \scirc h_n$ is $\alpha \scirc \alpha$-continuous,
$d(h(x),h(y)) \leq \alpha \scirc \alpha(d(x,y))$.

(b) Denote $h = \bcirc_{n \in \sboldbbN} h_n$.
Let $x,y \in X$. Then there are $m,n \in \bbN$ such that\break
$x,y \in \rfs{supp}(h_m) \cup \rfs{supp}(h_n) \cup
(X - \bigcup_{i \in \sboldbbN} \rfs{supp}(h_i))$.
If $x$ or $y$ belong to $X - \bigcup_{i \in \sboldbbN} \rfs{supp}(h_i)$,
or $x,y \in \rfs{supp}(h_m)$, or $x,y \in \rfs{supp}(h_n)$,
then either
$d(h(x),h(y)) = d(h_m(x),h_m(y)) \leq \alpha(d(x,y))$,
or
$d(h(x),h(y)) = d(h_n(x),h_n(y)) \leq \alpha(d(x,y))$.

So we may assume that $x \in \rfs{supp}(h_m)$
and $y \in \rfs{supp}(h_n)$.
Let $z \in [x,y] \cap \rfs{bd}(\rfs{supp}(h_m))$.
Then $z \in X$ and $z \not\in \rfs{supp}(h_n)$.
Hence
$h_m(z) = h_n(z) = z$.
So
\newline\rule{3pt}{0pt}
\renewcommand{\arraystretch}{1.5}
\addtolength{\arraycolsep}{-0pt}
$
\begin{array}[b]{ll}
\renewcommand{\arraystretch}{1.0}
&
d(h(x),h(y)) \leq d(h(x),h(z)) + d(h(z),h(y)) =
d(h_m(x),h_m(z)) + d(h_n(z),h_n(y))
\\
\leq
&
\alpha(d(x,z)) + \alpha(d(z,y)) \leq 2 \alpha(d(x,y)).
\end{array}
$
\renewcommand{\arraystretch}{1.0}
\addtolength{\arraycolsep}{0pt}
\newline
\rule{10pt}{0pt}\hfill\kern-05pt\myqed

\begin{theorem}\label{t4.18}\label{metr-bldr-t5.17}
Let $X,Y \in K_{\srfs{NRM}}^{\calO}$.
Suppose that $X$ has Property MV1,
and let $\tau \in H(X,Y)$ be such that
$(\rfs{UC}(X))^{\tau} \subseteq \rfs{BUC}(Y) \subseteq
(\rfs{BUC}(X))^{\tau}$.
Then $\tau$ is boundedness preserving.
\end{theorem}

\noindent
{\bf Proof }
Suppose otherwise. Let $\vecx \subseteq X$ be a bounded sequence
such that $\tau(\vecx)$ is unbounded.
We may assume that either $\vecx$ is a Cauchy sequence or
$\vecx$ is spaced.

{\bf Case 1 }
$\vecx$ is a Cauchy sequence.
Applying MV1 to the bounded set $\rfs{Rng}(\vecx)$ we obtain
$r = r_{\srfs{Rng}(\vecx)} > 0$ and
$\alpha = \alpha_{\srfs{Rng}(\vecx)} \in \rfs{MC}$.
Set $x^* = \lim^{\oversE} \vecx$,
and choose $\delta > 0$
such that $\delta, \alpha(\delta) < \dgfrac{r}{4}$,
and $m$ such that $d(x_m,x^*) < \delta$.
Let $\fnn{h}{[0,1] \times X}{X}$ be the isotopy provided by MV1
when $x$ and $s$ are taken to be $x_m$ and $r$,
and let $\barh = h^{\srfs{cl}}_{[0,1] \times \oversE}$.
(See Definition \ref{metr-bldr-d4.6}).
From the fact that $h$ is $\alpha$-continuous it follows that
$\fnn{\barh}{\rfs{cl}^{[0,1] \times \oversE}([0,1] \times X)}
{\rfs{cl}^{\oversE}(X)}$
and $\barh$ is
$\alpha$-continuous.
Since $\barh_1$ is $\alpha$-continuous,
$d(\barh_1(x^*),\barh_1(x_m)) \leq \alpha(d(x^*,x_m)) <
\alpha(\delta) < \dgfrac{r}{4}$.
So\break
$d(x^*, \barh_1(x^*)) \geq d(x_m,\barh_1(x_m)) -
d(x_m,x^*) - d(\barh_1(x_m),\barh_1(x^*)) >
r - \dgfrac{r}{4} - \dgfrac{r}{4} = \dgfrac{r}{2}$.
That is, $d(x^*, \barh_1(x^*)) > \dgfrac{r}{2}$.
For $n \in \bbN$ define $L_n = h(x_n,[0,1])$.

{\bf Claim 1.}
$\limti{n} d(\tau(L_n),0) = \infty$.
{\bf Proof }
Suppose otherwise.
Then there are a $\onetoonen$ sequence $\setm{n_k}{k \in \bbN}$
and a sequence $\setm{t_k}{k \in \bbN} \subseteq [0,1]$ such that
$\setm{\tau(h(x_{n_k},t_k))}{k \in \bbN}$ is bounded.
We may assume that
$\setm{t_k}{k \in \bbN}$
converges to $t^*$.
Since $h_{t^*} \in \rfs{UC}(X)$,
$(h_{t^*})^{\tau} \in \rfs{BUC}(Y)$.
In particular, $(h_{t^*})^{\tau} \in \rfs{BDD.P}(Y)$.
It follows that
$\setm{\tau(h_{t^*}(x_{n_k}))}{k \in \bbN} =
(h_{t^*})^{\tau}(\setm{\tau(x_{n_k})}{k \in \bbN})$ is unbounded.
Let $I_k$ be the interval whose endpoints are $t_k$ and $t^*$
and $L'_k = h(I_k \times \sngltn{x_{n_k}})$.
By the $\alpha$-continuity of $h$, $\limti{k} \rfs{diam}(L'_k) = 0$.
Proceeding as in the proof of Case 1 of Theorem \ref{t4.10},
we construct a $\onetoonen$ sequence $\setm{k_i}{i \in \bbN}$ and
$g \in \rfs{UC}(X)$ such that
$g(h(t_{k_i},x_{n_{k_i}})) = h(t^*,x_{n_{k_i}})$. 
The fact that $g \in \rfs{UC}(X)$ implies that
$g^{\tau} \in \rfs{BUC}(Y)$, so in particular,
$g^{\tau}$ boundedness preserving.
However, $g^{\tau}$ takes the bounded sequence
$\tau(h(t_{k_i},x_{n_{k_i}}))$ to the unbounded
sequence $\tau(h(t^*,x_{n_{k_i}}))$. A contradiction, so Claim~1 is
proved.\smallskip

Let $u_n = h(1,x_n)$ and $U_n = B^Y(\tau(L_n),1)$.
There is a subsequence
$\setm{U_{n_k}}{k \in \bbN}$ of $\setm{U_n}{n \in \bbN}$
such that for every $k \in \bbN$,
$U_{n_k} \subseteq \dgfrac{B(0,d(0,U_{n_{k + 1}}))}{2}$.
For every $k \in \bbN$, let $g_k \in \rfs{UC}(Y)$ be such that
$\rfs{supp}(g_k) \subseteq U_{n_k}$ and
$g_k(\tau(x_{n_k})) = \tau(u_{n_k})$.
Let $g = \bcirc_{k \in \sboldbbN} g_{2k}$ and $f = g^{\tau\inverse}$.

Clearly, $g \in \rfs{BUC}(Y)$. So $f$ must belong to $\rfs{BUC}(X)$.
Note that
$\lim_{n \in \sboldbbN} u_n = \barh_1(x^*) \neq x^* =
\lim_{n \in \sboldbbN} x_n$.
So since
$f(x_{n_{2k}}) = u_{n_{2k}}$ and
$f(x_{n_{2k + 1}}) = x_{n_{2k + 1}}$,
$\setm{f(x_{n_k})}{k \in \bbN}$ is not convergent in $\overE$.
However, $\setm{x_{n_k}}{k \in \bbN}$ is convergent in $\overE$.
Hence $f$ takes a Cauchy sequence to a sequence which is not a Cauchy
sequence.
So $f \not\in \rfs{BUC}(X)$, a contradiction.\kern-2pt
\smallskip

{\bf Case 2 } $\vecx$ is spaced.
Let $r_0 > 0$ be such that $\vecx$ is $5r_0$-spaced.
Applying MV1 to the bounded set $\rfs{Rng}(\vecx)$ we obtain
$r_1 = r_{\srfs{Rng}(\vecx)} > 0$ and
$\alpha = \alpha_{\srfs{Rng}(\vecx)} \in \rfs{MC}$.
Let $s = \min(r_0,r_1)$.
For every $n \in \bbN$ let 
$\fnn{h_n}{[0,1] \times X}{X}$ be the function assured by MV1 for
$x_n$ and $s$. Recall that for $t \in [0,1]$,
$h_{n,t}(x)$ is the homeomorphism of $X$ taking every $x \in X$
to $h_n(t,x)$.
Set $L_n = h_n([0,1] \times \sngltn{x_n})$.

{\bf Claim 2.}
$\limti{n} d(\tau(L_n),0) = \infty$.
{\bf Proof }
Suppose otherwise.
Then there are a $\onetoonen$ sequence $\setm{n_k}{k \in \bbN}$
and a sequence $\setm{t_k}{k \in \bbN} \subseteq [0,1]$ such that
$\setm{\tau(h_{n_k}(t_k,x_{n_k}))}{k \in \bbN}$ is bounded.
Clearly, for every distinct $m,n \in \bbN$ and $q,t \in [0,1]$,
{\thickmuskip=2mu \medmuskip=1mu \thinmuskip=1mu 
\hbox{$d(\rfs{supp}(h_{m,q}),\rfs{supp}(h_{n,t})) \geq r_0$.}
}
So by \ref{p4.17}(a),
$f \eqdf \bcirc_{k \in \sboldbbN} h_{n_k,t_k} \in \rfs{UC}(X)$.
So $f^{\tau} \in \rfs{BUC}(Y) \subseteq \rfs{BDD.P}(Y)$.
We shall reach a contradiction by showing that $f^{\tau}$
takes an unbounded sequence to a bounded sequence.
$\setm{\tau(x_{n_k})}{k \in \bbN}$ is unbounded, whereas
$f^{\tau}(\setm{\tau(x_{n_k})}{k \in \bbN}) =
\setm{\tau(h_{n_k}(t_k,x_{n_k}))}{k \in \bbN}$ is bounded.
Claim 2 is thus proved.\smallskip

Let $u_n = h_n(1,x_n)$, $v_n = h_n(\dgfrac{1}{n},x_n)$ 
and $U_n = B^Y(\tau(L_n),1)$.
There is a subsequence
$\setm{U_{n_k}}{k \in \bbN}$ of $\setm{U_n}{n \in \bbN}$
such that for every $k \in \bbN$,
$U_{n_k} \subseteq \dgfrac{B(0,d(0,U_{n_{k+1}}))}{2}$.
For every $k \in \bbN$, let $g_k \in \rfs{UC}(Y)$ be such that
$\rfs{supp}(g_k) \subseteq U_{n_k}$,
$g_k(\tau(x_{n_k})) = \tau(x_{n_k})$ and
$g_k(\tau(v_{n_k})) = \tau(u_{n_k})$.
Let \hbox{$g = \bcirc_{k \in \sboldbbN} g_{k}$}
and $f = g^{\tau\inverse}$.

Clearly, $g \in \rfs{BUC}(Y)$. So $f$ must belong to $\rfs{BUC}(X)$.
By the $\alpha$-continuity of all $h_n$'s,
$\limti{k} d(x_{n_k},v_{n_k}) = 0$, whereas for every $k \in \bbN$,
$d(f(x_{n_k}),f(v_{n_k})) = d(x_{n_k},u_{n_k}) = s$. So
$f \not\in \rfs{BUC}(X)$, a contradiction.
\smallskip\hfill\myqed

\rule{0pt}{0pt}\kern-3pt
Recall the convention that $X$ and $Y$
\hbox{denote open subsets of the normed spaces $E$ and $F$.}

\begin{cor}\label{c4.19}\label{metr-bldr-c5.18}
Let $X,Y \in K_{\srfs{NRM}}^{\calO}$.
Suppose that $X$ is BUD.AC, and $X$ has Property MV1.
Let $\tau \in H(X,Y)$ be such that
$(\rfs{UC}(X))^{\tau} \subseteq \rfs{BUC}(Y)$
and
$(\rfs{BUC}(Y))^{\tau\inverse} \subseteq \rfs{BUC}(X)$.
Then $\tau \in \rfs{BUC}(X,Y)$.
\end{cor}

\noindent
{\bf Proof }
Combine Lemma \ref{l4.14} and Theorem \ref{t4.18}.
\hfill\myqed

\kern1.7mm

The following Theorem is the final result for groups
of type $\rfs{BUC}(X)$.

\begin{theorem}\label{metr-bldr-t5.19}
Let $X,Y \in K_{\srfs{NFCB}}^{\calO}$.
Suppose that $X$ and $Y$ are BUD.AC, and $X$ and $Y$ have Property MV1.
Let $\iso{\varphi}{\rfs{BUC}(X)}{\rfs{BUC}(Y)}$.
Then there is $\tau \in \rfs{BUC}^{\pm}(X,Y)$ such that
$\tau$ induces $\varphi$.
\end{theorem}

\noindent
{\bf Proof }
Combine Corollaries \ref{metr-bldr-c2.26} and \ref{metr-bldr-c5.18}.
\hfill\myqed

\subsection{Groups of homeomorphisms which are uniformly
continuous on every bounded positively distanced set.}
\label{ss5.4}
\label{ss5.4-BPD-homeomorphisms}

We next deal with the group $\rfs{BPD.UC}(X)$ and with some related
groups.
Recall that $\rfs{BPD.UC}(X)$ is the group of all homeomorphisms
$f$ such that $f$ and $f\inverse$ take every subset of $X$ whose
distance from the boundary of $X$ is positive to a set whose distance
from the boundary of $X$ is positive,
and such that $f$ and $f\inverse$ are uniformly continuous
on every such set.
The generalization of $\rfs{BPD.UC}(X)$ to arbitrary moduli of
continuity is denoted by $H_{\itGamma}^{\srfs{NBPD}}(X)$.
That is, $\rfs{BPD.UC}(X)$ is the group
$H_{\itGamma}^{\srfs{NBPD}}(X)$ when $\itGamma = \rfs{MC}$.
These groups are explained in the next definition.
The final reconstruction result for such groups appears in
Theorem~\ref{metr-bldr-t5.30},
and this result is obtained for countably generated $\itGamma$'s and
for $\itGamma = \rfs{MC}$.
The main intermediate result for countably generated $\itGamma$'s
appears in Theorem~\ref{metr-bldr-t5.23}(b),
and it says that if
$(\rfs{LIP}_{00}(X))^{\tau} \subseteq H_{\itGamma}^{\srfs{NBPD}}(X)$,
then $\tau \in H_{\itGamma}^{\srfs{NBPD}}(X,Y)$.
The intermediate result fot $\itGamma = \rfs{MC}$
appears in Theorem~\ref{metr-bldr-t5.29}.
The analogous statement here is:
if $(\rfs{UC}_{00}(X))^{\tau} \subseteq \rfs{BPD.UC}(Y)$,
then $\tau \in \rfs{BPD.UC}(X,Y)$.
The groups $\rfs{LIP}_{00}(X)$ and $\rfs{UC}_{00}(X)$ are defined in
\ref{metr-bldr-d5.22}.

For open subsets of a Banach space we can also conclude that
$\tau\inverse \in \rfs{BPD.UC}(X,Y)$.
That is,
if $(\rfs{BUC}(X))^{\tau} \subseteq \rfs{BPD.UC}(Y)$,
then $\tau\inverse \in \rfs{BPD.UC}(Y,X)$.
This is done in Theorem~\ref{t4.27}(a).

A weaker variant of uniform continuity pops up,
and is also dealt with.
Groups arising from this variant are defined in
\ref{metr-bldr-d5.20}(c) and are denoted by
$H_{\itGamma}^{\srfs{WBPD}}(X)$.
The final result for such groups is stated in Theorem
\ref{metr-bldr-t5.35}.
The main intermediate results for such groups appear in
Theorem \ref{metr-bldr-t5.23}(a) and
Proposition \ref{metr-bldr-p5.34}.

\smallskip

We next define the groups
$H_{\itGamma}^{\srfs{BPD}}(X)$, $H_{\itGamma}^{\srfs{NBPD}}(X)$ and
$H_{\itGamma}^{\srfs{WBPD}}(X)$.

\begin{defn}\label{d-gamma.1}\label{metr-bldr-d5.20}
\begin{rm}
(a) Define
$$
H_{\itGamma}^{\srfs{BPD}}(X,Y) =
\setm{f \in \rfs{BPD.P}(X,Y)}{\kern1pt\mbox{for every }
\mbox{BPD set }
A \subseteq X,\ f \nrestriction A \mbox{ is }
\itGamma\mbox{-continuous}}.
$$

   \index{N@hbpd00@@$H_{\itGamma}^{\srfs{BPD}}(X,Y) =
          \setm{f \in \rfs{BPD.P}(X,Y)}{\kern1pt
	  \mbox{for every BPD set } A \subseteq X,
          f \nrestriction A \mbox{ is } \itGamma\mbox{-continuous}}$}

(b) Let $\itGamma$ be a modulus of continuity and $\fnn{f}{X}{Y}$.
We say that $f$ is {\it nearly $\itGamma$-continuous on BPD sets},
if for every $\rfs{BPD}$ set $A \subseteq X$
there are $\alpha \in \itGamma$
and $r > 0$ such that $f \nrestriction A$ is $(r,\alpha)$-continuous.
See Definition \ref{metr-bldr-d4.2}(b).\ 
We denote by $H_{\itGamma}^{\srfs{NBPD}}(X,Y)$
the set of all $h \in \rfs{BPD.P}(X,Y)$
such that $h$ is nearly $\itGamma$-continuous on BPD sets.
   \index{nearly $\itGamma$-continuous on BPD sets}
   \index{N@hnbpd@@$H_{\itGamma}^{\srfs{NBPD}}(X,Y) =
\setm{h \in BPD.P(X,Y)}{h \mbox{ is nearly $\itGamma$-continuous}
\mbox{ on BPD sets}}$}

(c) Let $\alpha \in \rfs{MC}$,
and $\fnn{f}{X}{Y}$ be a function between metric spaces.
Recall that according to Definition \ref{d1.6}(a), 
$f$ is locally $\sngltn{\alpha}$-continuous,
if for every $x \in X$ there is $U \in \rfs{Nbr}^X(x)$ such that
$f \nrestriction U$ is $\alpha$-continuous.
Let $\fnn{f}{X}{Y}$ be a function between metric spaces and
$\itGamma$ be a modulus of continuity.
Call $f$ {\it weakly $\itGamma$-continuous},
if there is $\alpha \in \itGamma$ such that
$f$ is locally $\sngltn{\alpha}$-continuous.
If $f \in H(X,Y)$ and both $f$ and $f\inverse$ are
weakly $\itGamma$-continuous, the $f$ is said to be
{\it weakly $\itGamma$-bicontinuous}.

Let $X$ and $Y$ be open subsets of normed spaces $E$ and $F$
respectively, $\itGamma$ be a modulus of continuity
and $\fnn{f}{X}{Y}$.
Call $f$ {\it weakly $\itGamma$-continuous on BPD sets},
if for every $\rfs{BPD}$ set $A \subseteq X$,
$f \nrestriction A$ is weakly $\itGamma$-continuous.
We denote by $H_{\itGamma}^{\srfs{WBPD}}(X,Y)$
the set of all $h \in \rfs{BPD.P}(X,Y)$
such that $h$ is weakly $\itGamma$-continuous on BPD sets.

   \index{weakly $\itGamma$-continuous function}
   \index{weakly $\itGamma$-bicontinuous function}
   \index{weakly $\itGamma$-continuous on BPD sets}
   \index{N@hwbpd@@$H_{\itGamma}^{\srfs{WBPD}}(X,Y) =
          \setm{h \in BPD.P(X,Y)}{h
	  \mbox{ is weakly $\itGamma$-continuous} \mbox{ on BPD sets}}$}

(d) Let $X$ a subset of a metric space $E$.
$X$ has
{\it the discrete path property for BPD sets},
if for every BPD subset $A \subseteq X$
there are $d > 0$ and $K \geq 1$ such that for every
$x,y \in A$ and $r > 0$ there are $n \in \bbN$ and
$x = x_0,\ldots,x_n = y \in X$ such that
$n \leq K \mcdot \frac{d(x,y)}{r}$,
and for every $i < n$,
$\delta(x_i) > d$ and $d(x_i,x_{i + 1}) \leq r$.
   \index{discrete path property for BPD sets}
\smallskip\hfill\proofend
\end{rm}
\end{defn}

Note that
$H_{\itGamma}^{\srfs{BPD}}(X)$, $H_{\itGamma}^{\srfs{NBPD}}(X)$ and
$H_{\itGamma}^{\srfs{WBPD}}(X)$ are groups.
It is easy to check that for $X$'s which are open subsets
of a finite-dimensional normed space,
$X$ has the discrete path property for BPD sets
iff $X$ is connected.
For infinite-dimensional normed spaces
neither of the above implications is true.
In any case, ``well-behaved'' open subsets of a normed space have the
discrete path property for BPD sets.
For example, an open ball has this property.
We first observe the following easy facts.
Part (a) follows from Proposition~\ref{metr-bldr-p4.3}(a),
and the proof of (b) is left to the reader.

\begin{prop}\label{metr-bldr-p5.21}
\num{a} $\rfs{BPD.UC}(X) = H_{\srfs{MC}}^{\srfs{NBPD}}(X)$.

\num{b} Suppose that $X$ has the discrete path property for BPD sets.
Then
$H_{\itGamma}^{\srfs{BPD}}(X) = H_{\itGamma}^{\srfs{NBPD}}(X)$.
\end{prop}


\begin{defn}\label{d-gamma.2}\label{n-gamma.3}\label{metr-bldr-d5.22}
\begin{rm}
(a) $X$ is {\it BPD-arcwise-connected ($X$ is BPD.AC)},
if for every BPD set $A \subseteq X$ there are
$C,D > 0$ such that for every $x,y \in A$ there is a rectifiable arc
$L \subseteq X$ connecting $x$ and $y$ such that
$\rfs{lngth}(L) \leq D$ and $\delta(L) \geq C$.

   \index{bpd-arcwise-connected@@BPD-arcwise-connected}
   \index{bpdac@@BPD.AC. Abbreviation of BPD-arcwise-connected}

(b) In some of the subsequent lemmas it will be
convenient to regard a sequence as a function whose domain
is an infinite subset of $\bbN$. So if $\sigma \subseteq \bbN$
is infinite, then the object $\setm{x_i}{i \in \sigma}$
is considered to be a sequence.
The notions of a subsequence, a convergent sequence etc.\ are 
easily modified to fit into this setting.

(c) Let
$\rfs{LIP}_{00}(X) = \setm{h \in \rfs{LIP}(X)}
{\rfs{supp}(h) \mbox{ is a BPD set}}$
and
$\rfs{UC}_{00}(X) = \setm{h \in \rfs{UC}(X)}
{\rfs{supp}(h) \mbox{ is a BPD set}}$.

   \index{N@lip10@@$\rfs{LIP}_{00}(X) =
   \setm{f \in \rfs{LIP}(X)}{\rfs{supp}(f)
   \mbox{ is a BPD set}}$}
   \index{N@uc08@@$\rfs{UC}_{00}(X) =
   \setm{f \in \rfs{UC}(X)}{\rfs{supp}(f)
   \mbox{ is a BPD set}}$}

(d) For $x \in X$
let $\delta_1^X(x) = \max(\norm{x},\dgfrac{1}{\delta^X(x)})$.
We abbreviate $\delta_1^X(x)$ by $\delta_1(x)$.

   \index{N@d04@@$\delta_1^X(x)$}

(e) Let $A \subseteq \bbN$ and $n \in \bbN$.
Define $A^{\geq n} = \setm{m \in A}{m \geq n}$.
The notations $A^{> n}$, $A^{\leq n}$,
$A^{< n}$ etc.\ are defined analogously.
\hfill\proofend
   \index{N@AAAA@@$A^{\geq n} = \setm{m \in A}{m \geq n}$}
   \index{N@AAAA@@$A^{> n} = \setm{m \in A}{m > n$}}
   \index{N@AAAA@@$A^{\leq n} = \setm{m \in A}{m \leq n$}}
   \index{N@AAAA@@$A^{< n} = \setm{m \in A}{m < n$}}
\end{rm}
\end{defn}

Note that if $X$ is BPD.AC, then $X$ is connected.
Note that a subset $A \subseteq X$ is BPD iff
$\sup(\setm{\delta_1^X(x)}{x \in A}) < \infty$.

\begin{theorem}\label{t-gamma.4}\label{metr-bldr-t5.23}
Let $\itGamma$ be a countably generated modulus of continuity.
Suppose that $X$ and $Y$ are open subsets of
normed spaces $E$ and $F$ respectively,
$X$ is BPD.AC and $\tau \in H(X,Y)$.

\num{a} If
$(\rfs{LIP}_{00}(X))^{\tau} \subseteq H_{\itGamma}^{\srfs{WBPD}}(Y)$,
then $\tau \in H_{\itGamma}^{\srfs{WBPD}}(X,Y)$.

\num{b} If
$(\rfs{LIP}_{00}(X))^{\tau} \subseteq H_{\itGamma}^{\srfs{NBPD}}(Y)$,
then $\tau \in H_{\itGamma}^{\srfs{NBPD}}(X,Y)$.
\end{theorem}

Variants of the argument appearing in Claims 3
will be used in several other proofs.

\begin{lemma}\label{l-gamma.5}\label{metr-bldr-l5.24}
Suppose that $X$ is BPD.AC, $\tau \in H(X,Y)$ and
$(\rfs{LIP}_{00}(X))^{\tau} \subseteq \rfs{BPD.P}(Y)$. Then 
$\tau \in \rfs{BPD.P}(X,Y)$.
\end{lemma}

\noindent
{\bf Proof } Let $X,Y$ and $\tau$ be as in the lemma.

{\bf Claim 1.} Suppose that $u \in X$, $0 < r < s$,
$B(u,s) \subseteq X$ and $\vecx \subseteq B(u,r)$.
Then $\tau(\vecx)$ is BPD in $Y$.
{\bf Proof } Suppose by contradiction that $\tau(\vecx)$
is not BPD in $Y$.
Let $a \in (0,1)$ be such that $\tau(B(u,ar))$ is BPD in $Y$.
Let $\fnn{\eta}{[0,\infty)}{[0,\infty)}$
be the piecewise linear function
with breakpoints at $ar$ and $\dgfrac{(r + s)}{2}$ such that
$\eta(ar) = r$ and for every $t \geq \dgfrac{(r + s)}{2}$,
$\eta(t) = t$.
Let
$h = \rfs{Rad}^E_{\eta,u} \nrestriction X$.
(See Definition \ref{metr-bldr-d3.17}(b)).
Then $h \in \rfs{LIP}_{00}(X)$.
Let $\vecv = h\inverse(\vecx)$.
Clearly, $\vecv \subseteq B(u,ar)$.
So $\tau(\vecv)$ is BPD in $Y$.
Obviously,
$h^{\tau}(\tau(\vecv)) = \tau(\vecx)$.
Hence $h^{\tau}$ takes a BPD set to a set which is not BPD.
That is, $h^{\tau} \not\in \rfs{BPD.P}(Y)$, a contradiction.

{\bf Claim 2.} If $\vecx$ is a BPD sequence in $X$ and
$\vecx$ is a Cauchy sequence, $\tau(\vecx)$ is a BPD sequence in $Y$.
{\bf Proof} 
Suppose by contradiction that $\vecx$ is a counter-example.
Let $x^* = \lim^{\oversE}(\vecx)$.
Clearly, $x^* \in \overfs{int}(X)$.
Let $u \in X$ and $r > 0$ be such that
$x^* \in B^{\oversE}(u,r)$ and $B^E(u,2r) \subseteq X$.
Let $\vecy$ be a final segment of $\vecx$ such that
$\vecy \subseteq B(u,r)$.
Then $\vecy$ is a counter-example to Claim~1.
This proves Claim 2.

Suppose by contradiction that $\tau \not\in \rfs{BPD.P}(X,Y)$.
Let $\vecx$ be a BPD $\onetoonen$ sequence such that
$\tau(\vecx)$ is not BPD.
We may assume that
$\limti{n} \delta_1(\tau(x_n)) = \infty$.
Hence for every subsequence $\vecy$ of $\vecx$,
$\tau(\vecy)$ is not BPD.

It follows from Claim 2 that $\vecx$ has no Cauchy subsequences.
Let $x^* \in X - \rfs{Rng}(\vecx)$ and
$A = \rfs{Rng}(\vecx) \cup \sngltn{x^*}$.
Let $C$ an $D$ be as assured by the property BPD.AC.
For every $n \in \bbN$ let $L_n \subseteq X$ be a rectifiable arc
connecting $x^*$ and $x_n$
such that $\delta(L_n) \geq C$ and $\rfs{lngth}(L_n) \leq D$.
Note that $\bigcup_{n \in \sboldbbN} L_n$ is a BPD set.
Let $\fnn{\gamma_n}{[0,1]}{L_n}$
be a parametrization of $L_n$ such that $\gamma_n(0) = x^*$,
$\gamma_n(1) = x_n$, and for every $t \in [0,1]$,
$\rfs{lngth}(\gamma_n([0,t])) = t \cdot \rfs{lngth}(L_n)$.

For every infinite $\sigma \subseteq \bbN$ and $t \in [0,1]$ let
$A[\sigma,t] = \setm{\gamma_n(t)}{n \in \sigma}$.
We regard $A[\sigma,t]$ as a sequence whose domain is $\sigma$.
Clearly, for every $t \in [0,1]$,
$A[\bbN,t] \subseteq \rfs{cl}(B(x^*,tD))$.
So by the continuity of $\tau$, there is $t_0 > 0$ such that
for every $t \in [0,t_0]$, and
$\sigma \subseteq \bbN$, $\tau(A[\sigma,t])$ is a BPD set.
For every infinite $\sigma \subseteq \bbN$ let
$s_{\sigma} = \inf(\setm{t \in [0,1]}{\tau(A[\sigma,t])
\mbox{ is not a BPD set}})$. So $s_{\sigma} > 0$.

For $\sigma,\eta \subseteq \bbN$ let 
$\sigma \almostcontained \eta$ mean that $\sigma - \eta$ is finite.

{\bf Claim 3.} There is an infinite $\sigma \subseteq \bbN$ such that
for every infinite $\eta \subseteq \sigma$, $s_{\eta} = s_{\sigma}$.
{\bf Proof }
Suppose by contradiction that no such $\sigma$ exists.
Clearly if $\eta \almostcontained \sigma$, then
$s_{\eta} \geq s_{\sigma}$. We define by transfinite induction on
$\nu < \omega_1$ an infinite subset $\sigma_{\nu} \subseteq \bbN$ 
such that for every $\nu < \mu$:
$\sigma_{\mu} \almostcontained \sigma_{\nu}$ and
$s_{\sigma_{\mu}} > s_{\sigma_{\nu}}$.
If $\sigma_{\nu}$ has been defined, let 
$\sigma_{\nu + 1} \subseteq \sigma_{\nu}$ be such that 
$s_{\sigma_{\nu + 1}} > s_{\sigma_{\nu}}$.
If $\mu$ is a limit ordinal, and $\sigma_{\nu}$ has been defined
for every $\nu < \mu$, let $\sigma_{\mu}$ be an infinite set such
that for every $\nu < \mu$,
$\sigma_{\mu} \almostcontained \sigma_{\nu}$. 
By the induction hypthesis, if $\nu < \mu$, then
$s_{\sigma_{\nu + 1}} > s_{\sigma_{\nu}}$.
Hence $s_{\sigma_{\mu}} \geq s_{\sigma_{\nu}} > s_{\sigma_{\nu}}$.
So the induction hypothesis holds.
The set $\setm{s_{\sigma_{\nu}}}{\nu < \omega_1}$ is a subset of
$\bbR$ order isomorphic to $\omega_1$, a contradiction. Claim 3 is
proved.\smallskip

Let $\sigma$ be as assured by Claim 3 and write $s = s_{\sigma}$.

{\bf Claim 4.} $A[\sigma,s]$ does not have Cauchy subsequences.
{\bf Proof } Suppose by contradiction that $\eta \subseteq \sigma$ is
infinite, and $A[\eta,s]$ is a Cauchy sequence.
Since $A[\bbN,1] = \vecx$ does not contain Cauchy subsequences, $s < 1$.
Let $\hatx = \lim A[\eta,s]$.
Since $A[\eta,s]$ is a BPD sequence $\hatx \in \overfs{int}(X)$.
So there are $u \in X$ and 
$r > 0$ such that 
$\hatx \in B^{\oversE}(u,r)$ and $B^E(u,3r) \subseteq X$.
We may assume that $A[\eta,s] \subseteq B(u,r)$.
For every $i$ and $t$,
$\norm{\gamma_i(t) - \gamma_i(s)} \leq (t - s) \mcdot D$.
So for every $t \in (s,s + \dgfrac{r}{D})$,
$A[\eta,t] \subseteq B(u,2r)$.
By the definition of $\sigma$, $s_{\eta} = s_{\sigma} = s$.
So there is $t \in (s,s + \dgfrac{r}{D})$ such that
$\tau(A[\eta,t])$ is not a BPD subset of $Y$.
But $A[\eta,t] \subseteq B(u,2r)$ and $B(u,3r) \subseteq X$.
This contradicts Claim 1.
So Claim~4 is proved.

By Proposition \ref{p4.7}(a) and Claim 4, we may assume that 
there is $d > 0$ such that $A[\sigma,s]$ is $d$-spaced.
Let $r = \dgfrac{\min(C,d)}{4}$.
$\delta(A[\sigma,s]) \geq C$,
and so $B^E(A[\sigma,s],r) \subseteq X$ and
$\delta(B^E(A[\sigma,s],r)) > 0$.
Also for every distinct $m,n \in \sigma$,
$d(B(\gamma_m(s),r),B(\gamma_n(s),r)) \geq \dgfrac{d}{2}$.
Let $t_1 \in (s - \frac{r}{2D},s)$.
Since $t_1 < s$, it follows that $\tau(A[\sigma,t_1])$ is a BPD set.
Let $t_2 \in [s,s + \frac{r}{2D})$ be such that
$\tau(A[\sigma,t_2])$ is not a BPD set.

By Lemma \ref{l2.6}(b),
there is $K \geq 1$ such that for every normed space~$E$,
$u \in E$, $r > 0$ and
$x,y \in B^E(u,\dgfrac{r}{2})$ there is $h \in H(E)$ such that
$h$ is $K$-bilipschitz,
$\rfs{supp}(h) \subseteq B^E(u,r)$ and $h(x) = y$.

Clearly, for every $n \in \sigma$,
$\gamma_n(t_1),\gamma_n(t_2) \in B(\gamma_n(s),\dgfrac{r}{2})$.
So by the above fact, there is  $h_n \in H(X)$ such that
$h_n$ is $K$-bilipschitz,
$\rfs{supp}(h_n) \subseteq B(\gamma_n(s),r)$ and
$h_n(\gamma_n(t_1)) = \gamma_n(t_2)$.

By Proposition \ref{p4.17}(b),
$h \eqdf \bcirc_{n \in \sigma} h_n \in \rfs{LIP}(X)$.
Since $\rfs{supp}(h) \subseteq B^E(A[\sigma,s],r)$,
and $\delta(B^E(A[\sigma,s],r)) > 0$,
$h \in \rfs{LIP}_{00}(X)$.
Hence $h^{\tau} \in \rfs{BPD.P}(Y)$.
However, 
$\tau(A[\eta,t_1])$ is a BPD set,
$\tau(A[\eta,t_2])$ is not a BPD set,
and $h^{\tau}(\tau(A[\eta,t_1])) = \tau(A[\eta,t_2])$.
A contradiction.\break
\rule{0pt}{0pt}\medskip\hfill\myqed

\begin{prop}\label{p-gamma.6}\label{metr-bldr-p5.25}
For a compact metric space $C$ and $t > 0$
let $\nu_C(t)$ denote the minimal cardinality of a cover of $C$
consisting of subsets of $C$ whose diameter $\leq t$.

Let $\vecC = \setm{C_i}{i \in \bbN}$ be a sequence of compact
subsets of a metric space $X$, and let
$\fnn{\nu}{(0,\infty)}{\bbN}$.
Suppose that for every $i \in \bbN$,
$\nu_{C_i} \leq \nu$.
Suppose further that there is no infinite set
$\eta \subseteq \bbN$ and a sequence
$\setm{c_i}{i \in \eta}$ such that for every $i \in \eta$,
$c_i \in C_i$,
and $\setm{c_i}{i \in \eta}$ is a Cauchy sequence.
Then there is a subsequence $\vecD$ of
$\vecC$ such that $\vecD$ is spaced.
\end{prop}

\noindent
{\bf Proof } Suppose that $\vecC$ has no spaced subsequences,
and we show that there are an infinite set $A \subseteq \bbN$
and a Cauchy sequence $\vecc = \setm{c_i}{i \in A}$ such that for every
$i \in A$, $c_i \in C_i$.
There are a subsequence $\vecC^1$ of $\vecC$
and $r \in \bbR \cup \sngltn{\infty}$
such that $\lim_{i,j \rightarrow \infty}d(C^1_i,C^1_j) = r$.
Since $\vecC$ has no spaced subsequences, $r = 0$.
We may assume that $\vecC = \vecC^1$.

For $\vecp \subseteq \bbN$ let $T_{\vecp}$ be the tree of finite
sequences $\vecn$ such that for every $i < \rfs{lngth}(\vecn)$,
$n_i < p_i$.
Let $S_{\vecp} = \prod_{i \in \sboldbbN} \bbN^{< p_i}$.

Let $p_i = \prod_{j \leq i} \nu(\dgfrac{1}{j})$, $T = T_{\vecp}$
and $S = S_{\vecp}$.
Then for every $i \in \bbN$
there is $\setm{C_{i,\vecn}}{\vecn \in T}$
such that for every $\vecn \in T$,
$C_{i,\vecn}$ is closed
and $\rfs{diam}(C_{i,\vecn}) \leq \dgfrac{1}{\rfs{lngth}(\vecn)}$;
for every $\ell \in \bbN$,\break
$C_i =
\bigcup \setm{C_{i,\vecn}}
{\vecn \in T \mbox{ and } \rfs{lngth}(\vecn) = \ell}$;
and for every $\vecm,\vecn \in T$: if $\vecm$ is an initial
segment of $\vecn$, then $C_{i,\vecn} \subseteq C_{i,\vecm}$.

By Ramsey Theorem,
there are a sequence of infinite subsets of $\bbN$,
$A_0 \supseteq A_1 \supseteq\ldots\ $ and $\vecq,\vecr \in S$
such that for every $\ell$ and $i,j \in A_{\ell}$:
if $i < j$, then
$d(C_{i,\vecq \rharp \sboldbbN^{\leq \ell}},
C_{j,\vecr \rharp \sboldbbN^{\leq \ell}}) = d(C_i,C_j)$.

Let $A \subseteq \bbN$ be an infinite set such that for every
$i$, $A - A_i$ is finite.
For every $i \in A$ let
$D_i = \bigcap_{j \in \sboldbbN} C_{i,\vecq \rharp \sboldbbN^{\leq j}}$
and
$E_i = \bigcap_{j \in \sboldbbN} C_{i,\vecr \rharp \sboldbbN^{\leq j}}$.
Clearly, $D_i,E_i$ are singletons,
denote them by $x_i$ and $y_i$ respectively.
We check that
$\lim_{i \rightarrow \infty, i < j} d(x_i,y_j) = 0$.
Let $\varepsilon > 0$.
Then there is $N_1$ such that for every $i,j > N_1$,
$d(C_i,C_j) < \dgfrac{\varepsilon}{3}$.
Let $N_2$ be such that $\dgfrac{1}{N_2} < \dgfrac{\varepsilon}{3}$,
$N_3$ be such that $A^{\geq N_3} \subseteq A_{N_2}$
and $N = \max(N_1,N_3)$.
Let $i < j$ and $i,j \in A^{\geq N}$.
So $i,j \in A_{N_2}$. Hence
$d(C_{i,\vecq \rharp \sboldbbN^{\leq N_2}},
C_{i,\vecr \rharp \sboldbbN^{\leq N_2}}) =
d(C_i,C_j) < \dgfrac{\varepsilon}{3}$.
It follows that
$$
d(x_i,y_j) \leq
\rfs{diam}(C_{i,\vecq \rharp \sboldbbN^{\leq N_2}}) +
d(C_i,C_j) + \rfs{diam}(C_{j,\vecr \rharp \sboldbbN^{\leq N_2}}) <
\dgfrac{\varepsilon}{3} + \dgfrac{\varepsilon}{3} +
\dgfrac{\varepsilon}{3} = \varepsilon.
$$
We have proved that
$\lim_{i \rightarrow \infty, i < j} d(x_i,y_j) = 0$.
Let $\varepsilon > 0$.
Choose $N$ such that for every $i,j \in A^{\geq N}$: if $i < j$,
then $d(x_i,y_j) < \dgfrac{\varepsilon}{2}$.
Suppose that $i_1,i_2 \in A^{\geq N}$ and let $j$ be such that
$i_1,i_2 < j \in A$.
Then
$d(x_{i_1},x_{i_2})
\leq d(x_{i_1},y_j) + d(y_j,x_{i_2}) < \varepsilon$.
So $\setm{x_i}{i \in A}$ is a Cauchy sequence.
\rule{0pt}{1pt}\hfill\myqed

\begin{lemma}\label{l-gamma.7}\label{metr-bldr-l5.26}
There is $K_{\srfs{arc}}(\ell,t) > 0$
such that for every normed space $E$, $L,r > 0$,
and a rectifiable arc $\gamma \subseteq E$ with endpoints
$x,y$:
if $\rfs{lngth}(\gamma) \leq L$, then there is
$h \in H(E)$ such that:
\begin{itemize}
\addtolength{\parskip}{-11pt}
\addtolength{\itemsep}{06pt}
\item[\num{1}] 
$\rfs{supp}(h) \subseteq B(\gamma,r)$;
\item[\num{2}] 
$h \nrestriction B(x,\dgfrac{r}{2}) =
\rfs{tr}_{y - x} \nrestriction B(x,\dgfrac{r}{2})$;
\item[\num{3}] 
$h$ is $K_{\srfs{arc}}(L,r)$-bilipschitz.
\vspace{-05.7pt}
\end{itemize}
   \index{N@karc@@$K_{\srfs{arc}}(\ell,t)$}
\end{lemma}

\noindent
{\bf Proof }
Let $n = [\frac{L}{\dgfrac{r}{2}}] + 1$.
Suppose that $\fnn{\gamma}{[0,1]}{X}$.
There are
$0 = t_0, t_1,\ldots,t_n = 1$
such that for every $i < n$,
$\rfs{lngth}(\gamma \nrestriction [t_i,t_{i + 1}]) < \dgfrac{r}{2}$.
Let $x_i = \gamma(t_i)$.
Then for every $z \in [x_i,x_{i + 1}]$,
$d(z,\gamma \nrestriction [t_i,t_{i + 1}]) < \dgfrac{r}{4}$.
So $\bigcup_{i < n} B([x_i,x_{i + 1}],\dgfrac{3r}{4}) \subseteq
B(\gamma,r)$.
By Lemma \ref{l2.6}(c), there are $h_1,\ldots,h_n \in H(E)$
such that for every $i = 1,\ldots,n$:
\begin{list}{}
{\setlength{\leftmargin}{35pt}
\setlength{\labelsep}{08pt}
\setlength{\labelwidth}{20pt}
\setlength{\itemindent}{-00pt}
\addtolength{\topsep}{-04pt}
\addtolength{\parskip}{-02pt}
\addtolength{\itemsep}{-05pt}
}
\item[(1.1)] 
$\rfs{supp}(h_i) \subseteq B([x_{i - 1},x_i],\dgfrac{3r}{4})$;
\item[(1.2)] 
$h_i \nrestriction B(x_{i - 1},\frac{2}{3} \mcdot \frac{3r}{4}) =
\rfs{tr}_{x_i - x_{i - 1}} \rest
B(x_{i - 1},\frac{2}{3} \mcdot \frac{3r}{4})$.
\item[(1.3)] 
$h_i$ is $K_{\srfs{seg}}(\dgfrac{r}{2},\dgfrac{3r}{4})$-bilipschitz.
\vspace{-02.0pt}
\end{list}

Let $h = h_n \scirc \ldots \scirc h_1$.
Then $h$ satisfies requirements (1) and (2) in the lemma.
Also, $h$ is
$K_{\srfs{seg}}(\dgfrac{r}{2},\dgfrac{3r}{4})^n$-bilipschitz.
Since $n = [\dgfrac{2L}{r}] + 1$, we may define
$K_{\srfs{arc}}(\ell,t) =
K_{\srfs{seg}}(\dgfrac{t}{2},\dgfrac{3t}{4})^{[\frac{2 \ell}{t}] + 1}$.
\smallskip\hfill\myqed

If $L$ is a rectifiable arc
let $\fnn{\gamma_L}{[0,1]}{L}$ be a parametrization of $L$
such that for every $t \in [0,1]$,
$\rfs{lngth}(\gamma_L \nrestriction [0,t]) = t\cdot \rfs{lngth}(L)$.

\begin{lemma}\label{l-gamma.8}\label{metr-bldr-l5.27}
Let $X$ be an open subset of a normed space $E$.
For $n \in \bbN$ let $L_n \subseteq X$ be a rectifiable arc with
$\rfs{lngth}(L_n) \leq M$ and $\delta(L_n) \geq d > 0$.
Let $\gamma_n = \gamma_{L_n}$ and $x_n = \gamma_n(0)$.
Suppose that $\setm{x_n}{n \in \bbN}$ is spaced and $\onetoonet$
and that there is $x^* \in X$ such that for every $n \in \bbN$,
$\gamma_n(1) = x^*$.
\underline{Then} there are $\hat{x} \in X$, $r > 0$,
an infinite $\eta \subseteq \bbN$
and $t \in (0,1]$ such that:
\begin{itemize}
\addtolength{\parskip}{-11pt}
\addtolength{\itemsep}{06pt}
\item[\num{1}] 
$B(\hat{x},r) \subseteq X$, $B(\hat{x},r)$ is a BPD set,
and for every $n \in \eta$, $x_n \not\in \rfs{cl}^E(B(\hat{x},r))$;
\item[\num{2}]
for every $n \in \eta$, $\gamma_n(t) \in B(\hat{x},r)$;
\item[\num{3}]
$\setm{\gamma_n \nrestriction [0,t]}{n \in \eta}$ is spaced.
\vspace{-05.7pt}
\end{itemize}
\end{lemma}

\noindent
{\bf Proof } For $\eta \subseteq \bbN$ and $t \in [0,1]$ let
$A[\eta,t] = \setm{\gamma_n(t)}{n \in \bbN}$.
We regard $A[\eta,t]$ both as a set and as a sequence.
For every infinite $\eta \subseteq \bbN$ \ let
$$
s_{\eta} =
\inf(\setm{s \in [0,1]}{ A[\eta,s]
\mbox{ contains a Cauchy sequence}}).
$$
Since for every $n \in \bbN$, $\gamma_n(1) = x^*$, $s_{\eta}$ is
well defined.
Clearly, if $\eta \subseteq \sigma$, then $s_{\eta} \geq s_{\sigma}$.

As in  \ref{l-gamma.5}, there is an infinite $\sigma \subseteq \bbN$
such that for every infinite $\eta \subseteq \sigma$,
$s_{\eta} = s_{\sigma}$.
Let $s = s_{\sigma}$. We show that if $t \in [0,s)$,
then
\begin{itemize}
\addtolength{\parskip}{-11pt}
\addtolength{\itemsep}{06pt}
\item[$(*)$] there is no infinite set $\eta \subseteq \sigma$ and a
sequence
$\setm{t_i}{i \in \eta}$ such that for every $i \in \eta$,
$t_i \in [0,t]$,
and $\setm{\gamma_i(t_i)}{i \in \eta}$
is a Cauchy sequence.
\vspace{-05.7pt}
\end{itemize}
Suppose otherwise.
We may assume that
$\setm{t_i}{i \in \eta}$ is a convergent sequence.
Let $t^*$ be the limit of this sequence. So $t^* < s$.
Let $I_i$ be the interval whose endpoints are $t_i$ and $t^*$.
Recall that
$\rfs{lngth}(\gamma_i \nrestriction I_i) =
|t^* - t_i| \cdot \rfs{lngth}(\gamma_i) \leq |t^* - t_i| \cdot M$.
So $\lim_{\,i \,\in\, \eta} d(\gamma_i(t_i),\gamma_i(t^*)) = 0$.
Hence $\setm{\gamma_i(t^*)}{i \in \eta}$ is a Cauchy sequence.
This contradicts the definition of $s$.

Suppose by contradiction that there is an infinite
$\eta \subseteq \sigma$
such that $A[\eta,s]$ is spaced. Let $e > 0$ be such that
$A[\eta,s]$ is $e$-spaced. Then for every
$t \in [s,s + \frac{e}{3M}]$,
$A[\eta,t]$ is spaced. So $s_{\eta} > s_{\sigma}$. This contradicts
the definition of $\sigma$.

It follows that $A[\sigma,s]$ contains a Cauchy sequence.
Hence we may assume that $A[\sigma,s]$ is a Cauchy sequence.
Let $\bar{x} = \lim^{\oversE} A[\sigma,s]$.
Since $\setm{x_n}{n \in \sigma}$ is 1\,-\,1,
we may assume that for every $n$,
$\barx \neq x_n$.
Since $\delta(L_n) \geq d > 0$,
$d^{\oversE}(\bar{x},E - X) \geq d > 0$.
Since $\setm{x_n}{n \in \bbN}$ is spaced,
there is $0 < r < d$ such that
$\setm{x_n}{n \in \sigma} \cap B^{\oversE}(\barx,r) = \emptyset$.
Let $t = s - \frac{r}{2M}$.
There is $i_0$ such that for every $i_0 \leq i \in \sigma$,
$d(\gamma_i(s),\bar{x}) < \dgfrac{r}{4}$.
We may assume that $i_0 = 0$.
So for every $i \in \sigma$, \ %
$d(\gamma_i(t),\bar{x}) \leq
d(\gamma_i(t),\gamma_i(s)) + d(\gamma_i(s),\bar{x}) <
\rfs{lngth}(\gamma_i \nrestriction [t,s]) + \dgfrac{r}{4} \leq
(s - t) \cdot M + \dgfrac{r}{4} \leq \dgfrac{3r}{4}$.
Let $\hatx \in E \cap B^{\oversE}(\barx,\dgfrac{r}{8})$.
So for every $i \in \sigma$,
$d(\gamma_i(t),\hat{x}) < \dgfrac{7r}{8}$.

By $(*)$ and Proposition \ref{p-gamma.6},
there is an infinite $\eta \subseteq \sigma$ such that
$\setm{\gamma_i \nrestriction [0,t]}{i \in \eta}$ is spaced.
Also, since $\delta(\hat{x}) \geq d - \dgfrac{r}{8}$,
$\delta(B(\hat{x},\dgfrac{7r}{8})) \geq d - r > 0$.
So $B(\hat{x},\dgfrac{7r}{8})$ is a BPD set.
Hence $\hat{x}$, $r$, $\eta$ and $t$ are as required in the lemma.
\rule{0pt}{1pt}\hfill\myqed

\begin{prop}\label{p-gamma.9}\label{metr-bldr-p5.28}
Let $\itGamma$ be a countably generated modulus of continuity,
and let $a > 0$.
Then there is $\setm{\alpha_n}{n \in \bbN} \subseteq \itGamma$ such that
\begin{itemize}
\addtolength{\parskip}{-11pt}
\addtolength{\itemsep}{06pt}
\item[\num{1}]
For every $\alpha \in \itGamma$ there is $n \in \bbN$ such that
$\alpha \preceq \alpha_n$. That is, $\setm{\alpha_n}{n \in \bbN}$
generates $\itGamma$,
\item[\num{2}]
for every $m < n$,
$\alpha_m \nrestriction [0,a] \leq \alpha_n \nrestriction [0,a]$.
\vspace{-05.7pt}
\end{itemize}
\end{prop}

\noindent
{\bf Proof } Let $\setm{\beta_n}{n \in \bbN}$ be a generating set
for $\itGamma$ such that for every $m < n$, $\beta_m \preceq \beta_n$.
We define by induction $K_n > 0$ and $\alpha_n \in \itGamma$.
We assume by induction that $\alpha_n = K_n \beta_n$.
Let $K_0 = 1$ and $\alpha_0 = \beta_0$.
Suppose that $K_n$ and $\alpha_n$ have been defined.
Let $i \leq n$.
Since $\beta_i \preceq \beta_{n + 1}$ and $\alpha_i = K_i \beta_i$,
it follows that
$M_i \eqdf
\sup_{x \in [0,a]} \frac{\alpha_i(x)}{\beta_{n + 1}(x)} < \infty$.
Let $K_{n + 1} = \max(M_0,\ldots,M_n) + 1$
and $\alpha_{n + 1} = K_{n + 1} \beta_{n + 1}$.
Obviously, $\setm{\alpha_n}{n \in \bbN} \subseteq \gG$
and $\setm{\alpha_n}{n \in \bbN}$ is as required.
\medskip\hfill\myqed

\noindent
{\bf Proof of Theorem \ref{metr-bldr-t5.23} }
(a) Let $\itGamma$, $X$, $Y$ and $\tau$ be as in Part (a).
We have that\break
$ \rfs{LIP}_{00}(X) \subseteq H_{\itGamma}^{\srfs{WBPD}}(X)$
and
$H_{\itGamma}^{\srfs{WBPD}}(Y) \subseteq \rfs{BPD.P}(Y)$,
hence
$(\rfs{LIP}_{00}(X))^{\tau}  \subseteq \rfs{BPD.P}(Y)$.
So by Lemma \ref{l-gamma.5}, $\tau \in \rfs{BPD.P}(X,Y)$.

Using the notation of Definition \ref{d2.7}(a),
$\rfs{LIP}_{00}(X) = \rfs{LIP}(X;\calU)$,
where $\calU$ is the set of all open BPD subsets of $X$.
Clearly,
$H_{\itGamma}^{\srfs{WBPD}}(Y) \subseteq H^{\srfs{LC}}_{\itGamma}(Y)$
so we have that
$(\rfs{LIP}(X;\calU))^{\tau} \subseteq H^{\srfs{LC}}_{\itGamma}(Y)$.
Hence by Theorem \ref{metr-bldr-t3.27},
$\tau$ is locally $\itGamma$-continuous.

Suppose by contradiction that there is an open BPD set $U \subseteq X$
such that for no $\alpha \in \itGamma$, $\tau \nrestriction U$ is
locally $\sngltn{\alpha}$-continuous.
Let $\setm{\alpha_n}{n \in \bbN}$ generate $\itGamma$.
We may assume that for every $m < n \in \bbN$,
$\alpha_m \preceq \alpha_n$.
For every $n \in \bbN$ let $\beta_n = \alpha_n \scirc \alpha_n$
and $x_n \in U$ be such that for every $V \in \rfs{Nbr}^X(x_n)$,
$\tau \nrestriction V$ is not $\beta_n$-continuous.
Let $\vecx = \setm{x_n}{n \in \bbN}$.

Suppose by contradiction that $\vecx$ has a Cauchy subsequence $\vecy$.
Let $\bary = \lim^{\oversE} \vecy$.
Since $U$ is a BPD set and $\rfs{Rng}(\vecy) \subseteq U$,
$\bary \in \overfs{int}(X)$.
Let $u \in X$ and $r > 0$ be such that $B^E(u,2r) \subseteq X$
and $\bary \in B^{\oversE}(u,r)$.
Since $\tau$ is locally $\itGamma$-continuous,
there are $V \in \rfs{Nbr}^X(u)$ and $\beta \in \itGamma$
such that $\tau \nrestriction V$ is $\beta$-continuous.
There is $h \in \rfs{LIP}(X)\sprtl{B(u,r)}$ such that
$h^{\srfs{cl}}_{\oversE}(\bary) \in \overfs{int}(V)$.
Since $h \in \rfs{LIP}_{00}(X)$,
$h^{\tau} \in H_{\itGamma}^{\srfs{WBPD}}(Y)$.

Recall that $\tau \in \rfs{BPD.P}(X,Y)$.
Since $B(u,r)$ is a BPD set in $X$,
$W \eqdf \tau(B(u,r))$
is a BPD set in $Y$.
So there is $\alpha \in \itGamma$ such that
$(h^{\tau})\inverse \nrestriction W$
is locally $\sngltn{\alpha}$-continuous.
Since $\lim \vecy = \bary$
and $h^{\srfs{cl}}_{\oversE}(\bary) \in \overfs{int}(V)$,
we may assume that $h(\vecy) \subseteq V$.
Let $K$ be such that $h$ is $K$-bilipschitz, and define
$\gamma(t) = Kt$. So $\gamma \in \itGamma$.
We show that for every $n \in \bbN$,
$\tau$ is $\alpha \scirc \beta \scirc \gamma$-bicontinuous at $y_n$.
Note that $\tau = (h^{\tau})\inverse \scirc \tau \scirc h$.
We have
\begin{itemize}
\addtolength{\parskip}{-11pt}
\addtolength{\itemsep}{06pt}
\item[(i)] 
$h$ is $\gamma$-bicontinuous at $y_n$.
\vspace{-05.7pt}
\end{itemize}
Since $h(y_n) \in B(u,r)$, we have
\begin{itemize}
\addtolength{\parskip}{-11pt}
\addtolength{\itemsep}{06pt}
\item[(ii)] 
$\tau$~is $\beta$-bicontinuous at $h(y_n)$.
\vspace{-05.7pt}
\end{itemize}
Also, $\tau(h(y_n)) \in \tau(B(u,r)) = W$. So
\begin{itemize}
\addtolength{\parskip}{-11pt}
\addtolength{\itemsep}{06pt}
\item[(iii)] 
$(h^{\tau})\inverse$ is $\alpha$-bicontinuous at $\tau(h(y_n))$.
\vspace{-05.7pt}
\end{itemize}
It follows from (i)\,-\,(iii) that $\tau$ is
$\alpha \scirc \beta \scirc \gamma$-bicontinuous at~$y_n$.
Clearly, $\alpha \scirc \beta \scirc \gamma \in \itGamma$,
so there is $n$ such that
$\alpha \scirc \beta \scirc \gamma \preceq \beta_n$.
Hence $\tau$ is $\beta_n$-bicontinuous at $y_n$. This contradicts
the choice of $y_n$.
So $\vecx$ does not have Cauchy subsequences.

We may thus assume that $\vecx$ is spaced.
Let $x^* \in U$. Since $X$ is BPD.AC, there are $M,d > 0$ and
rectifiable arcs $\setm{L_n}{n \in \bbN}$ such that
for every $n \in \bbN$, $L_n$ connects $x_n$ with $x^*$,
$\delta(L_n) \geq d$ and $\rfs{lngth}(L_n) \leq M$.
Applying Lemma~\ref{l-gamma.8} to $x^*$ and $\setm{L_n}{n \in \bbN}$
we obtain $\hat{x} \in X$, $r > 0$,
an infinite $\eta \subseteq \bbN$ and $t \in (0,1]$
as assured by that lemma.
So for the parametrization $\gamma_n$ of $L_n$ defined in
Lemma~\ref{l-gamma.8} the following holds:
\begin{list}{}
{\setlength{\leftmargin}{39pt}
\setlength{\labelsep}{08pt}
\setlength{\labelwidth}{20pt}
\setlength{\itemindent}{-00pt}
\addtolength{\topsep}{-04pt}
\addtolength{\parskip}{-02pt}
\addtolength{\itemsep}{-05pt}
}
\item[(1.1)]
$B(\hat{x},r) \subseteq X$, $B(\hat{x},r)$ is a BPD set,
and for every $n \in \eta$, $x_n \not\in \rfs{cl}^E(B(\hat{x},r))$;
\item[(1.2)]
for every $n \in \eta$, $\gamma_n(t) \in B(\hat{x},r)$;
\item[(1.3)]
$\setm{\gamma_n \nrestriction [0,t]}{n \in \eta}$ is spaced.
\vspace{-02.0pt}
\end{list}
We may assume that $\eta = \bbN$.
{\thickmuskip=3mu \medmuskip=2mu \thinmuskip=1mu 
For every $n \in \bbN$ let $t_n$ be the least $t'$ such that
$\gamma_n(t') \in \rfs{cl}^E(B(\hat{x},r))$.
Let $\gammaprime_n = \gamma_n \nrestriction [0,t_n]$
and $y_n = \gamma_n(t_n)$.
So $d(y_n,\hat{x}) = r$
and $\rfs{Rng}(\gammaprime_n) \cap B(\hat{x},r) = \emptyset$.}

Since $\tau$ is locally $\itGamma$-continuous,
there is $\alpha^* \in \itGamma$ and $r_1 < r$ such
that $\tau \nrestriction B(\hat{x},r_1)$ is
$\alpha^*$-continuous.
Let
$z_n =
\hat{x} + \frac{r_1}{2} \cdot
\dgfrac{(y_n - \hat{x})}{{\norm{y_n - \hat{x}}}}$
and
$L^*_n = \rfs{Rng}(\gammaprime_n) \cup [y_n,z_n]$.
So $L^*_n$ is a rectifiable arc.
Clearly, there are $M^*,d^*,D^* > 0$ such that for every
distinct $m,n \in \bbN$,
\begin{list}{}
{\setlength{\leftmargin}{39pt}
\setlength{\labelsep}{08pt}
\setlength{\labelwidth}{20pt}
\setlength{\itemindent}{-00pt}
\addtolength{\topsep}{-04pt}
\addtolength{\parskip}{-02pt}
\addtolength{\itemsep}{-05pt}
}
\item[(2.1)] $\rfs{lngth}(L^*_m) \leq M^*$;
\item[(2.2)] $\delta(L^*_m) \geq d^*$;
\item[(2.3)] $d(L^*_m,L^*_n) \geq D^*$.
\vspace{-02.0pt}
\end{list}
Let $r^* > 0$ be such that
$r^* < \dgfrac{d^*}{2}, \dgfrac{D^*}{3}, \dgfrac{r_1}{2}$.
For every $n \in \bbN$ we apply Lemma \ref{l-gamma.7} with
$L = M^*$, $r = r^*$, $\gamma = L^*_n$, $x = x_n$ and
$y = z_n$.
We obtain $h_n \in H(X)$ such that:
\begin{list}{}
{\setlength{\leftmargin}{39pt}
\setlength{\labelsep}{08pt}
\setlength{\labelwidth}{20pt}
\setlength{\itemindent}{-00pt}
\addtolength{\topsep}{-04pt}
\addtolength{\parskip}{-02pt}
\addtolength{\itemsep}{-05pt}
}
\item[(3.1)] $\rfs{supp}(h_n) \subseteq B(L^*_n,r^*)$;
\item[(3.2)] $h_n \nrestriction B(x_n,\dgfrac{r^*}{2}) =
\rfs{tr}_{z_n - x_n} \nrestriction B(x_n,\dgfrac{r^*}{2})$;
\item[(3.3)] $h_n$ is $K_{\srfs{arc}}(M^*,r^*)$-bilipschitz.
\vspace{-02.0pt}
\end{list}
Clearly,
$\setm{h_n}{n \in \bbN}$ and $\setm{h_n\inverse}{n \in \bbN}$
satisfy the conditions of
Proposition \ref{p4.17}(b)
with $\alpha(x) = K_{\srfs{arc}}(M^*,r^*) \mcdot x$.
Define $h = \bcirc_{n \in \sboldbbN} h_n$ and $g = h\inverse$. So by
Proposition \ref{p4.17}(b),
$h$ and $g$ are $2K(M^*,r^*)$-Lipschitz.
Also $\delta(\rfs{supp}(h)) \geq d^* - r^* > 0$. So
$h,g \in \rfs{LIP}_{00}(X)$.
Since $\tau \in \rfs{BPD.P}(X,Y)$, $\tau(U)$ is a BPD subset of $Y$.
We shall thus reach a contradiction by proving the following statement.
\begin{itemize}
\addtolength{\parskip}{-11pt}
\addtolength{\itemsep}{06pt}
\item[$(*)$] 
There is no $\alpha \in \itGamma$ such that
$g^{\tau} \nrestriction \tau(U)$
is locally $\sngltn{\alpha}$-continuous.
\vspace{-05.7pt}
\end{itemize}

Let $\alpha \in \itGamma$. Choose $n$ such that
$\alpha,\alpha^* \preceq \alpha_n$
and set $u = \tau(z_n)$.
For $s > 0$ define $U_s = B(u,s)$,
$T_s = \tau\inverse(U_s)$ and $S_s = h\inverse(T_s)$.
There is $s > 0$ such that:
\begin{list}{}
{\setlength{\leftmargin}{39pt}
\setlength{\labelsep}{08pt}
\setlength{\labelwidth}{20pt}
\setlength{\itemindent}{-00pt}
\addtolength{\topsep}{-04pt}
\addtolength{\parskip}{-02pt}
\addtolength{\itemsep}{-05pt}
}
\item[(4.1)]
$\alpha \nrestriction [0,2s] \leq \alpha_n \nrestriction [0,2s]$;
\item[(4.1)]
$T_s \subseteq B(z_n,\dgfrac{r^*}{2})$;
\item[(4.1)]
$\alpha^* \nrestriction [0,\rfs{diam}(T_s)] \leq
\alpha_n \nrestriction [0,\rfs{diam}(T_s)]$.
\vspace{-02.0pt}
\end{list}
Let $s' < s$. We show that $h^{\tau} \nrestriction B(u,s')$ is not
$\alpha$-continuous.
Since $S_{s'}$ is a neighborhood of $x_n$, there are
$x^1,x^2 \in S_{s'}$ such that
\begin{list}{}
{\setlength{\leftmargin}{39pt}
\setlength{\labelsep}{08pt}
\setlength{\labelwidth}{20pt}
\setlength{\itemindent}{-00pt}
\addtolength{\topsep}{-04pt}
\addtolength{\parskip}{-02pt}
\addtolength{\itemsep}{-05pt}
}
\item[(5.1)]
$d(\tau(x^1),\tau(x^2)) > \beta_n(d(x^1,x^2))$.
\vspace{-02.0pt}
\end{list}
For $i = 1,2$ let  $z^i = h(x^i)$ and $u^i = \tau(z^i)$.
So $z^1,z^2 \in T_{s'}$ and so $u^1,u^2 \in U_{s'}$.
By (4.2), the choice of $z_n$ and the choice of $r^*$,
$T_{s'} \subseteq B(z_n,\dgfrac{r^*}{2}) \subseteq B(\hat{x},r_1)$.
So $\tau \nrestriction T_{s'}$ is $\alpha^*$-continuous.
\newline
Hence $\alpha^*(d(z^1,z^2)) \geq d(u^1,u^2)$.
\newline
By (4.3), $\alpha_n(d(z^1,z^2)) \geq \alpha^*(d(z^1,z^2))$.
\newline
So $\alpha_n(d(z^1,z^2)) \geq d(u^1,u^2)$. Hence
\begin{list}{}
{\setlength{\leftmargin}{39pt}
\setlength{\labelsep}{08pt}
\setlength{\labelwidth}{20pt}
\setlength{\itemindent}{-00pt}
\addtolength{\topsep}{-04pt}
\addtolength{\parskip}{-02pt}
\addtolength{\itemsep}{-05pt}
}
\item[(5.2)]
$d(z^1,z^2) \geq (\alpha_n)\inverse(d(u^1,u^2))$.
\vspace{-02.0pt}
\end{list}
Since $T_{s'} \subseteq B(z_n,\dgfrac{r^*}{2})$
and by Property (3.2) of $h_n$, $h\inverse \nrestriction U_{s'}$
is an isometry.
So
\begin{list}{}
{\setlength{\leftmargin}{39pt}
\setlength{\labelsep}{08pt}
\setlength{\labelwidth}{20pt}
\setlength{\itemindent}{-00pt}
\addtolength{\topsep}{-04pt}
\addtolength{\parskip}{-02pt}
\addtolength{\itemsep}{-05pt}
}
\item[(5.3)]
$d(z^1,z^2) = d(x^1,x^2)$.
\vspace{-02.0pt}
\end{list}
By (5.1) and (5.3),
\begin{list}{}
{\setlength{\leftmargin}{39pt}
\setlength{\labelsep}{08pt}
\setlength{\labelwidth}{20pt}
\setlength{\itemindent}{-00pt}
\addtolength{\topsep}{-04pt}
\addtolength{\parskip}{-02pt}
\addtolength{\itemsep}{-05pt}
}
\item[(5.4)]
$d(\tau(x^1),\tau(x^2)) > \beta_n(d(z^1,z^2))$.
\vspace{-02.0pt}
\end{list}
Combining (5.2) and (5.4) we obtain
\begin{list}{}
{\setlength{\leftmargin}{39pt}
\setlength{\labelsep}{08pt}
\setlength{\labelwidth}{20pt}
\setlength{\itemindent}{-00pt}
\addtolength{\topsep}{-04pt}
\addtolength{\parskip}{-02pt}
\addtolength{\itemsep}{-05pt}
}
\item[(5.5)]
$d(\tau(x^1),\tau(x^2)) > \beta_n((\alpha_n)\inverse(d(u^1,u^2)))$.
\vspace{-02.0pt}
\end{list}
But $\beta_n = \alpha_n \scirc \alpha_n$. So
\begin{list}{}
{\setlength{\leftmargin}{39pt}
\setlength{\labelsep}{08pt}
\setlength{\labelwidth}{20pt}
\setlength{\itemindent}{-00pt}
\addtolength{\topsep}{-04pt}
\addtolength{\parskip}{-02pt}
\addtolength{\itemsep}{-05pt}
}
\item[(5.6)]
$d(\tau(x^1),\tau(x^2)) > \alpha_n(d(u^1,u^2))$.
\vspace{-02.0pt}
\end{list}
By Clause (4.1) in the definition of $s$,
and since $u_1,u_2 \in B(u,s)$,
\begin{list}{}
{\setlength{\leftmargin}{39pt}
\setlength{\labelsep}{08pt}
\setlength{\labelwidth}{20pt}
\setlength{\itemindent}{-00pt}
\addtolength{\topsep}{-04pt}
\addtolength{\parskip}{-02pt}
\addtolength{\itemsep}{-05pt}
}
\item[(5.7)]
$d(\tau(x^1),\tau(x^2)) > \alpha(d(u^1,u^2))$.
\vspace{-02.0pt}
\end{list}
But $\tau(x^i) = (h\inverse)^{\tau}(u^i) = g^{\tau}(u^i)$. So
\begin{list}{}
{\setlength{\leftmargin}{39pt}
\setlength{\labelsep}{08pt}
\setlength{\labelwidth}{20pt}
\setlength{\itemindent}{-00pt}
\addtolength{\topsep}{-04pt}
\addtolength{\parskip}{-02pt}
\addtolength{\itemsep}{-05pt}
}
\item[(5.8)]
$d(g^{\tau}(u^1),g^{\tau}(u^2)) > \alpha(d(u^1,u^2))$.
\vspace{-02.0pt}
\end{list}
We have proved $(*)$, and this contradicts the fact that
$g^{\tau} \in H_{\itGamma}^{\srfs{WBPD}}(Y)$. So Part (a) is proved.
\smallskip

(b) Let $\itGamma$, $X$, $Y$ and $\tau$ be as in Part (b).
As in the proof of Part (a), we conclude that
$\tau \in \rfs{BPD.P}(X,Y)$
and $\tau$ is locally $\itGamma$-continuous.

Suppose by contradiction that there is an open BPD set $U \subseteq X$
such that for no $\alpha \in \itGamma$ and $r > 0$,
$\tau \nrestriction U$ is $(r,\alpha)$-continuous.
By Proposition \ref{p-gamma.9}, there is a set
$\setm{\alpha_n}{n \in \bbN}$ which generates $\itGamma$
and such that
$\alpha_m \nrestriction [0,1] \leq \alpha_n \nrestriction [0,1]$
for every $m < n$.
For every $n \in \bbN$ let $\beta_n = \alpha_n \scirc \alpha_n$,
and $x_n, x'_n \in U$
be such that $d(x_n,x'_n) < \dgfrac{1}{n}$
and $d(\tau(x_n),\tau(x'_n)) > \beta_n(d(x_n x'_n))$.
Let
$\vecx = \setm{x_n}{n \in \bbN}$.

Suppose by contradiction that $\setm{x_{n_i}}{i \in \bbN}$ is a Cauchy
subsequence $\vecx$. Denote $y_i = x_{n_i}$ and $y_i' = x_{n_i}'$.
Let $\bary = \lim^{\oversE} \vecy$.
Since $U$ is a BPD set and $\rfs{Rng}(\vecy) \subseteq U$,
$\bary \in \overfs{int}(X)$.
Let $u \in X$ and $r > 0$ be such that $B^E(u,2r) \subseteq X$
and $\bary \in B^{\oversE}(u,r)$.
Since $\tau$ is locally $\itGamma$-continuous,
there are $V \in \rfs{Nbr}^X(u)$ and $\beta \in \itGamma$
such that $\tau \nrestriction V$ is $\beta$-continuous.
There is $h \in \rfs{LIP}(X)\sprtl{B(u,r)}$ such that
$h^{\srfs{cl}}_{\oversE}(\bary) \in \overfs{int}(V)$.
Since $h \in \rfs{LIP}_{00}(X)$,
$h^{\tau} \in H_{\itGamma}^{\srfs{NBPD}}(Y)$.

Recall that $\tau \in \rfs{BPD.P}(X,Y)$.
Since $B(u,r)$ is a BPD set in $X$,
it follows that $W \eqdf \tau(B(u,r))$ is a BPD set in $Y$.
So there are $\alpha \in \itGamma$ and $s > 0$
such that
$h^{\tau} \nrestriction W$ is $(s,\alpha)$-continuous,
and $(h^{\tau})\inverse \nrestriction W$ is $(s,\alpha)$-continuous.
Since $\lim \vecy = \lim \vecy\fprime = \bary$
and $h^{\srfs{cl}}_{\oversE}(\bary) \in \overfs{int}(V)$,
we may assume that $h(\vecy),h(\vecy\fprime) \subseteq V$.

From the fact $h \in \rfs{LIP}(X)$ it follows that
$\lim_{i \rightarrow \infty} d(h(y_i),h(y_i')) = 0$.
Set $u_i = h(y_i)$ and $u_i' = h(y_i')$.
Since $h(\vecy),h(\vecy\fprime) \subseteq V$
and $\tau \nrestriction V$ is $\beta$-continuous,
it follows that
$\lim_{i \rightarrow \infty} d(\tau(u_i),\tau(u_i')) = 0$.
We may thus assume that for every $i \in \bbN$,
$d(\tau(u_i),\tau(u_i')) < s$.

Let $K$ be such that $h$ is $K$-bilipschitz,
define $\gamma(t) = Kt$ and
$\rho = \alpha \scirc \beta \scirc \gamma$.
So $\gamma \in \itGamma$ and hence $\rho \in \itGamma$.
We show that for every $i \in \bbN$
\begin{itemize}
\addtolength{\parskip}{-11pt}
\addtolength{\itemsep}{06pt}
\item[$(\dagger)$] 
$d(\tau(y_i),\tau(y_i')) \leq \rho(d(y_i,y_i'))$.
\vspace{-05.7pt}
\end{itemize}
Note that $\tau(y_i) = (h^{\tau})\inverse \scirc \tau \scirc h(y_i)$,
and the same holds for $y_i'$.
So
\begin{itemize}
\addtolength{\parskip}{-11pt}
\addtolength{\itemsep}{06pt}
\item[(1)]
$d(h(y_i),h(y_i') \leq \gamma(d(y_i,y_i')$.
\vspace{-05.7pt}
\end{itemize}
Now, $h(y_i),h(y_i') \in V$
and $\tau \nrestriction V$ is $\beta$-continuous, so
\begin{itemize}
\addtolength{\parskip}{-11pt}
\addtolength{\itemsep}{06pt}
\item[(2)]
$d(\tau(h(y_i)),\tau(h(y_i'))) \leq \beta(\gamma(d(y_i,y_i')))$. 
\vspace{-05.7pt}
\end{itemize}
Since $d(\tau(u_i),\tau(u_i')) < s$ and $\tau(u_i),\tau(u_i') \in W$,
it follows that
\begin{itemize}
\addtolength{\parskip}{-11pt}
\addtolength{\itemsep}{06pt}
\item[(3)]
$d((h^{\tau})\inverse(\tau(u_i)),
(h^{\tau})\inverse(\tau(u_i'))) \leq \alpha(d(\tau(u_i),\tau(u_i'))$.
\vspace{-05.7pt}
\end{itemize}
Obviously, (1)\,-\,(3) imply $(\dagger)$.

Denote $\hatbeta_i = \beta_{n_i}$.
There is $j$ such that $\rho \preceq \hatbeta_j$.
Let $\ell \in \bbN$ be such that
$\rho \rest[0,\dgfrac{1}{\ell}] \leq
\hatbeta_j \rest[0,\dgfrac{1}{\ell}]$.
Let $i = \max(j,\ell)$.
So $d(y_i,y_i') \leq \dgfrac{1}{n_i} \leq \dgfrac{1}{\ell}$.
From~$(\dagger)$ and the fact
$\hatbeta_j \nrestriction [0,1] \leq \hatbeta_i \nrestriction [0,1]$
we conclude that
$d(\tau(y_i),\tau(y_i')) \leq \rho(d(y_i,y_i')) \leq
\hatbeta_i(d(y_i,y_i'))$.
That is,
$$
d(\tau(x_{n_i}),\tau(x_{n_i}')) \leq
\beta_{n_i}(d(x_{n_i},x_{n_i}')).
$$
This contradicts the way that $x_{n_i}$ and $x_{n_i}'$ were chosen.
So $\vecx$ has no Cauchy subsequences.

We may thus assume that $\vecx$ is spaced.
Let $x^* \in U$. Since $X$ is BPD.AC, there are $M,d > 0$ and
rectifiable arcs $\setm{L_n}{n \in \bbN}$ such that
for every $n \in \bbN$, $L_n$ connects $x_n$ with $x^*$,
$\delta(L_n) \geq d$ and $\rfs{lngth}(L_n) \leq M$.
From Lemma~\ref{l-gamma.8} we obtain
$\hat{x} \in X$, $r > 0$,
an infinite $\eta \subseteq \bbN$
and $t \in (0,1]$
such that for the parametrization $\gamma_n$ of $L_n$ defined in
Lemma~\ref{l-gamma.8} the following holds:
$B(\hat{x},r)$ is a BPD subset of $X$,
for every $n \in \eta$, $x_n \not\in \rfs{cl}^E(B(\hat{x},r))$
and $\gamma_n(t) \in B(\hat{x},r)$
and the set of arcs
$\setm{\gamma_n \nrestriction [0,t]}{n \in \eta}$ is spaced.
We may assume that $\eta = \bbN$.

For every $n \in \bbN$ let $t_n$ be the least $t'$ such that
$\gamma_n(t') \in \rfs{cl}^E(B(\hat{x},r))$.
Let
$\gammaprime_n =
\gamma_n \nrestriction [0,t_n]$ and $y_n = \gamma_n(t_n)$.
So $d(y_n,\hat{x}) = r$
and $\rfs{Rng}(\gammaprime_n) \cap B(\hat{x},r) = \emptyset$.

Since $\tau$ is locally $\itGamma$-continuous,
there is $\alpha^* \in \itGamma$ and $r_1 < r$ such
that $\tau \nrestriction B(\hat{x},r_1)$ is
$\alpha^*$-continuous.
Let
$z_n =
\hat{x} + \frac{r_1}{2} \cdot
\dgfrac{(y_n - \hat{x})}{\norm{y_n - \hat{x}}}$
and
$L^*_n = \rfs{Rng}(\gammaprime_n) \cup [y_n,z_n]$.
So $L^*_n$ is a rectifiable arc.
Clearly, there are $M^*,d^*,D^* > 0$
such that for any distinct $m,n \in \bbN$,
$\rfs{lngth}(L^*_m) \leq M^*$,
$\delta(L^*_m) \geq d^*$
and $d(L^*_m,L^*_n) \geq D^*$.
Let $r^* > 0$ be such that
$r^* < \dgfrac{d^*}{2},\, \dgfrac{D^*}{3},\, \dgfrac{r_1}{2}$.

For every $n \in \bbN$ we apply Lemma \ref{l-gamma.7} with
$L = M^*$, $r = r^*$, $\gamma = L^*_n$, $x = x_n$ and
$y = z_n$.
We obtain $h_n \in H(X)$ such that
$\rfs{supp}(h_n) \subseteq B(L^*_n,r^*)$,
$h_n \nrestriction B(x_n,\dgfrac{r^*}{2}) =
\rfs{tr}_{z_n - x_n} \nrestriction B(x_n,\dgfrac{r^*}{2})$
and $h_n$ is $K_{\srfs{arc}}(M^*,r^*)$-bilipschitz.

The families
$\setm{h_n}{n \in \bbN}$ and $\setm{h_n\inverse}{n \in \bbN}$
satisfy the conditions of Proposition~\ref{p4.17}(b)
with $\alpha(x) = K_{\srfs{arc}}(M^*,r^*) \mcdot x$.
Let $h = \bcirc_{n \in \sboldbbN} h_n$ and $g = h\inverse$.
So by Proposition \ref{p4.17}(b),
$h$ is $2K_{\srfs{arc}}(M^*,r^*)$-bilipschitz.
Also, $\delta(\rfs{supp}(h)) \geq d^* - r^* > 0$,
and hence $h,g \in \rfs{LIP}_{00}(X)$.
Since $\tau \in \rfs{BPD.P}(X,Y)$, $\tau(U)$ is a BPD subset of $Y$.
From the fact
$(\rfs{LIP}_{00}(X))^{\tau} \subseteq H_{\itGamma}^{\srfs{NBPD}}(Y)$
it follows that for some $\alpha \in \itGamma$ and $r > 0$,
$g^{\tau} \nrestriction \tau(U)$ is $(r,\alpha)$-bicontinuous.
We shall thus reach a contradiction by proving the following statement.
\begin{itemize}
\addtolength{\parskip}{-11pt}
\addtolength{\itemsep}{06pt}
\item[$(*)$]
There are no $r > 0$ and $\alpha \in \itGamma$ such that
$g^{\tau} \nrestriction \tau(U)$ is $(r,\alpha)$-continuous.
\vspace{-05.7pt}
\end{itemize}

Let $r > 0$ and $\alpha \in \itGamma$.
For $n \in \bbN$ set $z'_n = h(x'_n)$,
$u_n = \tau(z_n)$ and $u'_n = \tau(z'_n)$.
Choose $m \in \bbN$ and $b \in (0,1)$ such that
$\alpha \nrestriction [0,b],
\alpha^* \nrestriction [0,b] \leq \alpha_m \nrestriction [0,b]$.
So for every $n \geq m$,
\begin{itemize}
\addtolength{\parskip}{-11pt}
\addtolength{\itemsep}{06pt}
\item[(1)] 
$\alpha \nrestriction [0,b] \leq \alpha_n \nrestriction [0,b]$;
\item[(2)] 
$\alpha^* \nrestriction [0,b] \leq \alpha_n \nrestriction [0,b]$.
\vspace{-05.7pt}
\end{itemize}
There is $n \geq m$ such that:
\begin{itemize}
\addtolength{\parskip}{-11pt}
\addtolength{\itemsep}{06pt}
\item[(3)] 
$\dgfrac{1}{n} < b$;
\item[(4)] 
$\alpha^*(\dgfrac{1}{n}) < r$;
\item[(5)] 
$\alpha^*(\dgfrac{1}{n}) < b$;
\item[(6)] 
$\dgfrac{1}{n} < \dgfrac{r^*}{2}$.
\vspace{-05.7pt}
\end{itemize}

By the choice of $z_n$ and $r^*$,
$B(z_n,r^*) \subseteq B(\hat{x},r_1)$.
So $\tau \nrestriction B(z_n,r^*)$ is
$\alpha^*$-continuous.
Since $d(x_n,x'_n) \leq \dgfrac{1}{n} < \dgfrac{r^*}{2}$ and by the
definition of $h_n$ and $h$,
\begin{itemize}
\addtolength{\parskip}{-11pt}
\addtolength{\itemsep}{06pt}
\item[(7)] 
$d(x_n,x'_n) = d(z_n,z'_n)$.
\vspace{-05.7pt}
\end{itemize}
Hence $z'_n \in B(z_n,r^*)$, and so
\begin{itemize}
\addtolength{\parskip}{-11pt}
\addtolength{\itemsep}{06pt}
\item[(8)] 
$d(u_n,u'_n) \leq \alpha^*(d(z_n,z'_n))$.
\vspace{-05.7pt}
\end{itemize}
{\thickmuskip=3.0mu \medmuskip=2mu \thinmuskip=1mu 
By (3) and (7) $d(z_n,z'_n) \leq \dgfrac{1}{n} < b$,
so by (2) and (8),
$d(u_n,u'_n) \leq \alpha_n(d(z_n,z'_n))$.
\hbox{It follows that}
}
\begin{itemize}
\addtolength{\parskip}{-11pt}
\addtolength{\itemsep}{06pt}
\item[(9)] 
$d(z_n,z'_n) \geq \alpha_n\inverse(d(u_n,u'_n))$.
\vspace{-05.7pt}
\end{itemize}
By (7) and (9), $d(x_n,x'_n) \geq \alpha_n\inverse(d(u_n,u'_n))$.
By the definition of $\beta_n$, $x_n$
and $x'_n$,
$d(\tau(x_n),\tau(x'_n)) > \alpha_n \scirc \alpha_n(d(x_n,x'_n))$.
So
\begin{itemize}
\addtolength{\parskip}{-11pt}
\addtolength{\itemsep}{06pt}
\item[(10)] 
$d(\tau(x_n),\tau(x'_n)) > \alpha_n(d(u_n,u'_n))$.
\vspace{-05.7pt}
\end{itemize}
Note that $\tau(x_n) = g^{\tau}(u_n)$
and $\tau(x'_n) = g^{\tau}(u'_n)$. So
\begin{itemize}
\addtolength{\parskip}{-11pt}
\addtolength{\itemsep}{06pt}
\item[(11)] 
$d(g^{\tau}(u_n),g^{\tau}(u'_n)) > \alpha_n(d(u_n,u'_n))$.
\vspace{-05.7pt}
\end{itemize}
Since $d(z_n,z'_n) \leq \dgfrac{1}{n}$, by (8) and (5),
$d(u_n,u'_n) \leq b$.
So by (1),
$\alpha_n(d(u_n,u'_n)) \geq \alpha(d(u_n,u'_n))$.
It now follows from (11) that
\begin{itemize}
\addtolength{\parskip}{-11pt}
\addtolength{\itemsep}{06pt}
\item[(12)] 
$d(g^{\tau}(u_n),g^{\tau}(u'_n)) > \alpha(d(u_n,u'_n))$.
\vspace{-05.7pt}
\end{itemize}
By (8), $d(u_n,u'_n) \leq \alpha^*(\dgfrac{1}{n})$.
So by (4), 
\begin{itemize}
\addtolength{\parskip}{-11pt}
\addtolength{\itemsep}{06pt}
\item[(13)] 
$d(u_n,u'_n) < r$.
\vspace{-05.7pt}
\end{itemize}
Facts (12), (13) mean that $g^{\tau} \nrestriction \tau(U)$ is not
$(r,\alpha)$-continuous.
This was proved for arbitrary $r$ and $\alpha$,
namely, we have proved $(*)$.
We have a contradiction to the fact that
$g^{\tau} \in H_{\itGamma}^{\srfs{NBPD}}(Y)$. So Part (b) is proved.
\hfill\myqed

\begin{question}\label{metr-bldr-q5.30}
\begin{rm}
Does Theorem \ref{metr-bldr-t5.23} remain true when the assumption
that $\itGamma$ is countably generated
is dropped or replaced by the assumption that $\itGamma$ is generated
by a set whose cardinality is $\leq \kappa(X)$?
\smallskip
\end{rm}
\end{question}

Note that the use of the countable generatedness of $\itGamma$
in the proof of \ref{metr-bldr-t5.23} was essential.

\begin{theorem}\label{t4.23}\label{metr-bldr-t5.29}
Let $X,Y \in K_{\srfs{NRM}}^{\calO}$.
Suppose that $X$ is BPD.AC. Let $\tau \in H(X,Y)$ be such that
$(\rfs{UC}_{00}(X))^{\tau} \subseteq \rfs{BPD.UC}(Y)$. Then
$\tau \in \rfs{BPD.UC}(X,Y)$.
\end{theorem}

\noindent
{\bf Proof }
By definition, $\rfs{BPD.UC}(Y) \subseteq \rfs{BPD.P}(Y)$,
hence by Lemma \ref{l-gamma.5}, $\tau \in \rfs{BPD.P}(X,Y)$.

Suppose by contradiction that $\tau \not\in \rfs{BPD.UC}(X,Y)$.
Then there are $d > 0$ and $\vecx,\vecy \subseteq X$ such that 
$\rfs{Rng}(\vecx) \cup \rfs{Rng}(\vecy)$ is a BPD set,
$\lim_{n \raro \infty} d(x_n,y_n) = 0$, and for every
$n \in \bbN$, $d(\tau(x_n),\tau(y_n)) \geq d$.

Suppose by contradiction that $\vecx$ has a Cauchy subsequence.
We may then assume that $\vecx$ is a Cauchy sequence.
Let $\barx = \lim^{\oversE} \vecx$.
Since $\rfs{Rng}(\vecx)$ is a BPD set,
$\barx \in \overfs{int}(X)$.
Let $u \in X$ and $r > 0$ be such that $B^E(u,2r) \subseteq X$
and $\barx \in B^{\oversE}(u,r)$.

$\rfs{BPD.UC}(X) \subseteq \rfs{LUC}(X)$
and $\rfs{UC}_{00}(X) = \rfs{UC}(X,\calU)$, where $\calU$ is the
set of all open BPD subsets of $X$.
So by Theorem \ref{metr-bldr-t4.8}(b), $\tau \in \rfs{LUC}(X,Y)$.
So there is $V \in \rfs{Nbr}^X(u)$
such that $\tau \nrestriction V$ is uniformly continuous.
There is $h \in \rfs{LIP}(X)\sprtl{B(u,r)}$ such that
$h^{\srfs{cl}}_{\oversE}(\bary) \in \overfs{int}(V)$.
Since $h \in \rfs{UC}_{00}(X)$,
$h^{\tau} \in \rfs{BPD.UC}(X)$.

Recall that $\tau \in \rfs{BPD.P}(X,Y)$.
Since $B(u,r)$ is a BPD set in $X$,
$W \eqdf \tau(B(u,r))$
is a BPD set in $Y$.
So $h^{\tau} \nrestriction W$ is bi-UC.
Since
$\lim \vecx = \lim \vecy = \barx$
and $h^{\srfs{cl}}_{\oversE}(\barx) \in \overfs{int}(V)$,
we may assume that $h(\vecx),h(\vecy) \subseteq V$.
Since $h$ is uniformly continuous and
$\tau \nrestriction V$ is uniformly continuous,
\begin{itemize}
\addtolength{\parskip}{-11pt}
\addtolength{\itemsep}{06pt}
\item[(1)] 
$\limti{i} d(\tau(h(x_i)),\tau(h(y_i))) = 0$.
\vspace{-05.7pt}
\end{itemize}
Note that
$(h^{\tau})\inverse(\tau(h(x_i))) = \tau(x_i)$,
and the same holds for $y_i$.
\hbox{So for every $i$,}
\begin{itemize}
\addtolength{\parskip}{-11pt}
\addtolength{\itemsep}{06pt}
\item[(2)] 
$d((h^{\tau})\inverse(\tau(h(x_i))),(h^{\tau})\inverse(\tau(h(y_i))))
\geq d$.
\vspace{-05.7pt}
\end{itemize}
(1) and (2) contradict the fact that $h^{\tau} \nrestriction W$
is bi-UC. So $\vecx$ has no Cauchy subsequences.

We may thus assume that there is $s > 0$ such that $\vecx$
is $s$-spaced.
Let $r = \dgfrac{\min(s,\delta(\vecx))}{3}$.
We may assume that for every $n \in \bbN$, $d(y_n,x_n) < \dgfrac{r}{3}$.
Let $r_n = 2d(y_n,x_n)$. Hence $B^E(x_n,r_n) \subseteq X$,
and $\limti{n}\rfs{diam}(B^E(x_n,r_n)) = 0$.
Also, for every distinct $m,n \in \bbN$,
$d(B^E(x_m,r_m),B^E(x_n,r_n)) \geq \dgfrac{s}{3}$.

For every $n \in \bbN$, let $z_n \in [x_n,y_n]$ be such that
$d(\tau(z_n),\tau(x_n)) \leq \dgfrac{d}{(n + 2)}$,
and $h_n \in \rfs{UC}(X)$ be such that
$\rfs{supp}(h_n) \subseteq B(x_n,r_n)$,
$h_n(x_n) = x_n$
and $h_n(z_n) = y_n$.
By Proposition \ref{p4.9},
$h \eqdf \bcirc_{n \in \sboldbbN} h_n \in \rfs{UC}(X)$.
Also $\delta(\rfs{supp}(h)) \geq \dgfrac{r}{3}$.
So $h \in \rfs{UC}_{00}(X)$.
Hence $h^{\tau} \in \rfs{BPD.UC}(Y)$.
$\vecx \cup \vecy \cup \vecz$
is a BPD set.
So since $\tau \in \rfs{BPD.P}(X,Y)$,
it follows that
$\tau(\vecx) \cup \tau(\vecy) \cup \tau(\vecz)$
is a BPD set.
However,
$h^{\tau} \nrestriction (\tau(\vecx) \cup \tau(\vecy) \cup \tau(\vecz))$
is not UC. This is so, because
$\limti{n} d(\tau(x_n),\tau(z_n)) = 0$,
whereas for every $n \in \bbN$,
$d(h^{\tau}(\tau(x_n)),h^{\tau}(\tau(z_n))) = 
d(\tau(x_n),\tau(y_n)) \geq d$.
A contradiction.\bigskip\hfill\myqed

\begin{theorem}\label{metr-bldr-t5.30}
Let $\itGamma,\itDelta$ be moduli of continuity.
Suppose that $\itGamma$ is
countably generated or $\itGamma = \rfs{MC}$,
and that the same holds for $\itDelta$.
Let
$X,Y \in K_{\srfs{NFCB}}^{\calO}$,
and assume that $X$ and $Y$ are BPD.AC.
Suppose that
$\iso{\varphi}
{H_{\itGamma}^{\srfs{NBPD}}(X)}{H_{\itDelta}^{\srfs{NBPD}}(Y)}$.
\underline{Then} $\itGamma = \itDelta$ and there is
$\tau \in (H_{\itGamma}^{\srfs{NBPD}})^{\pm}(X,Y)$
such that $\tau$ induces $\varphi$.
\end{theorem}

\noindent
{\bf Proof }
Let $\calU$ denote the set of all open BPD subsets of $X$.
Note that
\begin{itemize}
\addtolength{\parskip}{-11pt}
\addtolength{\itemsep}{06pt}
\item[(1)] 
$\rfs{LIP}_{00}(X) \leq H_{\itGamma}^{\srfs{NBPD}}(X) \leq \rfs{IXT}(X)$
and
$\rfs{LIP}_{00}(X) = \rfs{LIP}(X,\calU)$.
\vspace{-05.7pt}
\end{itemize}
Hence
by Corollary \ref{metr-bldr-c2.26}, there is $\tau \in H(X,Y)$ such that
$\tau$ induces $\varphi$.
Suppose that $\itDelta$ is countably generated.
Clearly,
\begin{itemize}
\addtolength{\parskip}{-11pt}
\addtolength{\itemsep}{06pt}
\item[(2)] 
$H_{\itDelta}^{\srfs{NBPD}}(Y) \subseteq H_{\itDelta}^{\srfs{LC}}(Y)$.
\vspace{-05.7pt}
\end{itemize}
By (1) and (2),
$(\rfs{LIP}(X,\calU))^{\tau} \subseteq H_{\itDelta}^{\srfs{LC}}(Y)$.
By Theorem \ref{metr-bldr-t3.27},
$\tau$ is locally $\itDelta$-bicontinuous.
Suppose by contadiction that $\alpha \in \itDelta - \itGamma$.
Let $B$ be an open ball in $E$ such that $B$ is a BPD subset of $X$
and such that for some $\beta \in \itDelta$, $\tau \nrestriction B$ is
$\beta$-bicontinuous.
There is $g \in H(X)\sprt{B}$ such that $g$ is $\alpha$-bicontinuous,
and for every $\gamma \in \itGamma$, $g$ is not
$\gamma$-bicontinuous.
So $g \not\in H_{\itGamma}^{\srfs{NBPD}}(X)$,
but $g^{\tau} \in H_{\itDelta}^{\srfs{NBPD}}(Y)$, a contradiction.
So $\itDelta \subseteq \itGamma$.
An identical argument shows that $\itGamma \subseteq \itDelta$.
Hence $\itGamma = \itDelta$.
Applying Theorem \ref{metr-bldr-t5.23} to $\tau$ and $\tau\inverse$,
we conclude that
$\tau \in (H_{\itGamma}^{\srfs{NBPD}})^{\pm}(X,Y)$.

Suppose next that $\itGamma = \itDelta = \rfs{MC}$.
Since $\rfs{UC}_{00}(X) \leq H_{\srfs{MC}}^{\srfs{NBPD}}(X)$,
we have
$(\rfs{UC}_{00}(X))^{\tau} \subseteq H_{\srfs{MC}}^{\srfs{NBPD}}(X)$,
and the same holds for $Y$.
Hence Theorem \ref{metr-bldr-t5.29}
may be applied to $\tau$ and $\tau\inverse$.
We conclude that $\tau \in \rfs{BPD.UC}^{\pm}(X,Y)$.
That is, $\tau \in (H_{\srfs{MC}}^{\srfs{NBPD}})^{\pm}(X,Y)$.
\smallskip\hfill\myqed

We now turn to the group $H_{\srfs{MC}}^{\srfs{WBPD}}(X)$.
We shall reach the same final result as for the groups of type
$H_{\srfs{MC}}^{\srfs{NBPD}}(X)$. But here we need the extra
assumption that $X$ is fillable. This notion is defined below.

\begin{defn}\label{metr-bldr-d5.32}
\begin{rm}
Let $X$ be a topological space and $G \leq H(X)$.
A sequence $\vecx \subseteq X$ is called a {\it $G$-filling} of $X$,
if the following holds.
For every sequence $\setm{U_i}{i \in \bbN}$ such that for every $i$,
$U_i \in \rfs{Nbr}(x_i)$,
there is sequence $\setm{g_i}{i \in \bbN} \subseteq G$ such that
$\bigcup_{i \in \sboldbbN} g_i(U_i) = X$.
We say that $X$ is {\it $G$-fillable} if $X$ has a $G$-filling.
   \index{filling. $G$-filling}
   \index{fillable. $G$-fillable}
\end{rm}
\end{defn}

The trivial verification of the following observation is left
to the reader.

\begin{prop}\label{metr-bldr-p5.33}
Let $E$ be a normed space.

\num{a} If $E$ is separable and $X \subseteq E$ is open,
then $X$ is $\rfs{LIP}_{00}(X)$-fillable.

\num{b} If $r > 0$, then
$B^E(0,r)$ is $\rfs{LIP}_{00}(X)$-fillable.
\end{prop}

The following observation gives some answer for the groups
of type $H_{\srfs{MC}}^{\srfs{WBPD}}(X)$.

\begin{prop}\label{metr-bldr-p5.34}
Suppose that $X$ is BPD.AC,
$\rfs{UC}_{00}(X) \leq G \leq H_{\srfs{MC}}^{\srfs{WBPD}}(X)$
and $X$
is $G$-fillable.
Let $\tau \in H(X,Y)$ be such that
$G^{\tau} \subseteq H_{\srfs{MC}}^{\srfs{WBPD}}(Y)$.
Then $\tau \in H_{\srfs{MC}}^{\srfs{WBPD}}(X,Y)$.
\end{prop}

\noindent
{\bf Proof }
Let $\calU$ be the set of all open BPD subsets of $X$.
Then $\rfs{UC}_{00}(X) = \rfs{UC}(X,\calU)$.
Note that $H_{\srfs{MC}}^{\srfs{WBPD}}(Y) \subseteq \rfs{LUC}(Y)$.
So $(\rfs{UC}(X,\calU))^{\tau} \subseteq \rfs{LUC}(Y)$.
By Theorem \ref{metr-bldr-t4.8}(b),
$\tau \in \rfs{LUC}^{\pm}(X,Y)$.
Similarly,
$(\rfs{LIP}_{00}(X))^{\tau} \subseteq 
(\rfs{UC}_{00}(X))^{\tau} \subseteq 
H_{\srfs{MC}}^{\srfs{WBPD}}(Y) \subseteq
\rfs{BPD.P}(Y)$.
So by Lemma \ref{metr-bldr-l5.24},
$\tau \in \rfs{BPD.P}(X,Y)$.

Let $\vecx$ be a $G$-filling for $X$.
For every $i \in \bbN$ let $U_i \in \rfs{Nbr}(x_i)$ and $\alpha_i$
be such that $\tau \nrestriction U_i$ is $\alpha_i$-bicontinuous.
Let $\setm{g_i}{i \in \bbN} \subseteq G$  be such that
$\bigcup \setm{g_i(U_i)}{i \in \bbN} = X$.

Let $A \subseteq X$ be a BPD set.
We show that $\tau \nrestriction A$ is weakly $\rfs{MC}$-bicontinuous.
Since $\tau \in \rfs{BPD.P}(X,Y)$, $\tau(A)$ is a BPD set.
For every $i \in \bbN$ let $\beta_i$ be such that
$g_i \nrestriction A$ is locally $\sngltn{\beta_i}$-bicontinuous
and $\gamma_i$ be such that
$g_i^{\tau} \nrestriction \tau(A)$
is locally $\sngltn{\gamma_i}$-bicontinuous.
Next note that
$$\tau \nrestriction g_i(U_i) =
(g_i^{\tau} \nrestriction \tau(U_i)) \scirc (\tau \nrestriction U_i)
\scirc ((g_i)\inverse \nrestriction g_i(U_i)).
$$
Hence $\tau \nrestriction (g_i(U_i) \cap A)$ is locally
$\sngltn{\gamma_i \scirc \alpha_i \scirc \beta_i}$-bicontinuous.

There is $\rho \in \rfs{MC}$ such that for every $i \in \bbN$,
$\gamma_i \scirc \alpha_i \scirc \beta_i  \preceq \rho$.
\hbox{Hence for every $i \in \bbN$,}
$\tau \nrestriction (g_i(U_i) \cap A)$ is locally
$\sngltn{\rho}$-bicontinuous,
and from the fact $\bigcup_{i \in \sboldbbN} (g_i(U_i) \cap A) = A$
we conclude that
$\tau \nrestriction A$ is $\sngltn{\rho}$-bicon\-tinuous.
So $\tau \in H_{\srfs{MC}}^{\srfs{WBPD}}(X,Y)$.
\hfill\myqed

\begin{theorem}\label{metr-bldr-t5.35}
Let $\itGamma,\itDelta$ be moduli of continuity.
Suppose that $\itGamma$ is
countably generated or $\itGamma = \rfs{MC}$,
and that the same holds for $\itDelta$.
Let $X,Y \in K_{\srfs{NFCB}}^{\calO}$. Assume that
\begin{itemize}
\addtolength{\parskip}{-11pt}
\addtolength{\itemsep}{06pt}
\item[\num{1}] $X$ and $Y$ are BPD.AC;
\item[\num{2}]
If $\itGamma = \rfs{MC}$,
then $X$ is $H_{\srfs{MC}}^{\srfs{WBPD}}(X)$-fillable,
and the same holds for $\itDelta$ and $Y$.
\vspace{-05.7pt}
\end{itemize}
Suppose that
$\iso{\varphi}
{H_{\itGamma}^{\srfs{WBPD}}(X)}{H_{\itDelta}^{\srfs{WBPD}}(Y)}$.
\underline{Then} $\itGamma = \itDelta$ and there is
$\tau \in H_{\itGamma}^{\srfs{WBPD}}(X,Y)$ such that $\tau$ induces
$\varphi$.
\end{theorem}

\noindent
{\bf Proof } The proof is very similar to the proof of
Theorem \ref{metr-bldr-t5.30}.
\smallskip\hfill\myqed

In some cases we reach a final reconstruction result of the following
strong form.
\begin{itemize}
\addtolength{\parskip}{-11pt}
\addtolength{\itemsep}{06pt}
\item[(1)] 
If $\iso{\varphi}{\calP(X)}{\calQ(Y)}$,
then either $\calP(X) = \calQ(X)$
and there is $\tau \in \calQ^{\pm}(X,Y)$
such that $\tau$ induces $\varphi$,
or
$\calP(Y) = \calQ(Y)$
and there is $\tau \in \calP^{\pm}(X,Y)$
such that $\tau$ induces $\varphi$.
\vspace{-05.7pt}
\end{itemize}
In other cases we are able to reach only the following weaker
conclusion.
\begin{itemize}
\addtolength{\parskip}{-11pt}
\addtolength{\itemsep}{06pt}
\item[(2)] 
If $\iso{\varphi}{\calP(X)}{\calP(Y)}$, then there is
$\tau \in \calP^{\pm}(X,Y)$ such that $\tau$ induces $\varphi$.
\vspace{-05.7pt}
\end{itemize}
Roughly speaking, in order to prove results of the first form,
we need to prove the following intermediate claim.
\begin{itemize}
\addtolength{\parskip}{-11pt}
\addtolength{\itemsep}{06pt}
\item[(3)] 
If $\tau \in H(X,Y)$ and $(\calP(X))^{\tau} \subseteq \calP(Y)$,
then $\tau \in \calP^{\pm}(X,Y)$,
\vspace{-05.7pt}
\end{itemize}
and in order to prove a result of the second form,
the following intermediate claim suffices.
\begin{itemize}
\addtolength{\parskip}{-11pt}
\addtolength{\itemsep}{06pt}
\item[(4)] 
If $\tau \in H(X,Y)$ and $(\calP(X))^{\tau} \subseteq \calP(Y)$,
then $\tau \in \calP(X,Y)$.
\vspace{-05.7pt}
\end{itemize}
For example, Theorem \ref{metr-bldr-t4.8} which deals with the group
$\rfs{LUC}(X)$ has the stronger form~(3),
and Theorem \ref{metr-bldr-t5.6}
which deals with the group $\rfs{UC}(X)$
has the weaker form (4).

The strong intermediate claim is not always true.
Example \ref{metr-bldr-e5.8} shows that Statement~(3)
is false for $\rfs{UC}(X)$,
and Statement~(3) is also false for $\rfs{BPD.UC}(X)$,
as is shown in Example \ref{metr-bldr-e5.37}(a).
However, if $X$ is an open subset of a Banach space,
and $X$ fulfills some additional requirements, then the implication
$$
(\rfs{BPD.UC}(X))^{\tau} \subseteq \rfs{BPD.UC}(Y) \Rightarrow
\tau \in \rfs{BPD.UC}^{\pm}(X,Y)
$$
is true. This will be proved in Theorem \ref{metr-bldr-t5.40}(a).
Later, in Theorem \ref{metr-bldr-t7.7} we shall prove an
analogous statement for $\rfs{UC}(X)$. Namely,
if $X$ satisfies certain additional requirements,
then
$(\rfs{UC}(X))^{\tau} \subseteq \rfs{UC}(Y)$ $\Rightarrow$
$\tau \in \rfs{UC}^{\pm}(X,Y)$.
\smallskip

We need yet another notion of weak uniform arcwise connectedness.
This will be the additional assumption
in Theorem~\ref{metr-bldr-t5.40}(a).

\begin{defn}\label{d4.24}\label{metr-bldr-d5.36}
\begin{rm}
Let $E$ be a metric space, $X \subseteq E$ and $x \in \rfs{bd}^E(X)$.
We say that $X$ is {\it locally arcwise connected at $x$},
if for every $\varepsilon > 0$ there is $\delta > 0$ such that for every
$y,z \in X$: if $d(x,y),d(x,z) < \delta$, then there is an arc
$L \subseteq X$ connecting $y$ and $z$
such that $\rfs{diam}(L) < \varepsilon$.
We then call $x$ a {\it simple boundary point} of~$X$.
We say that $X$ is
{\it locally arcwise connected at its boundary with respect
to $E$ $($BR.LC.AC with respect to $E$$)$},
if every boundary point of $X$ is simple.
\hfill\proofend

   \index{locally arcwise connected at a boundary point}
   \index{simple boundary point}
   \index{locally arcwise connected at the boundary. Abbreviated by
          BR.LC.AC}
   \index{brlcac@@BR.LC.AC. Abbreviation of
          locally arcwise connected at the boundary}
\end{rm}
\end{defn}

An equivalent formulation of simplicity is as follows.
For every $\varepsilon > 0$ there is $\delta > 0$ such that
for every $y,z \in X \cap B(x,\delta)$ there is an arc $L$
connecting $y$ and $z$ such that
$L \subseteq X \cap B(x,\varepsilon)$.
Note that being locally arcwise connected at $x \in \rfs{bd}(X)$
implies but is not equivalent to the fact that $X \cup \sngltn{x}$ is
locally arcwise connected at $x$.

The following example shows that the completeness requirement in
Lemma \ref{metr-bldr-l5.38} cannot be dropped.

\begin{example}\label{bddly-lip-bldr-e5.31}\label{metr-bldr-e5.37}
Let $E$ be an incomplete normed space,
$K \subseteq [0,\dghalf)$ be a closed nowhere dense perfect set
containing $0$,
$X' = B^E(0,2) - \overB^E(0,1)$, $u \in S^E(0,1)$,
\newline
$C = \setm{(1 + t) \ncdot u}{t \in K - \sngltn{0}}$,
$X = X' - C$,
$Y' = B^E(0,1)$,
$D = \setm{(1 - t) \ncdot u}{t \in K - \sngltn{0}}$,
and $Y = Y' - D$.

\num{a} $X$ and $Y$ are BPD.AC, BR.LC.AC and UD.AC.

\num{b} There is $\tau \in H(X,Y)$ such that
$(\rfs{BPD.UC}(X))^{\tau} \subseteq \rfs{BPD.UC}(Y)$
and $\tau\inverse \not\in \rfs{BPD.UC}(Y,X)$.

\num{c} There is $\tau \in H(X,Y)$ such that
\begin{itemize}
\addtolength{\parskip}{-11pt}
\addtolength{\itemsep}{06pt}
\item[\num{1}]
$(\rfs{BPD.UC}(X))^{\tau} \subseteq \rfs{BPD.UC}(Y)$,
\item[\num{2}]
$\tau\inverse \not\in \rfs{BPD.P}(Y,X)$,
\item[\num{3}]
for every BPD set $A \subseteq X$, $\tau \nrestriction A$
is bilipschitz.
\vspace{-05.7pt}
\end{itemize}
\end{example}

\noindent
{\bf Proof } 
(a) This part is trivial, so we leave its verification to the reader.
In any case, Part~(a) shows that the fact that the boundaries
of $X$ and $Y$ are well-behaved, does not by itself imply
that $\tau\inverse \in \rfs{BPD.UC}(Y,X)$.

(b) This part follows from (c). So it suffices to prove (c).

(c) Note the following facts:
$C \subseteq B^E(0,\dgfrac{3}{2}) - \overB^E(0,1)$,
$D \subseteq B^E(0,1) - \overB^E(0,\dghalf)$,
$u \in \rfs{acc}(C)$ and $u \in \rfs{acc}(D)$.

{\thickmuskip=3.5mu \medmuskip=2.5mu \thinmuskip=2mu 
Let $y \in B^{\oversE}(0,\dghalf) - B^E(0,\dghalf)$.
By Proposition \ref{metr-bldr-p2.25}(b),
there is $\rho \in \rfs{LIP}(\overE)\sprtl{B^{\oversE}(0,\dghalf)}$}
such that $\rho(0) = y$
and $\rho(E - \sngltn{0}) = E$.
So $\rho(D) = D$ and hence $\rho(Y - \sngltn{0}) = Y$.
Let $\fnn{\eta}{X}{Y - \sngltn{0}}$ be defined by
$\eta(x) = (2 - \norm{x}) \mcdot \frac{x}{\norm{x}}$
and $\tau = \rho \scirc \eta$. Clearly, $\tau \in H(X,Y)$,
and it is easy to check that $\tau$ satisfies Clause (3).

Let $r > 0$ be such that
$B^{\oversE}(y,r) \subseteq B^{\oversE}(0,\dghalf)$
and $M = B^{\oversE}(y,r) \cap E$.
Then $M$ is a BPD subset of $Y$.
However, $\tau\inverse(M)$ contains a set of the form
$B^E(0,2) - \overB(0,2 - \varepsilon)$, where $\varepsilon > 0$.
So $\tau\inverse(M)$ is not a BPD subset of $X$.
Hence Clause (4) is fulfilled.

We show that $\tau$ fulfills Clause (1).
It is easy to check that
$(\rfs{BPD.P}(X))^{\tau} \subseteq \rfs{BPD.P}(Y)$.
So it remains to show that if $h \in \rfs{BPD.UC}(X)$ and
$M \subseteq Y$ is a BPD set, then $h^{\tau} \nrestriction M$ is bi-UC.

Since $\rho$ is bilipschitz it suffices to show that
for every $h \in \rfs{BPD.UC}(X)$
and $M \subseteq Y - \sngltn{0}$:
if $d(M,D \cup S(0,1)) > 0$,
then $h^{\eta} \nrestriction M$ is bi-UC.
(Indeed we show that $h^{\eta} \nrestriction M$ is bi-UC,
even for $M$'s which satisfy $M \subseteq Y - \sngltn{0}$
and $d(M,D) > 0$).

{\bf Claim 1.} Let $Z,W$ be metric spaces, $z \in Z$,
and $\fnn{f}{Z}{W}$.
Suppose that $f$ is continuous at $z$,
and for every $r > 0$, $f \nrestriction (Z - B(z,r))$ is UC.
Then $f$ is UC.
{\bf Proof } Let $\varepsilon > 0$. There is $\delta_1 > 0$
such that $\rfs{diam}(f(B(z,\delta_1))) < \varepsilon$.
Let $\delta_2 > 0$ be such that for every
$x,y \in Z - B(z,\dgfrac{\delta_1}{2})$: if $d(x,y) < \delta_2$,
then $d(f(x),f(y)) < \varepsilon$.
Let $\delta = \min(\dgfrac{\delta_1}{2},\delta_2)$.
Suppose that $d(x,y) < \delta$.
Either $x,y \in Z - B(z,\dgfrac{\delta_1}{2})$
or $x,y \in B(z,\delta_1)$.
In either case $d(f(x),f(y)) < \varepsilon$.
Claim~1 is proved.

{\bf Claim 2.} Let $h \in \rfs{BPD.P}(X)$ and
$\vecx \subseteq X$ be such that \hbox{$\limti{n} \norm{x_n} = 2$.}
Then $\limti{n} \norm{h(x_n)} = 2$.
{\bf Proof } Suppose by way of contradiction that this is not true,
and let $\vecx$ be a counter-example.
Since $h \in \rfs{BPD.P}(X)$, for every subsequence $\vecx\fprime$ of
$\vecx$, $h(\vecx\fprime)$ is not a BPD sequence.
It follows easily that either $\vecx$ has a subsequence $\vecx\fprime$
such that $\limti{n} \norm{h(x'_n)} = 1$,
or $\vecx$ has a subsequence $\vecx\fprime$
which converges to a member of $C$.
Taking $\vecx$ to be $\vecx\fprime$
we may assume that either (i) $\limti{n} \norm{h(x_n)} = 1$
or (ii) for some $\hatu \in C$, $\lim h(\vecx) = \hatu$.

Suppose that (i) happens. Then
for every $n \in \bbN$ there are $u_n \in C$, $s_n > r_n > 0$
and an arc $L_n \subseteq X$ such that the following hold. 
\begin{itemize}
\addtolength{\parskip}{-11pt}
\addtolength{\itemsep}{06pt}
\item[(1)] 
$h(x_n) \in L_n$
and $L_n$ intersects both $S(u_n,r_n)$ and $S(u_n,s_n)$.
\item[(2)] 
Set $S_n = S(u_n,r_n) \cup S(u_n,s_n)$. Then $\delta^X(S_n) > 0$.
(Hence $S_n \subseteq X$).
\item[(3)] 
Define $d_n = \sup(\setm{d(z,S(0,1))}{z \in L_n \cup S_n})$.
Then $\limti{n} d_n = 0$.
\item[(4)] 
$(B(u_n,s_n) - \overB(u_n,r_n)) \cap C \neq \emptyset$.
\vspace{-05.7pt}
\end{itemize}

Suppose that (ii) happens. Then
for every $n \in \bbN$ there are $s_n > r_n > 0$
and an arc $L_n \subseteq X$ such that the following hold. 
\begin{itemize}
\addtolength{\parskip}{-11pt}
\addtolength{\itemsep}{06pt}
\item[(5)] 
$h(x_n) \in L_n$, and $L_n$ intersects both $S(\hatu,r_n)$
and $S(\hatu,s_n)$.
\item[(6)] 
Set $S_n = S(\hatu,r_n) \cup S(\hatu,s_n)$.
Then $\delta^X(S_n) > 0$. (Hence $S_n \subseteq X$).
\item[(7)] 
Define $d_n = \sup(\setm{d(z,\hatu)}{z \in L_n \cup S_n})$.
Then $\limti{n} d_n = 0$.
\item[(8)] 
$(B(\hatu,s_n) - \overB(\hatu,r_n)) \cap C \neq \emptyset$.
\vspace{-05.7pt}
\end{itemize}

In both Case (i) and Case (ii) denote $A_n = L_n \cup S_n$
and $B_n = h\inverse(A_n)$.
Let $\vecz$ be a sequence such that $z_n \in B_n$
for every $n \in \bbN$.
By (3) and (7), $\limti{n} \delta^X(h(z_n)) = 0$.
From the fact that $h \in \rfs{BPD.P}(X)$ it follows that
$\limti{n} \delta^X(z_n) = 0$.
There is a subsequence $\setm{n_i}{i \in \bbN}$
such that either
$\limti{n} d(z_{n_i},S(0,2)) = 0$
or
$\limti{n} d(z_{n_i},S(0,1) \cup C) = 0$.\break
We may assume that $n_i = i$ for every $i$.
Suppose by contradiction that the latter happens.
Now, $x_n,z_n \in B_n$, $B_n$ is connected
and $\limti{n} d(x_n,S(0,2)) = 0$.
We also have that $d(S(0,2),S(0,1) \cup C) > 0$.
Choose $y_n \in B_n$ such that
$\norm{y_n - x_n} = \dgfrac{\norm{z_n - x_n}}{2}$.
Then $d(\setm{y_n}{n \in \bbN},\rfs{bd}(X)) > 0$, a contradiction.
So $\limti{n} d(z_n,S(0,2)) = 0$.

Let $e_n = \sup(\setm{d(z,S(0,2))}{z \in B_n})$.
It follows that $\limti{n} e_n = 0$.
Let $n$ be such that $e_n \leq \dgquarter$.
Denote $S = S_n$ and $T = h\inverse(S)$.
Since $S$ is a BPD set, $T$ is a BPD set.
Let $d = d(T,S(0,2))$.
It is obvious that $X - S$ has three connected components,
and neither of them is a BPD set. So the same holds for $T$.
However, since $e_n \leq \dgquarter$,
$T \subseteq B(0,2) - \overB(0,\dgfrac{3}{2})$ and so
$X \cap B(0,\dgfrac{3}{2})$ is contained in a component of $X - T$,
and $B(0,2) - \overB(0,2 - d)$ is also contained in a
component of $X - T$.
It follows that one of the components of $X - T$ is contained in
$W \eqdf \overB(0,2 - d) - B(0,\dgfrac{3}{2})$.
But $W$ is a BPD subset of $X$. A contradiction, so Claim 2 is proved.

\kern1.5mm

Let $h \in \rfs{BPD.UC}(X)$ and denote $g = h^{\eta}$.

{\bf Claim 3.}
$0 \in \rfs{Dom}(g^{\srfs{cl}})$ and $g^{\srfs{cl}}(0) = 0$.
{\bf Proof }
Let $\vecx \subseteq B(0,1) - \sngltn{0}$
be such that $\lim \vecx = 0$.
Then $\limti{n} \norm{\eta\inverse(x_n)} = 2$.
Note that $h \in \rfs{BPD.P}(X)$. Applying Claim 2 to $h$,
we conclude that
$\limti{n} \norm{h(\eta\inverse(x_n))} = 2$.
Hence
$\limti{n} \norm{\eta(h(\eta\inverse(x_n)))} = 0$.
That is, $\limti{n} \norm{g(x_n)} = 0$.
So Claim 3 is proved.
\smallskip

Let $M \subseteq Y - \sngltn{0}$ be such that $d(M,D) > 0$.
Let $r > 0$ and
$N = \eta\inverse(M - B(0,r))$.
Then $d(N,S(0,2)) \geq r$. So $\eta \nrestriction N$ is bilipschitz,
hence (i) $\eta\inverse \nrestriction (M - B(0,r))$ is bilipschitz.
$N$ is a BPD subset of $X$.
So (ii) $h \nrestriction N$ is bi-UC.
Also, $h(N)$ is a BPD subset of $X$.
In particular, $d(h(N),S(0,2)) > 0$.
So (iii) $\eta \nrestriction h(N)$ is bilipschitz.
\vspace{1.5mm}
\newline
\rule{7pt}{0pt}
\renewcommand{\arraystretch}{1.5}
\addtolength{\arraycolsep}{-3pt}
$
\begin{array}{ll}
&
g \nrestriction (M - B(0,r)) =
\eta \scirc h \scirc \eta\inverse \nrestriction (M - B(0,r))
\\
=
\rule{5pt}{0pt}
&
(\eta \nrestriction h(\eta\inverse(M - B(0,r)))) \scirc
(h \nrestriction \eta\inverse(M - B(0,r))) \scirc
(\eta\inverse \nrestriction (M - B(0,r)))
\\
=
\rule{5pt}{0pt}
&
\eta \nrestriction h(N) \scirc
(h \nrestriction N) \scirc
(\eta\inverse \nrestriction (M - B(0,r))).
\end{array}
$
\vspace{1.7mm}
\renewcommand{\arraystretch}{1.0}
\addtolength{\arraycolsep}{3pt}
\newline
By (i)\,-\,(iii), $g \nrestriction (M - B(0,r))$ is bi-UC.
By Claim 3 and Claim 1, $g^{\srfs{cl}} \nrestriction M$ is UC.
Applying the same argument to $h\inverse$ we conclude that
$(g^{\srfs{cl}})\inverse \nrestriction g(M)$ is UC.
So $g \nrestriction M$ is bi-UC.
That is, $h^{\eta} \nrestriction M$ is bi-UC.
It has been already argued that this implies that
$h^{\tau} \in \rfs{BPD.UC}(Y)$.
\rule{18pt}{0pt}\hfill\myqed

\begin{lemma}\label{l4.25}
\label{metr-bldr-l5.38}
Suppose that $X$ is an open subset of a Banach space $E$.

\num{a} $\rfs{BUC}(X) \subseteq \rfs{BPD.UC}(X)$.

\num{b}
Suppose that $X$ is BR.LC.AC, $\tau \in H(X,Y)$
and $(\rfs{BUC}(X))^{\tau} \subseteq \rfs{BPD.P}(Y)$.
Then $\tau\inverse \in \rfs{BPD.P}(Y,X)$.
\end{lemma}

\noindent
{\bf Proof }
(a) Let $h \in \rfs{BUC}(X)$.
Suppose that $x \in \rfs{bd}(X)$, $\vecx \subseteq X$
and $\lim \vecx = x$.
Then $h(\vecx)$ is a Cauchy sequence. Let $y = \lim h(\vecx)$.
Clearly, $y \in \rfs{bd}^{\oversE}(X) \cup \overfs{int}(X)$
and $y \not\in X$.
Since $E$ is complete, $\overfs{int}(X) = X$. Hence $y \in \rfs{bd}(X)$.
We have shown that $\rfs{Dom}(h^{\srfs{cl}}) = \rfs{cl}(X)$
and that $h^{\srfs{cl}}(\rfs{bd}(X)) \subseteq \rfs{bd}(X)$.
Applying the same argument to $h\inverse$ one concludes that
$(\dagger)$ $h^{\srfs{cl}}(\rfs{bd}(X)) = \rfs{bd}(X)$.
It is trivial that $h^{\srfs{cl}} \in \rfs{BUC}(\rfs{cl}(X))$.

Suppose by contradiction that $A$ is a BPD set and $h(A)$ is not
a BPD set. By definition, $h$ is boundedness preserving.
So $h(A)$ is bounded and hence $\delta(h(A)) = 0$.
Let $\vecx \subseteq h(A)$ and $\vecy \subseteq \rfs{bd}(X)$
be such that $\limti{n} d(x_n,y_n) = 0$.
By $(\dagger)$,
$(h^{\srfs{cl}})\inverse(\vecy) \subseteq \rfs{bd}(X)$.
So for every $n$, $d(h\inverse(x_n),h\inverse(y_n)) \geq \delta(A) > 0$.
Hence
$(h^{\srfs{cl}})\inverse \nrestriction (\rfs{Rng}(\vecx)
\cup \rfs{Rng}(\vecy))$
is not uniformly continuous. A contradiction.

(b) Let $X,E,Y$ and $\tau$ be as in Part (b), and suppose that
$Y$ is an open subset of the normed space $F$.
Then $F$ is a Banach space.
To see this note that an open ball $B$ of $F$ is homeomorphic to
an open subset of $E$. So $B$ is completely metrizable. But
$F \cong B$, so $F$ is completely metrizable. So $F$ is a
dense $G_{\delta}$ subset of $\overF$, and so is every coset of $F$ in
$\overF$. Since $\overF$ has no disjoint dense $G_{\delta}$
subsets, $F = \overF$.
Suppose be contradiction $\tau\inverse \not\in \rfs{BPD.P}(Y,X)$.
Then there is a $\onetoonen$ sequence $\vecx \subseteq X$ such that
$\vecx$ is not a BPD sequence, but $\tau(\vecx)$ is a BPD sequence.
We may assume that $\limti{n} \delta_1^X(x_n) = \infty$.

Since $\tau(\vecx)$ is a BPD set,
it does not have convergent subsequences in $F$,
hence $\tau(\vecx)$ does not have Cauchy subsequence.
So we may assume that there is $d > 0$ such that $\tau(\vecx)$ is 
$d$-spaced.

{\bf Claim 1.} $\vecx$ is not a Cauchy sequence. {\bf Proof }
Suppose otherwise, and let $x^* = \lim \vecx$.
Then $x^* \in \rfs{bd}(X)$,
for if $x^* \not\in X$, then $\vecx$ is a BPD sequence.

By the simplicity of $x^*$, we can find a subsequence $\vecy$ of
$\vecx$, arcs $\setm{L_n}{n \in \bbN}$ and open sets
$\setm{U_n}{n \in \bbN}$ such that
$y_{2n},y_{2n + 1} \in L_n \subseteq U_n \subseteq \rfs{cl}(U_n)
\subseteq X$, for any distinct $m,n \in \bbN$, $d(U_m,U_n) > 0$,
and $\limti{n} \rfs{diam}(U_n) = 0$.

Let $z_n \in L_n - \sngltn{y_{2n}}$ be such that
$\limti{n} d(\tau(y_{2n}),\tau(z_n)) = 0$. It follows that
\begin{itemize}
\addtolength{\parskip}{-11pt}
\addtolength{\itemsep}{06pt}
\item[(1)] 
$\tau(\vecy) \cup \tau(\vecz)$ is a BPD set.
\vspace{-05.7pt}
\end{itemize}
Let $h_n \in \rfs{UC}(X)$ be such that
$h_n(y_{2n}) = y_{2n}$,
$h_n(z_n) = y_{2n + 1}$ and \hbox{$\rfs{supp}(h_n) \subseteq U_n$.}
By Proposition \ref{p4.9},
$h \eqdf \bcirc_{n \in \sboldbbN} h_n \in \rfs{UC}(X)$.
However,
\begin{itemize}
\addtolength{\parskip}{-11pt}
\addtolength{\itemsep}{06pt}
\item[(2)] 
$h^{\tau} \nrestriction (\tau(\vecy) \cup \tau(\vecz))$ is not UC.
\vspace{-05.7pt}
\end{itemize}
To see this recall that $\limti{n} d(\tau(y_{2n}),\tau(z_n)) = 0$.
However,
$d(h^{\tau}(\tau(y_{2n})),h^{\tau}(\tau(z_n))) = 
d(\tau(y_{2n}),\tau(y_{2n + 1})) \geq d$.
Facts (1) and (2) mean that $h^{\tau} \not\in \rfs{BPD.UC}(Y)$.
A contradiction, so Claim 1 is proved.

{\bf Claim 2.} It is not true that $\limti{n} \delta(x_n) = 0$.
{\bf Proof } Suppose otherwise. By Claim~1, we may assume that there is
$e_1 > 0$ such that $\vecx$ is $e_1$-spaced. For every $n \in \bbN$
let $b_n \in \rfs{bd}(X)$ be such that $d(x_n,b_n) \leq 2\delta(x_n)$,
and $[x_n,b_n) \subseteq X$.

For every $n \in \bbN$ let
$\vecx^n = \setm{x^n_i}{i \in \bbN} \subseteq [x_n,b_n)$ be a sequence
converging to $b_n$. By Claim 1, $\tau(\vecx^n)$ is not a BPD set.
It follows that there is a sequence $\setm{i_n}{n \in \bbN}$ such that
$\setm{\tau(x^n_{i_n})}{n \in \bbN}$ is not a BPD set.
Let $y_n = x^n_{i_n}$. Since $\vecx$ is $e_1$-spaced and
$\limti{n} d(x_n,y_n) = 0$, we may assume that there is $e > 0$
such that $\setm{[x_n,y_n]}{n \in \bbN}$ is $e$-spaced.

Let $\setm{U_n}{n \in \bbN}$ be a sequence of open subsets of $X$
such that
$[x_n,y_n] \subseteq U_n$,
$\limti{n} \rfs{diam}(U_n) = 0$
and for any distinct $m,n \in \bbN$, $d(U_m,U_n) > 0$.
Let $h_n \in \rfs{UC}(X)$ be such that $\rfs{supp}(h_n) \subseteq U_n$
and $h_n(x_n) = y_n$.
By Proposition~\ref{p4.9},
$h \eqdf \bcirc_{n \in \sboldbbN} h_n \in \rfs{UC}(X)$, but
$h^{\tau} \not\in \rfs{BPD.P}(Y)$. This is so, because
$h^{\tau}(\tau(\vecx)) = \tau(\vecy)$,
and $\tau(\vecx)$ is a BPD set, whereas $\tau(\vecy)$ is not.
A contradiction. This proves Claim 2.\medskip

From Claims 1 and 2 and the fact that $\vecx$ is not a BPD sequence,
it follows that $\vecx$ is unbounded.
So we may assume that $\setm{\norm{x_n}}{n \in \bbN}$
is a strictly increasing sequence converging to $\infty$.
Recall also that $\tau(\vecx)$ is $d$-spaced.
We now deal with two cases.

{\bf Case 1 } $E - X$ is bounded. We may assume that
$E - X \subseteq B(0,\dgfrac{\norm{x_0}}{2})$.
Set $x_{-1} = 0$.
Choose $y_n \in (x_{2n},x_{2n + 1}]$ such that
$d(\tau(x_{2n}),\tau(y_n)) < \dgfrac{1}{(n + 1)}$.
Define $r_n =
\dgfrac{\min(\norm{x_{2n} - x_{2n - 1}},
\norm{x_{2n + 2} - x_{2n + 1}})}{2}$
and let $h_n \in \rfs{UC}(X)$ be such that
$h_n(x_{2n}) = x_{2n}$,
$h_n(y_n) = x_{2n + 1}$ and
$\rfs{supp}(h_n) \subseteq B([x_{2n},x_{2n + 1}],r_n)$.
Clearly, $\rfs{supp}(h_m) \cap \rfs{supp}(h_n) = \emptyset$
for every $n \neq m$ and hence
$h \eqdf \bcirc_{n \in \sboldbbN} h_n \in \rfs{BUC}(X)$.
Since $\limti{n} d(\tau(x_{2n}),\tau(y_n)) = 0$,
it follows that
$\tau(\vecx) \cup \tau(\vecy)$ is a BPD set.
But $h^{\tau} \nrestriction (\tau(\vecx) \cup \tau(\vecy))$ is not UC.
So $h^{\tau} \not\in \rfs{BPD.UC}(Y)$.
A contradiction, so Case 1 does not occur.

{\bf Case 2 } $E - X$ is unbounded.
We define by induction on $n \in \bbN$: $u_n \in \rfs{Rng}(\vecx)$,
$v_n \in X$, $h_n \in \rfs{UC}(X)$ and $r_n > 0$.
Let $r_{-1} = 0$. Suppose that $r_{n - 1}$ has been defined.
Let $u_n \in \rfs{Rng}(\vecx) - \rfs{cl}(B(0,r_{n - 1}))$ and
Let
$b_n \in \rfs{bd}(X) - \rfs{cl}(B(0,r_{n - 1}))$.
We may assume that there is an arc
$L_n \subseteq (X \cup \sngltn{b_n}) - \rfs{cl}(B(0,r_{n - 1}))$
connecting $u_n$ and $b_n$.
Let\break
$\vecv^n \eqdf \setm{v^{n,i}}{i \in \bbN} \subseteq L_n - \sngltn{b_n}$
be a sequence converging to $b_n$. So $\vecv^n$ is a Cauchy sequence.
So by Claim 1, $\tau(\vecv^n)$ is not a BPD set.
Hence there is
$v_n \in L_n - \sngltn{b_n}$ such $\delta_1(\tau(v_n)) > n$.

Let $r_n$ be such that $L_n \subseteq B(0,r_n)$ and
$h_n \in \rfs{UC}(X)$ be such that
$h_n(u_n) = v_n$ and
$\rfs{supp}(h_n) \subseteq B(0,r_n) - \rfs{cl}(B(0,r_{n - 1}))$.
Clearly,
$\rfs{supp}(h_m) \cap \rfs{supp}(h_n) = \emptyset$ for every $m \neq n$,
and hence $h \eqdf \bcirc_{n \in \sboldbbN} h_n \in \rfs{BUC}(X)$.
However, since $\tau(\vecu)$ is a BPD sequence,
$\tau(\vecv)$ is not a BPD sequence,
and $h^{\tau}(\tau(\vecu)) = \tau(\vecv)$,
$h^{\tau} \not\in \rfs{BPD.P}(Y)$.
A contradiction, so Case 2 does not happen.
It follows that $\tau\inverse \in \rfs{BPD.P}(Y,X)$.
\smallskip\hfill\myqed

If $X$ is BPD.AC, and we remove from $X$ a spaced set,
then the resulting open set is also BPD.AC.
This is proved in the next proposition.
Although this fact is quite trivial, a complete proof
requires much writing.

\begin{prop}\label{p4.26}\label{metr-bldr-p5.39}
\num{a} Let $E$ be a normed space which is not $1$-dimensional.
Let $u,v,w \in E$ be such that
$\norm{u - w} = \norm{v - w} = r > 0$.
Then there is an arc $L \subseteq X$
connecting $u$ and $v$ such that $L \cap B(w,r) = \emptyset$, and
$\rfs{lngth}(L) \leq 8r$.

\num{b} Suppose that $X$ is BPD.AC, and is not $1$-dimensional.
If $A \subseteq X$ is spaced, then $X - A$ is BPD.AC.
\end{prop}

\noindent
{\bf Proof }
(a) We may assume that $E$ is 2-dimensional, $w = 0$ and $r = 1$.
Let $z \in S(0,1)$ be such that $\ell \eqdf \setm{u + tz}{t \in \bbR}$
is a supporting line for $B(0,1)$.
Represent $v$ as $v = au + bz$,
and choose $z$ in such a way that $b > 0$.
Let $L_1 = [u,u + 2z]$, $L_2 = [2z + u,2z - u]$,
$L_3 = [2z - u,-u]$ and $L_0 = L_1 \cup L_2 \cup L_3$.
Since $\ell$ is a supporting line of $B(0,1)$ it follows that
$L_1$ and $L_3$ are disjoint from $B(0,1)$.
Suppose that $w \in L_2$. So $w = 2z + tu$, where $\abs{t} \leq 1$.
We may assume that $t \geq 0$.
Then $\norm{w} \geq 2\norm{z} - t\norm{u} \geq 1$.
So $L_2 \cap B(0,1) = \emptyset$.
Recall that $v = au + bz \in S(0,1)$.
From the fact that $\ell$ supports $B(0,1)$ it follows that
$a \leq 1$. Then $1 = \norm{v} \geq b - a \geq b - 1$.
So $b \leq 2$.
Let $\lambda = \min(\dgfrac{1}{\abs{a}},\dgfrac{2}{b})$
and $L_v = [v,\lambda v]$.
Clearly, $L_v \cap B(0,1) = \emptyset$.
Either  $\lambda v = u + b_1 z$, where $b_1 \in [0,2]$,
or
$\lambda v = -u + b_1 z$, where $b_1 \in [0,2]$,
or
$\lambda v = a_1 u + 2z$, where $a_1 \in [-1,1]$.
Hence $\lambda v \in L_1 \cup L_3 \cup L_2 = L_0$.
The set $L_0 \cup L_v$  is disjoint from $B(0,1)$
and contains an arc $L$ connecting $u$ and $v$.
Obviously, for $i = 1,\ldots,3$, $\rfs{lngth}(L_i) = 2$
and
$\rfs{lngth}(L_v) = \norm{\lambda v} - \norm{v} \leq
2\norm{z} + \norm{u} - 1 = 2$.
So $\rfs{lngth}(L) \leq 8$.

(b) We prove Claim 1 stated below, and leave it to the reader to
verify that (b) is implied by Claim 1.

{\bf Claim 1.}
For every $r,C,D > 0$ there are \hbox{$r_1,C_1,D_1 > 0$}
such that for every
normed space $E$, an open subset $X \subseteq E$ and an $r$-spaced
subset $A \subseteq X$ the following holds.
If $x,x^* \in X - A$ are such that
$d(\sngltn{x,x^*},A) \geq r$,
and $L \subseteq X$ is an arc connecting $x$ and $x^*$ such that
$\delta^X(L) \geq C$ and
$\rfs{lngth}(L) \leq D$,
then there is an arc $M \subseteq X - A$
connecting $x$ and $x^*$ such that
$d(M,A) \geq r_1$,
$\delta^X(M) \geq C_1$ and
$\rfs{lngth}(M) \leq D_1$.

{\bf Proof } Let $D_1 = 8D$, $C_1 = \dgfrac{C}{2}$ and
$r_1 = \dgfrac{\min(r,C)}{64}$.
Let $E,\,X,\,A,\,x,\,x^*$ and $L$ be as in the claim and
$\fnn{\gamma}{[0,1]}{L}$ be a parametrization of $L$ which satisfies
$\rfs{lngth}(\gamma \nrestriction [0,t]) = t \cdot \rfs{lngth}(L)$
for every $t \in [0,1]$.
For every $a \in A$
let $T_a = \setm{t \in [0,1]}{\gamma(t) \in B(w,2r_1)}$.
Clearly, $T_a$ is an open subset of $(0,1)$,
and $\rfs{cl}(T_a) \cap \rfs{cl}(T_b) = \emptyset$
for any distinct $a,b \in A$.
Define $T = \bigcup \setm{T_a}{a \in A}$,
and let $\calI$ be a set of pairwise disjoint open intervals
of $(0,1)$ such that $\bigcup \calI = T$.
For an open interval $I$ in $(0,1)$ denote by $s_I$ and $t_I$
the left and right endpoints of $I$,
and if $I \in \calI$ denote by $a_I$
that member of $A$ such that $I \subseteq T_a$.
Clearly, $s_I,t_I \in S(a_I,2r_1)$.
For every $I \in \calI$ let $L_I = \gamma([s_I,t_I])$
and $M_I$ be a rectifiable arc connecting $a_I$
and $b_I$ such that $M_I \cap B(a_I,2r_1) = \emptyset$
and $\rfs{lngth}(M_I) \leq 16r_1$. 
The existence of $M_I$ is assured by Part (a).
Let $\calI_0 = \setm{I \in \calI}{d(L_I,a_I) \leq r_1}$.
Let
$$\hbox{$
M = L - \bigcup_{I \in \calI_0} L_I \cup
\bigcup_{I \in \calI_0} M_I.
$}$$
Certainly, $M$ is an arc whose endpoints are $x$ and $x^*$.
It is trivial that if $I \in \calI_0$,
then $\rfs{lngth}(L_I) \geq 2r_1$,
and so for every $I \in \calI_0$,
$\dgfrac{\rfs{lngth}(M_I)}{\rfs{lngth}(L_I)} \leq 8$.
It follows that $M$ is rectifiable
and that $\rfs{lngth}(M) \leq 8 \mcdot \rfs{lngth}(L) \leq 8D$.

Let $w \in M$. If $w \in L - \bigcup_{I \in \calI_0} L_I$,
then $d(w,A) \geq 2 r_1$. If there is $I \in \calI_0$
such that $w \in M_I$, then $d(w,a_I) \geq 2r_1$
and for every $b \in A - \sngltn{a_I}$,
$$
d(w,b) \geq d(b,a_I) - d(w,a_I) \geq r - 8r_1 - 2r_1 \geq
64r_1 - 10r_1 = 54r_1.
$$
It follows that $d(M,A) \geq r_1$.

It remains to show that $\delta^X(M) \geq \dgfrac{C}{2}$.
Obviously,
$\delta^X(L - \bigcup_{I \in \calI_0} L_I) \geq \delta^X(L) \geq C$.
Let $I \in \calI_0$ and be such that $w \in M_I$.
Then
$$
d(w,E - X) \geq d(a_I,E - X) - d(w,a_I) \geq
C - 8r_1 - 2r_1 = C - 10r_1 \geq C - 16r_1 \geq \dgfrac{C}{2}.
$$
It follows that $\delta^X(M) \geq \dgfrac{C}{2}$.
We have proved Claim 1.
\bigskip\hfill\myqed

We are ready to prove that for open subsets of Banach spaces,
if $(\rfs{BUC}(X))^{\tau} \subseteq \rfs{BPD.UC}(Y)$,
then $\tau\inverse \in \rfs{BPD.UC}(Y,X)$.
This is the contents of Part (a) of the next theorem.
The main argument lies though in Part (b), and once it is proved,
(a) follows easily. So we shall start with the proof of (b).

\begin{theorem}\label{t4.27}\label{metr-bldr-t5.40}
Let $E$ be a Banach space and $X$ be an open subset of $E$.

\num{a} Suppose that $X$ is BPD.AC and BR.LC.AC,
and that $\tau \in H(X,Y)$ is such that
$(\rfs{BUC}(X))^{\tau} \subseteq \rfs{BPD.UC}(Y)$.
Then $\tau\inverse \in \rfs{BPD.UC}(Y,X)$.

\num{b} Suppose that $X$ is BPD.AC,
and that $\tau \in H(X,Y)$ is such that
$(\rfs{LIP}_{00}(X))^{\tau} \subseteq \rfs{BPD.UC}(Y)$.
Assume further that $\tau\inverse \in \rfs{BPD.P}(Y,X)$.
Then $\tau\inverse \in \rfs{BPD.UC}(Y,X)$.
\end{theorem}

\noindent
{\bf Proof }
(b)
We shall see that the proof of (b) can be reduced to an instance
of Lemma~\ref{l-gamma.5}.
Suppose by contradiction that
$\tau\inverse \not\in \rfs{BPD.UC}(Y,X)$.
So there are sequences $\vecx\fprime,\vecy\fprime$
in $Y$ and $e > 0$ such that
$\rfs{Rng}(\vecx\fprime) \cup \rfs{Rng}(\vecy\fprime)$
is a BPD subset of $Y$,
$\limti{n} d(x'_n,y'_n) = 0$,
and $d(\tau\inverse(x'_n),\tau\inverse(y'_n)) > e$
for every $n \in \bbN$.
We may assume that $\vecx\fprime$ is either a Cauchy sequence or
$\vecx\fprime$ is spaced.
However, $\vecx\fprime$ cannot be a Cauchy sequence because in that
case its limit belongs to $Y$,
and this violates the continuity of $\tau\inverse$.
So we may assume that $\vecx\fprime$ is spaced.
Set $\vecx = \tau\inverse(\vecx\fprime)$
and $\vecy = \tau\inverse(\vecy\fprime)$.
From the fact that $\tau\inverse \in \rfs{BPD.P}(Y,X)$ it follows
that $\rfs{Rng}(\vecx)$ is a BPD set.
We may assume that $\vecx$ is either spaced or is a Cauchy sequence.
But if it is a Cauchy sequence then its limit belongs to $X$,
and by the continuity of $\tau$ at $x$, $\vecx\fprime$ is a Cauchy
sequence, which we have already excluded.
So we may assume that $\vecx$ is spaced.
Let $d > 0$ be such that $\vecx$ is $d$-spaced.
Then for every $n \in \bbN$ there is at most one $m$ such that
$\norm{y_n - x_m} < \dgfrac{d}{2}$.
It follows that there is an infinite set $\eta \subseteq \bbN$,
such that $\norm{y_n - x_m} \geq \rfs{min}(e,\dgfrac{d}{2})$
for every $m,n \in \eta$.
We may thus assume that $d(\rfs{Rng}(\vecx),\rfs{Rng}(\vecy)) > 0$.

We denote $\rfs{Rng}(\vecx), \rfs{Rng}(\vecy),
\rfs{Rng}(\vecx\fprime)$ and $\rfs{Rng}(\vecy\fprime)$
by $A,B,A'$ and $B'$ respectively.
Let $\whatX = X - A$, $\whatY = Y - A'$
and $\hattau = \tau \nrestriction \whatX$.
So $\hattau \in H(\whatX,\whatY)$. We shall prove that
\begin{list}{}
{\setlength{\leftmargin}{39pt}
\setlength{\labelsep}{08pt}
\setlength{\labelwidth}{20pt}
\setlength{\itemindent}{-00pt}
\addtolength{\topsep}{-04pt}
\addtolength{\parskip}{-02pt}
\addtolength{\itemsep}{-05pt}
}
\item[(i)]
$\whatX$ is BPD.AC,
\item[(ii)]
$(\rfs{LIP}_{00}(\whatX))^{\hattau} \subseteq \rfs{BPD.P}(\whatY)$,
\item[(iii)]
$B$ is a BPD subset of $\whatX$,
whereas $\hattau(B)$ is not a BPD subset of $\whatY$.
\vspace{-02.0pt}
\end{list}
Facts (i)\,-\,(iii) contardict Lemma \ref{l-gamma.5}.

(i) By Proposition \ref{p4.26}(b), $\whatX$ is BPD.AC.

(ii) Let $h \in \rfs{LIP}_{00}(\whatX)$. Then $h$ is extendible, and
$h^{cl} \nrestriction \rfs{bd}(\whatX) = \rfs{Id}$.
So $h^{cl}(A) = A$. Hence 
$h^* \eqdf h^{cl} \nrestriction X \in H(X)$
and clearly, $h^* \in \rfs{LIP}_{00}(X)$.
So $(h^*)^{\tau} \in \rfs{BPD.UC}(Y)$.
We show that if $C$ is a BPD subset of $\whatY$,
then $h^{\hattau}(C)$ is a BPD subset of $\whatY$.
Clearly, $h^{\hattau} = (h^*)^{\tau} \nrestriction \whatY$.
Obviously, $C \cup A'$ is a BPD subset of $Y$, and hence
$(h^*)^{\tau} \nrestriction (C \cup A')$ is bi-UC.
So since $d(C,A') > 0$, $d((h^*)^{\tau}(C),(h^*)^{\tau}(A')) > 0$.
Since $(h^*)^{\tau}(A') = A$, it follows that
$(\dagger)$ $d((h^*)^{\tau}(C),A') > 0$.
Since $(h^*)^{\tau} \in \rfs{BPD.P}(Y)$,
and $C$ is a BPD subset of $Y$, we also have
$(\dagger\dagger)$ $(h^*)^{\tau}(C)$ is a BPD subset of $Y$.
From $(\dagger)$ and $(\dagger\dagger)$ it follows that
$(h^*)^{\tau}(C)$ is a BPD subset of $\whatY$.
That is, $h^{\hattau}(C)$ is a BPD subset of $\whatY$.
We have shown that for every $h \in \rfs{LIP}_{00}(\whatX)$,
$h^{\hattau}$ is BPD.P. The same holds for $h\inverse$,
so $(\rfs{LIP}_{00}(\whatX))^{\hattau} \subseteq \rfs{BPD.P}(\whatY)$.

(iii) Since $\tau\inverse \in \rfs{BPD.P}(Y,X)$ and $B'$ is a BPD
subset of $Y$, we have that $B$ is a BPD subset of $X$.
From the fact that $d(A,B) > 0$ we conclude that
$B$ is a BPD subset of $X - A = \whatX$.
On the other hand,
$d(A',B') = d(\rfs{Rng}(\vecx\fprime),\rfs{Rng}(\vecy\fprime)) = 0$,
so $B'$ is not a BPD subset of $\whatY$.

Facts (i)\,-\,(iii) contradict Lemma \ref{l-gamma.5},
so $\tau\inverse \in \rfs{BPD.UC}(Y,X)$.
Part (b) is thus proved.\kern8pt
\smallskip

(a) Let $X,Y,\tau$ be as in (a).
Then $(\rfs{BUC}(X))^{\tau} \subseteq \rfs{BPD.P}(Y)$.
So by Lemma~\ref{l4.25}(b), $\tau\inverse \in \rfs{BPD.P}(Y,X)$.
We also have that
$(\rfs{LIP}_{00}(X))^{\tau} \subseteq \rfs{BPD.UC}(Y)$.
So by Part (b) of this theorem, $\tau\inverse \in \rfs{BPD.UC}(Y,X)$.
\bigskip\hfill\myqed

\newpage

\section{Groups of extendible homeomorphisms and the\\
reconstruction of the closure of open sets}
\label{s6}

\subsection{General description.}
\label{ss6.1}

This chapter deals
with the homeomorphism groups of closed sets which are
the closure of an open subset of a normed space
and with groups of extendible homeomorphisms.
Under appropriate assumptions on the open sets $X$ and $Y$
we prove that if
$\iso{\varphi}{H(\rfs{cl}(X))}{H(\rfs{cl}(Y))}$,
then there is $\tau \in H(\rfs{cl}(X),\rfs{cl}(Y))$
such that $\tau$ induces $\varphi$.
Under the same assumptions we also prove that
if $\iso{\varphi}{\rfs{EXT}(X)}{\rfs{EXT}(Y)}$,
then there is $\tau \in \rfs{EXT}^{\pm}(X,Y)$
such that $\tau$ induces $\varphi$.
The definitions of
$\rfs{EXT}(X,Y)$ and $\rfs{EXT}(X)$ appear in
\ref{metr-bldr-d4.6}(b) and \ref{metr-bldr-d5.1}(a).

The results about $H(\rfs{cl}(X))$ appear in
Theorems \ref{metr-bldr-t6.20} and \ref{metr-bldr-t6.22},
and those about $\rfs{EXT}(X)$
appear in Theorems \ref{metr-bldr-t6.5},
\ref{metr-bldr-t6.10} and \ref{metr-bldr-t6.15}.
These theorems cover open subsets of a normed space whose boundary
may be quite complicated.
So they go far beyond the class of open sets whose closure
is a manifold with a boundary.
Nevertheless, the statements
\begin{itemize}
\addtolength{\parskip}{-11pt}
\addtolength{\itemsep}{06pt}
\item[] 
Every $\iso{\varphi}{H(\rfs{cl}(X))}{H(\rfs{cl}(Y))}$
is induced by some $\tau \in H(\rfs{cl}(X),\rfs{cl}(Y))$
\vspace{-05.7pt}
\end{itemize}
and
\begin{itemize}
\addtolength{\parskip}{-11pt}
\addtolength{\itemsep}{06pt}
\item[] 
Every $\iso{\varphi}{\rfs{EXT}(X)}{\rfs{EXT}(Y)}$
is induced by some $\tau \in \rfs{EXT}^{\pm}(X,Y)$
\vspace{-05.7pt}
\end{itemize}
are not true for every pair of open subsets of a normed space,
not even in the finite-dimensional case.
Example \ref{metr-bldr-e5.9} exhibits two different trivial reasons
why the above statements are not true in their full generality.

The proofs of the theorems about $\rfs{EXT}(X)$ and about
$H(\rfs{cl}(X))$ are essentially identical.
Moreover, for finite-dimensional normed spaces the question about the
faithfulness of $\setm{H(\rfs{cl}(X))}{X \mbox{ is open}}$
is a special case of the question about the EXT-detremined-ness of 
$\setm{X}{X \mbox{ is open}}$.
To see this, notice the following facts.
\begin{itemize}
\addtolength{\parskip}{-11pt}
\addtolength{\itemsep}{06pt}
\item[(1)] 
If $U$ is a regular open subset of $\bbR^n$,
then $\rfs{EXT}(U) = H(\rfs{cl}(U))$.
\item[(2)] 
If $X \subseteq \bbR^n$ is open and
$\whatX = \rfs{int}(\rfs{cl}(X))$,
then $\whatX$ is regular open and $\rfs{cl}(X) = \rfs{cl}(\whatX)$.
\vspace{-05.7pt}
\end{itemize}
Suppose now that $\iso{\varphi}{H(\rfs{cl}(X))}{H(\rfs{cl}(Y))}$.
By (2), $\iso{\varphi}{H(\rfs{cl}(\whatX))}{H(\rfs{cl}(\whatY))}$,
and by~(1),
$\iso{\varphi}{\rfs{EXT}(\whatX)}{\rfs{EXT}(\whatY)}$.
So if it can be proved that there is
$\tau \in \rfs{EXT}^{\pm}(\whatX,\whatY)$
such that $\tau$ induces $\varphi$,
then this $\tau$ indeed belongs to $H(\rfs{cl}(X),\rfs{cl}(Y))$.

Theorems \ref{metr-bldr-t6.5} and \ref{metr-bldr-t6.15}
prove the EXT-determined-ness of certain classes.
In \ref{metr-bldr-t6.5} it is assumed that the members of the
EXT-determined class are BR.LC.AC, see \ref{metr-bldr-d5.36}.
This property is a weakening of uniformly-in-diameter
arcwise connectedness.
It may happen though that every point in the boundary of such a set
is fixed under $\rfs{EXT}(X)$.
In \ref{metr-bldr-t6.15}, on the other hand, the EXT-determined-ness
is derived from the property that the $\rfs{EXT}(X)$-orbit
of every member of $\rfs{bd}(X)$ contains an arc,
but $X$ need not be BR.LC.AC.

In Corollary \ref{c4.33}(a) we prove that if $X$ and $Y$ satisfy
certain weak assumptions on arcwise connectedness,
and $(\rfs{EXT}(X))^{\tau} = \rfs{EXT}(Y)$,
then $\tau \in \rfs{EXT}(X,Y)$.
A statement of the form:
``$(\rfs{EXT}(X))^{\tau} \subseteq \rfs{EXT}(X)$ \kern1mm $\Rightarrow$
\kern1mm $\tau \in \rfs{EXT}(X,Y)$''
is also proved, but only under rather restrictive
assumptions on $X$ and $Y$. See Corollary \ref{metr-bldr-c6.5}(b).

Suppose that $X$ is an open subset of $\kern1pt\bbR^n$.
Then $\rfs{EXT}(X) = \rfs{BUC}(X)$.
If in addition, $X$ is bounded,
then $\rfs{EXT}(X) = \rfs{UC}(X)$.
So for finite-dimensional bounded $X$'s Corollary \ref{metr-bldr-c5.7}
which deals with $\rfs{BUC}(X)$
is indeed about $\rfs{EXT}(X)$.
However, Theorems \ref{metr-bldr-t6.10} and \ref{metr-bldr-t6.15}
are stronger than \ref{metr-bldr-c5.7} even for
finite-dimensional bounded $X$'s.

Groups of completely locally uniformly continuous homeomorphisms
are dealt with in Theorem \ref{metr-bldr-t6.17}.
(See definition \ref{metr-bldr-d5.3}(f)).
The $\itGamma$-continuous version of these groups is the subject
of Chapters \ref{s8}\,-\,\ref{s12}.

At the end of this chapter
in items \ref{metr-bldr-r6.23}\,-\,\ref{metr-bldr-t6.28}, we discuss
two generalizations of\break
the results.
The first generalization deals with subsets $Z$ of a normed space
such that $Z \subseteq \rfs{cl}(\rfs{int}(Z))$.
The second generalization deals with sets which are the closure of
an open subset in a normed manifold.

Recall that if not otherwise stated, then $X$ and $Y$ denote
respectively open subsets of the normed spaces $E$ and $F$.

\subsection{Groups of extendible homeomorphisms.}
\label{ss6.2}
\label{ss6.2-extendible-homeomorphism}

The following definition contains some notions related
to arcwise connectedness. These notions are used
in the statement of Theorem~\ref{metr-bldr-t6.5}
which deals with EXT-determined-ness.
In the next definition only, $E$ denotes a general metric space.

\begin{defn}\label{d4.30.1}\label{metr-bldr-d6.2.1}
\begin{rm}
Let $E$ be a metric space and $X \subseteq E$.

(a) A set of pairwise disjoint sets
is called a {\it pairwise disjoint family}.
Let $\calA$ be a pairwise disjoint family of subsets of $X$.
$\calA\kern1pt$ is {\it completely discrete} with respect to $E$,
if for every $x \in E$ there is $U \in \rfs{Nbr}(x)$ such that
$\setm{A \in \calA}{A \cap U \neq \emptyset}$ is finite.
A set $A \subseteq X$ is {\it completely discrete} with respect to $E$,
if $A$ does not have accumulation points in $E$.
The mention of $E$ in the above definition is omitted, since
$E$ is usually understood from the context.
A sequence $\vecx \subseteq X$ is a {\it completely discrete sequence},
if it is $\onetoonen$, and its range is completely discrete.

   \index{pairwise disjoint family. A set of pairwise disjoint sets}
   \index{completely discrete family of sets. A set $\calA$ of
   pairwise disjoint sets such that\newline\indent
   $\forall B((\forall A \in \calA)(\abs{B \cap A} \leq 1)
   \rightarrow \rfs{acc}(B) = \emptyset)$}
   \index{completely discrete set}
   \index{completely discrete sequence}

\num{b} $X$ is said to be {\it boundedly arcwise connected (BD.AC)},
if for every bounded $A \subseteq X$ there is $d > 0$ such that
for every $x,y \in A$ there is a rectifiable arc $L \subseteq X$
connecting $x$ and $y$ such that $\rfs{lngth}(L) \leq d$.

   \index{boundedly arcwise connected. Abbreviated by BD.AC}
   \index{bd.ac@@BD.AC. Abbreviation of boundedly arcwise connected}

\num{c}
$X$ is said to be a {\it wide set},
if for every infinite completely discrete set
$A \subseteq X$ there is an infinite $B \subseteq A$,
a set $\setm{y_b}{b \in B}$
and a set of arcs $\setm{L_b}{b \in B}$ such that:\break
$\setm{y_b}{b \in B}$ is bounded;
for every $b \in B$, \ $y_b,b \in L_b \subseteq X$;
and $\setm{L_b}{b \in B}$ is completely discrete.

   \index{wide set}

(d) Let $\vecx \subseteq X$ be a completely discrete sequence.
Let
$x^* \in \rfs{cl}(X)$,
$\setm{L_n}{n \in \bbN}$ be a sequence of arcs
and $\vecy \subseteq X$.
Assume that
\begin{itemize}
\addtolength{\parskip}{-11pt}
\addtolength{\itemsep}{06pt}
\item[(1)]
$L_n \subseteq X$ for every $n \in \bbN$,
\item[(2)]
$L_n$ connects $x_n$ with $y_n$ for every $n \in \bbN$,
\item[(3)]
$\lim \vecy = x^*$,
\item[(4)]
$L_m \cap L_n = \emptyset$
for any distinct $m,n \in \bbN$,
\item[(5)]
for every $r > 0$, $\setm{L_n - B^E(x^*,r)}{n \in \bbN}$
is completely discrete.
\vspace{-05.7pt}
\end{itemize}
Then $\fortpl{\vecx}{x^*}{\setm{L_n}{n \in \bbN}}{\vecy}$
is called a {\it joining system for $\vecx$ with respect to~$E$}.

(e) $X$ is {\it jointly arcwise connected (JN.AC) with respect to $E$},
if for every completely discrete sequence $\vecx \subseteq X$
there is a subsequence $\vecx\fprime$ of $\vecx$
such that $\vecx\fprime$ has a joining system.

   \index{joining system}
   \index{jointly arcwise connected  }
   \index{jnac@@JN.AC. Abbreviation of jointly arcwise connected}
\end{rm}
\end{defn}

In (a)-(d) of the next proposition we infer joint arcwise connectedness
from various simpler properties of $X$.
Part (e) is a trivial observation, so we do not prove it.

\begin{prop}\label{p4.31}\label{metr-bldr-p5.43}\label{metr-bldr-p6.3}
\num{a} Suppose that $\vecx \subseteq X$ is a Cauchy sequence and
$\lim^{\oversE} \vecx \in \overfs{int}(X) - X$.
Then $\vecx$ has subsequence $\vecx\fprime$
such that $\vecx\fprime$ has a joining system.

\num{b} Suppose that $X$ is an open subset
of a finite-dimensional normed space. Then
$X$ is JN.AC iff $X$ is bounded.

\num{c} Suppose that $X$ is an open subset of a Banach space
and $X$ is BD.AC. Then every bounded completely discrete sequence
$\vecx \subseteq X$ has a subsequence $\vecx\fprime$
such that $\vecx\fprime$ has a joining system.
In particular, if in addition $X$ is bounded, then $X$ is JN.AC.

\num{d}
If $X$ is an open subset of a Banach space,
$X$ is wide and $X$ is BD.AC,
then $X$ is JN.AC.

\num{e} Let $X$ be a bounded subset of a finite-dimensional
normed space. Then $X$ is BR.LC.AC iff $X$ is UD.AC.
\end{prop}

\noindent
{\bf Proof }
(a) Let $\barx = \lim^{\oversE} \vecx$.
Let $u \in E$ and $r > 0$ be such that $B(u,r) \subseteq E$
and  $\barx \in B^{\oversE}(u,r)$. Let $v \in B(u,r)$.
There is a subsequence $\vecy$ of $\vecx$ such that
$\vecy \subseteq B(u,r)$ and
$\setm{[y_n,v)}{n \in \bbN}$ is a pairwise disjoint family.
Let $v_n \in [y_n,v)$ be such that $\lim \vecv = v$.
Then 
$\fortpl{\vecy}{v}{\setm{[y_n,v_n]}{n \in \bbN}}{\vecv}$
is a joining system for $\vecy$.

(b) If $X$ is a bounded open subset of a finite-dimensional space,
then $X$ does not contain an infinite completely discrete
set. So $X$ is JN.AC.

Suppose that $X$ is an unbounded open subset of a finite-dimensional
space, Let $\vecx \subseteq X$ be a 1\,-\,1 sequence such that
$\limti{n} \norm{x_n} = \infty$. Then $\vecx$ is completely discrete,
and it is trivial that $\vecx$ has no joining system.

(c) Let $X$ be as in Part (c). Let $\vecx \subseteq X$ be
completely discrete. Since $X$ is an open subset of a Banach space,
we may assume that $\vecx$ is spaced.
Let $u \in X$. For every $n \in \bbN$ let $L_n \subseteq X$
be a rectifiable arc connecting $x_n$ with $u$ such that
$\rfs{lngth}(L_n) \leq d$.
Let $\gamma_n(t)$ be the parametrization of $L_n$
satisfying $\gamma_n(0) = u$, $\gamma_n(1) = x_n$ and
$\rfs{lngth}(\gamma_n([0,t])) = t \cdot\rfs{lngth}(L_n)$.

For every $\sigma \subseteq \bbN$ and $t \in [0,1]$ set
$A[\sigma,t] = \setm{\gamma_n(t)}{n \in \sigma}$,
and if $\sigma$ is infinite
define $t_{\sigma} = \inf(\setm{t}{A[\sigma,t] \mbox{ is spaced}})$.
There is an infinite $\sigma$ such that for every infinite
$\eta \subseteq \sigma$, $t_{\eta} = t_{\sigma}$.
It is easy to see that there is no infinite $\eta \subseteq \sigma$
such that $A[\eta,t_{\sigma}]$ is spaced.
So there is $\eta \subseteq \sigma$ such that $A[\eta,t_{\sigma}]$
is a Cauchy sequence.
Then $A[\eta,1]$ is a subsequence of $\vecx$ and
$\fortpl{A[\eta,1]}{\lim A[\eta,t_{\eta}]}
{\setm{\gamma_n([t_{\eta},1])}{n \in \eta}}{A[\eta,t_{\eta}]}$
is a joining system for $A[\eta,1]$.

(d) This part follows easily from (c).
\smallskip\hfill\myqed

In the next theorem, Part (a) is a special case of (b).
It seems wothwhile to state (a) separately,
because the class considered there
is more understandable than the class dealt with in (b).

\begin{theorem}\label{t4.33}\label{metr-bldr-t6.5}
\num{a} Let $K_{\srfs{BCX}}^{\calO}$
denote the class of all $X \in K_{\srfs{BNC}}^{\calO}$
such that $X$ is wide, BR.LC.AC
and BD.AC.
Suppose that $X,Y \in K_{\srfs{BCX}}^{\calO}$
and $\iso{\varphi}{\rfs{EXT}(X)}{\rfs{EXT}(Y)}$.
Then there is
$\tau \in \rfs{EXT}^{\pm}(X,Y)$ such that
$\tau$ induces $\varphi$.

Note that $K_{\srfs{BCX}}^{\calO}$ contains the class of all
bounded members of $K_{\srfs{BNC}}^{\calO}$
which are BR.LC.AC and BD.AC.
   \index{N@kobcx@@$K_{\srfs{BCX}}^{\calO}$}

\num{b} Let $K_{\srfs{NMX}}^{\calO}$
denote the class of all $X \in K_{\srfs{NRM}}^{\calO}$
such that $X$ is BR.LC.AC and JN.AC.
Let $X,Y \in K_{\srfs{NMX}}^{\calO}$.
Suppose that $\iso{\varphi}{\rfs{EXT}(X)}{\rfs{EXT}(Y)}$.
Then there is
$\tau \in \rfs{EXT}^{\pm}(X,Y)$ such that
$\tau$ induces $\varphi$.
   \index{N@konmx@@$K_{\srfs{NMX}}^{\calO}$}
\end{theorem}

The proof of Theorem \ref{metr-bldr-t6.5}
appears after Corollary~\ref{metr-bldr-c6.5}.
\smallskip

\noindent
{\bf Remark }
(a) By Proposition \ref{metr-bldr-p6.3}(c),
$K_{\srfs{BCX}}^{\calO} \subseteq K_{\srfs{NMX}}^{\calO}$.
So \ref{metr-bldr-t6.5}(b) is a special case of \ref{metr-bldr-t6.5}(a).

(b) Note that all members of $K_{\srfs{BCX}}^{\calO}$
which are subsets of a finite-dimensional normed space are bounded.
This is so, since for finite-dimensional spaces, wideness implies
boundedness.
Yet $K_{\srfs{BCX}}^{\calO}$ contains unbounded subsets of
infinite-dimensional Banach spaces.

(c) There is a regular open subset $X \subseteq \bbR^3$
such that $X \in K_{\srfs{BCX}}^{\calO}$
and $g^{\srfs{cl}} \nrestriction \rfs{bd}(X) = \rfs{Id}$
for every $g \in \rfs{EXT}(X)$.
This is maybe somewhat unexpected,
since it means that $\rfs{bd}(X)$ is recoverable from $\rfs{EXT}(X)$
even though every member of $\rfs{EXT}(X)$
is the identity on $\rfs{bd}(X)$.
See Example~\ref{e4.34}(d).
\smallskip

Recall that
$
\rfs{UC}_0(X) =
\setm{f \in \rfs{UC}(X)}{\rfs{Dom}(f^{\srfs{cl}}) = \rfs{cl}(X)
\mbox{ and } f^{\srfs{cl}} \rest\kern1pt \rfs{bd}(X) = \rfs{Id}}.
$

\begin{prop}\label{p4.29}\label{metr-bldr-p6.1}
Suppose that $X$ is BR.LC.AC,
and let $\tau \in H(X,Y)$ be such that
$(\rfs{UC}_0(X))^{\tau} \subseteq \rfs{EXT}(Y)$.
Let $x \in \rfs{bd}(X)$, $y \in \rfs{bd}(Y)$
and $\vecx \subseteq X$ be such that
$\lim \vecx = x$ and $\lim \tau(\vecx) = y$. 
Then $\tau \cup \sngltn{\pair{x}{y}}$ is continuous.
\end{prop}

\noindent
{\bf Proof }
Let $\vecu \subseteq X$ be such that $\lim \vecu = x$.
Suppose by contradiction that $\tau(\vecu)$ does not converge to $y$.
We may assume that $y$ is not a limit point of $\tau(\vecu)$.

We now repeat the construction appearing in the proof of Case 1 in
Theorem \ref{t4.10}.
Using the fact that $X$ is BR.LC.AC, we construct by induction
on $i \in \bbN$, $n_i \in \bbN$ and $L_i \subseteq X$ such that:
(i) $L_i$ is an arc connecting $x_{n_i}$ and $u_{n_i}$;
(ii) $\limti{i} \rfs{diam}(L_i) = 0$; and
(iii) for every $i \in \bbN$, $d(L_i,\bigcup_{j \neq i} L_j) > 0$.
For every $i \in \bbN$ let $U_i \subseteq X$ be an open set such that
$L_i \subseteq U_i$, $\limti{i} \rfs{diam}(U_i) = 0$,
and for every $i \neq j$, $d(U_i,U_j) > 0$.

Let $h_i \in \rfs{UC}(X)$ be such that
$\rfs{supp}(h_i) \subseteq U_{2i}$
and $h_i(x_{n_{2i}}) = u_{n_{2i}}$. By Proposition \ref{p4.9},
$h \eqdf \bcirc_{i \in \sboldbbN} h_i \in \rfs{UC}(X)$.
It is also obvious that $h \in \rfs{UC}_0(X)$.
However, $h^{\tau}$ is not exendible, since
$\tau(\vecx)$ is convergent,
whereas $h^{\tau}(\tau(\vecx))$ is not convergent. A contradiction.
\bigskip\hfill\myqed

Our next goal is to show that if
$(\rfs{EXT}(X))^{\tau} \subseteq \rfs{EXT}(Y)$, then for every
$y \in \rfs{bd}(Y)$ there is a sequence $\vecy$ converging to $y$
such that $\tau\inverse(\vecy)$ is a convergent sequence.
This holds automatically when $X$ is bounded and finite-dimensional,
but in that case extendibility is equivalent to uniform continuity,
and so Theorem \ref{t4.4} already answers our question.
In the general case we have to make an additional
arcwise connectedness assumption on $X$.

   \index{N@luc02@@$\rfs{LUC}_{01}(X) =
          \setm{h \in \rfs{LUC}(X)}{(\exists U)( U \mbox{ is $E$-open, }
          U \supseteq \rfs{bd}(X) \mbox{ and }
          \rfs{supp}(h) \subseteq X - U)}$}
For a metric space $E$ and $X \subseteq E$ define
{\thickmuskip=3.0mu \medmuskip=2.0mu \thinmuskip=2.0mu 
$$\rfs{LUC}_{01}(X) =
\setm{h \in \rfs{LUC}(X)}{\mbox{\,there is an $E$-open set }
U \supseteq \rfs{bd}(X)
\mbox{ such that } h \nrestriction (U \cap X) = \rfs{Id}}.
$$
}
\begin{lemma}\label{l4.32}\label{metr-bldr-l6.4}
Assume that $X$ is JN.AC, $\tau \in H(X,Y)$
and $(\rfs{LUC}_{01}(X))^{\tau} \subseteq \rfs{EXT}(Y)$,
and let $y \in \rfs{bd}(Y)$.

\num{a}
Suppose that $\vecx \subseteq X$ is completely discrete,
$\fortpl{\vecx}{x^*}{\setm{L_n}{n \in \bbN}}{\vecx\fprime}$
is a joining system for $\vecx$ and $\lim \tau(\vecx) = y$.
Then there is a sequence $\vecu \subseteq X$ such that
$\lim \vecu = x^*$ and $\lim \tau(\vecu) = y$.

\num{b}
There is a sequence $\vecu \subseteq X$
such that $\vecu$ converges to a member of $\rfs{bd}(X)$
and $\lim \tau(\vecu) = y$.
\end{lemma}

\noindent
{\bf Proof }
(a) Suppose that $\vecx$ is completely discrete,
$\fortpl{\vecx}{x^*}{\setm{L_n}{n \in \bbN}}{\vecx\fprime}$
is a joining system for $\vecx$,
and $\tau(\vecx)$ converges to $y$.
We may assume that $x^* \not\in \setm{x_n}{n \in \bbN}$.
Hence since $\vecx$ is completely discrete, $d \eqdf d(\vecx,x^*) > 0$.
Also assume that $L_n(0) = x_n$ and $L_n(1) = x'_n$.

{\bf Claim 1.}
For every $r > 0$ there is a sequence
$\vecu^{\fr} \subseteq B(x^*,r) \cap X$
such that $\tau(\vecu^{\fr})$ converges to $y$.
{\bf Proof }
Let $r \in (0,d)$. For every $n \in \bbN$ we define $v_n$.
If $n$ is even and $d(x'_n,x^*) \leq \dgfrac{r}{2}$, let
$t_n = \min \setm{t \in [0,1]}{d(L_n(t),x^*) = \dgfrac{r}{2}}$ and
$v_n = L_n(t_n)$.
Otherwise, let $v_n = x_n$.
Let $\vecv = \setm{v_n}{n \in \bbN}$.
Let $L'_n$ be the subarc of $L_n$ connecting $x_n$ with $v_n$. 
Clearly $L'_n \cap B(x^*,\dgfrac{r}{2}) = \emptyset$,
and hence by Definition \ref{metr-bldr-d6.2.1}(d)(5),
$\setm{L'_n}{n \in \bbN}$
is completely discrete.
It is easy to see that there is a completely discrete family of open
sets $\setm{U_n}{n \in \bbN}$ such that for every $n \in \bbN$,
$L'_n \subseteq U_n \subseteq \rfs{cl}(U_n) \subseteq X$.
Let $h_n \in \rfs{UC}(X)$ be such that $\rfs{supp}(h_n) \subseteq U_n$
and $h_n(x_n) = v_n$. It is easy to see that
$h \eqdf \bcirc \setm{h_n}{n \in \bbN} \in \rfs{LUC}_{01}(X)$.
Hence $h^{\tau} \in \rfs{EXT}(Y)$.

The facts that $\tau(\vecx)$ is convergent in $\rfs{cl}(Y)$
and that
$h^{\tau} \kern-1.0pt\in\kern-1.0pt \rfs{EXT}(Y)$
imply that
$h^{\tau}(\tau(\vecx))$ is also convergent in $\rfs{cl}(Y)$.
\hbox{Note that $h^{\tau}(\tau(\vecx)) = \tau(\vecv)$.}
So $\tau(\vecv)$ is convergent in $\rfs{cl}(Y)$.
Recall that for every $n \in \bbN$, $v_{2n+1} = x_{2n+1}$.
So \hbox{$\lim \tau(\vecv) = \lim \tau(\vecx) = y$.}
Let $N_r \in \bbN$ be such that for every $n > N_r$,
$d(x'_n,x^*) \leq \dgfrac{r}{2}$
and define
$\vecu^{\fr} = \setm{v_{2n}}{2n > N_r}$.
Then $\vecu^{\fr} \subseteq B(x^*,r) \cap X$ and hence
$\vecu^{\fr}$ is as required in Claim 1.
\smallskip

Let $r_n = \dgfrac{1}{n}$. For every $n \in \bbN$ let $k_n$ be such that
$d(y,\tau(u^{r_n}_{k_n})) < \dgfrac{1}{n}$.
Then $\vecu \eqdf \setm{u^{r_n}_{k_n}}{n \in \bbN}$
converges to $x^*$ and $\lim \tau(\vecu) = y$.
\smallskip

(b)
Suppose by contradiction that $y$ is a counter-example to the claim of
Part~(b).
Let $\vecy \subseteq Y$ be a
$\onetoonen$ sequence converging to $y$
and $\vecz = \tau^{-1}(\vecy)$.
If $\vecz$ has a convergent subsequence, then this subsequence
converges to a member of $\rfs{bd}(X)$, so $y$ is not a counter-example.
Hence $\vecz\kern1.3pt$ is completely discrete.

Since $X$ is JN.AC, there is a subsequence $\vecx$ of $\vecz$
such that $\vecx$ has a joining system
$\fortpl{\vecx}{x^*}{\setm{L_n}{n \in \bbN}}{\vecx\fprime}$.
By Part (a) there is a sequence $\vecu \subseteq X$ such that
$\lim \vecu = x^*$ and $\lim \tau(\vecu) = y$.
If $x^* \in X$,
then $y = \lim \tau(\vecu) = \tau(x^*) \in Y$, a contradiction.
So $x^* \in \rfs{bd}(X)$.
This means that $y$ is not a counter-example to (b).
A contradiction, so (b) is proved.
\smallskip\hfill\myqed

The fact $(\rfs{EXT}(X))^{\tau} \subseteq \rfs{EXT}(X)$
does not imply that $\tau \in \rfs{EXT}(X,Y)$.
To deduce that $\tau \in \rfs{EXT}(X,Y)$, we need to assume that
$(\rfs{EXT}(X))^{\tau} = \rfs{EXT}(X)$.
This is shown in Part~(a) of the next corollary.
In (b) we show that if $\rfs{EXT}(X)$ acts transitively
on $\rfs{bd}(X)$,
then the assumption $(\rfs{EXT}(X))^{\tau} \subseteq \rfs{EXT}(X)$
does suffice.

\begin{cor}\label{c4.33}\label{metr-bldr-c6.5}
\num{a} Suppose that $X$ is BR.LC.AC, and $Y$ is JN.AC.
Let $\tau \in H(X,Y)$ be such that
$(\dagger)$
$(\rfs{UC}_0(X))^{\tau} \subseteq \rfs{EXT}(Y)$
and
$(\dagger\dagger)$
$(\rfs{LUC}_{01}(Y))^{\tau^{-1}} \subseteq \rfs{EXT}(X)$.
Then $\tau \in \rfs{EXT}(X,Y)$.

\num{b} Suppose that $X$ is BR.LC.AC, $X$ is JN.AC,
and that the boundary of X has the following transitivity property:
$(*)$ for every $x,y \in \rfs{bd}(X)$ there is $h \in \rfs{EXT}(X)$
such that $h^{\srfs{cl}}(x) = y$.
Let $\tau \in H(X,Y)$ be such that
$(\rfs{EXT}(X))^{\tau} \subseteq \rfs{EXT}(Y)$.
Then $\tau \in \rfs{EXT}(X,Y)$.
\end{cor}

\noindent
{\bf Proof }
The two parts of the corollary will be proved
by combining Lemma \ref{l4.32}(b),
and Propositions~\ref{p4.29} and \ref{p4.28}(a).

(a) Let $x \in \rfs{bd}(X)$. By Lemma \ref{l4.32}(b)
applied to $\tau^{-1}$,
there is \hbox{$\vecx \subseteq X$} converging to $x$
such that $\tau(\vecx)$ converges to a point in $\rfs{bd}(Y)$.
Let
$y = \lim \tau(\vecx)$. By Proposition \ref{p4.29},
$\tau \cup \sngltn{\pair{x}{y}}$ is continuous.
So by Proposition \ref{p4.28}(a), $\tau$ is extendible.

(b) By Lemma \ref{l4.32}(b) applied to $\tau$,
there are $x_0 \in \rfs{bd}(X)$ and $\vecx \subseteq X$
converging to $x_0$
such that $\tau(\vecx)$ converges to a member of $\rfs{bd}(Y)$.
Let
$x \in \rfs{bd}(X)$. There is
$h \in \rfs{EXT}(X)$ such that $h(x_0) = x$.
Since $h^{\tau} \in \rfs{EXT}(Y)$, 
$h^{\tau}(\tau(\vecx))$ converges to a member of $\rfs{bd}(Y)$.
But $\tau(h(\vecx)) = h^{\tau}(\tau(\vecx))$.
It follows that for every $x \in \rfs{bd}(X)$
there is a sequence $\vecu$ converging to $x$
such that $\tau(\vecu)$ is convergent.
By Propositions~\ref{p4.29} and \ref{p4.28}(a),
$\tau \in \rfs{EXT}(X,Y)$.
\smallskip\hfill\myqed

\noindent
{\bf Proof of Theorem \ref{metr-bldr-t6.5} }
(a) This is a special case of Part (b),
because by Proposition~\ref{metr-bldr-p6.3}(d),
a BD.AC wide open subset of a Banach space is JN.AC.

(b) $\rfs{LIP}_{00}(X) \subseteq \rfs{EXT}(X)$
and $\rfs{LIP}_{00}(X) = \rfs{LIP}(X,\calS)$,
where $\calS$ is the set of all open BPD subsets of $X$.
The same holds for $Y$.
So by Theorem \ref{t2.4}(b), there is $\tau \in H(X,Y)$ such that
$\tau$ induces $\varphi$.
From the fact that $\rfs{UC}_0(X) \subseteq \rfs{EXT}(X)$
we conclude that  $(\rfs{UC}_0(X))^{\tau} \subseteq \rfs{EXT}(Y)$.
So \ref{metr-bldr-c6.5}(a) can be applied to $\tau$ and $\tau\inverse$.
We conclude that $\tau \in \rfs{EXT}^{\pm}(X,Y)$. This proves (b).
\smallskip\hfill\myqed

Part (a) of the next example is designed to show that the condition
$(\dagger)$ of \ref{metr-bldr-c6.5}(a) is needed.
Indeed, for $Y,X$ and $\tau\inverse$
of (a), $(\dagger\dagger)$ holds but the conclusion of
\ref{metr-bldr-c6.5}(a) does not.
Part (b) shows that assumption $(\dagger\dagger)$
in Corollary~\ref{metr-bldr-c6.5}(a) cannot be omitted.
The example is infinite-dimensional. Indeed, for finite-dimensional
normed spaces $(\dagger)$ {\it does} suffice. This follows from
Theorem~\ref{metr-bldr-t5.6} and Proposition~\ref{metr-bldr-p6.3}(e).
Part (c) shows that the transitivity assumption $(*)$
in Corollary~\ref{metr-bldr-c6.5}(b) is indeed needed.
Part (d) shows that there is $X \in K^{\calO}_{\srfs{BCX}}$ such that
$\rfs{EXT}(X)$ fixes $\rfs{bd}(X)$ pointwise.
The set $X$ is a regular open subset of $\bbR^3$,
therefore $\rfs{EXT}(X) = H(\rfs{cl}(X))$.

Let $\rfs{Cmp}(X)$ denote the set of connected components of a
topological space $X$.

   \index{N@cmp@@$\rfs{Cmp}(X)$. The set of connected components of $X$}

\begin{example}\label{e4.34}
\num{a} There are bounded regular open connected sets
$X$ and $Y$ in $\bbR^2$ and
$\tau \in H(X,Y)$ such that $X$ and $Y$ are BR.LC.AC,
$(\rfs{EXT}(X))^{\tau} \subseteq \rfs{EXT}(Y)$,
but $\tau^{-1} \not\in \rfs{EXT}(Y,X)$.
Note that by Proposition~\ref{metr-bldr-p5.43}(b),
$X$ and $Y$ are JN.AC.

\num{b} There are regular open bounded domains $X$ and $Y$ in an
infinite-dimensional Banach space and $\tau \in H(X,Y)$ such that
$X$ and $Y$ are BR.LC.AC and JN.AC,
\newline
$(\rfs{EXT}(X))^{\tau} \subseteq \rfs{EXT}(Y)$,
but $\tau \not\in \rfs{EXT}(X,Y)$.

\num{c} There are bounded domains $X$ and $Y$ in an
infinite-dimensional Banach space and $\tau \in H(X,Y)$ such that
$X$ and $Y$ are BR.LC.AC and JN.AC,
$\rfs{bd}(X)$ has two connected components, $\rfs{bd}(Y)$
is connected,
$\rfs{EXT}(X)$ and $\rfs{EXT}(Y)$
act very transitively on $\rfs{bd}(X)$
and $\rfs{bd}(Y)$ respectively,
$(\rfs{EXT}(X))^{\tau} \subseteq \rfs{EXT}(Y)$,
but $\tau \not\in \rfs{EXT}(X,Y)$.

\num{d} There is $X \in K^{\calO}_{\srfs{BCX}}$ such that
$X$ is a regular open bounded subset of $\bbR^3$,
and $g^{\srfs{cl}} \nrestriction \rfs{bd}(X) = \rfs{Id}$
for every $g \in \rfs{EXT}(X)$.
\end{example}

\noindent{\bf Proof }
(a) Let $X' \subseteq \bbR^2$ be the open square whose vertices are
$(0,0)$, $(1,0)$, $(0,1)$ and $(1,1)$,
and $Y' \subseteq \bbR^2$ be the open triangle whose vertices are
$(0,0)$, $(0,1)$ and $(1,1)$. Let $\tau' \in H(X',Y')$ be defined by
$\tau'((x,y)) = (xy,y)$. Let $A = [(0,0),(1,0)]$.

Clearly, $\tau' \in \rfs{EXT}(X',Y')$,
$(\tau')^{cl} \nrestriction (\rfs{cl}(X) - A) \in
H(\rfs{cl}(X') - A\,,\,\rfs{cl}(Y') - \sngltn{(0,0)})$
and
$(\tau')^{cl}(A) = \sngltn{(0,0)}$.
Also, if $g \in \rfs{EXT}(X',X')$ and $g^{cl}(A) = A$,
then $g^{\tau'} \in \rfs{EXT}(Y')$.

For $n > 1$ and $1 \leq k < n$ let
$x_{n,k} = (\dgfrac{k}{2^n},\dgfrac{1}{2^n})$,
$B_{n,k} = \rfs{cl}(B(x_{n,k},\dgfrac{1}{8^n}))$
and $\calB = \setm{B_{n,k}}{n > 1,\  1 \leq k < n}$.
Note that $\calB$ is a pairwise disjoint family of closed balls
contained in $X'$ and $\rfs{cl}(\bigcup \calB) - \bigcup \calB = A$.
Let $X = X' - \bigcup \calB$,
$Y = \tau'(X)$ and $\tau = \tau' \nrestriction X$.
Clearly, for every $g \in \rfs{EXT}(X)$, $g^{cl}(A) = A$. It follows
that $X$, $Y$ and $\tau$ are as required.
{\thickmuskip=4.3mu \medmuskip=3mu \thinmuskip=2mu 
Note also that for every $x,y \in A - \dbltn{(0,0)}{(1,0)}$
there is  $g \in \rfs{EXT}(X)$ such that $g(x) = y$.
\kern-8pt\smallskip
}

(b) Let $E$ be the Hilbert space $\ell_2$, $Y'$ be the open cylinder
defined by
$$\hbox{$
Y' =
\setm{(x_0,x_1,\ldots)}
{\abs{x_0} < 3 \mbox{ and } \sum\limits_{i = 1}^{\infty} x_i^2 < 9}
$}$$
and $X' = Y' - \overB^E(0,1)$.
Let $\iso{\tau_1}{X'}{Y' - \sngltn{0}}$ be such that
$\tau_1 \nrestriction (Y' - B^E(0,2)) = \rfs{Id}$.
Let $\iso{\tau_2}{Y' - \sngltn{0}}{Y'}$ be such that
$\tau_2 \nrestriction (Y' - B^E(0,2)) = \rfs{Id}$
and $\tau' = \tau_2 \scirc \tau_1$.
The existence of $\tau_2$ follows from the facts that a point in
$\bbR^{\sboldbbN}$ is a strongly negligible set, and that
$\ell_2 \cong \bbR^{\sboldbbN}$.
See \cite{BP} Chapter IV Definition 5.1 and
Chapter V Proposition~2.2(c) and Theorem 6.4.

Note that $\tau'$ cannot be continued to a continuous function defined
on $S(0,1)$.
Hence $\tau' \not\in \rfs{EXT}(X',Y')$.
It is trivial that $\rfs{bd}(Y')$ is homeomorphic to a sphere,
and that $\rfs{bd}(X')$ has two components:
$\rfs{bd}(Y')$ and $S(0,1)$.
It can be easily checked that for every $h \in \rfs{EXT}(X')$: if
$h^{cl}(S(0,1)) = S(0,1)$, then $h^{\tau'} \in \rfs{EXT}(Y')$.
However, there is $h \in \rfs{EXT}(X')$ such that
$h^{cl}(S(0,1)) = \rfs{bd}(Y')$.
This implies that $(\rfs{EXT}(X'))^{\tau} \not\subseteq \rfs{EXT}(Y')$,
contrary to what is required in this example.

For a pairwise disjoint family $\calC$ of subsets a topological space
$Z$ define
$$
\rfs{acc}^Z(\calC) =
\setm{z \in Z}{\mbox{for every }U \in \rfs{Nbr}^Z(z),
\ \setm{C \in \calC}{U \cap C \neq \emptyset}
\mbox{ is in infinite}}.
$$
To define $X$
we construct a pairwise disjoint family $\calF$ of closed sets
such that\break
(i) $\bigcup \calF \subseteq Y' - \overB(0,2)$
and (ii) $\rfs{acc}(\calF) \subseteq \rfs{bd}(Y')  \cup \bigcup \calF$.
We shall then define $X$, $Y$ and $\tau$ to be respectively
$X' - \bigcup \calF$,
$\tau'(X)$ and $\tau' \nrestriction X$.
It follows from (ii) that $X$ is open,
and the construction of $\calF$ will ensure that
$S(0,1)$ is the unique connected component of
$\rfs{bd}(X)$ which is clopen in $\rfs{bd}(X)$ and which is also
strongly connected, (a notion to be defined later).
It will thus follow that for every
$h \in \rfs{EXT}(X)$, $h^{cl}(S(0,1)) = S(0,1)$,
and this in turn implies that 
$(\rfs{EXT}(X))^{\tau} \subseteq \rfs{EXT}(Y)$.

Let $\setm{e_i}{i \in \bbN}$ be the standard basis of $\ell_2$,
denote by $T$ the set of finite sequences of natural numbers,
let $\fnn{f}{T}{\bbN - \sngltn{0}}$ be a $\onetoonen$ function,
and for $\eta \in T$ define $d_{\eta} = e_{f(\eta)}$.
Let $\sLambda$ denote the empty sequence
and $T^* = T - \sngltn{\sLambda}$.
The relation ``$\nu$ is a proper initial segment of $\eta$''
is denoted by $\eta \iseg \nu$.
Suppose that $\eta = \nu \concat \ontpl{i}$,
$\zeta = \nu \concat \ontpl{j}$ and  $i \neq j$.
In that case we say that $\nu = \rfs{pred}(\eta)$,
$\eta \in \rfs{Suc}(\nu)$ and $\zeta \in \rfs{Brthr}(\eta)$.

Let $<^T$ be the relation on $T$ defined by
$\nu <^T \eta$ if either $\eta \iseg \nu$
or there is $n \in \rfs{Dom}(\nu) \cap \rfs{Dom}(\eta)$
such that $\nu \nrestriction \bbN^{<n} = \eta \nrestriction \bbN^{<n}$
and $\nu(n) < \eta(n)$.
It is easy to check that $<{^T}$ is a dense linear ordering
with maximum $\sLambda$ and with no minimum.
Denote by $T_n$ the set of all $\eta \in T$ such that
$\rfs{Dom}(\eta) = \bbN^{< m}$ for some $m \leq n$.
Then $T_n$ is well-ordered by~$<^T$.

We define a line segment $L_{\eta}$
for every $\eta \in T^*$.
If $\eta = \nu \concat \ontpl{m}$, then $L_{\eta}$ has the form
$[d_{\nu} + a_{\eta} \mcdot e_0,d_{\eta} + a_{\eta} \mcdot e_0]$,
where $2 < a_{\eta} < 3$.
So for $L_{\eta}$ to be defined we need to define $a_{\eta}$.
We define $a_{\eta}$ by induction.
Let $\setm{\eta_n}{n \in \bbN}$ be a $\onetoonen$ enumeration of $T$
such that for every $n \in \bbN$ and $\nu \iseg \eta_n$
there is $m < n$ such that $\eta_m = \nu$.
Define
$S_n =
\setm{\eta_m \sconcat \ontpl{i}}{m < n \mbox{ and } i \in \bbN}$.
We define by induction on $n$ the set
$\setm{a_{\nu}}{\nu \in S_n}$.
So at stage $n$ we need to define the set
$\setm{a_{\eta_n \sconcat \ontpl{i}}}{i \in \bbN}$.
Since $\setm{\nu}{\nu \iseg \eta_n} \subseteq \setm {\eta_m}{m < n}$
for every $n$, it follows that $\eta_0 = \sLambda$.
Let $\sngltn{a_{\ontpl{i}}}_{i \in \bbN}$
be a strictly increasing sequence converging to $3$
such that $a_{\ontpl{0}} = \dgfrac{5}{2}$.
So
$$
L_{\ontpl{i}} =
[d_{\tLambda} + a_{\ontpl{i}} \mcdot e_0,
d_{\ontpl{i}} + a_{\ontpl{i}} \mcdot e_0].
$$
Let $n > 0$ and suppose that $a_{\nu}$
has been defined for every $\nu \in S_n$.
Let $\barzero = \pair{0}{\ldots}$ denote the infinite sequence of $0$'s.
It is convenient to define $a_{\barzero} = 2$.
We assume by induction that
\begin{itemize}
\addtolength{\parskip}{-11pt}
\addtolength{\itemsep}{06pt}
\item[(1)]
$2 < a_{\nu} < 3$ for every $\nu \in S_n$,
\item[(2)]
{\thickmuskip=3mu \medmuskip=2mu \thinmuskip=1mu 
$\setm{a_{\eta_m \sconcat \ontpl{i}}}{i \in \bbN}$
is a strictly increasing sequence converging to $a_{\eta_m}$
for every $0 < m < n$,}
\item[(3)]
if $\nu,\rho \in S_n$ and $\nu <^T \rho$,
then
$a_{\nu} < a_{\rho}$.
\vspace{-05.7pt}
\end{itemize}
Note that for $n = 1$ the induction hypotheses hold.
Clearly, $S_n \subseteq T_{n + 1}$,
so $\setm{a_{\nu}}{\nu \in S_n}$ is well-ordered.
Obviously, $\eta_n \in S_n$. If $\eta_n = \trpl{0}{\ldots}{0}$,
then $\eta_n = \rfs{min}(S_n)$.
In this case set $\rho_n = \barzero$.
Otherwise, write $\eta_n$ as
$\nu \concat \ontpl{k} \concat \trpl{0}{\ldots}{0}$,
where $k > 0$, and the sequence of $0$'s at the end of $\eta_n$
may be the empty sequence.
Define $\rho_n = \nu \concat \ontpl{k - 1}$.
It is easy to check that in this case
$\rho_n$ is the predecessor of $\eta_n$ in $S_n$.
Choose $\setm{a_{\eta_n \sconcat \ontpl{i}}}{i \in \bbN}$
to be a strictly increasing sequence converging to $a_{\eta_n}$
such that
$a_{\eta_n \sconcat \ontpl{0}} =
\dgfrac{(a_{\rho_n} + a_{\eta_n})}{2}$.
It is left to the reader to verify that the induction hypotheses hold.

Let $\calL = \setm{L_{\eta}}{\eta \in T^*}$,
set $a_{\tLambda} = 3$,
for $\eta \in T$ define $c_{\eta} = d_{\eta} + a_{\eta} e_0$
and let\break
$C = \setm{c_{\eta}}{\eta \in T}$.
Note that $c_{\tLambda} \in \rfs{bd}(Y')$.
For $\eta = \nu \concat \ontpl{i} \in T^*$
define $b_{\eta} = d_{\nu} + a_{\eta} e_0$.
So $L_{\eta} = [b_{\eta},c_{\eta}]$.

We first establish some facts about the distance between the members of
$\calL$.

{\bf Claim 1.} If $\nu \neq \rfs{pred}(\eta)$,
$\eta \neq \rfs{pred}(\nu)$
and $\rfs{pred}(\nu) \neq \rfs{pred}(\eta)$,
then $d(L_{\nu},L_{\eta}) > 1$.

{\bf Proof }
$L_{\nu}$ and $L_{\eta}$ can be written as
$L_{\nu} = a_{\nu} e_0 + [b,c]$
and
$L_{\eta} = a_{\eta} e_0 + [d,e]$,
where $b,c,d,e \in \setm{e_i}{i \in \bbN^{\geq 1}}$ and
$\dbltn{b}{c} \cap \dbltn{d}{e} = \emptyset$.
So $(d(L_{\nu},L_{\eta}))^2 =
(a_{\nu} - a_{\eta})^2 + 4 \mcdot \frac{1}{4} > 1$.

{\bf Claim 2.}
Suppose that $\nu = \rfs{pred}(\eta)$ or $\eta = \rfs{pred}(\nu)$
or $\nu \in \rfs{Brthr}(\eta)$
and write $L_{\nu} = a_{\nu} e_0 + [b,c]$
and $L_{\eta} = a_{\eta} e_0 + [b,d]$,
where $b,c,d \in \setm{e_i}{i \in \bbN^{\geq 1}}$.
Let $x \in L_{\nu}$ and write
$x = a_{\nu} e_0 + b + e$. Then
$d(x,L_{\eta}) > \frac{\sqrt{3}}{2} \norm{e}$.

{\bf Proof }
Clearly, $e$ can be written as $e = t(c - b)$
and so
$$
d(x,L_{\eta})^2 = (a_{\eta} - a_{\nu})^2 + d(b + e,[b,d])^2 >
d(b + e,[b,d])^2 = d(t(c - b),[0,d - b])^2.
$$
Also, 
$$\hbox{$
d(t(c - b),[0,d - b]) \geq d(t(c - b),\setm{s(d - b)}{s \in \bbR}) =
\norm{t(c - b)} \mcdot \sin\frac{\pi}{3} =
\frac{\sqrt{3}}{2} \norm{e}.
$}$$
So $d(x,L_{\eta}) > \frac{\sqrt{3}}{2} \norm{e}$.
This proves Claim 2.
\smallskip

If we define $X_0 = X' - \bigcup \calL$, $Y_0 = \tau'(X_0)$
and $\tau_0 = \tau' \nrestriction X_0$, then all the requirements of
Part (b) are fulfilled except that $X_0$ is not regular open.
To achieve that $X$ be regular open,
we replace every $L_{\eta}$ by a set $F_{\eta}$
such that $F_{\eta} = \rfs{cl}(\rfs{int}(F_{\eta}))$.
This will ensure that $X$ is regular open.
The verification of the following trivial fact is left to the reader.

{\bf Claim 3.} $C$ is $\sqrt{2}$-spaced.

Let $\eta,\nu \in T^*$.
For distinct $x,y \in \ell_2$ define
$H_{x,y} = (\setm{t(y - x)}{t \in \bbR})^{\perp}$.
Let $\theta$ be such that $\tan\theta = \dgfrac{1}{8}$
and define the ``closed double cone'' of $x,y$ to be
$$
\rfs{dcone}(x,y) =
\setm{z \in [x,y] + H_{x,y}}
{d(z,[x,y]) \leq d(z,\dbltn{x}{y}) \mcdot \sin\theta}.
$$
Note that $\rfs{dcone}(x,y)$ is the union of two cones with vertices
$x,y$. The common base of the two cones is
$B(\dgfrac{(x + y)}{2},r) \cap (\dgfrac{(x + y)}{2} + H_{x,y})$,
where $r = \half \norm{y - x} \mcdot \tan\theta$
and the opening angle of the cones is $\theta$.
The verification of the following fact is omitted.

{\bf Claim 4.}
There is $K > 1$ such that for every distinct $x,y,u,v \in \ell_2$
and $\varepsilon > 0$:
if $u,v \not\in \rfs{dcone}(x,y)$
and
$d(u,\rfs{dcone}(x,y)), d(v,\rfs{dcone}(x,y)) \leq
\varepsilon$,
then there is a rectifiable arc $J$ connecting $u,v$
such that
$J \subseteq
\setm{z}{d(z,\rfs{dcone}(x,y)) \leq \varepsilon} - \rfs{dcone}(x,y)$,
$d(J,\dbltn{x}{y}) = d(\dbltn{u}{v},\dbltn{x}{y})$
and $\rfs{lngth}(J) \leq K \norm{u - v}$.

Note that in order to prove Claim~4 it suffices to consider the affine
subspace of $\ell_2$ generated by $x,y,u,v$. So the proof can be
carried out in a 3-dimensional Euclidean space.

Define $F_{\eta} = \rfs{dcone}(b_{\eta},c_{\eta})$,
$\calF = \setm{F_{\eta}}{\eta \in T^*}$, $\whatF = \bigcup \calF$,
$X = X' - \whatF$, $Y = Y' - \whatF$
and $\tau = \tau' \nrestriction X$.
Clearly, $\tau \in H(X,Y)$.
Since $\tau'$ cannot be continued to a continuous function defined
on $S(0,1)$,
$\tau$ too cannot be continued to a continuous function defined
on $S(0,1)$.
Hence $\tau \not\in \rfs{EXT}(X,Y)$.
The next claim contains the central fact about $\calF$.

{\bf Claim 5.} Let $\eta \in T^*$ and $r > 0$.
Then $d(F_{\eta} - B(c_{\eta},r), \whatF - F_{\eta}) > 0$.
{\bf Proof } Let $\eta = \nu \concat \ontpl{i}$.
If $i > 0$ define
$\delta_{\eta} =
\rfs{min}(
a_{\nu \sconcat \ontpl{i + 1}} - a_{\nu \sconcat \ontpl{i}},
a_{\nu \sconcat \ontpl{i}} - a_{\nu \sconcat \ontpl{i - 1}})$
and if $i = 0$ define
$\delta_{\eta} =
a_{\nu \sconcat \ontpl{i + 1}} - a_{\nu \sconcat \ontpl{i}}$.
Let
$\varepsilon_{\eta,r} =
\min(\frac{3}{4},\frac{3r}{4},\frac{\delta_{\eta}}{3})$.
Let $\zeta \in T^* - \sngltn{\eta}$.
We show that
$d(F_{\eta} - B(c_{\eta},r),F_{\zeta}) \geq \varepsilon_{\eta,r}$.
If
$\zeta \not\in
\rfs{Brthr}(\eta) \cup \rfs{Suc}(\eta) \cup \sngltn{\rfs{pred}(\eta)}$,
then by Claim~1, $d(L_{\eta},L_{\zeta}) > 1$.
So
$d(F_{\eta},F_{\zeta}) >
1 - 2 \mcdot \eighth \frac{\sqrt{2}}{2} > \dgfrac{3}{4}$.

Suppose that $\zeta \in \rfs{Suc}(\eta)$.
Recall that $c_{\eta} = a_{\eta} e_0 + d_{\eta}$.
Let $x \in F_{\eta} - B(c_{\eta},r)$ and let
$y$ be the nearest to $x$ in $L_{\eta}$.
Then $y = a_{\eta} e_0 + d_{\eta} + e$,
where $e$ has the form $e = s(d_{\nu} - d_{\eta})$.
Since $\norm{x - y} \leq \dgfrac{\norm{e}}{8}$
and $\norm{x - (a_{\eta} e_0 + d_{\eta})} \geq r$,
we have that and $\norm{e} \geq \dgfrac{8r}{9}$.
Take a point $z \in F_{\zeta}$,
let $w$ be the nearest point to $z$ in $L_{\zeta}$
and suppose that $\norm{w - (a_{\zeta} e_0 + d_{\eta})} = t$. 
Then $\norm{y - w} > \sqrt{\norm{e}^2 + t^2}$
and hence
$\norm{y - z} > \sqrt{\norm{e}^2 + t^2} - \dgfrac{t}{8}$.
The minimal value of the function
$g(t) = \sqrt{\norm{e}^2 + t^2} - \dgfrac{t}{8}$
is $\geq \norm{e} - \dgfrac{\norm{e}}{56}$.
This implies that
$d(y,F_{\zeta}) \geq \norm{e} - \dgfrac{\norm{e}}{56}$.
Since $\norm{x - y} \leq \dgfrac{\norm{e}}{8}$,
it follows that
$d(x,F_{\zeta}) \geq
\norm{e} - \dgfrac{\norm{e}}{56} - \dgfrac{\norm{e}}{8} =
\dgfrac{6 \norm{e}}{7}$.
Hence
$d(x,F_{\zeta}) \geq \frac{6}{7} \mcdot \frac{8}{9} r \geq
\dgfrac{3r}{4}$.

Assume that
$\zeta \in \rfs{Brthr}(\eta) \cup \sngltn{\rfs{pred}(\eta)}$.
Define $f = a_{\eta} e_0 + d_{\nu}$.
Let $x \in F_{\eta}$
and suppose first that
$\norm{x - f} \leq \dgfrac{\delta_{\eta}}{2}$.
If $w \in L_{\zeta}$ and $\norm{w - (a_{\zeta} e_0 + d_{\nu})} = t$,
then $d(f,w) \geq \sqrt{\delta_{\eta}^2 + t^2}$,
So the distance between $f$ and a general point in $F_{\zeta}$ is
$\geq \sqrt{\delta_{\eta}^2 + t^2} - \dgfrac{t}{8}$.
So $d(f,F_{\zeta}) \geq \delta_{\eta} - \dgfrac{\delta_{\eta}}{56}$
and hence 
$d(x,F_{\zeta}) \geq
\delta_{\eta} - \dgfrac{\delta_{\eta}}{56} - \dgfrac{\delta_{\eta}}{2}
> \dgfrac{\delta_{\eta}}{3}$.

Suppose that $x \in F_{\eta}$ and
$\norm{x - f} \geq \dgfrac{\delta_{\eta}}{2}$.
Let $y$ be the nearest point to $x$ in $L_{\eta}$
and $\delta = \norm{y - f}$.
Then $d(y,F_{\zeta}) \geq \delta - \dgfrac{\delta}{56}$
and hence
$d(x,F_{\zeta}) \geq
\delta - \dgfrac{\delta}{56} - \dgfrac{\delta}{8} =
\dgfrac{6\delta}{7}$.
Also, $\delta \geq \frac{8}{9} \mcdot \frac{\delta_{\eta}}{2}$.
So
$d(x,F_{\zeta}) \geq \frac{6}{7} \mcdot \frac{8}{9} \mcdot
\frac{\delta_{\eta}}{2} > \frac{\delta_{\eta}}{3}$.
The proof of Claim~5 is complete.
\smallskip

{\bf Claim 6.}
(i) $\calF$ is a pairwise disjoint family and $\rfs{acc}^E(\calF) = C$.

(ii) Let $\eta \in T$, $\setm{F_n}{n \in \bbN} \subseteq \calF$
be a $\onetoonen$ sequence, $x_n \in F_n$ and
$\limti{n} x_n = c_{\eta}$.
Then
$\setm{F_n}{n \in \bbN} -
\setm{F_{\eta \sconcat \ontpl{i}}}{i \in \bbN}$
is finite.

{\bf Proof }
By Claim 5, $(F_{\eta} - \sngltn{c_{\eta}}) \cap F_{\zeta} = \emptyset$
for every distinct $\eta,\zeta \in T^*$.
Since $c_{\eta} \neq c_{\zeta}$ for any $\eta \neq \zeta$,
it follows that $\calF$ is pairwise disjoint.

We show that $C \subseteq \rfs{acc}(\calF)$.
Recall that $C = \setm{c_{\eta}}{\eta \in T}$,
where $c_{\eta} = d_{\eta} + a_{\eta} e_0$
and $a_{\tLambda} = 3$.
We start with $c_{\tLambda}$.
By the construction,
$a_{\ontpl{n}} \mcdot e_0 + d_{\tLambda} \in L_{\ontpl{n}} \subseteq
F_{\ontpl{n}}$
and 
$c_{\tLambda} = 3e_0 + d_{\tLambda} =
\limti{n} a_{\ontpl{n}} \mcdot e_0 + d_{\tLambda}$.
So $c_{\tLambda} \in \rfs{acc}(\calF)$.
Suppose now that $\eta \neq \sLambda$.
Then
$a_{\eta \sconcat \ontpl{n}} \mcdot e_0 + d_{\eta} \in
L_{\eta \sconcat \ontpl{n}} \subseteq F_{\eta \sconcat \ontpl{n}}$
and 
$c_{\eta} = a_{\eta} e_0 + d_{\eta} =
\limti{n} a_{\eta \sconcat \ontpl{n}} \mcdot e_0 + d_{\eta}$.
So $c_{\eta} \in \rfs{acc}(\calF)$.
We have shown that $C \subseteq \rfs{acc}(\calF)$.
\smallskip

Let $\setm{\nu_i}{i \in \bbN} \subseteq T^*$ be a $\onetoonen$ sequence,
$x_i \in F_{\nu_i}$, and suppose that $\sngltn{x_i}_{i \in \bbN}$
is convergent.
Let $x = \limti{i} x_i$.
We shall show that for some $\eta \in T$, $x = c_{\eta}$
and that
$\setm{F_{\nu_i}}{i \in \bbN} -
\setm{F_{\eta \sconcat \ontpl{i}}}{i \in \bbN}$
is finite.
This will imply both that $\rfs{acc}(\calF) \subseteq C$ and (ii).
We color the unordered pairs of $\bbN$ in three colors.
The pair $\dbltn{i}{j}$ has Color 1 if $\nu_i \in \rfs{Brthr}(\nu_j)$,
and $\dbltn{i}{j}$ has Color 2 if $\nu_i = \rfs{pred}(\nu_j)$
or $\nu_j = \rfs{pred}(\nu_i)$.
The remaining unordered pairs have Color 3.
By Ramsey Theorem we may assume the $\bbN$ is monochromatic.
Color 2 has no infinite monochromatic sets, and
if $\bbN$ has Color 3,
then by the first paragraph in proof of Claim~5
the sequence $\sngltn{x_i}_{i \in \bbN}$ is $\frac{3}{4}$-spaced.
It follows that for some $\eta \in T$,
$\setm{\nu_i}{i \in \bbN} \subseteq \rfs{Suc}(\eta)$.

Let $y_i$ be the nearest point to $x_i$ in $L_{\nu_i}$,
and write $y_i = a_{\nu_i} \mcdot e_0 + d_{\eta} + f_i$,
where $f_i = t_i(d_{\eta \sconcat\ontpl{n_i}} - d_{\eta})$
for some $t_i \in [0,1]$.
We may assume that 
$\sngltn{f_i}_{i \in \bbN}$ is convergent and let $f = \limti{i} f_i$.
Suppose by way of contradiction that $f \neq 0$.
Let $n$ be such that for every $i,j \geq n$,
$\norm{x_i - x_j} < \varepsilon$,
where $\varepsilon$ is to be chosen later,
and $\frac{4}{5} \norm{f} < \norm{f_i} < 2 \norm{f}$.
Let $i,j \geq n$ be distinct.
Then $\norm{x_i - y_i} \leq \dgfrac{\norm{f_i}}{8} \leq 
\dgfrac{\norm{f}}{4}$
and $\norm{x_j - y_j} \leq \dgfrac{\norm{f}}{8} \leq
\dgfrac{\norm{f}}{4}$.
So
$$\hbox{$
\norm{y_i - y_j} \leq \norm{y_i - x_i} + \norm{x_i - x_j} +
\norm{x_j - y_j} \leq \dgfrac{\norm{f_i}}{8} +
\dgfrac{\norm{f_j}}{8} + \varepsilon <
\dgfrac{\norm{f}}{2} + \varepsilon.
$}$$
On the other hand, by Claim 2,
$$\hbox{$
\norm{y_i - y_j} \geq d(y_i,L_{\nu_j}) \geq 
\frac{\sqrt{3}}{2} \norm{f_i} \geq
\frac{2\sqrt{3}}{5} \norm{f}.
$}$$
If $\varepsilon$ is sufficiently small, then the last two inequalities
are contradictory. So $f = 0$.
Now, $\norm{x_i - y_i} \leq \dgfrac{\norm{f_i}}{8}$.
So $\limti{i} \norm{x_i - y_i} = 0$ and hence
$$
\limti{i} x_i = \limti{i} y_i =
\limti{i} a_{\nu_i} e_0 + d_{\eta} + f_i =
\limti{i} a_{\nu_i} e_0 + d_{\eta} = a_{\eta} e_0 + d_{\eta} =
c_{\eta} \in C.
$$
We have proved that $\rfs{acc}(\calF) \subseteq C$.
We have also shown that
if $\setm{F_n}{n \in \bbN} \subseteq \calF$
is a $\onetoonen$ sequence, $x_n \in F_n$ and
$\limti{n} x_n = c_{\eta}$,
then
$\setm{F_n}{n \in \bbN} \cap
\setm{F_{\eta \sconcat \ontpl{i}}}{i \in \bbN}$
is infinite. Obviously, this implies (ii).
This completes the proof of Claim~6.
\smallskip

Denote $\whatF \cup \sngltn{c_{\tLambda}}$ by
$\wtildeF$.
Since every member of $\calF$ is closed
and
$\rfs{acc}(\calF) = C \subseteq \wtildeF$,
it follows that
$\wtildeF$ is closed.
Recall that $c_{\tLambda} \in \rfs{bd}(Y')$ and hence
$c_{\tLambda} \not\in X'$.
It follows that $X = X' - \wtildeF$,
so $X$ is open.
Clearly, $F = \rfs{cl}(\rfs{int}(F))$ for every $F \in \calF$.
So $\wtildeF = \rfs{cl}(\rfs{int}(\wtildeF))$.
This implies that $E - \wtildeF$ is regular open,
and hence $X = X' \cap (E - \wtildeF)$ is regular open.
An identical argument shows that $Y$ is regular open in~$E$.
\smallskip

{\bf Claim 7.} Let $K$ be the constant mentioned in Claim~4.
Then for every $x,y \in Y$ there is a rectifiable arc $J \subseteq Y$
connecting $x$ and $y$ such that $\rfs{lngth}(J) \leq 2K \norm{x - y}$.
Similarly, let $K_1 = \max(2K,\pi)$.
Then for every $x,y \in X$
there is a rectifiable arc $J \subseteq X$
connecting $x$ and $y$ such that $\rfs{lngth}(J) \leq K_1 \norm{x - y}$.

{\bf Proof } Let $x,y \in Y$.
By Claim~3, $C$ is spaced,
so for every $\varepsilon > 0$ there is
$z \in B(y,\varepsilon)$ such that $[x,z] \cap C = \emptyset$.
Choose such a $z$ for a small $\varepsilon$ which will be determined
later.
Since $Y$ is open, we may choose $z$ such that $[z,y] \subseteq Y$,
and since $Y'$ is convex, $[x,z] \subseteq Y'$.
Since $[x,z] \cap C = \emptyset$ and $\rfs{acc}(\calF) = C$,
$\calF_0 \eqdf \setm{F \in \calF}{F \cap [x,z] \neq \emptyset}$
is finite. The fact that $C$ is spaced implies that
$r \eqdf d([x,z],C) > 0$.
Let $\calF_0 = \fsetn{F_0}{F_{n - 1}}$,
$F_i = F_{\eta_i}$, $b_i = b_{\eta_i}$,
$c_i = c_{\eta_i}$,
$F_i \cap [x,z] = [x_{i,0},x_{i,1}]$ and
$$\hbox{$
\delta_i =
\half\min(d(F_i - B(c_i,\dgfrac{r}{2}), \whatF - F_i),r,
\delta^{Y'}(\bigcup \calF_0)).
$}$$
By Claim~5, $\delta_i > 0$.
Let $\hatx_{i,j} \in [x,z]$ be such that
$\norm{\hatx_{i,j} - x_{i,j}} \leq \delta_i$
and $[\hatx_{i,j},x_{i,j}) \cap F_i = \emptyset$.
By Claim~4, there is a rectifiable arc $J_i$ connecting
$\hatx_{i,0}$ and $\hatx_{i,1}$ such that
$\rfs{lngth}(J_i) \leq K \norm{\hatx_{i,0} - \hatx_{i,1}}$,
$J_i \subseteq \setm{z \in \ell_2}{d(z,F_i) \leq \delta_i} - F_i$
and
$d(J_i,\dbltn{b_i}{c_i}) =
d(\dbltn{\hatx_{i,0}}{\hatx_{i,1}},\dbltn{b_i}{c_i})$.
Since $d(\dbltn{\hatx_{i,0}}{\hatx_{i,1}},\dbltn{b_i}{c_i}) \geq r$,
it follows that $d(J_i,c_i) \geq r$.
Let $u \in J_i$ and $v$ be the nearest point to $u$ in $F_i$.
Then $\norm{c_i - v} \geq \norm{c_i - u} - \norm{u - v} \geq
\dgfrac{r}{2}$.
So $v \in F_i - B(c_i,\dgfrac{r}{2})$,
and hence $d(v,\whatF - F_i) \geq 2 \delta_i$.
From the fact that $\norm{u - v} \leq \delta_i$ it follows that
$u \not\in \whatF - F_i$, so $J_i \cap \whatF = \emptyset$.
Also, since for every $u \in J_i$,
$d(u,\bigcup \calF_0) < \delta^{Y'}(\bigcup \calF_0)$,
we have that $J_i \subseteq Y'$.
Let
$J' = [x,z] \cup \bigcup_{i < n} J_i -
\bigcup_{i < n} [\hatx_{i,0},\hatx_{i,1}]$
and $J = J' \cup [z,y]$.
It is easily seen that $J'$ and $J$ are rectifiable arcs,
and it follows that $J \subseteq Y' - \whatF = Y$.
From the fact that
$\rfs{lngth}(J_i) \leq K \norm{\hatx_{i,0} - \hatx_{i,1}}$,
it follows that $\rfs{lngth}(J') \leq K \norm{z - x}$.
Recall that $\norm{y - z} < \varepsilon$.
So if $\varepsilon$ is sufficiently small,
then $\rfs{lngth}(J) < 2K \norm{y - x}$.

The proof of the analogous fact for $X$ is almost identical.
We have proved Claim 7.
\smallskip

We now show that $X$ and $Y$ are BR.LC.AC and JN.AC.
Claim 7 implies that $Y$ is UD.AC and BD.AC.
It follows directly from the definitions that
if $F$ is any metric space, $Z \subseteq F$ and $Z$ is UD.AC,
then  $Z$ is BR.LC.AC with respect to $F$.
Hence $Y$ is BR.LC.AC with respect to $\ell_2$.
The bounded arcwise connectedness of $Y$ and Proposition~\ref{p4.31}(c)
imply that $Y$ is JN.AC.
The same arguments apply to $X$, hence $X$ too is BR.LC.AC and JN.AC.

Our next goal is to show
$(*)$ $h(S(0,1)) = S(0,1)$ for every $h \in \rfs{EXT}(X)$.
It may very well be true that
$(\dagger)$ $S(0,1)$ is the only clopen component
of $\rfs{bd}(X)$.
This would imply~$(*)$, but we do not know to prove this.
So instead we prove
$(\dagger\dagger)$
$S(0,1)$ is the only clopen component
of $\rfs{bd}(X)$ which is strongly connected in $\rfs{bd}(X)$.
This also implies $(*)$.

Let $Z$ be a connected space. We say that $Z$ is
{\it strongly connected}
if for every $z \in Z$ and $U \in \rfs{Nbr}(z)$,
there is $V \in \rfs{Nbr}(z)$ such that $V \subseteq U$
and $Z - V$ is connected.
Clearly, $S(0,1)$ is strongly connected.

For $\eta \in T^*$ let $S_{\eta} = \rfs{bd}^{\ell_2}(F_{\eta})$.
It is easy to see that
$\rfs{bd}(X) =
S(0,1) \cup S(0,3) \cup \bigcup_{\eta \in T^*} S_{\eta}$.
Obviously, $S(0,1)$ is a component of $\rfs{bd}(X)$,
and $S(0,1)$ is clopen in $\rfs{bd}(X)$.
Let $\calK$ denote the set of components of $\rfs{bd}(X)$
which are clopen in $\rfs{bd}(X)$ and which are
different from $S(0,1)$.
Let $\eta \in T$ and $T' \subseteq T$.
We say that $T'$ is $\eta$-large
if $\eta \in T' \subseteq T^{\srisegeq \eta}$,
and for every $\nu \in T'$,
$\setm{i}{\nu \concat \ontpl{i} \not\in T'}$ is finite.
Define $S_{\tLambda} = S(0,3)$ and for $T' \subseteq T$
set $S_{T'} = \bigcup_{\nu \in T'} S_{\nu}$.

{\bf Claim 8.}
For every $K \in \calK$ there are a finite set $\sigma \subseteq T$
and a family $\setm{T_{\nu}}{\nu \in \sigma}$ such that
$T_{\nu}$ is $\nu$-large for every $\nu \in \sigma$,
and $K = \bigcup_{\nu \in \sigma} S_{T_{\nu}}$.

{\bf Proof }
Note that $S_{\eta}$ is connected for every $\eta \in T$.
Hence for every
$K \in \calK$ and $\eta \in T$, either $S_{\eta} \subseteq K$
or $S_{\eta} \cap K = \emptyset$.
Also,
for every $\eta \in T$ and an infinite $\sigma \subseteq \bbN$,
$S_{\eta} \cap
\rfs{acc}(\setm{S_{\eta \sconcat \ontpl{i}}}{i \in \sigma})
\neq \emptyset$.
This implies that
$(\dagger)$
if $K \in \calK$ and $S_{\eta} \cap K \neq \emptyset$,
then $S_{\eta} \subseteq K$
and $\setm{i}{S_{\eta \sconcat \ontpl{i}} \not\subseteq K}$
is finite.
The fact that the members of $\calK$ are closed implies that
$(\dagger\dagger)$
if $K \in \calK$
and $\setm{i}{S_{\eta \sconcat \ontpl{i}} \subseteq K}$ is infinite,
then $S_{\eta} \subseteq K$.
Facts $(\dagger)$ and $(\dagger\dagger)$ imply that Claim~8 is true.
\smallskip

Let $K \in \calK$
and suppose that $\sigma \subseteq T$
and $\setm{T_{\nu}}{\nu \in \sigma}$ are as assured by Claim~8.
So there are $\eta \in T^*$ and an infinite $T' \subseteq T$
such that $S_{\eta} \subseteq K = S_{T'}$.
By Claim~5, $d(S_{\eta} - B(c_{\eta},r),\whatF - F_{\eta}) > 0$
for every $r > 0$.
Since $S_{\eta}$ and $S_{\tLambda}$ are closed and disjoint,
it follows that $d(S_{\eta},S_{\tLambda}) > 0$,
and from the facts that $K \subseteq \whatF \cup S_{\tLambda}$
and $S_{\eta} \subseteq F_{\eta}$ we conclude that
$d(S_{\eta} - B(c_{\eta},r),K - S_{\eta}) > 0$.
So $S_{\eta} - B(c_{\eta},r)$ is clopen in $K$.
This implies that $K$ is not strongly connected.
We have shown that $S(0,1)$ is the unique clopen strongly connected
component of $\rfs{bd}(X)$.
Hence $h(S(0,1)) = S(0,1)$ for every $h \in \rfs{EXT}(X)$.
It follows that $(\rfs{EXT}(X))^{\tau} \subseteq \rfs{EXT}(Y)$.
This completes the proof of (b).
\smallskip

(c) Let $S \subseteq \ell_2$ be a two-dimensional sphere with radius 1
and center at 0. Let\break
$X = B(0,3) - S$ and $Y = B(0,3)$.
Then there is $\tau \in H(X,Y)$ such that
$\tau \nrestriction (B(0,3) - \overB(0,2)) = \rfs{Id}$.
It is trivial that $X$ and $Y$ are BR.LC.AC and JN.AC,
and it is easy to see that
$(\rfs{EXT}(X))^{\tau} \subseteq \rfs{EXT}(Y)$
and $\tau \not\in H(X,Y)$.

(d) We construct a set $X$ with the following properties:
\begin{itemize}
\addtolength{\parskip}{-11pt}
\addtolength{\itemsep}{06pt}
\item[(1)] 
$X$ is a regular open bounded subset of $\bbR^3$,
\item[(2)] 
there is $K > 1$
such that for every $x,y \in X$
there is a rectifiable arc $L \subseteq X$
such that $\rfs{lngth}(L) \leq K \norm {x - y}$,
\item[(3)] 
for every $g \in \rfs{EXT}(X)$,
$g^{\srfs{cl}} \nrestriction \rfs{bd}(X) = \rfs{Id}$.
\vspace{-05.7pt}
\end{itemize}
It is easy to verify that if $X$ satisfies (1)-(3), then it fulfills the
requirements of the example.

We turn to the construction of $X$.
Let $\whatR_n$ be the $n$-fold solid torus
and $\whatT_n$ denote its boundary.
A subset $A \subseteq \bbR^3$ is {\it $K$-bypassable}, if for every
$x,y \in \bbR^3 - A$ there is a rectifiable arc $L \subseteq \bbR^3 - A$
connecting $x$ and $y$ such that $\rfs{lngth}(L) \leq K \norm{x - y}$
and $d(z,A) \leq d(x,A), d(y,A)$ for every $z \in L$.
Obviously, there is $K > 1$ such that for every $n$
there is a $K$-bypassable $F \subseteq \bbR^3$
such that $F \cong \whatR_n$.
Let $D$ be a countable dense subset of $B(0,1)$,
$E$ be a countable dense subset of $S(0,1)$
and $\setm{\dbltn{a_n}{b_n}}{n \in \bbN}$ be a list
of all $2$-element subsets of $D$
and all singletons from $D \cup E$.
Also assume that $a_0 = b_0 \in E$.
We define by induction a finite family of open sets $\calU_n$
and a finite family of closed sets $\calF_n$
such that for any distinct $A \in \calU_n \cup \calF_n$
and $F \in \calF_n$,
$\rfs{cl}(A) \subseteq B(0,1)$,
$\rfs{cl}(A) \cap F = \emptyset$
and $F$ is $K$-bypassable.
Let $\calU_0 = \calF_0 = \emptyset$.
Suppose that $\calU_n$ and $\calF_n$ have been defined.

{\bf Case 1 } $a_n \neq b_n$.
If $\dbltn{a_n}{b_n} \cap \bigcup \calF_n \neq \emptyset$
define $\calU_{n + 1} = \calU_n$ and $\calF_{n + 1} = \calF_n$.
Suppose otherwise.
Define $\calF_{n + 1} = \calF_n$.
Since $\calF_n$ is a finite pairwise disjoint family of closed
$K$-bypassable sets there is a rectifiable arc
$L_n \subseteq B(0,1) - \bigcup \calF_n$
connecting $a_n$ and $b_n$
such that $\rfs{lngth}(L_n) \leq K \norm{a_n - b_n}$.
Let $r = d(L_n,S(0,1) \cup \bigcup \calF_n)$
and $\calU_{n + 1} = \calU_n \cup \sngltn{B(L_n,\dgfrac{r}{2})}$.

{\bf Case 2 } $a_n = b_n$.
If  $a_n \in D$ let $c_n \in \bigcup_{F \in \calF_n} \rfs{bd}(F))$
be such that
$\norm{c_n - a_n} = d(a_n,\bigcup_{F \in \calF_n} \rfs{bd}(F))$
and $H_n \in \calF_n$ be such that $c_n \in \rfs{bd}(H_n)$.
If $a_n \in E$ let $c_n = a_n$ and $H_n = S(0,1)$.
Let
$F_n \subseteq B(c_n,\frac{1}{n + 1}) \cap
(B(0,1) - \bigcup \calF_n  - \bigcup_{U \in \calU_n} \rfs{cl}(U))$
be such that $F_n \cong \whatR_n$ and $F_n$ is $K$-bypassable.
Define $\calF_{n + 1} = \calF_n \cup \sngltn{F_n}$.
Let \hbox{$r_n = d(H_n,S(0,1) \cup \bigcup \calF_{n + 1} - H_n)$}
and
{\thickmuskip=3.5mu \medmuskip=2.5mu \thinmuskip=1.5mu 
$U_{n,0} = B(0,1) \cap
(B(H_n,\dgfrac{r_n}{2}) - \rfs{cl}(B(H_n,\dgfrac{r_n}{4})))$.
Let $x_n \in B(0,1) \cap (B(c_n,\dgfrac{r_n}{2}) - H_n)$,\break
}
$s_n \in (0,\dgfrac{r_n}{2})$ be such that
$U_{n,1} \eqdf B(x_n,s_n)$ is disjoint from $H_n$
and $\calU_{n + 1} = \calU_n \cup \dbltn{U_{n,0}}{U_{n,1}}$.
This concludes the inductive construction.

Let $X = B(0,1) - \rfs{cl}(\bigcup_{n \in \bbN} \calF_n)$.
Since every two members of $D \cap X$
lie in the same member of $\rfs{Cmp}(X)$
and $D \cap X$ is dense in $X$,
it follows that $X$ is connected.

Set $A = \setm{n}{a_n = b_n}$, for every $n \in A$
let $\iso{f_n}{\whatR_n}{F_n}$ and $T_n = f_n(\whatT_n)$
and define $T = S(0,1) \cup \bigcup_{n \in A} T_n$.
The verification of the following facts is left to the reader.
\begin{itemize}
\addtolength{\parskip}{-11pt}
\addtolength{\itemsep}{06pt}
\item[(1)] 
$\rfs{bd}(X) = \rfs{cl}(T)$
and $T \subseteq \rfs{cl}(\rfs{int}(\bbR^3 - X))$.
\item[(2)] 
For every $n \in A$, $T_n \in \rfs{Cmp}(\rfs{bd}(X))$,
and $S(0,1) \in \rfs{Cmp}(\rfs{bd}(X))$.
\item[(3)] 
For every
$C \in \rfs{Cmp}(\rfs{bd}(X)) - \setm{T_n}{n \in A} - \sngltn{S(0,1)}$,
$\bbR^3 - C$ is connected.
\vspace{-05.7pt}
\end{itemize}
Fact (1) implies that $X$ is regular open.
It follows from (3) and Alexander's Duality Theorem,
that for every
$C \in \rfs{Cmp}(\rfs{bd}(X)) - \setm{T_n}{n \in A} - \sngltn{S(0,1)}$
and $n \in \bbN$, $C \not\cong \whatT_n$.
Let $x \in T$.
Then there is a sequence $\setm{k_n}{n \in \bbN} \subseteq A$
such that $\limti{k_n} T_n = x$.
Hence $x$ has the following property:
\begin{itemize}
\addtolength{\parskip}{-11pt}
\addtolength{\itemsep}{06pt}
\item[] 
There is a sequence $\setm{C_n}{n \in \bbN}$ of members of
$\rfs{Cmp}(\rfs{bd}(X))$ such that $C_n \cong T_{k_n}$ and
$\limti{n} C_n = x$.
\vspace{-05.7pt}
\end{itemize}
However, if $y \in \rfs{bd}(X) - \sngltn{x}$,
then $y$ does not have this property.
Since $\rfs{bd}(X)$ is invariant under $\rfs{EXT}(X)$,
it follows that $g(x) = x$ for every $g \in \rfs{EXT}(X)$.
That is, $g \nrestriction T = \rfs{Id}$ for every $g \in \rfs{EXT}(X)$.
Since $T$ is dense in $\rfs{bd}(X)$,
it follows that
$g \nrestriction \rfs{bd}(X) = \rfs{Id}$
for every $g \in \rfs{EXT}(X)$.
\smallskip\hfill\myqed

{\bf Remark }
Recall that in Corollary \ref{c4.33}(b) it was assumed that
for every $x,y \in \rfs{bd}(X)$ there is $h \in \rfs{EXT}(X)$
such that $h^{\srfs{cl}}(x) = y$.
In Part (c) of the above example
$\rfs{bd}(X)$ has two connected components $K_0,K_1$,
both $K_0$ and $K_1$ are not a singleton,
and for every $i = 0,1$ and $x,y \in K_i$
there is $h \in \rfs{EXT}(X)$ such that $h^{\srfs{cl}}(x) = y$.
The space $Y$ in the above example
has the property that $\rfs{bd}(Y)$ is connceted,
$\rfs{bd}(Y)$ is not a singleton,
and for every $x,y \in \rfs{bd}(Y)$,
there is $h \in \rfs{EXT}(X)$ such that $h^{\srfs{cl}}(x) = y$.
These transitivity properties of $\rfs{bd}(X)$ and $\rfs{bd}(Y)$,
though quite strong, do not imply the conclusion of \ref{c4.33}(b).
\smallskip

In Theorem \ref{metr-bldr-t6.5} it was shown that if
$\iso{\varphi}{\rfs{EXT}(X)}{\rfs{EXT}(Y)}$, then $\varphi$ is
induced by some $\tau \in \rfs{EXT}^{\pm}(X,Y)$.
But Theorem \ref{metr-bldr-t6.5} applies only to sets $X$
with finitely many connected components.
To see this let $X$ be BR.LC.AC and JN.AC
as was assumed in \ref{metr-bldr-t6.5},
and suppose by contradiction that
$X$ has infinitely many connected components.
Let $\vecz$ be a sequence of members of $X$ which
lie in distinct components of $X$.
Let $\fortpl{\vecx}{x^*}{\setm{L_n}{n \in \bbN}}{\vecx\fprime}$
be a joining system for some subsequence $\vecx$ of $\vecz$.
Then $\vecx\fprime$ is a convergent sequence,
but each member of $\rfs{Rng}(\vecx\fprime)$ lies in a different
component of $X$. This contradicts the fact that $X$ is BR.LC.AC.
So $X$ has only finitely many connected components.

Our next goal is to extend \ref{metr-bldr-t6.5} to sets $X$ that
may have infinitely many connected components.
We have four test cases $X$ for which $\rfs{EXT}(X)$ seems to be
sufficiently well behaved to imply a reconstruction theorem for
$\rfs{EXT}(X)$,
but which are not covered by Theorem~\ref{metr-bldr-t6.5}.
The first example which is defined below,
has infinitely many components.
The three others appear in Example~\ref{metr-bldr-e6.12},
and they are connected.

\begin{example}\label{metr-bldr-e6.7}
\begin{rm}
Let $E$ be a Banach space.
We define
$$\hbox{$
R_1^E = \bigcup_{n \in \bbN}
(B^E(0,1 - \frac{1}{2n + 3}) - \overB^E(0,1 - \frac{1}{2n + 2}))
$}$$
The set $R_1^E$ is the union of a sequence of pairwise disjoint
open rings converging to $S^E(0,1)$.
\end{rm}
\end{example}

We shall prove a reconstruction theorem
for a class which contains~$R_1^E$.
The definition of this class is rather technical, but it contains
quite complicated sets. This class will be denoted by
$K_{\srfs{BX}}^{\calO}$.
For simplicity, we consider only subsets of a Banach spaces
and not subsets of general normed spaces.
Hence only \ref{metr-bldr-t6.5}(a) is extended.
That is, $K_{\srfs{BCX}}^{\calO} \subseteq K_{\srfs{BX}}^{\calO}$.

\begin{defn}\label{metr-bldr-d6.8}
\begin{rm}
(a) Recall that $\rfs{Cmp}(X)$ denotes the set of connected components
of a topological space $X$.
For $x,y \in X$, $x \simeq^X y$ denotes that $x$ and $y$
lie in the same connected component of $X$.
The notation $\vecx \simeq^X \vecy$ means that $x_n \simeq^X y_n$
for every $n \in \bbN$.

   \index{N@AAAA@@
          {\thickmuskip=2mu \medmuskip=1mu \thinmuskip=1mu 
          $x \simeq^X y$. This means that $x$ and $y$ lie in the
          same connected component of $X$}}
   \index{N@AAAA@@$\vecx \simeq^X \vecy$.
          This means that for every $n$, $x_n \simeq^X y_n$}

(b) Let $X$ be a metric space.
We say that
{\it $X$ is boundedly component-wise arcwise connected (BD.CW.AC)},
if for every bounded set $A \subseteq X$ there is $d = d_A$
such that for every $x,y \in A$: if $x \simeq^X y$,
then there is a rectifiable arc $L \subseteq X$ connecting $x$ and $y$
such that $\rfs{lngth}(L) \leq d$.

   \index{boundedly component-wise arcwise connected, $X$ is}
   \index{bdcwac@@BD.CW.AC. Abbreviation of
          boundedly component-wise arcwise connected}

(c) Let $X \in K_{\srfs{NRM}}^{\calO}$ and $x \in \rfs{bd}(X)$.
We say that
$X$ is {\it component-wise locally arcwise connected} at $x$,
if for every $\varepsilon > 0$ there is $\delta > 0$ such that
for every
$y,z \in B(x,\delta) \cap X$: if $y \simeq^X z$,
then there is an arc $L \subseteq B(x,\varepsilon) \cap X$
connecting $y$ and $z$.
We say that
$X$ is {\it component-wise locally arcwise connected at its boundary
(BR.CW.LC.AC)},
if $X$ is component-wise locally arcwise connected
at every $x \in \rfs{bd}(X)$.

   \index{component-wise locally arcwise connected at $x$}
   \index{component-wise locally arcwise connected at the boundary}
   \index{brcwlac@@BR.CW.LC.AC. Abbreviation of $X$ is
          component-wise locally arcwise connected at
          \newline\indent
          its boundary}

(d) Let $X \in K_{\srfs{NRM}}^{\calO}$.
Call $X$ a {\it component-wise wide space},
if for every $r > 0$,\break
$\bigcup \setm{C \in \rfs{Cmp}(X)}{C \cap B(0,r) \neq \emptyset}$
is wide.

   \index{component-wise wide}

(e) Let $X \subseteq E$. A point $x \in \rfs{bd}(X)$
is called a {\it multiple boundary point} of $X$, if for every
$C \in \rfs{Cmp}(X)$, $x \in \rfs{bd}(X - C)$,
and $x$ is a {\it double boundary point} of $X$,
if there are distinct $C_1,C_2 \in \rfs{Cmp}(X)$
such that $x \in \rfs{bd}(C_1) \cap \rfs{bd}(C_2)$.
   \index{multiple boundary point}
   \index{double boundary point}

(f) A subspace $X \subseteq E$
is {\it locally movable at its multiple boundary},
if for every $\vecx \subseteq X$ which converges in $E$ to a
multiple boundary point and $U \in \rfs{Nbr}^{\srfs{cl}(X)}(\lim \vecx)$
there is a subsequence $\vecx\fprime$ of $\vecx$ and
$g \in \rfs{EXT}(X)$ such that:
$g(\vecx\fprime) \simeq^X \vecx\fprime$,
$g^{\srfs{cl}}(\lim \vecx) \neq \lim \vecx$
and $\rfs{supp}(g) \subseteq U$.

   \index{locally movable at the multiple boundary}

(g) Let $K_{\srfs{BX}}^{\calO}$ be the class of all
$X \in K_{\srfs{BNC}}^{\calO}$ such that:
\begin{itemize}
\addtolength{\parskip}{-11pt}
\addtolength{\itemsep}{06pt}
\item[(1)] 
$X$ is component-wise wide,
BR.CW.LC.AC and BD.CW.AC,
\item[(2)] 
$X$ is locally movable at its multiple boundary.
\vspace{-05.7pt}
\end{itemize}
   \index{N@kobx@@$K_{\srfs{BX}}^{\calO}$}
\end{rm}
\end{defn}

\begin{prop}\label{metr-bldr-p6.11}
\num{a} Let $R_1^E$ be as defined in Example \ref{metr-bldr-e6.7}.
Then $R_1^E \in K_{\srfs{BX}}^{\calO}$.

\num{b} $K_{\srfs{BCX}}^{\calO} \subseteq K_{\srfs{BX}}^{\calO}$.
\end{prop}

\noindent
{\bf Proof } The proofs of both parts are trivial. Anyway, we indicate
the proof of (b).
Suppose that $X \in K_{\srfs{BCX}}^{\calO}$.
It is easily seen that the multiple boundary of $X$ is empty,
hence $X$ is locally movable at its multiple boundary.
The fact that $X$ is wide implies
that it is component-wise wide.
Similarly, since $X$ is BR.LC.AC and BD.AC,
it is BR.CW.LC.AC and BD.CW.AC.
So $X \in K_{\srfs{BX}}^{\calO}$.
\hfill\myqed

\begin{prop}\label{metr-bldr-p6.9}
\num{a} Let $X \in K_{\srfs{BX}}^{\calO}$.
Then for every $C \in \rfs{Cmp}(X)$,
$C$ is BR.LC.AC and JN.AC.

\num{b} Let $X,Y \in K_{\srfs{BX}}^{\calO}$ and $\tau \in H(X,Y)$
be such that $(\rfs{EXT}(X))^{\tau} = \rfs{EXT}(Y)$.
Let $C \in \rfs{Cmp}(X)$, $D = \tau(C)$
and $\eta = \tau \nrestriction C$.
Then $D \in \rfs{Cmp}(Y)$ and $\eta \in \rfs{EXT}^{\pm}(C,D)$.
\end{prop}

\noindent
{\bf Proof }
(a) The fact that $X$ is component-wise wide implies that $C$ is wide.
The fact that $X$ is BD.CW.AC implies that $C$ is BD.AC.
So by Proposition \ref{metr-bldr-p6.3}(d), $C$ is JN.AC.

Let $x \in \rfs{bd}(C)$. The fact that $X$ is component-wise locally
arcwise connected at $x$ implies that $C$ is locally arcwise connected
at $x$. So $C$ is BR.LC.AC.

(b) It is trivial that $C$ is an open subset of $E$
and that $D \in \rfs{Cmp}(Y)$.
So by Part (a), $C$ is JN.AC and BR.LC.AC, and the same holds for $D$.
We wish to apply Corollary \ref{metr-bldr-c6.5}(a) to $\eta$,
so we need to check that
$(\rfs{UC}_0(C))^{\eta} \subseteq \rfs{EXT}(D)$
and that $(\rfs{LUC}_{01}(D))^{\eta^{-1}} \subseteq \rfs{EXT}(C)$.
Let $g \in \rfs{UC}_0(C)$.
Set $h = g \cup \rfs{Id} \nrestriction (X - C)$.
Then $h \in \rfs{UC}_0(X) \subseteq \rfs{EXT}(X)$.
So $h^{\tau} \in \rfs{EXT}(Y)$.
Hence $g^{\eta} = h^{\tau} \nrestriction D \in \rfs{EXT}(D)$.
A similar argument shows that
$(\rfs{LUC}_{01}(D))^{\eta^{-1}} \subseteq \rfs{EXT}(C)$.
By Corollary \ref{metr-bldr-c6.5}(a),
$\eta \in \rfs{EXT}(C,D)$.
The same argument can be applied to $\eta\inverse$.
Hence $\eta \in \rfs{EXT}^{\pm}(C,D)$.
\hfill\myqed

\begin{theorem}\label{metr-bldr-t6.10}
Let $X,Y \in K_{\srfs{BX}}^{\calO}$ and
$\iso{\varphi}{\rfs{EXT}(X)}{\rfs{EXT}(Y)}$.
Then there is
$\tau \in \rfs{EXT}^{\pm}(X,Y)$
such that $\tau$ induces~$\varphi$.
\end{theorem}

\noindent
{\bf Proof } By Theorem \ref{t2.4}(b), there is
$\tau \in H(X,Y)$ such that $\tau$ induces $\varphi$.

{\bf Claim 1.} Let $\vecx,\vecu \subseteq X$.
Suppose that $\vecx,\vecu,\tau(\vecx),\tau(\vecu)$
are convergent sequences and $\lim \vecx = \lim \vecu \in \rfs{bd}(X)$.
Then $\lim \tau(\vecx) = \lim \tau(\vecu)$.
{\bf Proof }
Let $x = \lim \vecx$,
$y = \lim \tau(\vecx)$ and $v = \lim \tau(\vecu)$,
and suppose by contradiction that $y \neq v$.
Clearly, $y,v \in \rfs{bd}(Y)$.
Assume first that either $y$ or $v$ is a multiple boundary point of
$Y$, and assume without loss of generality that $y$ is such a point.
Since $Y$ is locally movable at its multiple boundary,
there are $h \in \rfs{EXT}(Y)$ and a subsequence $\vecy\fprime$ of 
$\tau(\vecx)$ such that
$h^{\srfs{cl}}(y) \neq y$, $h(\vecy\fprime) \simeq^Y \vecy\fprime$
and for some $W \in \rfs{Nbr}^{\srfs{cl}(Y)}(v)$,
$h \nrestriction (W \cap Y) = \rfs{Id}$.
By removing an initial segment of $\tau(\vecu)$ we may assume that
$\tau(\vecu) \subseteq W$.
So $h$, $\vecy\fprime$ and $W$ satisfy
\begin{itemize}
\addtolength{\parskip}{-09pt}
\addtolength{\itemsep}{06pt}
\item[$(*)$] 
$h \in \rfs{EXT}(Y)$, $\vecy\fprime$ is a subsequence of $\tau(\vecx)$,
$W \in \rfs{Nbr}^{\srfs{cl}(Y)}(v)$, $h^{\srfs{cl}}(y) \neq y$,
$\tau(\vecu) \subseteq W$
and $h^{\srfs{cl}} \nrestriction W = \rfs{Id}$.
\vspace{-05.7pt}
\end{itemize}

Now assume that $y,v$ are not multiple boundary points of $Y$.
Then there are $C_1,C_2 \in \rfs{Cmp}(X)$
such that all but finitely members of $\vecx$ belong to $C_1$,
and all but finitely members of $\vecu$ belong to $C_2$.
From Proposition~\ref{metr-bldr-p6.9}(b)
and the fact that $\lim \tau(\vecx) \neq \lim \tau(\vecu)$
it follows that, $C_1 \neq C_2$.
So $x$ is a double boundary point of $X$.
Let $D_1 = \tau(C_1)$
and set $\whatD = Y - D_1$.
Then by \ref{metr-bldr-p6.9}(b), $D_1 \in \rfs{Cmp}(Y)$,
and since $y$ is not a multiple boundary point of $Y$,
it follows that $y \in \rfs{bd}(D_1) - \rfs{cl}(\whatD)$.
Let $V \in \rfs{Nbr}^F(y)$ be such that 
$\rfs{cl}(V) \cap \rfs{cl}(\whatD) = \emptyset$,
and let $U \in \rfs{Nbr}^E(x)$ be such that
$\tau(U \cap C_1) \subseteq V$.
Since $X$ is locally movable at its multiple boundary,
there is $k \in \rfs{EXT}(X)$ and a subsequence $\vecz\fprime$
of $\vecx$ such that $k^{\srfs{cl}}(x) \neq x$,
$\rfs{supp}(k) \subseteq U$
and $k(\vecz\fprime) \simeq \vecz\fprime$.
Let
$h = (k \nrestriction C_1)^{\tau} \cup 
\rfs{Id} \nrestriction (Y - D_1)$.
Then $h \nrestriction D_1 \in \rfs{EXT}(D_1)$.
Also,
$$
\rfs{supp}(h) = \rfs{supp}(h \nrestriction D_1) =
\tau(\rfs{supp}(k \nrestriction C_1)) \subseteq \tau(U \cap C_1)
\subseteq V.
$$
So
$\rfs{supp}((h \nrestriction D_1)^{\srfs{cl}}) \subseteq
\rfs{cl}(V)$.
From the fact that
$\rfs{cl}(V) \cap \rfs{cl}(\whatD) = \emptyset$,
it follows that $h \in \rfs{EXT}(Y)$.
Let $\vecy\fprime = \tau(\vecz\fprime)$.
Then
$h^{\srfs{cl}}(y) \neq y$ and $h(\vecy\fprime) \simeq^Y \vecy\fprime$.
Clearly, $v \in \rfs{cl}(\tau(C_2))$
and $\tau(C_2) \subseteq \whatD$.
So $v \in \rfs{cl}(\whatD)$,
and hence for some $W \in \rfs{Nbr}^{\srfs{cl}(Y)}(v)$,
$h \nrestriction (W \cap Y) = \rfs{Id}$.
By removing an initial segment of $\tau(\vecu)$ we may assume that
$\tau(\vecu) \subseteq W$.
It follows that $h$, $\vecy\fprime$ and $W$ satisfy $(*)$.
So in both cases: when $\dbltn{u}{v}$ contains
a multiple boundary point, and when it does not,
we have found $h,\vecy\fprime$ and $W$ satisfying $(*)$.

Let $g = h^{\tau\inverse}$
and $\vecx\fprime = \tau\inverse(\vecy\fprime)$.
So $g \in \rfs{EXT}(X)$ and $g \nrestriction \vecu = \rfs{Id}$.
Since $\vecu \cup \vecx\fprime$ converges to $x$
and $g \in \rfs{EXT}(X)$, $\lim g(\vecx\fprime) = x$.
Since $h(\vecy\fprime) \simeq^Y \vecy\fprime$,
it follows that $g(\vecx\fprime) \simeq^X \vecx\fprime$.
Since $X$ is BR.CW.AC,
there is $\setm{f_k}{k \in \bbN} \subseteq \rfs{UC}(X)$
and subsequences $\sngltn{n_k}_{k \in \bbN}$
and $\sngltn{m_k}_{k \in \bbN}$
such that: (i) for every $k$, $f_k(x_{n_k}') = g(x_{n_k}')$,
$\rfs{cl}(\rfs{supp}(f_k)) \subseteq X$
and $f_k \nrestriction \setm{x_{m_k}'}{k \in \bbN} = \rfs{Id}$,
(ii) $\limti{k} \rfs{diam}(\rfs{supp}(f_k)) = 0$,
(iii) for any $\ell \neq k$, 
$\rfs{supp}(f_{\ell}) \cap \rfs{supp}(f_k) = \emptyset$.
Let $f = \bcirc_{k \in \bbN} f_k$.
So $f \in \rfs{UC}_0(X) \subseteq \rfs{EXT}(X)$,
and hence $f^{\tau}$ must belong to $\rfs{EXT}(Y)$.
Let us see that this does not happen.
Recall that $\lim \vecy\fprime = y$.
However, $\lim_k f^{\tau}(y'_{n_k}) = \lim_k h(y'_{n_k}) =
h(y) \neq y$,
and on the other hand,
$\lim_k f^{\tau}(y'_{m_k}) = \lim_k y'_{m_k} = y$.
So $\vecy\fprime$ is convergent, but $f^{\tau}(\vecy\fprime)$ is not,
and hence $f^{\tau} \not\in \rfs{EXT}(Y)$.
A contradiction, so Claim 1 is proved.\kern-5pt
\smallskip

{\bf Claim 2.} Let $\vecx \subseteq X$ be a convergent sequence in $E$.
Then there is a subsequence $\vecx\fprime$ of $\vecx$
such that $\tau(\vecx\fprime)$ is convergent in $F$.
{\bf Proof } Let $x = \lim \vecx$.
We may assume that $x \in \rfs{bd}(X)$.
If for some $C \in \rfs{Cmp}(X)$,
$\setm{n}{x_n \in C}$ is infinite,
then by Proposition \ref{metr-bldr-p6.9}(b), there is a subsequence
as required in the claim.

Hence we may assume that for every $m \neq n$,
$x_m \not\simeq^X x_n$, and so $x$ is a multiple boundary point.
For every $n$ let $y_n = \tau(x_n)$,
and $C_n$ and $D_n$ be such that $x_n \in C_n \in \rfs{Cmp}(X)$
and $y_n \in D_n \in \rfs{Cmp}(Y)$.

Suppose by contradiction that $\setm{D_n}{n \in \bbN}$
is completely discrete.
Let $\vecu \in \prod_{n \in \bbN} C_n$.
Define $\vecv = \tau(\vecu)$.
There is $k \in \rfs{EXT}(Y)$
such that for every $n$, $k(y_{2n}) = v_{2n}$
and
$k(y_{2n + 1}) = y_{2n + 1}$.
Let $g = k^{\tau\inverse}$.
Then $g \in \rfs{EXT}(X)$.
Since $\vecx$ is convergent, $g(\vecx)$ is convergent.
For every~$n$, $g(x_{2n}) = u_{2n}$
and $g(x_{2n + 1}) = x_{2n + 1}$.
So $\limti{n} u_{2n} = \limti{n} x_{2n + 1} = x$.
This implies that $\limti{n} C_{2n} = x$.
Hence for every $f \in \rfs{EXT}(X)$:
if $\setm{n \in \bbN}{f(x_{2n}) \simeq^X x_{2n}}$ is infinite,
then $f(x) = x$.
Clearly, $x$ is a multiple boundary point.
So the above fact is in contradiction with the the fact
that $X$ is locally movable at its multiple boundary.
It follows that $\setm{D_n}{n \in \bbN}$ is not completely discrete.
By choosing a subsequence of $\vecx$ we may assume that there is
$\vecv \in \prod_{i \in \bbN} D_n$
such that $\vecv$ is convergent in $F$. Let $v = \lim \vecv$.

Suppose by way of contradiction that
$\vecy$ does not contain a convergent subsequence.
We show that if $\vecy$ is unbounded,
then there is another counter-example to Claim~2
in which $\vecy$ is bounded.
Let $r$ be such that $v \in B^F(0,r)$.
Then for every $n$, $D_n \cap B^F(0,r) \neq \emptyset$.
Since $Y$ is component-wise wide,
there are a subsequence $\vecy\fprime$ of $\vecy$, $s > 0$
and a completely discrete sequence of arcs $\setm{L_n}{n \in \bbN}$
such that for every $n$,
$L_n \subseteq D_n$
and $L_n$ connects $y_n'$ with a member of $B^F(0,s)$.
We may assume that $\vecy\fprime = \vecy$.

Denote the endpoint of $L_{2n}$ which is not $y_{2n}$ by $\hatw_n$.
Let $\hatk \in \rfs{EXT}(Y)$ be such that for every $n$,
$\hatk(y_{2n}) = \hatw_n$ and $\hatk(y_{2n + 1}) = y_{2n + 1}$
and set $\hatg = \hatk^{\tau\inverse}$.
Then $\hatg \in \rfs{EXT}(X)$ and hence $\lim \hatg(\vecx)$ exists.
So $\limti{n} \hatg(x_{2n + 1}) = \limti{n} \hatg(x_{2n}) = x$.
Since $\hatk(y_{2n}) = \hatw_n$, it follows that
$\hatg(x_{2n}) = \tau\inverse(\hatw_n)$.
That is, $\tau(\hatg(x_{2n})) = \hatw_n$.
So $\setm{\tau(g(x_{2n}))}{n \in \bbN}$ is bounded and completely
discrete.
By replacing $\vecx$ by $\setm{\hatg(x_{2n})}{n \in \bbN}$
we obtain a counter-example to Claim~2 in which $\vecy$ is bounded.
Since $E$ is a Banach space,
we may also assume that $\vecy$ is spaced.

Since $Y$ is BR.CW.AC, there are $d$ and rectifiable arcs
$L_n \subseteq D_n$ such that $L_n$ connects $y_n$ with $v_n$
and $\rfs{lngth}(L_n) \leq d$.
Let $\gamma_n(t)$ be a parametrization of $L_n$ such that
$\gamma_n(1) = y_n$, 
$\gamma_n(0) = v_n$, 
and for every $t$,
$\rfs{lngth}(\gamma_n([0,t])) = t \mcdot \rfs{lngth}(L_n)$.
For every infinite $\sigma \subseteq \bbN$ let
$s_{\sigma} =
\inf(\setm{t \in [0,1]}{\setm{\gamma_n([t,1])}{n \in \sigma}
\mbox{ is spaced}})$.
Let $\sigma$ be an infinite set such that for every infinite
$\eta \subseteq \sigma$, $s_{\eta} = s_{\sigma}$.
Then $\setm{\gamma_n(s_{\sigma})}{n \in \sigma}$ contains a Cauchy
sequence, and for every $t > s_{\sigma}$,
$\setm{\gamma_n([t,1])}{n \in \sigma}$ is spaced.
Set $s = s_{\sigma}$.
It can be assumed
that $\setm{\gamma_n(s)}{n \in \sigma}$ is a Cauchy sequence,
that $\sigma = \bbN$
and that $s = 0$.
So $\gamma_n(1) = y_n$ for every $n \in \bbN$,
$\setm{\gamma_n(0)}{n \in \bbN}$ is a Cauchy sequence,
and $\setm{\gamma_n([t,1])}{n \in \bbN}$ is spaced
for every $t \in (0,1]$.
Let $w_n = \gamma_n(0)$ and $w = \lim \vecw$.

For every $t > 0$ let $\vecw^t = \setm{\gamma_{2n}(t)}{n \in \bbN}$.
Let $\vecy^{\fzero} = \setm{y_{2n}}{n \in \bbN}$
and $\vecy^{\fone} = \setm{y_{2n + 1}}{n \in \bbN}$.
For every $t > 0$ there is $k_t \in \rfs{EXT}(Y)$
such that $k_t(\vecy^{\fzero}) = \vecw^t$
and $k_t(\vecy^{\fone}) = \vecy^{\fone}$.
This follows from the fact that for $t > 0$,
$\setm{\gamma_n([t,1])}{n \in \bbN}$ is completely discrete.
We check that for every $t \in (0,1]$, $\lim \tau\inverse(\vecw^t) = x$.
Let $h_t = k_t^{\tau\inverse}$.
Then $h_t(x_{2n + 1}) = x_{2n + 1}$
and $h_t(x_{2n}) = \tau\inverse(w^t_n)$.
Clearly, $h_t \in \rfs{EXT}(X)$,
so $h_t$ takes $\vecx$ to a convegent sequence.
But $h_t(x_{2n}) = x_{2n}$,
hence $\lim h_t(\vecx) = \lim_n x_{2n} = x$.
So $\lim_n \tau\inverse(w^t_n) = x$.

Note that for every $\varepsilon > 0$
there are $t_{\varepsilon} > 0$ and $m_{\varepsilon}$
such that for every $t \leq t_{\varepsilon}$
and $n \geq m_{\varepsilon}$,
$\norm{w^t_n - w} < \varepsilon$.
Also, $x_{2n} \simeq^X \tau\inverse(w_n^t)$
for every $n$ and $t$.
It follows that there are sequences
$\vecz$ and $\sngltn{n_k}_{k = 1}^{\infty}$
such that $\lim \vecz = x$, $\lim \tau(\vecz) = w$,
and for every $k$, $z_k \simeq^X x_{2n_k}$.
To see this, take $z_k$ to be $\tau\inverse(w_{n_k}^{t_k})$,
where $\sngltn{t_k}_{k = 1}^{\infty}$ is any sequence converging to $0$
and $n_k$ is such that
$n_k \geq m_{\dgfrac{1}{k}}$
and $\norm{\tau\inverse(w_{n_k}^{t_k}) - x} < \dgfrac{1}{k}$.

From the facts $X$ is BR.CW.AC, $z_k \simeq^X x_{2n_k}$
and $\lim \vecz = \lim_k x_{2n_k}$,
we conclude that there is $g \in \rfs{EXT}(X)$
such that for infinitely many $k$'s,
$g(x_{2n_k}) = z_k$.
We now check that $g^{\tau} \not\in \rfs{EXT}(Y)$, and this is
of course a contradiction.
Using the fact that $\tau(\vecx) = \vecy$,
it is evident that $g^{\tau}$ takes an infinite subsequence of $\vecy$
to an infinite subsequence of $\tau(\vecz)$.
However, $\vecy$ is spaced, and $\tau(\vecz)$ is converges to $w$,
that is, $g^{\tau}$ takes a spaced sequence to a convergent sequence.
Hence $g^{\tau} \not\in \rfs{EXT}(Y)$. A contradiction.
This proves Claim 2.
\smallskip

We prove that $\tau \in \rfs{EXT}(X,Y)$.
Suppose by contradiction that $\vecx \subseteq X$ is a convergent
sequence and $\tau(\vecx)$ is not a convergent sequence.
By Claim 2, there is a subsequence $\vecx^{\fzero}$ of
$\vecx$ such that $\tau(\vecx^{\fzero})$ is convergent.
Since $\tau(\vecx)$ is not convergent,
there is a subsequence $\vecx^{\ftwo}$ of $\vecx$ such that
$d(\tau(\vecx^{\ftwo}),\tau(\vecx^{\fzero})) > 0$.
By Claim 2, there is a subsequence $\vecx^{\fone}$ of
$\vecx^{\ftwo}$ such that $\tau(\vecx^{\fone})$ is convergent.
But $\lim \tau(\vecx^{\fzero}) \neq \lim \tau(\vecx^{\fone})$.
This contradicts Claim 1.
So $\tau \in \rfs{EXT}(X,Y)$.
The assumptions on $X$, $Y$ and $\tau$ were symmetric
with respect to $X$ and~$Y$.
So $\tau \in \rfs{EXT}^{\pm}(X,Y)$.
\rule{1pt}{0pt}\hfill\myqed

\begin{remark}\label{r6.12}
\begin{rm}
The requirement that $X$ be locally movable at its multiple boundary,
which appears in Theorem \ref{metr-bldr-t6.10} is stronger than
what is really needed in the proof of that theorem.
However, the exact assumption needed in that proof is longer
and more complicated, so we include it only as a remark.
Thus in Theorem \ref{metr-bldr-t6.10} the assumption that
$X$ is locally movable at its multiple boundary can be replaced by the
following weaker requirement. The proof remains essentially unchanged.

Let $X \subseteq E$. Then
\begin{itemize}
\addtolength{\parskip}{-09pt}
\addtolength{\itemsep}{06pt}
\item[(1)] 
For every $\vecx \subseteq X$ which is converges to a
multiple boundary point and $z \in \rfs{bd}(X) - \sngltn{\lim \vecx}$,
there is a subsequence $\vecx\fprime$ of $\vecx$ and
$g \in \rfs{EXT}(X)$ such that:
$g(\vecx\fprime) \simeq^X \vecx\fprime$,
$g^{\srfs{cl}}(\lim \vecx\fprime) \neq \lim \vecx\fprime$
and for some $U \in \rfs{Nbr}^E(z)$,
$g \nrestriction (U \cap X) = \rfs{Id}$.
\item[(2)] 
For every $\vecx \subseteq X$ which converges to a
double boundary point and $U \in \rfs{Nbr}^E(\lim \vecx)$
there is a subsequence $\vecx\fprime$ of $\vecx$ and
$g \in \rfs{EXT}(X)$ such that:
$g(\vecx\fprime) \simeq^X \vecx\fprime$,
$g^{\srfs{cl}}(\lim \vecx) \neq \lim \vecx$
and $\rfs{supp}(g) \subseteq U$.
\hfill\proofend
\vspace{-05.7pt}
\end{itemize}
\end{rm}
\end{remark}

The requirement that $X$ be locally movable at its multiple boundary
which appears in Definition \ref{metr-bldr-d6.8}(g)
cannot be entirely omitted.
This is demonstrated by the following trivial example.

\begin{example}\label{metr-bldr-e6.18}
\begin{rm}
There are regular open subsets $X,Y \subseteq \bbR^2$ which satisfy
Clause~1 in the definition of $K_{\srfs{BX}}^{\calO}$
such that
$\rfs{EXT}(X) \cong \rfs{EXT}(Y)$
and $\rfs{cl}(X) \not\cong \rfs{cl}(Y)$.
\end{rm}
\end{example}

\noindent
{\bf Proof }
Let $u \in \bbR^2$ and $F_0,\ldots,F_3 \subseteq \bbR^2$
be closed solid triangles
such that for every $i \neq j$,
$F_i \cap F_j = \sngltn{u}$.
For $i = 1,2,3$ let
$\setm{D_{i,j}}{j < i}$ be a set of pairwise disjoint
closed balls such that
$D_{i,j} \subseteq \rfs{int}(F_i)$ for every $j < i$.
Let
$X = \bigcup_{i < 4} \rfs{int}(F_i) -
\bigcup \setm{D_{i,j}}{i = 1,2,3,\ j < i}$.

Let $v,w \in \bbR^2$
and $G_0,\ldots,G_3 \subseteq \bbR^2$
be closed solid triangles such that
$G_0 \cap G_1 = \sngltn{v}$,
$G_2 \cap G_3 = \sngltn{w}$
and $G_i \cap G_{\ell} = \emptyset$
for every $i \in \dbltn{0}{1}$ and $\ell \in \dbltn{2}{3}$.
For $i = 1,2,3$ let $\setm{E_{i,j}}{i = 1,2,3,\ j < i}$
be a set of pairwise disjoint closed balls such that
$E_{i,j} \subseteq \rfs{int}(G_i)$
for every $j < i$.
Let
$Y = \bigcup_{i < 4} \rfs{int}(G_i) -
\bigcup \setm{E_{i,j}}{i = 1,2,3,\ j < i}$.
Then $X$ and $Y$ are as required in the example.
\smallskip\hfill\myqed

For open subsets of finite-dimensional spaces we have
Theorem \ref{metr-bldr-t5.2} which says that the class of bounded
sets which are the closure of an open UD.AC subset of a
Euclidean space is faithful.
We shall next define another faithful class of spaces which are not
required to be UD.AC.
This class, denoted by $K_{\srfs{IMX}}^{\calO}$,
is defined in \ref{metr-bldr-d6.13}(b).
Loosely speaking,
we replace the assumption that $X$ is UD.AC by the assumption
that the orbit of every $x \in \rfs{bd}(X)$
under $\rfs{EXT}(X)$ contains an arc.
This gives rise to a rather large class.
See Proposition~\ref{metr-bldr-p6.14}.

The next example contains finite and infinite dimensional sets
which belong to $K_{\srfs{IMX}}^{\calO}$ but do not belong
to any of the previously defined EXT-determined classes.
The three examples are connected.
The first example is a subset of $\bbR^2$ which is not UD.AC.
The second set is infinite-dimensional.
It is quite similar to the set $R_1^E$ defined in \ref{metr-bldr-e6.7},
yet it does not belong to $K^{\calO}_{\srfs{BX}}$.
Note the second example is BD.AC, and the first two examples are
regular open.

\begin{example}\label{metr-bldr-e6.12}
\begin{rm}
(a) Let
$R_2 = \setm{(r,\theta)}{\theta \in (\pi,\infty) \mbox{ and }
1 - \frac{1}{\theta - \dgfrac{\pi}{2}} < r <
1 - \frac{1}{\theta + \dgfrac{\pi}{2}}}$.
($R_2$ is described in polar coordinates).
So $R_2$ is an open spiral strip converging to $S(0,1)$.
Note that $R_2$ is connected, $R_2$ is not UD.AC
and $R_2 \not\in K_{\srfs{BX}}^{\calO}$.

(b) Let $E  = \ell_2$
and $R_1^E$ be as in Example~\ref{metr-bldr-e6.7}.
So the set $R_1^E$ is the union of a sequence of pairwise disjoint
open rings converging to $S^E(0,1)$.
We connect each two consecutive rings by an open tube whose closure
is disjoint from the closure of any other ring.
The set of tubes is to be spaced.
{\thickmuskip=2mu \medmuskip=1mu \thinmuskip=1mu 
Let $\sngltn{e_n}_{n \in \bbN}$ be the standard basis of $\ell_2$
and $L_n = [(1 - \frac{1}{2n + 3}) e_n,(1 - \frac{1}{2n + 4}) e_n]$.
}
So each $L_n$ connects two consecutive rings in $R_1^E$.
For some $d > 0$, $\setm{L_n}{n \in \bbN}$ is $d$-spaced.
Let $s_n = \frac{1}{2n + 4} - \frac{1}{2n + 5}$
and $r_n = \min(\dgfrac{d}{3},s_n)$
and $R_3 = R_1^E \cup
\bigcup_{n \in \sboldbbN} B(L_n,r_n)$.
It follows that $R_3$ is connected,
$R_3$ is not UD.AC and $R_3 \not\in K_{\srfs{BX}}^{\calO}$.
However, $R_3$ is JN.AC.

\num{c} Let $E$ be a normed space with dimension $> 2$ and
$F$ be a subspace of $E$ with co-dimension $1$.
Let $R_4^E = B^E(0,2) - \overB^F(0,1)$.
\end{rm}
\end{example}

\begin{defn}\label{metr-bldr-d6.13}
\begin{rm}
(a) Let $\fnn{h}{[a,b] \times Z_1}{Z_2}$ and $t \in [a,b]$.
We denote by $h_t$ the function $g(z) = h(t,z)$.
Let $X \in K_{\srfs{NRM}}^{\calO}$ and $x \in \rfs{bd}(X)$.
We say that $x$ is {\it isotopically movable} with respect to $X$,
if for every $r > 0$ there is a continuous function
$\fnn{h}{[0,1] \times \rfs{cl}(X)}{\rfs{cl}(X)}$
such that $h_0 = \rfs{Id}$, $h_1(x) \neq x$,
and for every $t \in [0,1]$,
$h_t \nrestriction X \in \rfs{EXT}(X)$
and $\rfs{supp}(h_t) \subseteq B(x,r)$.
We say that
$X$ is {\it isotopically movable at its boundary (BR.IS.MV)},
if every $x \in \rfs{bd}(X)$
is isotopically movable with respect to $X$.

   \index{isotopically movable with respect to $X$}
   \index{isotopically movable at the boundary}
   \index{brismv@@BR.IS.MV. Abbreviation of
   isotopically movable at the boundary}

(b) Let $K_{\srfs{IMX}}^{\calO}$
be the class of all open subsets $X$ of a normed space
such that $X$ is JN.AC and BR.IS.MV.
   \index{N@koimx@@$K_{\srfs{IMX}}^{\calO}$.
          Class of open finite-dimensional BR.IS.MV open sets}
\hfill\proofend
\end{rm}
\end{defn}

The next observation and Proposition \ref{metr-bldr-p6.3}
show that $K_{\srfs{IMX}}^{\calO}$ is a large class.
Let $E$ be a normed space and
$X \subseteq E \times (0,\infty)$ be open
and $Z = \setm{z \in E}{\exists a ((z,a) \in X)}$.
The {\it body of revolution} of $X$ is defined as follows.
$$
\rfs{revb}(X) =
\setm{(z,u,v)}
{(z,\sqrt{u^2 + v^2}\kern2pt) \in X}.
$$
So $\rfs{revb}(X)$ is an open subset of $E \times \bbR^2$.
If $\inf(\setm{a}{(z,a)) \in X}) > 0$ for every $z \in Z$,
then $\rfs{revb}(X)$ is called a {\it hollow body of revolution}.
Clearly if $\rfs{revb}(X)$ is hollow,
then $\rfs{revb}(X) \cong X \times S^1$.

\begin{prop}\label{metr-bldr-p6.14}
Let $X,Y \in K^{\calO}_{\srfs{NRM}}$.
\begin{itemize}
\addtolength{\parskip}{-11pt}
\addtolength{\itemsep}{06pt}
\item[\num{1}] 
If $Y$ is BR.IS.MV, then $X \times Y$ is BR.IS.MV.
\item[\num{2}] 
If $X$ and $Y$ are JN.AC, then $X \times Y$ is JN.AC.
\item[\num{3}] 
If $X \subseteq \bbR^n$,
$Y \subseteq \bbR^m$, $X,Y$ are bounded, and $Y$ is BR.IS.MV,
then $X \times Y \in K_{\srfs{IMX}}^{\calO}$.
\item[\num{4}] 
If $X$ and $Y$ are JN.AC
and $Y$ is  BR.IS.MV, then $X \times Y \in K_{\srfs{IMX}}^{\calO}$.
\item[\num{5}] 
If $X \subseteq \bbR^n$, $X$ is bounded
and $\rfs{revb}(X)$ is hollow,
then $\rfs{revb}(X) \in K_{\srfs{IMX}}^{\calO}$.
\vspace{-05.7pt}
\end{itemize}
\end{prop}

\noindent
{\bf Proof } The proof is trivial.
For Parts (3) and (5) see \ref{metr-bldr-p6.3}(b).
\smallskip\hfill\myqed

\noindent
{\bf Remark } The class $K_{\srfs{IMX}}^{\calO}$ does not contain
any of the classes $K_{\srfs{NMX}}^{\calO}$,
$K_{\srfs{BCX}}^{\calO}$
and $K_{\srfs{BX}}^{\calO}$ defined in \ref{metr-bldr-t6.5}
and \ref{metr-bldr-d6.8}(g).
Recall that
$K_{\srfs{BCX}}^{\calO} \subseteq
K_{\srfs{BX}}^{\calO}, K_{\srfs{NMX}}^{\calO}$.
Example \ref{metr-bldr-e6.7} belongs to $K_{\srfs{BCX}}^{\calO}$
but not to $K_{\srfs{IMX}}^{\calO}$.

\begin{theorem}\label{metr-bldr-t6.15}
Suppose that $X,Y \in K_{\srfs{IMX}}^{\calO}$
and $\iso{\varphi}{\rfs{EXT}(X)}{\rfs{EXT}(Y)}$.
Then there is $\tau \in \rfs{EXT}^{\pm}(X,Y)$ such that
$\tau$ induces $\varphi$.
\end{theorem}

\noindent
{\bf Proof }
By Theorem \ref{t2.4}(b), there is $\tau \in H(X,Y)$ such that
$\tau$ induces $\varphi$.

{\bf Claim 1.}
For every $x \in \rfs{bd}(X)$ there is a sequence $\vecx$ converging
to $x$ such that $\tau(\vecx)$ converges to a member of
$\rfs{bd}(Y)$.
{\bf Proof } This claim follows from
Lemma \ref{metr-bldr-l6.4}(b) applied to $\tau\inverse$.
\smallskip

{\bf Claim 2.} Let $x \in \rfs{bd}(X)$ and $\vecx,\vecu \subseteq X$.
Suppose that $\lim \vecx = \lim \vecu = x$ and that
$\tau(\vecx)$ and $\tau(\vecu)$ are convergent.
Then $\lim \tau(\vecx) = \lim \tau(\vecu)$.
{\bf Proof }
Set $\vecy = \tau(\vecx)$,\break
$\vecv = \tau(\vecu)$,
$y = \lim \vecy$, $v = \lim \vecv$, and suppose by contradiction that
$y \neq v$.
Obviously, $y,v \in \rfs{bd}(Y)$.
Let $r = \dgfrac{\norm{y - v}}{2}$.
We may assume that $\vecv \subseteq B(v,r)$ and that
$\vecy \cap B(v,r) = \emptyset$.
Let\break
$\fnn{h}{[0,1] \times \rfs{cl}(Y)}{\rfs{cl}(Y)}$
be an isotopy
as assured by the fact that
$v$ is isotopically movable with respect to $Y$,
and
such that for every $t \in [0,1]$, $\rfs{supp}(h_t) \subseteq B(v,r)$.

For every $t \in [0,1]$
let $u_{n,t} = \tau\inverse(h(t,v_n))$.
We first prove the following fact.
$(*)$ For every $t \in [0,1]$,
$\limti{n} u_{n,t} = x$.
Let $t \in [0,1]$. Let $\barh = h_t \nrestriction Y$
and $\barg = \barh^{\tau\inverse}$.
Then $\barg \in \rfs{EXT}(X)$.
Also $\barg \nrestriction \vecx = \rfs{Id}$.
So $\barg^{\srfs{cl}}(x) = x$.
Hence
$\limti{n} u_{n,t} = \limti{n} \barg(u_n) =
\barg(\limti{n} u_n) = \barg(x) = x$.
So $(*)$ is proved.

Let $L_n = h([0,1] \times \sngltn{v_n})$ and $K_n = \tau\inverse(L_n)$.
We prove that $\limti{n} K_n = x$.
Suppose by contradiction that this is not true.
Then there are $d > 0$, $\vect \subseteq [0,1]$
and a $\onetoonen$ sequence $\setm{n_i}{i \in \bbN}$ such that
$d(x,u_{n_i,t_i}) \geq d$ for every $i \in \bbN$.
We may assume that $\vect$ is convergent.
Set $t^* = \lim \vect$. Let $I_i$ be the closed interval whose
endpoints are $t_i$ and $t^*$
and\break
$J_i = h(I_i \times \sngltn{v_{n_i}})$.
Then $\limti{i} J_i = h(t^*,v)$.
Since for every $t \in [0,1]$, $h_t \nrestriction Y \in \rfs{EXT}(Y)$
and $v \in \rfs{bd}(Y)$,
it follows that $h(t^*,v) \in \rfs{bd}(Y)$.
The fact $\vecv \subseteq Y$ implies that $J_i \subseteq Y$,
and hence $h(t^*,v) \not\in J_i$
for every $i \in \bbN$.
Since $J_i$ is compact, $d(J_i,h(t^*,v)) > 0$.
We may thus replace $\sngltn{n_i}_{i \in \bbN}$ by a subsequence
and obtain that
\hbox{$\max(\setm{d(z,h(t^*,v))}{z \in J_{i + 1}}) <$} $d(J_i,h(t^*,v))$
for every $i \in \bbN$.
There is a sequence $\sngltn{V_i}_{i \in \bbN}$ of open sets
such that for any distinct $i,j \in \bbN$,
$J_i \subseteq V_i \subseteq \rfs{cl}(V_i) \subseteq Y \cap B(v,r)$,
$V_i \cap V_j = \emptyset$
and $\limti{i} V_i = \limti{i} J_i$.
From the fact that $J_i$ connects $h(v_{n_i},t^*)$ and $h(v_{n_i},t_i)$,
it follows that there is $h_i \in \rfs{UC}(Y)\sprtm{V_i}$
such that $h_i(h(v_{n_i},t^*)) = h(v_{n_i},t_i)$.
Let $\hath = \bcirc_{i \in \bbN} h_i$.
Then by Proposition~\ref{metr-bldr-p4.5},
$\hath \in \rfs{UC}_0(Y) \subseteq \rfs{EXT}(Y)$.
Clearly, $\rfs{supp}(\hath) \subseteq B(v,r)$
and so $\hath \nrestriction \vecy = \rfs{Id}$.
Let $\hatg = \hath^{\tau\inverse}$.
So $\hatg \in \rfs{EXT}(X)$.
Since $\hath \nrestriction \vecy = \rfs{Id}$,
it follows that $\hatg \nrestriction \vecx = \rfs{Id}$
and hence $\hatg^{\srfs{cl}}(x) = x$.
Clearly, for every~$i$,
$\hatg(u_{n_i,t^*}) = u_{n_i,t_i}$,
and from~$(*)$ it follows that $\limti{i} u_{n_i,t^*} = x$.
So
$$
\limti{i} u_{n_i,t_i} =
\limti{i} \hatg(u_{n_i,t^*}) =
\hatg^{\srfs{cl}}(\limti{i} u_{n_i,t^*}) =
\hatg^{\srfs{cl}}(x) = x. 
$$
This contradicts the fact that
$d(x,u_{n_i,t_i}) \geq d$.
So $\limti{n} K_n = x$.

There is an infinite set $\sigma \subseteq \bbN$
such that $K_i \cap K_j = \emptyset$
for any distinct $i,j \in \sigma$.
Let $\setm{U_i}{i \in \sigma}$ be such that
$K_i \subseteq U_i \subseteq X$,
$U_i$ is open,
$U_i \cap U_j = \emptyset$ for every $i \neq j$
and $\lim_{i \in \sigma} U_i = x$.
Let $\eta \subseteq \sigma$ be such that $\eta$ and $\sigma - \eta$
are infinite.
For every $i \in \eta$ let $g_i \in \rfs{UC}(X)\sprtm{U_i}$
be such that $g_i(u_i) = u_{i,1}$.
Let $\barg = \bcirc_{i \in \eta} g_i$ and $\barh = \barg^{\tau}$.
By Proposition \ref{metr-bldr-p4.5},
$\barg \in \rfs{UC}_0(X) \subseteq \rfs{EXT}(X)$,
hence it follows that $\barh \in \rfs{EXT}(Y)$.

For every $i \in \eta$, $\barh(v_i) = h(v_i,1)$,
so $\lim_{i \in \eta} \barh(v_i) = h(v,1)$.
For every $i \in \sigma - \eta$, $\barh(v_i) = v_i$,
so $\lim_{i \in \sigma - \eta} \barh(v_i) = v$.
Recall that $h(v,1) \neq v$.
Also, $\limti{i} v_i = v$. So $\vecv$ is convergent and
$\barh(\vecv)$ is not convergent. Hence $\barh \not\in \rfs{EXT}(Y)$.
A contradiction, so Claim 2 is proved.\smallskip

Suppose by contradiction that $x \in \rfs{bd}(X)$ and
$x \not\in \rfs{Dom}(\tau^{\srfs{cl}})$.
By Claim~1,
there is a sequence $\vecx \subseteq X$ such that $\lim \vecx = x$
and $\tau(\vecx)$ is convergent. Set $y = \lim \tau(\vecx)$.
There are a $\onetoonen$ sequence $\vecu \subseteq X$ and $d > 0$
such that $\lim \vecu = x$ and $d(\tau(\vecu),y) \geq d$.
Define $\vecv = \tau(\vecu)$.
By Claim 2, $\vecv$ does not have a convergent subsequence.
That is, $\vecv$ is completely discrete. 
Since $Y$ is JN.AC, there is a subsequence $\vecw$ of $\vecv$
such that $\vecw$ has a joining system.
Let
$\fortpl{\vecw}{w^*}{\setm{L_n}{n \in \bbN}}{\vecw\fprime}$
be a joining system for $\vecw$.
We may assume that $w^* \not\in \rfs{Rng}(\vecw)$.

We show that it can be assumed that $w^* \neq y$.
Suppose that $w^* = y$. Let $r = d(\vecw,y)$.
Since $Y$ is BR.IS.MV and $y \in \rfs{bd}(Y)$,
there is $h \in \rfs{EXT}(Y)$ such that $\rfs{supp}(h) \subseteq B(y,r)$
and $h^{\srfs{cl}}(y) \neq y$.
So $h \nrestriction \vecw = \rfs{Id}$.
It follows that 
$\fortpl{\vecw}{h^{\srfs{cl}}(y)}{\setm{h(L_n)}{n \in \bbN}}
{h(\vecw\fprime)}$
is a joining system for $\vecw$.
So we may assume that $w^* \neq y$.

Recall that $Y$ is JN.AC. So we may apply Lemma \ref{metr-bldr-l6.4}(b)
to $\tau\inverse$.
Recall also that
\hbox{$\lim \tau\inverse(\vecw) = \lim \tau\inverse(\vecv) = x$.}
Hence there is $\vecz \subseteq Y$
such that $\lim \vecz = w^*$ and $\lim \tau\inverse(\vecz) = x$.
We now have two sequences: $\vecx$ and $\tau\inverse(\vecz)$,
both converge to $x$,
and $\tau(\vecx)$ and $\tau(\tau\inverse(\vecz))$ are convergent,
but not to the same point.
This contradicts Claim 2, so $x \in \rfs{Dom}(\tau^{\srfs{cl}})$.

We have shown that $\tau \in \rfs{EXT}(X,Y)$,
and an identical argument shows that
$\tau\inverse \in \rfs{EXT}(Y,X)$.
That is, $\tau \in \rfs{EXT}^{\pm}(X,Y)$.
\medskip\hfill\myqed

\subsection{Completely locally uniformly continuous homeomorphism\\
groups.}
\label{ss6.3}
\label{ss6.3-CMP.LUC-homeomorphisms}

Having obtained the results about $\rfs{EXT}(X)$ and $\rfs{LUC}(X)$,
only little extra work is needed to prove
CMP.LUC\,-\,determined-ness. See Definition \ref{metr-bldr-d5.3}(f).
This faithfulness result will complete the picture on groups of type
$H_{\itGamma}^{\srfs{CMP.LC}}(X)$ discussed in
Chapters \ref{s8}\,-\ref{s12}.

The following is a strengthening of property BR.LC.AC.

\begin{defn}\label {metr-bldr-d6.16}
$X$ is
{\it locally uniformly\,-\,in\,-\,diameter arcwise\,-\,connected
(LC.UD.AC)},
if for every $x \in \rfs{bd}(X)$ there is $U \in \rfs{Nbr}(x)$
such that for every $\varepsilon > 0$ there is $\delta > 0$ such that
for every $u,v \in U$: if $d(u,v) < \delta$, then there is an arc
$L \subseteq X$ connecting $u$ and $v$
such that $\rfs{diam}(L) < \varepsilon$.
\begin{rm}
\vspace{-2.0mm}
\end{rm}
\end{defn}

\begin{theorem}\label{metr-bldr-t6.17}
\num{a} Let $X,Y \in K_{\srfs{NRM}}^{\calO}$.
Suppose that $X$ and $Y$ are LC.UD.AC and JN.AC.
Let $\iso{\varphi}{\rfs{CMP.LUC}(X)}{\rfs{CMP.LUC}(Y)}$.
Then there is
$\tau \in \rfs{CMP.LUC}^{\pm}(X,Y)$ such that
$\tau$ induces $\varphi$.

\num{b} Suppose that $X$ is LC.UD.AC and $Y$ is JN.AC,
and let $\tau \in H(X,Y)$ be such\break
that
$(\rfs{UC}_0(X))^{\tau} \subseteq \rfs{CMP.LUC}(Y)$ \ and \ %
$(\rfs{LUC}_{01}(Y))^{\tau^{-1}} \subseteq \rfs{CMP.LUC}(X)$.
Then $\tau \in \rfs{CMP.LUC}(X,Y)$.
\end{theorem}

\noindent
{\bf Proof }
We shall see that (b) implies (a). So we start by proving Part (b).

(b) It is trivial that $X$ is BR.LC.AC.
We first show that $\tau \in \rfs{EXT}(X,Y)$.
By definition, $\rfs{CMP.LUC}(X) \subseteq \rfs{EXT}(X)$.
So $(\rfs{UC}_0(X))^{\tau} \subseteq \rfs{EXT}(Y)$
and
$(\rfs{LUC}_{01}(Y))^{\tau^{-1}} \subseteq \rfs{EXT}(X)$.
By Corollary \ref{metr-bldr-c6.5}(a), $\tau \in \rfs{EXT}(X,Y)$.

We show that $\tau \in \rfs{LUC}(X,Y)$.
Let $\calS$ be the set of BPD-subsets of $X$.
Then $\rfs{UC}(X,\calS) \subseteq \rfs{UC}_0(X)$
and $\rfs{CMP.LUC}(Y) \subseteq \rfs{LUC}(Y)$.
So
$(\rfs{UC}(X,\calS))^{\tau} \subseteq \rfs{LUC}(Y)$.
By Theorem \ref{metr-bldr-t4.8}(b),
$\tau \in \rfs{LUC}^{\pm}(X,Y)$.

Let $x^* \in \rfs{bd}(X)$. We show that there is
$U \in \rfs{Nbr}(x^*)$ such that $\tau \nrestriction (U \cap X)$ is UC.
The proof is very much a repetition of the proof of Part 1 of
Theorem \ref{metr-bldr-t4.8}(c).

Suppose by contradiction that for every $U \in \rfs{Nbr}^X(x^*)$,
$\tau \nrestriction U$ is not UC.
The following claim is an easy consequence of the fact that
$\tau \nrestriction B(x^*,r) \cap X$ is not UC. Its proof is left
to the reader.

{\bf Claim 1.} For every $r > 0$ there are sequences
$\vecx,\vecy$ and $d,e > 0$ such that:
\begin{itemize}
\addtolength{\parskip}{-11pt}
\addtolength{\itemsep}{06pt}
\item[(1)] 
$\rfs{Rng}(\vecx) \cup
\rfs{Rng}(\vecy) \subseteq B^X(x^*,\dgfrac{r}{2})$;
\item[(2)] 
$\limti{n} \norm{x_n - y_n} = 0$;
\item[(3)] 
either (i) for any distinct $m,n \in \bbN$,
$d(\uopair{x_m}{y_m},\uopair{x_n}{y_n}) \geq e$,
or (ii) $\vecx$ is a Cauchy sequence;
\item[(4)] 
$d(\rfs{Rng}(\vecx) \cup \rfs{Rng}(\vecy),x^*) > e$;
\item[(5)] 
for every $n \in \bbN$,
$\norm{\tau(x_n) - \tau(y_n)} \geq d$.
\vspace{-05.7pt}
\end{itemize}

Let $U \in \rfs{Nbr}(x^*)$ be as assured by the LC.UD.AC-ness of $X$.
There is $a > 0$ and a function $\fnn{\eta}{(0,a]}{\bbR}$
such that $\lim_{t \rightarrow 0} \eta(t) = 0$
and for every $u,v \in U \cap X$, if $\norm{u - v} \leq t$,
then there is an arc $L \subseteq X$ connecting $u$ and $v$ sucu that
$\rfs{diam}(L) \leq \eta(t)$.

Let $e_{-1} > 0$ be such that $B^E(x^*,e_{-1}) \subseteq U$.
It is easy to define by induction on $i \in \bbN$,
$r_i > 0$, sequences $\vecx\farsu{i},\vecy\farsu{i}$
and $d_i,e_i > 0$
such that: (i) $\vecx\farsu{i},\vecy\farsu{i}, d_i,e_i$
satisfy the conclusion of Claim 1 for $r_i$; and
(ii) for every $i \in \bbN$, $r_i = \dgfrac{e_{i - 1}}{8}$.
Clearly $e_{i + 1} \leq \dgfrac{e_i}{4}$.
By deleting initial segments from the $\vecx\farsu{i}$'s and
$\vecy\farsu{i}$'s,
we may further assume that for every $i,n \in \bbN$,
$\eta(\norm{x^i_n - y^i_n}) < \dgfrac{e_i}{8}$.
We may further assume that either for every $i \in \bbN$
Clause (3)(i) of Claim 1 holds,
or for every $i \in \bbN$
Clause (3)(ii) of Claim 1 holds.

{\bf Case 1 } Clause (3)(i) of Claim 1 holds.
Let $\setm{\pair{i(k)}{n(k)}}{k \in \bbN}$
be a $\onetoonen$ enumerarion of $\bbN^2$.
Then
$\limti{k} \norm{x^{i(k)}_{n(k)} - y^{i(k)}_{n(k)}} = 0$.
Set $u_k = x^{i(k)}_{n(k)}$, $v_k = y^{i(k)}_{n(k)}$
and let $L_k \subseteq X$ be an arc connecting $u_k$ and $v_k$ such that
$\rfs{diam}(L_k) \leq \eta(\norm{u_k - v_k})$.
Let $B_k = B(L_k,\dgfrac{e_{i(k) + 1}}{4})$.
Then
$$
\rfs{diam}(B_k) \leq \rfs{diam}(L_k) + \dgfrac{e_{i(k) + 1}}{2} \leq
\eta(\norm{u_k - v_k}) + \dgfrac{e_{i(k) + 1}}{2} \leq
\dgfrac{e_{i(k)}}{8} + \dgfrac{e_{i(k) + 1}}{2} \leq
\dgfrac{e_{i(k)}}{4}.
$$
It follows that if $i(k) = i(\ell)$,
then $d(B_k,B_{\ell}) \geq \dgfrac{e_{i(k)}}{2}$.
Suppose that $i(k) < i(\ell)$.
Then $\norm{u_{\ell} - u_k} \geq \dgfrac{7 e_{i_{\ell}}}{8}$,
$\rfs{diam}(B_k) \leq \dgfrac{e_{i(k)}}{4} \leq \dgfrac{e_{i(\ell)}}{4}$
and
$\rfs{diam}(B_{\ell}) \leq \dgfrac{e_{i(\ell)}}{4}$.
So $d(B_k,B_{\ell}) \geq \dgfrac{3 e_{i_{\ell}}}{8}$.
Obviously, $\limti{k} \rfs{diam}(B_k) = 0$.
Let $w_k \in L_k - \sngltn{u_k}$
be such that $\norm{\tau(w_k) - \tau(u_k)} < \frac{1}{k + 1}$.
By Lemma \ref{l2.6}(d),
there is $h_k \in \rfs{LIP}(X)$ such that
$\rfs{supp}(h_k) \subseteq B_k$, $h_k(u_k) = u_k$ and $h_k(w_k) = v_k$.
By Propostion \ref{metr-bldr-p4.5},
$h \eqdf \bcirc_{k \in \sboldbbN} h_k \in \rfs{UC}(X)$
and indeed $h \in \rfs{UC}_0(X)$.

Let us see that for every $V \in \rfs{Nbr}(\tau^{\srfs{cl}}(x^*))$,
$h^{\tau} \nrestriction (V \cap Y)$ is not UC.
For $i \in \bbN$ define $\sigma_i = \setm{k}{i(k) = i}$.
So if $k \in \sigma_i$,
then $L_k \subseteq B(x^*,\eta(2r_k))$.
Since $\limti{i} \eta(2r_i) = 0$, and since $\tau^{\srfs{cl}}$ is
continuous at $x^*$, there is $i$ such that for every $k \in \sigma_i$,
$\tau(L_k) \subseteq V$.

For every $k \in \sigma_i$, $\tau(u_i),\tau(w_i) \in V$.
Clearly, $\lim_{k \in \sigma_i} \norm{\tau(u_k) - \tau(w_k)} = 0$.
However, for every $k \in \sigma_i$,
$\norm{h^{\tau}(\tau(u_k)) - h^{\tau}(\tau(w_k))} =
\norm{\tau(u_i)) - \tau(v_i))} \geq d_i$.
So $h^{\tau} \nrestriction (V \cap Y)$ is not UC.
Hence $h^{\tau} \not\in \rfs{CMP.LUC}(Y)$
even though $h \in \rfs{UC}_0(X)$, a contradiction.

{\bf Case 2 } Clause (3)(ii) of Claim 1 holds.
Let $\barz_i = \lim \vecx\farsu{i}$.
Clearly,
$\barz_i \in B^{\oversE}(x^*,r_i) - B^{\oversE}(x^*,e_i)$.
So $\setm{\barz_i}{i \in \bbN}$ is $\onetoonen$
and $\limti{i} \barz_i = x^*$.
Also, $\barz_i \in \overE - E$. This is so,
because if $\barz_i \in E$,
then either $\barz_i \in X$ and $\tau$ is not continuous at $\barz_i$,
or $\barz_i \in \rfs{bd}^E(X)$ and
$\barz_i \not \in \rfs{Dom}(\tau^{\srfs{cl}})$.
Both situations are impossible.
For every $i$ and $n$ let $L_{i,n} \subseteq X$ be an arc connecting
$x^i_n$ and $y^i_n$ such that
$\rfs{diam}(L_{i,n}) \leq \eta(\norm{x^i_n - y^i_n})$.
Note that for every $i$, $\limti{n} L_{i,n} = \barz_i$.
From the facts $\barz_i \not\in E$ and $L_{i,n} \subseteq E$
we conclude that $d(\barz_i,L_{i,n}) > 0$.

It follows easily that there is a sequence
$\setm{\pair{i(k)}{n(k)}}{k \in \bbN}$ such that
\begin{itemize}
\addtolength{\parskip}{-11pt}
\addtolength{\itemsep}{06pt}
\item[(1)] 
for every $i \in \bbN$, $\setm{k}{i(k) = i}$ is infinite,
\item[(2)] 
for every $k \in \bbN$,
$c_k \eqdf d(L_{i(k),n(k)}\kern1pt,
\kern1pt\bigcup_{m \neq k} L_{i(m),n(m)}) > 0$.
\vspace{-05.7pt}
\end{itemize}
It is also clear from the construction that
\begin{itemize}
\addtolength{\parskip}{-11pt}
\addtolength{\itemsep}{06pt}
\item[(3)] 
$\limti{k} \rfs{diam}(L_{i(k),n(k)}) = 0$.
\vspace{-05.7pt}
\end{itemize}

Set $L_k = L_{i(k),n(k)}$, $u_k = x^{i(k)}_{n(k)}$,
$v_k = y^{i(k)}_{n(k)}$ and $B_k = B(L_k,\dgfrac{c_k}{3})$.
Clearly, for every $\ell \neq k$, $d(B_{\ell},B_k) \geq c_k$
and $\limti{k} \rfs{diam}(B_k) = 0$.
From this point on the proof proceeds exactly as in Case 1.
So in Case 2 too, a contradiction is reached.

It follows that there is $U \in \rfs{Nbr}(x^*)$
such that $\tau \nrestriction (U \cap X)$ is UC,
and this implies that $\tau^{\srfs{cl}}$ is UC at $x^*$.
Recall that we have already shown before that $\tau \in \rfs{EXT}(X,Y)$
and that $\tau \in \rfs{LUC}(X,Y)$. So $\tau \in \rfs{CMP.LUC}(X,Y)$.
\smallskip

(a)
Let $\iso{\varphi}{\rfs{CMP.LUC}(X)}{\rfs{CMP.LUC}(Y)}$.
Clearly,
$\rfs{LIP}^{\srfs{LC}}(X) \leq \rfs{CMP.LUC}(X) \leq H(X)$,
and the same holds for $Y$.
So By Theorem \ref{metr-bldr-t2.8}(a),
there is $\tau \in H(X,Y)$ such that $\tau$ induces $\varphi$.
Hence $(\rfs{CMP.LUC}(X))^{\tau} = \rfs{CMP.LUC}(Y)$.
Obviously, $\rfs{UC}_0(X) \subseteq \rfs{CMP.LUC}(X)$
and $\rfs{LUC}_{01}(Y) \subseteq \rfs{CMP.LUC}(Y)$.
So Part (b) of this lemma can be applied.
Hence $\tau \in \rfs{CMP.LUC}(X,Y)$.
Similarly, $\tau\inverse \in \rfs{CMP.LUC}(Y,X)$. That is,
$\tau \in \rfs{CMP.LUC}^{\pm}(X,Y)$.
\hfill\myqed

\subsection{The reconstruction of $\rfs{cl}(X)$ from $H(\rfs{cl}(X))$.}
\label{ss6.4}

The next two theorems \ref{t4.35} and \ref{metr-bldr-t6.22}
deal with the reconstruction of $F$ from $H(F)$,
when $F$ is the closure of an open subset of a normed space.
The sets to which these theorems apply
may have rather complicated boundaries.
It is not true though that for every $F,K$
which are the closure of an open subset of a normed space,
$H(F) \cong H(K)$ implies that $F \cong K$.
See Example \ref{metr-bldr-e5.9}.

Recall that if $A \subseteq E$ has a nonempty interior,
then
$\rfs{ENI}(A) \eqdf
\setm{h(x)}{x \in \rfs{int}^E(A) \mbox{ and}\break
h \in H(A)}$.
For $f \in  \rfs{UC}_0(X)$,
define
$f^{\srfs{eni}} = f^{\srfs{cl}} \nrestriction \rfs{ENI}(\rfs{cl}(X))$.
{\thickmuskip=2.5mu \medmuskip=1.5mu \thinmuskip=1mu 
Hence $f^{\srfs{eni}} \in H(\rfs{ENI}(\rfs{cl}(X)))$.
}
Also define
$\rfs{UC}_0^{\srfs{eni}}(X) =
\setm{f^{\srfs{eni}}}{f \in  \rfs{UC}_0(X)}$.

   \index{N@AAAA@@$f^{\srfs{eni}} =
          f^{\srfs{cl}} \nrestriction \rfs{ENI}(\rfs{cl}(X))$}
   \index{N@uceni@@$\rfs{UC}_0^{\srfs{eni}}(X) =
          \setm{f^{\srfs{eni}}}{f \in  \rfs{UC}_0(X)}$}

Parts (a) and (b) of the next proposition
are analogous to Proposition~\ref{metr-bldr-p6.1}
and\break
Lemma~\ref{metr-bldr-l6.4}(a).
The proofs of (a) and (b) are essentially identical to the proofs
of their counterparts, so they are omitted.
Part (c) is analogous to Lemma~\ref{metr-bldr-l6.4}(b),
but (c) is stated for $\eta\inverse$ rather than for $\eta$.

\begin{prop}\label{metr-bldr-p6.19}
\num{a} Let $X$ be BR.LC.AC
and $\tau \in H(\rfs{ENI}(\rfs{cl}(X)),\rfs{ENI}(\rfs{cl}(Y)))$.
Assume that
$(\rfs{UC}_0^{\srfs{eni}}(X))^{\tau} \subseteq
\rfs{EXT}(\rfs{ENI}(\rfs{cl}(Y)))$.
Let $x \in \rfs{bd}(X) - \rfs{ENI}(\rfs{cl}(X))$,
$y \in \rfs{bd}(Y)$ and $\vecx \subseteq X$ be such that
$\lim \vecx = x$ and $\lim \tau(\vecx) = y$. 
Then $(\tau \nrestriction X) \cup \sngltn{\pair{x}{y}}$ is continuous.

\num{b} Let $X$ be JN.AC
and $\tau \in H(\rfs{ENI}(\rfs{cl}(X)),\rfs{ENI}(\rfs{cl}(Y)))$
be such that
$(\rfs{LUC}_{01}(X))^{\tau} \subseteq H(\rfs{ENI}(\rfs{cl}(Y)))$.
Let $y \in \rfs{bd}(Y) - \rfs{ENI}(\rfs{cl}(Y))$.
Suppose that $\vecx \subseteq X$ is completely discrete,
$\fortpl{\vecx}{x^*}{\setm{L_n}{n \in \bbN}}{\vecx\fprime}$
is a joining system for $\vecx$ and $\lim \tau(\vecx) = y$.
Then there is a sequence $\vecu \subseteq X$ such that
$\lim \vecu = x^*$ and $\lim \tau(\vecu) = y$.

\num{c}
Let $X,Y \in K^{\calO}_{\srfs{NRM}}$.
Assume that $Y$ is JN.AC.
Set $K = \rfs{cl}(X)$ and $M = \rfs{cl}(Y)$,
and let $\eta \in H(\rfs{ENI}(K),\rfs{ENI}(M))$ be such that
for every $h \in H(M)$,
$((h \nrestriction \rfs{ENI}(M))^{\eta\inverse})^{\srfs{cl}} \in H(K)$.
Then for every $x \in K - \rfs{ENI}(K)$
there is a sequence $\vecx \subseteq X$ converging to $x$
such that $\eta(\vecx) \subseteq Y$,
and $\eta(\vecx)$ is convergent in $M$.
\end{prop}

\noindent
{\bf Proof }
(c)
Let $x \in K - \rfs{ENI}(K)$.
Let $\vecx\fprime \subseteq X$ be a sequence converging to $x$. 
For every $n \in \bbN$ let $r_n = \min(\delta(x'_n),d(x'_n,x))$.
So $B^E(x'_n,r_n)$ is a nonempty open subset of $\rfs{ENI}(K)$.
\hbox{Clearly,
$\rfs{bd}(Y) \cap \rfs{ENI}(M)$} is nowhere dense in $\rfs{ENI}(M)$.
So there is $x_n \in B^E(x'_n,r_n)$ such that
$\eta(x_n) \not\in \rfs{bd}(Y) \cap \rfs{ENI}(M)$.
That is, $\eta(x_n) \in Y$.
So $\vecx \subseteq X$, $\lim \vecx = x$ and $\eta(\vecx) \subseteq Y$.

Define $\vecy = \eta(\vecx)$.
Suppose that $\vecy$ has a subsequence $\vecy\fprime$ such that
$\vecy\fprime$ is convergent in $\rfs{cl}(Y)$.
Then $\eta\inverse(\vecy\fprime)$
is as required in the proposition.
Suppose that such a $\vecy\fprime$ does not exist.
Hence $\vecy$ is completely discrete.

Let $\fortpl{\vecy}{y^*}{\setm{L_n}{n \in \bbN}}{\vecy\fprime}$
be a joining system for $\vecy$.
By \ref{metr-bldr-p6.19}(b) applied to $\vecy$ and $\eta\inverse$,
there is $\vecv \subseteq Y$
such that $\lim \vecv = y^*$ and $\lim \eta\inverse(\vecv) = x$.
It is obvious that\break
$y^* \in \rfs{bd}(Y) - \rfs{ENI}(\rfs{cl}(Y))$.

As in the beginning of the proof, there is a sequence
$\vecv\fprime \subseteq Y$ such that $\lim \vecv\fprime = y^*$,
$\eta\inverse(\vecv\fprime) \subseteq X$ and
$\lim \eta\inverse(\vecv\fprime) = \lim \eta\inverse(\vecv) = x$.
So $\eta\inverse(\vecv\fprime)$ is as required.
\smallskip\hfill\myqed

The following theorem is analogous to
Theorem~\ref{metr-bldr-t6.5}(b). The proofs are essentially the same.

\begin{theorem}\label{t4.35}\label{metr-bldr-t6.20}
{\thickmuskip=2mu \medmuskip=1mu \thinmuskip=1mu 
Let $X,Y \in K^{\calO}_{\srfs{NMX}}$. (See \ref{metr-bldr-t6.5}(b)).
If
$\iso{\varphi}{H(\rfs{cl}(X))}{H(\rfs{cl}(Y))}$,
}
then there is $\iso{\tau}{\rfs{cl}(X)}{\rfs{cl}(Y)}$,
such that $\tau$ induces $\varphi$.
\end{theorem}

\noindent
{\bf Proof }
Let $K = \rfs{cl}(X)$ and $M = \rfs{cl}(Y)$.
From Theorem \ref{t2.15}(c) it follows that there is
$\eta \in H(\rfs{ENI}(K),\rfs{ENI}(M))$
such that $\eta$ induces $\varphi$.

For every $x \in \rfs{bd}(X) - \rfs{ENI}(\rfs{cl}(X))$
let $\vecx \subseteq X$ be such that $\lim \vecx = x$ and
$\eta(\vecx)$ is convergent in $M$.
The existence of $\vecx$
is assured by Proposition \ref{metr-bldr-p6.19}(c).
Let $y_x = \lim\, \eta(\vecx)$.
Since $\rfs{Rng}(\eta) \supseteq Y$, $y_x \in \rfs{bd}(Y)$.
Since $\eta$ induces $\varphi$, for every $g \in H(K)$,\break
$((g \nrestriction \rfs{ENI}(K))^{\eta})^{\srfs{cl}} \in
\rfs{EXT}(\rfs{ENI}(M))$.
In particular,
$(\rfs{UC}_0^{\srfs{eni}}(X))^{\eta} \subseteq
\rfs{EXT}(\rfs{ENI}(M))$.
Hence by Proposition~\ref{metr-bldr-p6.19}(a),
$\eta \nrestriction X \cup \sngltn{\pair{x}{y_x}}$ is continuous.
Also, for every $x \in \rfs{bd}(X) \cap \rfs{ENI}(\rfs{cl}(X))$,
$\eta \nrestriction X \cup \sngltn{\pair{x}{\eta(x)}}$ is continuous.
We thus have
\begin{itemize}
\addtolength{\parskip}{-11pt}
\addtolength{\itemsep}{06pt}
\item[(1)] 
for every $x \in \rfs{bd}(X) - \rfs{ENI}(\rfs{cl}(X))$,
$\eta \nrestriction X \cup \sngltn{\pair{x}{y_x}}$ is continuous,
\item[(2)] 
for every $x \in \rfs{bd}(X) \cap \rfs{ENI}(\rfs{cl}(X))$,
$\eta \nrestriction X \cup \sngltn{\pair{x}{\eta(x)}}$ is continuous.
\vspace{-05.7pt}
\end{itemize}
So by Proposition \ref{metr-bldr-p4.7}(a),
$\eta \cup
\setm{\pair{x}{y_x}}{x \in \rfs{bd}(X) - \rfs{ENI}(\rfs{cl}(X))}$
is continuous.
So $\eta$ can be extended to a continuous function $\tau$
from $\rfs{cl}(X)$ to $\rfs{cl}(Y)$.

Similarly,
$\eta\inverse$ can be extended to a continuous function $\rho$ from
$\rfs{cl}(Y)$ to $\rfs{cl}(X)$.
{\thickmuskip=3.3mu \medmuskip=2.2mu \thinmuskip=1.1mu 
It follows easily that $\tau$ is $\onetoonen$ and that
$\tau\inverse = \rho$.
So $\tau \in H(\rfs{cl}(X),\rfs{cl}(Y))$.
}
Since $\eta$ induces $\varphi$ and $\rfs{Dom}(\eta)$ is dense
in $\rfs{Dom}(\tau)$,
it follows that $\tau$ induces $\varphi$.
\bigskip\hfill\myqed

\begin{prop}\label{metr-bldr-p6.21}
\num{a} Let $X \in K^{\calO}_{\srfs{NRM}}$, $K = \rfs{cl}(X)$,
$U \subseteq \rfs{ENI}(K)$ be open in~$K$, $L \subseteq U$ be an arc
and $x,y$ be the endpoints of $L$.
Then there is $h \in H(K)\sprt{U}$ such that $h(x) = y$.

\num{b} Let $Z$ be a topological space $z \in Z$
and $\setm{h_i}{i \in \bbN} \subseteq H(Z)$
be such that for every $i \neq j$,
$\rfs{supp}(h_i) \cap \rfs{supp}(h_j) = \emptyset$
and $\limti{i} \rfs{supp}(h_i) = z$.
Then $\bcirc_{i \in \bbN} h_i \in H(Z)$.
\end{prop}

\noindent
{\bf Proof }
(a) Let $\fnn{\gamma}{[0,1]}{L}$ be a paramtrization of $L$
such that $\gamma(0) = x$ and $\gamma(1) = y$.
There are $n \in \bbN$, $\setm{U_i}{i < n}$
and $0 = t_0 < \ldots t_n = 1$
such that for every $i < n$:
$U_i$ is open in $K$,
$U_i$ is homeomorphic to an open ball of a normed space,
$U_i \subseteq U$
and $\gamma([t_i,t_{i + 1}]) \subseteq U_i$.
So for every $i < n$ there is
$h_i \in H(K)\sprtm{U_i}$ such that $h_i(z_i) = z_{i + 1}$.
Clearly, $h_{n - 1} \scirc \ldots \scirc h_0$ is as required.

(b) The proof of Part (b) is trivial.
\smallskip\hfill\myqed

The following theorem is analogous to Theorem
\ref{metr-bldr-t6.15}. The proofs are essentially the same.

\begin{theorem}\label{metr-bldr-t6.22}
Let $X,Y \in K_{\srfs{IMX}}^{\calO}$ and
$\iso{\varphi}{H(\rfs{cl}(X))}{H(\rfs{cl}(Y))}$.
Then there is $\tau \in H(\rfs{cl}(X),\rfs{cl}(Y)))$
such that $\tau$ induces~$\varphi$.
\end{theorem}

\noindent
{\bf Proof }
Set $K = \rfs{cl}(X)$ and $M = \rfs{cl}(Y)$.
Then by Theorem \ref{t2.15}(c),
there is
$\eta \in H(\rfs{ENI}(K),\rfs{ENI}(M))$ such that $\eta$
induces $\varphi$.
So for every $g \in H(K)$,
$((g \nrestriction \rfs{ENI}(K))^{\eta})^{\srfs{cl}} = \varphi(g) \in
H(M)$.
We shall prove that $\eta^{\srfs{cl}} \in H(K,M)$.

{\bf Claim 1.} Let $x \in K - \rfs{ENI}(K)$
and $\vecx,\vecu \subseteq X$.
Suppose that
$\lim \vecx = \lim \vecu = x$ and that
$\eta(\vecx)$ and $\eta(\vecu)$ are convergent in~$M$.
Then $\lim \eta(\vecx) = \lim \eta(\vecu)$.
{\bf Proof }
Let $\vecy = \eta(\vecx)$, $\vecv = \eta(\vecu)$,
$y = \lim \vecy$, $v = \lim \vecv$, and suppose by contradiction that
$y \neq v$.
Obviously, $y,v \in \rfs{bd}(Y)$.
Let $r = \dgfrac{\norm{y - v}}{2}$.
We may assume that $\vecv \subseteq B(v,r)$ and that
$\vecy \cap B(v,r) = \emptyset$.
Let
$\fnn{h}{[0,1] \times \rfs{cl}(Y)}{\rfs{cl}(Y)}$
be an isotopy
as assured by the fact that
$v$ is isotopically movable with respect to $Y$,
and
such that for every $t \in [0,1]$, $\rfs{supp}(h_t) \subseteq B(v,r)$.

For every $t \in [0,1]$
let $u_{n,t} = \eta\inverse(h(t,v_n))$.
We prove the following fact.
$(*)$ For every $t \in [0,1]$,
$\limti{n} u_{n,t} = x$.
Let $t \in [0,1]$. Let $\barh = h_t \nrestriction \rfs{ENI}(M)$
and $\barg = \barh^{\eta\inverse}$.
Then $\barg \in \rfs{EXT}(\rfs{ENI}(K))$.
Clearly, $\barg \nrestriction \vecx = \rfs{Id}$
and so $\barg^{\srfs{cl}}(x) = x$.
Hence
$\limti{n} u_{n,t} = \limti{n} \barg(u_n) =
\barg(\limti{n} u_n) = \barg(x) = x$.
So $(*)$ is proved.

Let $L_n = h([0,1] \times \sngltn{v_n})$ and $K_n = \eta\inverse(L_n)$.
We prove that $\limti{n} K_n = x$.
Suppose by contradiction that this is not true.
Then there are $d > 0$, $\vect \subseteq [0,1]$
and a $\onetoonen$ sequence $\setm{n_i}{i \in \bbN}$ such that
for every $i \in \bbN$, $d(x,u_{n_i,t_i}) \geq d$.
We may assume that $\vect$ is convergent.
Let $t^* = \lim \vect$. Let $I_i$ be the closed interval whose
endpoints are $t_i$ and $t^*$
and\break
$J_i = h(I_i \times \sngltn{v_{n_i}})$.
Then $\limti{i} J_i = h(t^*,v)$.
Since for every $t \in [0,1]$, $h_t \nrestriction Y \in \rfs{EXT}(Y)$
and $v \in \rfs{bd}(Y)$, it follows that $h(t^*,v) \in \rfs{bd}(Y)$.
The fact that $v_{n_i} \in Y$ implies that $J_i \subseteq Y$.
Hence for every $i \in \bbN$, $h(t^*,v) \not\in J_i$.
We may thus assume that for any
$i \neq j$, $J_i \cap J_j = \emptyset$.

There is a sequence $\sngltn{V_i}_{i \in \bbN}$
of pairwise disjoint open sets such that for every $i \in \bbN$,
$J_i \subseteq V_i \subseteq \rfs{cl}(V_i) \subseteq Y \cap B(v,r)$
and $\limti{i} V_i = h(t^*,v)$.
Let $h_i \in \rfs{UC}(Y)\sprtm{V_i}$ be such that
$h_i(h(v_{n_i},t^*)) = h(v_{n_i},t_i)$
and
$\tildeh = \bcirc_{i \in \bbN} h_i$.
Then $\tildeh \in \rfs{UC}_0(Y)$.
Hence $\hath \eqdf \tildeh^{\srfs{eni}} \in \rfs{EXT}(\rfs{ENI}(M))$.
Let $\hatg = \hath^{\eta\inverse}$.
So $\hatg \in \rfs{EXT}(\rfs{ENI}(K))$.
Clearly, $\hatg \nrestriction \vecx = \rfs{Id}$
and hence $\hatg^{\srfs{cl}}(x) = x$.
Also, for every~$i \in \bbN$, \,$\hatg(u_{n_i,t^*}) = u_{n_i,t_i}$.
It follows from $(*)$ that \hbox{$\limti{i} u_{n_i,t^*} = x$}
and so
$$\limti{i} u_{n_i,t_i} =
\limti{i} \hatg(u_{n_i,t^*}) =
\hatg^{\srfs{cl}}(\limti{i} u_{n_i,t^*}) =
\hatg^{\srfs{cl}}(x) = x. 
$$
This contradicts the fact that $d(x,u_{n_i,t_i}) \geq d$,
so $\limti{n} K_n = x$.

Recall that $x \in K - \rfs{ENI}(K)$, and note that
$K_i = \eta\inverse(L_i) \subseteq \eta\inverse(Y) \subseteq
\rfs{ENI}(K)$.
So $x \not\in K_i$.
Hence there is an infinite set $\sigma \subseteq \bbN$
such that for any distinct $i,j \in \sigma$,
$K_i \cap K_j = \emptyset$.
There is a sequence $\setm{U_i}{i \in \sigma}$ of pairwise disjoint
sets such that $K_i \subseteq U_i \subseteq \rfs{ENI}(K)$,\break
$U_i$ is open in $\rfs{ENI}(K)$
and $\lim_{i \in \sigma} U_i = x$.
Let $\rho \subseteq \sigma$ be such that $\rho$ and $\sigma - \rho$
are infinite.

By Proposition \ref{metr-bldr-p6.21}(a),
for every $i \in \rho$ there is $g_i \in H(K)\sprtm{U_i}$
such that $g_i(u_i) = u_{i,1}$.
By Proposition \ref{metr-bldr-p6.21}(b),
$\hatg \eqdf \bcirc_{i \in \rho} g_i \in H(K)$.
Let
$\barg = \hatg \nrestriction \rfs{ENI}(K)$
and $\barh = \barg^{\eta}$.
Then $\barg^{\srfs{cl}} = \hatg \in H(K)$.
From the fact that $\eta$ induces $\varphi$
it follows that $\barh^{\srfs{cl}} \in H(M)$.

For every $i \in \rho$, $\barh(v_i) = h(v_i,1)$.
So $\lim_{i \in \rho} \barh(v_i) = h(v,1)$.
For every $i \in \sigma - \rho$, $\barh(v_i) = v_i$.
So $\lim_{i \in \sigma - \rho} \barh(v_i) = v$.
Recall that $h(v,1) \neq v$
and that $\limti{i} v_i = v$. So $\vecv$ is convergent and
$\barh(\vecv)$ is not convergent.
Hence $\barh^{\srfs{cl}} \not\in H(M)$.
A contradiction, so Claim 1 is proved.
\smallskip

Suppose by contradiction that $x \in K - \rfs{ENI}(K)$ and
$x \not\in \rfs{Dom}((\eta \nrestriction X)^{\srfs{cl}})$.
Recall that $Y \in K^{\calO}_{\srfs{IMX}}$ and hence $Y$ is JN.AC.
So by Proposition \ref{metr-bldr-p6.19}(c),
For every $x \in K - \rfs{ENI}(K)$
there is a sequence $\vecx \subseteq X$ converging to $x$
such that $\eta(\vecx) \subseteq Y$,
and $\eta(\vecx)$ is convergent in~$M$.
Set $y = \lim \eta(\vecx)$.
Obviously, $y \in \rfs{bd}(Y)$.
Since $x \not\in \rfs{Dom}((\eta \nrestriction X)^{\srfs{cl}})$,
there are a $\onetoonen$ sequence $\vecu \subseteq X$ and $d > 0$
such that $\lim \vecu = x$ and $d(\eta(\vecu),y) \geq d$.
Define $\vecv = \eta(\vecu)$.
Then by Claim 1, $\vecv$ does not have a convergent subsequence.
That is, $\vecv$ is completely discrete. 
Since $Y$ is JN.AC, there is a subsequence $\vecw$ of $\vecv$
such that $\vecw$ has a joining system.
Let\break
$\fortpl{\vecw}{w^*}{\setm{L_n}{n \in \bbN}}{\vecw\fprime}$
be a joining system for $\vecw$.
We may assume that $w^* \not\in \rfs{Rng}(\vecw)$.

It can be assumed that $w^* \neq y$.
For suppose that $w^* = y$. Let $r = d(\vecw,y)$.
Since $Y$ is BR.IS.MV and $y \in \rfs{bd}(Y)$,
there is $h \in \rfs{EXT}(Y)$ such that $\rfs{supp}(h) \subseteq B(y,r)$
and $h^{\srfs{cl}}(y) \neq y$.
So $h \nrestriction \vecw = \rfs{Id}$.
It follows that 
$\fortpl{\vecw}{h^{\srfs{cl}}(y)}{\setm{h(L_n)}{n \in \bbN}}
{h(\vecw\fprime)}$
is a joining system for $\vecw$,
and if we redefine $w^*$ to be $h^{\srfs{cl}}(y)$, then $w^* \neq y$.

Recall that $Y$ is JN.AC. So we may apply Lemma \ref{metr-bldr-p6.19}(b)
to $\eta\inverse$.
Recall also that
$\lim \eta\inverse(\vecw) = \lim \eta\inverse(\vecv) = x$.
Hence there is $\vecz \subseteq Y$
such that $\lim \vecz = w^*$ and $\lim \eta\inverse(\vecz) = x$.
The two sequences $\vecx$ and $\eta\inverse(\vecz)$
converge to $x$,
however, $\eta(\vecx)$ and $\eta(\eta\inverse(\vecz))$ are convergent,
but they do not converge to the same point.
This contradicts Claim 1,
so
$\rfs{Dom}((\eta \restriction X)^{\srfs{cl}}) \supseteq
K - \rfs{ENI}(K)$.
Since $\rfs{Dom}(\eta) = \rfs{ENI}(K)$,
we have that $\rfs{Dom}(\eta^{\srfs{cl}}) = K$.

We have shown that $\eta \in \rfs{EXT}(\rfs{ENI}(X),\rfs{ENI}(Y))$.
An identical argument shows that
$\eta\inverse \in \rfs{EXT}(\rfs{ENI}(Y),\rfs{ENI}(X))$.
Hence $\eta^{\srfs{cl}} \in H(K,M)$.
Since $\eta$ induces $\varphi$, $\eta^{\srfs{cl}}$ induces $\varphi$.
\rule{15pt}{0pt}\medskip\hfill\myqed

\subsection{Generalizations to manifolds and to nearly open sets.}
\label{ss6.5}
\label{ss6.5-generalization-to-manifolds-nearly-open}

The results of this chapter are true in two other settings,
which are more general than the present setting.
The proofs remain exactly the same as they were.

\begin{remark}\label{metr-bldr-r6.23}
\begin{rm}
(a) Let $Z$ be a subset of the normed space $E$.
$Z$ is a {\it nearly open set},
if $Z \subseteq \rfs{cl}^E(\rfs{int}^E(Z))$.
The results of this chapter can be extended to the
class of nearly open subsets of a normed space. Let
\newline\centerline{
$K^{\calN\calO}_{\srfs{NRM}} =
\setm{\pair{X}{Z}}{X \in K^{\calO}_{\srfs{NRM}}
\mbox{ and } X \subseteq Z \subseteq \rfs{cl}(X)}$.
}
Note that
$\setm{\pair{X}{\rfs{cl}(X)}}{X \in K^{\calO}_{\srfs{NRM}}} \subseteq
K^{\calN\calO}_{\srfs{NRM}}$.

   \index{nearly open set.
   $Z$ is nearly open, if $Z \subseteq \rfs{cl}(\rfs{int}(Z))$}

   \index{N@knonrm@@
   $K^{\calN\calO}_{\srfs{NRM}} =
   \setm{\pair{X}{Z}}{X \in K^{\calO}_{\srfs{NRM}}
   \mbox{ and } X \subseteq Z \subseteq \rfs{cl}(X)}$}

(b) The analogy with $K^{\calO}_{\srfs{NRM}}$ is as follows.
Let $\pair{X}{Z} \in K^{\calN\calO}_{\srfs{NRM}}$.
The group
$$
\rfs{EXT}^Z(X) =
\setm{h \kern1pt\nrestriction\kern1pt X\kern2pt}
{\kern2pt h \in H(Z) \mbox{ and } h(X) = X}
$$
is the analogue of
$\rfs{EXT}^E(X)$,
and the group $H(Z)$ is the analogue of
$H(\rfs{cl}(X))$.

(c) Suitable reformulations of Theorem \ref{metr-bldr-t6.5},
Corollary \ref{c4.33} and Theorems \ref{metr-bldr-t6.15},
\ref{metr-bldr-t6.17}, \ref{metr-bldr-t6.20} and \ref{metr-bldr-t6.22}
are true for $K^{\calN\calO}_{\srfs{NRM}}$.
\hfill\proofend
\end{rm}
\end{remark}

We demonstrate the generalization discussed in
Remark \ref{metr-bldr-r6.23} by describing the analogues of
Theorem~\ref{metr-bldr-t6.5}(b) and Theorem \ref{metr-bldr-t6.20}.
The faithful class captured by this generalization
contains $2^{2^{\aleph_0}}$ subsets of $\bbR^3$.

Let $K^{\calN\calO}_{\srfs{NMX}}$ be the class of all
$\pair{X}{Z} \in K^{\calN\calO}_{\srfs{NRM}}$ such that
$X$ is BR.LC.AC with respect to $Z$,
and $X$ is JN.AC with respect to $Z$.
Evidently, this is the analogue of $K^{\calO}_{\srfs{NMX}}$
defined in \ref{metr-bldr-t6.5}(b).
   \index{N@knonmx@@$K^{\calN\calO}_{\srfs{NMX}} =
   \setm{\pair{X}{Z} \in  K^{\calN\calO}_{\srfs{NRM}}}
   {X \mbox{ is BR.LC.AC and JN.AC with repect to } Z}$}
Let us first see that $K^{\calN\calO}_{\srfs{NMX}}$ is a large class.
Write $X = (0,1)^3$, that is, $X$ is an open cube in $\bbR^3$.
We construct sets $Z$ such that
$\pair{X}{Z} \in K^{\calN\calO}_{\srfs{NMX}}$,
and in fact, we show that
$\abs{\setm{Z}{\pair{X}{Z} \in K^{\calN\calO}_{\srfs{NMX}}}} =
2^{2^{\aleph_0}}$.
We skip the easy proof of Part (b) of the next example.

\begin{example}\label{metr-bldr-e6.24}
\begin{rm}
Let $X = (0,1)^3$.

(a) For $x,y \in \bbR$\kern1pt\ let $L_{x,y} = [(x,0,0),(x,y,0)]$.
Let
$\emptyset \neq A \subseteq [0,1]$ and $\fnn{\rho}{A}{[0,1)}$.\break
(We do not assume the $\rho$ is continuous).
Let
$Z_{\rho} = X \cup \bigcup_{x \in A} L_{x,\rho(x)}$.
Then $\pair{X}{Z_{\rho}} \in K^{\calN\calO}_{\srfs{NMX}}$.

(b) Let $F$ be a closed nonempty subset of $\rfs{bd}^{\bbR^3}(X)$.
Then
$\pair{X}{X \cup F} \in K^{\calN\calO}_{\srfs{NMX}}$.
\vspace{-2.3mm}
\end{rm}
\end{example}

\noindent
{\bf Proof }
(a)
Let $X,A,\rho$ and $Z$ be as above.
It is trivial that $X$ is BR.LC.AC with respect to $Z$.
We show that $X$ is JN.AC with respect to $Z$.
Let $\vecu = \sngltn{u_n}_{n \in \bbN} \subseteq X$ be
a completely discrete sequence with respect to $Z$.
It may be assumed that $\vecu$ is convergent in $\bbR^3$,
and we denote its limit by $\hatu$.
So $\hatu \in \rfs{cl}^{\bbR^3}(X) - Z$.
Write $u_n = (x_n,y_n,z_n)$ and $\hatu = (\hatx,\haty,\hatz)$.

{\bf Case 1 } Assume that $\hatz = 0$.
Suppose first that there is $a \in A$ such that
$\setm{n}{x_n = a}$ is infinite.
So we may assume that $x_n = a$ for every $n \in \bbN$.
It follows that for some $b > \rho(a)$, $\lim \vecu = (a,b,0)$.
Hence $\vecu$ has a subsequence $\vecv$
such that
$[v_m,(a,\rho(a),0)] \cap [v_m,(a,\rho(a),0)] = \sngltn{(a,\rho(a),0)}$
for any $m \neq n$.
Choose $w_n \in [v_n,(a,b,0))$ such that $\limti{n} w_n = (a,b,0)$
and define $L_n = [v_n,w_n]$. It is easy to see that
$\fortpl{\vecv}{(a,\rho(a),0)}{\vecL}{\sngltn{w_n}_{n \in \bbN}}$
is a joining system for $\vecv$

Suppose next that for every $a \in A$, $\setm{n}{x_n = a}$ is finite.
Choose any $a \in A$ and remove from $\vecu$ all $u_n$'s such that
$x_n = a$. Then $a \neq x_n$ for every $n \in \bbN$.
We may also assume that $z_0 < \dghalf$
and that $\sngltn{z_n}_{n \in \bbN}$ is strictly decreasing.
Let $y'_n = \max(1 - z_n,y_n)$, $u'_n = (x_n,y'_n,z_n)$
and $L^0_n = [u_n,u'_n]$.
We show that $\vecL^0 \eqdf \sngltn{L^0_n}_{n \in \bbN}$
is completely discrete with respect to $Z$.
Since $\sngltn{z_n}_{n \in \bbN}$ is $\onetoone$,
$\vecL^0$ is a pairwise disjoint sequence,
that is, $L^0_m \cap L^0_n = \emptyset$ for any $m \neq n$.
If $(x,y,z) \in \rfs{acc}^{\bbR^3}(\vecL^0)$,
then $x = \hatx$, $z = 0$ and $y \geq \haty$,
and since $(\hatx,\haty,0) \not\in Z$,
it follows that $(\hatx,y,0) \not\in Z$.
The sequence $\sngltn{y'_n}_{n \in \bbN}$ converges to $1$,
so we may assume that it is strictly increasing.
Let $v_n = (x_n,y'_n,\dghalf)$ and $L^1_n = [u'_n,v_n]$.
It is trivial that
$\vecL^1 \eqdf \sngltn{L^1_n}_{n \in \bbN}$
is a pairwise disjoint sequence.
If $(x,y,z) \in \rfs{acc}^{\bbR^3}(\vecL^1)$, then $y = 1$
and so $(x,y,z) \not\in Z$.
So $\vecL^1$ is completely discrete with respect to $Z$.
Suppose that $m < n$. Then $L^0_m \cap L^1_n = \emptyset$,
since the $y$-coordinate of any member of $L^0_m$ is $\leq y'_m$,
and the $y$-coordinate of any member of $L^1_n$ is equal to $y'_n$
which is $> y'_m$.
Similarly, $L^1_m \cap L^0_n = \emptyset$, since members of 
$L^1_m$ and $L^0_n$ differ in their $z$-coordinate.
We conclude that
$(L^0_m \cup L^1_m) \cap (L^0_n \cap L^1_n) = \emptyset$
for any $m \neq n$.

Let $w_n = (a,y'_n,\dghalf)$ and $L^2_n = [v_n,w_n]$.
The sequence $\vecL^2 \eqdf \sngltn{L^2_n}_{n \in \bbN}$
is a pairwise disjoint sequence, since members of $L^2_m$ and $L^2_n$
differ in their $y$-coordinate.
Also,\break
$L^2_n \cap (L^0_m \cup L^1_m) = \emptyset$ for any $m \neq n$.
This follows from the fact that the only point in
$L^0_m \cup L^1_m$ whose $z$-coordinate is $\dghalf$ is $v_m$
and $v_m \not\in L^2_n$.
The $y$-coordinate of any member of $\rfs{acc}^{\bbR^3}(\vecL^2)$
is $1$, so $\rfs{acc}^{\bbR^3}(\vecL^2) \cap Z = \emptyset$
and hence $\vecL^2$ is completely discrete with respect to $Z$.
Let $w^* = (a,\rho(a),0)$, choose $w'_n \in [w_n,w^*)$
such that $\limti{n} w'_n = w^*$ and define $L^3_n = [w_n,w'_n]$.
Clearly, $\vecL^3 \eqdf \sngltn{L^3_n}_{n \in \bbN}$
is a pairwise disjoint sequence.
Since $\limti{n} w_n = (a,1,0)$,
it follows that
$\rfs{acc}^{\bbR^3}(\vecL^3) = [w^*,(a,1,0)]$.
So $\rfs{acc}^Z(\vecL^3) = \sngltn{w^*}$.
It follows that for every $r > 0$,
$\setm{L^3_n - B(w^*,r)}{n \in \bbN}$ is completely discrete
with respect to $Z$.
Note that $w_m$ is the only point in $\bigcup_{i \leq 2} L^i_m$
whose $x$-coordinate is $a$.
So since for $n \neq m$, $w_m \not\in L^3_n$,
$L^3_n \cap (\bigcup_{i \leq 2} L^i_m) = \emptyset$.
Define $L_n = \bigcup_{i \leq 3} L^i_m$,
$\vecw\fprime = \sngltn{w'_n}_{n \in \bbN}$
and $\vecL = \sngltn{L_n}_{n \in \bbN}$.
It follows that $\vecL$ is a pairwise disjoint sequence and that
for every $r > 0$,
$\setm{L_n - B(w^*,r)}{n \in \bbN}$ is completely discrete
with respect to $Z$.
So
$\fortpl{\vecu}{w^*}{\vecL}{\vecw\fprime}$
is a joining system for $\vecu$.
\smallskip

The case that $\hatz \neq 0$ is divided into several subcases.
Their proofs are similar to the proof of Case 1, but they are simpler.
\hfill\myqed

\begin{theorem}\label{metr-bldr-t6.25}
For $\ell = 1,2$ let
$\pair{X_{\ell}}{Z_{\ell}} \in K^{\calN\calO}_{\srfs{NMX}}$.

\num{a} If $\iso{\varphi}{Z_1}{Z_2}$,
then there is $\tau \in H(Z_1,Z_2)$
such that $\tau$ induces $\varphi$.

\num{b} If $\iso{\varphi}{\rfs{EXT}^{Z_1}(X_1)}{\rfs{EXT}^{Z_2}(X_2)}$,
then there is $\tau \in \rfs{EXT}^{Z_1,Z_2}(X_1,X_2)$
such that $\tau$ induces $\varphi$.
\end{theorem}

\noindent
{\bf Proof } The proof of Part (a) is identical to the proof of
Theorem \ref{metr-bldr-t6.20}. The proof of Part~(b) is identical
to the proof of Theorem \ref{metr-bldr-t6.5}.
\smallskip\hfill\myqed

\begin{remark}\label{r4.6}\label{metr-bldr-e6.26}
\begin{rm}
The second generalization is motivated by the following example.
Let $E = \bbR \times S^{\sboldbbR^2}(0,1)$,
$Y = [0,1] \times S^{\sboldbbR^2}(0,1)$
and
$X = (0,1) \times S^{\sboldbbR^2}(0,1)$.
$X$ is a normed manifold. So its reconstruction from subgroups of
$H(X)$ is included in Theorem \ref{t2.15}(a).
The local $\itGamma$-continuity of conjugating homeomorphisms of
$X$ is proved in \ref{metr-bldr-t3.47}(a),
\ref{metr-bldr-c3.48}(a) and \ref{metr-bldr-c4.10}.
The space $Y$ however, is not covered by any of the above theorems
because it is not a normed manifold.
Also, $Y$ is not the closure of an open subset of a normed space.
So the theorems proved so far in Chapter \ref{s6} do not apply to $Y$.
However, $Y$ is a well-behaved space
and is very similar to the spaces which have been already dealt with.
\hfill\proofend
\end{rm}
\end{remark}

The above remark calls for the setting
in which $E$ is as a normed manifold,
$X$ is an open subset of $E$ and $Y = \rfs{cl}^E(X)$.
This setting will yield reconstruction results for $Y$.

\begin{defn}\label{metr-bldr-d6.27}
\begin{rm}
(a) Let $\trpl{X}{\itPhi}{d}$ be such that
$\pair{X}{\itPhi}$ is a normed manifold,
$\pair{X}{d}$ is a metric space,
and there is $K$ such that for every $\varphi \in \itPhi$,
$\varphi$ is $K$-bilipschitz.
Then $\trpl{X}{\itPhi}{d}$ is called a
{\it normed Lipschitz manifold}.

   \index{normed Lipschitz manifold}

(b) Let
$K^{\calO}_{\srfs{NLPM}} = \setm{Y}{Y \mbox{ is an open subset of a
normed Lipschitz manifold}}$.
\hfill\proofend

   \index{N@konlpm@@$K^{\calO}_{\srfs{NLPM}} =
          \setm{Y}{Y \mbox{ is an open subset of a normed Lipschitz
          manifold}}$}
\end{rm}
\end{defn}

Chapter \ref{s6} in its entirety can be proved for
$K^{\calO}_{\srfs{NLPM}}$.

\begin{theorem}\label{metr-bldr-t6.28}
In Definitions \ref{metr-bldr-d6.2.1}, \ref{metr-bldr-d6.8},
\ref{metr-bldr-d6.13}, \ref{metr-bldr-d6.16}
and in Remark \ref{metr-bldr-r6.23} change every mention of
$K^{\calO}_{\srfs{NRM}}$ to a mention of $K^{\calO}_{\srfs{NLPM}}$.
Then the variants obtained in this way from
Theorem \ref{metr-bldr-t6.5}
and Theorems \ref{metr-bldr-t6.10},
\ref{metr-bldr-t6.15}, \ref{metr-bldr-t6.17}, \ref{metr-bldr-t6.20},
\ref{metr-bldr-t6.22} and \ref{metr-bldr-t6.25} are true.
\end{theorem}

\noindent
{\bf Proof } The proofs of all the above theorems are identical to
the proofs of their counterparts.
\rule{0pt}{0pt}\hfill\myqed

\newpage

\noindent
\section{Groups which are not of the same type are not\\ isomorphic}
\label{s7}

In the previous chapters we considered several properties of
homeomorphisms. For instance, UC homeomorphisms, LUC homeomorphisms,
extendible homeomorphisms
and homeomorphisms which are uniformly continuous on every bounded
positively distanced set.
In this chapter we prove that for properties $\calP$ and $\calQ$ as
above, if $\calP(X) \cong \calQ(Y)$,
then either $\calP(X) = \calQ(X)$ or $\calP(Y) = \calQ(Y)$.
But before we deal with these questions,
we prove some additional facts about the group $\rfs{UC}(X)$.

\subsection{The group $\rfs{UC}(X)$ revisited.}
\label{ss7.1}

We have seen in Theorem \ref{t4.10}
that if $X,Y \in K^{\calO}_{\srfs{NRM}}$, $X$ is UD.AC
and $(\rfs{UC}(X))^{\tau} \subseteq \rfs{UC}(Y)$,
then $\tau$ is uniformly continuous.
We next reconsider the problem of deducing that $\tau^{-1}$ is
uniformly continuous from the fact that
$(\rfs{UC}(X))^{\tau} \subseteq \rfs{UC}(Y)$.
Recall that the implication
\begin{itemize}
\addtolength{\parskip}{-11pt}
\addtolength{\itemsep}{06pt}
\item[$(\dagger)$]\hfill
$(\rfs{UC}(X))^{\tau} \subseteq \rfs{UC}(Y)$ $\Rightarrow$
$\tau\inverse$ is uniformly continuous\hfill
\vspace{-05.7pt}
\end{itemize}
is not true for every $X,Y \in K^{\calO}_{\srfs{NRM}}$.
Counter-examples appear in \ref{metr-bldr-e5.8} and \ref{e4.34}(a).
Yet, $(\dagger)$~holds when $X$ and $Y$ are well-behaved.
Theorem \ref{t4.36} below deals with finite-dimensional spaces
for which $(\dagger)$ is true.
The infinite-dimensional case is considered in \ref{t4.41}.
The result of \kern1pt\ref{t4.41}
is needed in the proof of Corollary \ref{metr-bldr-c7.11}(d) and (e).

\begin{theorem}\label{t4.36}\label{metr-bldr-t7.2}
Let $X,Y \in K^{\calO}_{\srfs{NRM}}$.
Suppose that
$X$ is finite-dimensional and bounded,
$X$ is UD.AC, $\abs{\rfs{Cmp}(\rfs{bd}(X))} \leq \aleph_0$
and $(*)$ for every $C \in \rfs{Cmp}(\rfs{bd}(X))$,
distinct $x,y \in C$, and $z \in \rfs{bd}(X) - \dbltn{x}{y}$,
there is $f \in \rfs{UC}(X)$ such that either $f^{cl}(x) = y$
and $f^{cl}(z) = z$, or $f^{cl}(z) = y$ and $f^{cl}(x) = x$.
Suppose that for every $C \in \rfs{Cmp}(\rfs{bd}(Y))$, $\abs{C} > 1$.
Let $\tau \in H(X,Y)$ be such that
$(\rfs{UC}(X))^{\tau} \subseteq \rfs{UC}(Y)$.
Then $\tau^{-1}$ is uniformly continuous.
\end{theorem}

\noindent
{\bf Proof }
By Theorem \ref{t4.10}, $\tau$ is uniformly continuous, and hence
$\tau^{cl}$ maps $\rfs{cl}(X)$ onto $\rfs{cl}(Y)$.
It thus suffices to show that $\tau^{cl}$ is injective.
Suppose otherwise. For $x \in \rfs{bd}(X)$
let $C_x$ denote the connected component of $\rfs{bd}(X)$
containing $x$.
It follows from $(*)$ that if for some $z \neq x$,
$\tau^{cl}(x) = \tau^{cl}(z)$,
then for every $y \in C_x$, $\tau^{cl}(y) = \tau^{cl}(x)$.
The argument is as follows.
Suppose indeed that $z \neq x$, $\tau^{cl}(x) = \tau^{cl}(z)$
and $y \in C_x - \dbltn{x}{z}$.
Let $f \in \rfs{UC}(X)$ be as assured by $(*)$.
We assume first that $f^{\srfs{cl}}(x) = x$ and $f^{\srfs{cl}}(z) = y$,
Let $\vecx,\vecy \subseteq X$ converge respectively to
$x$ and $y$,
and let $\vecx\fprime = f\inverse(\vecx)$
and $\vecz = f\inverse(\vecy)$.
Then
$$
\tau^{\srfs{cl}}(y) =
\lim \tau(\vecy) = \lim f^{\tau} \scirc \tau \scirc f\inverse(\vecy) = 
\lim f^{\tau}(\tau(\vecz)) = (f^{\tau})^{\srfs{cl}}(\lim \tau(\vecz)) =
(f^{\tau})^{\srfs{cl}}(\tau^{\srfs{cl}}(z)).
$$
Similarly,
{\thickmuskip=3mu \medmuskip=1mu \thinmuskip=1mu 
$$
\tau^{\srfs{cl}}(x) =
\lim \tau(\vecx) = \lim f^{\tau} \scirc \tau \scirc f\inverse(\vecx) = 
\lim f^{\tau}(\tau(\vecx\fprime)) =
(f^{\tau})^{\srfs{cl}}(\lim \tau(\vecx\fprime)) =
(f^{\tau})^{\srfs{cl}}(\tau^{\srfs{cl}}(\lim \vecx\fprime)).
$$
}
\kern-3pt
Since $f^{\srfs{cl}}(x) = x$ and $\lim \vecx = x$, we have that
$\lim \vecx\fprime = x$. So
$$
\tau^{\srfs{cl}}(x) =
(f^{\tau})^{\srfs{cl}}(\tau^{\srfs{cl}}(\lim \vecx\fprime)) =
(f^{\tau})^{\srfs{cl}}(\tau^{\srfs{cl}}(x)) =
(f^{\tau})^{\srfs{cl}}(\tau^{\srfs{cl}}(z)) =
\tau^{\srfs{cl}}(y).
$$

The same argument applies to the case that
$f^{\srfs{cl}}(z) = z$ and $f^{\srfs{cl}}(x) = y$.
It follows that for any distinct $C,D \in \rfs{Cmp}(\rfs{bd}(X))$,
either $\tau^{cl}(C) = \tau^{cl}(D)$ and $\tau^{cl}(C)$ is a singleton,
or $\tau^{cl}(C) \cap \tau^{cl}(D) = \emptyset$.

Let $x$ and $y$ be distinct members of $\rfs{bd}(X)$ such that
$\tau^{cl}(x) =  \tau^{cl}(y)$, and $C$ be the component of
$\tau^{cl}(x)$ in \rfs{bd}(Y).
The family $\setm{\tau^{cl}(C_u) \cap C}{u \in \rfs{bd}(X)}$
is a partition of $C$ into more than one and at most countably
many closed sets.
This contradicts the theorem of Sierpi\'{n}ski that a continuum cannot
be partitioned into countably many nonempty closed sets.
See \cite{En} Theorem~6.1.27.\bigskip\hfill\myqed

We do not know whether in the above theorem, the requirement that
$\rfs{bd}(X)$ has at most countably many components can be dropped.
Here is an easy example of a bounded regular open subset
$X \subseteq \bbR^3$
such that $X$ is UD.AC, $X$ satisfies $(*)$ of
Theorem \ref{metr-bldr-t7.2}, 
every connected component of $\rfs{bd}(X)$ has cardinality $> 1$,
and $\rfs{bd}(X)$ has $2^{\aleph_0}$ connected components.
\smallskip

\begin{example}\label{editor-e7.2}
\begin{rm}
Let $C \subseteq [0,1]$ be the Cantor set.
Let $K = C \times \sngltn{1}$.
So $K \subseteq B^{\kern0.5pt\sboldbbR^2}(0,2)$
and $B^{\kern0.5pt\sboldbbR^2}(0,2) - K$ is connected.
Let
$A = \setm{a_n}{n \in \bbN} \subseteq B^{\kern0.5pt\sboldbbR^2}(0,2)$
be such that $\rfs{cl}(A) - A = K$, and every member of $A$ is an
isolated point in~$A$.
Let $r_n > 0$ and $D_n = \overB(a_n,r_n)$.
Assume that $D_n \subseteq B(0,2) \cap \setm{(x,y)}{x > 0}$
and $B(a_m,2r_m) \cap B(a_n,2r_n) = \emptyset$
for any $m \neq n$,
and that
$\rfs{cl}(\bigcup_{n \in \sboldbbN} D_n) -
\bigcup_{n \in \sboldbbN} D_n = K$.
Let $U = B(0,3) - \rfs{cl}(\bigcup_{n \in \sboldbbN} D_n)$.
Let $X \subseteq \bbR^3$ be the set obtained by rotating $U$
about the $x$-axis.
Note that if $x,y \in U$, then there is an arc $L \subseteq U$
connecting $x$ and $y$
such that $\rfs{lngth}(L) \leq 2\pi \cdot \norm{x - y}$.
It follows easily that $X$ is as required.
\hfill\proofend
\end{rm}
\end{example}

We next deal with infinite-dimensional open sets for which
the fact that $(\rfs{UC}(X))^{\tau} \subseteq \rfs{UC}(Y)$ implies that
$\tau^{-1}$ is uniformly continuous.

\begin{defn}\label{d4.37}\label{metr-bldr-d7.3}
\begin{rm}
(a) For $A \subseteq X$
define $\Delta^{X,E}(A) = \sup_{a \in A} d(a,E - X)$.
As usual, we abbreviate $\Delta^{X,E}(A)$ by $\Delta(A)$.

   \index{N@d06@@$\Delta(A)$. $\sup_{a \in A} d(a,E - X)$}

(b) Let $h \in H(X)$. We say that $h$ is {\it strongly extendible},
if for every
$\varepsilon > 0$
there is $\tilde{h} \in H(E)$
such that $\tilde{h}$ extends $h$ and
$\rfs{supp}(\tildeh) \subseteq B(\rfs{supp}(h),\varepsilon)$.
Define
$\rfs{UC}_{\rme}(X) \eqdf
\setm{h \in \rfs{UC}(X)}{h \mbox{ is strongly extendible}}$.

   \index{strongly extendible}
   \index{N@uce@@$\rfs{UC}_{\rme}(X) =
          \setm{h \in \rfs{UC}(X)}{h \mbox{ is strongly extendible}}$}

(c) A {\it simple arc} is a space homeomorphic to $[0,1]$.
For a simple arc $L$ and $x,y \in L$
let $[x,y]^L$ denote the subarc of $L$ whose endpoints are $x$ and $y$.
Let $\alpha \in \rfs{MBC}$
and $\fnn{\eta}{(0,\infty)}{(0,\infty)}$ be such that
$\eta$ is monotonic and $\limtz{t} \eta(t) = 0$.
Let $X$ be a metric space and $L \subseteq X$
be a simple arc.
We say that $L$ is an {\it $\pair{\alpha}{\eta}$-track}
if for every $x,y \in L$
there is $h \in \rfs{UC}(X)$ such that $h$ is $\alpha$-bicontinuous,
$h(x) = y$
and $\fs{supp}(h) \subseteq B([x,y]^L,r)$,
where $r = \eta(\rfs{diam}([x,y]^L)$.
   \index{track. $\pair{\alpha}{\eta}$-track}
If in the above definition we require that $h \in \rfs{UC}_{\rme}(X)$,
then $L$ is called an {\it $\pair{\alpha}{\eta}$-e-track}.
   \index{e-track. $\pair{\alpha}{\eta}$-e-track}

(d) We define the notion of a track system for $\vecx$.
Let $\vecx \subseteq X$ be a completely discrete sequence,
$y^* \in \rfs{bd}(X)$, $\vecy \subseteq X$
and $\vecL = \setm{L_n}{n \in \bbN}$ be a sequence of simple arcs
such that $\lim \vecy = y^*$, $L_n \subseteq X$,
$L_n$ connects $x_n$ with $y_n$
and $\bigcup_{n \in \bbN} L_n$ is bounded.
Assume that
\begin{itemize}
\addtolength{\parskip}{-11pt}
\addtolength{\itemsep}{06pt}
\item[(1)] 
there are $\alpha$ and $\eta$ such that $L_n$ is an
$\pair{\alpha}{\eta}$-track for every $n \in \bbN$,
\item[(2)] 
there are $\beta \in \rfs{MC}$ and for every $n$ a parametrization
$\gamma_n$ of $L_n$ such that $\rfs{Dom}(\gamma_n) = [0,1]$,
$\gamma_n$ is $\beta$-UC for every $n \in \bbN$ and
$\gamma_n(0) = y_n$ and $\gamma_n(1) = x_n$.
\vspace{-05.7pt}
\end{itemize}
Then $T = \fortpl{\vecx}{y^*}{\vecL}{\vecy}$
is called a {\it track system for~$\vecx$},
and $\gamma_n$ is called a {\it legal parametrization} of $L_n$
in~$T$.
Note that Clause (2) just means that $\setm{\gamma_n}{n \in \bbN}$ is
equicontinuous.
If in Clause (1) we require that $L_n$ be an e-track,
then $T$ is called an {\it e-track system}.

   \index{track system}
   \index{e-track system}
   \index{legal parametrization}

Let $T = \fortpl{\vecx}{y^*}{\vecL}{\vecy}$
be a track system.
If for every $r > 0$,
$\setm{L_n - B(y^*,r)}{n \in \bbN}$ is completely discrete,
then $T$ is called a {\it completely discrete track system}.
If for every $r > 0$, $\setm{L_n - B(y^*,r)}{n \in \bbN}$
is spaced, then $T$ is called a {\it spaced track system}.

   \index{completely discrete track system}
   \index{spaced track system}

(f)
$X$ is {\it jointly track connected (JN.TC)},
if for every completely discrete bounded sequence $\vecx \subseteq X$:
if $\limti{n} \delta(x_n) = 0$,
then $\vecx$ has a subsequence $\vecy$
such that $\vecy$ has a track system.
$X$ is {\it jointly e-track connected (JN.ETC)},
if the above subsequence $\vecy$ is required to have an e-track system.

   \index{jntc@@JN.TC}
   \index{jnetc@@JN.ETC}

\end{rm}
\end{defn}

\begin{remark}\label{r4.38}\label{metr-bldr-r7.4}
\begin{rm}
We explain the notion of a track system by an example.
Let $X$ be the unit ball of the Hilbert space $\ell_2$
and $S$ be the unit sphere. Let $\vecx$ be a completely discrete
sequence in $X$ such that $\delta(\vecx) = 0$. We construct a track
system for a subsequence of $\vecx$.
Let $e_0 = (1,0,0,\ldots)$.
Take a subsequence $\vecy$ of $\vecx$ such that
$\sngltn{e_0} \cup \rfs{Rng}(\vecy)$ is an independent set.
For $n \in \bbN$ let $z_n = \norm{y_n} e_0$,
$S_n = S(0,\norm{y_n}) \cap \rfs{span}(\dbltn{y_n}{e_0})$
and $L_n$ be any of the two subarcs of $S_n$ connecting
$y_n$ with $z_n$.
Then
$T =
\fortpl{\vecy}{e_0}{\sngltn{L_n}_{n \in \bbN}}
{\sngltn{z_n}_{n \in \bbN}}$
is a track system for $\vecy$. Indeed, $T$ is an e-track system.

The property JN.ETC is needed in the proof that
$\rfs{UC}(X) \not\cong \rfs{EXT}(X)$.
\end{rm}
\end{remark}

\begin{prop}\label{metr-bldr-p7.5}
\num{a} Let $\setm{h_n}{n \in \bbN} \subseteq \rfs{UC}_{\rme}(X)$,
and suppose that $\setm{\rfs{supp}(h_n)}{n \in \bbN}$ is spaced.
Then $\bcirc_{n \in \sboldbbN} h_n \in \rfs{UC}_{\rme}(X)$.

\num{b} Let $x,y \in E$ be such that $\norm{x} = \norm{y}$ and
$\norm{x - y} = d > 0$. Let $L = \setm{tx}{t \geq 0}$.
Then $d(y,L) \geq \dgfrac{d}{2}$.

\num{c}
If $T = \fortpl{\vecx}{y^*}{\vecL}{\vecy}$ a track system,
then the following hold.
\begin{itemize}
\addtolength{\parskip}{-11pt}
\addtolength{\itemsep}{06pt}
\item[\num{i}] 
For every $t \in (0,1)$,
$T_t \eqdf \fortpl{\sngltn{\gamma_n(t)}_{n \in \bbN}}{y^*}
{\sngltn{\gamma_n([0,t])}_{n \in \bbN}}{\vecy}$
is a track system,
and if $T$ is completely discrete, so is $T_t$.
\item[\num{ii}] 
$\limti{n} \Delta(L_n) = 0$.
\vspace{-05.7pt}
\end{itemize}

\num{d} Let
$\fortpl{\vecx}{y^*}{\vecL}{\vecy}$
be a completely discrete track system.
Then there is an infinite $\sigma \subseteq \bbN$ such that
$\fortpl{\vecx \nrestriction \sigma}{y^*}{\vecL \nrestriction \sigma}
{\vecy \nrestriction \sigma}$
is a spaced track system.

\num{e} Let $T = \fortpl{\vecx}{y^*}{\vecL}{\vecy}$ be a track system.
Let $\gamma_n$ be legal parametrization of $L_n$\break
in $T$.
Then there are $t \in [0,1)$, $z^* \in \rfs{bd}(X)$
and an infinite $\sigma \subseteq \bbN$
such that\break
$\fortpl{\vecx \nrestriction \sigma}{z^*}
{\setm{\gamma_n([t,1])}{n \in \sigma}}
{\setm{\gamma_n(t)}{n \in \sigma}}$
is a spaced track system.

\num{f} Let $T = \fortpl{\vecx}{y^*}{\vecL}{\vecy}$ be a
completely discrete track system
and $C \in \rfs{Cmp}(\rfs{bd}(X))$ be such that $d(\vecx,C) = 0$.
Then $y^* \in C$.

\num{g}
{\thickmuskip=2mu \medmuskip=1mu \thinmuskip=1mu 
Let $T = \fortpl{\vecx}{y^*}{\vecL}{\vecy}$ be a track system,
$h \in \rfs{UC}(X)$ and
$T' \eqdf \fortpl{h(\vecx)}{h^{\srfs{cl}}(y^*)}{h(\vecL)}{h(\vecy)}$.
}
Then $T'$ is a track system.
\end{prop}

\noindent
{\bf Proof }
(a) The proof is trivial and is left to the reader.

(b) We may assume that $\norm{x} = 1$.
Let $tx \in L$. If $|1 - t| \leq \dgfrac{d}{2}$, then use the
triangle with vertices $x$, $tx$ and $y$ to conclude that
$\norm{y - tx} \geq \norm{y - x} - \norm{x - tx}
\geq \dgfrac{d}{2}$;
and if $|1 - t| \geq \dgfrac{d}{2}$, then use the triangle
with vertices $0$, $tx$ and $y$ to conclude that
$\norm{y - tx} \geq |\kern 1pt\norm{y - 0} - \norm{tx - 0}\kern 1pt| =
|1 - t| \geq \dgfrac{d}{2}$.

(c) The first part of (c) follows from the definition of a track system.
To prove the second part,
suppose by way of contradiction that for some $d > 0$,
$\setm{n}{\Delta(L_n) > d}$ is infinite.
Let $\alpha$ and $\eta$
be as assured by the fact that $T$ is a track system.
Since $\lim \vecy = y^* \in \rfs{bd}(X)$,
there is $n$ such that $\alpha(\delta(y_n)) < d$ and $\Delta(L_n) > d$.
Choose $z \in L_n$ such that $\delta(z) > d$
and $w \in \rfs{bd}(X)$ such that $\alpha(\norm{y_n - w} < d$.
Since $L_n$ is an $\pair{\alpha}{\eta}$-track,
there is $h \in H(X)$
such that $h$ is $\alpha$-bicontinuous and $h(y_n) = z$.
Then
$$
\norm{h(y_n) - h(w)} = \norm{z - h(w)} \geq d(z,\rfs{bd}(X)) >
d > \alpha(\norm{y_n - w}),
$$
and this contradicts the $\alpha$-continuity of $h$.
\smallskip

(d)
For every $r > 0$,
$\setm{L_i - B(y^*,r)}{i \in \bbN}$ is completely discrete.
So by Proposition~\ref{metr-bldr-p5.25},
for every $r > 0$ and an infinite $\eta \subseteq \bbN$
there is an infinite $\nu \subseteq \eta$ such that
$\setm{L_i - B(y^*,r)}{i \in \nu}$ is spaced.
We define by induction $\rho_n \subseteq \bbN$.
Let $\rho_0 = \bbN$. For every $n \in \bbN$ let
$\rho_{n + 1}$ be an infinite subset of $\rho_n$ such that
$\setm{L_i - B(z^*,\frac{1}{n + 1})}{i \in \rho_{n + 1}}$
is spaced.
Let\break
$\sigma = \setm{\min(\rho_n \cap \bbN^{\geq n})}{n \in \bbN}$.
It is easy to see that for every $r > 0$,
$\setm{L_i - B(z^*,r)}{i \in \sigma}$ is spaced.
So
$\fortpl{\vecx \nrestriction \sigma}{y^*}
{\vecL \nrestriction \sigma}
{\vecy \nrestriction \sigma}$
is a spaced track system.

(e) For every infinite $\eta \subseteq \bbN$ and $t \in [0,1]$
let $A[\eta,t] = \setm{\gamma_n(t)}{n \in \eta}$.
Let\break
$s_{\eta} =
\sup(\setm{t}{A[\eta,t] \mbox{ is not completely discrete}})$.
Let $\rho \subseteq \bbN$ be an infinite set such that for every
infinite $\eta \subseteq \rho$, $s_{\eta} = s_{\rho}$.
Denote $s = s_{\rho}$.
Suppose by contradiction that $A[\rho,s]$
does not contain a Cauchy  sequence.
Then for some infinite $\eta \subseteq \rho$ and $d > 0$,
$A[\eta,s]$ is $d$-spaced.
There is $\varepsilon > 0$
such that for every $t > s - \varepsilon$,
$A[\eta,t]$ is spaced.
The existence of $\varepsilon$ follows from the equicontinuity of
$\setm{\gamma_n}{n \in \bbN}$,
that is, from the existence of $\beta$ appearing in Clause~(2)
of the definition of a track system.
So $s_{\eta} < s$. A contradiction.
So $A[\rho,s]$ contains a Cauchy sequence.
We may thus assume that $A[\rho,s]$ is a Cauchy sequence.
Let $z^* = \lim A[\rho,s]$.

Let $J_i = \gamma_i([s,1])$.
We show that there are no $r > 0$, an infinite $\eta \subseteq \rho$
and\break
$\vecu \in \prod_{i \in \eta} (J_i - B(z^*,r))$ such that
$\vecu$ is a Cauchy sequence.
Suppose otherwise.
Let $t_i \in [s,1]$ be such that $u_i = \gamma_i(t_i)$.
We may assume that $ \vect = \setm{t_i}{i \in \eta}$
is a Cauchy sequence.
Let $t^* = \lim \vect$.
Since $\rfs{Rng}(\vecu) \cap B(z^*,r) = \emptyset$,
it follows that $t^* \neq s$,
and since $\lim_{i \in \eta} d(\gamma_i(t_i),\gamma_i(t^*)) = 0$,
we have that $\setm{\gamma_i(t^*)}{i \in \eta}$ is a Cauchy sequence.
That is, $s_{\eta} > s$, a contradiction.
We have shown that
$\fortpl{\setm{x_n}{n \in \rho}}{z^*}
{\setm{\gamma_n([s,1])}{n \in \rho}}{A[\rho,s]}$
is a completely discrete track system.
By Part (d), there is an infinite $\sigma \subseteq \rho$ such that
$\fortpl{\setm{x_n}{n \in \sigma}}{z^*}
{\setm{\gamma_n([s,1])}{n \in \sigma}}{A[\sigma,s]}$
is a spaced track system.

(f) Suppose by contradiction that $y^* \not\in C$.
By (d), we may assume that $T$ is a spaced track system.
Let $\alpha,\eta$ be as assured by the fact that
$T$ is a track system.
Clearly, $a \eqdf d(y^*,C) > 0$. Choose $u \in C$,
and for every $n \in \bbN$ choose
$z_n \in (B(y^*,\dgfrac{a}{2}) - B(y^*,\dgfrac{a}{4})) \cap L_n$
and set $J_n = [x_n,z_n]^{L_n}$.
Then $b_1 \eqdf d(u,\bigcup_{n \in \bbN} J_n) > 0$,
and there is $b_2$ such that $\setm{J_n}{n \in \bbN}$ is $b_2$-spaced.
Set $b = \dgfrac{\min(b_1,b_2)}{3}$,
and let $c > 0$ be such that $c + \eta(c) < b$.
From the equicontinuity of $\sngltn{\gamma_n}_{n \in \bbN}$ it follows
that there is $k \in \bbN$ and $\setm{z_{n,i}}{n \in \bbN,\ i \leq k}$
such that for every $n \in \bbN$, $z_{n,0} = x_n$, $z_{n,k} = z_n$
and $z_{n,i} \in L_n$,
and $\rfs{diam}([z_{n,i},z_{n,i + 1}]^{L_n}) < c$
for every $i < k$.
So for every $n \in \bbN$ and $i < k$
there is $h_{n,i} \in \rfs{UC}(X)$ such that
$h_{n,i}$ is $\alpha$-bicontinuous, $h_{n,i}(z_{n,i}) = z_{n,i + 1}$
and $\rfs{supp}(h_{n,i}) \subseteq B([z_{n,i},z_{n,i + 1}]^{L_n},c)$.
Let $h_n = \bcirc_{i < k} h_{n,i}$.
Clearly, $h_n \in \rfs{UC}(X)$, and it is easily seen
that $\setm{\rfs{supp}(h_n)}{n \in \bbN}$ is $\frac{b_2}{3}$-spaced
and $d(u,\rfs{supp}(h_n)) > \dgfrac{b_1}{2} > 0$.
It follows that $h \eqdf \bcirc_{i < k} h_{n,i} \in \rfs{UC}(X)$,
$h^{\srfs{cl}}(u) = u$ and $h(x_n) = z_n$ for every $n \in \bbN$.
Since $h(u) = u$, it follows that $h(C) = C$.
However, $d(\vecx,C) = 0$
and $d(h(\vecx),h(C)) = d(\vecz,C) > \dgfrac{a}{2} > 0$.
This contradicts the fact that $h$ is uniformly continuous.

(g) By Proposition \ref{metr-bldr-p4.3}(c),
there is $\gamma \in \rfs{MBC}$ such that $h$ is $\gamma$-bicontinuous.
Let $\alpha$, $\eta$ and $\beta$ be as in the definition of
a track system.
Define $\alpha' = \gamma \scirc \alpha \scirc \gamma$,
$\eta' = \gamma \scirc \eta \scirc \gamma$
and $\beta' = \gamma \scirc \beta$.
Then $\alpha',\eta'$ and $\beta'$ demonstrate that $T'$ is
a track system.
\rule{7pt}{0pt}\hfill\myqed

\begin{prop}\label{p4.40}\label{metr-bldr-p7.6}
Let $Z$ be a metric space, and
$\setm{F_n}{n \in \bbN}$ and $\setm{K_n}{n \in \bbN}$ be sequences of
compact subsets of $Z$ such that:
\num{i} $\setm{F_n}{n \in \bbN}$ is spaced;
\num{ii} for every $\varepsilon > 0$ there is
$\ell_{\varepsilon} \in \bbN$ such that for
every $n \in \bbN$ and a subset $A \subseteq K_n$,
if $|A| \geq \ell_{\varepsilon}$,
then there are distinct $x,y \in A$
such that $d(x,y) < \varepsilon$; and
\num{iii} $\inf(\setm{d(F_n,K_n)}{n \in \bbN}) > 0$.
Then there is an infinite $\sigma \subseteq \bbN$ such that
$d(\bigcup \setm{F_n}{n \in \sigma},
\bigcup \setm{K_n}{n \in \sigma}) > 0$.
\end{prop}

\noindent
{\bf Proof }
Write $\bbN^+ = \setm{n \in \bbN}{n > 0}$.
   \index{N@n00@@$\bbN^+ = \setm{n \in \bbN}{n > 0}$}
We define by induction on $i \in \bbN^+$ a sequence of infinite subsets
of~$\bbN$, $\sigma_0 \supseteq \sigma_1 \supseteq
\kern-1pt\dots$\kern2.5pt.
Let $\sigma_0 = \bbN$. Suppose that $\sigma_i$ has been defined.
\hbox{We colour the increasing} pairs $\pair{m}{n}$
of members of $\sigma_i$ in four colours,
according to whether $d(F_m,K_n) < \dgfrac{1}{i}$ or not,
and according to whether $d(K_m,F_n) < \dgfrac{1}{i}$ or not.
By Ramsey~Theorem, there is a monochromatic infinite
$\sigma_{i + 1} \subseteq \sigma_i$.
If there is $i \in \bbN^+$ such that for every distinct
$m,n \in \sigma_i$,
$d(F_m,K_n) \geq \dgfrac{1}{i}$ and
\hbox{$d(K_m,F_n) \geq \dgfrac{1}{i}$,} then $\sigma \eqdf \sigma_i$
is as required.
Otherwise, for every $i \in \bbN$
either
(1) for every $m < n$ in $\sigma_i$, $d(F_m,K_n) < \dgfrac{1}{i}$,
or
(2) for every $m < n$ in $\sigma_i$, $d(K_m,F_n) < \dgfrac{1}{i}$.

Let $i \in \bbN$ and $\ell = \ell_{\dgfrac{1}{i}}$ be as assured by
clause (ii).
Let $k_0 < \ldots < k_{\ell}$ be members of $\sigma_i$.
Suppose that Case 1 occurs. For every $j < \ell$ let $x_j \in F_j$
and $y_j \in K_{\ell}$ be such that $d(x_j,y_j) < \dgfrac{1}{i}$.
Hence for some $j_1 < j_2 < \ell$, $d(y_{j_1},y_{j_2}) < \dgfrac{1}{i}$.
So $d(F_{j_1},F_{j_2}) < \dgfrac{3}{i}$.
The same argument is repeated in Case 2.
Hence for every $i \in \bbN^+$ there are distinct $j_1$ and $j_2$
such that $d(F_{j_1},F_{j_2}) < \dgfrac{3}{i}$,
contradicting the fact that $\setm{F_n}{n \in \bbN}$ is spaced.
\medskip\hfill\myqed

The properties that $X$ is required to fulfill in the next theorem
are quite restrictive.
However, they are shared by ``well-behaved'' open sets.
For example, if $X = B - \bigcup_{i < k} \overB_i$,
where $B$ is an open ball and $\fsetn{\overB_0}{\overB_{k - 1}}$
is a pairwise disjoint family of closed balls contained in $B$,
then $X$ fulfills the requirements of the theorem.
Part (b) of the theorem is a slight modification of its first part.
This modification is needed in the proof that $\rfs{UC}(X)$ and
$\rfs{EXT}(X)$ are not isomorphic unless they coincide.

\begin{theorem}\label{t4.41}\label{metr-bldr-t7.7}
\num{a} Let $X \in K^{\calO}_{\srfs{BNC}}$.
Suppose that the following hold.
\begin{itemize}
\addtolength{\parskip}{-11pt}
\addtolength{\itemsep}{06pt}
\item[\rm (1)] 
$X$ is bounded
and $X$ is UD.AC,
\item[\rm (2)] 
$\rfs{bd}(X)$ has finitely many connected components,
\item[\rm (3)] 
if $C \in \rfs{Cmp}(\rfs{bd}(X))$, $x,y \in C$ are distinct
and $z \in \rfs{bd}(X) - \dbltn{x}{y}$,
then there is $f \in \rfs{UC}(X)$ such that either $f^{cl}(x) = y$
and $f^{cl}(z) = z$, or $f^{cl}(z) = y$ and $f^{cl}(x) = x$,
\item[\rm (4)] 
$X$ is JN.TC,
\vspace{-05.7pt}
\end{itemize}
Let $Y \in  K^{\calO}_{\srfs{BNC}}$ and assume that
\begin{itemize}
\addtolength{\parskip}{-11pt}
\addtolength{\itemsep}{06pt}
\item[\rm (5)] 
if $C$ is a component of $\rfs{bd}(Y)$, then $\abs{C} > 1$.
\vspace{-05.7pt}
\end{itemize}
Let $\tau \in H(X,Y)$ be such that
$(\rfs{UC}(X))^{\tau} \subseteq \rfs{UC}(Y)$.
Then $\tau^{-1}$ is uniformly continuous.

\num{b} Modify Clause \num{3} of Part \num{a} by requiring that
$f \in \rfs{UC}_{\rme}(X)$,
and modify Clause \num{4} by requiring that $X$ is JN.ETC.
Let $\tau \in H(X,Y)$ be such that
$(\rfs{UC}_{\rme}(X))^{\tau} \subseteq \rfs{UC}(Y)$.
Then $\tau^{-1}$ is uniformly continuous.
\end{theorem}

\noindent
{\bf Proof }
The proofs of (a) and (b) are identical. We prove Part (a).
Recall that $X$ and $Y$ are subsets of the Banach spaces $E$
and $F$ respectively.
By Theorem \ref{metr-bldr-t5.6},
$\tau$ is uniformly continuous.
\smallskip

{\bf Claim 1.} Let $\vecx \subseteq X$ be a completely discrete
sequence such that $\tau(\vecx)$ is a Cauchy sequence.
Then there is a sequence $\vecx\fprime \subseteq X$
such that $\limti{n} \delta(x'_n) = 0$,
$\vecx\fprime$ is completely discrete, and
$\limti{n} \tau(\vecx\fprime) = \limti{n} \tau(\vecx)$.
{\bf Proof }
If $\delta(\vecx) = 0$, then we take $\vecx\fprime$ to be a subsequence
of $\vecx$ such that $\limti{n} \delta(x'_n) = 0$.
Suppose otherwise. Since $X \in K^{\calO}_{\srfs{BNC}}$,
we may assume that for some $d > 0$, $\vecx$ is $d$-spaced,
and since $X$ is bounded,
we may also assume that for every $n \in \bbN^+$,
$d(x_n,x_0) \leq d + \dgfrac{d}{8}$.
Without loss of generality, $x_0 = 0$.
For every $n \in \bbN^+$ let
$t_n = \min(\setm{t > 1}{t x_n \in \rfs{bd}(X)})$, $y_n = t_n x_n$,
$L_n = [x_n,y_n]$ and
$\gamma_n(t) = x_n + t(y_n - x_n)$, $t \in [0,1]$.
If $m \neq n$, then
$$
\norm{d \ncdot \frac{x_m}{\norm{x_m}} - d \ncdot \frac{x_n}{\norm{x_n}}}
\geq
\norm{x_m - x_n} - \norm{x_m - d \ncdot \frac{x_m}{\norm{x_m}}} -
\norm{x_n - d \ncdot \frac{x_n}{\norm{x_n}}} \geq \frac{3d}{4}.
$$
Hence by Proposition \ref{metr-bldr-p7.5}(b),
$d(L_m,L_n) \geq \dgfrac{3d}{8}$.

Define $\eta(t) = \delta(\setm{\gamma_n(t)}{n \in \bbN^+})$.
Since $\setm{\norm{x_n - y_n}}{n \in \bbN}$ is bounded,
$\eta$ is continuous.
Also, $\eta(1) = 0$.
Let $s = \min(\eta\inverse(0))$.
We may assume that for every $n \in \bbN^+$,
\hbox{$\delta(\gamma_n(s)) < \dgfrac{1}{n}$.}
It follows that for every $t \in (0,s)$,
the family $\setm{\gamma_n([0,t])}{n \in \bbN^+}$ is spaced,
and $\delta(\bigcup \setm{\gamma_n([0,t])}{n \in \bbN^+}) > 0$.
Also, since $X$ is bounded, $\setm{d(x_n,\gamma_n(t))}{n \in \bbN^+}$
is bounded.
So for every $t < s$ there is $h_t \in \rfs{UC}(X)$
such that for every $n \in \bbN^+$,
$h_t(x_{2n}) = \gamma_{2n}(t)$ and $h_t(x_{2n - 1}) = x_{2n - 1}$.
Let $z^* = \lim \tau(\vecx)$.
Let $t \in (0,s)$.
Clearly,
$\tau(\setm{\gamma_{2n}(t)}{n \in \bbN^+} \cup
\setm{x_{2n-1}}{n \in \bbN^+}) = (h_t)^{\tau}(\vecx)$,
and since 
$(h_t)^{\tau}\in \rfs{UC}(Y)$ and $\tau(\vecx)$ is a Cauchy sequence,
$\tau(\setm{\gamma_{2n}(t)}{n \in \bbN^+} \cup
\setm{x_{2n-1}}{n \in \bbN^+})$ is a Cauchy sequence.
Denote this sequence by $\vecu$.
$\tau(\setm{x_{2n-1}}{n \in \bbN^+})$ is a subsequence of $\vecu$
converging to $z^*$. So $\vecu$ converges to $z^*$, and hence
$\tau(\setm{\gamma_{2n}(t)}{n \in \bbN^+})$ converges to $z^*$.
Let $\vecs \subseteq (0,s)$ be a sequence converging to $s$.
For every $n \in \bbN^+$ let $k_n \geq n$ be such that
$d(\tau(\gamma_{2k_n}(s_n)),z^*) < \dgfrac{1}{n}$.
Let $x'_n = \gamma_{2k_n}(s_n)$.
So $\lim \tau(\vecx\fprime) = z^*$,
$\limti{n} \delta(x'_n) = 0$ and $\vecx\fprime$ is spaced.
Claim~1 is thus proved.\smallskip

{\bf Claim 2.}
Let $T = \fortpl{\vecy}{y^*}{\setm{L_n}{n \in \bbN}}{\vecz}$ be a
completely discrete track system in $X$,
and suppose that $\lim \tau(\vecy) = w^*$.
Then $\tau^{\srfs{cl}}(y^*) = w^*$.
{\bf Proof }
Suppose by contradiction that $\tau^{cl}(y^*) \neq w^*$.
Let $\gamma_n$ be a legal parametrization of $L_n$,
and $\beta \in \rfs{MC}$ be such that for every
$t_1,t_2 \in [0,1]$ and $n \in \bbN$,
$\gamma_n(t_1) - \gamma_n(t_2) \leq \beta(\abs{t_1 - t_2})$.

We now follow the proof of Lemma \ref{l-gamma.5}. 
For every infinite $\sigma \subseteq \bbN$ and $t \in [0,1]$ let
$A[\sigma,t] = \setm{\gamma_n(t)}{n \in \sigma}$ and
$s_{\sigma} = \inf(\setm{t \in [0,1]}{\tau(A[\sigma,t])
\mbox{ converges to } w^*})$.
Since $\tau^{cl}(y^*) \neq w^*$, there is $U \in \rfs{Nbr}(y^*)$
such that $d(w^*,\tau(U \cap X)) > 0$. Thus there is $t_0 > 0$
such that for every $t < t_0$, $d(w^*, \tau(A[\bbN,t])) > 0$. So for
every infinite $\sigma \subseteq \bbN$, $s_{\sigma} > 0$.
As in Lemma \ref{l-gamma.5},
there is an infinite $\sigma \subseteq \bbN$
such that for every infinite $\eta \subseteq \sigma$,
$s_{\eta} = s_{\sigma}$.
Write $s = s_{\sigma}$.

Suppose by contradiction that $d(A[\sigma,s],y^*) = 0$.
We may assume that $\lim A[\sigma,s] = y^*$.
Let $r > 0$. Then there is $m$ such that
$A[\sigma^{\hbox{\tiny $\geq$} m},s] \subseteq B(y^*,\dgfrac{r}{2})$.
By the definition of $s$, there is $t \geq s$
such that $\beta(t - s) < \dgfrac{r}{2}$
and $\lim \tau(A[\sigma,t]) = w^*$.
Then $A[\sigma^{\hbox{\tiny $\geq$} m},t] \subseteq B(y^*,r)$.
Hence for every $r,\varepsilon > 0$ there are $m \in \bbN$
and $t \in [s,s + \varepsilon)$ such that
$A[\sigma^{\hbox{\tiny $\geq$} m},t] \subseteq B(y^*,r)$
and $\lim \tau(A[\sigma,t]) = w^*$.
It follows that there is a sequence $\vecu \subseteq X$ such that
$\lim \vecu = y^*$ and $\lim \tau(\vecu) = w^*$,
and hence $\tau^{\srfs{cl}}(y^*) = w^*$. A contradiction,
so $d(A[\sigma,s],y^*) > 0$.

From the fact that
$\setm{L_n - B(y^*,r)}{n \in \bbN}$ is completely discrete
for every $r > 0$,
it follows that $A[\sigma,s]$ is completely discrete. So we may assume
that for some $d > 0$, $A[\sigma,s]$ is $d$-spaced.
Let $\alpha$ and $\eta$ be as assured by the fact that $T$
is a track system.
It follows from the equicontinuity of $\sngltn{\gamma_n}_{n \in \bbN}$
that there is $\delta > 0$ such that for every $n \in \bbN$
and $t_1,t_2 \in [0,1]$: if $0 < t_2 - t_1 < \delta$, then
$$
{\thickmuskip=7mu \medmuskip=8mu \thinmuskip=8mu 
\rfs{diam}(\gamma_n([t_1,t_2])) +
\eta\left(\rule{0pt}{10pt}\rfs{diam}(\gamma_n([t_1,t_2]))\right) <
\dgfrac{d}{3}.
}
$$
Choose $t_1 \in [s,s + \dgfrac{\delta}{2}) \cap [0,1]$
such that $\lim \tau(A[\sigma,t_1]) = w^*$
and $t_2 \in (s - \dgfrac{\delta}{2},s) \cap [0,1]$.
For every $n \in \sigma$
let
$x_n = \gamma_n(t_1)$, $u_n = \gamma_n(t_2)$
and $J_n = [x_n,u_n]^{L_n}$,
that is, $J_n = \gamma_n([t_2,t_1])$.
Since $\abs{t_1 - t_2} < \delta$, it follows that
$$
\rfs{diam}\left(\rule{0pt}{10pt}B(J_n,\eta(\rfs{diam}(J_n)))\right)
\leq \rfs{diam}(J_n) + \eta(\rfs{diam}(J_n)) \leq
\dgfrac{d}{3}.
$$
We may assume that $\sigma = \bbN$.
Since $T$ is a track system, there is $h_n \in H(X)$ such that
$h_n(x_n) = u_n$,
$\rfs{supp}(h_n) \subseteq B(J_n,\eta(\rfs{diam}(J_n)))$
and $h_n$ is $\alpha$-bicontinuous.
We check\break
that $\setm{\rfs{supp}(h_n)}{n \in \bbN}$ is
$\frac{d}{3}$-spaced.
Let $m \neq n$. Then $\gamma_m(s),\gamma_n(s) \in A[\sigma,s]$
and so\break
$\norm{\gamma_m(s) - \gamma_n(s)} \geq d$.
Since $\gamma_m(s) \in J_m$ and the same holds for $n$,
it follows that
$$
d\left(\rule{0pt}{10pt}B(J_n,\eta(\rfs{diam}(J_n))),
B(J_n,\eta(\rfs{diam}(J_n))\right)
\geq
d - \dgfrac{2d}{3} = \dgfrac{d}{3}.
$$
So $\setm{\rfs{supp}(h_n)}{n \in \bbN}$ is $\frac{d}{3}$-spaced.

By Proposition \ref{metr-bldr-p5.16}(a),
$h \eqdf \bcirc_{n \in \bbN} h_{2n} \in \rfs{UC}(X)$.
It follows that $h^{\tau} \in \rfs{UC}(Y)$.
Let $w_n = x_n$ if $n$ is odd, and $w_n = u_n$ if $n$ is even.
Hence $h^{\tau}(\tau(\vecx)) = \tau(\vecw)$.
By the choice of $t_1$, $\tau(\vecx)$ converges to $w^*$.
By the choice of $\sigma$ and $t_2$,
$\tau(\setm{u_{2n}}{n \in \bbN})$ does not converge to $w^*$.
So $\tau(\vecw)$ is not a Cauchy sequence.
This contradicts the fact that $h^{\tau} \in \rfs{UC}(Y)$.
We have thus proved Claim 2.\smallskip

{\bf Claim 3.} $\rfs{bd}(Y) \subseteq \rfs{Rng}(\tau^{cl})$.
{\bf Proof }
Suppose by contradiction that
$z^* \in \rfs{bd}(Y) - \rfs{Rng}(\tau^{cl})$.
Let $\vecz \subseteq Y$ converge to $z^*$.
So $\vecx \eqdf \tau^{-1}(\vecz)$ is completely discrete.
By Claim~1, we may assume that $\limti{n} \delta(x_n) = 0$.
Let $\vecy$ be a subsequence of $\vecx$ which has a track system.
By Proposition \ref{metr-bldr-p7.5}(e),
$\vecy$ has a completely discrete track system
$\fortpl{\vecy}{y^*}{\setm{L_n}{n \in \bbN}}{\vecy\fprime}$.
By Claim~2, $\tau^{cl}(y^*) = z^*$. A contradiction,
so Claim~3 is proved.\smallskip

{\bf Claim 4.}
If $C \in \rfs{Cmp}(\rfs{bd}(X))$,
then $\tau^{\srfs{cl}}(C)$ is closed in $F$.
{\bf Proof } Let $C \in \rfs{Cmp}(\rfs{bd}(X))$
and $v \in \rfs{cl}(\tau^{\srfs{cl}}(C))$.
Let $\vecx\fprime \subseteq C$ be such that
$\limti{n} \tau^{\srfs{cl}}(x'_n) = v$.
If $\vecx\fprime$ has a Cauchy subsequence $\vecy$,
then $\lim \vecy \in C$ and $\tau^{\srfs{cl}}(\lim \vecy) = v$.
Suppose that $\vecx\fprime$ does not have Cauchy subsequences,
that is, $\vecx\fprime$ is completely discrete.
There is $\vecx \subseteq X$ such that
$\limti{n} d(x_n,x'_n) = 0$ and
$\limti{n} \tau(x_n) = v$.
So $\vecx$ is completely discrete.\break
Since $X$ is JN.TC,
there are a subsequence $\vecy$ of $\vecx$ and a track system
$T = \fortpl{\vecy}{z^*}{\vecL}{\vecz}$.
By Proposition \ref{metr-bldr-p7.5}(e), we may assume that $T$ is
a spaced track system, and by \ref{metr-bldr-p7.5}(f), $z^* \in C$.
By Claim 2, $\tau^{\srfs{cl}}(z^*) = v$,
so $\tau^{\srfs{cl}}(C)$ is closed.
\smallskip

{\bf Claim 5.}
$\tau^{cl}$ is $\onetoonen$.
{\bf Proof }
By Clause (3), for every component $C \in \rfs{Cmp}(\rfs{bd}(X))$,
either $\tau^{\srfs{cl}} \nrestriction C$ is $\onetoone$ or
$\tau^{\srfs{cl}}(C)$ is a singleton;
and for any distinct $C,D \in \rfs{Cmp}(\rfs{bd}(X))$,
either $\tau^{\srfs{cl}}(C) = \tau^{\srfs{cl}}(D)$
and $\tau^{\srfs{cl}}(C)$ is a singleton,
or $\tau^{\srfs{cl}}(C) \cap \tau^{\srfs{cl}}(D) = \emptyset$.
The argument is as in the proof of Theorem \ref{t4.36}.

Suppose by contradiction that $\tau^{\srfs{cl}}$ is not $\onetoonen$.
Then there is $C \in \rfs{Cmp}(\rfs{bd}(X))$ and $y \in \rfs{bd}(Y)$
such that $\tau^{\srfs{cl}}(C_0) = \sngltn{y}$.
Let $D$ be the component of $y$ in $\rfs{bd}(Y)$.
Then $\abs{D} > 1$. By Claims 3 and 4,
$\setm{\tau^{\srfs{cl}}(C)}{C \in \rfs{Cmp}(\rfs{bd}(X))
\mbox{ and } \tau^{\srfs{cl}}(C) \subseteq D}$
is a partition of $D$ into finitely many and more than 1 closed sets.
This contradicts the connectivity of~$D$.
\smallskip

{\bf Claim 6.}
Let $T = \fortpl{\vecx}{y^*}{\vecL}{\vecy}$
be a track system in $X$.
Then for every $d > 0$ there is $h \in \rfs{UC}(X)$ such that
$h^{cl}(y^*) \neq y^*$ and $\rfs{supp}(h) \subseteq B(y^*,d)$.
{\bf Proof }
Let $\alpha$ and $\eta$ be as assured by the fact that $T$
is a track system.
We may assume that $y^* \not\in \rfs{Rng}(\vecx)$,
and hence we may also assume that $d < d(\vecx,y^*)$.
Let $a > 0$ be such that $2a + \eta(a) < d$
and $b$ be such that $\alpha(b) < a - b$.
Clearly, $b < a$.
Let $n$ be such that $\norm{y_n - y^*} < b$.
Then $\norm{x_n - y_n} \geq d - b > a$,
and hence there is $z \in L_n$
such that $\norm{z - y_n} = \rfs{diam}([z,y_n]^{L_n}) = a$.
Since $L_n$ is an $\pair{\alpha}{\eta}$-track,
there is $h \in H(X)$ such that $h$ is $\alpha$-bicontinuous,
$h(y_n) = z$ and $\rfs{supp}(h) \subseteq B([z,y_n]^{L_n},\eta(a))$.
Clearly,
$B([z,y_n]^{L_n},\eta(a)) \subseteq B(y^*,b + a + \eta(a)) \subseteq
B(y^*,d)$.
So $\rfs{supp}(h) \subseteq B(y^*,d)$.
Suppose by way of contradiction that $h(y^*) = y^*$.
Then
$\norm{z - y^*} = \norm{h(y_n) - h(y^*)} \leq \alpha(\norm{y_n - y^*} <
\alpha(b)$. 
However,
$\norm{z - y^*} \geq \norm{z - y_n} - \norm{y_n - y^*} \geq a - b$.
That is, 
$\alpha(b) > a - b$, a contradiction.
So $h(y^*) \neq y^*$. 
So Claim 6 is proved.
\smallskip

{\bf Claim 7.} There is no sequence $\vecx \subseteq X$ such that
$\vecx$ is completely discrete, and $\tau(\vecx)$ is a Cauchy sequence.
{\bf Proof } Suppose otherwise, and let $\vecx$ be a counter-example
to the claim.
By Claim 1, we may assume that $\limti{n} \delta(x_n) = 0$.
Since $X$ is JN.TC,
there are a subsequence $\vecy$ of $\vecx$, $y^*$,
$\vecL$ and $\vecz$ such that
$T = \fortpl{\vecy}{y^*}{\vecL}{\vecz}$
is a track system.
By Proposition \ref{metr-bldr-p7.5}(e),
we may assume that $T$ is a spaced track system.
Let $w = \lim \tau(\vecx)$. So $w = \lim \tau(\vecy)$.
By Claim 2, (i) $\tau^{\srfs{cl}}(y^*) = w$.
Since $y^* \in \rfs{bd}(X)$ and $\vecy \subseteq X$,
it follows that $y^* \not\in \rfs{Rng}(\vecy)$,
and since $\vecy$ is completely discrete, $d(\vecy,y^*) > 0$.
By Claim 6, there is $h \in \rfs{UC}(X)$ such that
(ii) $h^{\srfs{cl}}(y^*) \neq y^*$
and $\rfs{supp}(h) \subseteq B(y^*,d(\vecy,y^*))$.
So $h \nrestriction \vecy = \rfs{Id}$.
By Proposition \ref{metr-bldr-p7.5}(g),
$T' \eqdf \fortpl{h(\vecy)}{h^{\srfs{cl}}(y^*)}{h(\vecL)}{h(\vecz)}$
is a track system.
Since $T$ is spaced and $h \in \rfs{UC}(X)$
it follows that $T'$ is also spaced.
Recall that $h(\vecy) = \vecy$ and so $\lim h(\vecy) = w$.
So by Claim 2 applied to $T'$,
(iii) $\tau^{\srfs{cl}}(h^{\srfs{cl}}(y^*)) = w$.
Facts (i)\,-\,(iii) contradict the fact that $\tau^{\srfs{cl}}$ is
$\onetoonen$. This proves Claim 7.
\smallskip

Suppose by contradiction that $\tau^{-1}$ is not uniformly
continuous. Then there\break
are sequences $\vecx,\vecy \subseteq X$
and $d > 0$ such that for every $n \in \bbN$,
$d(x_n,y_n) \geq d$
and $\limti{n} d(\tau(x_n),\tau(y_n)) = 0$.
We may assume that each of the sequences $\vecx$, $\vecy$,
$\tau(\vecx)$ and $\tau(\vecy)$ is either spaced or is a Cauchy
sequence.\smallskip

{\bf Claim 8.}
The sequences $\vecx$, $\vecy$, $\tau(\vecx)$ and $\tau(\vecy)$ are
spaced.
{\bf Proof }
Suppose by contradiction that $\vecx$ is a Cauchy sequence.
Since $\tau^{\srfs{cl}}$ is uniformly continuous and
$\rfs{Dom}(\tau^{cl}) = \rfs{cl}(X)$,
it follows that $\tau(\vecx)$ is a Cauchy sequence.
Hence $\tau(\vecy)$ is also a Cauchy sequence.
If $\vecy$ is a Cauchy sequence, then $\tau^{cl}$ is not $\onetoone$,
contradicting Claim 5; and if $\vecy$ is completely discrete,
then Claim 7 is contradicted. So $\vecx$ is not a Cauchy sequence.
The same is true for $\vecy$. By Claim 7,
$\tau(\vecx)$ and $\tau(\vecy)$ are completely discrete.
Claim 8 is proved.\medskip

We call a pair of sequences $\pair{\vecu}{\vecv}$ in $X$
a counter-example, if $\vecu$ and $\vecv$ are spaced,
$\inf(\setm{d(u_n,v_n)}{n \in \bbN}) > 0$ and
$\limti{n} d(\tau(u_n),\tau(v_n)) = 0$.

{\bf Claim 9.}
There is a counter-example $\pair{\vecu}{\vecv}$ such that
$\delta(\vecu) = 0$.
{\bf Proof }
By Claim~8, there is a counter-example $\pair{\vecx}{\vecy}$.
If $\delta(\vecx) = 0$ or $\delta(\vecy) = 0$,
then there is nothing to prove. Suppose otherwise.
By Proposition \ref{p4.40}, we may assume that
$d(\vecx,\vecy) > 0$.
Let $d > 0$ be such that $\vecx$ is $d$-spaced and
$d(\vecx,\vecy) \geq d$.
By possibly interchanging $\vecx$ and
$\vecy$, we may also assume that there are
$e_1 \geq e_2 > 0$ such that 
$\limti{n} \norm{x_n} = e_1$ and $\limti{n} \norm{y_n} = e_2$.
Let $x'_n = \frac{e_1}{\norm{x_n}}x_n$.
Since $\delta(\vecx) > 0$, there is $a > 0$ such that for every
$n \in \bbN$, $B(x'_n,a) \subseteq X$.
We may further assume that $a < \dgfrac{d}{8}$,
and that for every $n \in \bbN$, $d(x_n,x'_n) < \dgfrac{a}{2}$.
So $(\bigcup \setm{B(x'_n,a)}{n \in \bbN}) \cap
\setm{y_n}{n \in \bbN} = \emptyset$,
and for every distinct $m,n \in \bbN$,
$d(B(x'_m,a),B(x'_n,a)) > \dgfrac{d}{2}$.
Let $x''_n = (1 + \dgfrac{a}{2})x'_n$.
It follows that there is $h \in \rfs{LIP}(X)$ such that for every
$n \in \bbN$, $h(x_n) = x''_n$ and
$\rfs{supp}(h) \subseteq \bigcup \setm{B(x'_n,a)}{n \in \bbN}$.
Since $h(\vecx) = \vecx''$ and $h(\vecy) = \vecy$,
it follows that $\pair{\vecx''}{\vecy}$ is a counter-example.
So we may assume that $e_1 > e_2$,
and that $\norm{x_n} = e_1$ for every $n \in \bbN$.
We still assume that $\vecx$ is $d$-spaced and that
$d(\vecx,\vecy) \geq d$.

We now proceed as in the proof of Claim 1.
For $n \in \bbN^+$ let
$t_n = \min(\setm{t > 1}{t x_n \in \rfs{bd}(X)})$, $z_n = t_n x_n$,
$L_n = [x_n,z_n]$ and
$\gamma_n(t) = x_n + t(z_n - x_n)$, $t \in [0,1]$.
By Propo\-sition~\ref{metr-bldr-p7.5}(b),
for every distinct $m,n \in \bbN$,
$d(L_m,L_n) \geq \dgfrac{d}{2}$,
and clearly, $d(L_m,\vecy) \geq e_1 - e_2$.

Let $s = \min(\setm{t}{\delta(\setm{\gamma_n(t)}{n \in \bbN^+}) = 0})$.
We may assume that for every $n \in \bbN^+$,
$\delta(\gamma_n(s)) < \dgfrac{1}{n}$.
It follows that for every $t \in (0,s)$,
the family $\setm{\gamma_n([0,t])}{n \in \bbN^+}$ is spaced,
$d(\bigcup \setm{\gamma_n([0,t])}{n \in \bbN^+},\vecy) > 0$
and $\delta(\bigcup \setm{\gamma_n([0,t])}{n \in \bbN^+}) > 0$.
Also, since $X$ is bounded, $\setm{d(x_n,\gamma_n(t))}{n \in \bbN^+}$
is bounded.
So for every $t < s$ there is $h_t \in \rfs{UC}(X)$
such that for every $n \in \bbN^+$,
$h_t(x_n) = \gamma_n(t)$ and $h_t(y_n) = y_n$.
Since $h_t^{\tau} \in \rfs{UC}(Y)$,
$\limti{n} d(\tau(x_n),\tau(y_n)) = 0$
and $h_t^{\tau}(\tau(x_n)) = \tau(\gamma_n(t))$,
it follows that $\limti{n} d(\tau(\gamma_n(t)),\tau(y_n)) = 0$.

Let $\vecs \subseteq (0,s)$ be a sequence converging to $s$.
For every $n \in \bbN^+$ let $k_n \geq n$ be such that
$d(\tau(\gamma_{k_n}(s_n)),\tau(y_n)) < \dgfrac{1}{n}$.
Define $x'_n = \tau(\gamma_{k_n}(s_n))$.
It follows that $d(\vecx\fprime,\vecy) > 0$,
$\limti{n} d(\tau(x'_n),\tau(y_n)) = 0$,
$\limti{n} \delta(x'_n) = 0$ and $\vecx\fprime$ is spaced.
Claim 9 is thus proved.
\smallskip

Let $T = \fortpl{\vecy}{y^*}{\vecL}{\vecz}$
be a track system, and $\gamma_n$ be
a legal para\-metrization of $L_n$.
We say that $T$ is good, if for every $t \in [0,1)$,
$\inf(\setm{d(y_n,\gamma_n([0,t]))}{n \in \bbN}) > 0$.
\smallskip

{\bf Claim 10 }
If $\fortpl{\vecy}{y^*}{\vecL}{\vecz}$
is a track system, and $\gamma_n$ is a legal parametrization of $L_n$,
then there is $s \in (0,1]$ and an infinite $\sigma \subseteq \bbN$
such that
$\fortpl{\setm{\gamma_n(s)}{n \in \sigma}}
{y^*}{\setm{\gamma_n([0,s])}{n \in \sigma}}{\vecz}$
is a good track system,
and
$\lim_{n \in \sigma} d(\tau(y_n),\tau(\gamma_n(s))) = 0$.
{\bf Proof }
For every infinite $\eta \subseteq \bbN$ let
$s_{\eta} =
\inf(\setm{t \in [0,1]}{\lim_{n \in \eta} d(y_n,\gamma_n(t)) = 0})$.
As in previous analogous arguments, there is an infinite
$\eta \subseteq \bbN$ such that for every infinite
$\zeta \subseteq \eta$, $s_{\zeta} = s_{\eta}$.
Let $s = s_{\eta}$ and $\vect\,$ be a sequence converging to $s$
such that for every $i \in \bbN$,
$\lim_{n \in \eta} d(y_n,\gamma_n(t_i)) = 0$.
Let $\sigma = \setm{n_i}{i \in \bbN} \subseteq \eta$ be an
increasing sequence such that
$\limti{i} d(y_{n_i},\gamma_{n_i}(t_i)) = 0$.
By the equicontinuity of $\sngltn{\gamma_n}_{n \in \bbN}$,
$\limti{i} d(\gamma_{n_i}(t_i),\gamma_{n_i}(s)) = 0$.
So $\lim_{n \in \sigma} d(y_n,\gamma_n(s)) = 0$. Hence since
$\tau$ is uniformly continuous,
$\lim_{n \in \sigma} d(\tau(y_n),\tau(\gamma_n(s))) = 0$.
Now suppose by contradiction that there is \hbox{$t < s$} such that
$\liminf_{n \in \sigma} d(y_n,\gamma_n([0,t])) = 0$.
So there is an increasing sequence
$\zeta = \setm{k_i}{i \in \bbN} \subseteq \sigma$
and $\vect \subseteq [0,t]$ such that
$\limti{i} d(y_{k_i},\gamma_{k_i}(t_i)) = 0$.
We may assume that $\vect$ converges, say to $s^*$.
Hence $s^* \leq t < s$,
and $\lim_{n \in \zeta} d(y_n,\gamma_n(s^*)) = 0$.
So $s_{\zeta} \leq s^* < s$ a contradiction.
{\thickmuskip=2mu \medmuskip=1mu \thinmuskip=1mu 
So for every $t \in [0,s)$,
$\liminf_{n \in \sigma} d(y_n,\gamma_n([0,t])) > 0$.
\hbox{Since $\lim_{n \in \sigma} d(y_n,\gamma_n(s)) = 0$,}\break
}
it follows that for every $t \in [0,s)$,
$\liminf_{n \in \sigma} d(\gamma_n(s),\gamma_n([0,t])) > 0$;
and the fact that $L_n$ is a simple arc implies that
$\gamma_n(s) \not\in \gamma_n([0,s))$.
So $\inf(\setm{d(\gamma_n(s),\gamma_n([0,t]))}{n \in \sigma}) > 0$.
Claim~10 is proved.\smallskip

{\bf Claim 11 }
There are a counter-example $\pair{\vecu}{\vecv}$
and a completely discrete track system
$\fortpl{\vecu}{u^*}{\vecJ}{\vecu\fprime}$
such that $\inf_{n \in \bbN} d(J_n,v_n) > 0$.
{\bf Proof }
By Claim~9, there is a counter-example $\pair{\vecx}{\vecy}$ such that
$\delta(\vecx) = 0$.
Let $T = \fortpl{\vecx}{x^*}{\vecL}{\vecx\fprime}$
be a completely discrete track system for $\vecx$.
By Claim~10, we may assume that $T$ is a good track system.

Suppose first that $d \eqdf \liminfti{n} d(L_n,y_n) > 0$.
Let $\setm{\ell_i}{i \in \bbN}$ be a subsequence of $\bbN$ such that
{\thickmuskip=2mu \medmuskip=1mu \thinmuskip=1mu 
$d(L_{\ell_i},y_{\ell_i}) \geq \dgfrac{d}{2}$.
Hence $\vecu = \setm{x_{\ell_i}}{i \in \bbN}$,
$u^* = x^*$,
$\vecv = \setm{y_{\ell_i}}{i \in \bbN}$ and
$\vecJ = \setm{L_{\ell_i}}{i \in \bbN}$,
are as required in the claim.
}

{\thickmuskip=1mu \medmuskip=1mu \thinmuskip=0mu 
Assume next that $\liminfti{n} d(L_n,y_n) = 0$.
\hbox{So we may assume that $\limti{n} d(L_n,y_n) = 0$.}
}
Let $\gamma_n$ be a legal paramet\-rization of $L_n$.
Hence there is $\vect \subseteq [0,1]$
such that\break
$\limti{n} d(\gamma_n(t_n),y_n) = 0$.
We may assume that $\vect$ is convergent.
Let $t = \lim \vect$.
It easily follows that $\limti{n} d(\gamma_n(t),y_n) = 0$.
Clearly $t < 1$, for otherwise $\limti{n} d(x_n,y_n) = 0$.
For every $n \in \bbN$ let $u_n = \gamma_n(t)$,
$v_n = x_n$ and $J_n = \gamma_n([0,t])$.

Since $\tau$ is uniformly continuous,
we have that $\limti{n} d(\tau(u_n),\tau(y_n)) = 0$.
Also, $\limti{n} d(\tau(v_n),\tau(y_n)) = 0$.
Hence
$\limti{n} d(\tau(u_n),\tau(v_n)) = 0$.
Since $\fortpl{\vecx}{x^*}{\vecL}{\vecx\fprime}$
is a good track system,
$\inf_{n \in \bbN} d(x_n,\gamma_n([0,t])) > 0$.
That is,
$\inf_{n \in \bbN} d(v_n,J_n) > 0$.
By Proposition \ref{metr-bldr-p7.5}(c)(i) applied to $T$ and $t$,
$\fortpl{\vecu}{x^*}{\vecJ}{\vecx\fprime}$ is a track system.
So $\vecu$, $\vecv$, $x^*$ and $\vecJ$ are as required.
Claim~11 is proved.\smallskip

{\bf Claim 12 }
There are a counter-example $\pair{\vecu}{\vecv}$
and a completely discrete track system
$\fortpl{\vecu}{u^*}{\vecJ}{\vecu\fprime}$
such that $d(\bigcup \setm{J_n}{n \in\bbN},\vecv) > 0$.
{\bf Proof }
Let $\pair{\vecu}{\vecv}$
and $\fortpl{\vecu}{u^*}{\vecJ}{\vecu\fprime}$
be as assured by the previous claim.
We show that there is an infinite
$\sigma \subseteq \bbN$ such that
$\pair{\vecu \nrestriction \sigma}{\vecv \nrestriction \sigma}$
and $\fortpl{\vecu \nrestriction \sigma}{u^*}
{\vecJ \nrestriction \sigma}{\vecu\fprime \nrestriction \sigma}$
are as required in the claim.
We shall apply Proposition~\ref{p4.40} with $F_n$ taken to be
$\sngltn{v_n}$ and $K_n$ taken to be $J_n$.
By our assumptions, clauses (i) and (iii) of \ref{p4.40} do hold.
We show that (ii) holds. Let $\gamma_n$ be a legal parametrization
of $J_n$. Suppose that $\varepsilon > 0$.
Then by the equicontinuity of $\sngltn{\gamma_n}_{n \in \bbN}$,
there is $\delta > 0$ such that for every $n \in \bbN$
and $t_1,t_2 \in [0,1]$: if $\abs{t_1 - t_2} < \delta$,
then $\norm{\gamma_n(t_1) - \gamma_n(t_2)} < \varepsilon$.
Define $\ell_{\varepsilon} = [\dgfrac{1}{\delta}] + 1$.
Then $\ell_{\varepsilon}$ fulfills the requirement of Clause~(ii) of
\ref{p4.40}.
The set $\sigma$ obtained from \ref{p4.40} is as required.
This proves Claim~12.\smallskip

{\bf Conclusion of the proof of the theorem:} 
Let $\pair{\vecx}{\vecy}$ and
$T = \fortpl{\vecx}{x^*}{\vecL}{\vecx\fprime}$
be as assured by Claim 12.
By Claim 7, $\tau(\vecy)$ is completely discrete.
So we may assume that $\tau(\vecy)$ is spaced.
Write $d_1 = d(\bigcup \setm{L_n}{n \in \bbN},\vecy)$.
Let $\gamma_n$ be a legal parametrization of $L_n$.
For every infinite $\sigma \subseteq \bbN$ let
$s_{\sigma} = \inf(\setm{t \in [0,1]}
{\lim_{n \in \sigma} d(\tau(\gamma_n(t)),\tau(y_n)) = 0})$.
Let $\sigma$ be such that for every infinite $\eta \subseteq \sigma$,
$s_{\eta} = s_{\sigma}$. Since $\tau(\vecx\fprime)$ is convergent
and $\tau(\vecy)$ is spaced, $s \eqdf s_{\sigma} > 0$.
As in previous analogous arguments,
$\setm{\gamma_n(s)}{n \in \sigma}$
is completely discrete.
So we may assume that that for some $d_2 > 0$,
$\setm{\gamma_n(s)}{n \in \sigma}$ is $d_2$-spaced.
Set $d = \min(d_1,d_2)$.
Let $\alpha,\eta$ be as assured by the fact that $T$ is a track system.
Let $a > 0$ be such that $a + \eta(a) < \dgfrac{d}{3}$.
By the equicontinuity of $\sngltn{\gamma_n}_{n \in \bbN}$,
there is $\delta > 0$ such that for every $n \in \bbN$
and $t_1,t_2 \in [0,1]$: if $\abs{t_1 - t_2} < \delta$,
then $\norm{\gamma_n(t_1) - \gamma_n(t_2)} < a$.
By the choice of $s$, there is $t_1 \in [s,s + \dgfrac{\delta}{2})$
such that $\lim_{n \in \sigma} d(\tau(\gamma_n(t_1)),\tau(y_n)) = 0$.
Also, choose $t_2 \in (s - \dgfrac{\delta}{2},s)$.
Then by the choice of $\sigma$ and $s$,
$\inf_{n \in \sigma} d(\tau(\gamma_n(t_2)),\tau(y_n)) > 0$.
For $n \in \sigma$ write $u_n = \gamma_n(t_1)$, $v_n = \gamma_n(t_2)$
and $J_n = \gamma_n([t_2,t_1])$.
Let $n \in \sigma$. Then since $L_n$
is an $\pair{\alpha}{\eta}$-track,
there is $h_n \in H(X)$ such that $h_n$ is $\alpha$-bicontinuous,
$h_n(u_n) = v_n$
and $\rfs{supp}(h_n) \subseteq B(J_n,\eta(\rfs{diam}(J_n)))$.
Since $\abs{t_1 - t_2} < \delta$, it follows that
$\rfs{diam}(J_n) < a$.
\hbox{So for every $x \in \rfs{supp}(h_n)$,
$\norm{x - \gamma_n(s)} < a + \eta(a) < \dgfrac{d_2}{3}$.}
This implies that
$d(\rfs{supp}(h_m),\rfs{supp}(h_n)) > \dgfrac{d_2}{3}$
for any $m \neq n$.
We conclude that $h \eqdf \bcirc_{n \in \sigma_1} h_n$ is well-defined
and belongs to $\rfs{UC}(X)$.
Clearly, $\rfs{supp}(h) \subseteq B(\bigcup_{n \in \bbN} L_n,\eta(a))$.
Since $d(\bigcup_{n \in \bbN} L_n,\vecy) = d_1$ and $\eta(a) < d_1$,
we have that $\rfs{supp}(h) \cap \rfs{Rng}(\vecy) = \emptyset$
and hence $h \nrestriction \vecy = \rfs{Id}$.
It follows that
$\inf_{n \in \sigma} d(h^{\tau}(\tau(y_n)),h^{\tau}(\tau(u_n)))
= \inf_{n \in \sigma} d(\tau(y_n),\tau(v_n))) > 0$.
But
$\lim_{n \in \sigma} d(\tau(y_n),\tau(u_n)) = 0$.
So $h^{\tau} \not\in \rfs{UC}(Y)$. A contradiction.
\hfill\myqed

\begin{remark}\label{metr-bldr-r7.8}
\begin{rm}
(a)
Clause (2) in Theorem \ref{metr-bldr-t7.7} can be relaxed.
In that case Clause (5) has to be strengthened.
Replace Clauses (2) and (5) by Clauses (2.1) and (5.1) stated below.
\begin{itemize}
\addtolength{\parskip}{-11pt}
\addtolength{\itemsep}{06pt}
\item[(2.1)] $\rfs{bd}(X)$ has countably many components.
\item[(5.1)] If $C$ is a component of $\rfs{bd}(Y)$,
then $C$ is not a singleton,
and either $C$ is arcwise connected or $C$ is locally connected.
\vspace{-05.7pt}
\end{itemize}
The proof of \ref{metr-bldr-t7.7} is changed only in one place.
In the proof of Claim 5, the component $D$ of $\rfs{bd}(Y)$ is
partitioned into countably many closed sets.
By (5.1), this is impossible. So a contradiction is reached.

There are spaces $X$ which satisfy (1), (2.1), (3) and (4),
but do not satisfy (2). However, such examples are rare.

(b) Let
$K^{\calO}_{\srfs{BLPM}} = \setm{Y}{Y \mbox{ is an open subset of a
Banach Lipschitz } \mbox{manifold}}$.
(See Definition \ref{metr-bldr-d6.27}).
In Theorem \ref{metr-bldr-t7.7} replace the assumption that
$X \in K^{\calO}_{\srfs{BNC}}$ By the assumption that
$X \in K^{\calO}_{\srfs{BLPM}}$.
Then Parts (a) and (b) of \ref{metr-bldr-t7.7} remain true,
and the proof remains as is.

   \index{N@koblpm@@
          $K^{\calO}_{\srfs{BLPM}} =
	  \setm{Y}{Y
	  \mbox{ is an open subset of a Banach Lipschitz manifold}}$}

(c) The sphere of a Banach space satisfies the assumptions of Part (b).
See Remark \ref{metr-bldr-r7.4}.
\end{rm}
\end{remark}

\begin{question}\label{metr-bldr-q7.9}
\begin{rm}
(a) Prove Theorem~\ref{metr-bldr-t7.7} for incomplete normed spaces.

(b) Let $E$ be a Banach space.
Let $\setm{\overB_n}{n \in \bbN}$ be a spaced set of closed balls
such that  for every $n$, $\overB_n \subseteq B^E(0,1)$.
Let $X = B^E(0,2) - \bigcup_{n \in \sboldbbN} \overB_n$.
Let $Y \in K^{\calO}_{\srfs{BNC}}$ and $\tau \in H(X,Y)$.
Suppose that $(\rfs{UC}(X))^{\tau} \subseteq \rfs{UC}(Y)$.
Is $\tau\inverse$ uniformly continuous?

Note that $X$ is not JN.TC, but it satisfies all the other
assumptions of Theorem~\ref{metr-bldr-t7.7}.
\end{rm}
\end{question}

\begin{prop}\label{metr-bldr-p7.10}
Suppose that $X$ is an open ball of a Banach space.
Then $X$ satisfies
Clauses \num{1}-\num{4} of Theorem $\ref{metr-bldr-t7.7}(b)$.
\end{prop}

\noindent
{\bf Proof } The proof is easy and is left to the reader.
\hfill\myqed

\subsection{The non-existence of isomorphisms between groups of\\
different types.}
\label{ss7.2}
\label{ss7.2-non-existence-of-isomorphisms}

In the previous chapters we considered groups of various types.
We now show that groups of different types cannot be isomorphic
unless they coincide. We shall deal with the groups
$\rfs{UC}(X)$, $\rfs{LUC}(X)$, $\rfs{BUC}(X)$, $\rfs{BPD.UC}(X)$ and
$\rfs{EXT}(X)$, and we add to this list the group $H(X)$.
Let $\calP,\calQ$ denote one of the above properties and
$\calP(X),\calQ(X)$ be the groups they define.
We describe the situation precisely.
It may happen that for distinct properties
$\calP$ and $\calQ$, there is $\varphi$ such that
$\iso{\varphi}{\calP(X)}{\calQ(Y)}$.
But in that case either $\calP(X) = \calQ(X)$ and $\varphi$ is induced
by a homeomorphism belonging to $\calQ^{\pm}(X,Y)$,
or
$\calP(Y) = \calQ(Y)$, and $\varphi$ is induced
by a homeomorphism belonging to $\calP^{\pm}(X,Y)$.
The situation with regard to such questions
is not sorted out completely,
and we only state results which follow directly from the theorems
that have been proved so far.
Only some of the possible consequences are stated and proved.

Let $X \in K^{\calO}_{\srfs{NRM}}$
and $h \in H(X)$. Recall that $h$ is said to be internally extendible,
if there is $\barh \in H(\overfs{int}(X))$
such that $\barh \supseteq h$.
Denote $\barh$ by $h^{\overfss{int}}$.
If
$\calP = \rfs{UC},\rfs{BUC},\rfs{BPD.UC}$,
then $\calP(X) \subseteq \rfs{IXT}(X)$.
See Definition \ref{metr-bldr-d2.23}(b).
For these $\calP$'s define $X^{\calP} = \overfs{int}(X)$ and
$\calP^{\srfs{BNO}}(X) = \setm{h^{\overfss{int}}}{h \in \calP(X)}$.
So $\pair{X^{\calP}}{\calP^{\srfs{BNO}}(X)} \in K_{\srfs{BO}}$.
See Definition \ref{d2.7}(b).
For $\calP = \rfs{LUC},\rfs{EXT}$, write $X^{\calP} = X$
and $\calP^{\srfs{BNO}}(X) = \calP(X)$.
So
$\pair{X^{\calP}}{\calP^{\srfs{BNO}}(X)} \in K_{\srfs{NO}}$.

\begin{cor}\label{c4.42}\label{metr-bldr-c7.11}
Let $X,Y \in K^{\calO}_{\srfs{NRM}}$.

\num{a} If $\iso{\varphi}{\rfs{LUC}(X)}{\calP(Y)}$,
then
$\calP(Y) = \rfs{LUC}(Y)$,
and there is
$\tau \in \rfs{LUC}^{\pm}(X,Y)$
such that $\tau$ induces $\varphi$.

\num{b} Let $X,Y \in K^{\calO}_{\srfs{NFCB}}$.
Assume that $X$ is BUD.AC and MV1, $Y$ is UD.AC
and that $\iso{\varphi}{\rfs{UC}(X)}{\rfs{BUC}(Y)}$.
Then $\rfs{BUC}(X) = \rfs{UC}(X)$,
and there is $\tau \in \rfs{BUC}^{\pm}(X,Y)$
such that $\tau$ induces $\varphi$.
($X$ may be unbounded, and $X$ need not be UC-equivalent to $Y$).

\num{c} Let $X,Y \in K^{\calO}_{\srfs{NFCB}}$.
Suppose that $X$ is BPD.AC,
$Y$ is UD.AC, and $Y$ has
the discrete path property for large distances.
Let
$\iso{\varphi}{\rfs{UC}(X)}{\rfs{BPD.UC}(Y)}$.
Then $\rfs{BPD.UC}(X) = \rfs{UC}(X)$,
and there is $\tau \in \rfs{BPD.UC}^{\pm}(X,Y)$
such that $\tau$ induces $\varphi$.

\num{d} Let $X,Y \in K^{\calO}_{\srfs{BNC}}$.
Suppose that $X$ is BPD.AC and BR.LC.AC.
Let
$\isobreak{\varphi}{\rfs{BUC}(X)}{\rfs{BPD.UC}(Y)}$.
Then $\rfs{BUC}(X) = \rfs{BPD.UC}(X)$,
and there is
$\tau \in \rfs{BPD.UC}^{\pm}(X,Y)$
such that $\tau$ induces $\varphi$.

\num{e} Suppose that $X,Y \in K^{\calO}_{\srfs{BNC}}$,
and $X$ or $Y$ are infinite-dimensional.
Then there is no $\iso{\varphi}{\rfs{UC}(X)}{\rfs{EXT}(Y)}$.
(Since $\rfs{EXT}(X) = \rfs{BUC}(X)$ whenever $X$
is finite-dimensional, such cases are included in part \num{c}).

\num{f} Suppose that $X,Y \in K^{\calO}_{\srfs{BNC}}$,
and $X$ or $Y$ are infinite-dimensional.
Then there is no $\iso{\varphi}{\rfs{UC}(X)}{H(Y)}$.
\end{cor}

\noindent
{\bf Proof }
(a) Since $\calP^{\srfs{BNO}}(Y) \cong \calP(Y)$,
there is $\iso{\barvarphi}{\rfs{LUC}(X)}{\calP^{\srfs{BNO}}(Y)}$.
$\pair{Y^{\calP}}{\calP^{\srfs{BNO}}(Y)} \in K_{\srfs{BNO}}$.
\kern-1ptAlso $\pair{X}{\rfs{LUC}(X)} \in K_{\srfs{BNO}}$.
So by Theorem~\ref{metr-bldr-t2.8}(b),
there is $\tau \in H(X,Y^{\calP})$ such that $\tau$ induces
$\barvarphi$.
Since $\pair{X}{\rfs{LUC}(X)}$ is transitive,
$\pair{Y^{\calP}}{\calP^{\srfs{BNO}}(Y)}$ is transitive.
Since $Y$ is an orbit of $\pair{Y^{\calP}}{\calP^{\srfs{BNO}}(Y)}$,
$Y^{\calP} = Y$.
Hence $\barvarphi = \varphi$, and hence $\tau$ induces $\varphi$.

Note that if $\calP = \rfs{UC},\rfs{LUC},\rfs{BUC},\rfs{BPD.UC}$,
then $\rfs{UC}_{00}(Y) \subseteq \calP(Y)$.
So $(\rfs{UC}_{00}(Y))^{\tau\inverse} \subseteq
(\calP(Y))^{\tau\inverse} \subseteq \rfs{LUC}(X)$.
Also,
$\rfs{UC}_{00}(Y) = \rfs{UC}(Y,\calU)$, where $\calU$ is the
set of all open BPD subsets of $Y$.
So by Theorem~\ref{metr-bldr-t4.8}(b),
{\thickmuskip=2mu \medmuskip=1mu \thinmuskip=1mu 
$\tau\inverse \in \rfs{LUC}^{\pm}(Y,X)$,
that is, $\tau \in \rfs{LUC}^{\pm}(X,Y)$.
So $\calP(Y) = (\rfs{LUC}(X))^{\tau} = \rfs{LUC}(Y)$.}
\smallskip

(b) By Corollary \ref{metr-bldr-c2.26}
there is $\tau \in H(X,Y)$ such that $\tau$
induces $\varphi$. So $(\dagger)$ $(\rfs{UC}(X))^{\tau} = \rfs{BUC}(Y)$.
{\thickmuskip=2mu \medmuskip=1mu \thinmuskip=1mu 
We show that $\tau \in \rfs{BUC}(X,Y)$.
By $(\dagger)$, $(\rfs{UC}(X))^{\tau} \subseteq \rfs{BUC}(Y)$
and $(\rfs{BUC}(Y))^{\tau\inverse} \subseteq \rfs{BUC}(X)$.
}
Recall that $X$ is BUD.AC and MV1.
So by Corollary \ref{metr-bldr-c5.18}, $\tau \in \rfs{BUC}(X,Y)$.

We show that $\tau\inverse \in \rfs{UC}(Y,X)$.
By $(\dagger)$, $\rfs{UC}_0(Y))^{\tau\inverse} \subseteq \rfs{UC}(X)$.
Recall that $Y$ is UD.AC. So by Theorem~\ref{metr-bldr-t5.6},
$\tau\inverse \in \rfs{UC}(Y,X)$,
and hence $\tau \in \rfs{BUC}^{\pm}(X,Y)$.
Then $\rfs{UC}(X) = (\rfs{BUC}(Y))^{\tau\inverse} = \rfs{BUC}(X)$.
\smallskip

(c) Let $\iso{\varphi}{\rfs{UC}(X)}{\rfs{BPD.UC}(Y)}$.
By Corollary \ref{metr-bldr-c2.26},
there is $\tau \in H(X,Y)$ such that $\tau$ induces $\varphi$.
So $(*)$ $(\rfs{UC}(X))^{\tau} = \rfs{BPD.UC}(Y)$.
By $(*)$, $(\rfs{UC}_{00}(X))^{\tau} = \rfs{BPD.UC}(Y)$.
Recall that $X$ is BPD.AC.
Hence by Theorem~\ref{t4.23}, $\tau \in \rfs{BPD.UC}(X,Y)$.

Obviously, $\rfs{UC}_0(Y) \subseteq \rfs{BPD.UC}(Y)$.
So by $(*)$,
$(\rfs{UC}_0(Y))^{\tau\inverse} \subseteq \rfs{UC}(X)$.
Recall that $Y$ is UD.AC.
Hence by Theorem \ref{t4.10}, $\tau\inverse \in \rfs{UC}(Y,X)$.
Since $Y$ has the discrete path property for large distances,
by Proposition~\ref{metr-bldr-p4.3}(b),
$\tau\inverse$ is uniformly continuous for all distances.
That is, for some $\alpha \in \rfs{MC}$, $\tau\inverse$ is
$\alpha$-continuous. In particular, $\tau\inverse$ is boundedness
preserving.
So $\tau\inverse \in \rfs{BPD.UC}(Y,X)$.
In summary, $\tau\inverse \in \rfs{BPD.UC}^{\pm}(Y,X)$.
It follows that
$\rfs{UC}(X) = (\rfs{BPD.UC}(Y))^{\tau\inverse} = \rfs{BPD.UC}(X)$.

(d) By Theorem~\ref{metr-bldr-t2.8}(a),
there is $\tau \in H(X,Y)$ such that $\tau$ induces $\varphi$.
This means that $(\rfs{BUC}(X))^{\tau} = \rfs{BPD.UC}(Y)$.
By Theorem~\ref{t4.23}, $\tau \in \rfs{BPD.UC}(X,Y)$,
and by Theorem~\ref{t4.27}(a),
$\tau\inverse \in \rfs{BPD.UC}(Y,X)$.
Hence $\tau\inverse \in \rfs{BPD.UC}^{\pm}(Y,X)$.
It follows that
$\rfs{BUC}(X) = (\rfs{BPD.UC}(Y))^{\tau\inverse} = \rfs{BPD.UC}(X)$.

(e) Suppose by contradiction that
$\iso{\varphi}{\rfs{UC}(X)}{\rfs{EXT}(Y)}$.
By Theorem \ref{t2.4}(a),
there is $\tau \in H(X,Y)$ such that $\tau$ induces $\varphi$.
So
$(\rfs{UC}(X))^{\tau} = \rfs{EXT}(Y)$.

Suppose that $Y$ is an open subset of the Banach space $F$.
Let $B$ be a ball in $F$ such that $\rfs{cl}^F(B) \subseteq Y$.
Clearly, for every $h \in \rfs{UC}_{\rme}(B)$ there
is $\tilde{h} \in \rfs{EXT}(Y)$ such that $\tilde{h}$ extends $h$.
Let $\eta = \tau\inverse \nrestriction B$ and $C = \eta(B)$.
Since $(\rfs{EXT}(Y))^{\tau\inverse} \subseteq \rfs{UC}(X)$,
$(\rfs{UC}_{\rme}(B))^{\eta} \subseteq \rfs{UC}(C)$. So also
$(\rfs{UC}_0(B))^{\eta} \subseteq \rfs{UC}(C)$.
So by Theorem \ref{t4.10}, $\eta$ is UC.
It follows that $C$ is bounded, and hence $\rfs{bd}(C)$
is not a singleton.
Clearly, $\rfs{bd}(C) = \eta^{cl}(\rfs{bd}(B))$,
and so $\rfs{bd}(C)$ is connected.
So no component of $\rfs{bd}(C)$ is a singleton.
By Proposition~\ref{metr-bldr-p7.10},
$B$ satisfies Clauses (1)\,-\,(4) of Theorem \ref{t4.41}(b).
By Theorem~\ref{t4.41}(b) applied to $B$, $C$ and $\eta$,
$\eta\inverse$ is UC.
In summary, $\eta \in \rfs{UC}^{\pm}(B,C)$.

Choose $h \in H(B) - \rfs{UC}(B)$ which is strongly extendible.
So there is
$\tilde{h} \in \rfs{EXT}(Y)$ extending $h$.
So $\tilde{h}^{\tau\inverse} \in \rfs{UC}(X)$.
Hence $h^{\eta} = \tilde{h}^{\tau\inverse} \nrestriction C \in
\rfs{UC}(C)$.
Since $\eta\inverse \in \rfs{UC}^{\pm}(C,B)$,
$h = (h^{\eta})^{\eta\inverse} \in \rfs{UC}(B)$.
A contradiction.

(f) The proof is identical to the proof of part (e).
\medskip\hfill\myqed

The following trivial examples show that the conclusions of
Corollary \ref{metr-bldr-c7.11}(b), (c) and (f) cannot be strengthened.

\begin{example}\label{bddly-lip-bldr-e7.12}
\begin{rm}
(a) There are regular open sets $X,Y \subseteq \bbR^2$ such that
\newline
(1) $\rfs{UC}(X) = \rfs{BUC}(X) \cong \rfs{BUC}(Y) \not\cong
\rfs{UC}(Y)$.
\newline
(2) $X$ is BUD.AC and MV1, and $Y$ is UD.AC.

(b) Let $X = (0,1)$. Then $\rfs{UC}(X) = \rfs{BPD.UC}(X)$.

(c) Let $E$ be a Banach space. Let $Y = B^E(0,1)$.
Let $\fnn{\tau}{E}{Y}$ be defined by
$\tau(x) = \frac{x}{1 + \norm{x}}$.
Then $\tau \in \rfs{BPD.UC}^{\pm}(E,Y)$,
$\rfs{BUC}(E) = \rfs{BPD.UC}(E)$ and 
$\rfs{BPD.UC}(Y) \not\cong \rfs{BUC}(Y)$.
\end{rm}
\end{example}

\noindent
{\bf Proof }
(a) For $n \in \bbN$ we define an open set $B_n$.
\newline\centerline{
$B_n = B(0,1) - \bigcup_{i < n} \overB((\dgfrac{i}{n},0),\frac{1}{3n})$.
}
So $B_n$ is obtained by removing from $B(0,1)$
$n$ pairwise disjoint closed balls each of which contained
in $B(0,1)$.
For every $n \in \bbN$ let
$X_n = (n,0) + \frac{1}{n + 4} \mcdot B_n$
and
$Y_n = (n,0) + \frac{1}{4} \mcdot B_n$.
Let $X = \bigcup_{n \in \bbN} X_n$ and
$Y = \bigcup_{n \in \bbN} Y_n$.
Note that for every $n \neq m$,
$d(X_n,X_m), d(Y_n,Y_m) \geq \dghalf$
and $X_n \cong Y_n \not\cong Y_m$.
Note that $\limti{n} \rfs{diam}(X_n) = 0$
and for every $n$, $\rfs{diam}(Y_n) = \dghalf$.
It is easy to check that $X$ and $Y$ have the required properties.

The proofs of Parts (b) and (c) are trivial.
\hfill\myqed

\begin{question}\label{metr-bldr-q5.4}
\begin{rm}
For $n > 1$, construct an open subset $X \subseteq \bbR^n$
such that $\rfs{UC}(X) = \rfs{BPD.UC}(X)$.
Note that if $X$ is such an example,
then every connected component of $X$ is an example.
Note that every example which is a connected set is bounded.
\vspace{-2.0mm}
\end{rm}
\end{question}

\newpage




\section{The group of locally $\Gamma$-continuous
homeo\-morphisms of the closure of an open set}
\label{s8}

\subsection{General description.}
\label{ss8.1}
Lipschitz equivalence between open subsets of $\bbR^n$ is relevant
in the theory of function spaces.
Suppose that $U,V$ be open subsets of $\bbR^n$.
The fact that $U,V$ are homeomorphic by a bilipschitz
homeomorphism or by a quasiconformal homeomorphism is equivalent to the
fact that certain Sobolev spaces of functions from $U$ to $\bbR$ and
from $V$ to $\bbR$ are isomorphic as lattice ordered vector spaces.
These results appear in \cite{GV1}, \cite{GV2} and \cite{GRo}.
We consider the analogous question for the
setting in which the Sobolev function spaces are replaced by
homeomorphism groups.

The simplest question of this kind is as follows.
Let $X \subseteq \bbR^n$ and $Y \subseteq \bbR^m$ be open sets.
Suppose that
$\iso{\varphi}{\rfs{LIP}(\rfs{cl}(X))}{\rfs{LIP}(\rfs{cl}(Y))}$.
Prove that there is $\tau \in \rfs{LIP}^{\pm}(X,Y)$ such that
$\tau$ induces $\varphi$.

We shall prove the above statement for bounded open subsets of $\bbR^n$
which have a well-behaved boundary.
In fact, we shall deal with a different group of homeomorphisms,
namely, the group $\rfs{LIP}^{\srfs{LC}}(\rfs{cl}(X))$
of locally bilipschitz homeomorphisms of $\rfs{cl}(X)$.
But for bounded subsets of $\bbR^n$ this group coincides
$\rfs{LIP}(\rfs{cl}(X))$.

The group of bilipschitz homeomorphisms is only a special case.
It is generalized to the setting
of $\itGamma$-bicontinuous homeomorphisms, where $\itGamma$ is any
principal modulus of continuity.
(See Property M6 in Definition \ref{nn1.3}).

The open sets for which we know to prove such results at this point,
have a very well-behaved boundary.
They are called locally $\itGamma$-LIN-bordered sets.
See Definition \ref{editor-d8.1}(c).
Essentially these are the open subsets of a normed space
whose closure is a manifold with a boundary.
For such sets we define
the group of completely locally $\itGamma$-bicontinuous homeomorphisms.
This group is denoted by $H^{\srfs{CMP.LC}}_{\itGamma}(X)$,
and is defined in Definition \ref{d-bddlip-1.6}.
We give here an equivalent definition.
Let $X$ be an open subset of a metric space $E$ and $\itGamma$ be
a modulus of continuity.
Define
$$
H^{\srfs{CMP.LC}}_{\itGamma}(X) =
\setm{g \in H(\rfs{cl}(X))}
{g \mbox{ is locally } \itGamma\mbox{-bicontinuous}
\mbox{ and } g(X) = X}.
$$
Suppose that $\itGamma,\itDelta$ are moduli of continuity and
$\itGamma$ is principal, $E,F$ are normed spaces, $X \subseteq E$,
$Y \subseteq F$ and $X,Y$ are locally $\itGamma$-LIN-bordered sets.
We shall prove that if
$\iso{\varphi}{H^{\srfs{CMP.LC}}_{\itGamma}(X)}
{H^{\srfs{CMP.LC}}_{\itDelta}(Y)}$,
then $\itGamma = \itDelta$
and there is $\iso{\tau}{\rfs{cl}(X)}{\rfs{cl}(Y)}$
such that $\tau(X) = Y$, $\tau$ is locally $\itGamma$-bicontinuous and
$\varphi(g) =
\tau \scirc g \scirc \tau\inverse$
for every $g \in H^{\srfs{CMP.LC}}_{\itGamma}(X)$.

The above statement is also true when $X$ and $Y$ are
open subsets of a normed Lipschitz manifold,
see Theorem~\ref{ams-bddly-lip-bldr-t1.4}(b).
The argument for manifolds is essentially identical,
so proofs will be given only for the class of open subsets of normed
spaces.

\subsection{Statement of the main theorems and the plan of the proof.}
\label{ss8.2}

We shall now define the class of open sets with a well-behaved
boundary.

\begin{defn}\label{editor-d8.1}
\begin{rm}
(a) Let $E$ be a normed space, $A \subseteq E$ and $r > 0$.
The set
$\rfs{BCD}^E(A,r) \eqdf
\break B^E(0,r) - A$
is called the
{\it boundary chart domain based on $E$ and $A$ with radius $r$}.
   \index{N@bcd00@@$\rfs{BCD}^E(A,r) = B^E(0,r) - A$.
          A boundary chart domain based on $E$ and $A$ with radius $r$}
   \index{boundary chart domain based on $E$ and $A$
          with radius $r$. A set of the form $\rfs{BCD}^E(A,r)$}
We say that $A \subseteq E$ is a {\it closed half space} of $E$,
if there is $\varphi \in E^*$
such that $A = \setm{x \in E}{\varphi(x) \geq 0}$.
   \index{closed half space. A set of the form
          $\setm{x \in E}{\varphi(x) \geq 0}$,
	  where $\varphi \in E^*$}
Suppose that $\rfs{dim}(E) > 1$,
and $A$ is either a closed subspace of $E$ 
different from $\sngltn{0}$ or a closed half space of $E$.
Then $\rfs{BCD}^E(A,r)$
is called a {\it a linear boundary chart domain}.
   \index{linear boundary chart domain. A set of the form
          $\rfs{BCD}^E(A,r)$, where $A$ is a
	  closed\\\indent subspace of $E$ different from $\sngltn{0}$
	  or a closed half space of $E$}

(b) Let $\trpl{Y}{\itPhi}{d}$ be a normed manifold, $X \subseteq Y$
be open,
$x \in \rfs{bd}(X)$ and
$\alpha \in \rfs{MBC}$.
We say that $X$ is {\it $\alpha$-linearly-bordered at $x$
$($$\alpha$-LIN-bordered$)$},
if there are a linear boundary chart domain
$\rfs{BCD}^{E}(A,r)$ and a function $\fnn{\psi}{B^E(0,r)}{Y}$
such that:
\begin{itemize}
\addtolength{\parskip}{-11pt}
\addtolength{\itemsep}{09pt}
\item[(i)]
$\iso{\psi}{B^E(0,r)}{\rfs{Rng}(\psi)}$,
\item[(ii)]
$\psi$ takes open subsets of $E$ to open subsets of $Y$
and closed subsets of $E$ to closed subsets of $Y$,
\item[(iii)]
$\psi(\rfs{BCD}^E(A,r)) = \rfs{Rng}(\psi) \cap X$,
\item[(iv)]
$\psi \nrestriction \rfs{BCD}^E(A,r)$ is $\alpha$-bicontinuous,
\item[(v)]
$\psi(0) = x$.
\vspace{-03.7pt}
\end{itemize}
$\trpl{\psi}{A}{r}$ is
called a {\it boundary chart element} for $x$.
   \index{lin-bordered0@@LIN-bordered. $\alpha$-LIN-bordered at $x$}
   \index{boundary chart element}

(c) Let $\itGamma \subseteq \rfs{MC}$.
We say that $X$ is {\it locally $\itGamma$-LIN-bordered},
if for every $x \in \rfs{bd}(X)$
there is $\alpha \in \itGamma$ such that
$X$ is $\alpha$-LIN-bordered at $x$.
   \index{locally-LIN-bordered. Locally $\itGamma$-LIN-bordered}
\hfill\proofend
\end{rm}
\end{defn}

The open sets that we had in mind when defining
LIN-borderedness are described below.
Take an open subset $U$ of $\bbR^n$ whose boundary is a
smooth submanifold. Let $K_1,\ldots,K_n$ be pairwise disjoint subsets
of $U$, and assume that for every $i$, $K_i$ is a compact
smooth submanifold of $\bbR^n$ which is not a singleton.
Then $U - \kern1pt\bigcup_{i = 1}^n K_i$
is $\itGamma^{\srfs{LIP}}$-LIN-bordered.

We remind the reader the definition of the group
$H^{\srfs{CMP.LC}}_{\itGamma}(X)$.

\begin{defn}\label{d-bddlip-1.6}
\begin{rm}
Suppose that $E,F$ are metric spaces, $X \subseteq E$
and $Y \subseteq F$, $\itGamma \subseteq \rfs{MC}$.
 
Let $\fnn{f}{X}{Y}$.
Then $f$ is {\it completely locally $\itGamma$-continuous},
if $f \in \rfs{EXT}^{E,F}(X,Y)$,
and for every $x \in \rfs{cl}^E(X)$ there are $\alpha \in \itGamma$ and
$T \in \rfs{Nbr}^E(x)$ such that
$f \nrestriction (T \cap X)$ is $\alpha$-continuous.
{\it Complete local $\itGamma$-bicontinuity} is defined analogously.

   \index{completely locally $\itGamma$-continuous}
   \index{completely locally $\itGamma$-bicontinuous}

$H^{\srfs{CMP.LC}}_{\itGamma}(X,Y;E,F)$
denotes the set of completely locally $\itGamma$-continuous
homeomorphisms between $X$ and $Y$.
We use the notation $H^{\srfs{CMP.LC}}_{\itGamma}(X,Y)$
as an abbreviation of $H^{\srfs{CMP.LC}}_{\itGamma}(X,Y;E,F)$.
The notations $(H^{\srfs{CMP.LC}}_{\itGamma})^{\pm}(X,Y)$
and $H^{\srfs{CMP.LC}}_{\itGamma}(X)$ are derived in the usual way.

   \index{N@hcmplc@@$H^{\srfs{CMP.LC}}_{\itGamma}(X) =
          \setm{g \in \rfs{EXT}(X)}{(\forall x \in \rfs{cl}^E(X))
          (\exists U \in \rfs{Nbr}^E(X))(\exists \alpha \in \itGamma)
          (g \nrestriction (U \cap X)\\\rule{27.3mm}{0pt}
          \mbox{ is } \alpha\mbox{-bicontinuous)}}$}
\end{rm}
\end{defn}

\begin{remark}\label{ams-bddly-lip-bldr-r1.3}
\begin{rm}
(a) Note that in the above definition,
if $E$ and $F$ are complete metric spaces,
then the requirement that $f \in \rfs{EXT}(X,Y)$ is not needed.

(b) In the above definition assume that $E,F$
are finite-dimensional normed spaces,
and $X,Y$ are bounded.
Let $g \in H(X,Y)$.
Then $g \in (H^{\srfs{CMP.LC}}_{\itGamma})^{\pm}(X,Y)$
iff there is $\alpha \in \itGamma$ such that $g^{\srfs{cl}}$ is
$\alpha$-bicontinuous.

(c) The motivation for dealing with groups of the type
$H^{\srfs{CMP.LC}}_{\itGamma}(X)$ is the finite-dimensional special
case described in (b).
However, the proof  of Theorem \ref{ams-bddly-lip-bldr-t1.4} below
covers other types of groups.
The following is such an example.
Let $E$ be a normed space, and $\overE$ be its completion.
Let $X \subseteq E$ be open. Denote
\newline\centerline{
$\overline{\rfs{bd}}(X) =
\rfs{cl}^{\overssE}(X) - \overline{\rfs{int}}(X)$.
}
See Definition \ref{metr-bldr-d2.23}(a).
Let
$\overline{\rfs{cl}}(X) = X \cup \overline{\rfs{bd}}(X)$.
Let
$\overH^{\srfs{CMP.LC}}_{\itGamma}(X) =
H^{\srfs{CMP.LC}}_{\itGamma}(X;\overline{\rfs{cl}}(X))$.

The proof of Theorem \ref{ams-bddly-lip-bldr-t1.4}
transfers to the group $\overH^{\srfs{CMP.LC}}_{\itGamma}(X)$
except for a slight change in the construction
of homeomorphisms in Chapter \ref{s11}.
\hfill\myqed
\end{rm}
\end{remark}

The next theorem is our main final goal. It is proved in
\ref{ams-bddly-lip-bldr-t6.1}(a).

\begin{theorem}\label{ams-bddly-lip-bldr-t1.4}
\num{a}
Let $\itGamma$ be a principal modulus of continuity
and $\itDelta$ be a modulus of continuity.
Let $E,F$ be normed spaces,
$X \subseteq E$ be a locally
$\itGamma$-LIN-bordered open set,
and $Y \subseteq F$ be a locally $\itDelta$-LIN-bordered open set.
Suppose that
$\iso{\varphi}{H^{\srfs{CMP.LC}}_{\itGamma}(X)}
{H^{\srfs{CMP.LC}}_{\itDelta}(Y)}$.
{\thickmuskip=4.5mu \medmuskip=2mu \thinmuskip=1mu 
Then $\itGamma = \itDelta$,
and there is $\tau \in (H^{\srfs{CMP.LC}}_{\itGamma})^{\pm}(X,Y)$
such that $\varphi(g) = g^{\tau}$ for every
$g \in H^{\srfs{CMP.LC}}_{\itGamma}(X)$.\kern-4pt
}

\num{b} In Part \num{a} assume that $E$ and $F$ are normed Lipschitz
manifolds. Then the Claim of Part \num{a} is true.
\end{theorem}

Part (a) is a special case of Part (b).
But we shall prove only (a),
since the setting of (b) is more complicated
and the proofs are essentially identical.
\smallskip

In the special case of bounded finite-dimensional spaces,
Theorem~\ref{ams-bddly-lip-bldr-t1.4} has a more natural formulation,
which we state in the next corollary.

\begin{cor}\label{editor-d8.5}
Let $\itGamma$ be a principal modulus of continuity,
$\itDelta$ be a modulus of continuity
and $\pair{X}{d}$ and $\pair{Y}{e}$
be compact metric Euclidean manifolds with boundary.
Assume that $\pair{X}{d}$ has an atlas
consisting of $\itGamma$-bicontinuous charts,
$\pair{Y}{e}$ has an atlas
consisting of $\itDelta$-bicontinuous charts
and $\iso{\varphi}{H_{\itGamma}(X)}{H_{\itDelta}(Y)}$.
Then $\itGamma = \itDelta$ and there is
$\iso{\tau}{X}{Y}$ such that $\tau$ is $\itGamma$-bicontinuous
and $\varphi(g) = g^{\tau}$ for every $g \in H_{\itGamma}(X)$.
\end{cor}

\noindent
{\bf Proof } The corollary follows from 
Theorem~\ref{ams-bddly-lip-bldr-t1.4}(b)
and Remark~\ref{ams-bddly-lip-bldr-r1.3}(b).
\smallskip\hfill\myqed

\noindent
{\bf Plan of the proof of Theorem~\ref{ams-bddly-lip-bldr-t1.4}(a).}

The proof of Theorem \ref{ams-bddly-lip-bldr-t1.4}(a)
has four main steps.
\newline
\rule{0pt}{0pt}
\hspace{5.7mm}
\vbox{
\begin{itemize}
\addtolength{\parskip}{-11pt}
\addtolength{\itemsep}{06pt}
\item[Step 1:]
There is $\tau \in H(X,Y)$ such that $\varphi(g) = g^{\tau}$
for every $g \in H^{\srfs{CMP.LC}}_{\itGamma}(X)$.
\item[Step 2:] $\itGamma = \itDelta$,
and $\tau$ is locally $\itGamma$-bicontinuous.
\item[Step 3:] $\tau \in \rfs{EXT}^{\pm}(X,Y)$.
\item[Step 4:] $\tau$ is completely locally $\itGamma$-bicontinuous.
\vspace{-10.1pt}
\end{itemize}
}
\rule{1pt}{0pt}The first two steps have already been accomplished.
Step 1 follows from Theorem~\ref{t2.4}
and Step 2 from Theorem~\ref{metr-bldr-t3.27}.
The exact statement of Step 3
is formulated in Theorem \ref{t-bddlip-1.8}.
The proof of this theorem takes all of Chapters \ref{s8}\,-\,\ref{s11},
and the conclusion of the proof appears at the end of Chapter \ref{s11}.
Chapter \ref{s12} is devoted to the proof of Step 4.

Theorem \ref{t-bddlip-1.8} has two variants.
Part (a) is indeed the main goal.
However, the strength of the argument is partially lost when dealing
only with groups of the type $H^{\srfs{CMP.LC}}_{\itGamma}(X)$.
Part (b) is stated in order to later reveal the full strength of the
argument.
See further explanation after the statement of
Theorem \ref{t-bddlip-1.8}.

\begin{definition}\label{ams-bddly-lip-bldr-d5.4}
\begin{rm}
(a)
Suppose that $X \subseteq E$ is open.
A subset $H \subseteq \rfs{EXT}^E(X)$ is {\it $E$-discrete},
if $\setm{\rfs{supp}(h)}{h \in H}$ is completely discrete
with respect to $E$.
(See Definition \ref{metr-bldr-d6.2.1}(a)).
Note that if $H$ is $E$-discrete,
then $\bcirc \setm{h}{h \in H} \in \rfs{EXT}^E(X)$.
   \index{discrete subset. $E$-discrete subset of $\rfs{EXT}^E(X)$}

(b)
A subgroup $G \leq \rfs{EXT}(X)$
is {\it closed under $E$-discrete composition},
if $\bcirc \setm{h}{h \in H} \in G$
for every $E$-discrete set $H \subseteq G$.
   \index{closed under $E$-discrete composition}

(c)
Let $E$ be a metric space, $X \subseteq E$ be open,
and $G \leq \rfs{EXT}(X)$. We say that $G$ is of
{\it boundary type $\itGamma$},
if for every $x \in \rfs{bd}(X)$:
\begin{itemize}
\addtolength{\parskip}{-11pt}
\addtolength{\itemsep}{06pt}
\item[(i)] there is $U \in \rfs{Nbr}^E(x)$
such that
$G\sprt{U \cap X} \supseteq
H^{\srfs{CMP.LC}}_{\itGamma}(X)\sprt{U \cap X}$,
\item[(ii)] for every $g \in G$,
there is $V \in \rfs{Nbr}^E(x)$ such that $g \nrestriction (V \cap X)$
is $\itGamma$-bicontinuous.
\vspace{-05.7pt}
\end{itemize}
   \index{boundary type. A group of boundary type $\itGamma$}
A subgroup $G \leq \rfs{EXT}(X)$ is {\it $\itGamma$-appropriate},
if $G$ is closed under $E$-discrete composition,
and $G$ is of boundary type $\itGamma$.
   \index{appropriate. A $\itGamma$-appropriate group}

(d)
Let
$H^{\srfs{BDR.LC}}_{\itGamma}(X) =
\setm{g \in \rfs{EXT}(X)}
{\mbox{for every } x \in \rfs{bd}(X),\ 
g \mbox{ is } \itGamma\mbox{-bicontinuous at } x}$.
Let $\itDelta$ be a modulus of continuity.
Define
$H^{\srfs{CMP.LC}}_{\itDelta,\itGamma}(X) =
H^{\srfs{LC}}_{\itDelta}(X) \cap H^{\srfs{BDR.LC}}_{\itGamma}(X)$.
   \index{N@hbdrlc00@@$H^{\srfs{BDR.LC}}_{\itGamma}(X) = 
          \setm{g \in \rfs{EXT}(X)} {\mbox{for every }
	  x \in \rfs{bd}(X),\ g \mbox{ is }
	  \itGamma\mbox{-bicontinuous at } x}$}
   \index{N@hbdrlc01@@ $H^{\srfs{CMP.LC}}_{\itDelta,\itGamma}(X) =
          H^{\srfs{LC}}_{\itDelta}(X) \cap
          H^{\srfs{BDR.LC}}_{\itGamma}(X)$}
\end{rm}
\end{definition}

\begin{example}\label{ams-bddly-lip-bldr-e1.6}
\begin{rm}
$H^{\srfs{CMP.LC}}_{\itGamma}(X)$ and $H^{\srfs{BDR.LC}}_{\itGamma}(X)$
are $\itGamma$-appropriate, and if $\itGamma \subseteq \itDelta$,
then $H^{\srfs{CMP.LC}}_{\itDelta,\itGamma}(X)$
is $\itGamma$-appropriate.
\end{rm}
\end{example}

\begin{theorem}\label{t-bddlip-1.8}
Let $\itGamma,\itDelta$ be countably generated moduli of continuity,
$E$ and $F$ be normed spaces and
$X \subseteq E$, $Y \subseteq F$ be open.
Suppose that $X$ is locally $\itGamma$-LIN-bordered,
and $Y$ is locally $\itDelta$-LIN-bordered
and let $\tau \in H(X,Y)$.

\num{a} If
$(H^{\srfs{CMP.LC}}_{\itGamma}(X))^{\tau} =
H^{\srfs{CMP.LC}}_{\itDelta}(Y)$,
then $\tau \in \rfs{EXT}^{\pm}(X,Y)$.

\num{b} Suppose that $G \leq \rfs{EXT}(X)$,
$H \leq \rfs{EXT}(Y)$ are respectively
$\itGamma$ and $\itDelta$ appropriate and
$G^{\tau} = H$. Then $\tau \in \rfs{EXT}^{\pm}(X,Y)$.
\end{theorem}

The proof of Theorem \ref{t-bddlip-1.8} appears at the end of
Chapter \ref{s11}.
\medskip

{\bf Explanation }
Suppose that
$(H^{\srfs{CMP.LC}}_{\itGamma}(X))^{\tau} =
H^{\srfs{CMP.LC}}_{\itDelta}(Y)$.
Then
\hbox{$\itGamma = \itDelta$.}
This is easily concluded in the following way.
Let $U \subseteq X$ be an open set such $\rfs{cl}(U) \subseteq X$
and $\rfs{cl}(\tau(U)) \subseteq Y$.
Since $\rfs{cl}(U) \subseteq X$,
$H^{\srfs{CMP.LC}}_{\itGamma}(X)\sprt{U} =
H^{\srfs{LC}}_{\itGamma}(X)\sprt{U}$.
Since $\rfs{cl}(\tau(U)) \subseteq Y$,
$H^{\srfs{CMP.LC}}_{\itDelta}(Y)\sprtl{\tau(U)} =
H^{\srfs{LC}}_{\itDelta}(Y)\sprtl{\tau(U)}$.
So
$(H^{\srfs{LC}}_{\itGamma}(X)\sprt{U})^{\tau} =
H^{\srfs{LC}}_{\itDelta}(Y)\sprtl{\tau(U)}$.
It now follows easily from Theorem~\ref{metr-bldr-t3.27} or from
Theorem~\ref{metr-bldr-t3.42}(b) that $\itGamma = \itDelta$.

When dealing with $H^{\srfs{BDR.LC}}_{\itGamma}(X)$,
the above argument is no longer valid.
Instead one has to infer that $\itGamma = \itDelta$ from the behavior
of $\tau$ at $\rfs{bd}(X)$. This is more difficult, and we  have a
proof only in special cases.
Part (b) of \ref{t-bddlip-1.8} prepares the ground for this argument.

\smallskip

As a consequence of Step 2, at the time that we reach Step 4, we
already know that $\itGamma = \itDelta$.
So the statement of Step 4 is as follows.

\begin{theorem}\label{ams-bddly-lip-bldr-t1.7}
Let $\itGamma$ be a principal modulus of continuity,
$X \subseteq E$ and
$Y \subseteq F$ be open subsets of
the normed spaces $E$ and $F$
and $\tau \in \rfs{EXT}^{\pm}(X,Y)$.
Suppose that $X$ and $Y$ are $\itGamma$-LIN-bordered
and
$(H^{\srfs{CMP.LC}}_{\itGamma}(X))^{\tau} =
H^{\srfs{CMP.LC}}_{\itGamma}(Y)$.
Then $\tau \in (H^{\srfs{CMP.LC}}_{\itGamma})^{\pm}(X,Y)$.
\end{theorem}

Chapter \ref{s12} is devoted to the proof of
Theorem \ref{ams-bddly-lip-bldr-t1.7}.
Actually, the main result of Chapter~\ref{s12} is
Theorem~\ref{ams-bddly-lip-bldr-t5.22},
and \ref{ams-bddly-lip-bldr-t1.7} is just a corollary of that theorem.
At the end of Chapter~\ref{s12} we prove
Theorem~\ref{ams-bddly-lip-bldr-t1.4}(a).
At that point it is only a matter of combining the intermediate
results from Chapters \ref{s11} and~\ref{s12}.
This is done in Theorem~\ref{ams-bddly-lip-bldr-t6.1}, 
and \ref{ams-bddly-lip-bldr-t1.4}(a) is the first part of that theorem.
\smallskip

Certain types of boundary points have to be treated differently
than others. These types are defined below.

\begin{defn}\label{d-bddlip-1.2+1}
\begin{rm}
If in Part (b) of Definition~\ref{editor-d8.1},
$A$ is a closed subspace of $E$ and $\rfs{dim}(A) = 1$,
or $\rfs{dim}(E) = 2$ and $A$ is a half space of $E$,
then we say that {\it $\rfs{bd}(X)$ is $1$-dimensional at $x$}.
   \index{dimension $1$ at $x$.
          $\rfs{bd}(X)$ is $1$-dimensional at $x$}

If in Part (b) of Definition~\ref{editor-d8.1},
$A$ is a closed subspace of $E$ and $\rfs{co-dim}(A) = 1$,
or $A$ is a half space of $E$,
then we say that {\it $\rfs{bd}(X)$ has co-dimension $1$ at $x$}.
   \index{co-dimension $1$ at $x$.
          $\rfs{bd}(X)$ has co-dimension $1$ at $x$}

If in Part (b) of \ref{editor-d8.1},
$A$ is a closed subspace of $E$ with co-dimension $1$,
then we say that {\it $X$ is two-sided at $x$.}
   \index{two-sided. $X$ is two-sided at $x$}
Hence $\rfs{Rng}(\psi) \cap X$ has two connected components.
Let $u,v \in \rfs{Rng}(\psi) \cap X$.
We say that $u,v \in X$ are {\it on the same side of $\rfs{bd}(X)$ with
respect to $\trpl{\psi}{A}{r}$},
if $u,v$ are in the same connected component of
$\rfs{Rng}(\psi) \cap X$.
We say that
$u,v \in X$ are {\it on different sides of $\rfs{bd}(X)$ with
respect to $\trpl{\psi}{A}{r}$},
if $u,v$ are in different connected components of
$\rfs{Rng}(\psi) \cap X$.
   \index{on the same side. $u,v$ are on the same side
          of $\rfs{bd}(X)$ with respect to $\trpl{\psi}{A}{r}$}
          \kern-6pt\index{on different sides.
          $u,v$ are on different sides of $\rfs{bd}(X)$
          with respect to $\trpl{\psi}{A}{r}$}

If in Part (b) of \ref{editor-d8.1}, (i) $\rfs{dim}(E) > 2$,
and (ii) $A$ is a closed subspace of $E$ of \hbox{dimension $> 1$}
or $A$ is a closed half space of $E$, then we say that $X$ is
{\it $\alpha$-simply-linearly-bordered $($$\alpha$-SLIN-bordered\,$)$
at $x$.}
   \index{slin-bordered@@SLIN-bordered.
          $\alpha$-simply-linearly-bordered at $x$
          ($\alpha$-SLIN-bordered at $x$)}
\hfill\proofend
\end{rm}
\end{defn}

Let $x \in \rfs{bd}(X)$.
Note that if
$\rfs{bd}(X)$ is $1$-dimensional at $x$,
and $\trpl{\psi}{A}{r}$ is any boundary chart element for $x$,
then either (i) $A$ is a $1$-dimensional subspace,
or (ii) $\rfs{dim}(E) = 2$ and $A$ is a closed half space.
Similarly, if $X$ is two-sided at $x$,
and $\trpl{\psi}{A}{r}$ is any boundary chart element for $x$,
then $A$ is a closed subspace with co-dimension $1$.

\begin{question}\label{editor-q8.10}
\begin{rm}
A subset $A \subseteq E$ is called a {\it closed half subspace} of $E$,
if there is a closed subspace $F$ of $E$
such that $F \neq \sngltn{0}$ and $A$ is a half space of $F$.
   \index{closed half subspace of a normed space}
Let $\rfs{BCD}^E(A,r)$ be a boundary chart domain.
We call $\rfs{BCD}^E(A,r)$ an
{\it almost linear boundary chart domain} if either
it is a linear boundary chart domain,
or $A$ is a closed half subspace of $E$.
   \index{almost linear boundary chart domain}
Let $\itGamma \subseteq \rfs{MC}$. Define the notion
{\it ``$X$ is locally $\itGamma$-almost-linearly-bordered''
(locally $\itGamma$-ALIN-bordered)}
in analogy with Definition~\ref{editor-d8.1}(c).

Are Theorems~\ref{t-bddlip-1.8} and \ref{ams-bddly-lip-bldr-t1.7}
true for locally ALIN-bordered sets?

In order to prove the analogues of
\ref{t-bddlip-1.8} and \ref{ams-bddly-lip-bldr-t1.7}
for locally ALIN-bordered sets,
only Lemma \ref{l-bddlip-1.10} needs to be generalized.
All other ingredients in the proof remain essentially the same.
\hfill\proofend
\end{rm}
\end{question}

Some ALIN-bordered sets are described below.
Take an open subset $U$ of $\bbR^n$ whose boundary is a
smooth submanifold.
Let $K_1,\ldots,K_n$ be pairwise disjoint subsets
of $U$, and assume that for every $i$, $K_i$ is a compact
manifold with a boundary which is not a singleton,
and $K_i$ is smoothly embedded in $\bbR^n$.
Then $U - \kern1pt\bigcup_{i = 1}^n K_i$
is $\itGamma^{\srfs{LIP}}$-ALIN-bordered.

\newpage

\section{The Uniform Continuity Constant}
\label{s9}

\subsection{Preliminary lemmas about the existence of
certain constants.}
\label{ss9.1}

In preparing the ground for the proof of Theorem \ref{t-bddlip-1.8},
we need to characterize the pairs of convergent sequences
$\vecx,\vecy$ in $X$ for which there is an
$\alpha$-bicontinuous homeomorphism $g \in H(X)$
and subsequences $\vecx\fprime,\vecy\fprime$ of $\vecx$ and $\vecy$
such that $g(\vecx\fprime) = \vecy\fprime$.
Stated more precisely,
let $z \in \rfs{bd}(X)$ and $\lim \vecx = \lim \vecy = z$,
and assume that for every $n \in \bbN$,
\begin{itemize}
\addtolength{\parskip}{-11pt}
\addtolength{\itemsep}{06pt}
\item[(1)] $\norm{x_n - z} \leq \alpha(\norm{y_n - z})$
\ \ and\ \ \ %
$\norm{y_n - z} \leq \alpha(\norm{x_n - z})$,
\item[(2)] $d(x_n,\rfs{bd}(X)) \leq \alpha(d(y_n,\rfs{bd}(X)))$
\ \ and\ \ \ %
$d(y_n,\rfs{bd}(X)) \leq \alpha(d(x_n,\rfs{bd}(X)))$.
\vspace{-05.7pt}
\end{itemize}
We shall prove that there are $g \in H(X)$ and
subsequences $\vecx\fprime$ and $\vecy\fprime$ of $\vecx$ and $\vecy$
respectively such that $g(\vecx\fprime) = \vecy\fprime$ and $g$ is
$N \ncdot \alpha \scirc \alpha \scirc \alpha \scirc
\alpha$-bicontinuous.
In fact, this is only an approximation of what we really prove.
The exact statement to be proved is the equivalence between
the conjunction of (1) and (2) above
and the fact that $\vecx \neweq^{N\alpha^4} \vecy$.
The relation $\neweq^{\alpha}$ is defined in
\ref{d3.1-bddly-lip-extending}(c),
and in Proposition~\ref{p3.3-bddly-lip-extending}(a)
we prove this equivalence.

The Uniform Continuity Constant Lemma \ref{l-bddlip-1.10}
is the main fact needed in the proof of the above.
It says that there is $K > 0$ for which
$\bfA \Rightarrow \bfB$, where $\bfA$ and $\bfB$ are the following
statements.
\newline
$(\bfA)$ \ $E$ is a normed vector space, $F$ is a closed subspace of $E$
with dimension $> 1$,
$\alpha \in \rfs{MBC}$,
$x,y \in E - F$,
$\norm{x} \leq \norm{y} \leq \alpha(\norm{x})$ and
$\alpha\inverse(d(x,F)) \leq d(y,F) \leq \alpha(d(x,F))$.
\newline
$(\bfB)$ \ There is an $K \ncdot \alpha \scirc \alpha$-bicontinuous
homeomorphism $g$ such that: $g(x) = y$, $g(F) = F$ and
$\rfs{supp}(g) \subseteq
B(0,2 \norm{y}) - \overB(0,\dgfrac{\norm{x}}{2})$.

This chapter is devoted to the proof of this lemma.
The geometric contents of the lemma is simple,
but a detailed proof seems to require much work.
When the claim of the lemma is restricted to pre-Hilbert spaces
and not to general normed spaces, the proof is easier.

We shall also need a statement analogous to 
$\bfA \Rightarrow \bfB$ for subspaces $F$ of $E$
with \hbox{$\rfs{dim}(F) = 1$.}
In this case Statements $\bfA$ and $\bfB$ need to be slightly
modified.
Chapter~\ref{s10} deals with this situation.
\smallskip

Before turning to the proof of The Uniform Continuity Constant Lemma
we quote some well-known basic facts from functional analysis,
and we also establish the existence of various types of homeomorphisms
which will be used in the proof of \ref{l-bddlip-1.10}.
These prepartions are carried out
in \ref{n-bddlip-n2.1}\,-\,\ref{p-bddlip-1.11-1}.
We start with some needed notation.

\begin{notation}\label{n-bddlip-n2.1}
\begin{rm}
(a) For $K \geq 1$ and $a,b \geq 0$ let $a \approx^K b$ mean that
$\dgfrac{a}{K} \leq b \leq Ka$.
   \index{N@AAAA@@$a \approx^K b$. This means
	  $\frac{1}{K} a \leq b \leq Ka$}
   \index{N@AAAA@@$\norm{\ }^1 \approx^K \norm{\ }^2$.
	  This means for every $u \in E$,
	  $\norm{u}^1 \approx^K \norm{u}^2$}
If $\norm{\ }^1,\norm{\ }^2$ are norms on a vector space $E$,
then $\norm{\ }^1 \approx^K \norm{\ }^2$
means that $\norm{u}^1 \approx^K \norm{u}^2$ for every $u \in E$.

   \index{N@AAAA@@$E = L \oplus^{\srfs{alg}} S$.
	  The algebraic direct sum}
(b) The notation $E = L \oplus^{\srfs{alg}} S$ means that
$L + S = E$ and $L \cap S = \sngltn{0}$. 
If $E = L \oplus^{\srfs{alg}} S$, then
$(x)_{L,S}, (x)_{S,L}$ denote the components of $x$ in $L$ and $S$
respectively.
In what follows we sometimes abbreviate $(x)_{L,S}$ by $(x)_L$
and $(x)_{S,L}$ by $(x)_S$.
   \index{N@AAAA@@$(x)_{L,S}$. The $L$-component of $x$ in
	  $L \oplus S$}
   \index{N@AAAA@@$(x)_L$. Abbreviation of $(x)_{L,S}$}
Suppose that $E = L \oplus^{\srfs{alg}} S$.
We define $\norm{u}^{L,S} = \norm{(u)_S} + \norm{(u)_L}$.
   \index{N@AAAA@@$\norm{u}^{L,S} =
	  \norm{(u)_S} + \norm{(u)_L}$}
The notation $E = L \oplus S$ means that $E = L \oplus^{\srfs{alg}} S$,
and that for some $K \geq 1$,
$\norm{\ }^{L,S} \approx^K \norm{\ }$.
In such a case $S$ is called a {\it complement} of $L$ in $E$.

(c) Let $L$ be a linear subspace of $E$.
Then $\rfs{co-dim}^E(L)$ denotes the co-dimension of $L$ in $E$.
This is abbreviated by $\rfs{co-dim}(L)$.
   \index{N@codim@@$\rfs{co-dim}^E(L)$. Co-dimension of $L$ in $E$.
	  Abbreviation $\rfs{co-dim}(L)$}

(d) Let $F$ and $H$ be linear subspaces of a normed space
$E$ and $M \geq 1$.
   \index{N@AAAA@@$H \perp^M F$.
	  This means for every $u \in H$,
          $d(u,F) \geq \frac{1}{M} \norm{u}$}
We define $H \perp^M F$
if $d(h,F) \geq \dgfrac{\norm{h}}{M}$ for every $h \in H$.

(e) Let $E = F \oplus^{\srfs{alg}} H$.
Then $\rfs{Proj}_{F,H}$ is denotes the function
$u \mapsto (u)_{F,H}$, \,$u \in E$.

(f)
   \index{N@b04@@$B(x;r,s) = \setm{y \in X}{r < d(x,y) < s}$}
Let $X$ be a metric space, $x \in X$ and $0 < r < s$.
The ring with center at $x$ and with radii $r,s$ is defined as
$$
B(x;r,s) = \setm{y \in X}{r < d(x,y) < s}.
$$
\end{rm}
\end{notation}

We quote without proof some basic and well-known facts from functional
analysis.

\begin{prop}\label{p-bddlip-bldr-1.9}
\num{a} For every $n > 0$ there is $M = M^{\srfs{aoc}}(n) \geq 1$
such that for every normed space $E$
and an $n$-dimensional subspace $L$ of $E$
there is a complement $S$ of $L$ in $E$ such that
$M \norm{x} \geq \norm{(x)_{L,S}} + \norm{(x)_{S,L}}$
\,for every $x \in E$.
A subspace $S$ satisfying the above is called an {\rm almost orthogonal
complement} of $L$.
   \index{N@maoc@@$M^{\srfs{aoc}}(n)$}
   \index{almost orthogonal complement}

\num{b} For every $n > 0$ there is $M = M^{\srfs{thn}}(n) \geq 1$
such that for every normed $n$-dimensional space $E$
there is a Hilbert norm $\norm{\ }^{\srbfs{H}}$ on $E$
such that $\norm{x} \leq \norm{x}^{\srbfs{H}} \leq M \norm{x}$\break
for every $x \in E$.
The norm $\norm{\ }^{\srbfs{H}}$ is called a
{\rm tight Hilbert norm} on $E$.
We denote $M^{\srfs{thn}}(2)$ by $M^{\srfs{thn}}$.

   \index{N@mthn@@$M^{\srfs{thn}}(n)$}
   \index{N@mthn@@$M^{\srfs{thn}} = M^{\srfs{thn}}(2)$}
   \index{tight Hilbert norm}
   \index{N@mhlb@@$M^{\srfs{hlb}}(n)$}
   \index{N@mhlb@@$M^{\srfs{hlb}} = M^{\srfs{hlb}}(2)$}
\num{c} For every $n > 0$ there is
$M = M^{\srfs{hlb}}(n) \geq 1$
such that for every normed space $E$ and an $n$-dimensional
linear subspace $L$ of $E$ there are a Euclidean norm
$\norm{\ }^{\srbfs{H}}$ on $L$ and a complement $S$ of $L$ such that
for every $x \in E$,
$$
\norm{(x)_{L,S}}^{\srbfs{H}} + \norm{(x)_{S,L}} \approx^M \norm{x}.
$$
Also, if $m < n$, then $M^{\srfs{hlb}}(m) \leq M^{\srfs{hlb}}(n)$.
A pair $\pair{\norm{\ }^{\srbfs{H}}}{S}$ satisfying the above is called
a {\rm tight Hilbert complementation} for $L$.
We denote $M^{\srfs{hlb}}(2)$ by $M^{\srfs{hlb}}$.
   \index{tight Hilbert complementation}

\num{d} Let $E = F \oplus H$ and $M \geq 1$.
Then $H \perp^M F$ iff $\norm{\rfs{Proj}_{H,F}} \leq M$.

\num{e} Let $E = F \oplus H$ and suppose that $H \perp^M F$.
Then $F \perp^{M + 1} H$.

\num{f} Let $E = F \oplus H$ and suppose that $H \perp^M F$.
Then $\norm{\ }^{F,H} \approx^{2 M + 1} \norm{\ }$.

\num{g} Let $E = F \oplus H$ and suppose that
$\norm{\ }^{F,H} \approx^M \norm{\ }$.
Then $H \perp^M F$.

\num{h} Let $\fnn{T}{E}{E}$ be a bounded linear projection with
a closed range. Then\break
$\rfs{ker}(T) \perp^{\norm{T} + 1} \rfs{Rng}(T)$.

\num{i} Let $x,y \in E - \sngltn{0}$ be such that
$\norm{x} \leq \norm{y}$. Let $z = \frac{\norm{x}}{\norm{y}} y$.
Then $\norm{y - z} \leq \norm{y - x}$
and $\norm{x - z} \leq 2 \norm{y - x}$.
\end{prop}

\begin{prop}\label{p-bddlip-1.11}
{\thickmuskip=3.3mu \medmuskip=2.2mu \thinmuskip = 1.1mu
Let $F$ be a closed subspace of a normed vector space $E$,
$x,y \in E - F$ and $\varepsilon > 0$.
Then there is a closed subspace $H$ of $E$ such that $F \subseteq H$,
$\rfs{span}(H \cup \dbltn{x}{y}) = E$,
$d(x,H) \geq \frac{1}{1 + \varepsilon} d(x,F)$
and $d(y,H) \geq \frac{1}{1 + \varepsilon} d(y,F)$.
}
\end{prop}

\noindent
{\bf Proof }
Let $\itDelta = 1 + \varepsilon$ and $\hatx \in F$ be such that
$\norm{x - \hatx} \leq \itDelta d(x,F)$.
Denote $x^{\nperp} = x - \hatx$.
Let $\psi$ be the linear functional on
$\rfs{span}(F \cup \sngltn{x})$ defined by
$\psi(x^{\nperp}) = \norm{x^{\nperp}}$
and $\psi(F) = \sngltn{0}$.
We check that $\norm{\psi} \leq \itDelta$.
Let $z \in \rfs{span}(F \cup \sngltn{x})$. If $z \in F$,
then $\abs{\psi(z)} = 0 \leq \itDelta \norm{z}$.
Suppose that $z = u + \lambda x^{\nperp}$, where $u \in F$ and
$\lambda \neq 0$. We may assume that $\lambda = 1$.
Then
$$
\abs{\psi(z)} = \norm{x^{\nperp}} \leq \itDelta d(x,F) \leq
\itDelta \norm{(u - \hatx) + x} =
\itDelta \norm{u + x^{\nperp}} =
\itDelta \norm{z}.
$$
Let $\varphi \in E^*$ be such that $\varphi$ extends $\psi$ and
$\norm{\varphi} = \norm{\psi}$.
Let $H_1 = \rfs{ker}(\varphi)$. So $F \subseteq H_1$.
Since $x = \hatx + x^{\nperp}$ and $\hatx \in H_1$,
$d(x,H_1) = d(x^{\nperp},H_1)$.
Let $u \in H_1$.
Then
$$
\norm{x^{\nperp} - u} \geq
\frac{\abs{\varphi(x^{\nperp} - u)}}{\itDelta} =
\frac{\norm{x^{\nperp}}}{\itDelta} \geq
\frac{d(x,F)}{\itDelta}.
$$
Hence
$d(x,H_1) = d(x^{\nperp},H_1) \geq
\frac{d(x,F)}{1 + \varepsilon}$.

Similarly,
there is a closed linear subspace $H_2$ with co-dimension $1$
such that\break
$d(y,H_2) \geq \frac{d(y,F)}{1 + \varepsilon}$.
Let $H = H_1 \cap H_2$. Then $H$ is as required.
\smallskip\hfill\myqed

The next proposition contains some additional basic and well-known
facts from functional analysis. The proofs are again omitted.

\begin{prop}\label{p-bddlip-1.12}
\num{a}
   \index{N@mprj@@$M^{\srfs{prj}}(n)$}
For every $n \in \bbN$ there is $M^{\srfs{prj}}(n)$ such that for every
normed space $E$ and a closed linear subspace $F \subseteq E$: if
$\rfs{co-dim}^E(F) = n$, then there is a projection
$\fnn{T}{E}{F}$ such that $\norm{T} \leq M^{\srfs{prj}}(n)$.

\num{b}
   \index{N@mort@@$M^{\srfs{ort}}(n)$}
   \index{N@mort@@$M^{\srfs{ort}} = M^{\srfs{ort}}(2)$}
For every $n \in \bbN$ there is $M = M^{\srfs{ort}}(n)$
such that for every normed space $E$ and a closed linear subspace
$F \subseteq E$: if $\rfs{co-dim}^E(F) \leq n$,
then there is a closed linear subspace $H \subseteq E$ such that 
$F \oplus H = E$ and $H \perp^M F$.
One can take $M^{\srfs{ort}}(n)$ to be $2^n - 1 + \varepsilon$
for any $\varepsilon > 0$.
We denote $M^{\srfs{ort}}(2)$ by $M^{\srfs{ort}}$.

\num{c} Let
   \index{N@mfdn@@$M^{\srfs{fdn}}(n)$}
   \index{N@mfdn@@$M^{\srfs{fdn}} = M^{\srfs{fdn}}(2)$}
$M^{\srfs{fdn}}(n) =
(1 + M^{\srfs{thn}}(n)) \cdot M^{\srfs{ort}}(n) + 1$.
Let $E$ be a normed space, $F \subseteq E$ be a closed subspace with
co-dimension $\leq n$ and $H$ be such that
$F \oplus H = E$ and $H \perp^{M^{\ssrfs{ort}}(n)} F$.
Let $\norm{\ }^{\sbfs{H}}$ be a Hilbert norm on $H$ such that
$\norm{\ }^{\sbfs{H}} \approx^{M^{\ssrfs{thn}}(n)}
\norm{\ } \,\nrestriction H$.
Define a new norm on $E$ by\,
$\norm{u}^{\sbfs{N}} = \norm{(u)_F} + \norm{(u)_H}^{\sbfs{H}}$.
Then $\norm{\ }^{\sbfs{N}} \approx^{M^{\ssrfs{fdn}}(n)} \norm{\ }$.
We denote $M^{\srfs{fdn}}(2)$ by $M^{\srfs{fdn}}$.
\end{prop}

\begin{defn}\label{d-bddlip-1.13}
\begin{rm}
(a) Let $H$ be a $2$-dimensional Hilbert space and $\theta \in \bbR$.
Then $\rfs{Rot}^H_{\theta}$ denotes the rotation by an angle of $\theta$
in $H$.
Let $E = F \oplus H$ be normed spaces.
Suppose that $H$ is a $2$-dimensional Hilbert space.
For $\theta \in \bbR$ let $\rfs{Rot}^{F,H}_{\theta} \in H(E)$
be defined by
$$\rfs{Rot}^{F,H}_{\theta}(u) = (u)_F + \rfs{Rot}^{H}_{\theta}((u)_H),
\ u \in E.
$$

   \index{N@rot00@@$\rfs{Rot}^H_{\theta}$.
	  In a $2$-dimensional Hilbert space $H$,
	  rotation by the angle $\theta$}
   \index{N@rot01@@$\rfs{Rot}^{F,H}_{\theta}$. For a
	  $2$-dimensional Hilbert space $H$ and a normed space $F$,
	  the operator on\\
	  \indent$H \oplus F$ which is
	  $\rfs{Rot}^H_{\theta}$ on $H$ and $\rfs{Id}$ on $F$}

(b) Let $h = \rfs{Rad}^E_{\eta,z}$ be a radial homeomorphism.
(See Definition~\ref{metr-bldr-d3.17}(b)).
We say that $h$ is {\it piecewise linearly radial},
if $\eta$ is piecewise linear.
   \index{piecewise linearly radial. A radial homeomorphism
          $\rfs{Rad}^E_{\eta}$ in which $\eta$ is piecewise linear}
\hfill\proofend
\end{rm}
\end{defn}

Part (a) of the following proposition is a variant Lemma \ref{l2.6}(c).

\begin{prop}\label{p-gamma.7}
\num{a} There is $M^{\srfs{seg}} > 1$ such that for every normed
   \index{N@mseg@@$M^{\srfs{seg}}$}
space~$E$, $x,y \in E$ and  $r > 0$,
there is $h \in H(E)$ such that
\begin{itemize}
\addtolength{\parskip}{-11pt}
\addtolength{\itemsep}{06pt}
\item[\num{1}] $\rfs{supp}(h) \subseteq B([x,y],r)$,
\item[\num{2}] $h(x) = y$,
\item[\num{3}] $h$ is
$M^{\srfs{seg}} \ncdot (\dgfrac{\norm{x - y}}{r} + 1)$-bilipschitz.
\vspace{-05.7pt}
\end{itemize}

\num{b} For every $t > 0$
there is $M^{\srfs{arc}}(t) > 1$ such that for every normed space $E$,
   \index{N@marc@@$M^{\srfs{arc}}(t)$}
a rectifiable arc $L \subseteq E$ with endpoints $x,y$
and  $r > 0$
there is $h \in H(E)$ such that
\begin{itemize}
\addtolength{\parskip}{-11pt}
\addtolength{\itemsep}{06pt}
\item[\num{1}] $\rfs{supp}(h) \subseteq B(L,r)$,
\item[\num{2}] $h(x) = y$,
\item[\num{3}] $h$ is
$M^{\srfs{arc}}(\frac{\srfs{lngth}(L)}{r})$-bilipschitz.
\vspace{-05.7pt}
\end{itemize}

\num{c} There is $M^{\srfs{rot}} \geq 1$ such that the following holds.
   \index{N@mrot@@$M^{\srfs{rot}}$}
Let $E = F \oplus H$ be normed spaces.
Suppose that $H$ is a $2$-dimensional Hilbert space, and that for every
$u \in E$,\break
$\norm{u} = \norm{(u)_F} + \norm{(u)_H}$.
Let $S$ be a closed subset of $E$,
$\fnn{\eta}{[0,\infty)}{\bbR}$,
and $K,r > 0$ be such that:
\num{i} $S \subseteq \barB(0,r)$;
\num{ii} for every $u \in S$ and $\theta \in \bbR$,
$\rfs{Rot}^{F,H}_{\theta}(u) \in S$;
\num{iii} $\eta$ is $K$-Lipschitz;
\num{iv} $\eta(s) = 0$ for every $s \geq r$.
Let
$\fnn{g}{E}{E}$ be defined by
$g(u) = \rfs{Rot}^{F,H}_{\eta(d(u,S))}(u)$.
Then $g \in H(E)$ and $g$ is
$(M^{\srfs{rot}} Kr + 1)$-bilipschitz.

\num{d} Suppose that $F,H$ are normed spaces, $E = F \oplus H$,
and $\norm{u + v} = \norm{u} + \norm{v}$
for every $u \in F$ and $v \in H$.
Let $\hatx \in F$, $x \in H$, $a > 1$,
$x' = \hatx + x$ and $x'' = \hatx + ax$.
Then there is $g \in H(E)$ such that
\begin{itemize}
\addtolength{\parskip}{-11pt}
\addtolength{\itemsep}{06pt}
\item[\num{1}] $g(x') = x''$,
\item[\num{2}] $g \nrestriction F = \rfs{Id}$,
\item[\num{3}] for every $u \in F$,
$\rfs{supp}(g) \subseteq
B(u;s,t)$,
where
$s = \dgfrac{\norm{x' - u}}{2}$ and $t = \dgfrac{3 \norm{x'' - u}}{2}$.
\item[\num{4}] $g$ is $2 M^{\srfs{seg}} a$\,-\,bilipschitz.
\vspace{-05.7pt}
\end{itemize}
\end{prop}

\noindent
{\bf Proof }
(a)
Set $\barx = \dgfrac{x}{\norm{x}}$
and $a = \norm{x - y}$. We may place the origin in such a way that
$x = (\dgfrac{r}{2}) \ncdot \barx$ and
$y = (\dgfrac{r}{2} + a) \ncdot \barx$.
We may assume that $r < a$.
Write $M =  M^{\srfs{aoc}}(1)$.
Let $L = \rfs{span}(\sngltn{x})$ and $S$ be a complement of $L$
such that $M \norm{u} \geq \norm{(u)_{L,S}} + \norm{(u)_{S,L}}$
for every $u \in E$.
So for every $u \in S$,
$\norm{u} \leq M \ncdot d(u,L)$.
Denote $(u)_{L,S} = \hatu$ and $(u)_{S,L} = u^{\sperp}$.
For every $u \in E$
let $\lambda_u$ be such that $\hatu = \lambda_u \barx$.
So $u = \lambda_u \barx + u^{\sperp}$.

Let $g(s,t) = g_s(t), \ s \geq 0, \ t \in \bbR$, be defined as follows.
For every $s \geq 0$,
$g_s(t)$ is a piecewise linear function satisfying the following.
\begin{itemize}
\addtolength{\parskip}{-11pt}
\addtolength{\itemsep}{06pt}
\item[(1)] The breakpoints of $g_s(t)$
are $0$, $\dgfrac{r}{2}$ and $a + r$.
\item[(2)] If $s \in [0,\frac{r}{2M}]$, then
$g_s(\dgfrac{r}{2}) =$
\lower0.7pt\hbox{{\large
$\frac{\frac{r}{2M} - s}{\frac{r}{2M}}$}$\cdot (a + \dgfrac{r}{2})
$,
}
and if $s \geq \frac{r}{2M}$,
then $g_s(\dgfrac{r}{2}) = \dgfrac{r}{2}$.
\item[(3)] If $t \leq 0$ or $t \geq a + r$, then $g_s(t) = t$.
\vspace{-05.7pt}
\end{itemize}
So $g_0(\dgfrac{r}{2}) = a + \dgfrac{r}{2}$
and $g_s = \rfs{Id}$ for every $s \geq \frac{r}{2M}$.
Define
$$
h(u) = u^{\sperp} + g_{d(u,L)}(\lambda_u) \cdot \barx.
$$
Clearly, $h(x) = y$.
Let $u \in E - B([x,y],r)$, and we prove that $h(u) = u$.
If $d(u,L) \geq \frac{r}{2M}$, then $g_{d(u,L)} = \rfs{Id}$.
So $h(u) = u^{\sperp} + \lambda_u \barx = u$.
Assume that $d(u,L) < \frac{r}{2M}$.
If $\lambda_u \leq 0$,
then for every $s$,
$g_s(\lambda_u) = \lambda_u$ and hence $h(u) = u$.
Assume that $\lambda_u > 0$.
Since $d(u,L) < \frac{r}{2M}$,
it follows that $\norm{u^{\sperp}} < \dgfrac{r}{2}$. Hence
$$
\abs{\lambda_u - (a + \dgfrac{r}{2})} = \norm{\hatu - y} \geq
\norm{u - y} - \norm{u^{\sperp}} > r - \dgfrac{r}{2} = \dgfrac{r}{2}.
$$
That is, either (i) $\lambda_u - (a + \dgfrac{r}{2}) > \dgfrac{r}{2}$ or
(ii) $\lambda_u - (a + \dgfrac{r}{2}) < -\dgfrac{r}{2}$.
Suppose by contradiction that (ii) happens.
Then $0 < \lambda_u < a$. If $\lambda_u \geq \dgfrac{r}{2}$,
then $\hatu = \lambda_u \barx \in [x,y]$, and hence
$d(u,[x,y]) \leq \norm{u - \hatu} = \norm{u^{\sperp}} < \dgfrac{r}{2}$.
So $u \in B([x,y],r)$, a contradiction.
If $\lambda_u < \dgfrac{r}{2}$, then
$d(u,[x,y]) \leq \norm{x - u} \leq
\norm{x - \hatu} + \norm{u^{\sperp}} <
\dgfrac{r}{2} + \dgfrac{r}{2} = r$.
So $u \in B([x,y],r)$.
A contradiction.
Hence $\lambda_u - (a + \dgfrac{r}{2}) > \dgfrac{r}{2}$.
So $\lambda_u > a + r$, and hence for every~$s$,
$g_s(\lambda_u) = \lambda_u$. So $h(u) = u$.
We have shown that $h \nrestriction (E - B([x,y],r)) = \rfs{Id}$.

For every $s \geq 0$ let $f_s = g_s\inverse$,
and let $f(s,t) = f_s(t)$.
Note that for every $u \in E$, $u^{\sperp} = (h(u))^{\sperp}$,
and hence $d(h(u),L) = d(u,L)$.
So if $w = h(u)$,
then $u = w^{\sperp} + f_{d(w,L)}(\lambda_w) \lcdot \barx$.
Hence $h\inverse$ exists and is continuous,
and so $h \in H(E)$.

We show that $h$ and $h\inverse$ are Lipschitz.
Note that for every $s$, the three slopes of $g_s$ are
$\leq \frac{a + \dgfrac{r}{2}}{\dgfrac{r}{2}}$.
Also, for every $s_1,s_2 \geq 0$  and $t \in \bbR$,
$\abs{g_{s_1}(t) - g_{s_2}(t)} \leq
\frac{a}{\frac{r}{2M}} \ncdot \abs{s_1 - s_2}$.
For $f_s$, the maximal slope is again
$\frac{a + \dgfrac{r}{2}}{\dgfrac{r}{2}}$ and
$\abs{f_{s_1}(t) - f_{s_2}(t)} \leq
\frac{a}{\frac{r}{2M}} \ncdot \abs{s_1 - s_2}$. Now
$$
h(u) - h(v) = u^{\sperp} - v^{\sperp} +
\left(\rule{0pt}{10pt}
g_{d(u,L)}(\lambda_u) - g_{d(u,L)}(\lambda_v)
\right) \barx +
\left(\rule{0pt}{10pt}
g_{d(u,L)}(\lambda_v) - g_{d(v,L)}(\lambda_v)
\right) \barx.
$$
Write $w = u - v$.
So
\vspace{1.5mm}
\newline
\rule{-2pt}{0pt}
\renewcommand{\arraystretch}{1.5}
\addtolength{\arraycolsep}{-0pt}
$
\begin{array}{ll}
&
\norm{h(u) - h(v)} \leq \norm{u^{\sperp} - v^{\sperp}} +
\abs{g_{d(u,L)}(\lambda_u) - g_{d(u,L)}(\lambda_v)} +
\abs{g_{d(u,L)}(\lambda_v) - g_{d(v,L)}(\lambda_v)}
\\
\leq
&
\norm{w^{\sperp}} +
\frac{a + \dgfrac{r}{2}}{\dgfrac{r}{2}} \norm{\hatu - \hatv} +
\frac{a}{\frac{r}{2M}} \ncdot (d(u,L) - d(v,L)) \leq
\norm{w^{\sperp}} + (\dgfrac{2a}{r} + 1) \ncdot \norm{\hatw} +
\norm{u - v}
\\
\leq
&
M \norm{u - v} + (\dgfrac{2a}{r} + 1) M \norm{u - v} +
\norm{u - v} \leq
(3M + 1)(\dgfrac{a}{r} + 1) \norm{u - v}.
\vspace{1.7mm}
\end{array}
$
\renewcommand{\arraystretch}{1.0}
\addtolength{\arraycolsep}{0pt}
\newline
\smallskip
An identical computation shows that
$h\inverse$ is $(3M + 1) (\dgfrac{a}{r} + 1)$-Lipschitz.
So $M^{\srfs{seg}} = 3M + 1$.

(b) Let $E$ be a normed space, $L \subseteq E$ be a rectifiable
arc with endpoints $x,y$ and $r > 0$.
Denote $\ell = \rfs{lngth}(L)$ and $n = [\dgfrac{\ell}{r}] + 1$.
There are $x_i \in L, \ i = 0,\ldots,n$ such that
$x_0 = x$, $x_n = y$ and for every $i < n$,
$\norm{x_i - x_{i + 1}} \leq r$.
For $i < n$ let $L_i = [x_i,x_{i + 1}]$.
Then $B(L_i,\dgfrac{r}{2}) \subseteq B(L,r)$.
By Part (a), there is $g_i \in H(E)$ such that
\begin{itemize}
\addtolength{\parskip}{-11pt}
\addtolength{\itemsep}{06pt}
\item[\num{1}] $\rfs{supp}(g_i) \subseteq B(L_i,\dgfrac{r}{2})$,
\item[\num{2}] $g_i(x_i) = x_{i + 1}$,
\item[\num{3}] $g_i$ is
$M^{\srfs{seg}} \ncdot
(\frac{\norm{x_i - x_{i + 1}}}{\dgfrac{r}{2}} + 1)$-bilipschitz.
\vspace{-05.7pt}
\end{itemize}
Since $\norm{x_i - x_{i + 1}} \leq r$ and by (3),
$g_i$ is $3 M^{\srfs{seg}}$-bilipschitz.
Let
$M^{\srfs{arc}}(t) = (3 M^{\srfs{seg}})^{[t] + 1}$.\break
Define $g = g_0 \scirc \ldots \scirc g_{n - 1}$.
It is easily seen that $g(x) = y$, $\rfs{supp}(g) \subseteq B(L,r)$
and $g$ is\break
$M^{\srfs{arc}}(\dgfrac{\ell}{r})$-bilipschitz.

(c) Suppose that some function $\fnn{f}{E}{E}$ has the property that
for some $a > 0$, $\norm{f(u) - f(v)} \leq M\,\norm{u - v}$
for every $u,v \in E$ such that $\norm{u - v} \leq a$.
Then $f$ is $M$-Lipschitz.
For the function $g$ we take $a$ to be $r$.
Let $u,v \in E$ be such that $\norm{u - v} \leq r$.
If $v \not\in B(0,3r)$, then $u \not\in B(0,2r)$. So $g(u) = u$ and
$g(v) = v$.
We may thus assume that $\norm{v} < 3r$.
Denote $(u)_H,(u)_F,(v)_H,(v)_F$ by $u_1,u_2,v_1,v_2$ respectively
and
$\theta(w) \eqdf \eta(d(w,S))$.
Then
\vspace{1.5mm}
\newline
\rule{13pt}{0pt}
\renewcommand{\arraystretch}{1.5}
\addtolength{\arraycolsep}{-3pt}
$
\begin{array}{lll}
\rule{0pt}{0pt}
g(v) - g(u)\kern5pt = &
\kern6pt
(\rfs{Rot}^H_{\theta(v)}(v_1) -
\rfs{Rot}^H_{\theta(u)}(v_1))\kern5pt + &
\\
\rule{0pt}{14pt}&
\kern6pt
(\rfs{Rot}^H_{\theta(u)}(v_1) -
\rfs{Rot}^H_{\theta(u)}(u_1)) \kern5pt +
&
(v_2 - u_2).
\end{array}
$
\renewcommand{\arraystretch}{1.0}
\addtolength{\arraycolsep}{3pt}
\newline
So
\vspace{1.5mm}
\newline
\rule{7pt}{0pt}
\renewcommand{\arraystretch}{1.5}
\addtolength{\arraycolsep}{0pt}
$
\begin{array}{llll}
\rule{0pt}{0pt}\hspace{-2.0mm}
\norm{g(v) - g(u)}& \kern0pt \leq &
\norm{(\rfs{Rot}^H_{\theta(v)}(v_1) -
\rfs{Rot}^H_{\theta(u)}(v_1))}\kern5pt + &
\\
\rule{0pt}{14pt}&&
\norm{(\rfs{Rot}^H_{\theta(u)}(v_1) -
\rfs{Rot}^H_{\theta(u)}(u_1)) \kern5pt +
&
\kern-9pt
(v_2 - u_2)}
\\
&
\kern2pt=
&
\rule{0pt}{14pt}
\norm{(\rfs{Rot}^H_{\theta(v)}(v_1) -
\rfs{Rot}^H_{\theta(u)}(v_1))}\kern5pt +&
\kern-6pt
\norm{v - u}.
\end{array}
$
\newline
\renewcommand{\arraystretch}{1.0}
\addtolength{\arraycolsep}{0pt}
\rule{0pt}{14pt}\kern-2pt
We estimate the first summand in the last expression.
\vspace{1.5mm}
\newline
\rule{4pt}{0pt}
\renewcommand{\arraystretch}{1.5}
\addtolength{\arraycolsep}{0pt}
$
\begin{array}{ll}
&
\norm{(\rfs{Rot}^H_{\theta(v)}(v_1) -
\rfs{Rot}^H_{\theta(u)}(v_1))} \leq
\abs{\theta(v) - \theta(u)} \cdot \norm{v_1} \leq
\abs{\theta(v) - \theta(u)} \cdot \norm{v}
\\
=
&
\abs{\eta(d(v,S)) - \eta(d(u,S))} \cdot \norm{v} \leq
K \cdot \abs{d(v,S) - d(u,S)} \cdot \norm{v}
\\
\leq
&
K \cdot \norm{v - u} \cdot \norm{v} \leq
3Kr \cdot \norm{v - u}.
\vspace{1.7mm}
\end{array}
$
\renewcommand{\arraystretch}{1.0}
\addtolength{\arraycolsep}{0pt}
\newline
It follows that
$\norm{g(v) - g(u)} \leq (3Kr + 1) \cdot \norm{v - u}$.

Note that
$g\inverse(u) = \rfs{Rot}^{F,H}_{-\eta(d(u,S))}(u)$.
Since \num{iii} and \num{iv} of Part (c) hold for $-\eta$,
we also have that $g\inverse$ is $(3Kr + 1)$-Lipschitz.
So $M^{\srfs{rot}} = 3$.

(d) Let $E,F,H,\hatx,x,a$ be as in Part (d).
It suffices to prove Part (d) for $\hatx = 0$.
This is so, since if $g$ satisfies the requirements of Part (d) for
$E,F,0,x,a$, then $g^{\srfs{tr}_{\kern-0.5pt\hatx}}$ satisfies those
requirements for $E,F,\hatx,x,a$.
So $x' = x$ and $x'' = ax$.
Let $L = [x,ax]$ and $r = \dgfrac{\norm{x}}{2}$. So\break
$$
\frac{\rfs{lngth}(L)}{r} + 1 \leq
\frac{(a - 1) \norm{x}}{\dgfrac{\norm{x}}{2}} + 1 =
2 (a - 1) + 1 \leq 2 a.
$$
It follows from Part (a) that there is $g \in H(E)$ such that
$\rfs{supp}(g) \subseteq B(L,r)$,
$g(x) = a x$
and $g$ is $2 a M^{\srfs{seg}}$-bilipschitz.
A trivial computation shows that $g$ fulfills Requirements (d)(2)
and (d)(3).
\hfill\myqed

\begin{prop}\label{p-bddlip-1.13}
For every $K \geq 1$ there is $M^{\srfs{bnd}}(K) \geq 1$
   \index{N@mbnd@@$M^{\srfs{bnd}}(K)$}
such that the following holds.
{\thickmuskip=3.3mu \medmuskip=1.5mu \thinmuskip=1mu 
Suppose that $E$ is a normed space and $F$ is a closed linear
subspace of $E$. Let $x \in E - F$ be such that
$d(x,F) > \dgfrac{\norm{x}}{K}$ and
$y \in F - \sngltn{0}$.
Then there is $g \in H(E)$ and $a,b > 0$ such that\kern-5pt
}
\begin{itemize}
\addtolength{\parskip}{-11pt}
\addtolength{\itemsep}{06pt}
\item[\num{1}] $g(x) = ax + by$,
\item[\num{2}] $\norm{g(x)} = \norm{x}$,
\item[\num{3}] $d(g(x),F) = \dgfrac{\norm{g(x)}}{K}$,
\item[\num{4}] $g \nrestriction F = \rfs{Id}$,
\item[\num{5}]
$\rfs{supp}(g) \subseteq B(0;\frac{\norm{x}}{2},\frac{3 \norm{x}}{2})$,
\item[\num{6}] $g$ is $M^{\srfs{bnd}}(K)$-bilipschitz.
\vspace{-05.7pt}
\end{itemize}
\end{prop}

\noindent
{\bf Proof } Let $x,y$ be as in the proposition.
We may assume that $\norm{y} = \norm{x}$.
Let $L_1 = [x,y]$. We find $D(K)$
such that $d(L_1,0) \geq D(K) \norm{x}$.
Let
$E_1 = \rfs{span}(\dbltn{x}{y})$
and $F_1 = \rfs{span}(\sngltn{y})$.
So
$\norm{x} \leq K d(x,F_1)$.
Set
$M = M^{\srfs{thn}}(2)$,
and let
$\norm{\ }^{\srbfs{H}}$ be a Hilbert norm on $E_1$ such that
$\norm{u} \leq \norm{u}^{\srbfs{H}} \leq M \norm{u}$
for every $u \in E_1$,
Hence
$\frac{\norm{x}^{\ssbfs{H}}}{M} \leq K d^{\srbfs{H}}(x,E_1)$.
Also,
$\norm{x}^{\srbfs{H}},\norm{y}^{\srbfs{H}} \geq
\frac{\norm{x}}{M}$.
Let $\alpha$ be the angle
between $x$ and $F_1$.
Hence
$\sin(\alpha) = \frac{d^{\ssbfs{H}}(x,E_1)}{\norm{x}^{\ssbfs{H}}} \geq
\frac{1}{MK}$.
It follows that
$$
d(0,L_1) \geq
\frac{d^{\srbfs{H}}(0,L_1)}{M} =
\frac{\sin(\dgfrac{\alpha}{2}) \norm{x}^{\srbfs{H}}}{M} \geq
\frac{\sin(\alpha)}{2M} \norm{x} \geq \frac{\norm{x}}{2M^2K}.
$$
So $D(K) = \frac{1}{2M^2K}$.

\rule{0pt}{15pt}\kern-0.5pt
Since $\dgfrac{d(x,F)}{\norm{x}} > \dgfrac{1}{K}$ and 
$\dgfrac{d(y,F)}{\norm{y}} = 0 < \dgfrac{1}{K}$,
there is $z_0 \in [x,y]$ such that
$\dgfrac{d(z_0,F)}{\norm{z_0}} = \dgfrac{1}{K}$.
Obviously, $\norm{z_0} \leq \norm{x}$.
Let $z = \frac{\norm{x}}{\norm{z_0}} \ncdot z_0$.
So
$\norm{z} = \norm{x}$ and $\dgfrac{d(z,F)}{\norm{z}} = \dgfrac{1}{K}$.
Obviously, for some $a,b > 0$, $z = ax + by$.
Let $L = [x,z]$.
For some $\lambda \geq 1$, $z = \lambda z_0$.
This implies that for every $u \in L$
there are $v \in [x,z_0]$ and $\mu \geq 1$ 
such that $u = \mu v$.
It follows that $d(L,0) \geq d([x,z_0],0)$,
and since $[x,z_0] \subseteq L_1$, we have that
$d(L,0) \geq d(L_1,0) \geq \frac{\norm{x}}{2 M^2 K}$.

Obviously, $\norm{x - z} \leq 2 \norm{x}$.
Let $r = \frac{\norm{x}}{4 M^2 K}$.
By Proposition \ref{p-gamma.7}(a),
there is $h \in H(E)$ such that
$h(x) = z$,
$\rfs{supp}(h) \subseteq B(L,r)$ and
$h$ is
$M^{\srfs{seg}} \ncdot (\frac{\norm{x - z}}{r} + 1)$-bilipschitz.
By the above,
$$
\frac{\norm{x - z}}{r} + 1 \leq
\frac{2 \norm{x}}{\dgfrac{\norm{x}}{(4 M^2 K)}} + 1 = 8 M^2 K + 1 \leq
9 M^2 K.
$$
So $h$ is
$9 M^{\srfs{seg}} M^2 K$-bilipschitz.

Recall that $d(z,F) = \dgfrac{\norm{z}}{K}$,
$d(x,F) > \dgfrac{\norm{x}}{K}$ and for some $u \in F$ and $c > 0$,
$x = u + cz$.
This implies that $d(L,F) = \dgfrac{\norm{z}}{K}$.
Hence $d(B(L,r),F) = \frac{\norm{x}}{K} - r =
\frac{\norm{x}}{K} - \frac{\norm{x}}{4 M^2 K} > 0$.
So $h \nrestriction F = \rfs{Id}$.

From the fact that $\norm{z} = \norm{x}$,
it follows that $L \subseteq \overB(0,\norm{x})$.
So
$B(L,r) \subseteq B(0,\norm{x} + r)$.
Hence
$\rfs{supp}(h) \subseteq B(0,(1 + \frac{1}{4 M^2 K}) \norm{x})$.
But $1 + \frac{1}{4 M^2 K} < \dgfrac{3}{2}$,
so $\rfs{supp}(h) \subseteq B(0,\frac{3}{2} \norm{x})$.
Clearly,
$h \nrestriction B(0,d(L,0) - r) = \rfs{Id}$.
Hence $h \nrestriction B(0,\frac{\norm{x}}{4 M^2 K}) = \rfs{Id}$.

Let $\eta \in H([0,\infty))$ be the piecewise linear function
such that: (i)
the breakpoints of $\eta$ are
$\frac{\norm{x}}{4 M^2 K}$ and $\norm{x}$;
(ii) $\eta(\frac{\norm{x}}{4 M^2 K}) = \dgfrac{\norm{x}}{2}$,
and $\eta(t) = t$ for every $t \geq \norm{x}$.
The slopes of the pieces of $\eta$ are $2 M^2 K$,
$\frac{4 M^2 K}{2 (4 M^2 K - 1)}$
and $1$. So $\eta$ is $2 M^2 K$-bilipschitz.

Let $k$ be the radial homeomorphism based on $\eta$.
Then by Proposition \ref{metr-bldr-p3.18},
$k$ is\break
$6 M^2 K$-bilipschitz.
Also, $k(B(0,\frac{\norm{x}}{4 M^2 K})) = B(0,\frac{\norm{x}}{2})$,
$k(B(0,\frac{3 \norm{x}}{2})) = B(0,\frac{3 \norm{x}}{2})$,
$k(F) = F$, $k(x) = x$ and $k(z) = z$.

Let $g = h^k$.
Then $g(x) = z$,
$\rfs{supp}(g) \subseteq B(0;\frac{\norm{x}}{2},\frac{3 \norm{x}}{2})$,
$g \nrestriction F = \rfs{Id}$,
and $g$ has bilipschitz constant
$(6 M^2 K)^2 \ncdot 9 M^{\srfs{seg}} M^2 K$.
So $M^{\srfs{bnd}}(K) = 324 M^{\srfs{seg}} M^6 K^3$.
\rule{1pt}{0pt}\hfill\myqed


\begin{prop}\label{p-bddlip-1.16}
There is $M^{\srfs{cmp}} \geq 1$
   \index{N@mcmp@@$M^{\srfs{cmp}}$}
such that the following holds.
Suppose that $E = F \oplus H$, $\rfs{dim}(H) \leq 2$
and $H \perp^{M^{\srs{ort}}} F$.
Let $x \in E - F$, $x = \hatx + x^{\nperp}$, $\hatx \in F$,
$\norm{x^{\nperp}} \leq \frac{4}{3} d(x,F)$
and $d(x,F) \leq \frac{1}{16} \norm{x}$.
Then there is $g \in H(E)$ such that
\begin{itemize}
\addtolength{\parskip}{-11pt}
\addtolength{\itemsep}{06pt}
\item[\num{1}] $g$ is $M^{\srfs{cmp}}$-bilipschitz,
\item[\num{2}] $g(x) = \hatx + (x)_H$,
\item[\num{3}] $g \nrestriction F = \rfs{Id}$,
\item[\num{4}]
$\rfs{supp}(g) \subseteq B(0;\frac{\norm{x}}{2},\frac{3 \norm{x}}{2})$.
\vspace{-05.7pt}
\end{itemize}
\end{prop}

\noindent
{\bf Proof }
%
Note that $\hatx + x^{\nperp} = x = (x)_F + (x)_H$.
So $(x)_H - x^{\nperp} = \hatx - (x)_F \in F$.
Hence
$d(x + \lambda ((x)_H - x^{\nperp}),F) = d(x,F)$
for every $\lambda \in \bbR$.
Consider the interval $L = [x,\hatx + (x)_H]$.
Then $L = \setm{x + \lambda ((x)_H - x^{\nperp})}{\lambda \in [0,1]}$
and so $d(L,F) = d(x,F)$.
It follows that
{\thickmuskip=3mu \medmuskip=2mu \thinmuskip=1mu 
$$\hbox{$
\rfs{lngth}(L) = \norm{(x)_H - x^{\nperp}} \leq
\norm{(x)_H} + \norm{x^{\nperp}} \leq
M^{\srfs{ort}}\cdot d(x,F) + \fourthirds d(x,F) =
(M^{\srfs{ort}} + \fourthirds) d(x,F)
$}$$
}\kern-0pt
and hence
$
\dgfrac{\rfs{lngth}(L)}{d(x,F)} + 1 \leq
M^{\srfs{ort}} + \dgfrac{4}{3} + 1 \leq
M^{\srfs{ort}} + 3.
$
\rule{0pt}{14pt}\kern-2pt
We shall now find $\min_{u \in L} \norm{u}$ and
$\max_{u \in L} \norm{u}$.
Let $u \in L$. Then for some $\lambda \in [0,1]$,
$u = x + \lambda ((x)_H - x^{\nperp})$.
Recall that $d((x)_H,F) \geq \dgfrac{\norm{(x)_H}}{M^{\srs{ort}}}$.
So
\vspace{1.5mm}
\newline
\rule{7pt}{0pt}
\renewcommand{\arraystretch}{1.5}
\addtolength{\arraycolsep}{-4pt}
$
\begin{array}{ll}
&
\norm{u} \geq \norm{x} - \norm{x^{\nperp}} - \norm{(x)_H} \geq
\norm{x} - \fourthirds d(x,F) - M ^{\srfs{ort}} \cdot d((x)_H,F)
\\
=
\rule{5pt}{0pt}
&
\norm{x} - \fourthirds d(x,F) - M ^{\srfs{ort}} \cdot d(x,F) =
\norm{x} - (M ^{\srfs{ort}} + \fourthirds) d(x,F)
\\
\geq
\rule{5pt}{0pt}
&
\norm{x} - (M ^{\srfs{ort}} + \fourthirds) \frac{\norm{x}}{16}
\geq
\frac{9}{16} \norm{x}.
\vspace{1.7mm}
\end{array}
$
\renewcommand{\arraystretch}{1.0}
\addtolength{\arraycolsep}{4pt}
\newline
For the maximum of $\norm{u}$ we have
{\thickmuskip=8mu \medmuskip=4mu \thinmuskip=2mu 
$$\hbox{$
\norm{u} \leq \norm{x} + \norm{x^{\nperp}} + \norm{(x)_H} \leq
\norm{x} + \frac{\dgfrac{4}{3}}{16} \norm{x} + 
\frac{M ^{\srs{ort}}}{16} \norm{x} < 1 \frac{7}{16} \norm{x}.
$}$$
}\kern0pt
It follows that
$B(L,d(x,F)) \subseteq B(L,\frac{\norm{x}}{16}) \subseteq
B(0;\frac{\norm{x}}{2},\frac{3 \norm{x}}{2})$.
So by Proposition \ref{p-gamma.7}(a),
there is $g \in H(E)$ such that
$\rfs{supp}(g) \subseteq B(L,d(x,F))$,
$g(x) = \hatx + (x)_H$
and $g$ is $M^{\srfs{seg}} \ncdot (M^{\srfs{ort}} + 3)$-bilipschitz.
It follows that $g$ satisfies requirements (3)\,-\,(4)
of the proposition. So we may define
$M^{\srfs{cmp}} = M^{\srfs{seg}} (M^{\srfs{ort}} + 3)$.
\hfill\myqed

\begin{defn}\label{d-bddly-lip-2.9-1}
\begin{rm}
(a) Let $\alpha \in \rfs{MBC}$ and $s,t \in [0,\infty)$.
Then $s \approx^{\alpha} t$ means that
$t \leq \alpha(s)$ and $s \leq \alpha(t)$.

   \index{N@AAAA@@$s \approx^{\alpha} t$.
          This means $t \leq \alpha(s)$ and $s \leq \alpha(t)$}

(b) Let $\alpha \in \rfs{MBC}$, $n \in \bbN$
and $\fnn{\rho}{[0,\infty)}{[0,\infty)}$ be continuous.
We say that $\rho$ is $(n,\alpha)$-continuous,
if there are $0 = a_0 < \ldots < a_{n - 1} < a_n = \infty$ such that
$$
\rho_i(t) \eqdf \rho(t + a_{i - 1}),\ \,t \in [0,a_i - a_{i - 1}),
$$
is $\alpha$-continuous for every $0 < i \leq n$.

   \index{continuous. $\rho$ is $(n,\alpha)$-continuous}

\end{rm}
\end{defn}

The four parts of the next proposition are trivial.
Their proofs are omitted.

\begin{prop}\label{p-bddlip-1.11-1}
\num{a}
Let $\alpha \in \rfs{MBC}$, $n \in \bbN$
and $\fnn{\rho}{[0,\infty)}{[0,\infty)}$.
If $\rho$ is $(n,\alpha)$-continuous,
then $\rho$ is $n \ncdot \alpha$-continuous.

\num{b} Let $\fnn{\rho}{[0,\infty)}{[0,\infty)}$
and $a > 0$. Define $\eta(s,t)$ as follows.
If $s \geq a$, then $\eta(s,t) = t$;
and if $s \in [0,a]$, then
$\eta(s,t) = (1 - \frac{s}{a}) \rho(t) + \frac{s}{a} t$.
Suppose that $\beta \in \rfs{MC}$ and $\rho$ is $\beta$-continuous.
Then $\eta_s(t) \eqdf \eta(s,t)$ is $\beta$-continuous
for every $s \in [0,\infty)$.
We denote $\eta(s,t)$ by $\eta_{(\rho,a)}(s,t)$.
\kern-8pt
   \index{N@eta@@$\eta_{(\rho,a)}(s,t)$}

\num{c} Let $\beta \in \rfs{MC}$, $a > 0$ and
$0 < m \leq \frac{\beta(a)}{a}$.
Then the function $f(t) = mt$, $t \in [0,a]$, is $\beta$-continuous.

\num{d} If $\beta \in \rfs{MC}$,
$M \geq 1$,
and $\gamma$ is the function defined by $\gamma(t) = \beta(M t)$,
then $\gamma \leq M \beta$.
\end{prop}

\subsection{The main construction.}
\label{ss9.2}

\begin{defn}\label{d-bddly-lip-2.9}
\begin{rm}
(a) Let $0 < a < 1$ and $b,M > 1$.
We say that
{\it $M$ is a Uniform Continuity constant for
$\pair{a}{b}$},
({\it $M$ is UC-constant} for $\pair{a}{b}$)
   \index{uc-constant@@UC-constant.
          $M$ is a UC-constant for $\pair{a}{b}$}
if the following holds.

Suppose that $E,F,\alpha,x,y$ satisfy the following assumptions.
\begin{itemize}
\addtolength{\parskip}{-11pt}
\addtolength{\itemsep}{06pt}
\item[A1] $E$ is a normed space, $F$ is a closed linear proper
subspace of $E$, $\rfs{dim}(F) > 1$, $\alpha \in \rfs{MBC}$
and $x,y \in E - F$,
\item[A2]
$\norm{x} \leq \norm{y}$ and $\norm{x} \approx^{\alpha} \norm{y}$,
\item[A3] $d(x,F) \approx^{\alpha} d(y,F)$,
\item[A4] if $\rfs{co-dim}^E(F) = 1$,
then $x,y$ are on the same side of $F$.
\vspace{-05.7pt}
\end{itemize}
Then there are $\bfig_1,\bfig_2 \in H(E)$ such that
\begin{itemize}
\addtolength{\parskip}{-11pt}
\addtolength{\itemsep}{06pt}
\item[B1] $\bfig_1,\bfig_2$ are $M \alpha$-bicontinuous,
\item[B2] $\bfig_2 \scirc \bfig_1(x) = y$,
\item[B3] $\bfig_1(F) = F$ and $\bfig_2(F) = F$,
\item[B4] for every $i = 1,2$, \,
$\rfs{supp}(\bfig_i) \subseteq B(0;a \norm{x},b \norm{y})$.
\vspace{-05.7pt}
\end{itemize}

(b) We define a relation $R(u,v,g;\alpha,a,b,F)$.
Let $F$ be a closed linear subspace of a normed space $E$,
$u,v \in E - F$, $g \in H(E)$, $0 < a < 1$, $b > 1$
and $\alpha \in \rfs{MBC}$.
The notation $R(u,v,g;\alpha,a,b,F)$ means that
\begin{itemize}
\addtolength{\parskip}{-11pt}
\addtolength{\itemsep}{06pt}
\item[R1] $g(u) = v$,
\item[R2] $g$ is $\alpha$-bicontinuous,
\item[R3] $g(F) = F$,
\item[R4] $g \nrestriction B(0;a \norm{u},b \norm{v}) = \rfs{Id}$.
\vspace{-05.7pt}
\end{itemize}
Let $M \geq 1$. Then $R(u,v,g;M,a,b,F)$ means that
$R(u,v,g;M \ncdot \rfs{Id}_{[0,\infty)},a,b,F)$ holds.

   \index{N@r00@@$R(u,v,g;\alpha,a,b,F)$}
   \index{N@r01@@$R(u,v,g;M,a,b,F)$}
\end{rm}
\end{defn}

The trivial proof of Part (b) in the next proposition is omitted.

\begin{prop}\label{editor-p.12}
\num{a}
$$
\left(
R(u,v,g;\alpha,a,b,F) \wedge R(v,w,h;\beta,c,d,F)
\rule{0pt}{12pt}
\right)
\Rightarrow
R(u,w,h \scirc g;\beta \scirc \alpha,ac,bd,F).
$$

\num{b}
$R(u,v,g;M,a,b,F) \Rightarrow R(v,u,g\inverse;M,\dgfrac{a}{M},Mb,F)$.
\end{prop}

\noindent
{\bf Proof }
(a)
It is obvious that $h \scirc g$ is $\beta \scirc \alpha$-bicontinuous,
$h \scirc g(u) = w$ and
$h \scirc g(F) = F$.

If $v = u$, then $ca \norm{u} <  c \norm{u} = c \norm{v}$.
So $h \nrestriction B(0,ca \norm{u}) = \rfs{Id}$.
If $v \neq u$, then $v \in \rfs{supp}(g)$.
This implies that $\norm{v} > a \norm{u}$
and hence $c \norm{v} > ca \norm{u}$.
So $h \nrestriction B(0,ca \norm{u}) = \rfs{Id}$.
Clearly, $ca \norm{u} < a \norm{u}$.
So $g \nrestriction B(0,ca \norm{u}) = \rfs{Id}$.
It follows that $h \scirc g \nrestriction B(0,ac \norm{u}) = \rfs{Id}$.

If $v = w$, then $bd \norm{w} = bd \norm{v} > b \norm{v}$.
So $\rfs{supp}(g) \subseteq B(0,bd \norm{w})$.
If $v \neq w$, then $v \in \rfs{supp}(g) \subseteq B(0,d \norm{w})$.
This implies that $\norm{v} < d \norm{w}$
and hence $b \norm{v} < bd \norm{w}$.
So $\rfs{supp}(g) \subseteq B(0,b \norm{v}) \subseteq B(0,bd \norm{w})$.
It follows that $\rfs{supp}(g) \subseteq B(0,bd \norm{w})$.
From the fact that $bd > d$ if follows that
$\rfs{supp}(h) \subseteq B(0,bd \norm{w})$.
So $\rfs{supp}(h \scirc g) \subseteq B(0,bd \norm{w})$.
We have shown that
$\rfs{supp}(h \scirc g \subseteq B(0;ac \norm{u},bd \norm{w})$.
So $R(u,w,h \scirc g;\beta \scirc \alpha,ac,bd,F)$ holds.
\rule{100pt}{0pt}\hfill\myqed

\begin{lemma}\label{l-bddlip-1.10} \ %
The Uniform Continuity Constant Lemma.

\num{a} There are  $0 < a < 1$, $b > 1$ and $\bfiM > 1$ such that
$\bfiM$ is a UC-constant for $\pair{a}{b}$.

\num{b} For every  $0 < a < 1$, $b > 1$ there is $\bfiM > 1$
such that $\bfiM$ is a UC-constant for $\pair{a}{b}$.
\end{lemma}

\noindent
{\bf Proof }
(a) The proof is long and has many steps.
The survey below may help guide the reader through the proof.
\medskip

\noindent
{\bf Plan of the proof}

Let $E$, $F$, $\alpha$, $\bfix_0,\bfiy_0$
satisfy conditions A1\,-\,A4 in the definition of a UC-constant.
We construct two bilipschitz homeomomorphisms
$\bfie$ and $\bfih$.
Set $\bfie(\bfix_0) = \bfix$ and $\bfiy = \bfih\inverse(\bfiy_0)$.
Next we construct $N \ncdot \alpha$-bicontinuous homeomorphisms
$\bfif_1$, $\bfif_2$ and $\bfiv \in E$
such that $\bfif_1(\bfix) = \bfiv$ and $\bfif_2(\bfiv) = \bfiy$.
Here $N$ is a fixed number independent of
$E$, $F$, $\alpha$, $\bfix_0$ and $,\bfiy_0$.
So we have
$$\hbox{
$\bfie(\bfix_0) = \bfix$, $\bfif_1(\bfix) = \bfiv$,
$\bfif_2(\bfiv) = \bfiy$ and $\bfih(\bfiy) = \bfiy_0$.
}$$
The homeomorphisms $\bfig_1 \eqdf \bfif_1 \scirc \bfie$
and $\bfig_2 \eqdf \bfih \scirc \bfif_2$ are the ones required
in the definition of a UC-constant.
To explain what each homeomorphism does,
we take the simpler situation in which $E$ is a pre-Hilbert space.
Let $E$ be a pre-Hilbert space
and $F$ be a closed linear subspace of $E$.
For any $z \in E$, denote $(z)_{F,F^{\perp}}$ by $\hatz$ and
$(z)_{F^{\perp},F}$ by $z^{\nperp}$.
The homeomorphism $\bfie$ is a composition of four actions.
So $\bfie = \bfie_4 \scirc \ldots \scirc \bfie_1$.
Similarly, $\bfih$ is a composition of two actions.
We shall define homeomorphisms $\bfih_1$ and $\bfih_2$,
and $\bfih$ will be the homeomorphism
$\bfih_1\inverse \scirc \bfih_2\inverse$.

The first action $\bfie_1$ is needed only if
$d(\bfix_0,F) > \dgfrac{\norm{\bfix_0}}{3}$.
Otherwise, $\bfie_1 = \rfs{Id}$.
If the former happen, then $\bfie_1(\bfix_0) = \bfix_1$,
where $d(\bfix_1,F) = \dgfrac{\norm{\bfix_1}}{3}$
and $\norm{\bfix_1} = \norm{\bfix_0}$.
A similar action is performed by a homeomorphism
$\bfih_1$ on $\bfiy_0$, and we denote 
$\bfih_1(\bfiy_0)$ by $\bfiy_1$.
We now have the points $\bfix_1$ and $\bfiy_1$
with the properties
$\norm{\bfix_1} = \norm{\bfix_0}$,
$d(\bfix_1,F) \leq \dgfrac{\norm{\bfix_1}}{3}$,
$\norm{\bfiy_1} = \norm{\bfiy_0}$
and $d(\bfiy_1,F) \leq \dgfrac{\norm{\bfiy_1}}{3}$.

Now, $\bfie_2$ takes $\bfix_1$ to
$\lambda \hatbfiy_1 + \bfix_1^{\nperp}$,
where $\lambda > 0$ and $\norm{\lambda \hatbfiy_1} = \norm{\hatbfix_1}$.
The action of $\bfie_2$ can be roughly described as a rotation
in the plane $F_1$ generated by $\hatbfix_1$ and $\hatbfiy_1$
and the identity on $F_1^{\perp}$.
It is at this stage that we need $F$ to be of dimension $\geq 2$.
Denote $\bfie_2(\bfix_1)$ by $\bfix_2$.

The homeomorphism $\bfie_3$ takes $\bfix_2$ to a vector $\bfix_3$
of the form $a \hatbfix_2 + b \bfix_2^{\nperp}$,
where $a,b > 0$, $\norm{\bfix_3} = \norm{\bfix_2}$,
$d(\bfix_3,F) \leq \dgfrac{\norm{\bfix_3}}{\itDelta}$,
and $\itDelta$ is a fixed number $> 1$ independent of
$E$, $F$, $\alpha$, $\bfix_0$ and $\bfiy_0$.
Similarly, $\bfih_2$ takes $\bfiy_1$ to a vector $\bfiy$
of the form $c \hatbfiy_1 + d \bfiy_1^{\nperp}$,
where $c,d > 0$, $\norm{\bfiy} = \norm{\bfiy_1}$
and $d(\bfiy,F) \leq \dgfrac{\norm{\bfiy}}{\itDelta}$.
Denote $\bfiy$ by $\bfiy_2$.

Note that the subspace
$K \eqdf \rfs{span}(\bfix_3^{\nperp},\bfiy^{\nperp})$
is orthogonal to $F$.
(This is not true when $E$ is a general normed space).
Set
$\bfix^{\svee} =
\frac{\norm{\sbfix_{\kern0.7pt3}^{\snperp}}}{\norm{\sbfiy^{\snperp}}}
\bfiy^{\nperp}$
and define $\bfix = \hatbfix_3 + \bfix^{\svee}$.
Clearly, $\bfix^{\nperp} = \bfix^{\svee}$.
The homeomorphism $\bfie_4$ takes $\bfix_3$ to $\bfix$.
The action of $\bfie_4$ can be roughly described as a rotation
in the plane $\hatbfix_3 + K$ and the identity on $K^{\perp}$.
Define $\bfix_4 = \bfix$.

We have the following situation:
$\bfix = \hatbfix + \bfix^{\nperp}$,
$\bfiy = \hat\bfiy + \bfiy^{\nperp}$,
$\hatbfix, \hat\bfiy \in F$ and for some $\lambda,\mu > 0$,
$\hat\bfiy = \lambda \hatbfix$
and $\bfiy^{\nperp} = \mu \bfix^{\nperp}$.
If $\norm{\hatbfiy} \geq \norm{\hatbfix}$
define $\bfiv = \hatbfiy + \bfix^{\nperp}$,
and if $\norm{\hatbfiy} < \norm{\hatbfix}$
define $\bfiv = \lambda \bfix$.
We shall define $\bfif_1$ such that
$\bfif_1(\bfix) = \bfiv$ and $\bfif_1$ is $N \ncdot \alpha$-bicontinuous
for some fixed $N$. 
If $\bfiv = \hatbfiy + \bfix^{\nperp}$, then $\bfif_1$ has the form
$\bfif_1(z) = z + a(z) \ncdot \hatbfix$,
and $a(z)$ tends to zero as $d(z,[\bfix,\bfiv\kern1pt])$ tends to
$\lambda$.
In the case that $\bfiv = \lambda \bfix$,
$\bfif_1$ is a piecewise linearly radial homeomorphism
and $\bfif_1$ is $N$-bilipschitz. This of course implies that
$\bfif_1$ is $N \ncdot \alpha$-bicontinuous.

No we have $\bfiv = \hatbfiy + \bfiv^{\nperp}$
and $\bfiy = \hatbfiy + \bfiy^{\nperp}$,
where for some $\nu > 0$, $\bfiy^{\nperp} = \nu \bfiv^{\nperp}$.
We shall define $\bfif_2$ which takes $\bfiv$ to $\bfiy$.
The homeomorphism $\bfif_2$ will have the form
$\bfif_2(z) = z + a(z) \ncdot \bfiv^{\nperp}$,
and it will be $N \ncdot \alpha$-bicontinuous.
\medskip

Along the construction described above,
but independenly of the particular choice of
$E$, $F$, $\alpha$, $\bfix_0,\bfiy_0$,
we shall define numbers
\begin{itemize}
\addtolength{\parskip}{-11pt}
\addtolength{\itemsep}{06pt}
\item[] 
$\bfiM_{1,i},a_{1,i},b_{1,i}$,\, for $i = 1,\ldots,4$;
\item[] 
$\bfiM_{2,i},a_{2,i},b_{2,i}$,\, for $i = 1,2$;
\item[] 
$\bfiM_{3,i},a_{3,i},b_{3,i}$,\, for $i = 1,2$.
\vspace{-05.7pt}
\end{itemize}
These numbers satisfy the following conditions.
\begin{itemize}
\addtolength{\parskip}{-11pt}
\addtolength{\itemsep}{06pt}
\item[C1]
for every $i = 1,\ldots,4$,\,
$R(\bfix_{i - 1},\bfix_i,\bfie_i;\bfiM_{1,i},a_{1,i},b_{1,i},F)$;
\item[C2]
for every $i = 1,2$,\,
$R(\bfiy_{i - 1},\bfiy_i,\bfih_i;\bfiM_{2,i},a_{2,i},b_{2,i},F)$;
\item[C3]
$R(\bfix_4,\bfiv,\bfif_1;\bfiM_{3,1} \ncdot \alpha,a_{3,1},b_{3,1},F)$;
\item[C4]
$R(\bfiv,\bfiy_2,\bfif_2;\bfiM_{3,2} \ncdot \alpha,a_{3,2},b_{3,2},F)$.
\vspace{-05.7pt}
\end{itemize}

We have thus the following conclusion.
There are $\bfiM_{i,j},a_{i,j}$, $b_{i,j}$
such that for every
$E$, $F$, $\alpha$, $\bfix_0,\bfiy_0$
satisfying conditions A1\,-\,A4 in the definition of a UC-constant,
there are
$\bfie_i \in H(E)$, $\bfix_i$, \,$i = 1,\ldots,4$; \ 
$\bfih_i \in H(E)$, $\bfiy_i$, \,$i = 1,2$;
$\bfif_1,\bfif_2 \in H(E)$ and $\bfiv$
such that C1\,-\,C4 hold.

We now find $a,b,\bfiM$ such that
$\bfiM$ is a UC-constant for $\pair{a}{b}$.
Let $E$, $F$, $\alpha$, $\bfix_0,\bfiy_0$
fulfill conditions A1\,-\,A4 in the definition of a UC-constant.
Then there are $\bfie_i$'s, $\bfif_i$'s, $\bfih_i$, etc.
which satisfy C1\,-\,C4.
Define $\bfie = \bfie_4 \scirc \ldots \scirc \bfie_1$,
$\bfih = \bfih_1\inverse \scirc \bfih_2\inverse$,
$\bfig_1 = \bfif_1 \scirc \bfie$
and $\bfig_2 = \bfih \scirc \bfif_2$.

\rule{0pt}{16pt}Let
$\bfiM_1 = \prod_{i = 1}^4 \bfiM_{1,i}$,
$\bfiA_1 = \prod_{i = 1}^4 a_{1,i}$ and
$\bfiB_1 = \prod_{i = 1}^4 b_{1,i}$.
Then by Proposition~\ref{editor-p.12}(a),
$R(\bfix_0,\bfix_4,\bfie;\bfiM_1,\bfiA_1,\bfiB_1,F)$ holds.
By \ref{editor-p.12}(b),
$R(\bfiy_1,\bfiy_0,\bfih_1\inverse;
\bfiM_{2,1},\frac{a_{2,1}}{\ssbfiM_{2,1}},\bfiM_{2,1} b_{2,1},F)$
and\break
$R(\bfiy_2,\bfiy_1,\bfih_2\inverse;\bfiM_{2,2},
\frac{a_{2,2}}{\ssbfiM_{2,2}},\bfiM_{2,2} b_{2,2},F)$
hold.
Let
$\bfiA_2 =\frac{a_{2,2}}{\ssbfiM_{2,2}} \frac{a_{2,1}}{\ssbfiM_{2,1}}$,
$\bfiB_2 = \bfiM_{2,2} b_{2,2} \bfiM_{2,1} b_{2,1}$
and
$\bfiM_2 = \bfiM_{2,2} \bfiM_{2,1}$.
Then by Proposition~\ref{editor-p.12}(a),
$R(\bfiy_2,\bfiy_0,\bfih;\bfiM_2,\bfiA_2,\bfiB_2,F)$ holds.
Let
$\bfiM\fprime = \bfiM_1 \bfiM_{3,1}$,
$\bfiA' = \bfiA_1 a_{3,1}$ and
$\bfiB\fprime = \bfiB_1 b_{3,1}$.
Note that if $\alpha \in \rfs{MC}$ and $M \geq 1$,
then $\alpha(M t) \leq M \alpha$.
So by Proposition~\ref{editor-p.12}(a),
\newline
(1)
\centerline{\kern-10mm
$R(\bfix_0,\bfiv,\bfig_1;
\bfiM\fprime \ncdot \alpha,\bfiA',\bfiB\fprime,F)$ holds.
}
Let
$\bfiM\fprime' = \bfiM_{3,2} \bfiM_2$,
$\bfiA'' = a_{3,2} \bfiA_2$ and
$\bfiB\fprime' = b_{3,2} \bfiB_2$.
Then by Proposition~\ref{editor-p.12}(a),
\newline
(2)
\centerline{\kern-9mm
$R(\bfiv,\bfiy_0,\bfig_2;
\bfiM\fprime' \ncdot \alpha,\bfiA'',\bfiB\fprime',F)$ holds.
}
Let $\bfiM = \max(\bfiM\fprime,\bfiM\fprime')$,
$a = \bfiA' \bfiA''$ and
$b = \bfiB\fprime \bfiB\fprime'$.
Then (1) and (2) imply that $B1\,-\,B4$
of Definition \ref{d-bddly-lip-2.9}(a) hold.
So $\bfiM$ is a UC-constant for $\pair{a}{b}$.
\smallskip

C1 is the conjunction of four requirements. Denote them by
C1.1,$\ldots$,C1.4.
Similarly, denote the two conjuncts of C2 by C2.1  and C2.2.
\medskip

\noindent
{\bf The construction}

\noindent
{\bf Part 1 } The construction of $\bfie_1$ and $\bfih_1$.
\\
Let $E, F, \alpha, \bfix_0, \bfiy_0$ satisfy conditions A1\,-\,A4
in the definition of a UC-constant.
Write $x = \bfix_0$ and $y = \bfiy_0$.

If $d(x,F) \leq \dgfrac{\norm{x}}{3}$,
let $\bfie_1 = \rfs{Id}$.
Otherwise let $u \in F - \sngltn{0}$ and $\bfie_1 \in H(E)$ be such that
\begin{itemize}
\addtolength{\parskip}{-11pt}
\addtolength{\itemsep}{06pt}
\item[\num{1}] $\bfie_1(x) \in \rfs{span}(\dbltn{x}{u})$,
\item[\num{2}] $\norm{\bfie_1(x)} = \norm{x}$,
\item[\num{3}] $d(\bfie_1(x),F) = \dgfrac{\norm{\bfie_1(x)}}{3}$,
\item[\num{4}] $\bfie_1 \nrestriction F = \rfs{Id}$; 
\item[\num{5}]
$\rfs{supp}(\bfie_1) \subseteq
B(0;\frac{\norm{x}}{2},\frac{3\norm{x}}{2})$,
\item[\num{6}] $\bfie_1$ is $M^{\srfs{bnd}}(3)$ bilipschitz.
\vspace{-05.7pt}
\end{itemize}
The existence of $\bfie_1$ is assured
by Proposition \ref{p-bddlip-1.13}.
Let $\bfix_1 = f(x)$, $\bfiM_{1,1} = M^{\srfs{bnd}}(3)$,
$a_{1,1} = \dghalf$ and $b_{1,1} = \dgfrac{3}{2}$.
Recall that $\bfix_0 = x$.
By (1)\,-\,(6),
$R(\bfix_0,\bfix_1,\bfie_1;\bfiM_{1,1},a_{1,1},b_{1,1},F)$
holds. So C1.1 is fulfilled.

Let $\bfih_1 \in H(E)$ have the same role for $y$ as $\bfie_1$ had
for $x$.
Let $\bfiy_1 = \bfih_1(y)$, $\bfiM_{2,1} = M^{\srfs{bnd}}(3)$,
$a_{2,1} = \dghalf$ and $b_{2,1} = \dgfrac{3}{2}$.
Recall that $\bfiy_0 = y$.
Then $R(\bfiy_0,\bfiy_1,\bfih_1;\bfiM_{2,1},a_{2,1},b_{2,1},F)$ holds.
So C2.1 is fulfilled.
\medskip

\noindent
{\bf Part 2 } The construction of $\bfie_2$.
\\
Since $\norm{\bfie_1(x)} = \norm{x}$ and $\norm{\bfih_1(y)} = \norm{y}$,
$\norm{\bfie_1(x)} \approx^{\alpha} \norm{\bfih_1(y)}$.
\vspace{1mm}
We check that
$$
d(\bfie_1(x),F) \approx^{\alpha} d(\bfih_1(y),F).
$$
If $\bfie_1 = \bfih_1 = \rfs{Id}$, then there is nothing to check.
Suppose that $\bfie_1 \neq \rfs{Id} \neq \bfih_1$.
Then
$d(\bfie_1(x),F) = \dgfrac{\norm{x}}{3}$ and
$d(\bfih_1(y),F) = \dgfrac{\norm{y}}{3}$.
So
$$\hbox{$
\frac{d(\ssbfih_1(y),F)}{d(\ssbfie_1(x),F)} =
\frac{\norm{y}}{\norm{x}} \leq
\frac{\alpha(\norm{x})}{\norm{x}} \leq
\frac{\alpha(\dgfrac{\norm{x}}{3})}{\dgfrac{\norm{x}}{3}} =
\frac{\alpha(d(\ssbfie_1(x),F))}{d(\ssbfie_1(x),F)}.
$}$$
Hence $d(\bfih_1(y),F) \leq \alpha(d(\bfie_1(x),F))$.
Since $\norm{x} \leq \norm{y}$,
$d(\bfih_1(x),F) \leq d(\bfie_1(y),F) \leq \alpha(d(\bfie_1(y),F))$.

Suppose that $\bfie_1 \neq \rfs{Id} = \bfih_1$.
Then $d(\bfih_1(y),F) \leq \dgfrac{\norm{y}}{3}$ and
\hbox{$d(\bfie_1(x),F) = \dgfrac{\norm{x}}{3}$}.
So
$d(\bfih_1(y),F) \leq \dgfrac{\alpha(\norm{x})}{3} \leq
\alpha(\dgfrac{\norm{x}}{3}) =
\alpha(d(\bfie_1(x),F))$.
Also,
$d(\bfie_1(x),F) \leq d(x,F) \leq \alpha(d(y,F)) =
\alpha(d(\bfih_1(y),F))$.
The argument in the case $\bfie_1 = \rfs{Id} \neq \bfih_1$ is identical.

Let $\bfie_1(x)$ take the role of $x$ and $\bfih_1(y)$
take the role of $y$. 
That is, $\bfie_1(x),\bfih_1(y)$
are renamed and are now denoted by $x$ and $y$.
Hence
$d(x,F) \leq \dgfrac{\norm{x}}{3}$ and
$d(y,F) \leq \dgfrac{\norm{y}}{3}$.
Let $\hatx,\haty \in F$ be such that
$\norm{x - \hatx} \leq (1 + \varepsilon)d(x,F)$
and $\norm{y - \haty} \leq (1 + \varepsilon)d(y,F)$.\break
$\varepsilon$ will be determined later.
Let $x^{\nperp} = x - \hatx$ and
$y^{\nperp} = y - \haty$.
$\bfie_2$ will take $x$ to a vector of the form
$\lambda \haty + x^{\nperp}$, where $\lambda > 0$.
It is in this part that $F$ needs to be of dimension $> 1$.
We may assume that:
\begin{itemize}
\addtolength{\parskip}{-11pt}
\addtolength{\itemsep}{06pt}
\item[{\smaller 2.1 }]
$x = \hatx + x^{\nperp}$ and $y = \haty + y^{\nperp}$;
\item[{\smaller 2.2 }]
$\hatx, \haty \in F$;
\item[{\smaller 2.3 }]
$\norm{x^{\nperp}} \leq (1 + \varepsilon)d(x,F)$ and
and $\norm{y^{\nperp}} \leq (1 + \varepsilon)d(y,F)$.
\item[{\smaller 2.4 }]
$d(x,F) \leq \dgfrac{\norm{x}}{3}$ and 
$d(y,F) \leq \dgfrac{\norm{y}}{3}$; 
\item[{\smaller 2.5 }]
$\norm{x} \approx^{\alpha} \norm{y}$ and
$d(x,F) \approx^{\alpha} d(y,F)$.
\item[{\smaller 2.6 }]
If $\rfs{co-dim}^E(F) = 1$,
then $x$ and $y$ are on the same side of $F$.
\vspace{-05.7pt}
\end{itemize}

{\thickmuskip=2.4mu \medmuskip=1.3mu \thinmuskip = 1.0mu
We define a functional $\psi$
on $\rfs{span}(F \cup \sngltn{x^{\nperp}})$:
$\psi(x^{\nperp}) = \norm{x^{\nperp}}$,
and $\psi(u) = 0$ for every $u \in F$.
Let $\varphi \in E^*$ be such that $\varphi$ extends $\psi$ and
$\norm{\varphi} = \norm{\psi}$.
Let $L = \rfs{span}(\sngltn{x^{\nperp}})$
and $H = \rfs{ker}(\varphi)$.
}
So $F \subseteq H$.
For every $u \in F$,
$$
\abs{\psi(u + x^{\nperp})} = \norm{x^{\nperp}} \leq
(1 + \varepsilon)d(x,F) = (1 + \varepsilon)d(x^{\nperp},F) \leq
(1 + \varepsilon) \norm{u + x^{\nperp}}.
$$
So $\norm{\varphi} = \norm{\psi} \leq 1 + \varepsilon$.

Let $u \in E$.
Define
$v = u - \varphi(u) \frac{x^{\snperp}}{\norm{x^{\snperp}}}$.
Then $(u)_H = v$
and $(u)_L = \varphi(u) \frac{x^{\snperp}}{\norm{x^{\snperp}}}$.
So
\vspace{1.5mm}
\newline
\rule{7pt}{0pt}
\renewcommand{\arraystretch}{1.5}
\addtolength{\arraycolsep}{-0pt}
$
\begin{array}{lll}
\norm{(u)_H}
&
=
&
\norm{v} =
\norm{u - \varphi(u)\frac{x^{\snperp}}{\norm{x^{\snperp}}}} \leq
\norm{u} +
\abs{\varphi(u)} \norm{\frac{x^{\snperp}}{\norm{x^{\snperp}}}}
=
\norm{u} + \abs{\varphi(u)}
\\
&
\leq
&
\norm{u} + \norm{\varphi}\norm{u} \leq (2 + \varepsilon) \norm{u}
\vspace{1.7mm}
\end{array}
$
\renewcommand{\arraystretch}{1.0}
\addtolength{\arraycolsep}{0pt}
\newline
and
$\norm{(u)_L} =
\norm{\varphi(u) \frac{x^{\snperp}}{\norm{x^{\snperp}}}} =
\abs{\varphi(u)} \leq \norm{\varphi} \norm{u} \leq
(1 + \varepsilon) \norm{u}$. So
$$
\norm{(u)_H} + \norm{(u)_L} \leq (3 + 2 \varepsilon) \norm{u}.
$$

Let $F_1$ be a $2$-dimensional subspace of $F$ such that
$\hatx,\haty \in F_1$. Such a subspace exists since $F$ is not
$1$-dimensional.
Let $H_1$ be an almost orthogonal complement of $F_1$ in $H$.
That is, $H_1 \oplus F_1 = H$, and for every $u \in H$,
$\norm{(u)_{F_1}} + \norm{(u)_{H_1}} \leq
M^{\srfs{aoc}}(2) \cdot \norm{u}$.
Let\break
$\norm{\ }^{\srbfs{H}}$ be a tight Hilbert norm on $F_1$.
So $\norm{\ }^{\srbfs{H}} \approx^{M^{\srs{thn}}} \norm{\ }^{F_1}$.

We define an equivalent norm $\norm{\ }^{\srbfs{N}}$ on $E$.
Let $u \in E$ and suppose that $u = u_1 + u_2 + u_3$,
where $u_1 \in F_1$, $u_2 \in H_1$
and $u_3 \in L$.
Define
$\norm{u}^{\srbfs{N}} \eqdf
\norm{u_1}^{\srbfs{H}} + \norm{u_2} + \norm{u_3}$.
Then
$\norm{u} \approx^{3 + 2 \varepsilon} \norm{u_1 + u_2} + \norm{u_3}$
and
$\norm{u_1 + u_2} \approx^{M^{\srs{hlb}}}
\norm{u_1}^{\srbfs{H}} + \norm{u_2}$.
Note that if $E = E_1 \oplus E_2$, for $\ell,i = 1,2$,
$\norm{\ }^{\ell,i}$ is a norm on $E_{\ell}$
and $\norm{\ }^{\ell,1} \approx^{M_{\ell}} \norm{\ }^{\ell,2}$,
then for every $u \in E$, \,
$\norm{(u)_{E_1}}^{1,1} + \norm{(u)_{E_2}}^{2,1} \approx^{\max(M_1,M_2)}
\norm{(u)_{E_1}}^{1,2} + \norm{(u)_{E_2}}^{2,2}$.
So
$
\norm{u_1 + u_2} + \norm{u_3}
\approx^{\max(M^{\srs{hlb}},1)}
\norm{u_1}^{\srbfs{H}} + \norm{u_2} + \norm{u_3}$.
That is,
$ \norm{u_1 + u_2} + \norm{u_3}
\kern3pt\approx^{M^{\srs{hlb}}}\kern3pt
\norm{u_1}^{\srbfs{H}} + \norm{u_2} + \norm{u_3}$.
Let $M^{\srfs{sp}} = (3 + 2 \varepsilon) M^{\srfs{hlb}}$.
Then
$\norm{u} \approx^{M^{\srs{sp}}}
\norm{u_1}^{\srbfs{H}} + \norm{u_2} + \norm{u_3} =
\norm{u}^{\srbfs{N}}$.
Let $d^{\srbfs{N}}$ denote the metric on $E$ obtained from
$\norm{\ }^{\srbfs{N}}$.

Let
$\hatz =
\frac{\norm{\hatx}^{\srs{H}}}{\norm{\haty}^{\srs{H}}} \, \haty$.
Then $\norm{\hatz}^{\srbfs{H}} = \norm{\hatx}^{\srbfs{H}}$.
The homeomorphism $\bfie_2$ will take $x$ to $\hatz + x^{\nperp}$.
Let $r = \norm{\hatx}^{\srbfs{H}}$,
$S_1 = S^{\srbfs{H}}(0,r)$ and
$S = \setm{u + \mu \cdot (x)^{\nperp}}{u \in S_1 \mbox{ and }
0 \leq \mu \leq 1}$.
Let $\theta_0$ be the angle from $\hatx$ to $\haty$.
That is,
$\rfs{Rot}^{F_1}_{\theta_0}(\hatx) = \hatz$.

Let $E_1 = H_1 + L$. Then $F_1 \oplus E_1 = E$.
We first define a function $\fnn{\eta}{[0,\infty)}{[0,\theta_0]}$,
and the homeomorphism $\bfie_2$ will be defined by means of $\eta$
as follows:
$$
\bfie_2(u) =
\rfs{Rot}^{F_1}_{\eta(d^{\srs{N}}(u,S))}((u)_{F_1}) + (u)_{E_1}.
$$
Define $\eta$ to be the piecewise linear function with one breakpoint
at $\dgfrac{r}{2}$, such that $\eta(0) = \theta_0$
and $\eta(s) = 0$ for every $s \geq \dgfrac{r}{2}$.

Note that $\hatx \in F_1$, $x^{\nperp} \in L$ and
$x = \hatx + x^{\nperp}$. So $(x)_{F_1} = \hatx$ and
$(x)_{E_1} = x^{\nperp}$.
Also, $x \in S$.
It follows that $\bfie_2(x) = \hatz + x^{\nperp}$.
Hence for some $\lambda > 0$,
$\bfie_2(x) = \lambda \haty + x^{\nperp}$.
Obviously, $\bfie_2(F_1) = F_1$. We verify that
\begin{equation}\tag{2.1}
\bfie_2(F) = F.
\end{equation}
Suppose that $u \in F$.
So $u = (u)_{F_1} + (u)_{E_1}$.
Hence $(u)_{E_1} \in F$.
For some angle $\beta$,
$\bfie_2(u) = \rfs{Rot}^{F_1}_{\beta}((u)_{F_1}) + (u)_{E_1}$.
Since $F_1 \subseteq F$,
$\rfs{Rot}^{F_1}_{\beta}((u)_{F_1}) \in F$.
So $\bfie_2(u) \in F$.

Note that
$d^{\srbfs{N}}(B^{\srbfs{N}}(u,s),S) = \dgfrac{r}{2}$.
Hence $\bfie_2 \nrestriction B^{\srbfs{N}}(0,\dgfrac{r}{2}) = \rfs{Id}$.
By {\smaller 2.3} and {\smaller 2.4},
\vspace{1.5mm}
\newline
\rule{7pt}{0pt}
\renewcommand{\arraystretch}{1.5}
\addtolength{\arraycolsep}{-0pt}
$
\begin{array}{lll}
r
&
=
&
\norm{\hatx}^{\srbfs{N}} \geq
\frac{\norm{\hatx}}{M^{\srs{sp}}} \geq
\frac{1}{M^{\srs{sp}}}
(\norm{x} - \norm{x^{\nperp}})
\geq
\frac{1}{M^{\srs{sp}}}
(\norm{x} - (1 + \varepsilon)d(x,F))
\\
&
\geq
&
\frac{1}{M^{\srs{sp}}}
(\norm{x} - (1 + \varepsilon) \frac{\norm{x}}{3})) =
\frac{1}{M^{\srs{sp}}} (\twothirds - \varepsilon) \norm{x} >
\frac{1}{2 M^{\srs{sp}}} \norm{x}
\vspace{1.7mm}
\end{array}
$
\renewcommand{\arraystretch}{1.0}
\addtolength{\arraycolsep}{0pt}
\newline
The last inequality holds when $\varepsilon$ is sufficiently small.
So
$\bfie_2 \nrestriction
B^{\srbfs{N}}(0,\frac{1}{4 M^{\srs{sp}}} \norm{x}) =
\rfs{Id}$.

Recall that $\norm{\ }^E \approx^{M^{\srs{sp}}} \norm{\ }^{\srbfs{N}}$.
So $B(0,\dgfrac{s}{M^{\srs{sp}}}) \subseteq B^{\srbfs{N}}(0,s)$
for every $s$.
It follows that 
$\bfie_2 \nrestriction B(0,\frac{1}{4 (M^{\srs{sp}})^2} \norm{x}) =
\rfs{Id}$.
Let $a_1 = \frac{1}{4 (M^{\srs{sp}})^2}$. We have shown that
\begin{equation}
\tag{2.2}
\bfie_2 \nrestriction B(0,a_1 \norm{x}) = \rfs{Id}.
\end{equation}

Now, $\rfs{supp}(\bfie_2) \subseteq
B^{\srbfs{N}}(0,\norm{x}^{\srbfs{N}} + \dgfrac{r}{2}) \subseteq
B(0,M^{\srfs{sp}}(\norm{x}^{\srbfs{N}} + \dgfrac{r}{2}))$
and
$\dgfrac{r}{2} = \dgfrac{\norm{\hatx}^{\srbfs{N}}}{2} \leq
\dgfrac{M^{\srfs{sp}} \norm{\hatx}}{2} \leq
\dgfrac{M^{\srfs{sp}} \ncdot \frac{4}{3} \norm{x}}{2} =
\frac{2}{3} M^{\srfs{sp}} \norm{x}$.
So 
$\rfs{supp}(\bfie_2) \subseteq
B(0, 2 (M^{\srfs{sp}})^2 \norm{x})$.
Define $b_1 = 2 (M^{\srfs{sp}})^2$. Then
\begin{equation}\tag{2.3}
\bfie_2 \nrestriction (E - B(0,b_1 \norm{x})) = \rfs{Id}.
\end{equation}

We next show that there is $M_1 > 0$ which is independent
of $x$, $F$ and $\theta_0$ such that $\bfie_2$ is $M_1$-bilipschitz.
Indeed, we shall find $M_1'$ such that for every $u,v \in E$:
if $\norm{u - v}^{\srbfs{N}} \leq \dgfrac{r}{2}$, then
$\norm{\bfie_2(u) - \bfie_2(v)}^{\srbfs{N}} \leq
M_1' \cdot \norm{u - v}^{\srbfs{N}}$.
This fact implies that $\bfie_2$ is $M_1'$-Lipschitz
in the metric~$d^{\srbfs{N}}$.

Obviously,
$\abs{\eta(t) - \eta(s)} \leq
\frac {\theta_0}{\dgfrac{r}{2}} \abs{t - s} \leq
\frac{2 \pi}{r} \abs{t - s}$ for every $s,t \in [0,\infty)$.
Define
$\theta(u) = \eta(d^{\srs{N}}(u,S))$.
So
$\abs{\theta(u) - \theta(v)} =
\abs{\eta(d^{\srs{N}}(u,S)) - \eta(d^{\srs{N}}(v,S))} \leq
\frac{2 \pi}{r} \norm{u - v}^{\srbfs{N}}$.

Clearly, $\norm{x^{\nperp}} < \dgfrac{\norm{x}}{2}$.
So $\norm{x} < 2 \norm{\hatx}$.
Hence $\norm{x^{\nperp}} < \norm{\hatx}$.
It follows that
$\norm{x^{\nperp}}^{\srbfs{N}} <
(M^{\srfs{sp}})^2 \norm{\hatx}^{\srbfs{N}}$.
Hence
$\max(\setm{\norm{u}^{\srbfs{N}}}{u \in S}) \leq
(1 + (M^{\srfs{sp}})^2) \cdot \norm{\hatx}^{\srbfs{N}} =
2 (1 + (M^{\srfs{sp}})^2) \cdot r$.

Let $u,v \in E$ be such that
$\norm{u - v}^{\srbfs{N}} \leq \dgfrac{r}{2}$.
If $\norm{u}^{\srbfs{N}} > 2 (1 + (M^{\srfs{sp}})^2) \cdot r + r$,
then
$\norm{v}^{\srbfs{N}} >
2 (1 + (M^{\srfs{sp}})^2) \cdot r + \dgfrac{r}{2}$.
So $\bfie_2(u) = u$ and $\bfie_2(v) = v$.
Suppose that
$\norm{u}^{\srbfs{N}} \leq 2 (1 + (M^{\srfs{sp}})^2) \cdot r + r$.
Denote
$M^{\srfs{sp1}} = 4 + 2 (M^{\srfs{sp}})^2$.
Then
$\norm{u}^{\srbfs{N}}, \norm{v}^{\srbfs{N}} < M^{\srfs{sp1}} \cdot r$.
\vspace{1.5mm}
\newline
\rule{7pt}{0pt}
\renewcommand{\arraystretch}{1.5}
\addtolength{\arraycolsep}{-0pt}
$
\begin{array}{lll}
\bfie_2(v) - \bfie_2(u)
&
=
&
(\rfs{Rot}^{F_1}_{\theta(v)}((v)_{F_1}) -
\rfs{Rot}^{F_1}_{\theta(u)}((v)_{F_1})) \kern4pt+
\\
\rule{5pt}{0pt}
&&
(\rfs{Rot}^{F_1}_{\theta(u)}((v)_{F_1}) -
\rfs{Rot}^{F_1}_{\theta(u)}((u)_{F_1})) \kern2pt + \kern2pt
((v)_{E_1} - (u)_{E_1}).
\vspace{1.7mm}
\end{array}
$
\renewcommand{\arraystretch}{1.0}
\addtolength{\arraycolsep}{0pt}
\newline
So
\newline
\rule{0pt}{0pt}\kern-6pt
\renewcommand{\arraystretch}{1.5}
\addtolength{\arraycolsep}{3pt}
$
\begin{array}{lll}
\rule{0pt}{16pt}
\norm{\bfie_2(v) - \bfie_2(u)}^{\srbfs{N}} \leq \kern-4pt
&
\norm{\rfs{Rot}^{F_1}_{\theta(v)}((v)_{F_1}) -
\rfs{Rot}^{F_1}_{\theta(u)}((v)_{F_1})}^{\srbfs{N}} +
&
\\&
\rule{0pt}{16pt}
\norm{(\rfs{Rot}^{F_1}_{\theta(u)}((v)_{F_1}) -
\rfs{Rot}^{F_1}_{\theta(u)}((u)_{F_1})) \kern4pt +
&
\rule{0pt}{0pt}\kern-5pt
((v)_{E_1} - (u)_{E_1})}^{\srbfs{N}} =
\\&
\rule{0pt}{16pt}
\norm{\rfs{Rot}^{F_1}_{\theta(v)}((v)_{F_1}) -
\rfs{Rot}^{F_1}_{\theta(u)}((v)_{F_1})}^{\srbfs{N}} +
&
\rule{0pt}{0pt}\kern-5pt
\norm{v - u}^{\srbfs{N}}.
\vspace{1.7mm}
\end{array}
$
\renewcommand{\arraystretch}{1.0}
\addtolength{\arraycolsep}{-3pt}
\newline
We deal with the first summand in the last expression.
\vspace{1.5mm}
\newline
\rule{7pt}{0pt}
\renewcommand{\arraystretch}{1.5}
\addtolength{\arraycolsep}{-3pt}
$
\begin{array}{ll}
&
\norm{\rfs{Rot}^{F_1}_{\theta(v)}((v)_{F_1}) -
\rfs{Rot}^{F_1}_{\theta(u)}((v)_{F_1})}^{\srbfs{N}} \leq
\abs{\theta(v) - \theta(u)} \cdot \norm{v}^{\srbfs{N}} \leq
\frac{2 \pi}{r} \norm{v - u}^{\srbfs{N}} \cdot \norm{v}^{\srbfs{N}}
\\
\leq
\rule{5pt}{0pt}
&
\frac{2 \pi}{r} \norm{v - u}^{\srbfs{N}} \cdot M^{\srfs{sp1}} \cdot r =
2 \pi \norm{v - u}^{\srbfs{N}} \cdot M^{\srfs{sp1}}.
\vspace{1.7mm}
\end{array}
$
\renewcommand{\arraystretch}{1.0}
\addtolength{\arraycolsep}{3pt}
\newline
It follows that for every $u,v \in E$,
$\norm{\bfie_2(v) - \bfie_2(u)}^{\srbfs{N}} \leq
(2 \pi M^{\srfs{sp1}} + 1) \cdot \norm{v - u}^{\srbfs{N}}$.

Obviously, for every $u \in E$,
$\bfie_2\inverse(u) = 
\rfs{Rot}^{F_1}_{-\eta(d^{\srs{N}}(u,S))}((u)_{F_1}) + (u)_{E_1}$.
So
$$
\norm{\bfie_2\inverse(v) - \bfie_2\inverse(u)}^{\srbfs{N}} \leq
(2 \pi M^{\srfs{sp1}} + 1) \cdot \norm{v - u}^{\srbfs{N}}.
$$
Let $M_1 = (2 \pi M^{\srfs{sp1}} + 1) \cdot (M^{\srfs{sp}})^2$.
Then
\begin{equation}\tag{2.4}
\hbox{$\bfie_2$ is $M_1$-bilipschitz in the norm $\norm{\ }^E$.}
\end{equation}

Set $\bfix_2 = \bfie_2(x)$ and recall that $\bfix_1 = x$.
Hence by (2.1)\,-\,(2.4),
$R(\bfix_1,\bfix_2,\bfie_2;M_1,a_1,b_1,F)$ holds.
That is, C1.2 is fulfilled with $\bfiM_{1,2} = M_1$, $a_{1,2} = a_1$
and $b_{1,2} = b_1$.

Since $\bfie_2$ is $M_1$-bilipschitz and $\bfie_2(0) = 0$,
it follows that $\norm{\bfie_2(x)} \approx^{M_1} \norm{x}$.
From the fact that $\bfie_2(F) = F$, it follows that
$d(\bfie_2(x),F) \approx^{M_1} d(x,F)$.
So
\begin{equation}
\tag{2.5}
\norm{\bfie_2(x)} \approx^{M_1 \cdot \alpha} \norm{y}
\mbox{ \,and\, }
d(\bfie_2(x),F) \approx^{M_1 \cdot \alpha} d(y,F).
\end{equation}

\noindent
{\bf Part 3 }
The construction of $\bfie_3$, $\bfih_2$ and $\bfie_4$.
\newline
Recall that $\bfix_2$ has the form $\lambda \haty + x^{\nperp}$.
Rename $\bfix_2$ and call it $x$,
and denote $\lambda \haty$ by $\hatx$.
We now have that
\begin{list}{}
{\setlength{\leftmargin}{31pt}
\setlength{\labelsep}{05pt}
\setlength{\labelwidth}{20pt}
\setlength{\itemindent}{-00pt}
\addtolength{\topsep}{-04pt}
\addtolength{\parskip}{-02pt}
\addtolength{\itemsep}{-05pt}
}
\item[{\smaller 3.1$^{\kern1pt*}$ }]
$x = \hatx + x^{\nperp}$ and $y = \haty + y^{\nperp}$,
\item[{\smaller 3.2$^{\kern1pt*}$ }]
$\hatx, \haty \in F$ and for some $\lambda > 0$,
$\hatx = \lambda \haty$,
\item[{\smaller 3.3$^{\kern1pt*}$ }]
$\norm{x^{\nperp}} \leq (1 + \varepsilon)d(x,F)$
and $\norm{y^{\nperp}} \leq (1 + \varepsilon)d(y,F)$,
\item[{\smaller 3.5$^{\kern1pt*}$ }]
$\norm{x} \approx^{M_1 \cdot \alpha} \norm{y}$ and
$d(x,F) \approx^{M_1 \cdot \alpha} d(y,F)$,
\item[{\smaller 3.6$^{\kern1pt*}$ }]
If $\rfs{co-dim}^E(F) = 1$,
then $x$ and $y$ are on the same side of $F$.
\vspace{-02.0pt}
\end{list}
Property {\smaller 3.4} which is analogous to {\smaller 2.4} is missing.
Only after applying $\bfie_3$ to $x$ and $\bfih_2$ to $y$,
we shall retain this property.

For the next step in the construction we choose some $\itDelta > 1$.
The value of $\itDelta$ will be determined later,
and it will be independent of $E,F,\alpha,\bfix_0$ and $\bfiy_0$.
The definition of $\bfie_3$ and $\bfih_2$ depends on $\itDelta$.

We first define $\bfie_3$.
If
$d(x,F) \leq \dgfrac{\norm{x}}{\itDelta}$, then define
$\bfie_3 = \rfs{Id}$.
Suppose that $d(x,F) > \dgfrac{\norm{x}}{\itDelta}$.
Then there are $\bfie_3 \in H(E)$ and $a,b > 0$ such that
\begin{itemize}
\addtolength{\parskip}{-11pt}
\addtolength{\itemsep}{06pt}
\item[\num{1}]
$\bfie_3(x) = a \hatx + bx$,
\item[\num{2}]
$\norm{\bfie_3(x)} = \norm{x}$,
\item[\num{3}]
$d(\bfie_3(x),F) = \dgfrac{\norm{\bfie_3(x)}}{\itDelta}$,
\item[\num{4}]
$\bfie_3 \nrestriction F = \rfs{Id}$,
\item[\num{5}]
$\rfs{supp}(\bfie_3) \subseteq
B(0;\frac{\norm{x}}{2},\frac{3\norm{x}}{2})$,
\item[\num{6}]
$\bfie_3$ is $M^{\srfs{bnd}}(\itDelta)$-bilipschitz.
\vspace{-05.7pt}
\end{itemize}
The existence of $\bfie_3$ follows from Proposition \ref{p-bddlip-1.13}.

Recall that $\bfix_2 = x$ and denote $\bfix_3 = \bfie_3(x)$.
Then
$R(\bfix_2,\bfix_3,\bfie_3;
M^{\srfs{bnd}}(\itDelta),\dghalf,\dgfrac{3}{2},F)$
holds.
That is, C1.3 is fulfilled with
$\bfiM_{1,3} = M^{\srfs{bnd}}(\itDelta)$,
$a_{1,3} = \dghalf$ and
$b_{1,3} = \dgthreehalves$.

There is $\bfih_2 \in H(E)$ which acts on $y$ in the way that $\bfie_3$
acts on $x$. That is,\break
if $d(y,F) \leq \dgfrac{\norm{y}}{\itDelta}$,
then $\bfih_2 = \rfs{Id}$,
and if $d(y,F) > \dgfrac{\norm{y}}{\itDelta}$,
then there are $c,d > 0$ such that (1)\,-\,(6) above hold
when $y,\bfih_2,c,d$ replace $x,\bfie_3,a,b$.
Recall that $\bfiy_1 = y$ and denote
$\bfiy_2 = \bfih_2(y)$. Then
$R(\bfiy_1,\bfiy_2,\bfih_2;
M^{\srfs{bnd}}(\itDelta),\dghalf,\dgfrac{3}{2},F)$
holds.
That is, C2.2 is fulfilled with
$\bfiM_{2,2} = M^{\srfs{bnd}}(\itDelta)$,
$a_{2,2} = \dghalf$ and $b_{2,2} = \dgthreehalves$.

Suppose that $\bfie_3 \neq \rfs{Id}$.
Then $(\star)$
$\bfie_3(x) = a \lambda \haty + b(\lambda \haty + x^{\nperp}) =
(a + b) \lambda \haty + b x^{\nperp}$.\break
By {\smaller 3.1$^{\kern1pt*}$\,-\,3.3$^{\kern1pt*}$},
$\norm{x^{\nperp}} \leq (1 + \varepsilon) d(x^{\nperp},F)$.
So from $(\star)$ it follows that
$\norm{b x^{\nperp}} \leq (1 + \varepsilon)d(\bfie_3(x),F)$.
Denote $(a + b) \lambda \haty$ by $\hatbfix_3$
and $b x^{\nperp}$ by $\bfix_3^{\nperp}$.
In
{\smaller 3.1}$^{\kern1pt*}$\,-\,{\smaller 3.3}$^{\kern1pt*}$
and in {\smaller 3.6}$^{\kern1pt*}$
replace $x$, $\hatx$ and\break
$x^{\nperp}$ by
$\bfix_3$, $\hatbfix_3$ and $\bfix_3^{\nperp}$,
and denote the resulting statements by
{\smaller 3.1}$^{\kern1pt*}(\bfix_3,y)$ etc..
Then\break
{\smaller 3.1}$^{\kern1pt*}
(\bfix_3,y)$\,-\,{\smaller 3.3}$^{\kern1pt*}(\bfix_3,y)$
and {\smaller 3.6}$^{\kern1pt*}(\bfix_3,y)$ hold.
Also,
\begin{equation}\tag{$\dagger$}
d(\bfix_3,F) \leq \dgfrac{\norm{\bfix_3}}{\itDelta}.
\end{equation}
If $\bfie_3 = \rfs{Id}$ and we define $\hatbfix_3$ to be $\hatx$
and $\bfix_3^{\nperp}$ to be $x^{\nperp}$, then again ($\dagger$) holds.

Applying the same argument to $\bfiy_2$ and defining
$\hatbfiy_2$ and $\bfiy_2^{\nperp}$ in analogy with
$\hatbfix_3$ and $\bfix_3^{\nperp}$ we conclude that
{\smaller 3.1}$^{\kern1pt*}
(x,\bfiy_2)$\,-\,{\smaller 3.3}$^{\kern1pt*}(x,\bfiy_2)$
and {\smaller 3.6}$^{\kern1pt*}(x,\bfiy_2)$ hold.
Also, ($\dagger$) holds for $\bfiy_2$.

From {\smaller 3.5$^{\kern1pt*}$} and from (6) applied to
$\bfie_3$ and $\bfih_2$ it follows that
$$\norm{\bfix_3} \approx^{M^{\srs{bnd}}(\itDelta)} \norm{x}
\approx^{M_1 \cdot \alpha} \norm{y}
\approx^{M^{\srs{bnd}}(\itDelta)} \norm{\bfiy_2}
$$
and hence ($\dagger$$\dagger$)
$\norm{\bfix_3} \approx^{M_1 (M^{\srs{bnd}}(\itDelta))^2 \cdot \alpha}
\norm{y}$.
Similarly,
($\dagger$$\dagger$$\dagger$)
$d(\bfix_3,F) \approx^{M_1 (M^{\srs{bnd}}(\itDelta))^2 \cdot \alpha}
d(\bfiy_2,F)$.

We now rename
$\bfix_3,\hatbfix_3,\bfix_3^{\nperp},
\bfiy_2,\hatbfiy_2,\bfiy_2^{\nperp}$
and denote them by $x,\hatx,x^{\nperp},y,\haty$ and $y^{\nperp}$.
We also denote $M_1 (M^{\srs{bnd}}(\itDelta))^2 \cdot \alpha$
by $\alpha_1$.
From the above we conclude that
\begin{itemize}
\addtolength{\parskip}{-11pt}
\addtolength{\itemsep}{06pt}
\item[{\smaller 3.1 }]
$x = \hatx + x^{\nperp}$ and $y = \haty + y^{\nperp}$,
\item[{\smaller 3.2 }]
$\hatx, \haty \in F$ and for some $\lambda > 0$,
$\hatx = \lambda \haty$,
\item[{\smaller 3.3 }]
$\norm{x^{\nperp}} \leq (1 + \varepsilon)d(x,F)$
and $\norm{y^{\nperp}} \leq (1 + \varepsilon)d(y,F)$,
\item[{\smaller 3.4 }]
$d(x,F) \leq \dgfrac{\norm{x}}{\itDelta}$
\kern4pt and \kern4pt
$d(y,F) \leq \dgfrac{\norm{y}}{\itDelta}$,
\item[{\smaller 3.5 }]
$\norm{x} \approx^{\alpha_1}
\norm{y}$
\kern4pt and \kern4pt
$d(x,F) \approx^{\alpha_1} d(y,F)$,
\item[{\smaller 3.6 }]
If $\rfs{co-dim}^E(F) = 1$,
then $x$ and $y$ are on the same side of $F$.
\vspace{-05.7pt}
\end{itemize}
Property {\smaller 3.1\,} follows from
{\smaller 3.1}$^{\kern1pt*}(\bfix_3,y)$
and {\smaller 3.1}$^{\kern1pt*}(x,\bfiy_2)$,
and the same is true for Properties~{\smaller 3.2}, {\smaller 3.3\,}
and {\smaller 3.6}.
Property {\smaller 3.4\,} is the conjunction of
($\dagger$) applied to $\bfix_3$ and to $\bfiy_2$
and {\smaller 3.6\,} is the conjunction of
($\dagger$$\dagger$) and ($\dagger$$\dagger$$\dagger$).

Set
$z^{\nperp} =
\norm{x^{\nperp}} \cdot \frac{y^{\snperp}}{\norm{y^{\snperp}}}$
and $z = \hatx + z^{\nperp}$.
We next define $\bfie_4$.
It will take $x$ to $z$.
So after applying $\bfie_4$ we shall reach the following situation:
$\bfix_4 = \hatbfix_4 + \bfix_4^{\nperp}$,
$\bfiy_2 = \hatbfiy_2 + \bfiy_2^{\nperp}$,
$\hatbfix_4 = \lambda \hatbfiy_2$ for some $\lambda > 0$
and $\bfix_4^{\nperp} = \mu \bfiy_2^{\nperp}$ for some $\mu > 0$.

There are two cases:
$\rfs{co-dim}^E(F) = 1$ and $\rfs{co-dim}^E(F) > 1$.

{\bf Case 1 } $\rfs{co-dim}^E(F) = 1$.
Since $x$ and $y$ are on the same side of $F$,
there are $\nu > 0$ and $u \in F$ such that
$z^{\nperp} = u + \nu x^{\nperp}$.
Let $L = [x,\hatx + z^{\nperp}]$.
We may assume that in {\smaller 3.3\kern1pt},
$\varepsilon \leq \dghalf$.
We show that $\dgfrac{\rfs{lngth}(L)}{d(L,F)} + 1 \leq 19$.
Clearly,
$\rfs{lngth}(L) = \norm{\hatx + z^{\nperp} - x} =
\norm{z^{\nperp} - x^{\nperp}} \leq 2 \norm{x^{\nperp}}$.
So
\begin{equation}
\tag{3.1}
\rfs{lngth}(L) \leq 2 \norm{x^{\nperp}}.
\end{equation}
Since for some $t$, $z^{\nperp} = t y^{\nperp}$,
we have that $\norm{z^{\nperp}} \leq (1 + \varepsilon) d(z^{\nperp},F)$.
So
$$
\norm{x^{\nperp}} = \norm{z^{\nperp}} \leq
(1 + \varepsilon) d(u + \nu x^{\nperp},F) =
(1 + \varepsilon) \nu d(x^{\nperp},F) \leq
(1 + \varepsilon) \nu \norm{x^{\nperp}}.
$$
Hence
$1 \leq (1 + \varepsilon) \nu$.
In the above argument we interchange the roles of 
$x^{\nperp}$ and $z^{\nperp}$.
That is, for some $u' \in F$,
$x^{\nperp} = u' + \frac{1}{\nu} z^{\nperp}$,
and hence $1 \leq (1 + \varepsilon) \frac{1}{\nu}$.
We conclude that
{\thickmuskip=3mu \medmuskip=2mu \thinmuskip=1mu 
$\frac{1}{1 + \varepsilon} \leq \nu \leq 1 + \varepsilon$.
Let $v \in L$. Then for some $t \in [0,1]$,
$v = \hatx + x^{\nperp} + t (z^{\nperp} - x^{\nperp}) =
\hatx + x^{\nperp} + t ((u + \nu x^{\nperp} - x^{\nperp})$.
}
So
\vspace{1.5mm}
\newline
\rule{7pt}{0pt}
\renewcommand{\arraystretch}{1.5}
\addtolength{\arraycolsep}{-0pt}
$
\begin{array}{lll}
d(v,F)
&
=
&
d((1 + t (\nu - 1)) x^{\nperp},F) =
\abs{1 + t(\nu - 1)} \cdot d(x^{\nperp},F) \geq
(1 - t \abs{\nu - 1}) \cdot d(x^{\nperp},F)
\\
&
\geq
&
(1 - \abs{\nu - 1}) \cdot d(x^{\nperp},F) \geq
(1 - (1 + \varepsilon - \frac{1}{1 + \varepsilon})) \cdot
d(x^{\nperp},F)
\\
&
=
&
(\frac{1}{1 + \varepsilon} - \varepsilon)
d(x^{\nperp},F) \geq
\frac{1}{6} d(x^{\nperp},F) \geq
\frac{1}{6(1 + \varepsilon)} \norm{x^{\nperp}} \geq
\frac{1}{9} \norm{x^{\nperp}}.
\vspace{1.7mm}
\end{array}
$
\renewcommand{\arraystretch}{1.0}
\addtolength{\arraycolsep}{0pt}
\newline
Hence
\begin{equation}
\tag{3.2}
d(L,F) \geq \dgfrac{\norm{x^{\nperp}}}{\kern1pt9}.
\end{equation}
It follows from (3.1) and (3.2) that
$\dgfrac{\rfs{lngth}(L)}{d(L,F)} + 1 \leq 19$.

Set $\itDelta = 8$. Then
$d(L,F) \leq d(x,F) \leq \dgfrac{\norm{x}}{8}$.
Hence
$\norm{x^{\nperp}} \leq \frac{3}{2} d(x,F) \leq \frac{3}{16} \norm{x}$.
So
$\rfs{lngth}(L) \leq \frac{3}{8} \norm{x}$.
Let $B = B(L,d(L,F))$.
Then
\begin{equation}
\tag{3.3}
\min_{v \in B} \norm{v} \geq \norm{x} - \rfs{lngth}(L) - d(L,F) \geq
\dgfrac{\norm{x}}{2}.
\end{equation}

Similarly,
\begin{equation}
\tag{3.4}
\max_{v \in B} \norm{v} \leq \norm{x} + \rfs{lngth}(L) + d(L,F)
\leq \dgfrac{3\norm{x}}{\kern1pt2}.
\end{equation}
{\thickmuskip=2.0mu \medmuskip=1mu \thinmuskip=1mu 
The endpoints of $L$ are $x$ and $\hatx + z^{\nperp}$,
so by Proposition \ref{p-gamma.7}(a), there is $\bfie_4 \in H(E)$
such that
}
\begin{itemize}
\addtolength{\parskip}{-11pt}
\addtolength{\itemsep}{06pt}
\item[(3.5)] \ $\rfs{supp}(\bfie_4) \subseteq B(L,d(L,F))$,
\item[(3.6)] \ $\bfie_4(x) =
\hatx + z^{\nperp}$,
\item[(3.7)] \ $\bfie_4$ is $19 M^{\srfs{seg}}$-bilipschitz.
\vspace{-05.7pt}
\end{itemize}
By (3.5), $\bfie_4 \nrestriction F = \rfs{Id}$.
By (3.3), (3.4) and (3.5),
$\rfs{supp}(\bfie_4) \subseteq
B(0;\frac{\norm{x}}{2},\frac{3 \norm{x}}{2})$.
Recall that $\bfix_3 = x$ and denote $\bfix_4 = \bfie_4(x)$.
It follows that
$R(\bfix_3,\bfix_4,\bfie_4;19 M^{\srfs{seg}},\dghalf,\dgfrac{3}{2},F)$
holds.

{\bf Case 2 } $\rfs{co-dim}^E(F) > 1$.
Let $\itUpsilon > 1$.
By Proposition \ref{p-bddlip-1.11},
there is a closed subspace $F_1$ of $E$ such that $F \subseteq F_1$,
\hbox{$\rfs{span}(F_1 \cup \dbltn{x}{y}) = E$,}
$d(x,F_1) \geq \frac{1}{\itUpsilon} d(x,F)$
and $d(y,F_1) \geq \frac{1}{\itUpsilon} d(y,F)$.
Obviously, either
$\rfs{co-dim}^E(F_1) = 1$ or
\hbox{$\rfs{co-dim}^E(F_1) = 2$.}
If $\rfs{co-dim}^E(F_1) = 1$,
let $F \subseteq F_2 \subseteq F_1$ be a closed subspace such that
$\rfs{co-dim}^E(F_2) = 2$.
Otherwise let $F_2 = F_1$.
It follows that
$\rfs{co-dim}^E(F_2) = 2$,
$d(x,F_2) \geq \frac{1}{\itUpsilon} d(x,F)$
and
$d(y,F_2) \geq \frac{1}{\itUpsilon} d(y,F)$.

In {\smaller 3.4}, choose $\itDelta = 24$.
Hence $d(x,F_2) \leq d(x,F) \leq \dgfrac{\norm{x}}{24}$.
In {\smaller 3.3}, choose $\varepsilon = \dgninth$,
and choose $\itUpsilon = 1 \ninth$.
So
$\norm{x^{\nperp}} \leq (1 + \varepsilon)d(x,F) \leq
\itUpsilon (1 + \varepsilon)d(x,F_2) \leq
(1 \ninth)^2 d(x,F_2) \leq \fourthirds d(x,F_2)$.
In summary,
\begin{equation}
\tag{3.8}
\hbox{
$\norm{x^{\nperp}} \leq \dgfrac{4 d(x,F_2)}{3}$
\hspace{2mm} and \hspace{2mm}
$d(x,F_2) \leq \dgfrac{\norm{x}}{24}$.
}
\end{equation}

Recall that
$z^{\nperp} =
\frac{\norm{x^{\snperp}}}{\norm{y^{\snperp}}}\, y^{\nperp}$
and
$z = \hatx + z^{\nperp}$.
We have that
$\norm{y^{\nperp}} \leq \fourthirds d(y,F_2)$.
This is shown in the same way
that the analogous fact was proved for $x$.
Obviously, $d(y,F_2) = d(y^{\nperp},F_2)$.
So $\norm{y^{\nperp}} \leq \fourthirds d(y^{\nperp},F_2)$.
Since $z^{\nperp}$ is a multiple of $y^{\nperp}$,
$\norm{z^{\nperp}} \leq \fourthirds d(z^{\nperp},F_2)$.
Also, $d(z,F_2) = d(z^{\nperp},F_2)$.
So $\norm{z^{\nperp}} \leq \fourthirds d(z,F_2)$.

Note that $z = x - x^{\nperp} + z^{\nperp}$.
So $\norm{z} \geq \norm{x} - \norm{x^{\nperp}} - \norm{z^{\nperp}} =
\norm{x} - 2 \norm{x^{\nperp}}$.
Also, $\norm{x^{\nperp}} \leq \fourthirds d(x,F_2) \leq
\fourthirds \cdot \frac{1}{24} \norm{x} = \frac{1}{18} \norm{x}$.
Hence
$
\frac{d(z,F_2)}{\norm{z}} \leq
\frac{\norm{z^{\nperp}}}{\norm{z}} \leq
\frac{\norm{x^{\nperp}}}{\norm{x} - 2 \norm{x^{\nperp}}} \leq
\frac{\dgfrac{\norm{x}}{18}}{\norm{x} - \dgfrac{\norm{x}}{9}}
=
\frac{1}{16}.
$
In summary,
\begin{equation}
\tag{3.9}
\hbox{
$\norm{z^{\nperp}} \leq \dgfrac{4 d(z,F_2)}{3}$
\hspace{2mm} and \hspace{2mm}
$d(z,F_2) \leq \dgfrac{\norm{z}}{16}$.
}
\end{equation}

Let $H$ be such that $E = F_2 \oplus H$
and $H \perp^{M^{\srs{ort}}} F_2$.
We apply Proposition~\ref{p-bddlip-1.16} to $x$ and to $z$.
Note that by (3.8) and 3.9), $x$ and to $z$
satisfy the assumptions of~\ref{p-bddlip-1.16}.
So there is $f_1 \in H(E)$
such that:
$f_1$ is $M^{\srfs{cmp}}$-bilipschitz,
$f_1(x) = \hatx + (x)_H$,
$f_1 \nrestriction F_2 = \rfs{Id}$
and
$\rfs{supp}(f_1) \subseteq B(0;\frac{\norm{x}}{2},\frac{3\norm{x}}{2})$.
Similarly, there is $h_1 \in H(E)$ such that:
$h_1$ is $M^{\srfs{cmp}}$-bilipschitz,
$h_1(z) = \hatx + (z)_H$,
$h_1 \nrestriction F_2 = \rfs{Id}$
and
$\rfs{supp}(h_1) \subseteq
B(0;\frac{\norm{z}}{2},\frac{3\norm{z}}{2})$.

We now translate what we have obtained for $f_1$ and $h_1$ to
statements of the form $R(.,.,f_1;\ldots)$
and $R(.,.,h_1;\ldots)$.
Since $f_1$ is $M^{\srfs{cmp}}$-bilipschitz
$f_1(x) = \hatx + (x)_H$ and $f_1(0) = 0$, it follows that
$\norm{x} \leq M^{\srfs{cmp}} \norm{\hatx + (x)_H}$.
So $\rfs{supp}(f_1) \subseteq
B(0;\half \norm{x},\frac{3 M^{\trfs{cmp}}}{2} \norm{\hatx + (x)_H})$.
This implies that
\begin{equation}\tag{3.10}
R(x,\hatx + (x)_H,f_1;
M^{\srfs{cmp}},\hbox{$\half$},\hbox{$\frac{3 M^{\trfs{cmp}}}{2}$},F)
\mbox{ holds.}
\end{equation}
Similarly,
\begin{equation}\tag{3.11}
R(z,\hatx + (z)_H,h_1;
M^{\srfs{cmp}},\hbox{$\half$},\hbox{$\frac{3 M^{\trfs{cmp}}}{2}$},F)
\mbox{ holds.}
\end{equation}

Let $\norm{\ }^{\srbfs{H}}$ be a tight equivalent Hilbert norm on $H$,
and define a new norm on $E$ by
$\norm{u}^{\srbfs{N}} = \norm{(u)_{F_2}} + 
\norm{(u)_H}^{\srbfs{H}}$. 
So $\norm{\ } \approx^{M^{\ssrfs{fdn}}} \norm{\ }^{\srbfs{N}}$.
This follows from Proposition~\ref{p-bddlip-1.12}(c).
Let $d^{\srbfs{N}}$ denote the metric induced
by $\norm{\ }^{\srbfs{N}}$ on $E$.

Set $x^* = (x)_H$, $z^* = (z)_H$
and
$z^{\ssharp} =
\frac{\norm{x^*}^{\ssbfs{N}}}{\norm{z^*}^{\ssbfs{N}}} z^*$.
We define a homeomorphism $g_{2,1}$
which takes $\hatx + x^*$ to $\hatx + z^{\ssharp}$.
A second homeomorphism $g_{2,2}$,
will take $\hatx + z^{\ssharp}$ to $\hatx + z^*$.
So
$$
x = \hatx + x^{\nperp}
\kern1pt\stackrel{f_1}{\rightarrow}\kern1pt \hatx + (x)_H
\kern1pt\stackrel{g_{2,1}}{\rightarrow}\kern1pt \hatx + z^{\ssharp}
\kern1pt\stackrel{g_{2,2}}{\rightarrow}\kern1pt \hatx + (z)_H
\kern1pt\stackrel{h_1\inverse}{\rightarrow} \hatx + z^{\nperp} = z.
$$
Finally, we shall define
$\bfie_4 \eqdf h_1\inverse \scirc g_{2,2} \scirc g_{2,1} \scirc f_1$.

Let $\theta$ be the angle from $x^*$ to $z^{\ssharp}$. That is,
$\theta \in [0,\pi]$ and
$\rfs{Rot}_{\theta}^H(x^*) = z^{\ssharp}$.
Let $\fnn{\eta}{[0,\infty)}{[0,\theta]}$ be the piecewise linear
function with one breakpoint at
$s_0 = \frac{\norm{x^*}^{\ssbfs N}}{2 M^{\srs{thn}}}$ such that
$\eta(0) = \theta$ and $\eta(s) = 0$ for every $s \geq s_0$.
Let $S_0$ be the circle in $\pair{H}{\norm{\ }^{\srbfs{H}}}$ with
center at $0$ and radius $\norm{x^*}^{\srbfs{H}}$,
and let $S = \hatx + S_0$.
Let $g_{2,1}$ be defined as follows.
For $u \in E$ set $u_1 = (u)_H$ and $u_2 = (u)_{F_2}$. Define
$$
g_{2,1}(u) = u_2 + \rfs{Rot}^H_{\eta(d^{\ssbfs{N}}(u,S))}(u_1).
$$
Since for every $u \in E$,
$d^{\srbfs{N}}(u,S) = d^{\srbfs{N}}(g_{2,1}(u),S)$,
it follows that
$g_{2,1} \in H(E)$.
Clearly,
$$
g_{2,1}(\hatx + x^*) = \hatx + z^{\ssharp}.
$$
Also, $\rfs{supp}(g_{2,1}) \subseteq B^{\srbfs{N}}(S,s_0)$.
If $u \in F_2$
then
$d^{\srbfs{N}}(u,S) = \norm{u - \hatx} + \norm{x^*}^{\srbfs{H}} > s_0$
and so $g_{2,1}(u) = u$.
That is, $g_{2,1} \nrestriction F_2 = \rfs{Id}$.
Since $F \subseteq F_2$,
$$
g_{2,1} \nrestriction F = \rfs{Id}.
$$
Note that $s_0 = \frac{\norm{x^*}^{\ssbfs{N}}}{2 M^{\srs{thn}}} \leq
\frac{\norm{x^*}}{2}$. 
So
$\rfs{supp}(g_{2,1}) \subseteq
B^{\srbfs{N}}(S,\frac{\norm{x^*}}{2})$.

Let
$u \in
B^{\srbfs{N}}(0,\norm{\hatx} - \frac{\norm{x^*}}{2})$.
So $\norm{u_2} \leq \norm{\hatx} - \frac{\norm{x^*}}{2}$.
Then
\vspace{1.5mm}
\newline
\rule{7pt}{0pt}
\renewcommand{\arraystretch}{1.5}
\addtolength{\arraycolsep}{-0pt}
$
\begin{array}{lll}
d^{\srbfs{N}}(u,S)
&
=
&
\norm{u_2 - \hatx} + d^{\srbfs{N}}(u_1,S_0) \geq
\norm{u_2 - \hatx} \geq \norm{\hatx} - \norm{u_2}
\\
&
\geq
&
\norm{\hatx} - (\norm{\hatx} - \frac{\norm{x^*}}{2}) =
\frac{\norm{x^*}}{2}.
\vspace{1.7mm}
\end{array}
$
\renewcommand{\arraystretch}{1.0}
\addtolength{\arraycolsep}{0pt}
\newline
It follows that
$g_{2,1} \nrestriction
B^{\srbfs{N}}(0,\norm{\hatx} - \frac{\norm{x^*}}{2}) =
\rfs{Id}$.

Let $r = \norm{\hatx} + 2 \norm{x^*}^{\srbfs{N}}$.
Suppose that $u \in E - B^{\srbfs{N}}(0,r)$.
Either $\norm{u_1} \geq \frac{3 \norm{x^*}^{\ssbfs{N}}}{2}$
or $\norm{u_2} \geq \norm{\hatx} + \frac{\norm{x^*}^{\ssbfs{N}}}{2}$.
If $v \in S$ then $v = \hatx + w$, where $w \in H$ and
$\norm{w}^{\srbfs{N}} = \norm{x^*}^{\srbfs{N}}$.
Hence
$\norm{u - v}^{\srbfs{N}} =
\norm{u_1 - w}^{\srbfs{N}} + \norm{u_2 - \hatx}$.
If 
$\norm{u_1} \geq \frac{3 \norm{x^*}^{\ssbfs{N}}}{2}$,
then
$\norm{u - v}^{\srbfs{N}} \geq
\norm{u_1 - w}^{\srbfs{N}} \geq
\frac{3 \norm{x^*}^{\ssbfs{N}}}{2} - \norm{x^*}^{\srbfs{N}} =
\frac{\norm{x^*}^{\ssbfs{N}}}{2}$.
So $u \not\in \rfs{supp}(g_{2,1})$.
If $\norm{u_2} \geq \norm{\hatx} + \frac{\norm{x^*}^{\ssbfs{N}}}{2}$,
then
$\norm{u - v}^{\srbfs{N}} \geq
\norm{u_2 - \hatx} \geq
\norm{\hatx} + \frac{\norm{x^*}^{\ssbfs{N}}}{2} - \norm{\hatx} =
\frac{\norm{x^*}^{\ssbfs{N}}}{2}$.
So $u \not\in \rfs{supp}(g_{2,1})$.
It follows that $\rfs{supp}(g_{2,1}) \subseteq B^{\srbfs{N}}(0,r)$.

By (3.8),
$\norm{x^{\nperp}} \leq \frac{1}{18}\, \norm{x}$,
and since $x = \hatx + x^{\nperp}$,
we have
$\frac{17}{18}\, \norm{x} \leq \norm{\hatx} \leq
\frac{19}{18}\, \norm{x}$.
Since $H \perp^{M^{\trfs{ort}}} F_2$,
$\norm{x^*} \leq M^{\srfs{ort}} d(x^*,F_2)$.
Also, $M^{\srfs{ort}} < 4$.
By the above and (3.8),
$\norm{x^*} \leq M^{\srfs{ort}} d(x^*,F_2) =
M^{\srfs{ort}} d(x,F_2) \leq \frac{4}{24}\, \norm{x}$.
Hence
$\norm{\hatx} - \frac{\norm{x^*}}{2} \geq
\frac{17 - 2}{24} \norm{x}$
and\break
$r = \norm{\hatx} + 2 \norm{x^*}^{\srbfs{N}} \leq
(1 + \frac{M^{\trfs{thn}}}{3}) \norm{x}$.
It follows that
$$
\rfs{supp}(g_{2,1}) \subseteq
B(0;\hbox{$\frac{\norm{x}}{2}$},2 M^{\srfs{thn}} \norm{x}).
$$

Next we find a Lipschitz constant for $g_{2,1}$.
By its definition, $\eta$ is
$\frac{\theta}
{\dgfrac{\norm{x^*}^{\ssbfs N}}{(2 M^{\srs{thn}})}}$\,-\,Lipschitz.
So $\eta$ is
$\frac{2 \pi M^{\srs{thn}}}{\norm{x^*}^{\ssbfs N}}$\,-\,Lipschitz.
Obviously,
$S \subseteq \hatx + \barB^{\sbfs{N}}(0,\norm{x^*}^{\sbfs{N}})$.
By \ref{p-gamma.7}(c),
$g_{2,1}$ is
$
(M^{\srfs{rot}} \cdot
\frac{2 \pi M^{\srs{thn}}}{\norm{x^*}^{\ssbfs{N}}} \cdot
\norm{x^*}^{\sbfs{N}} + 1)
$\,-\,Lipschitz in the norm $\norm{\ }^{\sbfs{N}}$.
That is, $g_{2,1}$ is
$(2 \pi M^{\srfs{rot}} \cdot M^{\srfs{thn}} + 1)$-Lipschitz
in the norm $\norm{\ }^{\sbfs{N}}$.
The same is true for $g_{2,1}\inverse$.
So $g_{2,1}$ is
$(2 \pi M^{\srfs{rot}} \cdot M^{\srfs{thn}} + 1)$-bilipschitz
in the norm $\norm{\ }^{\sbfs{N}}$.
Recall that
$\norm{\ } \approx^{M^{\ssrfs{fdn}}} \norm{\ }^{\srbfs{N}}$.
Write
$\whatM_{2,1} =
(M^{\ssrfs{fdn}})^2 (2 \pi M^{\srfs{rot}} \cdot M^{\srfs{thn}} + 1)$.
Then $g_{2,1}$ is $\whatM_{2,1}$-bilipschitz.
\smallskip

We may now write an $R(\ldots)$ statement for $g_{2,1}$.
Since $f_1$ is $M^{\srfs{cmp}}$-bilipschitz,
$f_1(x) = \hatx + (x)_H$ and $f_1(0) = 0$, it follows that
$\norm{x} \geq \frac{\norm{\hatx + (x)_H}}{M^{\trfs{cmp}}}$.
Similarly,
$g_{1,2} \scirc f_1(x) = \hatx + z^{\ssharp}$,
$g_{1,2} \scirc f_1(0) = 0$
and $g_{2,1} \scirc f_1$ is $\whatM_{2,1} M^{\srfs{cmp}}$-bilipschitz.
So
$\norm{x} \leq \whatM_{2,1} M^{\srfs{cmp}} \norm{\hatx + z^{\ssharp}}$.
It follows that
$$\hbox{$
\rfs{supp}(g_{2,1}) \subseteq
B(0;\frac{1}{2 M^{\trfs{cmp}}} \norm{\hatx + (x)_H},
2 M^{\trfs{thn}} \whatM_{2,1} M^{\srfs{cmp}}
\norm{\hatx + z^{\ssharp}}).
$}$$
Hence
\begin{equation}\tag{3.12}
R(\hatx + (z)_H,\hatx + z^{\ssharp},g_{2,1};
\whatM_{2,1},\hbox{$\frac{1}{2 M^{\trfs{cmp}}}$},
2 M^{\trfs{thn}} \whatM_{2,1} M^{\srfs{cmp}},F)
\mbox{ holds.}
\end{equation}

Our next goal is to define $g_{2,2}$.
Recall that $f_1(\hatx) = \hatx$ and
$f_1(\hatx + x^{\nperp}) = \hatx + x^*$.
Also, $f_1$ is $M^{\srfs{cmp}}$\,-\,bilipschitz.
So $\norm{x^*} \approx^{M^{\ssrfs{cmp}}} \norm{x^{\nperp}}$.
Similarly,
$\norm{z^*} \approx^{M^{\ssrfs{cmp}}} \norm{z^{\nperp}}$.
Also, $\norm{x^{\nperp}} = \norm{z^{\nperp}}$.
Let $M_{2,1} = (M^{\srfs{cmp}})^2$ and
$M_{2,2} = M_{2,1} \cdot (M^{\srfs{fdn}})^2$.
It follows that $\norm{x^*} \approx^{M_{2,1}} \norm{z^*}$.
By Proposition \ref{p-bddlip-1.12}(c),
$\norm{\ } \approx^{M^{\ssrfs{fdn}}} \norm{\ }^{\sbfs{N}}$,
and hence
$\norm{x^*}^{\sbfs{N}} \approx^{M_{2,2}} \norm{z^*}^{\sbfs{N}}$.
Since $\norm{z^{\ssharp}}^{\sbfs{N}} = \norm{x^*}^{\sbfs{N}}$, \ 
$\norm{{z^{\ssharp}}}^{\sbfs{N}} \approx^{M_{2,2}}
\norm{z^*}^{\sbfs{N}}$.
Let $a = \frac{\norm{z^*}^{\ssbfs{N}}}{\norm{x^*}^{\ssbfs{N}}}$.
So
\begin{itemize}
\addtolength{\parskip}{-11pt}
\addtolength{\itemsep}{06pt}
\item[(i)] $z^* = a z^{\ssharp}$,
\item[(ii)] $E = F_2 \oplus H$ and
$\norm{u + v}^{\sbfs{N}} = \norm{u} + \norm{v}^{\sbfs{H}}$
for every $u \in F_2$ and $v \in H$,
\item[(iii)] $\hatx \in F_2$ and $z^{\ssharp} \in H$,
\item[(iv)] $\dgfrac{1}{M_{2,2}} \leq a \leq M_{2,2}$.
\vspace{-05.7pt}
\end{itemize}

Assume first that $a \geq 1$.
Let $\hatx,z^{\ssharp},a,0$ take the roles of $\hatx,x,a$ and $u$
in Proposition~\ref{p-gamma.7}(d).
By (i)\,-\,(iii), the assumptions of \ref{p-gamma.7}(d) are fulfilled.
So relying also on (iv),
we conclude that there is $g_{2,2} \in H(E)$ such that
\num{1} $g_{2,2}(\hatx + z^{\ssharp}) = \hatx + z^*$;
\num{2} $g_{2,2} \nrestriction F_2 = \rfs{Id}$;
\num{3} $\rfs{supp}(g_{2,2}) \subseteq
B^{\sbfs{N}}(0;\frac{\norm{\hatx + z^{\sssharp}}^{\ssbfs{N}}}{2},
\frac{3 \norm{\hatx + z^*}^{\ssbfs{N}}}{2})$;
\num{4} $g_{2,2}$ is
$2 M^{\srfs{seg}} \cdot M_{2,2}$\,-\,bilipschitz
in the norm $\norm{\ }^{\sbfs{N}}$.

If $a < 1$ then we apply \ref{p-gamma.7}(d)
to $\hatx,z^*,\dgfrac{1}{a}$ and $0$
thus obtaining a homeomorphism $g_{2,2}' \in H(E)$ such that
$g_{2,2}'(\hatx + z^*) = \hatx + z^{\ssharp}$.
Define $g_{2,2} = (g_{2,2}')\inverse$.
Then (1), (2) and (4) remain true.
Instead of (3) we now have
$\rfs{supp}(g_{2,2}) \subseteq
B^{\sbfs{N}}(0;\frac{\norm{\hatx + z^*}^{\ssbfs{N}}}{2},
\frac{3 \norm{\hatx + z^{\sssharp}}^{\ssbfs{N}}}{2})$.
Note that by (i)\,-\,(iv), 
$\norm{\hatx + z^{\ssharp}}^{\sbfs{N}} \leq
M_{2,2} \norm{\hatx + z^*}^{\sbfs{N}}$.
So
$\rfs{supp}(g_{2,2}) \subseteq
B^{\sbfs{N}}(0;\frac{\norm{\hatx + z^{\sssharp}}^{\ssbfs{N}}}{2M_{2,2}},
\frac{3 M_{2,2} \norm{\hatx + z^*}^{\ssbfs{N}}}{2})$.
Recall that $z^* = (z)_H$. What we have shown implies that
\begin{equation}\tag{3.13}
R(\hatx + z^{\ssharp},\kern1pt\hatx + (z)_H,g_{2,2};\kern1pt
2(M^{\srfs{fdn}})^2 M^{\srfs{seg}} M_{2,2},\kern1pt
\hbox{$\frac{1}{2 M^{\trfs{fdn}} M_{2,2}}$},\kern1pt
2 M^{\srfs{fdn}} M_{2,2},\kern1pt F)
\mbox{ holds.}
\end{equation}
Note that in deducing (3.13) we used the fact that
$\norm{\ }^{\sbfs{N}} \approx^{M^{\trfs{fdn}}} \norm{\ }$.
This concludes the construction of $g_{2,2}$.
\smallskip

{\thickmuskip=3mu \medmuskip=2mu \thinmuskip=1mu 
Define
$\bfie_4 = h_1\inverse \scirc g_{2,2} \scirc g_{2,1} \scirc f_1$
and $\bfix_4 = z$.
Recall that $\bfix_3 = x$.
So $\bfix_4 = z = \bfie_4(\bfix_3)$.
We now
}
apply Proposition~\ref{editor-p.12}(a) and (b).
It follows from (3.10)\,-\,(3.13) and from
\ref{editor-p.12} that there are
$M'_2,A'_2,B'_2$ which do not depend on 
$E,F,\alpha,\bfix_0,\bfiy_0$
such that $R(\bfix_3,\bfix_4,\bfie_4;M'_2,A'_2,B'_2,F)$ holds.

In Case 1 too, we found $M'_1,A'_1,B'_1$
such that
$R(\bfix_3,\bfix_4,\bfie_4;M'_1,A'_1,B'_1,F)$ holds.
Define
$\bfiM_{1,4} = \max(M'_1,M'_2)$,
$a_{1,4} = \min(A'_1,A'_2)$
and $b_{1,4} = \max(B'_1,B'_2)$.
Then $\bfiM_{1,4},a_{1,4},b_{1,4}$
fulfill C1.4 in both Case 1 and Case 2.
\medskip

\noindent
{\bf Part 4 } The construction of $\bfif_1$.
\newline
We have shown that for $i = 1,\ldots,4$ there is $\bfiM_{1,i}$
which does not depend on $E,F,\alpha,\bfix_0,\bfiy_0$
such that $\bfie_i$ is $\bfiM_{1,i}$-bilipschitz.
We define $\bfie = \bfie_4 \scirc \ldots \scirc \bfie_1$.
Then $\bfie(\bfix_0) = \bfix_4 = z$ and $\bfie(0) = 0$.
Let $M_{3,1} = \prod_{i = 1}^4 \bfiM_{1,i}$.
So $\bfie$ is $M_{3,1}$-bilipschitz.
It follows that $\norm{z} \approx^{M_{3,1}} \norm{\bfix_0}$.
Similarly, for $i = 1,2$ there is $\bfiM_{2,i}$
such that $\bfih_i$ is $\bfiM_{2,i}$-bilipschitz.
We define $\bfih = \bfih_2 \scirc \bfih_1$.
Then $\bfih(\bfiy_0) = \bfiy_2 = y$ and $\bfih(0) = 0$.
Let $M_{3,2} = \bfiM_{2,1} \bfiM_{2,2}$.
So $\bfih$ is $M_{3,2}$-bilipschitz.
Let $M_{3,0} = M_{3,1} M_{3,2}$.
Then
\begin{itemize}
\addtolength{\parskip}{-11pt}
\addtolength{\itemsep}{06pt}
\item[{\smaller 4.1}]
$\norm{z} \approx^{M_{3,0}} \norm{\bfix_0}$.
\vspace{-05.7pt}
\end{itemize}
Since $\bfie(F) = F$, we have that
$d(z,F) \approx^{M_{3,1}} d(\bfix_0,F)$.
Similarly,
$d(y,F) \approx^{M_{3,2}} d(\bfiy_0,F)$.
Hence
\begin{itemize}
\addtolength{\parskip}{-11pt}
\addtolength{\itemsep}{06pt}
\item[{\smaller 4.2 }]
$\norm{z} \approx^{M_{3,0} \cdot \alpha} \norm{y}$
\,and\ \,$d(z,F) \approx^{M_{3,0} \cdot \alpha} d(y,F)$.
\vspace{-05.7pt}
\end{itemize}
The construction also implies that
\begin{itemize}
\addtolength{\parskip}{-11pt}
\addtolength{\itemsep}{06pt}
\item[{\smaller 4.3\ }]
$z = \hatz + z^{\nperp}$, $y = \haty + y^{\nperp}$,
where $\hatz, \haty \in F$, and for some $\lambda,\mu > 0$,
$\haty = \lambda \hatz$ and $y^{\nperp} = \mu z^{\nperp}$.\kern-0.1pt
\vspace{-05.7pt}
\end{itemize}

If Case 1 of Part 3 happens, let $\whatF = F$.
Suppose that Case 2 of Part 3 happens.
Let $F_2$ be as defined in Case 2 of Part 3.
So by (3.9), $\norm{z^{\nperp}} \leq \fourthirds d(z,F_2)$.
By Proposition~\ref{p-bddlip-1.11} applied to $F_2$ and taking
$x$ and $y$ to be $z^{\nperp}$, there is a closed subspace
$\whatF$ such that
$\norm{z^{\nperp}} \leq \threehalves d(z^{\nperp},\whatF)$,
$F_2 \subseteq \whatF$
and $\rfs{span}(E \cup \sngltn{z^{\nperp}}) = E$.
In both cases we have
\begin{itemize}
\addtolength{\parskip}{-11pt}
\addtolength{\itemsep}{06pt}
\item[{\smaller 4.4 }]
$F \subseteq \whatF$,
$\whatF \oplus \rfs{span}(\sngltn{z^{\nperp}}) = E$
and $\norm{z^{\nperp}} \leq 1 \half d(z^{\nperp},\whatF)$.
\vspace{-05.7pt}
\end{itemize}

{\bf Case 1 \ } $\norm{\haty} \geq \norm{\hatz}$.
In this case $\lambda \geq 1$.
Let $\bfiv = \haty + z^{\nperp}$.
We shall construct a homeomorphism $\bfif_1$
such that $\bfif_1(z) = \bfiv$.
(Recall that $z = \bfix_4)$.
Denote $\bfiv$ by $v$.
So $v = \lambda \hatz + z^{\nperp}$.
If $\lambda = 1$ let $\bfif_1 = \rfs{Id}$. Assume that $\lambda > 1$.

Let $H = \rfs{span}(\dbltn{\haty}{y^{\nperp}})$,
$H_1 = \rfs{span}(\sngltn{\haty})$ and
$H_2 = \rfs{span}(\sngltn{y^{\nperp}})$.
Let $F_3$ be a subspace of $\whatF$ such that for some
$\varphi \in \whatF^*$,
$\norm{\varphi} = 1$, $\varphi(\hatz) = \norm{\hatz}$ and
$F_3 = \rfs{ker}(\varphi)$.
It follows that $H_1 \oplus F_3 = \whatF$, $\whatF \oplus H_2 = E$
and $\whatF = H_1 \oplus H_2 \oplus F_3$.
Clearly, $\norm{\rfs{Proj}_{H_1,F_3}} = \norm{\varphi} = 1$.
So by Proposition~\ref{p-bddlip-bldr-1.9}(d),
$H_1 \perp^1 F_3$.

Let $S =
\setm{a \hatz + b z^{\nperp}}
{a \in \bbR, \ b \in [0,1]}$.
We define
$\fnn{\eta}{[0,\infty) \times [0,\infty)}{[0,\infty)}$.
For every $s$, $\eta_s(t) \eqdf \eta(s,t)$
is a piecewise linear function of $t$.
For $s \geq (\lambda - 1) \norm{\hatz}$, $\eta_s = \rfs{Id}$.
If $s < (\lambda - 1) \norm{\hatz}$, then $\eta_s(t)$ has breakpoints at
$\dgfrac{\norm{\hatz}}{2}$, $\norm{\hatz}$ and $2 \lambda \norm{\hatz}$;
$\eta_s(t) = t$ for every
$t \in
[0,\dgfrac{\norm{\hatz}}{2}) \cup [2 \lambda \norm{\hatz},\infty)$;
and
$\eta_s(\norm{\hatz}) =
(1 - \frac{s}{(\lambda - 1) \norm{\hatz}}) \ncdot \lambda \norm{\hatz} +
\frac{s}{(\lambda - 1) \norm{\hatz}} \ncdot \norm{\hatz}$.
Denote $(\lambda - 1) \norm{\hatz}$ by $a$.
Then in particular,
$\eta_0(\norm{\hatz}) = \lambda \norm{\hatz}$
and $\eta_a(\norm{\hatz}) = \norm{\hatz}$.

For $u \in E$ we denote $(u)_{H_1},(u)_{H_2},(u)_{F_3}$ by
$(u)_1,(u)_2$ and $(u)_3$ respectively,
and we abbreviate $(u)_i$ by $u_i$ when the notation
$(u)_i$ is too cumbersome.
Set $E^+ = \{t \hatz + w \,|\break t \geq 0, \ w \in H_2 \oplus F_3\}$.
Let $\bfif_1$ be defined by
$$
\bfif_1(u) =\kern2pt
\left\{
\renewcommand{\arraystretch}{1.5}
\addtolength{\arraycolsep}{4pt}
\begin{array}{ll}
\eta(d(u,S),\norm{u_1}) \,\frac{\hatz}{\norm{\hatz}} + u_2 + u_3
&
u \in E^+,
\\
u
&
u \in E - E^+.
\end{array}
\renewcommand{\arraystretch}{1.0}
\addtolength{\arraycolsep}{-4pt}
\right.
$$
\newline
Note that $\bfif_1 \nrestriction H_2 \oplus F_3 = \rfs{Id}$,
so $\bfif_1 \in H(E)$.
We shall define the constants mentioned in C3 and show that C3 holds.
Recall that
C3 \ $\equiv$\ \ 
$R(\bfix_4,\bfiv,\bfif_1;\bfiM_{3,1} \ncdot \alpha,a_{3,1},b_{3,1},F)$.
We verify R1\,-\,R4 in the definition of $R(\ldots)$.

{\bf R1: }
Clearly, $\bfif_1(\bfix_4) = \bfif_1(z) = v = \bfiv$.
\smallskip

{\bf R3: }
We verify that $\bfif_1(F) = F$.
For every $u \in E$ and in particular for every $u \in F$,\break
$\bfif_1(u) - u \in H_1 = \rfs{span}(\sngltn{\haty}) \subseteq F$.
So $\bfif_1(u) = u + (\bfif_1(u) - u) \in F$.
An identical argument shows that $\bfif_1\inverse(F) \subseteq F$.
Hence R3 holds.
\smallskip

{\bf R2: }
We find $\bfiM_{3,1}$
and prove that $\bfif_1$ is $\bfiM_{3,1} \ncdot \alpha$-bicontinuous.
Note that if $g \in H(E)$, $K \subseteq E$ is closed,
$\rfs{supp}(g) \subseteq K$
and $g \nrestriction K$ is $\beta$-continuous, then $g$ is
$2 \beta$-continuous. Since $\rfs{supp}(\bfif_1) \subseteq E^+$,
we may consider only points which belong to $E^+$.
Let $u,w \in E^+$. Then
\vspace{1.5mm}
\newline
\rule{-2pt}{0pt}
\renewcommand{\arraystretch}{1.5}
\addtolength{\arraycolsep}{-3pt}
$
\begin{array}{ll}
&
\norm{\bfif_1(w) - \bfif_1(u)} \kern1pt\leq\kern1pt
\abs{\eta(d(w,S),\norm{w_1}) - \eta(d(u,S),\norm{u_1})}
\kern1pt+\kern1pt
\norm{(w - u)_2} \kern1pt+\kern1pt
\norm{(w - u)_3}
\\
\leq
\rule{5pt}{0pt}
&
\abs{\eta(d(w,S),\norm{w_1}) - \eta(d(u,S),\norm{w_1})}
\kern1pt+\kern1pt
\abs{\eta(d(u,S),\norm{w_1}) - \eta(d(u,S),\norm{u_1})}
\\
+
\rule{5pt}{0pt}
&
\norm{(w - u)_2} \kern4pt+\kern4pt \norm{(w - u)_3}.
\vspace{1.7mm}
\end{array}
$
\renewcommand{\arraystretch}{1.0}
\addtolength{\arraycolsep}{3pt}
\newline
That is,
\vspace{1.5mm}
\newline
(4.1)
\rule{0pt}{0pt}
\renewcommand{\arraystretch}{1.5}
\addtolength{\arraycolsep}{-3pt}
$
\begin{array}[t]{ll}
&
\norm{\bfif_1(w) - \bfif_1(u)} \leq
\abs{\eta(d(w,S),\norm{w_1}) \kern1pt-\kern1pt \eta(d(u,S),\norm{w_1})}
\\
+\kern2pt
\rule{5pt}{0pt}
&
\abs{\eta(d(u,S),\norm{w_1}) - \eta(d(u,S),\norm{u_1})}
\kern4pt+\kern4pt
\norm{(w - u)_2} \kern4pt+\kern4pt \norm{(w - u)_3}.
\vspace{1.7mm}
\end{array}
$
\renewcommand{\arraystretch}{1.0}
\addtolength{\arraycolsep}{3pt}
\newline
The first summand in the right hand side of (4.1) has the form
$\abs{\eta(s_1,t) - \eta(s_2,t)}$.
If $s_1,s_2 \in [0,(\lambda - 1) \norm{\hatz}]$, then
$$
\abs{\eta(s_1,t) - \eta(s_2,t)} =
\frac{\abs{s_1 - s_2}}{(\lambda - 1) \norm{\hatz}} \ncdot
(\eta(0,t) - \eta((\lambda - 1) \norm{\hatz},t)) \leq
\frac{\lambda \norm{\hatz} - \norm{\hatz}}{(\lambda - 1) \norm{\hatz}}
\ncdot \abs{s_1 - s_2} = \abs{s_1 - s_2}.
$$

\noindent
The inequality between the first and last expression above is true
\kern-2pt\rule{0pt}{15pt}
for every\rule{3pt}{0pt}
$s_1,s_2 \in [0,\infty)$.
So \ \,%
$\abs{\eta(d(w,S),\norm{w_1}) - \eta(d(u,S),\norm{w_1})} \leq
\abs{d(w,S) - d(u,S)} \leq \norm{w - u}$.
That is,
\begin{equation}\tag{4.2}
\abs{\eta(d(w,S),\norm{w_1}) - \eta(d(u,S),\norm{w_1})}
\kern4pt\leq\kern4pt \abs{d(w,S) - d(u,S)} \kern4pt\leq\kern4pt
\norm{w - u}.
\end{equation}

The next computations are needed in order to estimate the
second summand in the right hand side of (4.1).
We find $A,B,C$ such that
$A \norm{z} \leq \norm{\hatz} \leq B \norm{z}$
and $\norm{z^{\nperp}} \leq C \norm{z}$.
There are different computations corresponding to Cases 1 and 2 of
Part 3.

In Case 1 of Part 3, $\itDelta = 8$ and $\varepsilon = \dghalf$.
So
$d(x,F) \leq \dgfrac{\norm{x}}{8}$ and
$\norm{x^{\nperp}} \leq 1 \half d(x,F)$.
Hence
$\norm{z^{\nperp}} = \norm{x^{\nperp}} \leq
\frac{3}{2} \cdot \eighth \norm{x} =
\frac{3}{16} \norm{x}$.
We have that $z = x - x^{\nperp} + z^{\nperp}$. Hence
$\norm{z} \geq \norm{x} - \norm{x^{\nperp}} - \norm{z^{\nperp}} =
\norm{x} - 2 \norm{x^{\nperp}}$.
Hence $\norm{z} \geq \norm{x} - \frac{3}{8} \norm{x} =
\frac{5}{8} \norm{x}$.
It follows that
$\norm{z^{\nperp}} \leq \frac{3}{16} \cdot \frac{8}{5} \norm{z}$.
That is,
\begin{equation}
\tag{4.4.1}
\norm{z^{\nperp}} \leq \hbox{$\frac{3}{10}$} \norm{z}.
\end{equation}
From the fact that 
$\hatz = z - z^{\nperp}$, we conclude
\begin{equation}
\tag{4.5.1}
\hbox{$\frac{7}{10}$} \norm{z} \leq \norm{\hatz} \leq
\hbox{$\frac{13}{10}$} \norm{z}.
\end{equation}

Recall that in Case 2 of Part 3,
$\itDelta = 24$ and $\varepsilon = \dgninth$.
We carry out a computation similar to the one in Case 1.
So $d(x,F) \leq \dgfrac{\norm{x}}{24}$ and
$\norm{x^{\nperp}} \leq 1 \ninth d(x,F)$.
So
$\norm{z^{\nperp}} = \norm{x^{\nperp}} \leq
\frac{10}{9 \cdot 24} \norm{x}
=
\frac{5}{108} \norm{x}$.
We have that
$\norm{z} \geq
\norm{x} - 2 \norm{x^{\nperp}} \geq \norm{x} - \frac{5}{54} \norm{x} =
\frac{49}{54} \norm{x}$
and hence
$\norm{z^{\nperp}} \leq \frac{5}{108} \cdot \frac{54}{49} \norm{z}$.
That is,
\begin{equation}
\tag{4.4.2}
\norm{z^{\nperp}} \leq \hbox{$\frac{5}{98}$} \norm{z}
\end{equation}
and hence
\begin{equation}
\tag{4.5.2}
\hbox{$\frac{93}{98}$} \norm{z} \leq
\norm{\hatz} \leq \hbox{$\frac{103}{98}$} \norm{z}.
\end{equation}
By (4.4.1) and (4.4.2),
\begin{equation}
\tag{4.4}
\norm{z^{\nperp}} \leq \hbox{$\frac{3}{10}$} \norm{z},
\end{equation}
and by (4.5.1) and (4.5.2),
\begin{equation}
\tag{4.5}
\hbox{$\frac{7}{10}$} \norm{z} \leq
\norm{\hatz} \leq \hbox{$\frac{13}{10}$} \norm{z}.
\end{equation}

Since $y$ also obeys {\smaller 3.3, 3.4},
in Case 1 of Part 3 we obtain that
$\frac{13}{16} \norm{y} \leq \norm{\haty} \leq
\frac{19}{16} \norm{y}$
and in Case 2,
$\frac{103}{108} \norm{y} \leq \norm{\haty} \leq
\frac{113}{108} \norm{y}$.
The following is thus true in both cases.
\begin{equation}
\tag{4.6}
\hbox{$\frac{13}{16}$} \norm{y} \leq \norm{\haty} \leq
\hbox{$\frac{19}{16}$} \norm{y}.
\end{equation}

{\thickmuskip=2mu \medmuskip=1mu \thinmuskip=1mu 
By {\smaller 4.3}, (4.6), {\smaller 4.2}, (4.5),
the monotonicity
of $\alpha$
and the fact that
$\alpha(A t) \leq A \alpha(t)$ for $A \geq 1$,
}
\vspace{1.5mm}
\newline
\rule{7pt}{0pt}
\renewcommand{\arraystretch}{1.5}
\addtolength{\arraycolsep}{-0pt}
$
\begin{array}{lll}
\lambda \norm{\hatz}
&
=
&
\norm{\haty} \leq
\frac{19}{16} \norm{y} \leq
\frac{19}{16} M_{3,0} \ncdot \alpha(\norm{z}) \leq
\frac{19}{16} M_{3,0} \ncdot \alpha(\frac{10}{7}\norm{\hatz})
\\
&
\leq
&
\frac{10}{7} \ncdot \frac{19}{16} M_{3,0} \ncdot \alpha(\norm{\hatz})
\leq
2 M_{3,0} \ncdot \alpha(\norm{\hatz}).
\vspace{1.7mm}
\end{array}
$
\renewcommand{\arraystretch}{1.0}
\addtolength{\arraycolsep}{0pt}
\newline
So
\begin{equation}
\tag{4.7}
\lambda \norm{\hatz} \leq 2 M_{3,0} \ncdot \alpha(\norm{\hatz}).
\end{equation}

Let $\rho = \eta_0$.
So $\rho$ is the piecewise linear function with breakpoints at
$\frac{\norm{\hatz}}{2}$, $\norm{\hatz}$ and $2 \lambda \norm{\hatz}$;
$\rho(t) = t$ for every
$t \in
[0,\dgfrac{\norm{\hatz}}{2}) \cup [2 \lambda \norm{\hatz},\infty)$;
and
$\rho(\norm{\hatz}) =
\lambda \norm{\hatz}$.
Clearly, $\rho \in H([0,\infty))$.
Using the notations of Proposition \ref{p-bddlip-1.11-1}(b),
$\eta = \eta_{(\rho,(\lambda - 1) \ncdot \norm{\hatz})}$.

We show that $\rho$ is $16 M_{3,0} \ncdot \alpha$-continuous.
The linear pieces of $\rho$ have the slopes:
$1$,
$\frac{\lambda \norm{\hatz} - \dgfrac{\norm{\hatz}}{2}}
{\dgfrac{\norm{\hatz}}{2}}$,
$\frac{2 \lambda \norm{\hatz} - \lambda \norm{\hatz}}
{2 \lambda \norm{\hatz} - \norm{\hatz}}$ and $1$.
That is, the slopes of the linear pieces of $\rho$ are
$1$, $2 \lambda - 1$ and $\frac{\lambda}{2 \lambda - 1}$.
We use the notations of Definition~\ref{d-bddly-lip-2.9-1}(b).
Let $a_0,\ldots,a_4$ denote
$0$, $\dgfrac{\norm{\hatz}}{2}$,
$\norm{\hatz}$, $2 \lambda \norm{\hatz}$ and $\infty$.
Then $\rho_1,\ldots,\rho_4$ are the functions
\vspace{1.5mm}
\newline
\rule{7pt}{0pt}
\renewcommand{\arraystretch}{1.5}
\addtolength{\arraycolsep}{-3pt}
$
\begin{array}{ll}
\rfs{Id} \nrestriction [0,\dgfrac{\norm{\hatz}}{2}],
&
\\
y = (2 \lambda - 1) t + \dgfrac{\norm{\hatz}}{2},
\rule{15pt}{0pt}
&
t \in [0,\dgfrac{\norm{\hatz}}{2}],
\\
y = \frac{\lambda}{2 \lambda - 1}t + \lambda \norm{\hatz},
\rule{5pt}{0pt}
&
t \in [0,(2 \lambda - 1) \norm{\hatz}],
\\
y = t + 2 \lambda \norm{\hatz},
&
t \in [0,\infty).
\vspace{1.7mm}
\end{array}
$
\renewcommand{\arraystretch}{1.0}
\addtolength{\arraycolsep}{3pt}
\newline
For $i = 1, 3, 4$, for every $t_1, t_2$,
$\abs{\rho_i(t_1) - \rho_i(t_2)} \leq \abs{t_1 - t_2} \leq
4 M_{3,0} \ncdot \alpha(\abs{t_1 - t_2})$.
Hence $\rho_i$ is $4 M_{3,0} \ncdot \alpha$-continuous.
We deal with $\rho_2$.
By (4.7),
$2 \lambda - 1 \leq 2 \lambda \leq
\frac{4 M_{3,0} \ncdot \alpha(\norm{\hatz})}{\norm{\hatz}}$.
So $\frac{2 \lambda - 1}{4 M_{3,0}} \leq
\frac{\alpha(\norm{\hatz})}{\norm{\hatz}}$.
Let $\rho_2^*(t)$ be the function
$$
y = \hbox{$\frac{2 \lambda - 1}{4 M_{3,0}}$} t,
\ \,t \in [0,\norm{\hatz}].
$$
Then by Proposition \ref{p-bddlip-1.11-1}(c),
$\rho_2^*(t)$ is $\alpha$-continuous.
Clearly,
$\rho_2(t) = 4 M_{3,0} \ncdot \rho_2^*(t) + \dgfrac{\norm{\hatz}}{2}$.
So $\rho_2$ is $4 M_{3,0} \ncdot \alpha$-continuous.
We have shown that $\rho$ is $(4,4 M_{3,0} \ncdot \alpha)$-continuous.
By Proposition~\ref{p-bddlip-1.11-1}(a), $\rho$ is
$16 M_{3,0} \ncdot \alpha$-continuous.
Denote $\gamma = 16 M_{3,0} \ncdot \alpha$.

We next deal with the second summand in the right hand side of
inequality (4.1).
It has the form $\abs{\eta(s,t_1) - \eta(s,t_2)}$.
Recall that $\eta = \eta_{(\rho,(\lambda - 1) \ncdot \norm{\hatz})}$.
Then by Proposition~\ref{p-bddlip-1.11-1}(b),
for every $s \in [0,\infty)$,
$\eta_s$ is $\gamma$-continuous. So
$$\abs{\eta(d(u,S),\norm{w_1}) - \eta(d(u,S),\norm{u_1})} \leq
\gamma(\abs{\kern0.7pt\norm{w_1} - \norm{u_1}\kern0.5pt}) =
\gamma(\norm{(w - u)_1}).
$$
That is,
\newline
(4.8)
\centerline{
$\abs{\eta(d(u,S),\norm{w_1}) - \eta(d(u,S),\norm{u_1})} \leq
\gamma(\norm{(w - u)_1})$.
}

We shall now bound the expressions $\norm{(w - u)_i}$
appearing in (4.1) and (4.8) by a multiple of $\norm{w - u)}$.
For $\baru \in E$ let $\baru_{1,3} = \baru_1 + \baru_3$.
Recall that $H_2 = \rfs{span}(\sngltn{y^{\nperp}})$.
By {\smaller 4.4},
$\norm{y^{\nperp}} \leq \frac{3}{2} d(y^{\nperp},\whatF)$.
Hence
$\norm{\baru_2} \leq \frac{3}{2} d(\baru_2,\whatF) \leq
\frac{3}{2} \norm{\baru}$.
From the fact that $\baru_{1,3} = \baru - \baru_2$, it follows that
$\norm{\baru_{1,3}} \leq \norm{\baru} + \norm{\baru_2} \leq
\frac{5}{2} \norm{\baru}$.
So we have
\begin{equation}\tag{$\star$}
\norm{\baru_{1,3}} \leq \hbox{$\frac{5}{2}$} \norm{\baru}.
\end{equation}
From the fact that $H_1 \perp^1 F_3$, it follows that
$\norm{\baru_1} \leq \norm{\baru_{1,3}}$,
and this implies that $\norm{\baru_3} \leq 2 \norm{\baru_{1,3}}$.
It follows that
\begin{equation}
\tag{4.9}
\norm{\baru_1} \leq \hbox{$\frac{5}{2}$} \norm{\baru}, \kern6pt
\norm{\baru_2} \leq \hbox{$\frac{3}{2}$} \norm{\baru}
\kern6pt \mbox{and} \kern6pt
\norm{\baru_3} \leq 5 \norm{\baru}.
\end{equation}
Substituting (4.2) and (4.8) into (4.1), we obtain
\begin{equation}
\tag{4.10}
\norm{\bfif_1(w) - \bfif_1(u)} \leq
\norm{w - u} + \gamma(\norm{(w - u)_1}) + \norm{(w - u)_2} +
\norm{(w - u)_3}.
\end{equation}
We substitute (4.9) into (4.10) and use
Proposition~\ref{p-bddlip-1.11-1}(d). So
$$
\norm{\bfif_1(w) - \bfif_1(u)} \leq
\hbox{$7 \half$} \norm{w - u} +
\hbox{$2 \half$} \gamma(\norm{w - u}).
$$
This means that $\bfif_1 \nrestriction E^+$ is
$(40 M_{3,0} \ncdot \alpha + 7 \half \rfs{Id})$-continuous.
Hence $\bfif_1 \nrestriction E^+$
is $50 M_{3,0} \ncdot \alpha$-continuous.
It follows that $\bfif_1$ is $100 M_{3,0} \ncdot \alpha$-continuous.

The computation which shows that for some $M$,
$\bfif_1\inverse$ is $M \ncdot \alpha$-continuous is
analogous. However, for $\bfif_1\inverse$
there is $M$ which does not depend on $E,F,\alpha,\bfix_0,\bfiy_0$
such that $\bfif_1\inverse$ is even $M$-Lipschitz.
For this $M$ it is also true that
$\bfif_1\inverse$ is $M \ncdot \alpha$-continuous.
We now carry out the computation for $\bfif_1\inverse$.
For $s \in [0, \infty)$ let $\theta_s = \eta_s\inverse$.
Denote $\theta(s,t) = \theta_s(t)$.
Note that for every $u \in E$, $d(\bfif_1(u),S) = d(u,S)$.
This implies that
$$
\bfif_1\inverse(u) =
\theta_{d(u,S)}(\norm{u_1}) \frac{u_1}{\norm{u_1}} + u_2 + u_3.
$$
The analogues (4.1$^*$) of (4.1) and (4.2$^*$) of (4.2)
obtained by replacing $\eta$ by
$\theta$ are still true.
Let $\mu = \theta_0$. So $\mu = \rho\inverse$
and $\theta = \eta_{(\mu,(\lambda - 1) \ncdot \norm{\hatz})}$.
The slopes of the linear pieces of $\mu$ are the inverses of the
slopes of the linear pieces of $\rho$.
Hence the slopes of the linear pieces of $\mu$ are:
$1$, $\frac{1}{2 \lambda - 1}$ and $\frac{2 \lambda - 1}{\lambda}$.
Clearly,
$1, \frac{1}{2 \lambda - 1}, \frac{2 \lambda - 1}{\lambda} \leq 2$.
So $\mu$ is $2$-Lipschitz.
By Proposition \ref{p-bddlip-1.11-1}(b), for every $s \in [0,\infty)$,
$\theta_s$ is $2$-Lipschitz.
Hence
$$
\abs{\theta(d(u,S),\norm{w_1}) - \theta(d(u,S),\norm{u_1})} \leq
2 \ncdot\abs{\kern0.7pt\norm{w_1} - \norm{u_1}\kern0.5pt} =
2 \ncdot \norm{(w - u)_1}.
$$
So (4.8) is replaced by
\begin{equation}\tag{4.11}
\abs{\theta(d(u,S),\norm{w_1}) - \theta(d(u,S),\norm{u_1})} \leq
2 \ncdot \norm{(w - u)_1}.
\end{equation}
Substituting (4.9) into (4.11) we get
\begin{equation}\tag{4.12}
\abs{\theta(d(u,S),\norm{w_1}) - \theta(d(u,S),\norm{u_1})} \leq
5 \norm{w - u}.
\end{equation}
Replace the first summand of the right hand side of (4.1$^*$)
by (4.2$^*$) and the second summand by (4.12).
Use (4.9) to estimate the last two summands of
(4.1$^*$).
So
$$
\norm{\bfif_1\inverse(w) - \bfif_1\inverse(u)} \leq
\norm{w - u} + 5 \norm{w - u} + \hbox{$\frac{3}{2}$} \norm{w - u} +
5 \norm{w - u}.
$$
That is,
$\norm{\bfif_1\inverse(w) - \bfif_1\inverse(u)} \leq
12 \half \norm{w - u}$.
From the fact that $\alpha \geq \rfs{Id}$, it follows that
$\bfif_1\inverse \nrestriction E^+$ is
$12 \half \ncdot \alpha$-continuous.
It follows that $\bfif_1\inverse$ is
$25 \ncdot \alpha$-continuous.
Hence
$\bfif_1$ is $100 M_{3,0} \ncdot \alpha$-bicontinuous.
So $\bfiM_{3,1} = 100 M_{3,0}$.

\kern1.3mm

{\bf R4: }
We next find $a_{3,1}$ such that
$\bfif_1 \nrestriction B(0, a_{3,1} \norm{z}) = \rfs{Id}$.
For every $t \leq \dgfrac{\norm{\hatz}}{2}$ and for every~$s$,
$\eta(s,t) = t$. So for every $u \in E$,
if $\norm{u_1} \leq \dgfrac{\norm{\hatz}}{2}$,
then $\bfif_1(u) = u$.
By (4.9) and the above, if $\norm{u} \leq \fifth \norm{\hatz}$,
then $\bfif_1(u) = u$.
By (4.5),
$\frac{7}{10} \norm{z} \leq \norm{\hatz}$.
So if $\norm{u} \leq \frac{7}{50} \norm{z}$,
then $\bfif_1(u) = u$.
Let $a_{3,1} = \frac{7}{50}$,
then $\bfif_1 \nrestriction B(0, a_{3,1} \norm{z}) = \rfs{Id}$.

We now find $b_{3,1}$ such that
$\rfs{supp}(\bfif_1) \subseteq B(0, b_{3,1} \norm{v})$.
We shall find $A_i$, \, $i = 1, 2, 3$, such that for every $u \in E$:
if $\norm{u_i} \geq A_i$, then $\bfif_1(u) = u$.
For every $t \geq 2 \lambda \hatz$ and every $s$,
$\eta(s,t) = t$.
So
\begin{equation}\tag{4.13}
\rule{0pt}{0pt}\kern-6mm
\hbox{If $\norm{u_1} \geq 2 \lambda \norm{\hatz}$,
then $\bfif_1(u) = u$.}\rule{8cm}{0pt}
\end{equation}
For every $s \geq (\lambda - 1) \norm{\hatz}$, $\eta_s = \rfs{Id}$.
So for every $u \in E$, if $d(u,S) \geq (\lambda - 1) \norm{\hatz}$,
then $\bfif_1(u) = u$.
Let $u \in E$. By the second part of (4.9),
$\norm{u} \geq \frac{2}{3} \norm{u_2}$.
Let $a > 0$.
If $\norm{u_2} \geq a + \norm{z^{\nperp}}$,
then for every $w \in S$, $\norm{(u - w)_2} \geq a$.
Hence
$\norm{u - w} \geq \frac{2}{3} \norm{(u - w)_2} \geq \frac{2}{3} a$.
Take $a = \frac{3}{2} (\lambda - 1) \norm{\hatz}$.
So if
$\norm{u_2} \geq
\frac{3}{2} (\lambda - 1) \norm{\hatz} + \norm{z^{\nperp}}$,
then for every $w \in S$,
$\norm{u - w} \geq (\lambda - 1) \norm{\hatz}$.
That is, if
$\norm{u_2} \geq
\frac{3}{2} (\lambda - 1) \norm{\hatz} + \norm{z^{\nperp}}$,
then $d(u,S) \geq (\lambda - 1) \norm{\hatz}$. Hence
\begin{equation}\tag{4.14}
\rule{0pt}{0pt}\kern-6mm
\hbox{If
$\norm{u_2} \geq
\frac{3}{2} (\lambda - 1) \norm{\hatz} + \norm{z^{\nperp}}$,
then $\bfif_1(u) = u$.}
\rule{4cm}{0pt}
\end{equation}
The third part of (4.9) says that
$\norm{\baru} \geq \fifth \norm{\baru_3}$
for every $\baru \in E$.
Let $u \in E$ be such that
$\norm{u_3} \geq 5 (\lambda - 1) \norm{\hatz}$.
For every $w \in S$, $(u - w)_3 = u_3$.
So
$\norm{(u - w)} \geq \fifth \norm{(u - w)_3} = \fifth \norm{u_3} \geq
(\lambda - 1) \norm{\hatz}$.
That is, $d(u,S) \geq (\lambda - 1) \norm{\hatz}$.
\begin{equation}\tag{4.15}
\rule{0pt}{0pt}\kern-6mm
\hbox{If $\norm{u_3} \geq 5 (\lambda - 1) \norm{\hatz}$,
then $\bfif_1(u) = u$.}\rule{5.5cm}{0pt}
\end{equation}
Combining (4.13) - (4.15) we conclude that
\begin{equation}\tag{4.16}
\rule{0pt}{0pt}\kern-6mm
\hbox{If
$\norm{u_1} + \norm{u_2} + \norm{u_3} \geq
(8 \half \lambda - 6 \half) \norm{\hatz} + \norm{z^{\nperp}}$,
then $\bfif_1(u) = u$.}\rule{2.1cm}{0pt}
\end{equation}
By {\smaller 4.4},
$\norm{z^{\nperp}} \leq 1 \half d(z^{\nperp},\whatF) =
1 \half d(z,\whatF) \leq 1 \half \norm{z}$,
and by (4.5),
$\norm{\hatz} \leq \frac{13}{10} \norm{z}$.
So
\begin{equation}\tag{4.17}
\hbox{$
(8 \half \lambda - 6 \half) \norm{\hatz} + \norm{z^{\nperp}} \leq
(\frac{13}{10} \ncdot (8 \half \lambda - 6 \half) + 1 \half) \norm{z}
\leq
10 \lambda \norm{z}$.}
\end{equation}
Note that
$z = \hatz + z^{\nperp} =
\frac{1}{\lambda} v - \frac{1}{\lambda} z^{\nperp} + z^{\nperp} =
\frac{1}{\lambda} v + (1 - \frac{1}{\lambda}) z^{\nperp}$.
Hence
$\norm{z} \leq
\frac{1}{\lambda} \norm{v} + \norm{z^{\nperp}}$.
By (4.4),
$\norm{z} \leq
\frac{1}{\lambda} \norm{v} + \frac{3}{10}\norm{z}$.
So
\begin{equation}\tag{4.18}
\norm{z} \leq \hbox{$\frac{10}{7 \lambda}$} \norm{v}.
\end{equation}
From (4.16), (4.17) and (4.18) we conclude that
\begin{equation}\tag{4.19}
\rule{0pt}{0pt}\kern-6mm
\hbox{If $\norm{u} \geq
\frac{100}{7} \ncdot \norm{v}$,
then $\bfif_1(u) = u$.}
\rule{7.5cm}{0pt}
\end{equation}
That is,
$\rfs{supp}(\bfif_1) \subseteq B(0,\frac{100}{7} \ncdot \norm{v})$.
So $b_{3,1} \eqdf \frac{100}{7}$ is as required in R4.
\smallskip

{\bf Case 2 } $\norm{\haty} < \norm{\hatz}$.
So $\lambda < 1$.
Let $\bfiv = v = \lambda z$,
and we construct $\bfif_1$ such that $\bfif_1(z) = v$.
By (4.6), $\norm{\haty} \geq \frac{13}{16} \norm{y}$,
and by (4.5), $\norm{\hatz} \leq \frac{13}{10} \norm{z}$.
So (i)
$\lambda = \frac{\norm{\haty}}{\norm{\hatz}} \geq
\frac{5}{8} \frac{\norm{y}}{\norm{z}}$.
By the construction of $\bfih_1$ and $\bfih_2$,
(ii) $\norm{y} = \norm{\bfiy_0}$.
By {\smaller 4.1}, (iii) $\norm{z} \approx^{M_{3,0}} \norm{\bfix_0}$.
Since $\bfix_0,\bfiy_0$ satisfy conditions A1\,-\,A4 appearing in the
definition of a UC-constant,
(iv) $\norm{\bfiy_0} \geq \norm{\bfix_0}$.
So by (i)-(iv),
\begin{equation}\tag{4.20}
\hbox{$\lambda \geq
\frac{5}{8} \frac{\norm{\sbfiy_0}}{M_{3,0} \norm{\sbfix_0}} \geq
\frac{1}{2 M_{3,0}}$}.
\end{equation}
Let $\fnn{\eta}{[0,\infty)}{[0,\infty)}$ be a piecewise linear function
with breakpoints
at $\frac{\lambda \norm{z}}{2}$, $\norm{z}$
and $2 \norm{z}$ such that
$\eta \rest
(\,[0,\frac{\lambda \norm{z}}{2}] \,\cup\, [2 \norm{z},\infty)\,) =
\rfs{Id}$
and $\eta(\norm{z}) = \lambda \norm{z}$.
Define $\bfif_1$ to be the piecewise linearly radial homeomorphism
based on $\eta$. (See Definition~\ref{d-bddlip-1.13}(b)).
Recall that $z = \bfix_4$, $\bfiv = v$.
We shall define $\bfiM\fprimei_{1,3}$, $a_{1,3}'$ and $b_{1,3}'$
such that
$R(\bfix_4,\bfiv,\bfif_1,\bfiM\fprimei_{1,3} \ncdot \alpha,
a_{1,3}',b_{1,3}',F)$
holds.

{\bf R1} and {\bf R3:} Obviously, $\bfif_1(\bfix_4) = \bfiv$
and $\bfif_1(F) = F$.

{\bf R2:} The slopes of the linear pieces of
$\eta$ are $1$,
$\frac{\half \lambda \norm{z}}{\norm{z} - \half \lambda \norm{z}}$,
and $\frac{2 \norm{z} - \lambda \norm{z}}{\norm{z}}$.
That is, they are
1, $\frac{\lambda}{2 - \lambda}$ and $2 - \lambda$.
Now, $\frac{\lambda}{2 - \lambda} \leq 1$
and by~(4.20), $\frac{1}{4 M_{3,0}} \leq \frac{\lambda}{2 - \lambda}$.
That is, $\frac{1}{4 M_{3,0}} \leq \frac{\lambda}{2 - \lambda} \leq 1$.
Also, $1 \leq 2 - \lambda \leq 2$.
Hence the slopes of all linear pieces of $\eta$ and $\eta\inverse$
are $\leq 4 M_{3,0}$.
So $\eta$ is $4 M_{3,0}$-bilipschitz.
By Proposition~\ref{metr-bldr-p3.18},
$\bfif_1$ is $12 M_{3,0}$-bilipschitz.
Since $\alpha \geq \rfs{Id}$,
$\bfif_1$ is $12 M_{3,0} \ncdot \alpha$-bicontinuous.
We may thus define $\bfiM\fprimei_{3,1} = 12 M_{3,0}$.

{\bf R4:} Obviously,
$\rfs{supp}(\bfif_1) \subseteq
B(0;\frac{\lambda \norm{z}}{2},2 \norm{z})$.
By (4.20),
$B(0,\frac{1}{4 M_{3,0}} \norm{z}) \subseteq
B(0,\frac{\lambda \norm{z}}{2})$.
So we may define $a_{3,1}' = \frac{1}{4 M_{3,0}}$.
Recall that $v = \lambda z$. So by (4.20),
$\norm{v} = \lambda \norm{z} \geq \frac{1}{2 M_{3,0}} \norm{z}$.
Hence $2 \norm{z} \leq 4M_{3,0} \norm{v}$.
It~follows that
$B(0,2 \norm{z}) \subseteq B(0,4 M_{3,0} \norm{v})$.
So we may take $b_{3,1}' = 4 M_{3,0}$.

We have shown that
$R(\bfix_4,\bfiv,\bfif_1,\bfiM\fprimei_{1,3} \ncdot \alpha,
a_{1,3}',b_{1,3}',F)$
holds.
Taking in account Case 1 and Case 2, we define
$\bfiM\fprime'_{3,1} = \max(\bfiM_{3,1},\bfiM\fprimei_{3,1})$,
$a''_{3,1} = \min(a_{3,1},a'_{3,1})$ and
$b''_{3,1} = \max(b_{3,1},b'_{3,1})$.
Then $\bfiM\fprime'_{3,1},a''_{3,1}, b''_{3,1}$ are as required
in C3.
\medskip

\noindent
{\bf Part 5 } The construction of $\bfif_2$.
\\
Let $v$ be as in Part 4. Remember that $v$ was defined in two
different ways.
In the case that $\norm{\hatz} \leq \norm{\haty}$,
$v = \haty + z^{\nperp}$
and in the case that $\norm{\hatz} > \norm{\haty}$, $v = \lambda z$.
Denote $v^{\nperp} = v - \haty$. The following holds.
\begin{itemize}
\addtolength{\parskip}{-11pt}
\addtolength{\itemsep}{06pt}
\item[{\smaller 5.1 }] \,$y = \haty + y^{\nperp}$,
$v = \haty + v^{\nperp}$,
$y^{\nperp} = \nu v^{\nperp}$, $\haty \in F$ and $\nu > 0$.
\vspace{-05.7pt}
\end{itemize}

If $\nu = 1$ let $\bfif_2 = \rfs{Id}$. Assume that $\nu \neq 1$.
The vector $y^{\nperp}$ is as in Part 4,
and in both Cases 1 and 2 of Part 4,
$v^{\nperp}$ is a multiple of~$y^{\nperp}$.
So the analogue of Clause {\smaller 4.4} in Part 4 holds
for $y^{\nperp}$ and $v^{\nperp}$. That is,
\begin{itemize}
\addtolength{\parskip}{-11pt}
\addtolength{\itemsep}{06pt}
\item[{\smaller 5.2 }]
\,$F \subseteq \whatF$,
$\whatF \oplus \rfs{span}(\sngltn{y^{\nperp}}) = E$
and
$\norm{y^{\nperp}} \leq 1 \half d(y^{\nperp},\whatF)$
and  equivalently
$\norm{v^{\nperp}} \leq 1 \half d(v^{\nperp},\whatF)$.
\vspace{-05.7pt}
\end{itemize}
Recall that $\bfig_1 = \bfif_1 \scirc \bfie$.
We shall next show that there is $\bfiN_1$ which does not depend on
$E,F,\alpha,\bfix_0,\bfiy_0$ such that
\begin{list}{}
{\setlength{\leftmargin}{31pt}
\setlength{\labelsep}{05pt}
\setlength{\labelwidth}{20pt}
\setlength{\itemindent}{-00pt}
\addtolength{\topsep}{-04pt}
\addtolength{\parskip}{-02pt}
\addtolength{\itemsep}{-05pt}
}
\item[\phantom{\sclub}(\sclub)]
for every $u \in E$,
$d(\bfig_1(u),F) \approx^{\ssbfiN_1} d(u,F)$.
In particular, $d(v,F) \approx^{\ssbfiN_1} d(\bfix_0,F)$,
\vspace{-02.0pt}
\end{list}
Recall that $M_{3,1} = \prod_{i = 1}^4 \bfiM_{1,i}$.
Then by C1, $d(\bfie(u),F) \approx^{M_{3,1}} d(u,F)$
for every $u \in E$.
In Case 1 of Part 4, $\bfif_1(u) - u \in F$ for every $u \in E$,
so $d(\bfif_1(u),F) = d(u,F)$.
So in Case 1 of Part 4,
$d(\bfig_1(u),F) \approx^{\ssbfiN_1} d(u,F)$ for every $u \in E$.

In Case 2 of Part 4, $\bfif_1$ is the piecewise linearly radial
homeomorphism based on $\eta$, and for any slope $a$ of a piece of
$\eta$, $\frac{1}{4 M_{3,0}} \leq a \leq 2 \leq 4 M_{3,0}$.
So for every $u \in E$,
$d(u,F) \approx^{4 M_{3,0}} d(\bfif_1(u),F)$.
Now, define $\bfiN_1 = 4 M_{3,1} M_{3,0}$.
Then in both Case 1 and Case 2 of Part 4,
$d(\bfig_1(u),F) \approx^{\ssbfiN_1} d(u,F)$
for every $u \in E$.
The fact $d(v,F) \approx^{\sbfiN_1} d(\bfix_0,F)$
is a special case of the above, since $v = \bfig_1(\bfix_0)$.

It is given that $d(\bfix_0,F) \approx^{\kern1pt\alpha} d(\bfiy_0,F)$.
Let $\bfiN_2 = \bfiM_{2,1} \bfiM_{2,2}$.
Then from C2 it follows that
$d(\bfiy_0,F) \approx^{\ssbfiN_2} d(\bfiy,F)$.
So
\begin{list}{}
{\setlength{\leftmargin}{31pt}
\setlength{\labelsep}{05pt}
\setlength{\labelwidth}{20pt}
\setlength{\itemindent}{-00pt}
\addtolength{\topsep}{-04pt}
\addtolength{\parskip}{-02pt}
\addtolength{\itemsep}{-05pt}
}
\item[(\sclub\sclub)]
$d(\bfix_0,F) \approx^{\ssbfiN_2 \ncdot \alpha} d(y,F)$.
\vspace{-02.0pt}
\end{list}
Let $\bfiN = \bfiN_1 \bfiN_2$ and $\beta = \bfiN \ncdot \alpha$.
It follows from (\sclub) and (\sclub\sclub) that
$d(v,F) \approx^{\beta} d(y,F)$.
By {\smaller 5.1}, $d(v,F) = d(v^{\nperp},F)$ and
$d(y,F) = d(y^{\nperp},F)$.
Hence
$$
d(y^{\nperp},F) \approx^{\kern1.0pt \beta} d(v^{\nperp},F).
$$
Clause~{\smaller 3.3 } in Part 3 says that
$\norm{y^{\nperp}} \leq (1 + \varepsilon)d(y,F)$.
In Cases 1 and 2 of Part 3, \,$\varepsilon$ was taken to be
$\dghalf$ and $\dgninth$ respectively. So
$\norm{y^{\nperp}} \leq \frac{3}{2} d(y^{\nperp},F)$.
Hence
$$
\norm{y^{\nperp}} \leq
\hbox{$\frac{3}{2}$} \ncdot \beta(d(v^{\nperp},F) \leq
\hbox{$\frac{3}{2}$} \ncdot \beta(\norm{v^{\nperp}}).
$$
Since $v^{\nperp}$ is a multiple of $y^{\nperp}$,
it follows that
$\norm{v^{\nperp}} \leq \frac{3}{2} d(v^{\nperp},F)$.
So
$$
\norm{v^{\nperp}} \leq
\hbox{$\frac{3}{2}$} \ncdot \beta(d(v^{\nperp},F) \leq
\hbox{$\frac{3}{2}$} \ncdot \beta(d(y^{\nperp},F) \leq
\hbox{$\frac{3}{2}$} \ncdot \beta(\norm{y^{\nperp}}).
$$
Let $\gamma = \dgfrac{3 \beta}{2}$. Hence
\begin{equation}\tag{5.1}
\norm{y^{\nperp}} \approx^{\kern1pt\gamma} \norm{v^{\nperp}}.
\end{equation}
From the fact that $y^{\nperp} = \nu v^{\nperp}$ and (5.1),
it follows that
\begin{equation}\tag{5.2}
\rule{0pt}{0pt}\kern-6mm
\hbox{
If $\nu > 1$, then
$\nu \mcdot \norm{v^{\nperp}} \leq \gamma(\norm{v^{\nperp}})$;
and
if $\nu < 1$, then
$\frac{1}{\nu} \mcdot \norm{y^{\nperp}} \leq \gamma(\norm{y^{\nperp}})$.
}
\rule{0.8cm}{0pt}
\end{equation}

Let $L = \setm{\haty + t y^{\nperp}}{t \in \bbR}$.
So $L$ is the straight line connecting $y$ and $v$.
Recall that $H_2 = \rfs{span}(\sngltn{y^{\nperp}})$.
By {\smaller 5.2}, $H_2 \perp^{1 \half} \whatF$.
So by Proposition~\ref{p-bddlip-bldr-1.9}(f),
$\norm{\ }^{\whatF,H_2} \approx^{2 \half} \norm{\ }$.
By {\smaller 5.1} and {\smaller 5.2}, $\haty \in \whatF$.
So for every $t \in \bbR$,
$\norm{\haty + t y^{\nperp}} \geq
\frac{2}{5} \ncdot (\norm{\haty} + \abs{t} \norm{y^{\nperp}}) \geq
\frac{2}{5} \ncdot \norm{\haty}$.
That is,
\begin{equation}\tag{5.3}
d(L,0) \geq \hbox{$\frac{2}{5}$} \ncdot \norm{\haty}.
\end{equation}

We show that
\begin{equation}\tag{5.4}
\norm{v^{\nperp}} \leq \hbox{$\frac{3}{7}$} \norm{\haty}.
\end{equation}
Let $\haty,\hatz$ be as in Part 4.
Suppose first that $\norm{\haty} \geq \norm{\hatz}$.
In this case $v^{\nperp} = z^{\nperp}$.
By (4.4), $\norm{z^{\nperp}} \leq \frac{3}{10} \norm{z}$.
Since $z = \hatz + z^{\nperp}$,
$\norm{z^{\nperp}} \leq \frac{3}{7} \norm{\hatz}$,
and since $\norm{\haty} \geq \norm{\hatz}$,
$\norm{z^{\nperp}} \leq \frac{3}{7} \norm{\haty}$.
That is,
if $\norm{\haty} \geq \norm{\hatz}$,
then $\norm{v^{\nperp}} \leq \frac{3}{7} \norm{\haty}$.
Next suppose that $\norm{\haty} < \norm{\hatz}$.
In this case
$\haty + v^{\nperp} = v = \lambda z = \lambda(\hatz + z^{\nperp}) =
\haty + \lambda z^{\nperp}$.
That is, $v^{\nperp} = \lambda z^{\nperp}$ and
$\haty = \lambda \hatz$.
Hence
$\frac{\norm{v^{\nperp}}}{\norm{\haty}} =
\frac{\norm{z^{\nperp}}}{\norm{\hatz}}$.
By (4.4),
$\frac{\norm{v^{\nperp}}}{\norm{\haty}} =
\frac{\norm{z^{\nperp}}}{\norm{\hatz}} \leq \frac{3}{7}$.
So, if $\norm{\haty} < \norm{\hatz}$,
then $\norm{v^{\nperp}} \leq \frac{3}{7} \norm{\haty}$.
We conclude that (5.4) holds in both cases.

Since $v = \haty + v^{\nperp}$,
it follows that $\norm{v} \leq \norm{\haty} + \norm{v^{\nperp}}$.
So by (5.4),
$\norm{v} \leq \frac{10}{7} \norm{\haty}$.
Similarly,
$\norm{\haty} \leq \norm{v} + \norm{v^{\nperp}} \leq
\norm{v} + \frac{3}{7} \norm{\haty}$.
So $\frac{4}{7} \norm{\haty} \leq \norm{v}$.
Hence
\begin{equation}\tag{5.5}
\hbox{$\frac{7}{10}$} \norm{v} \leq \norm{\haty} \leq
\hbox{$\frac{7}{4}$} \norm{v}.
\end{equation}
Fact (5.3) and the first inequality in (5.5) imply that
\begin{equation}\tag{5.6}
d(L,0) \geq \hbox{$\frac{14}{50}$} \norm{v}.
\end{equation}
In Case 1 of Part 3 we chose $\varepsilon = \half$ and $\itDelta = 8$.
So by {\smaller 3.3} and {\smaller 3.4},
$\norm{y^{\nperp}} \leq \frac{3}{2} d(y,F) \leq
\frac{3}{2} \cdot \frac{1}{8} \norm{y}$.
That is, $\norm{y^{\nperp}} \leq \frac{3}{16} \norm{y}$.
Since $y = \haty + y^{\nperp}$,
$\norm{\haty} \geq \frac{13}{16} \norm{y}$.
Hence in Case 1, $\norm{y^{\nperp}} \leq \frac{3}{13} \norm{\haty}$.
In Case 2 of Part 3 we follow the same computation with
$\varepsilon = \ninth$ and $\itDelta = \frac{1}{24}$.
We obtain that $\norm{y^{\nperp}} \leq \frac{10}{108} \norm{y}$
and hence
$\norm{y^{\nperp}} \leq \frac{10}{98} \norm{\haty}$.
So in both cases
\begin{equation}\tag{5.7}
\norm{y^{\nperp}} \leq \hbox{$\frac{3}{13}$} \norm{\haty}.
\end{equation}

We shall next define $g_4$.
The required $\bfif_2$ will be either $g_4$ or $g_4\inverse$.
Recall that $\whatF$ and $H_2$ were defined in Part 4,
and that $\nu$ was defined in {\smaller 5.1}.
For $u \in E$ set $u_1 \eqdf (u)_{\whatF}$ and $u_2 \eqdf (u)_{H_2}$.

If $\nu > 1$
let
$$\hbox{
$\barnu = \nu$, $\barv^{\nperp} = v^{\nperp}$,
$\bary^{\nperp} = y^{\nperp}$, $\barv = v$ and $\bary = y$,}
$$
and if $\nu < 1$
let
$$\hbox{
$\barnu = \frac{1}{\nu}$, $\barv^{\nperp} = y^{\nperp}$,
$\bary^{\nperp} = v^{\nperp}$, $\barv = y$ and $\bary = v$.}
$$
So $\barnu > 0$, $\bary^{\nperp} = \barnu \ncdot \barv^{\nperp}$
and by (5.2),
\begin{equation}\tag{5.8}
\barnu \leq
\frac{\gamma(\norm{\barv^{\nperp}})}{\norm{\barv^{\nperp}}}.
\end{equation}
Let $\rho \in H([0,\infty))$ be the piecewise linear function with
breakpoints at
$\dgfrac{\norm{\barv^{\nperp}}}{2}$,
$\norm{\barv^{\nperp}}$ and $2 \barnu \norm{\barv^{\nperp}}$
such that
$\rho \nrestriction
([0,\dgfrac{\norm{\barv^{\nperp}}}{2}] \cup
[2 \barnu \norm{\barv^{\nperp}},\infty)) =
\rfs{Id}$
and $\rho(\norm{\barv^{\nperp}}) = \barnu \norm{\barv^{\nperp}}$.
Define $\eta(s,t)$ to be the function
$$
\eta(s,t) =\kern2pt
\left\{
\renewcommand{\arraystretch}{1.5}
\addtolength{\arraycolsep}{4pt}
\begin{array}{ll}
(1 - \frac{s}{\dgfrac{\norm{\haty}}{5}}) \rho(t) +
\frac{s}{\dgfrac{\norm{\haty}}{5}} t
\rule{12pt}{0pt}&
s \in [0,\dgfrac{\norm{\haty}}{5}],
\\
t
&
s \geq \dgfrac{\norm{\haty}}{5}.
\end{array}
\renewcommand{\arraystretch}{1.0}
\addtolength{\arraycolsep}{-4pt}
\right.
$$
So $\eta = \eta_{(\rho,\dgfrac{\norm{\haty}}{5})}$ as defined in
Proposition \ref{p-bddlip-1.11-1}(b).
Let $\Ehat =
\setm{u \in E}{u_2 \geq 0}$
Define
$$
g_4(u) =\kern2pt
\left\{
\renewcommand{\arraystretch}{1.5}
\addtolength{\arraycolsep}{4pt}
\begin{array}{ll}
u_1 +
\eta(d(u,L),\norm{u_2}) \ncdot
\frac{\barv^{\nperp}}{\norm{\barv^{\nperp}}}
\rule{6pt}{0pt}&
u \in \Ehat,
\\
u
&
u \in E - \Ehat.
\end{array}
\renewcommand{\arraystretch}{1.0}
\addtolength{\arraycolsep}{-4pt}
\right.
$$
If $u_2 = 0$ then $g_4(u) = u$, so $g \nrestriction \whatF = \rfs{Id}$.
and hence $g_4 \in H(E)$.
Note that if $\nu > 1$, then $g_4(v) = y$,
and if $\nu < 1$, then $g_4\inverse(v) = y$.
Next we find $\bfiM_{3,2}, a_{3,2}, b_{3,2}$ independent of
$E,F,\alpha,\bfix_0,\bfiy_0$ such that
$R(v,y,g_4;\bfiM_{3,2} \ncdot \alpha,a_{3,2},b_{3,2},F)$ holds
or $R(y,v,g_4;\bfiM_{3,2} \ncdot \alpha,a_{3,2},b_{3,2},F)$ holds.
\smallskip

{\bf R3: }
Clearly, $g_4 \nrestriction F = \rfs{Id}$
and hence $g_4(F) = g_4\inverse(F) = F$.
\smallskip

{\bf R2: }
We shall next find $\bfiM_{3,2}$ such that $g_4$ is
$\bfiM_{3,2} \mcdot \alpha$-bicontinuous.
The slopes of the linear pieces of $\rho$ are:
$1$,
$\frac{\barnu \norm{\barv^{\nperp}} - \dgfrac{\norm{\barv^{\nperp}}}{2}}
{\norm{\barv^{\nperp}} - \dgfrac{\norm{\barv^{\nperp}}}{2}}$,
$\frac{2 \barnu \norm{\barv^{\nperp}} - \barnu \norm{\barv^{\nperp}}}
{2 \barnu \norm{\barv^{\nperp}} - \norm{\barv^{\nperp}}}$
and~$1$.
That is, the four slopes of $\rho$ are $1$, $2 \barnu - 1$,
$\frac{\barnu}{2 \barnu - 1}$ and $1$.
We apply Proposition \ref{p-bddlip-1.11-1}(a) to $\rho$
taking $a_0$ to~be~$0$, $a_1, a_2, a_3$ to be the breakpoints of $\rho$
and $a_4$ to be $\infty$.
Using the notation of Definition \ref{d-bddly-lip-2.9-1}(b),
the functions $\rho_1$, $\rho_3$ and $\rho_4$
are linear function with slopes $1$, $\frac{\barnu}{2 \barnu - 1}$
and $1$ respectively. So they are $1$-Lipschitz.
Clearly,
$\rho_2(t) = (2 \barnu - 1)t + c,
\ t \in [0,\dgfrac{\norm{\barv^{\nperp}}}{2})$.
By (5.8) and Proposition~\ref{p-bddlip-1.11-1}(c),
$\rho_2$ is $2 \ncdot \gamma$-continuous,
and so $\rho$ is $(4,2 \ncdot \gamma)$-continuous.
By Proposition~\ref{p-bddlip-1.11-1}(a),
\begin{equation}\tag{5.9}
\hbox{$\rho$ is $8 \mcdot \gamma$-continuous.}
\end{equation}

Let $u,w \in \Ehat$.
Then
\vspace{1.5mm}
\newline
(5.10)
\rule{7pt}{0pt}
\renewcommand{\arraystretch}{1.5}
\addtolength{\arraycolsep}{-3pt}
$
\begin{array}[t]{lll}
\norm{g_4(w) - g_4(u)} 
&
\kern6pt\leq\kern6pt
&
\norm{(w - u)_1} \kern4pt+\kern4pt
\abs{\eta(d(w,L),\norm{w_2}) - \eta(d(u,L),\norm{w_2})} \kern6pt+
\\
&
&
\abs{\eta(d(u,L),\norm{w_2}) - \eta(d(u,L),\norm{u_2})}.
\vspace{1.7mm}
\end{array}
$
\renewcommand{\arraystretch}{1.0}
\addtolength{\arraycolsep}{3pt}
\newline
Denote the three summands in the right hand of inequality (5.10) by
$D_1$, $D_2$ and $D_3$.
If $d(w,L),d(u,L) \in [0,\dgfrac{\norm{\haty}}{5})$, then
\vspace{1.5mm}
\newline
\rule{7pt}{0pt}
\renewcommand{\arraystretch}{1.5}
\addtolength{\arraycolsep}{0pt}
$
\begin{array}{lll}
D_2
&
\leq
&
\frac{\abs{d(w,L) - d(u,L)}}{\dgfrac{\norm{\haty}}{5}} \cdot
(\rho(\norm{w_2}) - \norm{w_2}) \leq
\frac{\norm{w - u}}{\dgfrac{\norm{\haty}}{5}} \cdot
(\rho(\norm{w_2}) - \norm{w_2})
\\
&
\leq
&
\frac{\norm{w - u}}{\dgfrac{\norm{\haty}}{5}} \cdot
(\barnu - 1) \norm{\barv^{\nperp}} \leq
\frac{\norm{w - u}}{\dgfrac{\norm{\haty}}{5}} \cdot
\barnu \cdot \norm{\barv^{\nperp}} \eqdf D_2'.
\vspace{1.7mm}
\end{array}
$
\renewcommand{\arraystretch}{1.0}
\addtolength{\arraycolsep}{0pt}
\newline
The above is true for every $u,w \in \Ehat$.
Since $\barnu \mcdot \barv^{\nperp} = v^{\nperp}$
or $\barnu \mcdot \barv^{\nperp} = y^{\nperp}$,
by (5.4) and (5.7),
$\frac{\barnu \kern0.8pt\cdot \norm{\barv^{\nperp}}}{\norm{\haty}} \leq
\frac{3}{7}$.
Hence,
$D_2' \leq \frac{15}{7} \mcdot \norm{w - u}$. That is,
\begin{equation}\tag{5.11}
\abs{\eta(d(w,L),\norm{w_2}) - \eta(d(u,L),\norm{w_2})} \leq
\hbox{$\frac{15}{7}$} \mcdot \norm{w - u}.
\end{equation}
By (5.9) and Proposition \ref{p-bddlip-1.11-1}(b),
$D_3 \leq
8 \mcdot \gamma(\abs{\kern0.8pt\norm{w_2} - \norm{u_2}\kern0.8pt})
\leq 8 \mcdot \gamma(\norm{(w - u)_2}) \eqdf D_3'$,\break
and by the second inequality in (4.9)
and Proposition \ref{p-bddlip-1.11-1}(d),
$D_3' \leq \frac{3}{2} \mcdot 8 \mcdot \gamma(\norm{w - u}$.
Hence
\begin{equation}\tag{5.12}
\abs{\eta(d(u,L),\norm{w_2}) - \eta(d(u,L),\norm{u_2})} \leq
12 \cdot \gamma(\norm{w - u}).
\end{equation}

Note that for every $\baru \in E$, $\baru_1$ of Part 5 is
$\baru_1 + \baru_3$ of Part 4.
So by the first and third inequalities in (4.9),
\begin{equation}\tag{5.13}
\norm{(w - u)_1} \leq \hbox{$7 \half$} \norm{w - u}.
\end{equation}

Substitute into (5.10) inequlities (5.13), (5.11) and (5.12).
We obtain the inequality
$\norm{g_4(w) - g_4(u)} \leq 9 \frac{9}{14} \norm{w - u} +
12 \mcdot \gamma(\norm{w - u})$.
Recall that $\gamma = \frac{3}{2} \beta$ and that
$\beta = \bfiN \alpha$.
Hence, since $\alpha \geq \rfs{Id}$,
$$
g_4 \nrestriction \Ehat \mbox{ is }
(18 \bfiN + 10) \mcdot \alpha\mbox{-continuous.}
$$

The computation which shows that for some $\bfiM$
independent of $E,F,\alpha,\bfix_0,\bfiy_0$,\break
$g_4\inverse \nrestriction \Ehat$
is $\bfiM \mcdot \alpha$-continuous is analogous.
But for $g_4\inverse$ there is $M$ which does not depend on
$E,F,\alpha,\bfix_0,\bfiy_0$
such that $g_4\inverse \nrestriction \Ehat$ is $M$-Lipschitz.
So we shall conclude that $g_4\inverse \nrestriction \Ehat$ is
$M \mcdot \alpha$-continuous.
This computation is analogous to the proof
that $\bfif_1\inverse$ is Lipschitz.

For $s \in [0, \infty)$ let $\theta_s = \eta_s\inverse$.
Denote $\theta(s,t) = \theta_s(t)$.
As in Part 4, for every $u \in E$,
$$
g_4\inverse(u) =
u_1 + \theta_{d(u,S)}(\norm{u_2})
\frac{\barv^{\nperp}}{\norm{\barv^{\nperp}}}.
$$
The analogue (5.10$^*$) of (5.10) and (5.11$^*$) of (5.11)
obtained by replacing $\eta$ by $\theta$ are true.
Let $\mu = \theta_0$. So $\mu = \rho\inverse$
and $\theta = \eta_{(\mu,\dgfrac{\norm{\haty}}{5})}$.
The slopes of the linear pieces of $\mu$ are the inverses of the
slopes of the linear pieces of $\rho$.
Hence the slopes are
$1$, $\frac{1}{2 \barnu - 1}$ and
$\frac{2 \barnu - 1}{\barnu}$.
The first two slopes are $\leq 1$ and the third is $\leq 2$.
So $\mu$ is $2$-Lipschitz.
By Proposition \ref{p-bddlip-1.11-1}(b), for every $s \in [0,\infty)$,
$\theta_s$ is $2$-Lipschitz.
Hence
$$
\abs{\theta(d(u,L),\norm{w_2}) - \theta(d(u,L),\norm{u_2})} \leq
2 \mcdot\abs{\kern0.7pt\norm{w_2} - \norm{u_2}\kern0.5pt} =
2 \mcdot \norm{(w - u)_2}.
$$
Applying the second inequality in (4.9) we conclude that
\begin{equation}\tag{5.14}
\abs{\theta(d(u,L),\norm{w_2}) - \theta(d(u,L),\norm{u_2})} \leq
3 \norm{w - u}.
\end{equation}
Substituting (5.13), (5.11$^*$) and (5.14)
into (5.10$^*$) we conclude that
$$
\norm{g_4\inverse(w) - g_4\inverse(u)} \leq
\hbox{$(7 \half + \frac{15}{7} + 3)$} \norm{w - u} \leq 13 \norm{w - u}.
$$
Since $13 \mcdot \rfs{Id} \leq (18 \bfiN + 10) \mcdot \alpha$,
$g_4\inverse \nrestriction \Ehat$
is $(18 \bfiN + 10) \mcdot \alpha$-continuous.
Hence
$g_4 \nrestriction \Ehat$
is $(18 \bfiN + 10) \mcdot \alpha$-bicontinuous
and so $g_4$ is $2 (18 \bfiN + 10) \mcdot \alpha$-bicontinuous.
So $\bfiM_{3,2} \eqdf 60 \bfiN$ is as required.
That is, $g$ and $g\inverse$
are $\bfiM_{3,2} \mcdot \alpha$-bicontinuous.
\smallskip

{\bf R4: }
We shall find $a'$ and $b'$ independent of
$E,F,\alpha,\bfix_0$ and $\bfiy_0$ such that
$\rfs{supp}(g_4) \subseteq B(0;a' \norm{\haty}, b' \norm{\haty})$.
Let $u \in B(0,\dgfrac{\norm{\haty}}{5})$.
By~(5.3), $d(u,L) > \dgfrac{\norm{\haty}}{5}$.
So for every $t \in [0,\infty)$,
$\eta(d(u,L),t) = t$.
In particular, $\eta(d(u,L),\norm{u_2}) = \norm{u_2}$.
Hence
$$
g_4(u) = u_1 +
\eta(d(u,L),\norm{u_2}) \mcdot
\frac{\barv^{\nperp}}{\norm{\barv^{\nperp}}} = u_1 + u_2 = u.
$$
That is,
$g_4 \nrestriction B(0,\dgfrac{\norm{\haty}}{5}) = \rfs{Id}$
and hence $a' = \dgfrac{1}{5}$.

Let $u \in E$. If $d(u,L) \geq \dgfrac{\norm{\haty}}{5}$
or $\norm{u_2} \geq 2 \barnu \norm{\barv^{\nperp}}$,
then
$\eta(d(u,L),\norm{u_2}) = \norm{u_2}$
and hence $g_4(u) = u$.
Recall that
$\barnu \cdot \barv^{\nperp} = v^{\nperp}$ or
$\barnu \cdot \barv^{\nperp} = y^{\nperp}$.
So if $\norm{u_2} \geq 2 \norm{\barv^{\nperp}}$
and
$\norm{u_2} \geq 2 \norm{\bary^{\nperp}}$,
then $g_4(u) = u$.
So by (5.4) and (5.7),
\begin{equation}\tag{5.15}\kern-6pt
\hbox{If $\norm{u_2} \geq \frac{6}{7} \norm{\haty}$,
then $g_4(u) = u$.}\rule{8.4cm}{0pt}
\end{equation}

Fact $(\star)$ in Part 4 (which precedes (4.9),
says that $\baru_{1,3} \leq 2 \half \baru$
for every $\baru \in E$.
But $\baru_{1,3}$ of Part 4 is $\baru_1$ of Part~5.
So $\norm{\baru_1} \leq \frac{5}{2} \norm {\baru}$
for every $\baru \in E$.
We show that
\begin{equation}\tag{5.16}\kern-7pt
\hbox{If $\norm{u_1} \geq 1 \half \norm{\haty}$,
then $g_4(u) = u$.}
\rule{8.2cm}{0pt}
\end{equation}
Suppose that $\norm{u_1} \geq 1 \half \norm{\haty}$
and let $w \in L$. Then $(u - w)_1 = u_1 - \haty$
and hence
$\norm{(u - w)_1} \geq \norm{u_1} - \norm{\haty} \geq
\half \norm{\haty}$.
So
$\norm{u - w} \geq \frac{2}{5} \norm{(u - w)_1} \geq
\frac{2}{5} \ncdot \half \norm{\haty} = \dgfrac{\norm{\haty}}{5}$.
Hence $d(u,L) \geq \dgfrac{\norm{\haty}}{5}$.
This implies that $g_4(u) = u$.
Suppose that
$\norm{u} \geq 3 \norm{\haty}$ and we show that $g_4(u) = u$.
Clearly,
$\norm{u_1} + \norm{u_2} \geq \norm{u} \geq 3 \norm{\haty}$.
So either $\norm{u_1} \geq 1 \half \norm{\haty}$
or
$\norm{u_2} \geq \frac{6}{7} \norm{\haty}$.
By (5.16) and (5.15), $g_4(u) = u$.
It follows that
$g_4 \nrestriction (E - B(0,3 \norm{\haty})) = \rfs{Id}$.
So $b' \eqdf 3$ is as desired.

Recall that (4.6) said that
$\norm{\haty} \leq \frac{19}{16} \norm{y}$,
and that (5.5) said that
$\frac{7}{10} \norm{v} \leq \norm{\haty}$.
It follows that
$\rfs{supp}(g_4) \subseteq
B(0;\fifth \ncdot \frac{7}{10} \norm{v},
3 \ncdot \frac{19}{16} \norm{y})$.
That is,
$\rfs{supp}(g_4) \subseteq
B(0;\hbox{$\frac{7}{50}$} \norm{v}, \hbox{$\frac{57}{16}$} \norm{y})$,
and the same is true for $g_4\inverse$.
Let
$a_{3,2} = \dgfrac{7}{50}$ and $b_{3,2} = \dgfrac{57}{16}$.
Then
$\rfs{supp}(g_2),\rfs{supp}(g_2\inverse) \subseteq
B(0;a_{3,2} \norm{v}, b_{3,2} \norm{y})$.
So R3 is proved.
\smallskip

{\bf The definition of $\bfif_2$: }
If $\nu > 1$ define $\bfif_2 = g_4$,
and if $\nu < 1$ define $\bfif_2 = g_4\inverse$.
\smallskip

{\bf R1: }
Clearly, $\bfif_2(v) = y$, and since $\bfiv = v$ and $\bfiy_2 = y$,
we have $\bfif_2(\bfiv) = \bfiy_2$.

We have found $\bfiM_{i,j}$'s, $a_{i,j}$'s and $b_{i,j}$'s
which fulfill C1\,-\,C4. It follows from the first part of the proof
of the lemma that there exist $\bfiM,a,b$ such that
$\bfiM$ is a UC-constant for $\pair{a}{b}$.

\kern1.3mm

(b) Let $\bfiM, a, b$ be as assured by Part (a),
and let $a' < 1$ and $b' > 1$.
We may assume that $a' > a$ and that $b' < b$.
Let $x,y \in E - F$ be as in the definition of a UC-constant.
Let $\bfig_1,\bfig_2$ be as assured in Part (a) for the numbers
$a$ and $b$. (See Definition \ref{d-bddly-lip-2.9}(a)).
Let $\eta \in H([0,\infty))$ be a piecewise linear function with
breakpoints at:
$a \cdot \norm{x},\norm{x}, \norm{y},
b \cdot \norm{y}, 2b \cdot \norm{y}$
and such that:
$\eta(0) = 0$;
$\eta(a \cdot \norm{x}) = a' \cdot \norm{x}$;
$\eta(\norm{x}) = \norm{x}$;
$\eta(\norm{y}) = \norm{y}$;
$\eta(b \cdot \norm{y}) = b' \cdot \norm{y}$;
$\eta \nrestriction [2b \cdot \norm{y},\infty) = \rfs{Id}$.
The slopes of the linear pieces of $\eta$ are:
$\frac{a'}{a}$, $\frac{1 - a'}{1 - a}$, $1$, $\frac{b' - 1}{b - 1}$,
$\frac{2b - b'}{b}$ and {1}.
These slopes depend only on $a, a', b, b'$ and not on $x$ and $y$.
Let $\bfiL$ be the maximum of all the above slopes and their inverses.
So $\eta$ is $\bfiL$-bilipschitz.

Let $k$ be the piecewise linearly radial homeomorphism based on $\eta$.
That is, for every $u \in E - \sngltn{0}$,
$k(u) = \eta(\norm{u}) \frac{u}{\norm{u}}$
and $k(0) = 0$.
By Proposition \ref{metr-bldr-p3.18}, $k$ is $3 \bfiL$-blipschitz.
For $i = 1,2$, let $\bfig\fprime_i = k \scirc \bfig_i \scirc k\inverse$.
Then $\bfig\fprime_i$
is
$(3 \bfiL \cdot \rfs{Id}) \scirc \alpha \scirc
(3 \bfiL \cdot \rfs{Id})$-bicontinuous.
So by Proposition~\ref{p-bddlip-1.11-1}(d),
$\bfig\fprime_i$ is
$9 \bfiL^2 \bfiM \kern1pt\cdot \alpha$-bicontinuous.
Define $\bfiM\fprime = 9 \bfiL^2 \bfiM$.
It is easy to verify that
Clauses B1\,-\,B4 in the definition of
a UC-constant (Definition~\ref{d-bddly-lip-2.9}(a)), are fulfilled by
$a'$, $b'$, $\bfig\fprime_1$, $\bfig\fprime_2$ and $\bfiM\fprime$.
Hence $\bfiM\fprime$ is a UC-constant for $\pair{a'}{b'}$.
\hfill\myqed

\newpage

\bigskip
\section{${\bf 1}$-dimensional boundaries}
\label{s10}
Chapter \ref{s9} dealt with the following situation.
$E$ is a normed space, $F$ is a closed subspace of $E$ with dimension
$\geq 2$, $x,y \in E - F$,
and $\norm{x} \approx^{\alpha} \norm{y}$
and $d(x,F) \approx^{\alpha} d(y,F)$.
It was shown that there is an
$M \ncdot \alpha \scirc \alpha$-bicontinuous
$g \in H(E)$ such that $g(x) = y$,
$g(F) = F$ and $\rfs{supp}(g)$ is contained the ring
$B(0;a \norm{x},b \norm{y})$.
When $F$ is $1$-dimensional, such a $g$ does not always exist.
The reason for this is that in order to move $x$ to $y$ we need to
rotate $x$ about an axis perpendicular to $F$.
See the construction of $g_1$ in Part 2 of the proof of
Lemma \ref{l-bddlip-1.10}(a).
When $F$ is $1$-dimensional, such a rotation does not exist.

Whereas Part 2 of the proof of Lemma \ref{l-bddlip-1.10}(a) fails
for a $1$-dimen\-sional subspace, Parts~1 and 3\,-\,5 remain without
change.
In these parts, the fact that $\rfs{dim}(F) \geq 2$ is not used.
By skipping Part~2 in the proof of Lemma \ref{l-bddlip-1.10}(a)
one obtains the following lemma.

Let $F,K$ be linear subspaces of a normed space $E$ and $u \in E$.
Then $u \perp F$ denotes the fact that $\norm{u} = d(u,F)$,
and $K \perp F$ means that $u \perp F$ for every $u \in K$.
   \index{N@AAAA@@ $u \perp F$. 
          This means: $F$ is linear subspace of a normed space $E$,
          $u \in E$ and $\norm{u} = d(u,F)$}
   \index{N@AAAA@@ $K \perp F$. 
          This means: $K,F$ are linear subspaces of $E$
          and for every $u \in K$, $u \perp F$}

\begin{lemma}\label{l3.1-temp}
Let $M$ be a UC-constant for $\pair{a}{b}$.
Then the following holds.
\newline
Let $E$ be a normed space and $F$ be a $1$-dimensional
linear subspace of $E$.
Let $\alpha \in \rfs{MBC}$ and
$x,y \in E - F$ be such that:
\begin{itemize}
\addtolength{\parskip}{-11pt}
\addtolength{\itemsep}{06pt}
\item[\num{i}]
$\norm{x} \leq \norm{y}$ and $\norm{x} \approx^{\alpha} \norm{y}$,
\item[\num{ii}] $d(x,F) \approx^{\alpha} d(y,F)$,
\item[\num{iii}] $x = \hatx + x^{\sperp}$, $y = \haty + y^{\sperp}$,
$\hatx, \haty \in F$, \ 
$x^{\sperp} \kern-1pt,\kern-1pt y^{\sperp} \perp F$,
and for some $\lambda > 0$,
$\hatx = \lambda \haty$,
\item[\num{iv}]
if $\rfs{dim}(E) = 2$, then $x,y$ are on the same side of $F$.
\vspace{-05.7pt}
\end{itemize}
Then there are $g_1,g_2 \in H(E)$ such that:
\begin{itemize}
\addtolength{\parskip}{-11pt}
\addtolength{\itemsep}{06pt}
\item[\num{1}] $g_1,g_2$ are $M \alpha$-bicontinuous,
\item[\num{2}] $g_1 \scirc g_2(x) = y$,
\item[\num{3}] $g_1(F) = F$ and $g_2(F) = F$,
\item[\num{4}] For every $i = 1,2$, \,
$\rfs{supp}(g_i) \subseteq B(0;a \norm{x},b \norm{y})$.
\vspace{-05.7pt}
\end{itemize}
\end{lemma}

\noindent
{\bf Proof } Parts 1, 3\,-\,5 of the proof of
Lemma \ref{l-bddlip-1.10}(a) constitue a proof of this Lemma.
\hfill\myqed

\begin{defn}\label{d3.2-temp}
\begin{rm}
Let $0 < a < 1$, $b > 1$ and $M \geq 1$.
We say that {\it ``$M$ is a $1$-dimensional
Uniform Continuity constant for $a$ and $b$''}
   \index{uc-constant@@UC-constant.
          $M$ is a 1UC-constant for $\pair{a}{b}$}
(abbreviated by ``$M$ is a {\it 1UC-constant} for $\pair{a}{b}$''),
if the following holds.
\newline
Suppose that $E,F,\alpha,x,y$ satisfy the following assumptions.
\begin{itemize}
\addtolength{\parskip}{-11pt}
\addtolength{\itemsep}{06pt}
\item[A1]
$E$ is a normed space and $F$ is a proper linear subspace of $E$
such that $\rfs{dim}(F) = 1$,
$\alpha \in \rfs{MBC}$ and
$x,y \in E - F$,
\item[A2]
$\norm{x} \leq \norm{y} \leq \alpha(\norm{x})$,
\item[A3]
$d(x,F) \approx^{\alpha} d(y,F)$,
\item[A4]
$\norm{x} \leq \alpha(d(x,F))$ and
$\norm{y} \leq \alpha(d(y,F))$,
\item[A5]
if $\rfs{dim}(E) = 2$, then $x,y$ are on the same side of $F$.
\vspace{-05.7pt}
\end{itemize}
Then there are $g_1,g_2,g_3 \in H(E)$ such that:
\begin{itemize}
\addtolength{\parskip}{-11pt}
\addtolength{\itemsep}{06pt}
\item[B1]
For every $i = 1,2,3$, \,
$g_i$ is $M \ncdot \alpha$-bicontinuous,
\item[B2]
$g_3 \scirc g_2 \scirc g_1(x) = y$,
\item[B3]
For every $i = 1,2,3$, \,
$g_i(F) =  F$,
\item[B4]
For every $i = 1,2,3$, \,
$\rfs{supp}(g_i) \subseteq B(0;a \norm{x},b \norm{y})$.
\vspace{-05.7pt}
\end{itemize}
\end{rm}
\end{defn}

{\bf Remark } Note that in the definition of a 1UC-constant
there is an extra assumption on $x$ and $y$ which did not appear in
the definition of a UC-constant. Namely,
Assumption A4 which says that
$\norm{x} \leq \alpha(d(x,F))$ and $\norm{y} \leq \alpha(d(y,F))$.
\smallskip

The rest of the chapter is devoted to the proof of the following lemma.

\begin{lemma}\label{l3.3-temp}
\num{a} There are $a,b,M$ such that
$M$ is a 1UC-constant for $a$ and~$b$.

\num{b} For every $0 < a < 1$ and $b > 1$ there is $M$ such that
$M$ is a 1UC-constant for $a$ and~$b$.
\end{lemma}

Items \ref{p3.4-temp}\,-\,\ref{p3.9-temp}
are needed in the proof of the above lemma.

\begin{prop}\label{p3.4-temp}
Let $F$ be a finite-dimensional linear subspace
of a normed space $E$ and $u \not\in F$.
Then there is a $1$-dimensional subspace
$L \subseteq \rfs{span}(F \cup \sngltn{u})$ such that
$L \perp F$.
\end{prop}

\begin{prop}\label{p3.5-temp}
Let $X$ be a metric space, $\alpha \in \rfs{MBC}$, $c > 0$,
$D,K \geq 1$, $g \in H(X)$, $\rfs{diam}(\rfs{supp}(g)) \leq D \alpha(c)$
and $g$ is $K \cdot \frac{\alpha(c)}{c}$-Lipschitz.
\newline
Then $g$ is $(D + K + 1) \cdot \alpha$-continuous.
\end{prop}

\noindent
{\bf Proof }
Note that if $\alpha \in \rfs{MC}$, then the function
$\frac{\alpha(t)}{t}$ is a decreasing function.
Let $x,y \in X$.
Suppose first that $d(x,y) \leq c$.
Then
\vspace{1.5mm}
\newline
\rule{7pt}{0pt}
\renewcommand{\arraystretch}{1.5}
\addtolength{\arraycolsep}{-3pt}
$
\begin{array}{ll}
&
d(g(x),g(y)) \leq K \frac{\alpha(c)}{c} d(x,y) \leq
K \frac{\alpha(d(x,y))}{d(x,y)} d(x,y) = K \alpha(d(x,y))
\\
\leq
\rule{5pt}{0pt}
&
(D + K + 1) \cdot \alpha(d(x,y)).
\vspace{1.7mm}
\end{array}
$
\renewcommand{\arraystretch}{1.0}
\addtolength{\arraycolsep}{3pt}
\newline
Next assume that $d(x,y) > c$.
If $x,y \in \rfs{supp}(g)$,
then
\vspace{1.5mm}
\newline
\rule{7pt}{0pt}
\renewcommand{\arraystretch}{1.5}
\addtolength{\arraycolsep}{-3pt}
$
\begin{array}{ll}
\phantom{\hbox{$\ \leq\ $}}
&
d(g(x),g(y)) \leq D \alpha(c) < D \alpha(d(x,y)) \leq
(D + K + 1) \cdot \alpha(d(x,y)).
\vspace{1.7mm}
\end{array}
$
\renewcommand{\arraystretch}{1.0}
\addtolength{\arraycolsep}{3pt}
\newline
If $x \not\in \rfs{supp}(g)$ and $y \in \rfs{supp}(g)$, then
\vspace{1.5mm}
\newline
\rule{7pt}{0pt}
\renewcommand{\arraystretch}{1.5}
\addtolength{\arraycolsep}{-3pt}
$
\begin{array}{ll}
&
d(g(x),g(y)) \leq d(x,y) + d(y,g(y)) \leq
\alpha(d(x,y)) + D \alpha(c)
\\
<
\rule{5pt}{0pt}
&
\rule{0mm}{14pt}
\alpha(d(x,y)) + D \alpha(d(x,y)) =
(D + K + 1) \alpha(d(x,y)).
\vspace{1.7mm}
\end{array}
$
\renewcommand{\arraystretch}{1.0}
\addtolength{\arraycolsep}{3pt}
\newline
The case that $x \in \rfs{supp}(g)$ and $y \not\in \rfs{supp}(g)$
is identical,
and the case that $x,y \not\in \rfs{supp}(g)$ is trivial.
\hfill\myqed

\begin{prop}\label{p3.6-temp}
There are $b > 1$, $0 < a < 1$ and $M > 1$
such that the following holds.
Suppose that:
\begin{itemize}
\addtolength{\parskip}{-11pt}
\addtolength{\itemsep}{06pt}
\item[\num{1}] $\alpha \in \rfs{MBC}$,
\item[\num{2}] $E$ is a normed space,
and $L$ is a $1$-dimensional linear subspace of $E$,
\item[\num{3}] $u \in E - L$
and $\norm{u} \leq \alpha(d(u,L))$,
\item[\num{4}] Let $u = \hatu + u^{\sperp}$, where $\hatu \in L$ and
$u^{\sperp} \perp L$,
and let $v = \frac{\norm{u}}{\norm{u^{\ssperp}}} u^{\sperp}$.
\vspace{-05.7pt}
\end{itemize}
Then there is $g \in H(E)$ such that:
\begin{itemize}
\addtolength{\parskip}{-11pt}
\addtolength{\itemsep}{06pt}
\item[\num{1}] $g(u) = v$,
\item[\num{2}] $g$ is $M \ncdot \alpha$-bicontinuous,
\item[\num{3}] $\rfs{supp}(g) \subseteq B(0;a \norm{u},b \norm{u})$,
\item[\num{4}] $g(L) = L$.
\vspace{-05.7pt}
\end{itemize}
Define $M^{\srfs{lift}} = M$, $a^{\srfs{lift}} = a$ and
$b^{\srfs{lift}} = b$.
Note that the conjunction of Clauses \num{1}\,-\,\num{4} is
the relation $R(u,v,g;M \ncdot \alpha,a,b,L)$
defined in Definition~\ref{d-bddly-lip-2.9}(b)).
   \index{N@mlit@@$M^{\srfs{lift}}$}
\end{prop}

\noindent
{\bf Proof } 
Let $A = [u,v]$.
Clearly, $d(u,L) = \norm{u^{\sperp}}$.
So $\norm{u^{\sperp}} \leq \norm{u}$.
We find an upper bound for $\norm{u - v}$.
\vspace{1.5mm}
\newline
\rule{7pt}{0pt}
\renewcommand{\arraystretch}{1.5}
\addtolength{\arraycolsep}{-3pt}
$
\begin{array}{ll}
&
\norm{u - v} \leq \norm{u - u^{\sperp}} + \norm{u^{\sperp} - v} =
\norm{\hatu} + (\norm{u} - d(u,L))
\\
\leq
\rule{5pt}{0pt}
&
(\norm{u} + \norm{u^{\sperp}}) + (\norm{u} - d(u,L))
\\
=
\rule{5pt}{0pt}
&
(\norm{u} + d(u,L)) + (\norm{u} - d(u,L)) =
2 \norm{u} \leq 2 \alpha(d(u,L)).
\vspace{1.7mm}
\end{array}
$
\renewcommand{\arraystretch}{1.0}
\addtolength{\arraycolsep}{3pt}
\newline
We show that $d(A,L) = d(u,L)$. For every $z \in A$
there are $\lambda \in [0,1]$ and $\mu \geq 1$ such that
$z = \lambda \hatu + \mu u^{\sperp}$.
So $d(z,L) = \mu \norm{u^{\sperp}} \geq \norm{u^{\sperp}} = d(u,L)$.
Since $u \in A$, we have that $d(A,L) = d(u,L)$.
We show that $d(A,0) \geq \dgfrac{\norm{u}}{4}$.
Let $w = \hatu + v$ and $C = [u,w] \cup [w,v]$.
We first show that $d(C,0) \geq \dgfrac{\norm{u}}{2}$.
If $z \in [v,w]$, then for some $t \in \bbR$,
$z = v + t \hatu$.
So
$\norm{z} \geq d(z,L) = d(v,L) = \norm{v} = \norm{u}$.
Hence $d([v,w],0) = \norm{u}$.

Note that
$[u,w] =
\setm{u +tv}{0 \leq t \leq 1 - \frac{\norm{u^{\ssperp}}}{\norm{v}}}$.
Let $z = u + tv \in [u,w]$.
If $t \leq \dghalf$, then
$\norm{u + tv} \geq \norm{u} - t \norm{v} \geq
\norm{u} - \dgfrac{\norm{u}}{2} = \dgfrac{\norm{u}}{2}$.
If $t \geq \dghalf$, then
\vspace{1.5mm}
\newline
\rule{7pt}{0pt}
\renewcommand{\arraystretch}{1.5}
\addtolength{\arraycolsep}{-3pt}
$
\begin{array}{ll}
&
\norm{u + tv} \geq d(u + tv,L) =
d(\hatu + u^{\sperp} + tv,L) = d(u^{\sperp} + tv,L)
\\
=
\rule{5pt}{0pt}
&
d((t + \frac{\norm{u^{\ssperp}}}{\norm{v}}) v,L) \geq
d(tv,L) \geq
\dgfrac{d(v,L)}{2} = \dgfrac{\norm{u}}{2}.
\vspace{1.7mm}
\end{array}
$
\renewcommand{\arraystretch}{1.0}
\addtolength{\arraycolsep}{3pt}
\newline
Hence $d([u,w],0) \geq \dgfrac{\norm{u}}{2}$.
It follows that $d(C,0) \geq \dgfrac{\norm{u}}{2}$.

We next prove that
$(*)$ for every $x \in A$ there are $z \in C$ and
$\mu \in [\dghalf,1]$ such that $x = \mu z$.
Recall that $w = \hatu + v$. The equation
$\mu w = \lambda u + (1 - \lambda) v$ has the solution
$\mu = \lambda = \frac{\norm{u}}{2 \norm{u} - \norm{u^{\ssperp}}}$.
So $\mu \in (0,1)$.
That is, there are $x \in A$, $z \in C$ and $\mu \in (0,1)$ such that
$x = \mu z$,
and hence
for every $x \in A$ there are $z \in C$ and $\mu \in (0,1)$ such that
$x = \mu z$.

Let $z \in [u,w]$. Then
$z = u + t \frac{\norm{u}}{\norm{u^{\ssperp}}} u^{\sperp}$, where
$0 \leq t \leq 1 - \frac{\norm{u^{\ssperp}}}{\norm{u}}$.
The equation\break
$\mu z = \lambda u + (1 - \lambda) v$ has the solution
$\lambda = \mu = \frac{1}{1 + t}$.
Since $t \in (0,1)$, $\mu \in [\half,1]$.
Let $z \in [v,w]$. Then
$z =
\frac{\norm{u}}{\norm{u^{\ssperp}}} u^{\sperp} + t(u - u^{\sperp})$,
where $0 \leq t \leq 1$.
The equation
$\mu z = \lambda u + (1 - \lambda) v$ has the solution
$\mu = \frac{\norm{u}}{\norm{u} + t (\norm{u} - \norm{u^{\ssperp}})}$,
$\lambda = t \mu$.
It follows that $\mu \in (\dghalf,1]$.
So $(*)$ is proved.
Hence $d(A,0) \geq \dgfrac{d(C,0)}{2} \geq \dgfrac{\norm{u}}{4}$.

Let $r = \dgfrac{d(u,L)}{8}$.
By Proposition \ref{p-gamma.7}(a),
there is $g \in H(E)$ such that
$g(u) = v$,
$\rfs{supp}(g) \subseteq B(A,r)$
and $g$ is
$M^{\srfs{seg}} \cdot (\frac{\srfs{lngth}(A)}{r} + 1)$-bilipschitz.
Hence Requirement (1) holds.
$$
M^{\srfs{seg}} \cdot (\frac{\rfs{lngth}(A)}{r} + 1) \leq
M^{\srfs{seg}} \cdot (\frac{16\, \alpha(d(u,L))}{d(u,L)} + 1) \leq
M^{\srfs{seg}} \cdot \frac{17\, \alpha(d(u,L))}{d(u,L)}.
$$
So
$g$ is
$17 M^{\srfs{seg}} \cdot \frac{\alpha(d(u,L))}{d(u,L)}$-bilipschitz.
Also,
$$
\rfs{diam}(B(A,r)) \leq \rfs{lngth}(A) + 2r \leq
2 \alpha(d(u,L)) + \dgfrac{d(u,L)}{4} \leq 3 \alpha(d(u,L)).
$$
We apply Proposition \ref{p3.5-temp} to $g$ and to $g\inverse$
with $c = d(u,L)$, $D = 3$ and $K = 17 M^{\srfs{seg}}$.
It follows that $g$ is
$(4 + 17 M^{\srfs{seg}}) \cdot \alpha$-bicontinuous.
So Requirement (2) holds
with $M = 4 + 17 M^{\srfs{seg}}$.
Since $d(A,L) = d(u,L)$ and $r < d(u,L)$,
it follows that $d(B(A,r),L) > 0$.
So $g \nrestriction L = \rfs{Id}$. Requirement~(4) thus holds.

We find the $a$ and $b$ of Requirement (3).
Let $r_0 = d(B(A,r),0)$.
So $g \nrestriction B(0,r_0) = \rfs{Id}$.
But
$r_0 = d(A,0) - r \geq
\dgfrac{\norm{u}}{4} - \dgfrac{d(u,L)}{8} \geq  \dgfrac{\norm{u}}{8}$.
So $a = \dgeighth$.
Let $r_1 = \sup_{x \in B(A,r)} \norm{x}$.
Then $\rfs{supp}(g) \subseteq B(0,r_1)$.
For every $x \in A$, $\norm{x} \leq \max(\norm{u},\norm{v}) = \norm{u}$.
So $r_1 \leq \norm{u} + r < 2 \norm{u}$.
Define $b = 2$.
Then $\rfs{supp}(g) \subseteq B(0;a \norm{u},b \norm{u})$.
So Requirement (3) is fulfilled with $a = \dgeighth$ and $b = 2$.
\rule{19pt}{0pt}\hfill\myqed

\begin{prop}\label{p3.7-temp}
Let $E$ be a $3$-dimensional Hilbert space,
$L$ be a $1$-dimensional sybspace of $E$,
$u,v \in E - L$ and $M \geq 1$.
Suppose that $\norm{u}, \norm{v} \leq M d(u,L)$
and $\norm{u}, \norm{v} \leq M d(v,L)$.
Then there is a rectifiable arc $A$ connecting $u$ and $v$ such that:
\begin{itemize}
\addtolength{\parskip}{-11pt}
\addtolength{\itemsep}{06pt}
\item[\num{1}] $\rfs{lngth}(A) \leq (4 + \pi)M \norm{u}$,
\item[\num{2}] $d(A,L) \geq \dgfrac{\norm{u}}{M}$,
\item[\num{3}] $\max(\setm{\norm{x}}{x \in A}) \leq M \norm{u}$.
\vspace{-05.7pt}
\end{itemize}
\end{prop}

\noindent
{\bf Proof }
Let $w_1 = u^{\perp}$, $w_2 = v^{\perp}$
and $w_3 = \frac{\norm{u^{\perp}}}{\norm{v^{\perp}}} v^{\perp}$.
Let $S$ be a subarc of $S(0, \norm{w_1}) \cap L^{\sperp}$
whose endpoints are $w_1$ and $w_3$ and such that
$\rfs{lngth}(S) \leq \pi \norm{w_1}$.
Define $A = [u,w_1] \cup S \cup [w_3,w_2] \cup [w_2,v]$.
Then $d(A,L) = \min(d(u,L),d(v,L)) \geq \dgfrac{\norm{u}}{M}$.
It is obvious that
$\max(\setm{\norm{x}}{x \in A}) = \max(\norm{u},\norm{v}) \leq
M \norm{u}$.
Now,
\vspace{1.5mm}
\newline
\rule{7pt}{0pt}
\renewcommand{\arraystretch}{1.5}
\addtolength{\arraycolsep}{0pt}
$
\begin{array}{lll}
\rfs{lngth}(A)
&
\leq
&
\norm{(u)_L} + \pi \norm{u^{\perp}} +
\abs{\norm{u^{\perp}} - \norm{v^{\perp}}} +
\norm{(v)_L}
\\
&
\leq
\rule{-3pt}{0pt}
&
\norm{u} + \pi \norm{u} + \norm{u} + \norm{v}
\leq (4 + \pi)M \norm{u}.
\vspace{1.7mm}
\end{array}
$
\renewcommand{\arraystretch}{1.0}
\addtolength{\arraycolsep}{0pt}
\newline
So $A$ is as required.
\hfill\myqed

\begin{prop}\label{p3.8-temp}
There are $M > 1$, $0 < a < 1$ and $b > 1$ such that the
following holds.
Suppose that:
\begin{itemize}
\addtolength{\parskip}{-11pt}
\addtolength{\itemsep}{06pt}
\item[\num{1}] $E$ is a normed space,
and $L$ is a $1$-dimensional linear subspace of $E$,
\item[\num{2}] $u,v \in E - L$, $\norm{u} = \norm{v}$, $u \perp L$
and $v \perp L$,
\item[\num{3}] If $E$ is $2$-dimensional,
then $u,v$ are on the same side of $ L$.
\vspace{-05.7pt}
\end{itemize}
Then there is $g \in H(E)$ such that
$R(u,v,g;M,a,b,L)$ holds.
(See Definition~\ref{d-bddly-lip-2.9}(b)).

We denote $M^{\srfs{perp}} = M$, $a^{\srfs{perp}} = a$ and
$b^{\srfs{perp}} = b$.
\end{prop}

\noindent
{\bf Proof }
If $E$ is $2$-dimensional, then $[u,v] \subseteq S(0,\norm{u})$.
So $d([u,v],L) = \norm{u}$ and\break
$\rfs{lngth}([u,v]) \leq \norm{u} + \norm{v} = 2 \norm{u}$.
By Proposition \ref{p-gamma.7}(a),
there is $g \in H(E)$ such that:
$g(u) = v$,
$\rfs{supp}(g) \subseteq B([u,v],\dgfrac{\norm{u}}{2})$,
and
$g$ is
$M^{\srfs{seg}} \cdot
\frac{2 \norm{u}}{\dgfrac{\norm{u}}{2}}$-bilipschitz.
So for $2$-dimensional $E$'s,
$M,a,b$ can be taken to be $4 M^{\srfs{seg}}$,
$\dghalf$ and $\dgfrac{3}{2}$.

Suppose that $\rfs{dim}(E) > 2$.
Let $F$ be a $3$-dimensional linear subspace of $E$ containing
$L$, $u$ and $v$,
and let $\norm{\ }^{\srbfs{H}}$ be a tight Hilbert norm on $F$.
Denote $N = M^{\srfs{thn}}(3)$.
(See Proposit\-ion~\ref{p-bddlip-bldr-1.9}(b)).
So for every $x \in F$,
$\norm{x} \leq \norm{x}^{\srbfs{H}} \leq N \norm{x}$.
Obviously,
$\norm{u}^{\srbfs{H}}, \norm{v}^{\srbfs{H}} \leq N d^{\srbfs{H}}(u,L)$,
and
$\norm{u}^{\srbfs{H}}, \norm{v}^{\srbfs{H}} \leq N d^{\srbfs{H}}(v,L)$.
By Proposition \ref{p3.7-temp},
there is a rectifiable arc $A$ in $F$ connecting $u$ and $v$ such that:
$\rfs{lngth}^{\srbfs{H}}(A) \leq
(4 + \pi)N \norm{u}^{\srbfs{H}}$,
$d^{\srbfs{H}}(A,L) \geq
\frac{1}{N} \norm{u}^{\srbfs{H}}$
and $\max(\setm{\norm{x}^{\srbfs{H}}}{x \in A}) \leq
N \norm{u}^{\srbfs{H}}$.
So
\begin{itemize}
\addtolength{\parskip}{-11pt}
\addtolength{\itemsep}{06pt}
\item[\num{1}] $\rfs{lngth}(A) \leq (4 + \pi)N^2 \norm{u}$,
\item[\num{2}] $d(A,L) \geq \frac{1}{N^2} \norm{u}$,
\item[\num{3}] $\max(\setm{\norm{x}}{x \in A}) \leq N^2 \norm{u}$.
\vspace{-05.7pt}
\end{itemize}
Let $r = \frac{1}{2 N^2} \norm{u}$,
By Proposition \ref{p-gamma.7}(b),
there is $g \in H(E)$ such that:
\begin{itemize}
\addtolength{\parskip}{-11pt}
\addtolength{\itemsep}{06pt}
\item[\num{4}] $\rfs{supp}(g) \subseteq B(A,r)$,
\item[\num{5}] $g(u) = v$,
\item[\num{6}] $g$ is
$M^{\srfs{arc}}(
\frac{(4 + \pi)N^2 \norm{u}}{\dgfrac{\norm{u}}{(2 N^2)}})$-bilipschitz.
\smallskip
\vspace{-05.7pt}
\end{itemize}
So $g$ is $M^{\srfs{arc}}(16 N^4)$-bilipschitz.

Since $d(B(A,r),L) \geq \frac{1}{2 N^2} \norm{u}$,
$g \nrestriction L = \rfs{Id}$.
Define
$M = M^{\srfs{arc}}(16 N^4)$,
$a = \frac{1}{2 N^2}$ and $b = N^2 + 1$.
Then $M,a,b$ are as required in the proposition.
\rule{10pt}{0pt}\hfill\myqed

\begin{prop}\label{p3.9-temp}
There are $M > 1$, $0 < a < 1$ and $b > 1$ such that the
following holds.
Suppose that
$E$ is a normed space, $u \in E - \sngltn{0}$,
$\alpha \in \rfs{MBC}$,
$1 \leq \lambda \leq \frac{\alpha(\norm{u})}{\norm{u}}$
and $v = \lambda u$.
Then there is a radial homeomorphism $g \in H(E)$ such that:
\hbox{$g(u) = v$,}\break
$g$ is $M \ncdot \alpha$-bicontinuous
and $\rfs{supp}(g) \subseteq B(0;a \norm{u},b \norm{v})$.
Note that this implies that\break
$R(u,v,g;M \ncdot \alpha,a,b,L)$ holds.
Denote $M,a,b$  by $M^{\srfs{dlt}},a^{\srfs{dlt}}$ and
$b^{\srfs{dlt}}$.
\end{prop}

\noindent
{\bf Proof } Let $\eta \in H([0,\infty))$ be the 
piecewise linear function which is determined by the following
equalities.
$\eta(0) = 0$,
$\eta(\dgfrac{\norm{u}}{2}) = \dgfrac{\norm{u}}{2}$,
$\eta(\norm{u}) = \lambda \norm{u}$,
and for every $t \geq \lambda \norm{u} + \norm{u}$, \ $\eta(t) = t$.
The slopes of the linear parts of $\eta$ are $1$,
$2 \lambda$ and $\dgfrac{1}{\lambda}$.
Since $1 \leq \lambda \leq \frac{\alpha(\norm{u})}{\norm{u}}$,
$\eta$ is\break
$2 \mcdot \frac{\alpha(\norm{u})}{\norm{u}}$-bilipschitz.
Let $g$ be the radial homeomorphism of $E$ based on $\eta$.
By Proposition~\ref{metr-bldr-p3.18},
$g$ is $6 \ncdot \frac{\alpha(\norm{u})}{\norm{u}}$-bilipschitz.
Also, $\lambda \norm{u} + \norm{u} \leq 2 \norm{v}$,
hence $\rfs{supp}(g) \subseteq B(0, 2 \norm{v})$.
By Propo\-sition~\ref{p3.5-temp},
$g$ is $(6 + 2 + 1) \ncdot \alpha$-bicontinuous.
So we may define $M = 9$,
$a = \dghalf$ and $b = 2$.
\hfill\myqed

\kern1.3mm

\noindent
{\bf Proof of Lemma \ref{l3.3-temp}}
(a) Let $E,F,x,y$ be as in the definition of a 1UC-constant
(Definition \ref{d3.2-temp}).
There are $\hatx$ and $x^{\sperp}$ such that $x = \hatx + x^{\sperp}$,
$\hatx \in F$ and $x^{\sperp} \perp F$.
Similarly, there are $\haty$ and $y^{\sperp}$ such that
$y = \haty + y^{\sperp}$, $\haty \in F$ and $y^{\sperp} \perp F$.
Let $x_1 = \frac{\norm{x}}{\norm{x^{\ssperp}}} x^{\sperp}$
and $y_1 = \frac{\norm{y}}{\norm{y^{\ssperp}}} y^{\sperp}$.
By Proposition~\ref{p3.6-temp}, there are $f_1,h_1 \in H(E)$
such that
$$
R(x,x_1,f_1;M^{\srfs{lift}} \kern-2pt\mcdot \alpha,a^{\srfs{lift}},
b^{\srfs{lift}},F)
\mbox{\ \ and\ \ \ }%
R(y,y_1,h_1;M^{\srfs{lift}} \kern-2pt\mcdot \alpha,a^{\srfs{lift}},
b^{\srfs{lift}},F).
$$
Let $y_2 = \frac{\norm{x_1}}{\norm{y_1}} y_1$.
Note that $\norm{x_1} = \norm{y_2}$,
$x_1 \perp F$ and $y_2 \perp F$, and if $E$ is $2$-dimensional then
$x_1$, $y_2$ are on the same side of $F$.
By Proposition \ref{p3.8-temp}, there is $f_2 \in H(E)$
such that
$$
R(x_1,y_2,f_2;
M^{\srfs{perp}}, a^{\srfs{perp}}, b^{\srfs{perp}},F).
$$
Since $\norm{y_2} = \norm{x}$ and $\norm{y_1} = \norm{y}$,
it follows that $\norm{y_2} \leq \norm{y_1} \leq \alpha(\norm{y_2})$.
So by Proposition~\ref{p3.9-temp},
there is $g_2 \in H(E)$ such that
$$
R(y_2,y_1,g_2;
M^{\srfs{dlt}} \kern-2pt\mcdot \alpha,
a^{\srfs{dlt}}, b^{\srfs{dlt}},F).
$$
Let $g_1 = f_2 \scirc f_1$ and $g_3 = h_1\inverse$.
Clearly, $g_3 \scirc g_2 \scirc g_1(x) = y$.
Let
$M = \max(M^{\srfs{lift}} M^{\srfs{perp}}, M^{\srfs{dlt}})$.
Note that $\norm{x_1} = \norm{x}$,
so $f_2 \nrestriction B(0,a^{\srfs{perp}} \norm{x}) = \rfs{Id}$.
Set
$a = \min(a^{\srfs{lift}}, a^{\srfs{perp}}, a^{\srfs{dlt}})$
and
$b = \max(b^{\srfs{lift}}, b^{\srfs{perp}}, b^{\srfs{dlt}})$.
It is obvious that Clauses B1\,-\,B4 in the definition of
a 1UC-constant hold for $M,a,b,x,y,g_1,g_2,g_3$ and $F$.

(b) Part (b) is deduced from Part (a) in the same way that
Part (b) of Lemma \ref{l-bddlip-1.10} is deduced from
Part (a) of that lemma.
\hfill\myqed

\newpage

\section{Extending the inducing homeomorphism to the boundary}
\label{s11}

A {\it sequence} means a function whose domain is an infinite subset
of $\bbN$.
If $\sigma \subseteq \bbN$ is infinite, then $\setm{x_i}{i \in \sigma}$
is abbreviated by $\vecx^{\srsp{\sigma}}$.
   \index{N@AAAA@@$\vecx^{\srsp{\sigma}}$. A sequence whose domain is
          $\sigma \subseteq \bbN$}
Suppose that $\vecx^{\srsp{\sigma}}$ is a sequence in $X$
and $g \in H(X,Y)$. Then $g(\vecx^{\srsp{\sigma}})$ denotes
the sequence $\setm{g(x_i)}{i \in \sigma}$.
For $n \in \bbN$ and an infinite $\sigma \subseteq \bbN$
let $\sigma^{\geq n} \eqdf \setm{k \in \sigma}{k \geq n}$.
For a sequence $\vecx$ let
$\vecx^{\kern1pt\geq n} \eqdf
\vecx \nrestriction \rfs{Dom}(\vecx)^{\kern1pt\geq n}$.
\kern-6pt\index{N@AAAA@@$\sigma^{\geq n} =
                \setm{k \in \sigma}{k \geq n}$}
\kern-6pt\index{N@AAAA@@
$\vecx^{\kern1pt\geq n} =
\vecx \nrestriction \rfs{Dom}(\vecx)^{\geq n}$}

Recall that if $\fnn{\alpha}{A}{A}$,
then $\alpha^{\sscirc n}$ denotes $\alpha \scirc \ldots \scirc \alpha$,
\,$n$ times.
Let $X,Y$ be open sets in metric spaces $E$ and $F$ respectively and
$\fnn{g}{X}{Y}$. If $x \in \rfs{Dom}(g^{\srfs{cl}})$,
then we sometimes abbreviate $g^{\srfs{cl}}(x)$ by $g(x)$.

\begin{defn}\label{d3.1-bddly-lip-extending}
\begin{rm}
(a) Let $X,Y$ be open sets in metric spaces $E$ and $F$ respectively.
Suppose that $x \in \rfs{cl}(X)$
and $g \in H(X,Y)$.
   \index{continuous. $\alpha$-continuous at $x \in \rfs{cl}(X)$}
We say that $g$ is {\it $\alpha$-continuous at $x$},
if there is $T \in \rfs{Nbr}(x)$ such that
$g \nrestriction (T \cap X)$ is $\alpha$-continuous.

Obviously, if $F$ is a complete metric space,
and $g$ is $\alpha$-continuous at $x$,
then $x \in \rfs{Dom}(g^{\srfs{cl}})$.
   \index{bicontinuous. $\alpha$-bicontinuous at $x \in \rfs{cl}(X)$}
We say that {\it $g$ is $\alpha$-bicontinuous at $x$},
if $g$ is $\alpha$-continuous at $x$,
$x \in \rfs{Dom}(g^{\srfs{cl}})$ and
$g\inverse$ is $\alpha$-continuous at $g^{\srfs{cl}}(x)$.
We say that {\it $g$ is $\itGamma$-bicontinuous at $x$},
if for some $\alpha \in \itGamma$,
$g$ is $\alpha$-bicontinuous at $x$.
   \index{bicontinuous. $\itGamma$-bicontinuous at $x \in \rfs{cl}(X)$}

(b) Suppose that $E$ is a metric space $X \subseteq E$ is open,
$b \in \rfs{bd}(X)$,
$\alpha \in \rfs{MBC}$ and $x,y \in X$. 
Recall that we denote
$\delta^{X,E}(x) = d(x,E - X)$. Superscripts $ ^E$ and $ ^X$ are omitted
when they are understood from the context.
   \index{N@d07@@$\delta^{X,E}(x) = d(x,E - X)$. Abbreviations:
          $\delta^X(x)$, $\delta(x)$}
   \index{N@AAAA@@$x \approx^{(\alpha,b)}_{(X,E)} y$. This means
          $d(x,b) \approx^{\alpha} d(y,b)$ and
          $\delta^X(x) \approx^{\alpha} \delta^X(y)$}
The notation
$x \approx^{(\alpha,b)}_{(X,E)} y$ means that
$$
d(x,b) \approx^{\alpha} d(y,b) \mbox{ \ \,and\ \ \,}
\delta^X(x) \approx^{\alpha} \delta^X(y).
$$
Suppose that $\vecx^{\srsp{\sigma}}$ and $\vecy^{\srsp{\sigma}}$
are sequences in $X$.
   \index{N@AAAA@@$\vecx^{\srsp{\sigma}} \approx^{(\alpha,b)}_{(X,E)}
          \vecy^{\srsp{\sigma}}$. This means: for every $n \in \sigma$,
$x_n \approx^{(\alpha,b)}_{(X,E)} y_n$}
   \index{N@AAAA@@$\approx^{(\alpha,b)}$.
          Abbreviation of \,$\approx^{(\alpha,b)}_{(X,E)}$}
Then $\vecx^{\srsp{\sigma}} \approx^{(\alpha,b)}_{(X,E)}
\vecy^{\srsp{\sigma}}$
means that for every
$n \in \sigma$, $x_n \approx^{(\alpha,b)}_{(X,E)} y_n$.
We abbreviate $\approx^{(\alpha,b)}_{(X,E)}$ by $\approx^{(\alpha,b)}$.
Note that the notation $\vecx \approx^{(\alpha,b)} \vecy$
entails that $\rfs{Dom}(\vecx) = \rfs{Dom}(\vecy)$.

(c) Let $X$ be a topological space, $A \subseteq H(X)$,
$\rho \subseteq \bbN$ be infinite
and $\vecx^{\srsp{\rho}},\vecy^{\srsp{\rho}}\kern1pt$
be sequences in $X$.
   \index{N@AAAA@@$\vecx^{\srsp{\rho}} \neweq^A \vecy^{\srsp{\rho}}$}
We define the relation
$\vecx^{\srsp{\rho}} \neweq^A \vecy^{\srsp{\rho}}$.
The relation $\vecx^{\srsp{\rho}} \neweq^A \vecy^{\srsp{\rho}}$
means that for every
infinite $\sigma,\eta \subseteq \rho$ there is $g \in A$ such that
$\setm{i \in \sigma}{g(x_i) = y_i}$
and
$\setm{i \in \eta}{g(x_i) = x_i}$
are infinite.
   \index{N@AAAA@@$\vecx^{\srsp{\rho}} \neweq^{\alpha}
          \vecy^{\srsp{\rho}}$}
If $\alpha \in \rfs{MBC}$,
then $\vecx^{\srsp{\rho}} \neweq^{\alpha} \vecy^{\srsp{\rho}}$
means that $\vecx^{\srsp{\rho}} \neweq^A \vecy^{\srsp{\rho}}$, where
$A = \setm{g \in H(X)}{g \mbox{ is } \alpha\mbox{-bicontinuous}}$.

(d) Let $E$ be a metric space, $X \subseteq E$ be open,
$\alpha \in \rfs{MBC}$ and $\itGamma$ be a modulus of continuity.
A sequence $\vecx$ in $X$ is called an {\it $\alpha$-abiding} sequence
if the following holds.
   \index{abiding sequence. $\alpha$-abiding sequence}
\begin{itemize}
\addtolength{\parskip}{-11pt}
\addtolength{\itemsep}{06pt}
\item[(i)]
$\vecx$ is convergent and $b \eqdf \lim \vecx \in \rfs{bd}(X)$;
\item[(ii)] There is $n = n(\vecx,\alpha) \in \bbN$ such that for every
$k \in \rfs{Dom}(\vecx)^{\geq n}$,
$d(x_n,b) \leq \alpha(\delta(x_n))$.
\vspace{-05.7pt}
\end{itemize}
A sequence $\vecx$ in $X$ is called a {\it $\itGamma$-\kern1ptevasive}
sequence if the following holds.
   \index{evasive sequence. $\itGamma$-\kern1ptevasive sequence}
\begin{itemize}
\addtolength{\parskip}{-11pt}
\addtolength{\itemsep}{06pt}
\item[(i)]
$\vecx$ is convergent and $b \eqdf \lim \vecx \in \rfs{bd}(X)$;
\item[(ii)]
For every subsequence $\vecy$ of $\vecx$ and $\alpha \in \itGamma$,
$\vecy$ is not $\alpha$-abiding.
\vspace{-05.7pt}
\end{itemize}
Equivalently, $\vecx$ is $\itGamma$-\kern1ptevasive,
iff (i) holds and for every $\alpha \in \itGamma$
there is $n \in \bbN$ such that
for every $m \in \rfs{Dom}(\vecx)^{\geq n}$,
$d(x_m,b) > \alpha(\delta(x_m))$.

(e) Let $X$ be an open subset of a normed space $E$,
and $x \in \rfs{bd}(X)$. Suppose that $X$ is two-sided at $x$,
and let $\trpl{\psi}{A}{r}$ be a boundary chart element for $x$.
Let $U,V \in \rfs{Nbr}(x)$
and $h \in \rfs{EXT}^{\pm}(U \cap X,V \cap X)$ be such that
$h^{\srfs{cl}}(x) = x$.
We say that $h$ is {\it side preserving} at $x$,
if there is $U^{\fprime} \in \rfs{Nbr}(x)$ such that
for every $u \in U^{\fprime} \cap X$, $u$ and $h(u)$ are on the
same side of $\rfs{bd}(X)$
with respect to $\trpl{\psi}{A}{r}$.
See Definition \ref{d-bddlip-1.2+1}.
We say that $h$ is {\it side reversing} at $x$,
if there is $U^{\fprime} \in \rfs{Nbr}(x)$ such that
for every $u \in U^{\fprime} \cap X$, $u$ and $h(u)$ are on different
sides of $\rfs{bd}(X)$
with respect to $\trpl{\psi}{A}{r}$.
Note that the properties of being side preserving or side reversing
does not depend on the choice of $\trpl{\psi}{A}{r}$.
   \index{side preserving at $x$}
   \index{side reversing at $x$}

(f) Let $X$ be an open subset of a normed space $E$,
and $x \in \rfs{bd}(X)$. Suppose that $\rfs{bd}(X)$ is
$1$-dimensional at $x$,
and let $\trpl{\psi}{A}{r}$ be a boundary chart element for $x$.
Let $L = \rfs{bd}(X) \cap \rfs{Rng}(\psi)$. So $L$ is an open arc.
Let $U,V \in \rfs{Nbr}(x)$
and $h \in \rfs{EXT}^{\pm}(U \cap X,V \cap X)$ be such that
$h^{\srfs{cl}}(x) = x$.
{\thickmuskip=2mu \medmuskip=1mu \thinmuskip=1mu 
We say that $h$ is {\it order preserving} at $x$,
if there is $U^{\fprime} \in \rfs{Nbr}(x)$ such that
for every $u \in U^{\fprime} \cap L$,
$u$ and $h^{\srfs{cl}}(u)$ are in the same connected component
of $L - \sngltn{x}$.
We say that $h$ is {\it order reversing} at $x$,
if there is $U^{\fprime} \in \rfs{Nbr}(x)$ such that
for every $u \in U^{\fprime} \cap X$, $u$ and $h(u)$ are in different
connected components of $L - \sngltn{x}$.
Note that the properties of being order preserving or order reversing
is independent of the choice of $\trpl{\psi}{A}{r}$.
}

   \index{order preserving at $x$}
   \index{order reversing at $x$}

Let $G \leq \rfs{EXT}(X)$.
We say that $\rfs{bd}(X)$ {\it is $G$-order-reversible at $x$,}
if there is $g \in G$ such that $g$ is order reversing at $x$,
and if $X$ is two-sided at $x$, then $g$ is side preserving.
If such a $g$ does not exist, then it is said that
$\rfs{bd}(X)$ {\it is $G$-order-irreversible at $x$}.
   \index{order-reversible. $\rfs{bd}(X)$ is $G$-order-reversible
          at $x$} \index{order-irreversible.
          $\rfs{bd}(X)$ is $G$-order-irreversible at $x$}
\end{rm}
\end{defn}

\begin{prop}\label{p3.2-bddly-lip-extending}
Let $E,F$ normed spaces. Suppose that $X \subseteq E$, $Y \subseteq F$
are open, $\alpha \in \rfs{MBC}$
and $g \in \rfs{EXT}^{\pm}(X,Y)$.
Let $b \in \rfs{bd}(X)$,
and suppose that $g$
is $\alpha$-bicontinuous at~$x$.

\num{a} There is $r_0 > 0$ such that for every $x \in B(b,r_0) \cap X$,
$\delta(x) \approx^{\alpha} \delta(g(x))$.

\num{b} Assume that $E = F$, $Y = X$ and $g(b) = b$.
Suppose that $\vecx$ is a sequence in $X$ converging to $b$.
Then for some $n \in \bbN$,
$\vecx^{\kern1pt\geq n} \approx^{(\alpha,b)}
g(\vecx)^{\kern1pt\geq n}$.

\num{c} Assume that $E = F$, $Y = X$ and $g(b) = b$.
Suppose that $X$ is two-sided at $b$.
Let $\trpl{\psi}{A}{r}$ be a boundary chart element for $b$.
Then there is $U \in \rfs{Nbr}(b)$
such that $U, g(U) \subseteq \rfs{Rng}(\psi)$,
and for every $u,v \in U \cap X$:
$u,v$  are on the same side of $\rfs{bd}(X)$ with respect to
$\trpl{\psi}{A}{r}$
\ iff \ 
$g(u),g(v)$  are on the same side
$\rfs{bd}(X)$ with respect to $\trpl{\psi}{A}{r}$.
\end{prop}

\noindent
{\bf Proof }
(a)
Let $r > 0$ be such that
$g \nrestriction (B(b,r) \cap X)$
is $\alpha$-continuous.
Choose $s > 0$ such that
$g\inverse \nrestriction (B(g(b),s) \cap Y)$
is $\alpha$-continuous,
and let $r_0$ be such that $r_0 < \dgfrac{r}{2}$
and $g(B(b,r_0) \cap X) \subseteq B(g(b),\dgfrac{s}{2}) \cap Y$.
Let $x \in B(b,r_0) \cap X$.
Suppose that $\varepsilon \in (0,\dgfrac{r}{2} - \norm{x - b})$.
Let $u \in \rfs{bd}(X)$ be such that
$\norm{u - x} < \delta(x) + \varepsilon$.
Since $\delta(x) < \norm{x - b} < r_0$,
it follows that
$$
\norm{u - b} \leq \norm{u - x} + \norm{x - b} <
\delta(x) + \dgfrac{r}{2} - \norm{x - b} + \norm{x - b} \leq
r_0  + \dgfrac{r}{2} < r.
$$
Hence
$g^{\srfs{cl}} \nrestriction \dbltn{x}{u}$ is $\alpha$-continuous.
So
$$
\delta(g(x)) \leq \norm{g^{\srfs{cl}}(x) - g^{\srfs{cl}}(u)} \leq
\alpha(\norm{x - u}) < \alpha(\delta(x) + \varepsilon).
$$
Since this argument is true for any
$\varepsilon \in (0,\dgfrac{r}{2} - \norm{x - b})$,
it follows that $\delta(g(x)) \leq \alpha(\delta(x))$.
We apply the analogous argument to $g(x)$.
This can be done, since $g(x) \in B(g(b),\dgfrac{s}{2}) \cap Y$.
So
$\delta(g\inverse(g(x))) \leq \alpha(\delta(g(x)))$.
That is, $\delta(x) \leq \alpha(\delta(g(x)))$.
We conclude that
$$\
\delta(x) \approx^{\alpha} \delta(g(x)).
$$

(b) Part (b) follows trivially from Part (a).

(c) There is $s \in (0,r)$ such that
$g(\psi(B(0,s))) \subseteq \rfs{Rng}(\psi)$.
Let
$U =\psi(B(0,s))$.
Let $u,v \in U \cap X$ be on the same side of $\rfs{bd}(X)$.
Let
$L = [\psi\inverse(u),\psi\inverse(v)]$.
Then $L \subseteq B(0,r) - A$
and thus $\psi(L) \subseteq X$.
So $g(\psi(L)) \subseteq X$.
Hence $\psi\inverse(g(\psi(L))) \subseteq B(0,r) - A$.
That is, there is an arc in $B(0,r) - A$ connecting
$\psi\inverse(g(u))$ and $\psi\inverse(g(v))$.
So $\psi\inverse(g(u))$ and $\psi\inverse(g(v))$
are on the same side of $A$. This means that $g(u)$ and $g(v)$ are
on the same side of $\rfs{bd}(X)$.\break
\rule{0pt}{0pt}\hfill\myqed

%

\begin{prop}\label{p3.3-bddly-lip-extending}
\num{a} There is $\bfiN > 1$ such that \num{a1} and \num{a2} below hold.

Let $\alpha,\beta \in \rfs{MBC}$,
$X$ be an open subset of a normed space.
Suppose that $b \in \rfs{bd}(X)$,\break
$X$ is $\beta$-LIN-bordered at $b$,
and $\rfs{bd}(X)$ is not $1$-dimensional at $b$.
Denote $\baralpha = \beta \scirc \alpha \scirc \beta$.

\num{a1}
Let $\vecx,\vecy$ be sequences in $X$ converging to $b$.
Suppose that $\vecx \approx^{(\alpha,b)} \vecy$.
Also assume that if $X$ is two-sided at $b$, then for every
$n \in \rfs{Dom}(\vecx)$,
$x_n$ and $y_n$ are on the same side of $\rfs{bd}(X)$.
Then
$\vecx
\neweq^
{\ssbfiN \kern1.5pt\ncdot\beta \scirc
\baralpha^{\sscirc 4} \scirc \beta}
\vecy$.

\num{a2}
Let $g \in \rfs{EXT}(X)$ be $\alpha$-bicontinuous at $b$.
Suppose that $g(b) = b$.
Suppose further that if $X$ is two-sided at $b$,
then $g$ is side preserving at $b$.
Let $\vecx$ be a sequence in $X$ converging to $b$.
Then
$\vecx
\neweq^
{\ssbfiN \kern1.5pt\ncdot\beta \scirc \baralpha^{\sscirc 4} \scirc
\beta}
g(\vecx)$.

\num{b}
Let $X$ be an open subset of a normed space and
$\beta \in \rfs{MBC}$.
Suppose that $b \in \rfs{bd}(X)$,
$X$ is $\beta$-LIN-bordered at $b$,
and $X$ is two-sided at $b$.
Let $g \in \rfs{EXT}(X)$ be such that $g(b) = b$,
and $g$ is side reversing at $b$.
Let $\vecx$ be a sequence in $X$ converging to $b$.
Then
$\vecx
\not\kern-3pt\neweq^{\srfs{EXT}(X)}
g(\vecx)$.
\end{prop}

\noindent
{\bf Proof }
(a)
Let $M$ be a UC-constant for
$\pair{\dghalf}{2}$, $\bfiM = M^2$ and $\bfiN =  \bfiM{\kern1.5pt^2}$.
(See Definition \ref{d-bddly-lip-2.9}(a)).
We shall prove that $\bfiN$ is as required in Part (a).

(a1) Let $X,b,\vecx,\vecy$ and $\alpha$ be as in (a1).
Let $\trpl{\psi}{A}{r}$ be a boundary chart element for $b$,
and assume that $\psi$ is $\beta$-bicontinuous.
We show that
$\vecx \neweq^
{\ssbfiN \kern1.5pt\ncdot\beta \scirc
\baralpha^{\sscirc 4} \scirc \beta}
\vecy$.
We may assume that $\vecx,\vecy \subseteq \rfs{Rng}(\psi)$.
Set $\vecw = \psi\inverse(\vecx)$
and $\vecz = \psi\inverse(\vecy)$.
Clearly, 
$\vecw \approx^{(\baralpha,0)} \vecz$.
Let $\sigma, \eta \subseteq \bbN$ be infinite.
We may assume that either for every $i \in \sigma$,
$\norm{w_i} \leq \norm{z_i}$
or for every $i \in \sigma$,
$\norm{z_i} < \norm{w_i}$.
Let us assume that the former happens.
The case that $\norm{z_i} < \norm{w_i}$ is dealt with in
a similar way.
Let $\setm{m_i}{i \in \bbN}$ and $\setm{m^1_i}{i \in \bbN}$
be respectively $\onetoone$ enumerations of $\sigma$ and $\eta$
and set $u_i = w_{m_i}$, $v_i = z_{m_i}$ and $u^1_i = w_{m^1_i}$.
So $\vecu \approx^{(\baralpha,0)} \vecv$.

We define by induction $i_n,j_n \in \bbN$ and $h_n \in H(B^E(0,r))$
such that:
\begin{itemize}
\addtolength{\parskip}{-11pt}
\addtolength{\itemsep}{06pt}
\item[(1)] $\norm{v_{i_0}} < \dgfrac{r}{2}$,
\item[(2)] $h_n(u_{i_n}) = v_{i_n}$,
\item[(3)] $h_n$ is
$\bfiM \kern2pt\ncdot \bar{\alpha} \scirc \baralpha$-bicontinuous,
\item[(4)] $\rfs{supp}(h_n) \subseteq
B(0;\half \kern1pt\norm{u_{i_n}},2 \kern1pt\norm{v_{i_n}})$,
\item[(5)] $\norm{u^1_{j_n}} < \dgfrac{\norm{u_{i_n}}}{2}$
and
$\norm{v_{i_{n + 1}}} < \dgfrac{\norm{u^1_{j_n}}}{2}$.
\item[(6)] $h_n(A) = A$.
\vspace{-05.7pt}
\end{itemize}
That the construction is possible follows from
Lemma~\ref{l-bddlip-1.10}(b).
Facts (4) and~(5) imply that
$\rfs{supp}(h_m) \cap \rfs{supp}(h_n) = \emptyset$
for any $m \neq n$.
So $h \eqdf \bcirc_n\kern2pt h_n$ is well defined.

Let $\gamma = \bfiM \kern2pt\ncdot \baralpha \scirc \baralpha$.
We verify that $h$ is $\gamma \scirc \gamma$-bicontinuous.
Let $u,v \in B^E(0,r)$. Then there are $m,n \in \bbN$ such that
$u,h_m(u) \in B^E(0,r) - \bigcup_{k \neq m} \rfs{supp}(h_k)$
and
$v,h_n(u) \in B^E(0,r) - \bigcup_{k \neq n} \rfs{supp}(h_k)$.
{\thickmuskip=3.5mu \medmuskip=1mu \thinmuskip=1mu 
If $m \neq n$, then $h(u) = h_m \scirc h_n(u)$ and 
$h(v) = h_m \scirc h_n(v)$,
and if $m = n$, then $h(u) = h_m(u)$ and  $h(v) = h_m(v)$.
Since $h_m \scirc h_n$ and $h_m$ are $\gamma \scirc \gamma$-continuous,
$\norm{h(u) - h(v)} \leq \gamma \scirc \gamma(\norm{u - v})$.
So $h$ is $\gamma \scirc \gamma$-continuous.
The same argument holds for $h\inverse$.
It follows that $h$ is $\gamma \scirc \gamma$-bicontinuous.
Since
$\gamma \scirc \gamma \leq
\bfiM^{\kern1pt2} \ncdot \baralpha^{\sscirc 4}$,
we have that
$h$ is $\bfiM^{\kern1pt2} \ncdot \baralpha^{\sscirc 4}$-bicontinuous.
By (4) and (5), $h(u^1_{j_n}) = u^1_{j_n}$ for every $n \in \bbN$.
Let
$g' = \psi \scirc h \scirc \psi\inverse \nrestriction \rfs{BCD}^E(A,r)$.
}
Then $\rfs{Dom}(g') = \rfs{Rng}(\psi) \cap X$.
Clearly, $g'$ is 
$\beta \scirc (\bfiM^{\kern1pt2} \ncdot \baralpha^{\sscirc 4}) \scirc
\beta$-bicontinuous, and hence $g'$ is
$\bfiM^{\kern1pt2} \ncdot \beta \scirc \baralpha^{\sscirc 4} \scirc
\beta$-bicontinuous.
Define $g = g' \cup \rfs{Id} \nrestriction (X - \rfs{Rng}(\psi))$.
From (1) and (4) it follows that $g \in H(X)$.
The fact that
$\bfiM^{\kern1pt2} \ncdot \beta \scirc \baralpha^{\sscirc 4} \scirc
\beta \in \rfs{MBC}$
implies that\break
$g$ too is
$\bfiM^{\kern1pt2} \ncdot \beta \scirc \baralpha^{\sscirc 4} \scirc
\beta$-bicontinuous.
Clearly, $x'_n \eqdf \psi(u^1_{j_n}) \in \setm{x_i}{i \in \eta}$
and $h(x'_n) = x'_n$.
For every $n \in \bbN$ there is $k(n) \in \sigma$ such that
$\psi(u_{i_n}) = x_{k(n)}$ and $\psi(v_{i_n}) = y_{k(n)}$.
From the fact that $h(u_{i_n}) = v_{i_n}$ it follows that
$g(x_{k(n)}) = y_{k(n)}$.
So $g$ fulfills the requirements which are needed in order to show that
$\vecx \neweq^
{\ssbfiN \kern1.5pt\ncdot\beta \scirc
\baralpha^{\sscirc 4} \scirc \beta}
\vecy$.

(a2)
It follows trivially from Proposition \ref{p3.2-bddly-lip-extending}(b)
and (a1) that $\bfiN$ is as required.

(b)
Suppose by contradiction that
$\vecx \neweq^{\srfs{EXT}(X)} g(\vecx)$.
Then $(*)$ there is $h \in \rfs{EXT}(X)$ such that
$\setm{i \in \bbN}{h(x_i) = g(x_i)}$ and
$\setm{i \in \bbN}{h(x_i) = x_i}$ are infinite.
Since $\lim \vecx = b$ and $g$ is side reversing, $(*)$ contradicts
Proposition \ref{p3.2-bddly-lip-extending}(c).
\hfill\myqed

\begin{prop}\label{p3.3.5-bddly-lip-extending}
There is $\bfiN > 1$ such that the following hold.

Let $\alpha,\beta \in \rfs{MBC}$,
and $X$ be an open subset of a normed space $E$.
Suppose that\break
$b \in \rfs{bd}(X)$,
$X$ is $\beta$-LIN-bordered at $b$,
and $\rfs{bd}(X)$ is $1$-dimensional at~$b$.
Denote $\baralpha = \beta \scirc \alpha \scirc \beta$.

\num{a}
Let $\vecx,\vecy$ be $\alpha$-abiding sequences in $X$ converging
to $b$ and $\vecx \approx^{(\alpha,b)} \vecy$.
Also assume that if $X$ is two-sided at $b$, then for every
$n \in \rfs{Dom}(\vecx)$,
$x_n$ and $y_n$ are on the same side of $\rfs{bd}(X)$.
Then
$\vecx
\neweq^
{\ssbfiN \kern1.5pt\ncdot\beta \scirc \baralpha^{\sscirc 6} \scirc
\beta}
\vecy$.

\num{b}
Let $g \in \rfs{EXT}(X)$ be $\alpha$-bicontinuous at $b$.
Suppose that $g(b) = b$.
Suppose further that if $X$ is two-sided at $b$,
then $g$ is side preserving at $b$.
Let $\vecx$ be an $\alpha$-abiding sequence in $X$ converging to $b$.
Then $\vecx
\neweq^
{\ssbfiN \kern1.5pt\ncdot\beta \scirc \baralpha^{\sscirc 6} \scirc
\beta}
g(\vecx)$.
\end{prop}

\noindent
{\bf Proof } (a) The proof follows the same steps as the proof of
Proposition~\ref{p3.3-bddly-lip-extending}(a1). But here
Lemma \ref{l3.3-temp} replaces the use of Lemma \ref{l-bddlip-1.10}
in the proof of \ref{p3.3-bddly-lip-extending}(a1).

(b) The proof follows the same steps as the proof of
Proposition \ref{p3.3-bddly-lip-extending}(a2).
\rule{1pt}{0pt}\hfill\myqed

\begin{prop}\label{p3.3.8-bddly-lip-extending}
\num{a}
There is $\bfiN > 1$ such that \num{a1} and \num{a2} below hold.

Let $\alpha,\beta \in \rfs{MBC}$,
and $X$ be an open subset of a normed space $E$.
Suppose that $b \in \rfs{bd}(X)$,
$X$ is $\beta$-LIN-bordered at $b$,
and $\rfs{bd}(X)$ is $1$-dimensional at~$b$.
Let $\trpl{\psi}{A}{r}$ be boundary chart element for $b$
with $\psi$ being $\beta$-bicontinuous.
If $A$ is a subspace of $E$, let $F = A$.
If $\rfs{dim}(E) = 2$ and $A$ is a half space of $E$, let
$F = \rfs{bd}(A)$. (So $F$ is a $1$-dimensional subspace of $E$).
Denote $\baralpha = \beta \scirc \alpha \scirc \beta$.

\num{a1}
Let $\vecx,\vecy \subseteq \rfs{Rng}(\psi)$ be sequences
which converge to $b$, and set $\vecu = \psi\inverse(\vecx)$ and
$\vecv = \psi\inverse(\vecy)$.
Assume that
\begin{itemize}
\addtolength{\parskip}{-11pt}
\addtolength{\itemsep}{06pt}
\item[\num{i}] $\vecx \approx^{(\alpha,b)} \vecy$,
\item[\num{ii}] for every $n \in \rfs{Dom}(\vecx)$ there are
$\hatu_n,u_n^{\sperp},\hatv_n,v_n^{\sperp}$ and $\lambda_n$
such that
$u_n = \hatu_n + u_n^{\sperp}$,\break
$v_n = \hatv_n + v_n^{\sperp}$,
$\hatu_n, \hatv_n \in F$, \
$u_n^{\sperp} \kern-1pt,\kern-1pt v_n^{\sperp} \perp F$,
$\lambda_n > 0$
and $\hatv_n = \lambda_n \hatu_n$,
\item[\num{iii}] if $X$ is two-sided at $b$, then for every
$n \in \rfs{Dom}(\vecx)$,
$x_n$ and $y_n$ are on the same side of $\rfs{bd}(X)$.
\vspace{-05.7pt}
\end{itemize}
\underline{Then} \ %
$\vecx \neweq^
{\ssbfiN \kern1.5pt\ncdot\beta \scirc \baralpha^{\sscirc 4} \scirc
\beta}
\vecy$.

\num{a2}
Let $\itGamma$ be a modulus of continuity,
$\alpha,\beta \in \itGamma$,
and $\vecx$ be a $\itGamma$-evasive sequence in $X$ converging to $b$.
Let $g \in \rfs{EXT}(X)$ be $\alpha$-bicontinuous at $b$,
and assume that:
$g(b) = b$,\break
$g$ is order preserving at $b$,
and if $X$ is two-sided at $b$
then $g$ is side preserving at $b$.\break
\underline{Then} \ %
$\vecx \neweq^
{\ssbfiN \kern1.5pt\ncdot\beta \scirc \baralpha^{\sscirc 4} \scirc
\beta}
g(\vecx)$.
\smallskip

In Parts \num{b}\,-\,\num{d} below we assume that
$\itGamma$ is a modulus of continuity,
$\beta \in \itGamma \cap \rfs{MBC}$,
$X$ is an open subset of a normed space $E$,
$b \in \rfs{bd}(X)$,
$X$ is $\beta$-LIN-bordered at $b$,
and $\rfs{bd}(X)$ is $1$-dimensional at~$b$.
We also assume that
$G \leq \rfs{EXT}^{\pm}(X)$, and $G$ is of boundary type $\itGamma$.

\num{b}
Let $\trpl{\psi}{A}{r}$ be boundary chart element for $b$
with $\psi$ being $\beta$-bicon\-tinuous.
If $A$ is a subspace of $E$ set $F = A$,
and if $\rfs{dim}(E) = 2$ and $A$ is a half space of $E$, set
$F = \rfs{bd}(A)$. (So $F$ is a $1$-dimensional subspace of $E$).
Let $\vecx,\vecy \subseteq \rfs{Rng}(\psi)$ be sequences
which converge to $b$, and set $\vecu = \psi\inverse(\vecx)$ and
$\vecv = \psi\inverse(\vecy)$.
Assume that
\begin{itemize}
\addtolength{\parskip}{-11pt}
\addtolength{\itemsep}{06pt}
\item[\num{i}] $\vecx,\vecy$ are $\itGamma$-evasive,
\item[\num{ii}] for every $n \in \rfs{Dom}(\vecx)$ there are
$\hatu_n,u_n^{\sperp},\hatv_n,v_n^{\sperp}$ and $\lambda_n$
such that
$u_n = \hatu_n + u_n^{\sperp}$,\break
$v_n = \hatv_n + v_n^{\sperp}$,
$\hatu_n, \hatv_n \in F$, \
$u_n^{\sperp} \kern-1pt,\kern-1pt v_n^{\sperp} \perp F$,
$\lambda_n < 0$
and $\hatv_n = \lambda_n \hatu_n$.
\vspace{-05.7pt}
\end{itemize}
\underline{Then} \ $\vecx \not\kern-3pt \neweq^G \vecy$.

\num{c}
Let $\vecx$ be a $\itGamma$-evasive sequence in $X$ converging to $b$.
Let $g \in G$.  Suppose that $g(b) = b$,
and $g$ is order reversing at $b$.
Then $\vecx \not\kern-3pt\neweq^{G} g(\vecx)$.

\num{d}
Let $\vecx$ be a sequence in $X$ converging to $b$.
Let $g \in G$ be such that $g(b) = b$
and $g$ is order preserving at $b$.
Assume further that if $X$ is two-sided at $b$,
then $g$ is side preserving.
Then $\vecx \neweq^{\itGamma} g(\vecx)$.
\end{prop}

\noindent
{\bf Proof }
(a1) The proof follows the same steps as the proof of
Proposition~\ref{p3.3-bddly-lip-extending}(a1). But here
Lemma \ref{l3.1-temp} replaces the use of Lemma \ref{l-bddlip-1.10}
in the proof of \ref{p3.3-bddly-lip-extending}(a1).

(a2) Let $\trpl{\psi}{A}{r}$ be boundary chart element for $b$
such that $\psi$ is $\beta$-bicontinuous.
If $A$ is a half space set $F = \rfs{bd}(A)$.
Otherwise, set $F = A$.
Let $B$ be an open ball with center at $b$ such that
$g^{\srfs{cl}} \nrestriction (\rfs{cl}(B) \cap \rfs{cl}(X))$ is
$\alpha$-bicontinuous, and
$\rfs{cl}(B), \,g^{\srfs{cl}}(\rfs{cl}(B) \cap \rfs{cl}(X)) \subseteq
\rfs{Rng}(\psi)$.  Let $U = \psi\inverse(B \cap X)$
and $h = (g \nrestriction (B \cap X))^{\psi\inverse}$.

We may assume that $\vecx \subseteq B$
and that $\rfs{Dom}(\vecx) = \bbN$.
Set $\vecu = \psi\inverse(\vecx)$,
and for every $n \in \bbN$ let $u_n = \hatu_n + u_n^{\sperp}$,
where $\hatu_n \in F$ and $u_n^{\sperp} \perp F$.
Denote $h(\vecu)$ by $\vecv$,
and for every $n \in \bbN$ let $v_n = \hatv_n + v_n^{\sperp}$,
where $\hatv_n \in F$ and $v_n^{\sperp} \perp F$.
Let $s > 0$ be such that
$B(0,s) \cap (E - A) \subseteq U,h(U)$.
We may assume that
$u_n, \hatu_n, u_n^{\sperp},v_n, \hatv_n, v_n^{\sperp} \in B(0,s)$
for every $n \in \bbN$.
In order to apply (a1), we need to show that
$\hatv_n = \lambda_n \hatu_n$, where $\lambda_n > 0$.
From Proposition \ref{p3.2-bddly-lip-extending}(a)
and the facts:\break
$\vecx$ is $\itGamma$-evasive, $\beta \in \itGamma$ and
$\psi$ is $\beta$-bicontinuous, it follows that
$\vecu$ is $\itGamma$-evasive.

Define $\baralpha = \beta \scirc \alpha \scirc \beta$.
Then $h$ is $\baralpha$-bicontinuous.
This implies that $\vecv$ too is $\itGamma$-evasive.
So $\lim \frac{d(u_n,F)}{\norm{u_n}} = 0$
and $\lim \frac{d(v_n,F)}{\norm{v_n}} = 0$.
We may thus assume that
$d(u_n,F) < \dgfrac{\norm{u_n}}{2}$
and $d(v_n,F) < \dgfrac{\norm{v_n}}{2}$
for every $n \in \bbN$.
It follows that for every $n$, $\hatu_n \neq 0$.
Let $\lambda_n$ be such that $\hatv_n = \lambda_n \hatu_n$.
It is trivial that $h^{\srfs{cl}}$ is $\baralpha$-bicontinuous,
and that $h^{\srfs{cl}} \nrestriction (F \cap B(0,s))$
is order preserving,
that is, for every $u \in F \cap B(0,s)$,
$u$ and $h^{\srfs{cl}}(u)$ are on the same side of $0$.
It follows that for every $n$ there is $\mu_n > 0$
such that $h^{\srfs{cl}}(\hatu_n) = \mu_n \hatu_n$.
Suppose by contradiction that for infinitely many $n$'s,
$\lambda_n \leq 0$.
Take such an~$n$.
Then
$$
\norm{h^{\srfs{cl}}(u_n) - h^{\srfs{cl}}(\hatu_n)} \leq
\baralpha(\norm{u_n - \hatu_n}) = \baralpha(\norm {u^{\sperp}_n}).
$$
But
$$
h^{\srfs{cl}}(u_n) - h^{\srfs{cl}}(\hatu_n) = v_n - \mu_n \hatu_n =
v^{\sperp}_n + \lambda_n \hatu_n - \mu_n \hatu_n =
v^{\sperp}_n - (\mu_n - \lambda_n) \hatu_n.
$$
So
\vspace{1.5mm}
\newline
\rule{7pt}{0pt}
\renewcommand{\arraystretch}{1.5}
\addtolength{\arraycolsep}{-3pt}
$
\begin{array}{ll}
&
\norm{h^{\srfs{cl}}(u_n) - h^{\srfs{cl}}(\hatu_n)} =
\norm{v^{\sperp}_n - (\mu_n - \lambda_n) \hatu_n} \geq
(\mu_n - \lambda_n) \norm{\hatu_n} - \norm{v^{\sperp}_n}
\\
\geq
\rule{5pt}{0pt}
&
\mu_n \norm{\hatu_n} - \norm{v^{\sperp}_n} =
\norm{h(\hatu_n)} - \norm{v^{\sperp}_n} \geq
\baralpha\inverse(\norm{\hatu_n}) - \norm{v^{\sperp}_n} \geq
\baralpha\inverse(\dgfrac{\norm{u_n}}{2}) -
\baralpha(\norm{u^{\sperp}_n}).
\vspace{1.7mm}
\end{array}
$
\renewcommand{\arraystretch}{1.0}
\addtolength{\arraycolsep}{3pt}
\newline
Note that
$\baralpha(\norm {u^{\sperp}_n}) =
\baralpha(\norm{u_n - \hatu_n}) \geq
\norm{h^{\srfs{cl}}(u_n) - h^{\srfs{cl}}(\hatu_n)}$.
It follows that
$$
\baralpha\inverse(\dgfrac{\norm{u_n}}{2}) -
\baralpha(\norm{u^{\sperp}_n}) \leq
\baralpha(\norm {u^{\sperp}_n}).
$$
So
\rule{0pt}{15pt}
$\norm{u_n} \leq
2 \baralpha \scirc \baralpha(\norm{u^{\sperp}_n}) +
2 \baralpha(\norm{u^{\sperp}_n}) \leq
4 \baralpha \scirc \baralpha(\norm{u^{\sperp}_n})$.
That is, $\norm{u_n} \leq 4 \baralpha \scirc \baralpha(d(u_n,F))$.
Since $4 \baralpha \scirc \baralpha \in \itGamma$,
and since the above holds for infinitely many $n$'s,
$\vecu$ is not $\itGamma$-evasive. A contradiction.
Hence for all but finitely many $n$'s, $\lambda_n > 0$.
Recall that $\vecv = h(\vecu)$. So $\vecv = \psi\inverse(g(\vecx))$.
Obviously, $\vecx \approx^{(\alpha,b)} g(\vecx)$.
Hence by Part (a), \,
$\vecx
\neweq^
{\ssbfiN \kern1.5pt\ncdot\beta \scirc \baralpha^{\sscirc 4} \scirc
\beta}
g(\vecx)$.

(b) Suppose by contradiction that there are infinite
$\sigma,\eta \subseteq \rfs{Dom}(\vecx)$ and $g \in G$
such that for every $i \in \sigma$, $g(x_i) = y_i$,
and for every $i \in \eta$, $g(x_i) = x_i$.
Let $h = g^{\psi\inverse}$.
So for some $\gamma \in \itGamma$,
$h$ is $\gamma$-bicontinuous at $0$.
Let $Y = E - A$. Then $\vecu$ is $\itGamma$-evasive with respect to
$Y$ and $E$.
Note that for every $i \in \sigma$, $h(u_i) = v_i$,
and for every $i \in \eta$, $h(u_i) = u_i$.
We abbreviate $h^{\srfs{cl}}$ by $h$.
Denote $h(\hatu_i) = \mu_i \hatu_i$. Assume by contradiction that
for infinitely many $i$'s in $\eta$, $\mu_i \leq 0$.
Since $\vecu$ is $\itGamma$-evasive, there is $n$ such that
for every $i \in \eta^{\geq n}$,
$\norm{u_i^{\sperp}} \leq \quarter \norm{u_i}$.
Let $i \in \eta^{\geq n}$, and assume that $\mu_i \leq 0$.
Then
\vspace{1.5mm}
\newline
\rule{7pt}{0pt}
\renewcommand{\arraystretch}{1.5}
\addtolength{\arraycolsep}{-0pt}
$
\begin{array}{lll}
\gamma(\delta(u_i))
&
=
&
\gamma(\norm{u_i^{\sperp}}) =
\gamma(\norm{u_i^{\sperp}}) =
\gamma(\norm{u_i - \hatu_i}) \geq
\norm{h(u_i) - h(\hatu_i)}
\\
&
=
&
\norm{u_i - \mu_i \hatu_i} =
\norm{u_i^{\sperp} + \hatu_i - \mu_i \hatu_i} \geq
(1 - \mu_i) \norm{\hatu_i} - \norm{u_i^{\sperp}}
\\
&
\geq
&
\norm{\hatu_i} - \norm{u_i^{\sperp}} \geq
\frac{3}{4} \norm{u_i} - \quarter \norm{u_i} = \half \norm{u_i}.
\vspace{1.7mm}
\end{array}
$
\renewcommand{\arraystretch}{1.0}
\addtolength{\arraycolsep}{0pt}
\newline
So $\vecu$ is not $\itGamma$-evasive, a contradiction.
It follows that there is $i$ such that $\mu_i > 0$.
This implies that $h$ is order preserving at $0$.
In particular, for every $i \in \sigma$, $\mu_i > 0$.
We claim that $\vecv$ is $\itGamma$-evasive.
\newpage
This is so, since (i) $\vecv = h(\vecu)$,
(ii) $\gamma \in \itGamma$, (iii) $h$ is $\gamma$-continuous
and (iv) $\vecu$ is $\itGamma$-evasive.
Let $n$ be such that for every $i \in \rfs{Dom}(\vecv)^{\geq n}$,
$\delta(v_i) \leq \quarter \norm{v_i}$.
Let $i \in \sigma^{\geq n}$. Then
\vspace{1.5mm}
\newline
\rule{7pt}{0pt}
\renewcommand{\arraystretch}{1.5}
\addtolength{\arraycolsep}{-0pt}
$
\begin{array}{lll}
\gamma^{\sscirc 2}(\delta(v_i))
&
\geq
&
\gamma(\delta(u_i)) =
\gamma(\norm{u_i - \hatu_i}) \geq
\norm{h(u_i) - h(\hatu_i)} = \norm{v_i - \mu_i \hatu_i}
\\
&
=
&
\norm{v_i^{\sperp} + (\lambda_i - \mu_i) \hatu_i} \geq
\abs{\lambda_i - \mu_i} \norm{\hatu_i} - \norm{v_i^{\sperp}} \geq
\abs{\lambda_i} \norm{\hatu_i} - \norm{v_i^{\sperp}}
\\
&
=
&
\norm{\hatv_i} - \norm{v_i^{\sperp}} \geq
\frac{3}{4} \norm{v_i} - \quarter \norm{v_i} = \half \norm{v_i}.
\vspace{1.7mm}
\end{array}
$
\renewcommand{\arraystretch}{1.0}
\addtolength{\arraycolsep}{0pt}
\newline
So $\vecv \nrestriction \sigma^{\geq n}$ is
$2 \cdot \gamma \scirc \gamma$-abiding.
This contradicts the fact that $\vecv$ is $\itGamma$-evasive.

(c) Let $\trpl{\psi}{A}{r}$ be boundary chart element for $b$
such that $\psi$ is $\beta$-bicontinuous.
Since $g \in G$ there is $\alpha \in \itGamma$ and
$U \in \rfs{Nbr}^E(b)$ such that
$g \nrestriction (U \cap X)$ is $\alpha$-bicontinuous.
We may assume that $\vecx \subseteq \rfs{Rng}(\psi) \cap U$.
Let $h = g^{\psi\inverse}$
and $\gamma = \beta \scirc \alpha \scirc \beta$.
Then $h$ is $\gamma$-bicontinuous.
Let $\vecu = \psi\inverse(\vecx)$ and $\vecv = h(\vecu)$.
So $\vecv = \psi\inverse(g(\vecx))$.
Also, let $u_i = \hatu_i + u_i^{\sperp}$ and
$v_i = \hatv_i + v_i^{\sperp}$, where $\hatu_i,\hatv_i \in F$
and $ u_i^{\sperp}, v_i^{\sperp} \perp F$.
Since $\vecu$ is $\itGamma$-evasive, and $\vecv = h(\vecu)$,
$\vecv$ is $\itGamma$-evasive. We may thus assume that
for every $i \in \rfs{Dom}(\vecu)$,
$\norm{u_i^{\sperp}} \leq \dgfrac{\norm{u_i}}{4}$
and $\norm{v_i^{\sperp}} \leq \dgfrac{\norm{v_i}}{4}$.
Let $\lambda_i$ be such that $\hatv_i = \lambda \hatu_i$.
Suppose by contradiction that for infinitely many $i$'s,
$\lambda_i \geq 0$.
We abbreviate $h^{\srfs{cl}}$ by $h$.
Let $\mu_i$ be such that $h(\hatu_i) = \mu_i \hatu_i$.
Since $g$ is order reversing at $b$,
$h$ is order reversing at $0$. So $\mu_i < 0$.
Let $i$ be such that $\lambda_i \geq 0$. Then
\vspace{1.5mm}
\newline
\rule{7pt}{0pt}
\renewcommand{\arraystretch}{1.5}
\addtolength{\arraycolsep}{-0pt}
$
\begin{array}{lll}
\gamma(\norm{u_i^{\sperp}})
&
\geq
&
\norm{h(u_i) - h(\hatu_i)} =
\norm{v_i^{\sperp} + \lambda_i \hatu_i - \mu_i \hatu_i} \geq
(\lambda_i - \mu_i) \norm{\hatu_i} - \norm{v_i^{\sperp}}
\\
&
\geq
&
\abs{\mu_i} \norm{\hatu_i} - \norm{v_i^{\sperp}} =
\norm{\hatv_i} - \norm{v_i^{\sperp}} \geq \dgfrac{\norm{v_i}}{2}.
\vspace{1.7mm}
\end{array}
$
\renewcommand{\arraystretch}{1.0}
\addtolength{\arraycolsep}{0pt}
\newline
But $\norm{u_i^{\sperp}} = \delta(u_i) \leq \gamma(\delta(v_i))$.
So
$2 \ncdot \gamma \scirc \gamma(\delta(v_i)) \geq \norm{v_i}$.
That is, $\vecv$ is not $\itGamma$-evasive,\break
a contradiction.
It follows that for all but finitely many $i$'s, $\lambda_i < 0$.
By Part (b), $\vecx \not\kern-3pt\neweq^{G} g(\vecx)$.

(d) Let $\sigma,\eta$ be infinite subsets of $\rfs{Dom}(\vecx)$.
Either (i) there is an infinite $\rho \subseteq \sigma$ and
$\gamma \in \itGamma$ such that $\vecx \nrestriction \rho$ is
$\gamma$-abiding; or (ii) there is an infinite $\rho \subseteq \sigma$
such that $\vecx \nrestriction \rho$ is $\itGamma$-evasive.

Suppose that case (i) happens. To get an $f \in G$ such that
$\setm{i \in \rho}{f(x_i) = g(x_i)}$
and $\setm{i \in \eta}{f(x_i) = x_i}$ are infinite,
follow the construction in
Proposition \ref{p3.3-bddly-lip-extending}(a).
However, Lemma \ref{l-bddlip-1.10} which was used in
\ref{p3.3-bddly-lip-extending}(a)
is replaced here by Lemma \ref{l3.3-temp}.
In case (ii), follow the the proof of (a2) in this proposition. 
\rule{0pt}{1pt}\hfill\myqed

Recall that we deal with the setting where we have a normed space
$E$ and an open subset $X \subseteq E$.
In this setting, when we write $\rfs{cl}(A)$ we mean
$\rfs{cl}^E(A)$.
If we wish to denote the closure of $A$ with respect to other sets,
e.g.\ the closure of $A$ with respect to $X$, then we write
$\rfs{cl}^X(A)$.

\begin{prop}\label{p3.4-bddly-lip-extending}
For a topological space $X$ and a subgroup $G \leq H(X)$,
we define the property $P_{\srfs{cmpct}}(\vecx)$ of sequences
$\vecx$ in $X$ as follows.
\begin{list}{}
{\setlength{\leftmargin}{11pt}
\setlength{\labelsep}{05pt}
\setlength{\labelwidth}{11pt}
\setlength{\itemindent}{-00pt}
\addtolength{\topsep}{-04pt}
\addtolength{\parskip}{-02pt}
\addtolength{\itemsep}{-05pt}
}
\item[]
$P_{\srfs{cmpct}}(\vecx) \equiv$
For every infinite $\sigma \subseteq \rfs{Dom}(\vecx)$
and a sequence\kern-1pt
{\thickmuskip=2.0mu \medmuskip=1mu \thinmuskip=1mu 
\hbox{
$\setm{U_i}{i \in \sigma} \in \prod \setm{\rfs{Nbr}(x_i)}{i \in \sigma}$
}
}
consisting of pairwise disjoint sets,
there is a sequence
$\setm{g_i}{i \in \sigma} \in \prod \setm{G\sprtm{U_i}}{i \in \sigma}$
such that $\bcirc \setm{g_i}{i \in \sigma} \not\in G$.
\vspace{-02.0pt}
\end{list}
Let $E$ be a normed space and $X \subseteq E$ be open.
Let $\itGamma$ be a countably generated modulus of continuity
and $G \leq \rfs{EXT}(X)$ be $\itGamma$-appropriate.
(See Definition \ref{ams-bddly-lip-bldr-d5.4}(c)).
Let $\vecx$ be a $\onetoonen$ sequence in $X$.
Then $\rfs{cl}(\rfs{Rng}(\vecx))$ is compact
iff
$P_{\srfs{cmpct}}(\vecx)$ holds.
\end{prop}

\noindent
{\bf Proof }
Suppose first that 
$\rfs{cl}(\rfs{Rng}(\vecx))$ is not compact.
Then there is an infinite $\sigma \subseteq \rfs{Dom}(\vecx)$ such that
either $\setm{x_i}{i \in \sigma}$ is spaced,
or $\setm{x_i}{i \in \sigma}$ is a non-convergent Cauchy sequence.
For every $i \in \sigma$ let
$r_i = \third \inf \setm{\norm{x_j - x_i}}{j \in \sigma - \sngltn{i}}$
and $U_i = B^X(x_i,r_i)$.
Hence $d(U_i,U_j) \geq \dgfrac{r_i}{3}$
for any $i \neq j$ in $\sigma$.
It is easily seen that that
$\setm{U_i}{i \in \bbN}$ is $\rfs{cl}(X)$-discrete.
Let
$\setm{g_i}{i \in \sigma} \in \prod \setm{G\sprtm{U_i}}{i \in \sigma}$.
So $\setm{g_i}{i \in \sigma}$ is $\rfs{cl}(X)$-discrete.
Since $G$ is $\itGamma$-appropriate,
$\bcirc \setm{g_i}{i \in \sigma} \in G$.
So $\neg P_{\srfs{cmpct}}(\vecx)$ holds.

Suppose that $\vecx$ is $\onetoonen$ and that
$\rfs{cl}(\rfs{Rng}(\vecx))$ is compact.
Let $\setm{\alpha_i}{i \in \bbN} \subseteq \itGamma$ be a generating
sequence for $\itGamma$.
That is, for every $\alpha \in \itGamma$ there is $i \in \bbN$
such that $\alpha \preceq \alpha_i$.
We also assume that for every $i \in \bbN$,
$\setm{j}{\alpha_j = \alpha_i}$ is infinite.
Let $\sigma \subseteq \rfs{Dom}(\vecx)$ be infinite,
and for every $i \in \sigma$ let $U_i \in \rfs{Nbr}^X(x_i)$.
Assume that for every $i \neq j$, $U_i \cap U_j = \emptyset$.
Since $\rfs{cl}(\rfs{Rng}(\vecx))$ is compact,
$\setm{x_i}{i \in \sigma}$ contains a $\onetoonen$ convergent
subsequence $\setm{x_{i_n}}{n \in \bbN}$.
Denote $y_n = x_{i_n}$ and $V_n = U_{i_n} \cap B(y_n,\frac{1}{n + 1})$.
For every $n$ let $g_{i_n} \in G\sprtm{V_n}$ be such that
$g_{i_n} \nrestriction V_n$ is not
$\alpha_n$-continuous. It is easy to see that such $g_{i_n}$ exists.
For $i \in \sigma - \setm{i_n}{n \in \bbN}$ let $g_i = \rfs{Id}$.
Let $y = \lim_n y_n$ and $g = \bcirc \setm{g_i}{i \in \sigma}$.
Then there is no $\alpha \in \itGamma$ and $U \in \rfs{Nbr}(y)$
such that $g \nrestriction (U \cap X)$ is $\alpha$-continuous.
We justify this claim.
Let $\alpha \in \itGamma$. Then for some $i \in \bbN$,
$\alpha \preceq \alpha_i$. Let $r > 0$ be such that
$\alpha \nrestriction [0,r) \leq \alpha_i \nrestriction [0,r)$.
There is $n$ such that $\alpha_{i_n} = \alpha_i$,
$\rfs{diam}(V_n) < r$ and $V_n \subseteq U$.
There are $u,v \in V_n$ such that
$\norm{g_{i_n}(u) - g_{i_n}(v)} > \alpha_i(\norm{u - v})$.
Since $\norm{u - v} < r$,
we have $\alpha_i(\norm{u - v}) \geq \alpha(\norm{u - v})$.
So $\norm{g_{i_n}(u) - g_{i_n}(v)} > \alpha(\norm{u - v})$.
That is, $\norm{g(u) - g(v)} > \alpha(\norm{u - v})$.
Hence $g \nrestriction (U \cap X)$ is not $\alpha$-continuous.
It follows that $g \not\in G$. So $P_{\srfs{cmpct}}(\vecx)$ holds.
\hfill\myqed

{\bf Explanation }
For a topological space $\pair{X}{\tau^X}$ and $G \leq H(X)$ let
$\fnn{\rfs{Ap}}{G \times X}{X}$ be the application function
that is, $\rfs{Ap}(g,x) = g(x)$
and let
$M(X,G)$ be the structure
$\semisixtpl{X}{\tau^X}{G}{\in}{\scirc}{\rfs{Ap}}$.
   \index{N@m0000@@$M(X,G) = 
          \semisixtpl{X}{\tau^X}{G}{\in}{\scirc}{\rfs{Ap}}$}
Note that $P_{\srfs{cmpct}}(\vecx)$ is a property of $\vecx$ which
can be expressed in $M(X,G)$.
Hence if $\vecx \subseteq X$, $P_{\srfs{cmpct}}(\vecx)$
holds and $\iso{\psi}{M(X,G)}{M(Y,H)}$,
then $P_{\srfs{cmpct}}(\psi(\vecx))$ holds.
So in the case that $X$ is an open subset of a normed space $E$
and $G$ is $\itGamma$-appropriate,
and a similar fact holds for $Y$,
then the property ``$\rfs{cl}(\rfs{Rng}(\vecx))$ is compact''
is preserved under $\psi$.
In what follows we shall define additional properties of $\vecx$
which are expressible in $M(X,G)$. So they too are preserved
under isomorphisms between $M(X,G)$ and $M(Y,H)$.

\begin{defn}\label{p3.4.5-bddly-lip-extending}
\begin{rm}
Let $X$ be a topological space,
$G \leq H(X)$ and $\vecx$ be a sequence in $X$.

(a) Let $P_{\srfs{prerep}}(\vecx)$
be the following property of $\vecx$.
\begin{itemize}
\addtolength{\parskip}{-11pt}
\addtolength{\itemsep}{06pt}
\item[(i)]
$\rfs{Dom}(\vecx) = \bbN$ and $\vecx$ is $\onetoonen$,
\item[(ii)]
no subsequence of $\vecx$ is convergent in $X$,
\item[(iii)]
$P_{\srfs{cmpct}}(\vecx)$ holds.
\vspace{-05.7pt}
\end{itemize}
A sequence $\vecx$ which fulfills $P_{\srfs{prerep}}$
is called a {\it point pre-represntative}.
   \index{point pre-representative}

(b) Let 
$P_{\srfs{cnvrg}}(\vecx)$ and $P_{\srfs{pnt}}(\vecx)$
be the following properties.
\begin{list}{}
{\setlength{\leftmargin}{11pt}
\setlength{\labelsep}{05pt}
\setlength{\labelwidth}{11pt}
\setlength{\itemindent}{-00pt}
\addtolength{\topsep}{-04pt}
\addtolength{\parskip}{-02pt}
\addtolength{\itemsep}{-05pt}
}
\item[]
$P_{\srfs{cnvrg}}(\vecx) \equiv$
For every infinite
$\sigma \subseteq \rfs{Dom}(\vecx)$ and $g \in G$,
if
$\vecx \nrestriction \sigma \kern-2pt\neweq^{\kern-3pt G}\kern-2pt
g(\vecx) \nrestriction \sigma$,
then
$\vecx \kern-2pt\neweq^{\kern-3pt G}\kern-2pt
g(\vecx)$.
\item[]
$P_{\srfs{pnt}}(\vecx) \equiv
P_{\srfs{prerep}}(\vecx) \wedge P_{\srfs{cnvrg}}(\vecx)$.
\vspace{-02.0pt}
\end{list}
\end{rm}
\end{defn}

\begin{lemma}\label{l3.5-bddly-lip-extending}
Let $\itGamma$ be a countably generated modulus of continuity.
Suppose that $E$ is a normed space, $X \subseteq E$ is open,
$X$ is locally $\itGamma$-LIN-bordered,
and
$G \leq \rfs{EXT}(X)$ is $\itGamma$-appropriate.
Let $\vecx$ be a point pre-representative in $X$.
Then
$P_{\srfs{cnvrg}}(\vecx)$ holds,
iff $\vecx$ is convergent, and \num{i}, \num{ii}, \num{iii},
\num{iv} or \num{v} below happen.
Denote $b = \lim \vecx$.
\begin{itemize}
\addtolength{\parskip}{-11pt}
\addtolength{\itemsep}{06pt}
\item[\num{i}]
For some $\beta \in \itGamma$, $X$ is $\beta$-SLIN-bordered at $b$.
\item[\num{ii}]
For some $\beta \in \itGamma$, $X$ is $\beta$-LIN-bordered at $b$,
$X$ is two-sided at $b$, and $\rfs{bd}(X)$ is not $1$-dimensional
at $b$.
\item[\num{iii}]
$\rfs{bd}(X)$ is $1$-dimensional and $G$-order-reversible at $b$,
and for some $\alpha \in \itGamma$, $\vecx$ is\break
$\alpha$-abiding.
\item[\num{iv}]
$\rfs{bd}(X)$ is $1$-dimensional and $G$-order-reversible at $b$,
and $\vecx$ is $\itGamma$-evasive.
\item[\num{v}]
$\rfs{bd}(X)$ is $1$-dimensional and $G$-order-irreversible at $b$.
\vspace{-05.7pt}
\end{itemize}
\end{lemma}

\noindent
{\bf Proof }
We shall use the following trivial facts.

{\bf Claim 1} If $\vecy \neweq^A \vecz$,
then for every infinite $\sigma \subseteq \rfs{Dom}(\vecy)$,
$\vecy \nrestriction \sigma \neweq^A \vecz \nrestriction \sigma$.
\smallskip

{\bf Claim 2}
Suppose that $\vecy$ is a sequence in $X$ converging to a point in
$\rfs{bd}(X)$.
Assume further that $\rfs{bd}(X)$ is $1$-dimensional at $\lim \vecy$.
Then either $\vecy$ is $\itGamma$-evasive,
or for some $\alpha \in \itGamma$,
$\vecy$ has an $\alpha$-abiding subsequence.
\smallskip

{\bf Claim 3}
Suppose that $\vecy$ is a sequence in $X$ converging to a point in
$\rfs{bd}(X)$.
Assume further that $\rfs{bd}(X)$ is two-sided at $\lim \vecy$.
Let $g \in \rfs{EXT}(X)$ be such that
$g^{\srfs{cl}}(\lim \vecy) = \lim \vecy$,
and suppose that $g$ is side reversing. Then
$g(\vecy) \not\kern-2.5pt\neweq^{\srfs{EXT}(X)} \vecy$.
{\bf Proof } The Claim follows trivially from
Proposition \ref{p3.2-bddly-lip-extending}(c).
\smallskip

{\bf Claim 4}
Let $\vecy$ be a sequence in $X$ such that $\vecy$ is convergent
in $\rfs{cl}(X)$.  Suppose that $g \in \rfs{EXT}(X)$ and
$g^{\srfs{cl}}(\lim \vecy\kern1pt) \neq \lim \vecy$.
Then
$g(\vecy) \not\kern-3pt\neweq^{\srfs{EXT}(X)} \vecy$.
\smallskip

The following fact does require a proof.

{\bf Claim 5} Let $\vecx$ be a point pre-representative.
If $P_{\srfs{cnvrg}}(\vecx)$ holds,
then $\vecx$ is convergent.
{\bf Proof } Suppose that $\vecx$ is not convergent.
Let $\vecy,\vecz$ be convergent subsequences of $\vecx$
such that $\lim \vecy \neq \lim \vecz$.
Assume further that $(*)$
if $\rfs{bd}(X)$ is $1$-dimensional at $\lim \vecy$,
then either $\vecy$ is $\itGamma$-evasive,
or for some $\alpha \in \itGamma$,
$\vecy$ is $\alpha$-abiding.
Since $X$ is locally $\itGamma$-LIN-bordered,
there is $g \in G$ such that
\begin{itemize}
\addtolength{\parskip}{-11pt}
\addtolength{\itemsep}{06pt}
\item[(1)]
$g^{\fss cl}(\lim \vecy\kern1pt) = \lim \vecy$
and
$g^{\srfs{cl}}(\lim \vecz\kern1pt) \neq \lim \vecz$,
\item[(2)]
If $X$ is two-sided at $\lim \vecy$, then $g$ is side preserving,
\item[(3)]
If $\rfs{bd}(X)$ is $1$-dimensional at $\lim \vecy$,
then $g$ is order preserving.
\vspace{-05.7pt}
\end{itemize}
By Propositions \ref{p3.3-bddly-lip-extending}(a2),
\ref{p3.3.5-bddly-lip-extending}(b) and
\ref{p3.3.8-bddly-lip-extending}(a2) and by $(*)$,
$g(\vecy) \neweq^G \vecy$.
By Claim 4, $g(\vecz) \not\kern-3pt\neweq^G \vecz$, and by Claim 1,
$g(\vecx) \not\kern-3pt\neweq^G \vecx$.
Hence $\neg P_{\srfs{cnvrg}}(\vecx)$ holds.
This proves Claim 5.
\smallskip

Suppose that $\vecx$ satisfies Clause (i) in the statement of the lemma.
We show that $P_{\srfs{cnvrg}}(\vecx)$ holds.
Let $g \in G$. If $g^{\srfs{cl}}(b) \neq b$,
then by Claim 4, 
$g(\vecx\fprime) \not\kern-2.5pt\neweq^G \vecx\fprime$,
for every subsequence of $\vecx\fprime$ of $\vecx$.
If $g^{\srfs{cl}}(b) = b$,
then by Proposition \ref{p3.3-bddly-lip-extending}(a2),
$g(\vecx) \neweq^G \vecx$.
So $P_{\srfs{cnvrg}}(\vecx)$ holds.

Suppose that $\vecx$ satisfies Clause (ii) in the statement of
the lemma.
Let $g \in G$. If $g^{\srfs{cl}}(b) \neq b$,
then by Claim 4, 
$g(\vecx\fprime) \not\kern-2.5pt\neweq^G \vecx\fprime$,
for every subsequence of $\vecx\fprime$ of $\vecx$.
Suppose that $g^{\srfs{cl}}(b) = b$. If $g$ is side reversing,
then by Claim 3, 
$g(\vecx\fprime) \not\kern-2.5pt\neweq^G \vecx\fprime$,
for every subsequence of $\vecx\fprime$ of $\vecx$.
If $g$ is side preserving,
then by Proposition \ref{p3.3-bddly-lip-extending}(a2),
$g(\vecx) \neweq^G \vecx$. So $P_{\srfs{cnvrg}}(\vecx)$ holds.

Suppose that $\vecx$ satisfies Clause (iii) above.
Let $g \in G$. The case $g^{\srfs{cl}}(b) \neq b$,
is treated as in (i) and (ii).
Suppose that $g^{\srfs{cl}}(b) = b$.
If $X$ is two-sided at $x$ and $g$ is side reversing,
then by Claim 3, 
$g(\vecx\fprime) \not\kern-2.5pt\neweq^G \vecx\fprime$,
for every subsequence of $\vecx\fprime$ of $\vecx$.
Suppose that either $X$ is not two-sided at $b$,
or $X$ is two-sided at $b$ and $g$ is side preserving.
Then by Proposition~\ref{p3.3.5-bddly-lip-extending}(b),
$g(\vecx) \neweq^G \vecx$.
So $P_{\srfs{cnvrg}}(\vecx)$ holds.

Suppose that $\vecx$ satisfies Clause (iv).
As above, we may assume that
$g^{\srfs{cl}}(b) = b$,
and that if $X$ is two-sided at $b$, then $g$ is side preserving.
If $g$ is order reversing at $b$, then by
Proposition \ref{p3.3.8-bddly-lip-extending}(c),
$g(\vecx\fprime) \not\kern-2.5pt\neweq^G \vecx\fprime$,
for every subsequence of $\vecx\fprime$ of $\vecx$.
If $g$ is order preserving at $b$,
then by Proposition \ref{p3.3.8-bddly-lip-extending}(a2),
$g(\vecx) \neweq^G \vecx$.
So $P_{\srfs{cnvrg}}(\vecx)$ holds.

Suppose that $\vecx$ satisfies Clause (v).
We may assume that $g^{\srfs{cl}}(b) = b$,
and that if $X$ is two-sided at $b$, then $g$ is side preserving.
Since $\rfs{bd}(X)$ is $G$-order-irreversible at $b$,
$g$ must be order preserving at $b$.
Then by Proposition \ref{p3.3.8-bddly-lip-extending}(d),
$g(\vecx) \neweq^G \vecx$.
So $P_{\srfs{cnvrg}}(\vecx)$ holds.

We have shown that if $\vecx$ is point pre-representative,
$\vecx$ is convergent, and $\vecx$ satisfies one of the clauses
(i)\,-\,(v), then $P_{\srfs{cnvrg}}(\vecx)$ holds.

Let $\vecx$ be a point pre-representative,
and suppose that $P_{\srfs{cnvrg}}(\vecx)$ holds.
By Claim 5, $\vecx$ is convergent.
Suppose by contradiction that $\vecx$ does not satisfy any of the
clauses (i)\,-\,(v).
Let $b = \lim \vecx$. Then $\rfs{bd}(X)$ is $1$-dimensional
and $G$-order-reversible at $b$,
and (1) $\vecx$ is not $\itGamma$-evasive;
(2) there is no $\alpha \in \itGamma$ such that $\vecx$ is
$\alpha$-abiding.
There is $\gamma \in \itGamma$ and a subsequence $\vecy$ of $\vecx$
such that $\vecy$ is $\gamma$-abiding.
Since $\itGamma$ is countably generated,
there is a subsequence $\vecz$ of $\vecx$ such that $\vecz$ is
$\itGamma$-evasive.
Let $g \in G$ be such that $g$ is order reversing at $b$,
and if $X$ is two-sided at $x$, then $g$ is side preserving.
By Proposition \ref{p3.3.8-bddly-lip-extending}(c),
$g(\vecz) \not\kern-3pt\neweq^G \vecz$.
So $g(\vecx) \not\kern-3pt\neweq^G \vecx$.
By Proposition \ref{p3.3.5-bddly-lip-extending}(b),
$g(\vecy) \neweq^G \vecy$.
So $\neg P_{\srfs{cnvrg}}(\vecx)$ holds. A contradiction.
\smallskip\hfill\myqed

We represent points in $\rfs{bd}(X)$ by
sequences $\vecx$ in $X$ which satisfy $P_{\srfs{pnt}}(\vecx)$.
Such sequences are called point representatives.
By the above proposition, for every $x \in \rfs{bd}(X)$, there is
$\vecx$ such that $\lim \vecx = x$ and $P_{\srfs{pnt}}(\vecx)$ holds.
So every point of $\rfs{bd}(X)$ is represented.

We shall find a property
$\varphi_{\srfs{pnteq}}(\vecx,\vecy)$
which for point representatives $\vecx,\vecy$ expresses the fact
that $\lim \vecx = \lim \vecy$.
Let $\vecx$ be a point representative. The weak stabilizer of
$\vecx$ is defined as follows.
$$
\rfs{wstab}(\vecx) = \setm{g \in G}{g(\vecx) \neweq^G \vecx}.
$$
Define
$$
P_{\srfs{pnteq}}(\vecx,\vecy) \equiv
(\rfs{wstab}(\vecx) \subseteq \rfs{wstab}(\vecy))
\vee
(\rfs{wstab}(\vecy) \subseteq \rfs{wstab}(\vecx)).
$$

For an open subset $U$ of $X$ define
$
\rfs{opcl}(U) = U \cup (\rfs{bd}(X) - \rfs{acc}^{\srfs{cl}(X)}(X - U)).
$
   \index{N@opcl@@$\rfs{opcl}(U) =
         \rfs{int}^{\srfs{cl}(X)}(\rfs{cl}^{\srfs{cl}(X)}(U))$}
Then \rule{2.5pt}{0pt}$\rfs{opcl}(U)$ is open in $\rfs{cl}(X)$.
Also, if $V \in \rfs{Ro}(\rfs{cl}(X))$,
then $V = \rfs{opcl}(V \cap X)$.
Let $\calB = \setm{\rfs{opcl}(U)}{U \mbox{ is open in } X}$.
Hence $\rfs{Ro}(\rfs{cl}(X)) \subseteq \calB$,
and so $\calB$ is an open base for $\rfs{cl}(X)$.
Every open subset $U$ of $X$ will represent $\rfs{opcl}(U)$.
So the set of open subsets of $\rfs{cl}(X)$ which are represented,
forms an open base for $\rfs{cl}(X)$.

We next define the property $P_{\srfs{blng}}(\vecx,U)$.
For a point representative $\vecx$ and an open subset $U$ of $X$,
$P_{\srfs{blng}}(\vecx,U)$ will express the fact that
$\lim \vecx \in \rfs{opcl}(U)$.
Let
{\thickmuskip=2mu \medmuskip=1mu \thinmuskip=1mu 
$$\hbox{
$P_{\srfs{blng}}(\vecx, U) \equiv$
For every sequence $\vecy$:
if $P_{\srfs{pnt}}(\vecy)$ and
$P_{\srfs{pnteq}}(\vecx,\vecy)$,
then $\rfs{Rng}(\vecy) \kern1pt-\kern1pt U$ is finite.}$$}

\begin{prop}\label{p3.6-bddly-lip-extending}
Let $\itGamma$ be a countably generated modulus of continuity.
Suppose that $E$ is a normed space, $X \subseteq E$ is open,
and $X$ is locally $\itGamma$-LIN-bordered.
Let $G$ be a $\itGamma$-appropriate subgroup of $\rfs{EXT}(X)$.

{\thickmuskip=3.6mu \medmuskip=1mu \thinmuskip=1mu 
\num{a}
Suppose that $\vecx,\vecy$ are point representatives.
Then \kern2pt
$\lim \vecx = \lim \vecy$\ \ iff\ \ $P_{\srfs{pnteq}}(\vecx,\vecy)$
holds.}

\num{b}
Let $\vecx$ be a point representative, and $U \subseteq X$ be open.
Then
$\lim \vecx \in \rfs{opcl}(U)$\ \ iff\ \ $P_{\srfs{blng}}(\vecx, U)$
holds.
\end{prop}

\noindent
{\bf Proof } (a) Let $\vecx,\vecy$ be point representatives.
If $\lim \vecx \neq \lim \vecy$, then there is
$g \in G$ such that $g$ is the identity on some neighborhood of
$\lim \vecx$ and $g(\lim \vecy) \neq \lim \vecy$.
So $g  \in \rfs{wstab}(\vecx) - \rfs{wstab}(\vecy)$.
Similarly,
\hbox{$\rfs{wstab}(\vecy) \not\subseteq \rfs{wstab}(\vecx)$.}
So $\rfs{wstab}(\vecx)$ and $\rfs{wstab}(\vecy)$ are incomparable.

Suppose that $\lim \vecx = \lim \vecy$. Denote $b = \lim \vecx$.
If for some $\alpha \in \itGamma$,
$\rfs{bd}(X)$ is $\alpha$-SLIN-bordered at $b$,
then by Proposition~\ref{p3.3-bddly-lip-extending}(a2),
$\rfs{wstab}(\vecx) = \rfs{wstab}(\vecy) = \setm{g \in G}{g(b) = b}$.

Suppose that $X$ is two-sided at $b$ and $\rfs{bd}(X)$ is not
$1$-dimensional at $b$. Then\break
$\rfs{wstab}(\vecx) = \rfs{wstab}(\vecy) =
\setm{g \in G}
{g(b) = b \mbox{ and } g \mbox{ is side preserving at } b}$.
This follows from Proposition
\ref{p3.3-bddly-lip-extending}(a2) and (b).

Suppose that $\rfs{bd}(X)$ is $1$-dimensional at $b$.
If $\rfs{bd}(X)$ is $G$-order-irreversible at $b$,
and $X$ is not two sided at $b$,
then
$\rfs{wstab}(\vecx) = \rfs{wstab}(\vecy) = \setm{g \in G}{g(b) = b}$.
This follows from Proposition \ref{p3.3.8-bddly-lip-extending}(d).
Next assume that $\rfs{bd}(X)$ is $G$-order-irreversible at $b$,
and $X$ is two sided at $b$. Then
$\rfs{wstab}(\vecx) = \rfs{wstab}(\vecy) =
\setm{g \in G}
{g(b) = b \mbox{ and } g \mbox{ is side preserving at } b}$.
This follows from Propositions \ref{p3.3.8-bddly-lip-extending}(d)
and \ref{p3.3-bddly-lip-extending}(b).

Suppose that $G$-order-reversible at $b$.
Then  by Lemma \ref{l3.5-bddly-lip-extending}, $\vecx$ is
$\itGamma$-evasive, or there is $\alpha \in \itGamma$ such that
$\vecx$ is $\alpha$-abiding. The same holds for $\vecy$.
If both $\vecx$ and $\vecy$ are evasive or both are abiding,
then $\rfs{wstab}(\vecx) = \rfs{wstab}(\vecy)$.
This follows from Propositions \ref{p3.3.5-bddly-lip-extending}(b),
\ref{p3.3.8-bddly-lip-extending}(a2),
\ref{p3.3.8-bddly-lip-extending}(c) and
\ref{p3.3-bddly-lip-extending}(b).
Suppose that $\vecx$ is evasive and $\vecy$ is abiding.
Then $\rfs{wstab}(\vecx)$ consists of all $g \in G$ such $g(b) = b$,
$g$ is order preserving at $b$, and if $X$ is two sided at $b$, then
$g$ is side preserving at $b$.
$\rfs{wstab}(\vecy)$ consists of all $g \in G$ such that $g(b) = b$,
and if $X$ is two sided at $b$, then $g$ is side preserving at $b$.
So $\rfs{wstab}(\vecx) \subseteq \rfs{wstab}(\vecy)$.
We have shown that if $\lim \vecx = \lim \vecy$,
then $\rfs{wstab}(\vecx)$ and $\rfs{wstab}(\vecy)$ are comparable.

(b) Let $\vecx$ be a point representative, $U \subseteq X$ be open
in $X$ and $b = \lim \vecx$.
If $b \in \rfs{opcl}(U)$, then for every sequence $\vecy$ in $X$,
such that $\lim \vecy = b$
there is $n$ such that $\rfs{Rng}(\vecy^{\kern1pt\geq n}) \subseteq U$.
So $P_{\srfs{blng}}(\vecx,U)$ holds.
If $b \not\in \rfs{opcl}(U)$, then there is a sequence $\vecy$ in $X$,
which converges to $b$ and such that $\rfs{Rng}(\vecy)$ is disjoint
from $U$. There is a subsequence $\vecz$ of $\vecy$ such that
$P_{\srfs{pnt}}(\vecz)$ holds.
So $P_{\srfs{pnteq}}(\vecx,\vecz)$ holds.
Hence $\neg P_{\srfs{blng}}(\vecx,U)$ holds.
\smallskip\hfill\myqed

\noindent
{\bf Proof of Theorem \ref{t-bddlip-1.8} } Part (a) of
\ref{t-bddlip-1.8} is a special case of (b), so we prove (b).
Let $X,Y,G,H$ and $\tau$ be as in (b).
Then $\tau$  induces an isomorphism $\tilde{\tau}$ between
$M(X,G)$ and $M(Y,H)$.
Clearly,
Properties $P_{\srfs{pnt}}(\vecx)$, $P_{\srfs{pnteq}}(\vecx)$
and $P_{\srfs{blng}}(\vecx)$ are preserved by $\tilde{\tau}$.
This implies the bi-extendability of $\tau$.
\hfill\myqed

\newpage

\section{The complete $\Gamma$-bicontinuity of the inducing\\
homeomorphism}
\label{s12}

In the previous chapter we have shown
that if 
$(H^{\srfs{CMP.LC}}_{\itGamma}(X))^{\tau} =
H^{\srfs{CMP.LC}}_{\itDelta}(Y)$,
then\break
$\tau \in \rfs{EXT}^{\pm}(X,Y)$.
Further, by Theorem \ref{metr-bldr-t3.27},
$\tau$ is locally $\itGamma$-bicontinuous.
In this chapter we finally conclude that $\tau$ is completely locally
$\itGamma$-bicontinuous.
However, at this point we only know to show this for principal
$\itGamma$'s.

\subsection{$\itGamma$-continuity in directions
parallel to the \hbox{boundary of $X$.}}\label{ss12.1}

\begin{defn}\label{ams-bddly-lip-bldr-d5.1}
\begin{rm}
(a) Let $S$ be a set and $\calP$ be a partition of $S$.
That is, $\calP$ is a pairwise disjoint family whose union is $S$.
Denote $S$ by $S_{\calP}$.
For $T \subseteq S$ let
   \index{N@AAAA@@$\calP \kern-1pt\restriction\kern-1pt T =
	  \setm{P \cap T}{P \in \calP}$}
$\calP \nrestriction T \eqdf \setm{P \cap T}{P \in \calP}$.
   \index{N@AAAA@@$a \sim^{\calP} b$.
	  This means: there is $P \in \calP$ such that $a,b \in P$}
Let $a \sim^{\calP} b$ mean that there is $P \in \calP$
such that $a,b \in P$.
   \index{open sum partition with respect to $X$}
   \index{N@s01@@$S_{\calP} = \bigcup \calP$}
If $X$ is a topological space, and $S \subseteq X$ is an open set,
then $\calP$ is called an {\it open sum partition with respect to $X$}.

In the parts (b)\,-\,(d) assume that
$\pair{X}{d},\pair{Y}{e}$ are metric spaces, $\iso{\tau}{X}{Y}$,
$\alpha \in \rfs{MC}$ and $\itGamma \subseteq \rfs{MC}$.
Let $\calP$ be an open sum partition with respect to $X$
and $S = S_{\calP}$.

(b) Call $\tau$ an {\it $\pair{\alpha}{\calP}$-continuous}
function,
if for every $P \in \calP$ and $x_1,x_2 \in P$,\break
$e(\tau(x_1),\tau(x_2)) \leq \alpha(d(x_1,x_2))$,
   \index{continuous. $\pair{\alpha}{\calP}$-continuous}
and call $\tau$ an {\it $\pair{\alpha}{\calP}$-inversely-continuous},
{\thickmuskip=3.5mu \medmuskip=1mu \thinmuskip=1mu 
if for every $P \in \calP$ and $x_1,x_2 \in P$,
$d(x_1,x_2) \leq \alpha(e(\tau(x_1),\tau(x_2)))$.
   \index{inversely $\pair{K}{\calP}$-continuous}
We say that $\tau$ is {\it $\pair{\alpha}{\calP}$-bicontinuous},
if for every $P \in \calP$ and $x_1,x_2 \in P$,
$e(\tau(x_1),\tau(x_2)) \leq \alpha(d(x_1,x_2))$ and
$d(x_1,x_2) \leq \alpha(e(\tau(x_1),\tau(x_2)))$.
}
   \index{bicontinuous. $\pair{K}{\calP}$-bicontinuous}\kern-10pt

(c) We say that $\tau$ is
{\it $\pair{\alpha}{\calP}$-continuous at $x$},
if there is $T \in \rfs{Nbr}(x)$ such that $T \subseteq S$
and $\tau$ is $\pair{\alpha}{\calP \nrestriction T}$-continuous,
   \index{continuous.
          $\pair{\alpha}{\calP}$-continuous at $x$}
and $\tau$ is said to be {\it $\pair{\itGamma}{\calP}$-continuous
at $x$}, if there is $\alpha \in \itGamma$ such that
$\tau$ is $\pair{\alpha}{\calP}$-continuous at $x$.
The notions $\pair{\alpha}{\calP}$-inverse-continuity at $x$,
$\pair{\alpha}{\calP}$-bicontinuity at $x$,
$\pair{\itGamma}{\calP}$-inverse-continuity at $x$ and
$\pair{\itGamma}{\calP}$-bicontinuity at $x$ are defined analogously.
   \index{continuous. $\pair{\itGamma}{\calP}$-continuous at $x$}

(d) Call $\tau$ a {\it locally $\pair{\itGamma}{\calP}$-continuous}
function,
if for every $x \in S$,
$\tau$ is {\it $\pair{\itGamma}{\calP}$-continuous at $x$}.
   \index{locally $\pair{\alpha}{\calP}$-continuous}
The notions of local $\pair{\itGamma}{\calP}$-inverse-continuity
and local $\pair{\itGamma}{\calP}$-bicontinuity are defined analogously.
\hfill\proofend
\end{rm}
\end{defn}

The partitions $\calP$ that will be used here
have of the following form.
Let $F$ be a closed linear subspace of $E$. $\calP$ is the partiton
of $E$ into the cosets of $F$.

The next goal is to show that if 
$(H^{\srfs{CMP.LC}}_{\itGamma}(X))^{\tau} \subseteq
H^{\srfs{CMP.LC}}_{\itGamma}(Y)$,
then for every\break
$x \in \rfs{bd}(X)$ there is $\alpha \in \itGamma$
and a neighborhood of the identity in the group of translations parallel
to the boundary of $X$ such that for every $h$ in this neighborhood,
$h^{\tau}$ is $\alpha$-bicontinuous at $\tau^{\srfs{cl}}(x)$.

Recall that the notion of decayability was defined in
Definition \ref{metr-bldr-d3.1}(c).
We shall use it now again for the following situation.
Let $\rfs{BCD}^E(A,r)$ be a linear boundary chart domain,
$X = \rfs{cl}^{B(0,r)}(\rfs{BCD}^E(A,r))$,
$H = \setm{\rfs{tr}_v}{v \in \rfs{bd}^E(A)}$ and
$\lambda$ be the natural action of $H$ on $X$.
Then $\lambda$ is decayable.

When dealing with partial actions, it is often the case that
we wish to perform a composition
$g \scirc f$, where $\rfs{Rng}(f) \not\subseteq \rfs{Dom}(g)$.
Such a composition is considered to be legal.
The domain of the resulting function is
$f\inverse(\rfs{Rng}(f) \cap \rfs{Dom}(g))$.

\begin{prop}\label{ams-bddly-lip-bldr-p5.3}
\num{a}
Suppose that $\rfs{BCD}^E(A,r\fprime)$ be a linear boundary chart domain
and $L = \rfs{bd}^E(A)$. So $L$ is a closed subspace of $E$.
Let $L\fprime = L \cap B(0,r\fprime)$.
So $L\fprime = \rfs{bd}^{B(0,r\fprime)}(\rfs{BCD}^E(A,r\fprime))$.
Let $X = \rfs{BCD}^E(A,r\fprime) \cup L\fprime$
and $H = \setm{\rfs{tr}_v^E}{v \in L}$.
We equip $H$ with the norm topology of $L$.
Let $\lambda$ be defined as follows.
$$\hbox{
$\rfs{Dom}(\lambda) =
\setm{\pair{h}{z}}{h \in H \mbox{ and }z,h(z) \in X}$\ \ and\ \ %
$\lambda(h,z) = h(z)$.
}$$
Then $\lambda$ is a partial action of $H$ on $X$.

\num{b} Let $\rfs{BCD}^E(A,r\fprime)$ etc.\ be as in Part \num{a}
and $\alpha(t) = 2t$.
Then $\lambda$ is $\alpha$-decayable in $X$.

\num{c} Let $\rfs{BCD}^E(A,r\fprime)$ etc.\ be as in Part \num{a}.
Then for every $x \in L\fprime$, $x$ is a $\lambda$-limit-point.
\end{prop}

\noindent
{\bf Proof }
(a) This part is trivial.

(b) It suffices to check that $\lambda$ is $\alpha$-decayable at $0$.
We take $r_0$ to be $r\fprime$.
For $r \in (0,r\fprime)$ we take $V = V_{0,r}$ to be
$\setm{\rfs{tr}^E_v}{v \in B^L(0,\dgfrac{r}{4}}$.
So indeed $V \times B(0,ar) \subseteq \rfs{Dom}(\lambda)$.
(Recall that $a = \dgfrac{1}{2}$).
It thus suffices to show that for every normed space $E$, $r > 0$
and $v \in B(0,\dgfrac{r}{4})$ there is $g \in H(E)$ such that
\begin{itemize}
\addtolength{\parskip}{-11pt}
\addtolength{\itemsep}{06pt}
\item[(i)] $\rfs{supp}(g) \subseteq B(0,r)$,
\item[(ii)] for every $x \in E$, $g(x) - x \in \rfs{span}(\sngltn{v})$,
\item[(iii)] $g \nrestriction B(0,\dgfrac{r}{2}) =
\rfs{tr}_v \nrestriction B(0,\dgfrac{r}{2})$,
\item[(iv)] $g$ is $2$-bilipschitz.
\vspace{-05.7pt}
\end{itemize}
Let $\fnn{k}{[0,\infty)}{[0,\infty)}$ be the piecewise linear function
such that $k(t) = 1$ for $t \in [0,\dgfrac{r}{2}]$,
$k(t) = 0$ for $t \geq r$, and $k$ is linear in $[\dgfrac{r}{2},r]$.
So $(k \nrestriction [\dgfrac{r}{2},r])(t) = 2 - \dgfrac{2t}{r}$.
Let
$$
g(x) = x + k(\norm{x}) \ncdot v.
$$
It is trivial that (i)\,-\,(iii) hold.
We check that (iv) holds.
Let $x,y \in E$. If $\norm{x},\norm{y} \geq r$ or
$\norm{x},\norm{y} \leq \dgfrac{r}{2}$,
then $\norm{g(x) - g(y)} = \norm{x - y}$.
Let $u = g(x)$ and $w = g(y)$.
Assume first that
$\norm{x},\norm{y} \in [\dgfrac{r}{2},r]$.
Then
$u - w = (x - y) - \frac{2}{r}(\norm{x} - \norm{y}) \ncdot v$. \ So
$$\hbox{$
\norm{u - w} \leq
\norm{x - y} + \frac{2}{r} \norm{x - y} \ncdot \norm{v} <
(1 + \frac{2}{r} \ncdot \frac{r}{4}) \norm{x - y} =
\frac{3}{2} \norm{x - y}
$}$$
and
$$\hbox{$
\norm{u - w} \geq
\norm{x - y} - \frac{2}{r} \norm{x - y} \ncdot \norm{v} >
(1 - \frac{2}{r} \ncdot \frac{r}{4}) \norm{x - y} =
\frac{1}{2} \norm{x - y}.
$}$$
That is, $\norm{x - y} < 2 \norm{u - w}$.

Suppose that $\dgfrac{r}{2} < \norm{x} \leq r$
and $\norm{y} < \dgfrac{r}{2}$.
Let $z \in [x,y]$ be such that $\norm{z} = \dgfrac{r}{2}$.
Let $f \in \dbltn{g}{g\inverse}$. Then
\vspace{1.5mm}
\newline
\rule{7pt}{0pt}
\renewcommand{\arraystretch}{1.5}
\addtolength{\arraycolsep}{-0pt}
$
\begin{array}{lll}
\norm{f(x) - f(y)}
&
\leq
&
\norm{f(x) - f(z)} + \norm{f(z) - f(y)} \leq
2 \norm{x - z} + \norm{z - y}
\\
&
<
&
2(\norm{x - z} + \norm{z - y}) =
2 \norm{x - y}.
\vspace{1.7mm}
\end{array}
$
\renewcommand{\arraystretch}{1.0}
\addtolength{\arraycolsep}{0pt}
\newline
The case that $\dgfrac{r}{2} < \norm{x} \leq r$ and $\norm{y} > r$
is dealt with in a similar way.
The case that $\norm{x} < \dgfrac{r}{2}$ and $\norm{y} > r$ too
is dealt with in a similar way.

(c) It is trivial that every $x \in X$,
and in particular every $x \in L'$, is a $\lambda$-limit-point.
\hfill\myqed

\begin{defn}\label{ams-bddly-lip-bldr-p5.5}
\begin{rm}
Let $\pair{X}{d}$ be a metric space,
$\calP$ be an open sum partition with $S = S_{\calP}$,
$H$ be a topological group and $\lambda$ be a partial action of
$H$ on $X$. Denote the unit of $H$\break
by $e_H$,
and for $g \in H$ set $\hatg = g_{\lambda}$.
   \index{translation-like. $\pair{H}{\lambda}$
          is $\calP$-translation-like at $x$}

(a) Let $x \in S$. We say that $\pair{H}{\lambda}$ is
{\it $\calP$-translation-like at $x$},
if
for every $M \in \rfs{Nbr}(e_H)$ and $U \in \rfs{Nbr}(x)$
there are:
\begin{itemize}
\addtolength{\parskip}{-11pt}
\addtolength{\itemsep}{06pt}
\item[$(i)$]
$N \in \rfs{Nbr}(e_H)$ such that $N \subseteq M$,
\item[$(ii)$]
$T,B \in \rfs{Nbr}(x)$
such that $T \subseteq B \subseteq S \cap U$
and $N \times B \subseteq \rfs{Dom}(\lambda)$,
\item[$(iii)$]
$K > 0$;
\vspace{-05.7pt}
\end{itemize}
such that
for every $P \in \calP$ and distinct $x_0,x_1 \in P \cap T$ there are
$n \leq \frac{K}{d(x_0,x_1)}$ and \break
$g_1, \ldots, g_n \in N$ which satisfy:
\begin{itemize}
\addtolength{\parskip}{-11pt}
\addtolength{\itemsep}{06pt}
\item[\num{1}] $g_1 = e_H$,
\item[\num{2}] for every $i = 1, \ldots, n - 1$, \,
$\hatg_i(x_1) = \hatg_{i + 1}(x_0)$,
\item[\num{3}] $\hatg_n(x_1) \not\in B$.
\vspace{-02.0pt}
\end{itemize}

(b) Let $L \subseteq S$. We say that 
$\pair{H}{\lambda}$ is {\it $\calP$-translation-like in $L$},
if for every $x \in L$,
$\pair{H}{\lambda}$ is $\calP$-translation-like at $x$.
   \index{translation-like.
   $\pair{H}{\lambda}$ is $\calP$-translation-like in $L$}
\end{rm}
\end{defn}

The notion of a $\calP$-translation-like action will be used in
the following setting.
Let $\rfs{BCD}^E(A,r)$ be a linear boundary chart domain,
$X = \rfs{cl}^{B(0,r)}(\rfs{BCD}^E(A,r))$
and
\newline
$H = \setm{\rfs{tr}_v}{v \in \rfs{bd}^E(A)}$.
The natural partial action of $H$ on $X$ is translation-like.

\begin{prop}\label{ams-bddly-lip-bldr-p5.6}
Let $\rfs{BCD}^E(A,r)$ be a linear boundary chart domain,
$L = \rfs{bd}^E(A)$. So $L$ is a closed subspace of $E$.
Let $L\fprime = L \cap B(0,r)$.
So $L\fprime = \rfs{bd}^{B(0,r)}(\rfs{BCD}^E(A,r))$.
Let $X = \rfs{BCD}^E(A,r) \cup L\fprime$,
$\calP = \setm{X \cap (v + L)}{v \in X}$
and $H = \setm{\rfs{tr}_v^E}{v \in L}$.
We equip $H$ with the norm topology of $L$.
Let $\lambda$ be the following partial action of $H$ on $X$.
$$\hbox{
$\rfs{Dom}(\lambda) =
\setm{\pair{h}{z}}{h \in H \mbox{ and }z,h(z) \in X}$\ \ and\ \ %
$\lambda(h,z) = h(z)$.
}$$
Then $\lambda$ is $\calP$-translation-like in $X$.
\end{prop}

\noindent
{\bf Proof } The proof is trivial.\smallskip\hfill\myqed

The following lemma will be applied to the group of translations
in a direction parallel to the boundary of a linear boundary
chart domain.
This lemma captures the main argument in the proof of
Lemma \ref{ams-bddly-lip-bldr-l5.9}.

\begin{lemma}\label{ams-bddly-lip-bldr-l5.7}
Let $\pair{X}{d^X}$ and $\pair{Y}{d^Y}$ be metric spaces,
and $\iso{\tau}{X}{Y}$.
Let $\itGamma$ be a countably generated modulus of continuity,
and let $\alpha \in \rfs{MBC}$.
Let $S \subseteq X$ be open, and $\calP$ be a partition of $S$.
Let $H$ be a topological group and $\lambda$ be a partial action
of $H$ on $X$. Let $x \in S$.
Assume that:
\begin{itemize}
\addtolength{\parskip}{-11pt}
\addtolength{\itemsep}{06pt}
\item[\num{i}] $S \subseteq \rfs{Fld}(\lambda)$,
\item[\num{ii}] $\lambda$ is $\calP$-translation-like at $x$,
\item[\num{iii}] $\lambda$ is $\alpha$-decayable in $S$,
\item[\num{iv}] $x$ is a $\lambda$-limit-point,
\item[\num{v}] There is $U \in \rfs{Nbr}(x)$
such that for every $g \in H(X)$,
if $\rfs{supp}(g) \subseteq U$
and $g$ is $\alpha \scirc \alpha$-bicontinuous,
then $g^{\tau}$ is $\itGamma$-bicontinuous at $\tau(x)$.
\vspace{-05.7pt}
\end{itemize}
Then $\tau$ is inversely $\pair{\itGamma}{\calP}$-continuous at $x$.
\end{lemma}

\noindent
{\bf Proof } Suppose by contradiction that
$\tau$ is not inversely $\pair{\itGamma}{\calP}$-continuous at~$x$.
The conditions of Lemma \ref{metr-bldr-l3.11} hold for $x$, according
to the following correspondence.
The group $G$ of \ref{metr-bldr-l3.11} is $H(X)$ here,
and $N$ of \ref{metr-bldr-l3.11} is $S$ here.
Also, since $x$ is a $\lambda$-limit-point,
{\thickmuskip=4.0mu \medmuskip=1mu \thinmuskip=1mu 
$\kappa \eqdf
\min(\setm{\kappa(x,V_{\lambda}(x))}{V \in \rfs{Nbr}(e_H)}) \geq
\aleph_0$.
Hence $\itGamma$ is $(\leq\kern-3pt\kappa)$-generated.
It follows from \ref{metr-bldr-l3.11} that there are
$V \in \rfs{Nbr}(x)$, $M \in \rfs{Nbr}(e_H)$
and $\gamma \in \itGamma$ such that
$M \times V \subseteq \rfs{Dom}(\lambda)$,
and}
\begin{itemize}
\addtolength{\parskip}{-11pt}
\addtolength{\itemsep}{06pt}
\item[$(i)$] for every $h \in M$,
$(h_{\lambda})^{\tau} \nrestriction \tau(V)$ is $\gamma$-bicontinuous.
\vspace{-05.7pt}
\end{itemize}
For $g \in H$ denote $\hatg = g_{\lambda}$.
Since $\lambda$ is $\calP$-translation-like at $x$,
there are:
\begin{itemize}
\addtolength{\parskip}{-11pt}
\addtolength{\itemsep}{06pt}
\item[$(ii)$] $N \in \rfs{Nbr}(e_H)$ such that $N \subseteq M$,
\item[$(iii)$] $T,B \in \rfs{Nbr}(x)$
such that $T \subseteq B \subseteq S \cap V$,
\item[$(iv)$] $K > 0$;
\vspace{-05.7pt}
\end{itemize}
such that for every $P \in \calP$ and distinct $x_0,x_1 \in P \cap T$
there are
$n \leq \dgfrac{K}{d(x_0,x_1)}$ and
$e_H = g_1, \ldots, g_n \in N$ which satisfy:
$\hatg_i(x_1) = \hatg_{i + 1}(x_0)$
for every $i = 1, \ldots, n - 1$,
and $\hatg_n(x_1) \not\in B$.

Let $C = \tau(B)$ and $y = \tau(x)$.
Since $C$ is a neighborhood of $y$,
$d \eqdf d(y,Y - C) > 0$.
Let $t > 0$ be such that $\tau(B(x,t)) \subseteq B(y,\dgfrac{d}{2})$
and $B(x,t) \subseteq T$.
Set $K^* = \dgfrac{2K}{d}$.
By Clause~M2 in Definition \ref{nn1.3},
$K^{\raisedstar} \ncdot \gamma \in \itGamma$.
We have assumed that $\tau\inverse$ is not
$\pair{\itGamma}{\tau(\calP)}$-continuous at $y$.
Hence there are $P \in \calP$ and
$y_0,y_1 \in \tau(B(x,t) \cap P)$
such that
$d(\tau\inverse(y_0),\tau\inverse(y_1)) >
K^{\raisedstar} \gamma(d(y_0,y_1))$.
For $\ell = 0,1$ let $x_{\ell} = \tau\inverse(y_{\ell})$,
hence $x_0,x_1 \in B(x,t) \subseteq T$.
So there are $n \leq \dgfrac{K}{d(x_0,x_1)}$
and $ e_H = g_1, \ldots, g_n \in N$ such that
for every $i = 1, \ldots, n - 1$, $g_i(x_1) = g_{i + 1}(x_0)$
and $g_n(x_1) \not\in B$.
For $i = 2, \ldots, n$ let $x_i = g_i(x_1)$
and $y_i = \tau(x_i)$.
Since $y_0 \in \tau(B(x,t)) \subseteq B(y,\dgfrac{d}{2})$,
we have that $d(y_0,y) < \dgfrac{d}{2}$.
Note that
\begin{itemize}
\addtolength{\parskip}{-11pt}
\addtolength{\itemsep}{06pt}
\item[(1)] For every $i = 1, \ldots, n$, $g_i^{\tau}(y_0) = y_{i - 1}$
and $g_i^{\tau}(y_1) = y_i$,
\vspace{-05.7pt}
\end{itemize}
and recall that
\begin{itemize}
\addtolength{\parskip}{-11pt}
\addtolength{\itemsep}{06pt}
\item[(2)] $y_0,y_1 \in \tau(B(x,t)) \subseteq \tau(V)$,
\item[(3)] $g_1, \ldots, g_n \in N \subseteq M$.
\vspace{-05.7pt}
\end{itemize}
So by $(i)$ and (1)\,-\,(3),
$d(y_{i - 1},y_i) \leq \gamma(d(y_0,y_1))$
for every $i = 1, \ldots, n$.
Recall that $d(x_0,x_1) > K^* \gamma(d(y_0,y_1))$.
Also, $x_n \not\in B$ and hence $y_n \not\in C$.
\smallskip
So
\vspace{1.5mm}
\newline
\rule{5pt}{0pt}
\renewcommand{\arraystretch}{1.5}
\addtolength{\arraycolsep}{-3pt}
$
\begin{array}{ll}
&
d(y,Y - C) \leq d(y,y_n) \leq
d(y,y_0) + \sum_{1 = 1}^n d(y_{i - 1},y_i) <
\dgfrac{d}{2} + n \gamma(d(y_0,y_1))
\\
\leq
\rule{5pt}{0pt}
&
\dgfrac{d}{2} + \frac{K}{d(x_0,x_1)} \mcdot \gamma(d(y_0,y_1)) <
\dgfrac{d}{2} +
\frac{K}{K^* \gamma(d(y_0,y_1))} \mcdot \gamma(d(y_0,y_1))
\\
=
&
\dgfrac{d}{2} +
\frac{K}{\frac{2K}{d} \gamma(d(y_0,y_1))} \mcdot \gamma(d(y_0,y_1)) = d.
\vspace{1.7mm}
\end{array}
$
\renewcommand{\arraystretch}{1.0}
\addtolength{\arraycolsep}{3pt}
\newline
But $d(y,Y - C) = d$, a contradiction.
\smallskip\hfill\myqed

\begin{lemma}\label{ams-bddly-lip-bldr-l5.9}
Assume the following situation.
\begin{itemize}
\addtolength{\parskip}{-11pt}
\addtolength{\itemsep}{06pt}
\item[\num{1}]
$\itGamma,\itDelta$ are countably generated moduli of continuity.
\item[\num{2}]
$X \subseteq E$ and $Y \subseteq F$ are open subsets of
the normed spaces $E$ and $F$,
$X$ is $\itGamma$-LIN-\break
bordered and
$Y$ is $\itDelta$-LIN-bordered.
\item[\num{3}]
$\tau \in \rfs{EXT}^{\pm}(X,Y)$,
$G$ is a $\itGamma$-appropriate subgroup of $\rfs{EXT}(X)$,
$H$ is a $\itDelta$-appropriate
subgroup of $\rfs{EXT}(Y)$
and $G^{\tau} = H$.
\item[\num{4}]
$x \in \rfs{bd}(X)$,
$\trpl{\varphi}{A}{r}$ is a boundary chart element for $x$,
$\gamma \in \itGamma$
and $\varphi$ is $\gamma$-bicontinuous.
\item[\num{5}]
$y \in \rfs{bd}(Y)$,
$\trpl{\psi}{B}{s}$ is a boundary chart element for $y$,
$\delta \in \itDelta$
and $\psi$ is $\delta$-bicontinuous.
\item[\num{6}]
$\tau^{\srfs{cl}}(x) = y$
and
$\tau(\varphi(\rfs{BCD}^{E}(A,r))) \subseteq \psi(\rfs{BCD}^{F}(B,s))$.
\item[\num{7}] Set
$L = \rfs{bd}(A)$,
$\swhatX = \rfs{BCD}^{E}(A,r) \cup (L \cap B(0,r))$,
$\hattau = \psi\inverse \scirc \tau^{\srfs{cl}} \scirc \varphi$,
$\whatY = \hattau(\swhatX)$
and $\calP = \setm{(v + L) \cap \swhatX}{v \in \swhatX}$.
\vspace{-05.7pt}
\end{itemize}
\underline{Then}
$\hattau$ is inversely $\pair{\itDelta}{\calP}$-continuous at $0$.
\end{lemma}

\noindent
{\bf Proof }
We may assume that $X - \rfs{Rng}(\varphi) \neq \emptyset$.
From the fact that $G$ has boundary type $\itGamma$ it follows that
there is
$Z \in \rfs{Nbr}^E(x)$
such that
$G\sprt{Z \cap X} \supseteq
H^{\srfs{CMP.LC}}_{\itGamma}(X)\sprt{Z \cap X}$.
We may also assume that
$\varphi(\rfs{BCD}^{E}(A,r)) \subseteq Z$.

We wish to apply Lemma \ref{ams-bddly-lip-bldr-l5.7} to
$\swhatX,\whatY$ and $\hattau$.
More specifically, the roles of the objects mentioned
in \ref{ams-bddly-lip-bldr-l5.7} are taken by the following objects
here.
The role of $\itGamma$ in \,\ref{ams-bddly-lip-bldr-l5.7} is
taken by $\itDelta$ here,
the spaces $X,Y$ in \,\ref{ams-bddly-lip-bldr-l5.7} are
$\swhatX,\whatY$ here,
$\tau$ of \,\ref{ams-bddly-lip-bldr-l5.7} is $\hattau$,
$\alpha$ of \,\ref{ams-bddly-lip-bldr-l5.7} is the function $y = 2x$,
$S$ is $\swhatX$ and $\calP$ of \,\ref{ams-bddly-lip-bldr-l5.7}
is $\calP$ here.
The topological group $H$ appearing in \ref{ams-bddly-lip-bldr-l5.7}
is $\setm{\rfs{tr}^E_v}{v \in L}$ equipped with the norm
topology of $L$,
and $\lambda$ is the natural partial action of
$\setm{\rfs{tr}^E_v}{v \in L}$ on $\swhatX$.

Our next goal is to define the open set $U$ appearing in
Clause \num{v} of \,\ref{ams-bddly-lip-bldr-l5.7}.
We first check that
$\varphi(\swhatX) = \rfs{cl}(X) \cap \rfs{Rng}(\varphi)$
and that $\varphi(\swhatX)$ is open in $\rfs{cl}(X)$.
Clearly, $\swhatX \subseteq \rfs{cl}^E(\rfs{BCD}^{E}(A,r))$.
So if $u \in \swhatX$, then by the continuity of $\varphi$,
$\varphi(u) \in \rfs{cl}^E(\varphi(\rfs{BCD}^{E}(A,r))) \subseteq
\rfs{cl}(X)$.
That is $\varphi(\swhatX) \subseteq \rfs{cl}(X)$.
$\swhatX$ is closed in $B(0,r)$ and so
$B(0,r) - \swhatX$ is open in $B(0,r)$.
So $B(0,r) - \swhatX$ is open in $E$. Since $\varphi$ takes open
subsets of $E$ to open subsets of $E$,
$\varphi(B(0,r) - \swhatX)$ is open in $E$.
Also, $\varphi(B(0,r) - \swhatX) \cap X = \emptyset$.
So $\varphi(B(0,r) - \swhatX) \cap \rfs{cl}(X) = \emptyset$.
It follows that
$\rfs{Rng}(\varphi) \cap \rfs{cl}(X) = \varphi(\swhatX)$.
From the fact that $\rfs{Rng}(\varphi)$ is open in $E$
it follows that $\varphi(\swhatX)$ is open in $\rfs{cl}(X)$.

{\thickmuskip=3mu \medmuskip=1mu \thinmuskip=1mu 
Since $x = \varphi(0)$ and $0 \in \swhatX$,
it follows that $\varphi(\swhatX) \in \rfs{Nbr}^{\srfs{cl}(X)}(x)$.
So 
$d(x,\rfs{cl}(X) - \varphi(\swhatX)) > 0$.
}
Let $r\fprime \in (0,r)$ be such that
$\rfs{diam}(\varphi(\swhatX \cap B(0,r\fprime))) <
\dgfrac{d(x,\rfs{cl}(X) - \varphi(\swhatX))}{2}$.
The open set $U$ appearing in
Clause \num{v} of \,\ref{ams-bddly-lip-bldr-l5.7}
is $\swhatX \cap B^E(0,r\fprime)$.

We have to show that Clauses (i)\,-\,(v) of
\ref{ams-bddly-lip-bldr-l5.7} hold.
It follows from Proposition \ref{ams-bddly-lip-bldr-p5.3}(b)
that $\lambda$ is $\alpha$-decayable in $\swhatX$,
and it follows from Proposition \ref{ams-bddly-lip-bldr-p5.3}(c)
that $0$ is a $\lambda$-limit-point.
It follows from Proposition \ref{ams-bddly-lip-bldr-p5.6}
that $\lambda$ is $\calP$-translation-like at $0$.

We check that $U$ satisfies Clause \num{v}
of \ref{ams-bddly-lip-bldr-l5.7}.
Note that $\swhatX = \rfs{cl}^{B(0,r)}(\rfs{BCD}^{E}(A,r))$.
We shall also use the fact that if
$\rfs{cl}^E(A) \subseteq \rfs{Dom}(\varphi)$,
then $\rfs{cl}^E(\varphi(A)) = \varphi(\rfs{cl}^E(A))$.
This follows from the fact that $\varphi$ takes closed subsets of $E$
to closed subsets of $E$.

Let $\beta = \alpha \scirc \alpha$. So $\beta(t) = 4t$.
Let $g \in H(\swhatX)$ be $\beta$-bicontinuous
and $\rfs{supp}(g) \subseteq U$.
In order to prove that Clause~\num{v} is fulfilled,
it has to be shown that
$g^{\hattau}$ is $\itDelta$-bicontinuous at $\hattau(0^E)$.
Recall that
$\hattau = \psi\inverse \scirc \tau^{\srfs{cl}} \scirc \varphi$.
So
$g^{\hattau} = ((g^{\varphi})^{\tau^{\srfs{cl}}})^{\psi\inverse}$.
Set $\hath = g^{\varphi}$
and $ \rho = \gamma \scirc \beta \scirc \gamma$.
Since $g$ is $\beta$-bicontinuous
and $\varphi$ is $\gamma$-bicontinuous,
it follows that
$\hath$ is $\rho$-bicontinuous.
Also, $\beta \in \itGamma^{\srfs{LIP}} \subseteq \itGamma$
and $\gamma \in \itGamma$, so $\rho \in \itGamma$.

Note that $\rfs{Dom}(g) = \swhatX \subseteq \rfs{Dom}(\varphi)$.
So $\rfs{Dom}(\hath) = \varphi(\swhatX)$
and hence $\rfs{Dom}(\hath)$ is open in $\rfs{cl}(X)$.
It follows trivially from the definitions of $\swhatX$ and $U$
that $\rfs{cl}^E(U) \subseteq \swhatX$.
Hence $\rfs{cl}^E(\rfs{supp}(g)) \subseteq \whatX \subseteq
\rfs{Dom}(\varphi)$
and so
$\rfs{cl}^E(\varphi(\rfs{supp}(g))) =
\varphi(\rfs{cl}^E(\rfs{supp}(g)))$.
So
$$
\rfs{cl}^E(\rfs{supp}(\hath)) =
\rfs{cl}^E(\varphi(\rfs{supp}(g))) =
\varphi(\rfs{cl}^E(\rfs{supp}(g))) \subseteq
\varphi(\rfs{cl}^E(U) \subseteq \varphi(\whatX) = \rfs{Dom}(\hath).
$$
Let
$\barh =
\hath \cup \rfs{Id} \nrestriction (\rfs{cl}(X) - \rfs{Dom}(\hath))$.
We show that $\barh \in H_{\itGamma}(\rfs{cl}(X))$.
That is, $\barh \in H(\rfs{cl}(X))$
and $\barh$ is $\itGamma$-bicontinuous.
Let $u \in \rfs{cl}(X)$.
If $u \in \rfs{Dom}(\hath)$,
then since $\rfs{Dom}(\hath)$ is open in $\rfs{cl}(X)$
and $\hath$ is continuous, we have that $\barh$ is continuous at $u$.
If $u \not\in \rfs{Dom}(\hath)$,
then since $\rfs{cl}^E(\rfs{supp}(\hath)) \subseteq \rfs{Dom}(\hath)$,
it follows that $u \in \rfs{cl}(X) - \rfs{cl}^E(\rfs{supp}(\hath))$.
So there is $V \in \rfs{Nbr}^{\srfs{cl}(X)}(u)$
such that $\barh \nrestriction V = \rfs{Id}$.
Hence $\barh$ is continuous at $u$.
The same argument applies to $\barh\inverse$.
So $\barh \in H(\rfs{cl}(X))$.

We now show that $\barh$ is $\itGamma$-bicontinuous.
Recall that $X - \rfs{Rng}(\varphi) \neq \emptyset$
and hence $X -\rfs{Dom}(\hath) \neq \emptyset$.
Since $\varphi$ is $\gamma$-continuous,
it follows that $\rfs{Dom}(\hath)$ and hence
$\rfs{supp}(\hath)$ are bounded.
Set $c = d(\rfs{supp}(\hath),\rfs{cl}(X) - \rfs{Dom}(\hath))$
and $e = \rfs{diam}(\rfs{supp}(\hath))$.
Clearly, $e < \infty$.
We show that $c > 0$.
Recall that $\rfs{supp}(g) \subseteq U$,
and hence
$\rfs{supp}(\hath)  = \varphi(\rfs{supp}(g)) \subseteq
\varphi(U)$. Also, $x = \varphi(0) \in \varphi(U)$.
So
\vspace{1.5mm}
\newline
\rule{7pt}{0pt}
\renewcommand{\arraystretch}{1.5}
\addtolength{\arraycolsep}{-0pt}
$
\begin{array}{lll}
c
&
=
&
d(\rfs{supp}(\hath),\rfs{cl}(X) - \rfs{Dom}(\hath)) \geq
d(\varphi(U),\rfs{cl}(X) - \varphi(\swhatX))
\\
&
\geq
&
d(x,\rfs{cl}(X) - \varphi(\swhatX)) - \rfs{diam}(\varphi(U))
\\
&
\geq
&
d(x,\rfs{cl}(X) - \varphi(\swhatX)) -
\dgfrac{d(x,\rfs{cl}(X) - \swhatX)}{2}
=
\dgfrac{d(x,\rfs{cl}(X) - \swhatX)}{2} > 0.
\vspace{1.7mm}
\end{array}
$
\renewcommand{\arraystretch}{1.0}
\addtolength{\arraycolsep}{0pt}
\newline

Let $u,v \in \rfs{cl}(X)$.
If $u,v \in \rfs{supp}(\hath)$,
then $\norm{\barh(u) - \barh(v)} \leq \rho(\norm{u - v})$.
If $u,v \in \rfs{cl}(X) - \rfs{supp}(\hath)$,
then $\norm{\barh(u) - \barh(v)} = \norm{u - v}$.
Suppose that $u \in \rfs{supp}(\hath)$ and
$v \not\in \rfs{supp}(\hath)$.
If $v \in \rfs{Dom}(\hath)$,
then $\norm{\barh(u) - \barh(v)} \leq \rho(\norm{u - v})$.
Otherwise,
\vspace{1.5mm}
\newline
\rule{7pt}{0pt}
\renewcommand{\arraystretch}{1.5}
\addtolength{\arraycolsep}{-3pt}
$
\begin{array}{ll}
&
\norm{\barh(u) - \barh(v)} \leq
\norm{\barh(u) - u} + \norm{u - v} \leq e + \norm{u - v}
\\
=
\rule{5pt}{0pt}
&
\frac{e}{c} \mcdot c + \norm{u - v}
\leq 
\frac{e}{c} \mcdot \norm{u - v} + \norm{u - v}
= \frac{e + c}{c} \mcdot \norm{u - v}.
\vspace{1.7mm}
\end{array}
$
\renewcommand{\arraystretch}{1.0}
\addtolength{\arraycolsep}{3pt}
\newline
It follows that
$\barh$ is $(1 + \dgfrac{e}{c}) \mcdot \rho$-continuous.
The same argument applies to $\barh\inverse$.
Since
$(1 + \dgfrac{e}{c}) \mcdot \rho \in \itGamma$, it follows
that $\barh$ is $\itGamma$-bicontinuous.

Let $h = \barh \nrestriction X$.
Then $\rfs{supp}(h) \subseteq Z \cap X$.
Hence
$h \in H^{\srfs{CMP.LC}}_{\itGamma}(X)\sprt{Z \cap X}$.
It follows that $h \in G$.
By Assumption~(3) in the statement of the lemma,
$h^{\tau} \in H$.
So $h^{\tau} \in \rfs{EXT}(Y)$ and $h^{\tau}$
is $\itDelta$-bicontinuous at $y$.
That is, for some $\nu \in \itDelta$,
$h^{\tau}$ is $\nu$-bicontinuous at $y$.
So $(h^{\tau})^{\srfs{cl}}$ is $\nu$-bicontinuous at $y$.
Now, $\barh = h^{\srfs{cl}}$,
hence $\barh^{\tau^{\trfs{cl}}} = (h^{\tau})^{\srfs{cl}}$
and so $\barh^{\tau^{\trfs{cl}}}$
is $\nu$-bicontinuous at $y$.
Recall that $\psi$ is $\delta$-bicontinuous, where $\delta \in \Delta$.
Also, $\psi\inverse(y) = 0^F$.
It follows that
$(\barh^{\tau^{\trfs{cl}}})^{\psi\inverse}$
is $\delta \scirc \nu \scirc \delta$-bicontinuous at $0^F$.
That is, $(\barh^{\tau^{\trfs{cl}}})^{\psi\inverse}$
is $\itDelta$-bicontinuous at $0^F$.
Finally,
$g^{\varphi} = \hath \subseteq \barh$
and $y \in \rfs{Dom}((g^{\varphi})^{\tau^{\trfs{cl}}})$.
So
$((g^{\varphi})^{\tau^{\trfs{cl}}})^{\psi\inverse}$
is $\itDelta$-bicontinuous at $0^F$.
That is, $g^{\hattau}$ is $\itDelta$-bicontinuous at $\hattau(0^E)$.

We have checked that the conditions of Lemma
\ref{ams-bddly-lip-bldr-l5.7} hold.
So $\hattau$ is inversely $\pair{\itDelta}{\calP}$-continuous at $0$.
\smallskip\hfill\myqed

\subsection{$\itGamma$-continuity for submerged pairs
\hbox{and the star operation.}}\label{s12.2}
The next intermediate goal is to show that in the above setting,
$\hattau$ is inversely $\itDelta$-continuous at $0^E$,
(Lemma \ref{ams-bddly-lip-bldr-l5.20}(b)).
Unfortunately, we know to prove this only under additional
asumptions on $\itGamma$ and $\itDelta$.
The assumptions $\itGamma = \itDelta$ and $\itGamma$ is principal
suffice. (See Clause M6 in Definition \ref{nn1.3}).
The exact extra assumptions use the notion of star-closedness which
is defined in Definition \ref{ams-bddly-lip-bldr-d5.14}(d).
They are:
$\itGamma \subseteq \itDelta$\ \ and\ \ $\itDelta$ is
$\itGamma$-star-closed.

\begin{prop}\label{ams-bddly-lip-bldr-p5.10}
Recall that for $\rho \in H([0,\infty))$ and a normed space $E$,
the homeomorphism
$\rfs{Rad}_{\rho}^E \in H(E)$ was defined as follows:
for $u \neq 0$,
$\rfs{Rad}_{\rho}^E(u) = \rho(\norm{u}) \ncdot \frac{u}{\norm{u}}$
and $\rfs{Rad}_{\rho}^E(0) = 0$.
If $\alpha \in \rfs{MC}$ and $\rho$ is $\alpha$-continuous,
then $\rfs{Rad}_{\rho}^E$ is $5 \mcdot \alpha$-continuous,
\end{prop}

\noindent
{\bf Proof } Let $x,y \in E$ and $y \neq 0$.
Denote $z = \norm{x} \ncdot \frac{y}{\norm{y}}$.
Then $\norm{y - z} = \abs{\norm{y} - \norm{x}} \leq \norm{y - x}$.
So
$\norm{x - z} \leq \norm{x - y} + \norm{y - z} \leq 2 \norm{y - x}$.
Let $h = \rfs{Rad}_{\rho}^E$. Suppose that $x \neq z$. Then
\vspace{1.5mm}
\newline
\rule{7pt}{0pt}
\renewcommand{\arraystretch}{1.5}
\addtolength{\arraycolsep}{-3pt}
$
\begin{array}{ll}
&
\norm{h(y) - h(x)} \leq \norm{h(y) - h(z)} + \norm{h(z) - h(x)} \leq
\alpha(\norm{y - z}) +
\frac{\rho(\norm{x})}{\norm{x}} \ncdot \norm{x - z}
\\
\leq
\rule{5pt}{0pt}
&
\alpha(\norm{y - x}) +
\frac{\alpha(\norm{x})}{\norm{x}} \ncdot \norm{x - z} \leq
\alpha(\norm{y - x}) +
\frac{\alpha(\norm{\frac{x - z}{2}})}{\norm{\frac{x - z}{2}}} \ncdot
\norm{x - z}
\\
=
\rule{5pt}{0pt}
&
\alpha(\norm{y - x}) +
2 \alpha(\norm{\frac{x - z}{2}}) \leq
\alpha(\norm{y - x}) +
2 \alpha(\norm{x - z}) \leq
\alpha(\norm{y - x}) + 2 \alpha(2 \norm{y - x})
\\
\leq
\rule{5pt}{0pt}
&
\alpha(\norm{y - x}) + 4 \alpha(\norm{y - x}) =
5 \alpha(\norm{y - x}).
\vspace{0.0mm}
\end{array}
$
\renewcommand{\arraystretch}{1.0}
\addtolength{\arraycolsep}{3pt}
\newline
If $x = z$, then $\norm{h(y) - h(x)} \leq \alpha(\norm{y - x})$.
So $\rfs{Rad}_{\rho}^E$ is $5 \mcdot \alpha$-continuous.
\hfill\myqed

\begin{prop}\label{ams-bddly-lip-bldr-p5.11}
There is $M^{\srfs{rtn}}$ such that the following holds.
Let $\alpha \in \rfs{MBC}$
and $a > 0$.
Let $E$ be a normed space, $x,y \in E$
and $\norm{x} = \norm{y} = \alpha(a)$.
Then there is $g \in H(E)$ such that
$g(0) = 0$, $g(x) = y$,
$\rfs{supp}(g) \subseteq B(0,\alpha(a) + \dgfrac{a}{2})$,
and $g$ is $M^{\srfs{rtn}} \ncdot \alpha \scirc \alpha$-bicontinuous.
   \index{N@mrtn@@$M^{\srfs{rtn}}$}
\end{prop}

\noindent
{\bf Proof } Let $b = \alpha(a)$,
$c = \alpha(a) + \dgfrac{a}{2}$
and $N = M^{\srfs{hlb}}$.
(See Proposition \ref{p-bddlip-bldr-1.9}(c)).
Let $\rho \in H([0,\infty))$ be the piecewise linear function
with breakpoints at $b$ and $c$ such that
$\rho(b) = \frac{b}{2N}$ and $\rho(t) = t$ for every $t \geq c$.
The slope of $\rho$ on $[0,\alpha(a)]$ is $\frac{1}{2N} < 1$.
\smallskip
The slope of $\rho$ on $[\alpha(a), \alpha(a) + \dgfrac{a}{2}]$
is
$\frac{c - \frac{b}{2N}}{\dgfrac{\alpha(a)}{2}} =
\frac{2 \alpha(a) + a - \dgfrac{\alpha(a)}{N}}{a} \leq
\frac{3 \alpha(a)}{a}$.
The slope of $\rho$ on $[\alpha(a) + \dgfrac{\alpha}{2},\infty)$ is $1$.
So $\rho$ is $(3,3\alpha)$-continuous.
(See Definition \ref{d-bddly-lip-2.9-1}(b)).
By Proposition~\ref{p-bddlip-1.11-1}(a), $\rho$ is
$9 \alpha$-continuous.
By Proposition~\ref{ams-bddly-lip-bldr-p5.10},
$\rfs{Rad}_{\rho}^E$ is $45 \ncdot \alpha$-continuous.
Clearly, $(\rfs{Rad}_{\rho}^E)\inverse = \rfs{Rad}_{\rho\inverse}^E$.
The slope of $\rho\inverse$ on
$[0,\frac{\alpha(a)}{2N}]$ is $2N$.
The slope of $\rho\inverse$ on
$[\frac{\alpha(a)}{2N},\alpha(a) + \frac{\alpha}{2}]$
is  $\leq \frac{a}{\alpha(a)} \leq 1$.
So $\rho\inverse$ is $(3,2N \alpha)$-continuous.
It  follows that $(\rfs{Rad}_{\rho}^E)\inverse$ is
$30N \ncdot \alpha$-continuous.
Let $M_1 = \max(30N,45)$. Then $\rfs{Rad}_{\rho}^E$
is $M_1 \ncdot \alpha$-bicontinuous.
Let $h = \rfs{Rad}_{\rho}^E$.
Then
\begin{itemize}
\addtolength{\parskip}{-11pt}
\addtolength{\itemsep}{06pt}
\item[(1)] $\rfs{supp}(h) \subseteq B(0,\alpha(a) + \frac{a}{2})$,
\item[(2)] $h(x) = \frac{x}{2N}$,
\item[(3)] $h$ is $M_1 \ncdot \alpha$-bicontinuous.
\vspace{-05.7pt}
\end{itemize}

Let $L = \rfs{span}(\dbltn{x}{y})$.
By Proposition \ref{p-bddlip-bldr-1.9}(c),
there are a Euclidean norm
$\norm{\ }^{\srbfs{H}}$ on $L$ and a complement $S$ of $L$ such that
for every $u \in E$,
$\norm{(u)_L}^{\srbfs{H}} + \norm{(u)_S} \approx^{M^{\ssrfs{hlb}}}
\norm{u}$.
Denote $\boldnorm{u} = \norm{(u)_L}^{\srbfs{H}} + \norm{(u)_S}$.
We shall apply Proposition \ref{p-gamma.7}(c).
Let $\hatx = \frac{x}{2N}$,
$\haty = \frac{\sboldnorm{\hatx}}{\sboldnorm{y}} y$,
and $\theta$ be the angle from $\hatx$ to $\haty$
in $\pair{L}{\norm{\ }^{\srbfs{H}}}$.
So $\boldnorm{\haty} = \boldnorm{\hatx}$.
Let $S = \barB^L(0,\boldnorm{\hatx})$.
Let $\eta$ be the piecewise linear function with breakpoint at
$\boldnorm{\hatx}$ such that
$\eta(0) = \theta$ and $\eta(\boldnorm{\hatx}) = 0$.
So $\eta$ is $\frac{\theta}{\sboldnorm{\hatx}}$-Lipschitz.
Hence the conditions of Proposition \ref{p-bddlip-bldr-1.9}(c) hold
with $r = \boldnorm{\hatx}$
and $K = \dgfrac{\theta}{\raise1.4pt\hbox{$\sboldnorm{\hatx}$}}$.
Let $\bard$ denote the distance function obtained
from $\boldnorm{\ }$.
Let $g_1$ be defined by
$g_1(u) = \rfs{Rot}^{F,H}_{\eta(\bard(u,S))}(u)$.
Then $g_1 \in H(E)$ and $g_1$ is
$(M^{\srfs{rot}} \kern-3pt\cdot\kern-2.5pt Kr + 1)$-bilipschitz
with respect to $\bard$.
Note that $Kr = \theta \leq \pi$.
So $g_1$ is $M^{\srfs{rot}} (\pi + 1)$-bilipschitz
with respect to $\bard$.
Write
$M_2 = (M^{\srfs{hlb}})^2 M^{\srfs{rot}} (\pi + 1)$.
Hence
\begin{itemize}
\addtolength{\parskip}{-11pt}
\addtolength{\itemsep}{06pt}
\item[(4)] 
$g_1$ is $M_2$-bilipschitz in $\pair{E}{\norm{\ }}$.
\vspace{-05.7pt}
\end{itemize}

Let $u \in E - B(0,\norm{x})$.
Then $\boldnorm{u} \geq \dgfrac{\norm{u}}{M^{\ssrfs{hlb}}} \geq
\dgfrac{\norm{x}}{N}$.
So $\bard(u,S) \geq \boldnorm{\hatx}$. Hence $g_1(u) = u$.
That is,
\begin{itemize}
\addtolength{\parskip}{-11pt}
\addtolength{\itemsep}{06pt}
\item[(5)] 
$\rfs{supp}(g_1) \subseteq B(0,\norm{x})$.
\vspace{-05.7pt}
\end{itemize}
It is also obvious that
\begin{itemize}
\addtolength{\parskip}{-11pt}
\addtolength{\itemsep}{06pt}
\item[(6)] 
$g_1(\hatx) = \haty$.
\vspace{-05.7pt}
\end{itemize}

Let $\bary = \frac{y}{2N}$. Then $\norm{\bary} = \norm{\hatx}$.
Recall that $\boldnorm{\haty} = \boldnorm{\hatx}$.
Since $\boldnorm{\hatx} \approx^{M^{\ssrfs{hlb}}} \norm{\hatx}$, 
\newline
$\norm{\haty} \approx^{M^{\ssrfs{hlb}}} \norm{\bary}$.
That is,
$\frac{1}{M^{\ssrfs{hlb}}} \ncdot \norm{\haty} \leq \norm{\bary} \leq
M^{\srfs{hlb}} \ncdot \norm{\haty}$.
We construct $g_2$ which takes $\haty$ to $\bary$.
Let $\fnn{\rho}{[0,\infty)}{[0,\infty)}$
be the piecewise linear function with breakpoints
$\norm{\haty}$ and $\norm{x}$ such that $\rho(0) = 0$,
$\rho(\norm{\haty}) = \norm{\bary}$
and $\rho(t) = t$ for every $t \geq \norm{x}$.
Since $\norm{\haty},\norm{\bary} < \norm{x}$,
$\rho \in H([0,\infty))$.
The slopes of $\rho$ are $\frac{\norm{\bary}}{\norm{\haty}}$,
$\frac{\norm{x} - \norm{\bary}}{\norm{x} - \norm{\haty}}$ and $1$,
and the slopes of $\rho\inverse$ are
$\frac{\norm{\haty}}{\norm{\bary}}$,
$\frac{\norm{x} - \norm{\haty}}{\norm{x} - \norm{\bary}}$ and $1$.
Clearly,
\rule{0pt}{12pt}
$\frac{\norm{\bary}}{\norm{\haty}} \leq M^{\srfs{hlb}} = N$.
Note that
$\norm{\haty} \leq \norm{\haty}^{\srbfs{H}} =
\boldnorm{\haty} = \boldnorm{\hatx} \leq 
M^{\srfs{hlb}} \ncdot \norm{\hatx} =
N \ncdot \frac{\norm{x}}{2N} = \dgfrac{\norm{x}}{2}$.
So
$\frac{\norm{x} - \norm{\bary}}{\norm{x} - \norm{\haty}} =
\frac{(1 - \frac{1}{2N}) \norm{x}}{\norm{x} - \norm{\haty}} \leq
\frac{\norm{x}}{\norm{x} - \dgfrac{\norm{x}}{2}} = 2$.
Hence $\rho$ is $\max(N,2)$-Lipschitz.
\smallskip

As to the slopes of $\rho\inverse$, clearly,
$\frac{\norm{\haty}}{\norm{\bary}} \leq N$ and
$\frac{\norm{x} - \norm{\haty}}{\norm{x} - \norm{\bary}} \leq
\frac{\norm{x}}{(1 - \frac{1}{2N}) \norm{x}} \leq 2$.
So $\rho\inverse$ is $\max(N,2)$-Lipschitz.
Let $M_3 = 3 \max(N,2)$ and $g_2 = \rfs{Rad}_{\rho}^E$.
By Proposition \ref{metr-bldr-p3.18},
\begin{itemize}
\addtolength{\parskip}{-11pt}
\addtolength{\itemsep}{06pt}
\item[(7)] 
$g_2$ is $M_3$-bilipschitz.
\vspace{-05.7pt}
\end{itemize}
It follows trivially from the definitions of $\rho$ and $g_2$ that
\begin{itemize}
\addtolength{\parskip}{-11pt}
\addtolength{\itemsep}{06pt}
\item[(8)] 
$g_2(\haty) = \bary$,
\item[(9)] $\rfs{supp}(g_2) \subseteq B(0,\norm{x})$.
\vspace{-05.7pt}
\end{itemize}
Let $g = h\inverse \scirc g_2 \scirc g_1 \scirc h$.
Note that
\begin{itemize}
\addtolength{\parskip}{-11pt}
\addtolength{\itemsep}{06pt}
\item[(10)] 
$h\inverse(\bary) = h\inverse(\frac{y}{2N}) = y$.
\vspace{-05.7pt}
\end{itemize}
It follows from (1)\,-\,(10) that
$g$ is $M_1^2 M_2 M_3 \ncdot \alpha \scirc \alpha$-bicontinuous,
$g(x) = y$ and
$\rfs{supp}(g) \subseteq B(0, \alpha(a) + \dgfrac{a}{2})$.
Define $M^{\srfs{rtn}} = M_1^2 M_2 M_3$. Then $M^{\srfs{rtn}}$ is as
required.
\hfill\myqed

\begin{defn}\label{ams-bddly-lip-bldr-d5.12}
\begin{rm}
(a) Let $E$ be a metric space,
$x,y \in X \subseteq E$ and $\alpha \in \rfs{MC}$.
We say that $\pair{x}{y}$ is {\it $\alpha$-submerged}
in $X$ with respect to $E$,
if $\delta^X(x) \geq \norm{x - y} + \alpha\inverse(\norm{x - y})$.

   \index{submerged. $\pair{x}{y}$ is $\alpha$-submerged in $X$.
          This means
          $\delta^X(x) \geq \norm{x - y} +
	  \alpha\inverse(\norm{x - y})$}

{\thickmuskip=3mu \medmuskip=2mu \thinmuskip=1mu 
(b) Let $X \subseteq E$, $Y \subseteq F$ be open subsets of the metric
spaces $E,F$,
$V \subseteq X$, $x \in \rfs{bd}(X)$, $\alpha,\beta \in \rfs{MC}$,
$\itGamma,\itDelta \subseteq \rfs{MC}$
and $\tau \in \rfs{EXT}^{\pm}(X,Y)$.
We say that $\tau$ is
{\it $\beta$-continuous for $\alpha$-submerged pairs}
in $V$, if for every
$\alpha$-submerged pair $\pair{y}{z}$ in $V$,
$d^Y(\tau(y),\tau(z)) \leq \beta(d^X(y,z))$.}

We say that $\tau$ is
{\it $\beta$-continuous for $\alpha$-submerged pairs}
at $x$, ($\tau$ is $\scolonpair{\beta}{\alpha}$-continuous
at $x$),
if there is $U \in \rfs{Nbr}^E(x)$
such that $\tau$ is $\beta$-continuous for $\alpha$-submerged pairs
in $U \cap X$.
We say that $\tau$ is
{\it $\itDelta$-continuous for
$\itGamma$-submerged pairs at $x$,}
($\tau$ is
$\scolonpair{\itDelta}{\itGamma}$-continuous at~$x$),
if for every $\alpha \in \itGamma$ there is $\beta \in \itDelta$
such that $\tau$ is $\scolonpair{\beta}{\alpha}$-continuous at $x$.

   \index{continuous. $\beta$-continuous for $\alpha$-submerged pairs.
          Abbreviation: $\scolonpair{\beta}{\alpha}$-continuous}
   \index{continuous. $\itDelta$-continuous for
          $\itGamma$-submerged pairs.
       Abbreviation: $\scolonpair{\itDelta}{\itGamma}$-continuous}

(c) Let $X \subseteq E$, $Y \subseteq F$ be open subsets of the metric
spaces $E,F$,
$V \subseteq X$, \hbox{$\alpha,\beta \in \rfs{MC}$}
and $\tau \in \rfs{EXT}^{\pm}(X,Y)$.
We say that $\tau$ {\it is almost $\beta$-continuous for
$\alpha$-submerged pairs} in~$V$,
($\tau$ is $\scolonpair{\beta}{\alpha}$-almost-continuous
in $V$),
if for every $\alpha$-submerged pairs $\pair{y}{z_1},
\pair{y}{z_2}$\break
in $V$: if $d(y,z_1) = d(y,z_2)$, then
$d^Y(\tau(y),\tau(z_2)) \leq \beta(d^Y(\tau(y),\tau(z_1)))$.
\hfill\proofend
   \index{almost $\beta$-continuous for $\alpha$-submerged pairs.
          Abbreviation: $\scolonpair{\beta}{\alpha}$-almost-continuous}
\end{rm}
\end{defn}

Under assumptions similar to Lemma \ref{ams-bddly-lip-bldr-l5.9},
we prove the submerged continuity of $\tau\inverse$.

\begin{lemma}\label{ams-bddly-lip-bldr-l5.13}
Assume the following facts.
\begin{itemize}
\addtolength{\parskip}{-11pt}
\addtolength{\itemsep}{06pt}
\item[\num{1}]
$\itGamma,\itSigma$ are countably generated moduli of continuity,
and $\itOmega$ is the modulus of continuity
\indent\kern3pt generated by
$\itGamma \cup \itSigma$.
\item[\num{2}]
$X \subseteq E$ and $Y \subseteq F$ are open subsets of
the normed spaces $E$ and $F$,
$X$ is $\itGamma$-LIN-
\indent\kern3pt bordered and
$Y$ is $\itSigma$-LIN-bordered.
\item[\num{3}]
$\tau \in \rfs{EXT}^{\pm}(X,Y)$,
$G$ is a $\itGamma$-appropriate subgroup of $\rfs{EXT}(X)$,
$H$ is a $\itDelta$-appropriate
\indent\kern2pt subgroup of $\rfs{EXT}(Y)$
and $G^{\tau} = H$.
\item[\num{4}]
$x \in \rfs{bd}(X)$,
$\trpl{\varphi}{A}{r}$ is a boundary chart element for~$x$,
$\gamma \in \itGamma$
and $\varphi$ is $\gamma$-bicontinuous.
\item[\num{5}]
$y \in \rfs{bd}(Y)$,
$\trpl{\psi}{B}{s}$ is a boundary chart element for~$y$,
$\sigma \in \itSigma$
and $\psi$ is $\sigma$-bicontinuous.
\item[\num{6}]
$\tau^{\srfs{cl}}(x) = y$
and
$\tau(\varphi(\rfs{BCD}^{E}(A,r))) \subseteq \psi(\rfs{BCD}^{F}(B,s))$.
\item[\num{7}]
Set
$\wtildeX = \rfs{BCD}^{E}(A,r)$,
$\tildetau = \psi\inverse \scirc \tau \scirc \varphi$
and
$\wtildeY = \tildetau(\wtildeX)$.
\vspace{-05.7pt}
\end{itemize}
\underline{Then}
$\tildetau\inverse$ is $\scolonpair{\itOmega}{\itSigma}$-continuous at
$\tildetau(0)$.
\end{lemma}

\noindent
{\bf Proof }
There is
$Z \in \rfs{Nbr}^F(y)$
such that
$H\sprt{Z \cap Y} \supseteq
H^{\srfs{CMP.LC}}_{\itSigma}(Y)\sprt{Z \cap Y}$,
and we may assume that $\psi(\rfs{BCD}^{F}(B,s)) \subseteq Z$.
Set
$L = \rfs{bd}(A)$,
$\swhatX = \rfs{BCD}^{E}(A,r) \cup (L \cap B^E(0,r))$,
$\hattau = \psi\inverse \scirc \tau^{\srfs{cl}} \scirc \varphi$,
$\whatY = \hattau(\swhatX)$
and $\calP = \setm{(v + L) \cap \swhatX}{v \in \swhatX}$.
Note that $\hattau = \tildetau^{\srfs{cl}}_{B^E(0,r),B^F(0,s)}$.
By Lemma \ref{ams-bddly-lip-bldr-l5.9},
$\hattau$ is inversely $\pair{\itSigma}{\calP}$-continuous at~$0$.
Let $r_0 \in (0,r)$ and $\sigma \in \itSigma$ be such that
$\hattau \nrestriction (B^E(0,r_0) \cap \swhatX)$ is inversely
$\pair{\sigma}{\calP}$-continuous.
Let $L_0 \subseteq L$ be any ray whose endpoint is $0$.
For every
$u \in B(0,r_0) \cap \swhatX$ let $x_u$ be the intersection
point of the ray \hbox{$u + L_0$}
with the sphere $S(0,r_0)$.
Clearly, $\lim_{u \rightarrow 0} x_u = x_{0^E}$.
So
$\lim_{u \rightarrow 0} d^F(\hattau(u),\hattau(x_u)) =
d^F(\hattau(0^E),\hattau(x_{0^E})) > 0$.
Also,
$\lim_{u \rightarrow 0}^{\wtildeY} \delta^{\wtildeY}(\tildetau(u)) = 0$.
Hence there is $r_1 \in (0,r_0)$
such that for every $u \in B(0,r_1) \cap \swtildeX$,
$d^F(\tildetau(u),\tildetau(x_u)) > \delta^{\wtildeY}(\hattau(u))$.
Let $V = \tildetau(B(0,r_1) \cap \swtildeX)$.
So for every $v \in V$ and $t \in [0,\delta^{\wtildeY}(v)]$
there is
$y(v,t) \in
\tildetau([\tildetau\inverse(v),x_{\tildetau\inverse(v)}])$
such that $d^F(y(v,t),v) = t$.
Denote $\tildetau\inverse$ by $\tildeeta$.
By the inverse $\pair{\sigma}{\calP}$-bicontinuity of $\hattau$,
for every $v$ and $t$ as above
$d^E(\tildeeta(y(v,t)),\tildeeta(v)) \leq \sigma(d^F(y(v,t),v))$.

{\bf Claim 1}
Let $\alpha \in \itSigma \cap \rfs{MBC}$. Then there are
$W \in \rfs{Nbr}^F(0)$ and $\gamma \in \itGamma$ such that
$\tildeeta$
is $\scolonpair{\gamma}{\alpha}$-almost-continuous in
$W \cap \wtildeY$.
{\bf Proof }
Suppose by  contradiction this is not so.
Let $\setm{\gamma_i}{i \in \bbN}$ be a generating set for $\itGamma$,
and assume that for every $i$, $\setm{j}{\gamma_j = \gamma_i}$ is
infinite.
There is a sequence $\setm{\trpl{y_i}{u_i}{v_i}}{i \in \bbN}$
such that: (i) for every $i$, $\pair{y_i}{u_i}$ is\break
$\alpha$-submerged
in $\wtildeY$ and $\norm{u_i - y_i} = \norm{v_i - y_i}$;
(ii) $\lim_i y_i = 0^F$;
{\thickmuskip=2.1mu \medmuskip=1mu \thinmuskip=1mu 
(iii) $\delta^{\wtildeY}(y_{i + 1}) <
\dgfrac{\alpha\inverse(\norm{y_i - u_i})}{4}$;}
(iv) $\norm{\tildeeta(v_i) - \tildeeta(y_i)} >
\gamma_i(\norm{\hateta(u_i) - \hateta(y_i)})$.
Let
$r_i = \norm{u_i - y_i} + \frac{\alpha\inverse(\norm{u_i - y_i})}{2}$.
Note that from (iii) and the fact that $\pair{y_i}{u_i}$
is $\alpha$-submerged it follows that
$B(y_i,r_i) \cap B(y_i,r_j) = \emptyset$ for any $i \neq j$.
By Proposition \ref{ams-bddly-lip-bldr-p5.11},
there is $g_i \in H(\wtildeY)$ such that
$g_i(y_i) = y_i$, $g(u_i) = v_i$,
$\rfs{supp}(g_i) \subseteq B(y_i,r_i)$,
and $g_i$ is $M^{\srfs{rtn}} \ncdot \alpha \scirc \alpha$-bicontinuous.
Since $\rfs{supp}(g_i) \cap \rfs{supp}(g_j) = \emptyset$
for any $i \neq j$,
we have that $\tildeg = \bcirc_i g_i$ is well-defined,
and $\tildeg$
is $(M^{\srfs{rtn}})^2 \ncdot \alpha^{\sscirc 4}$-bicontinuous.
We shall reach a contradiction by showing
that $\tildeg$ is $\itSigma$-bicontinuous at $0^F$, whereas
$\tildeg^{\tildetau\inverse}$ is not $\itGamma$-bicontinuous at $0^E$.

Define $\tildeh = \tildeg^{\psi}$
and $h = \tildeh \cup \rfs{Id} \nrestriction (Y - \psi(\wtildeY))$.
We shall show that $h \in H$.
Recall that $y = \psi(0^F)$ and set $h_i = g_i^{\psi}$.
Then $\tildeh = \bcirc_{i \in \bbN} h_i$.
Recall that $\rfs{supp}(g_i) \subseteq B^F(y_i,r_i)$
and note that $\lim_{i \in \bbN} \overB^F(y_i,r_i) = 0^F$.
Since
$\sngltn{0^F} \cup \bigcup_{i \in \bbN} \overB^F(y_i,r_i) \subseteq
\rfs{Dom}(\psi)$, it follows that
$\lim_{i \in \bbN} \psi(\overB^F(y_i,r_i)) = y$.
Also, $\rfs{supp}(h_i) = \psi(\rfs{supp}(g_i))$.
Hence
$\rfs{cl}(\rfs{supp}(h_i)) =
\psi(\rfs{cl}(\rfs{supp}(g_i)))  \subseteq \psi(\overB^F(y_i,r_i))$
and so $\lim_{i \in \bbN} \rfs{cl}(\rfs{supp}(h_i)) = y$.
We thus conclude that:
(1)
$\rfs{cl}(\rfs{supp}(\tildeh)) = \sngltn{y} \cup
\bigcup_{i \in \bbN\kern1pt} \rfs{cl}(\rfs{supp}(h_i))$.
It also follows that:
(2) if $\vecz \subseteq Y$ and $\lim \vecz = y$,
then $\lim h(\vecz) = y$.
Note that: (3) for every $i \in \bbN$,
$\rfs{cl}(\rfs{supp}(h_i)) \subseteq \psi(\overB^F(y_i,r_i)) \subseteq
\psi(\wtildeY)$.
Let $z \in \rfs{cl}(Y)$.
If $z \not\in \rfs{cl}(\rfs{supp}(h))$,
then $h \cup \sngltn{\pair{z}{z}}$ is continuous.
If $z \in \rfs{cl}(\rfs{supp}(h))$,
then $z \in \rfs{cl}(\rfs{supp}(\tildeh))$.
So by (1) and (3), either $z = y$ or $z \in \psi(\wtildeY)$.
If $z = y$, then by (2), $h \cup \sngltn{\pair{z}{z}}$ is continuous.
If $z \in \psi(\wtildeY)$, then $h(z) = \tildeh(z)$.
From the facts: $\tildeh$ is continuous,
$h \nrestriction \wtildeY = \tildeh$
and $\psi(\wtildeY)$ is open in $F$,
it follows that $h$ is continuous at $z$.
We have shown that $h$ is extendible in $F$.
The same argument applies to $h\inverse$, so $h \in \rfs{EXT}(Y)$.
Clearly,
$\rfs{supp}(h) = \rfs{supp}(\tildeh) \subseteq \psi(\wtildeY) \subseteq 
\psi(\rfs{BCD}^{F}(B,s)) \subseteq Z$.
That is, (4) $\rfs{supp}(h) \subseteq Z$.

We now show that $h \in H^{\srfs{CMP.LC}}_{\itSigma}(Y)$.
Write $\baralpha = (M^{\srfs{rtn}})^2 \ncdot \alpha^{\sscirc 4}$
and $\beta = \sigma \scirc \baralpha \scirc \sigma$.
Then $\beta \in \itSigma$.
We have seen that $\tildeg$ is $\baralpha$-bicontinuous.
So since $\psi$ is $\sigma$-bicontinuous, it follows that
$\tildeh$ is $\beta$-bicontinuous.
This implies that $\tildeh^{\srfs{cl}}$ is $\beta$-bicontinuous.
We show that for every $z \in \rfs{cl}(Y)$,
$h$ is $\beta$-bicontinuous at $z$.
This is certainly true if $z \not\in \rfs{cl}(\rfs{supp}(h))$.
So suppose that
$z \in \rfs{cl}(\rfs{supp}(h))$.
Then $z \in \rfs{cl}(\rfs{supp}(\tildeh))$.
By (1) and (3), either $z \in \psi(\wtildeY)$ or $z = y$.
If $z \in \psi(\wtildeY)$,
then $\psi(\wtildeY) \in \rfs{Nbr}^F(z)$
and
$h \nrestriction \psi(\wtildeY) = \tildeh \nrestriction \psi(\wtildeY)$.
So $h$ is $\beta$-bicontinuous at $z$.

Assume that $z = y$.
Recall that $x = \varphi(0^E)$ and $y = \psi(0^F)$
and define $X_0 = X \cup \sngltn{x}$ and $Y_0 = Y \cup \sngltn{y}$.
Note that $\psi(\wtildeY) = \tau(\varphi(\rfs{BCD}^{E}(A,r)))$.
Since
$\varphi(\rfs{BCD}^{E}(A,r)) = \varphi(B^E(0,r)) \cap X$
and $\varphi(B^E(0,r))$ is open in $E$,
it follows that
$\varphi(\rfs{BCD}^{E}(A,r)) \cup \sngltn{x} \in
\rfs{Nbr}^{X_0}(x)$.
From the fact that $\tau \in \rfs{EXT}^{\pm}(X,Y)$
it follows that
$\tau(\varphi(\rfs{BCD}^{E}(A,r))) \cup \sngltn{y} \in
\rfs{Nbr}^{Y_0}(y)$.
That is, $\psi(\wtildeY) \cup \sngltn{y} \in \rfs{Nbr}^{Y_0}(y)$.
So there is $W \in \rfs{Nbr}^F(y)$
such that $W \cap Y = \psi(\wtildeY)$.
Thus $h \nrestriction W = \tildeh \nrestriction W$.
It follows that $h$ is $\beta$-bicontinuous at $y$.
So $h \in H^{\srfs{CMP.LC}}_{\itSigma}(Y)$.
By (4), $h \in H^{\srfs{CMP.LC}}_{\itSigma}(Y)\sprt{Z \cap Y}$.
Also recall that
$H\sprt{Z \cap Y} \supseteq
H^{\srfs{CMP.LC}}_{\itSigma}(Y)\sprt{Z \cap Y}$.
So $h \in H$.

We conclude that $h^{\tau\inverse} \in G$.
Now, $G$ is of boundary type $\itGamma$,
so $h^{\tau\inverse}$ is $\itGamma$-bicontinuous at~$x$.
Since $\varphi$ is $\itGamma$-bicontinuous and $\varphi(0^E) = x$,
we have that
$(h^{\tau\inverse})^{\varphi\inverse}$ is $\itGamma$-bicontinuous
at $0^E$.
The following steps show that
$(h^{\tau\inverse})^{\varphi\inverse} =
\tildeg^{\tildeeta}$.
{\thickmuskip=3.5mu \medmuskip=3mu \thinmuskip=1mu 
$$
h^{\tau\inverse} =
(\tildeh \cup \rfs{Id} \nrestriction
(Y - \rfs{Dom}(\tildeh)))^{\tau\inverse} =
(\tildeg^{\psi} \cup \rfs{Id} \nrestriction
(Y - \psi(\wtildeY)))^{\tau\inverse} =
(\tildeg^{\psi})^{\tau\inverse} \cup
\rfs{Id} \nrestriction (X - \varphi(\wtildeX)).
$$
}
Since $\rfs{Rng}(\varphi)$ is disjoint from $X - \varphi(\wtildeX)$,
$$
\left(\rule{0pt}{0pt}
(\tildeg^{\psi})^{\tau\inverse} \cup
\rfs{Id} \nrestriction
(X - \varphi(\wtildeX))\right)^{\varphi\inverse} =
\left(\rule{0pt}{10pt}
(\tildeg^{\psi})^{\tau\inverse}\right)^{\varphi\inverse}.
$$
That is,
$$
(h^{\tau\inverse})^{\varphi\inverse} =
\left(\rule{0pt}{0pt}
(\tildeg^{\psi})^{\tau\inverse}\right)
^{\varphi\inverse} =
\tildeg^{\tildeeta}.
$$
We conclude that
$\tildeg^{\tildeeta}$
is $\itGamma$-bicontinuous at~$0^E$.
\smallskip

We shall now show that
$\tildeg^{\tildeeta}$ is not $\itGamma$-continuous at~$0^E$
thus reaching a contradiction.
Let $T \in \rfs{Nbr}^E(0)$ and $\gamma\fprime \in \itGamma$.
Then there are $i \in \bbN$ and $a > 0$ such that
$\gamma\fprime \nrestriction
[0,a] \leq \gamma_i \nrestriction [0,a]$, \ \ %
$\tildeeta(u_i),\tildeeta(y_i) \in T$
\ and \ $\norm{\tildeeta(u_i) - \tildeeta(y_i)} \leq a$.
So
$$
\norm{g^{\tildeeta}(\tildeeta(u_i)) - g^{\tildeeta}(\tildeeta(y_i))} =
\norm{\tildeeta(v_i) - \tildeeta(y_i)}  >
\gamma_i(\norm{\tildeeta(u_i) - \tildeeta(y_i)}) \geq
\gamma\fprime(\norm{\tildeeta(u_i) - \tildeeta(y_i)}).
$$
This shows that $g^{\tildeeta}$ is not $\itGamma$-continuous at~$0^E$.
A contradiction, so Claim 1 is proved.
\smallskip

Let $W$ and $\gamma$ be as in Claim 1.
We may assume that $W \subseteq V$.
There is $U \in \rfs{Nbr}^F(0)$ such that for every
$u,v \in U \cap \wtildeY$: if $\pair{u}{v}$ is $\alpha$-submerged in
$\wtildeY$, then $B(u,\norm{v - u}) \subseteq W$.
Let $u,v \in U \cap \wtildeY$ be such that $\pair{u}{v}$ is
$\alpha$-submerged in $\wtildeY$.
Let
$w = y(u,\norm{v - u})$. Then $w \in U$.
Hence
$$
\norm{\tildeeta(v) - \tildeeta(u)} \leq
\gamma(\norm{\tildeeta(w) - \tildeeta(u)}) \leq
\gamma \scirc \sigma(\norm{w - u}) =
\gamma \scirc \sigma(\norm{v - u}).
$$
Clearly, $\gamma \scirc \sigma \in \itOmega$,
and we have just shown that
$\tildeeta$ is $\scolonpair{\gamma \scirc \sigma}{\alpha}$-continuous
at $0^F$.
\hfill\myqed

\begin{defn}\label{ams-bddly-lip-bldr-d5.14}
\begin{rm}
(a) Let $\alpha \in H([0,\infty))$.
For every $t \in [0,\infty)$
we define a sequence $\vect = \setm{t_n}{n \in \bbN}$.
Define $t_0 = t$ and for every $n \in \bbN$,
let $t_{n + 1}$ satisfy the equation
$$
t_{n + 1} + \alpha(t_{n + 1}) = t_n
$$
and define
$$
p_{\alpha,n}(t) = t_n \mbox{\ \ and\ \ }%
q_{\alpha,n}(t)= t_n - t_{n + 1}.
$$
Note that $p_{\alpha,0} = \rfs{Id}$.
   \index{N@p03@@$p_{\alpha,n}(t) =
          ((\rfs{Id} + \alpha)\inverse)^{\sscirc n}(t)$}
   \index{N@q00@@$\rule{1.9pt}{0pt}q_{\alpha,n}(t) =
          p_{\alpha,n}(t) - p_{\alpha,n + 1}(t)$}
\rule{0pt}{1pt}\kern-5pt\index{N@AAAA@@
$\beta \rsstar \alpha(t) = \sum_{n = 0}^{\infty}
   \beta(q_{\alpha,n}(t))$}

(b) Let $\alpha,\beta \in H([0,\infty))$.
We define the function
$\fnn{\beta \rsstar \alpha}
{[0,\infty)}{[0,\infty) \cup \sngltn{\infty}}$.
$$
\beta \rsstar \alpha(t) =
\mathop{\hbox{$\sum$}}\limits_{n = 0}^{\infty}
\beta(q_{\alpha,n}(t)).\kern-5pt
$$

(c) For $\alpha \in \rfs{MC}$ let
$\itGamma_{\alpha} =
\rfs{cl}_{\spreceq}(\setm{\alpha^{\sscirc n}}{n \in \bbN})$.
   \index{N@gamma00@@$\itGamma_{\alpha} =
          \rfs{cl}_{\spreceq}(\setm{\alpha^{\sscirc n}}{n \in \bbN})$}

(d) Let $\itGamma \subseteq \rfs{MC}$ and $\alpha \in \rfs{MC}$.
We say that $\itGamma$ is {\it $\alpha$-star-closed},
if for every $\beta \in \itGamma$ there is $\gamma \in \itGamma$
such that $\beta \rsstar \alpha \preceq \gamma$.
   \index{star-closed. $\itGamma$ is $\alpha$-star-closed}
Let $\itDelta \subseteq \rfs{MC}$.
We say that $\itGamma$ is {\it $\itDelta$-star-closed},
if there is $\delta \in \itDelta$
such that $\itGamma$ is $\delta$-star-closed.
   \index{star-closed. $\itGamma$ is $\itDelta$-star-closed}
\hfill\proofend
\end{rm}
\end{defn}

The next proposition contains some trivial observations about the
operation ``$\rsstar$''. For the continuation of the proof of the main
theorems we need only Parts (a)\,-\,(c) of the proposition.
The other parts are mentioned in order to familiarize the reader with
this operation.
Part~(a) was proved by Wieslaw Kubis.

\begin{prop}\label{ams-bddly-lip-bldr-p5.15}
Let $\alpha,\beta,\gamma \in H([0,\infty))$.

\num{a} for every $n \in \bbN$,
$\alpha^{\sscirc n} \rsstar \alpha \leq
n \alpha^{\sscirc n} + \rfs{Id}$.

\num{b}
If $\gamma \preceq \beta$,
then $\gamma \rsstar \alpha \preceq \beta \rsstar \alpha$.

\num{c} For every $n \in \bbN$, $q_{\alpha,n}$ and $p_{\alpha,n + 1}$
are strictly increasing functions.

\num{d} If $s < t$,
then $\beta \rsstar \alpha(s) \leq \beta \rsstar \alpha(t)$.

\num{e} Either $\beta \rsstar \alpha \nrestriction (0,\infty)$
is the constant function $f(t) = \infty$,
or $\beta \rsstar \alpha \in H([0,\infty))$.
\end{prop}

\noindent
{\bf Proof }
(a) Let $t \in [0,\infty)$. Denote $p_{\alpha,n}(t) = p_n$ and
$q_{\alpha,n}(t) = q_n$. Hence $q_n = \alpha(p_n)$
and $p_n + q_n = p_{n - 1}$.
Let $k \geq n \geq 1$. Then
$$
\alpha^{\sscirc n}(q_k) \leq \alpha^{\sscirc n}(p_{k - 1}) =
\alpha^{\sscirc (n - 1)}(q_{k - 1}) \leq \ldots \leq
\alpha^{\sscirc (n - (n - 1))}(p_{k - n)}) =
\alpha(p_{k - n}) = q_{k - n}.
$$
Note that $\sum_{i = 0}^{\infty} q_i = t$. Let $n \geq 1$. Then
\vspace{1.5mm}
\newline
\rule{7pt}{0pt}
\renewcommand{\arraystretch}{1.5}
\addtolength{\arraycolsep}{-3pt}
$
\begin{array}{ll}
&
\alpha^{\sscirc n} \rsstar \alpha(t) =
\sum\limits_{k = 0}^{\infty}
\alpha^{\sscirc n}(q_k) =
\sum\limits_{k < n} \alpha^{\sscirc n}(q_k) +
\sum\limits_{k \geq n} \alpha^{\sscirc n}(q_k)
\\
\leq
\rule{5pt}{0pt}
&
\sum\limits_{k < n} \alpha^{\sscirc n}(t) +
\sum\limits_{k \geq n} q_{k - n} = n \alpha^{\sscirc n}(t) + t.
\vspace{1.7mm}
\end{array}
$
\renewcommand{\arraystretch}{1.0}
\addtolength{\arraycolsep}{3pt}

(b) This part is immediate.

(c)
Note that $p_{\alpha,n + 1} + q_{\alpha,n} = p_{\alpha,n}$.
This equality
together with the facts that $\alpha$ is strictly increasing
and $p_{\alpha,0} = \rfs{Id}$,
implies by induction that $q_{\alpha,n}$ and $p_{\alpha,n + 1}$ are
strictly increasing for every $n \in \bbN$.

(d) Part (d) follows from the facts that $q_{\alpha,n}$
and $\beta$ are increasing functions.

(e)
Note that $q_{\alpha,k}(p_{\alpha,n}(t)) = q_{\alpha,k + n}(t)$.
Hence $\beta \rsstar \alpha(p_{\alpha,n}(t))$ is a tail of
$\beta \rsstar \alpha(t)$.
So for every $n$, $\beta \rsstar \alpha(p_{\alpha,n}(t)) < \infty$ iff
$\beta \rsstar \alpha(t) < \infty$.
Note also that $\lim_n p_{\alpha,n}(t) = 0$.
Suppose that for some $t$,
$\beta \rsstar \alpha(t) = \infty$ and let $s > 0$.
Then there is $n$ such that $p_{\alpha,n}(t) < s$.
So
$\infty = \beta \rsstar \alpha(p_{\alpha,n}(t)) \leq
\beta \rsstar \alpha(s)$.
Hence $\beta \rsstar \alpha \nrestriction (0,\infty)$ is the constant
function with value $\infty$.

Suppose that $\beta \rsstar \alpha \nrestriction (0,\infty)$ is not the
constant $\infty$.
So $\rfs{Rng}(\beta \rsstar \alpha) \subseteq [0,\infty)$.
Note that
$q_{\alpha,0} = \alpha \scirc p_{\alpha,0} =
\alpha \scirc (\rfs{Id} + \alpha)\inverse$.
So $\limti{t} q_{\alpha,0}(t) = \infty$.
For $\beta$ we have
$\lim_{t \rightarrow \infty} \beta(t) = \infty$.
It follows that
$\lim_{t \rightarrow \infty} \beta \rsstar \alpha(t)
\geq \limti{t} \beta(q_{\alpha,0}(t)) = \infty$.

The strict increasingness of $\beta$ and all the $q_{\alpha,n}$'s
together with the fact that
$\beta \rsstar \alpha(t) < \infty$ for every $t$,
implies that $\beta \rsstar \alpha$ is strictly increasing.

It remains to show that $\beta \rsstar \alpha$ is continuous.
Let $a \in (0,\infty)$,
and we show that\break
$\sum_n \beta(q_{\alpha,n}(t))$ is uniformly convergent in $[0,a]$.
Let $\varepsilon > 0$.
There is $n$ such that\break
$\sum_{k \geq n} \beta(q_{\alpha,k}(a)) < \varepsilon$.
From the increasingness of $\beta$ and all the $q_{\alpha,n}$'s
it follows that $\sum_{k \geq n} \beta(q_{\alpha,k}(t)) < \varepsilon$
for all $t \in [0,a]$.
So \,$\sum_n \beta(q_{\alpha,n}(t))$ is uniformly convergent
in $[0,a]$.
Hence $\beta \rsstar \alpha$ is continuous.
\hfill\myqed

\begin{question}\label{ams-bddly-lip-bldr-q5.16}
\begin{rm}

(a) Let $\alpha, \beta \in \rfs{MC}$.
Is it true that either $\beta \rsstar \alpha \nrestriction (0,\infty)$
is the constant function $\infty$, or $\beta \rsstar \alpha$
belongs to $\rfs{MC}$?

(b) Let $\alpha_1, \alpha_2, \beta \in \rfs{MC}$.
Is the following statement true?
If $\alpha_1 \preceq \alpha_2$,
then $\beta \rsstar \alpha_2 \preceq \beta \rsstar \alpha_1$.

(c) Let $\alpha \in \rfs{MC}$.
Is there $\beta \in \rfs{MC} - \itGamma_{\alpha}$
such that $\itGamma_{\beta}$ is $\alpha$-star-closed?
\end{rm}
\end{question}

\begin{prop}\label{ams-bddly-lip-bldr-p5.17}
Let $K > 0$, $r \in (0,1)$, $\alpha(t) = Kt$ and
$\beta(t) = t^r$.
Then there is $C$ such that $\beta \rsstar \alpha = C \ncdot \beta$.
\end{prop}

\noindent
{\bf Proof } 
Abbreviate $q_{\alpha,n}(t)$ by $q_n$.
Let $t \geq 0$.
Then
$$\hbox{$
q_n = \left(\frac{1}{(1 + K)^n} - \frac{1}{(1 + K)^{n + 1}}\right)
\mcdot t =
\frac{1}{(1 + K)^n} \mcdot \frac{Kt}{1 + K}
$}$$
and hence
$$\hbox{$
\beta \rsstar \alpha =
\sum\limits_{n = 0}^{\infty} \frac{1}{(1 + K^r)^n} \mcdot
(\frac{Kt}{1 + K})^r =
\frac{(1 + K)^r}{(1 + K)^r - 1} \mcdot
\frac{K^r}{(1 + K)^r} \mcdot t^r =
\frac{K^r}{(1 + K)^r - 1} \mcdot \beta(t).
$}$$
So $C = \dgfrac{K^r}{((1 + K)^r - 1)}$.
\smallskip\rule{0pt}{1pt}\hfill\myqed

Lemma \ref{ams-bddly-lip-bldr-l5.20}(b) is our next main step.
It is preceded by two propositions.
Part (a) of \ref{ams-bddly-lip-bldr-l5.20} is also a step in the proof
of \ref{ams-bddly-lip-bldr-l5.20}(b).
For $\alpha \in \rfs{MC}$, a normed space $E$ and $x,y \in E$ let
$\rfs{prt}_{\alpha}(x,y)$ be the point $z$ in the line segment
$[x,y]$ such that $\alpha(\norm{z - y}) = \norm{x - z}$.

\begin{prop}\label{ams-bddly-lip-bldr-p5.18}
Let $\alpha \in \rfs{MC}$ and $a > 0$.
Then there is $\varepsilon = \varepsilon_{\alpha,a}$
such that the following holds.
If $F$ is a normed space, $M$ is a closed subspace of $F$
or a closed half space of $F$,
$x \in F - M$ and $d(x,M) = a$,
then for every $y \in \rfs{bd}(M)$: if $d(x,y) < a + \varepsilon$,
then $\pair{x}{\rfs{prt}_{\alpha}(x,y)}$ is $2 \alpha$-submerged in
$F - M$.
\end{prop}

\noindent
{\bf Proof }
Let
$q(t) = q_{\alpha,0}(t)$ and $f(t) = q(t) + (2 \alpha)\inverse(q(t))$.
Then
$f(t) = q(t) + \alpha\inverse(\half q(t)) <
q(t) + \alpha\inverse(q(t))$.
In particular,
$f(a) < q(a) + \alpha\inverse(q(a)) = a$.
So there is $\varepsilon > 0$ such that for every $t$: if
$\abs{t - a} < \varepsilon$, then
{\thickmuskip=3.5mu \medmuskip=2mu \thinmuskip=1mu 
$f(t) < \frac{f(a) + a}{2}$.
\hbox{Let $y \in \rfs{bd}(M)$ be such that $d(x,y) < a + \varepsilon$.
Then\kern-4pt}}
\vspace{1.5mm}
\newline
\rule{7pt}{0pt}
\renewcommand{\arraystretch}{1.5}
\addtolength{\arraycolsep}{-3pt}
$
\begin{array}{ll}
&
\norm{x - \rfs{prt}_{\alpha}(x,y)} +
(2\alpha)\inverse(\norm{x - \rfs{prt}_{\alpha}(x,y)}) =
q(\norm{x - y}) + (2\alpha)\inverse(q(\norm{x - y}))
\\
=
\rule{5pt}{0pt}
&
f(\norm{x - y}) < \frac{f(a) + a}{2} < a = \delta^{F - M}(x).
\\
\end{array}
$
\renewcommand{\arraystretch}{1.0}
\addtolength{\arraycolsep}{3pt}
\smallskip\newline
So $\pair{x}{\rfs{prt}_{\alpha}(x,y)}$ is $2 \alpha$-submerged in
$F - M$.
\hfill\myqed

\begin{prop}\label{ams-bddly-lip-bldr-p5.19}
Let $\alpha \in \rfs{MC}$,
$F$ be a normed space, $M$ be a closed subspace of $F$
or a closed half space of $F$,
{\thickmuskip=2.0mu \medmuskip=1.0mu \thinmuskip=1mu 
$x \in F - M$ and $y \in M$.
Then there is a sequence $\setm{x_i}{i \in \bbN}$ such that:}
\begin{itemize}
\addtolength{\parskip}{-11pt}
\addtolength{\itemsep}{06pt}
\item[\num{i}] $x_0 = x$,
\item[\num{ii}] for every $i \in \bbN$,
$\pair{x_i}{x_{i + 1}}$ is $2 \alpha$-submerged in $F - M$,
\item[\num{iii}] for every $i \in \bbN$,
$\norm{x_i - x_{i + 1}} \leq q_{\alpha,i}(\norm{x - y})$,
\item[\num{iv}] $\lim_i x_i$ exists and $\lim_i x_i \in \rfs{bd}(M)$,
\item[\num{v}] $\norm{\lim_i x_i - y} \leq 2 \norm{x - y}$.
\vspace{-05.7pt}
\end{itemize}
Note that the convergence of $\setm{x_i}{i \in \bbN}$
follows from \num{iii}, and need not be required.
\end{prop}

\noindent
{\bf Proof } Write $p_{\alpha,i} = p_i$ and $q_{\alpha,i} = q_i$.
Note that $p_1 \scirc p_i = p_{i + 1}$
and that $q_0 \scirc p_i = q_i$.
Let $x_0 = x$ and $y_0 = y$.
We define by induction $x_i \in F - M$ and $y_i \in \rfs{bd}(M)$.
Suppose that $x_i,y_i$ have been defined.
Let $y_{i + 1} \in \rfs{bd}(M)$
be such that
$\norm{x_i - y_{i + 1}} \leq \norm{x_i - y_i}$
and
$\pair{x_i}{\rfs{prt}_{\alpha}(x_i,y_{i + 1})}$
is
$2 \alpha$-submerged in $F - M$.
The existence of such $y_{i + 1}$ is assured by
Proposition \ref{ams-bddly-lip-bldr-p5.18}.
Let $x_{i + 1} = \rfs{prt}_{\alpha}(x_i,y_{i + 1})$.
(Note that if for some $\bary \in M$,
$d(x,M) = \norm{x - \bary}$, then $y_i$ can be chosen to be
$\bary$ for every $i \geq 1$).

By the definitions, Clauses (i) and (ii) hold.
We prove Clause (iii).
We prove by induction on $i$ that
$\norm{x_i - x_{i + 1}} \leq q_i(\norm{x - y})$
and
$\norm{x_{i + 1} - y_{i + 1}} \leq p_{i + 1}(\norm{x - y})$.
It is trivial that the induction hypotheses hold for $i = 0$.
Suppose that the induction hypotheses hold for $i - 1$.
Then
\vspace{1.5mm}
\newline
\rule{7pt}{0pt}
\renewcommand{\arraystretch}{1.5}
\addtolength{\arraycolsep}{-3pt}
$
\begin{array}{l}
\norm{x_i - x_{i + 1}} = q_0(\norm{x_i - y_{i + 1}}) \leq
q_0(\norm{x_i - y_i}) \leq
q_0(p_i(\norm{x - y})) = q_i(\norm{x - y}).
\\
\norm{x_{i + 1} - y_{i + 1}} = p_1(\norm{x_i - y_{i + 1}}) \leq
p_1(\norm{x_i - y_i}) \leq p_1(p_i(\norm{x - y})) =
p_{i + 1}(\norm{x - y}).
\vspace{1.7mm}
\end{array}
$
\renewcommand{\arraystretch}{1.0}
\addtolength{\arraycolsep}{3pt}
\newline
So Clause (iii) holds.

We prove Clause (iv).
Obviously,
$\sum_{i = 0}^{\infty} q_i(\norm{x - y}) = \norm{x - y}$.
Since
$\norm{x_i - x_{i + 1}} \leq q_i(\norm{x - y})$,
it follows that
$\sum_{i = 0}^{\infty} \norm{x_i - x_{i + 1}}$
is convergent.
So $\setm{x_i}{i \in \bbN}$ is convergent.
Let $\barx = \lim_i x_i$.
The facts $\lim_i p_i(\norm{x - y}) = 0$
and $\norm{x_i - y_i} \leq p_i(\norm{x - y})$
imply that $\lim_i \norm{x_i - y_i} = 0$.
Since $y_i \in \rfs{bd}(M)$,
it follows that $\barx \in \rfs{bd}(M)$.

We prove Clause (v).
$$
\norm{\barx - x} \leq
\mathop{\hbox{$\sum$}}\limits_{i = 0}^{\infty} \norm{x_i - x_{i + 1}}
\leq
\mathop{\hbox{$\sum$}}\limits_{i = 0}^{\infty} q_i(\norm{x - y}) =
\norm{x - y}.
$$
So $\norm{\barx - y} \leq \norm{\barx - x} + \norm{x - y)} \leq
2 \norm{x - y}$.
\hfill\myqed

\begin{lemma}\label{ams-bddly-lip-bldr-l5.20}
Assume that Clauses \num{1}-\num{7} of
Lemma \ref{ams-bddly-lip-bldr-l5.13} hold. That is,
\begin{itemize}
\addtolength{\parskip}{-11pt}
\addtolength{\itemsep}{06pt}
\item[\num{1}]
$\itGamma,\itSigma$ are countably generated moduli of continuity,
and $\itOmega$ is the modulus of continuity
\indent\kern3pt generated by
$\itGamma \cup \itSigma$.
\item[\num{2}]
$X \subseteq E$ and $Y \subseteq F$ are open subsets of
the normed spaces $E$ and $F$,
$X$ is $\itGamma$-LIN-
\indent\kern3pt bordered and
$Y$ is $\itSigma$-LIN-bordered.
\item[\num{3}]
$\tau \in \rfs{EXT}^{\pm}(X,Y)$,
$G$ is a $\itGamma$-appropriate subgroup of $\rfs{EXT}(X)$,
$H$ is a $\itDelta$-appropriate
\indent\kern2pt subgroup of $\rfs{EXT}(Y)$
and $G^{\tau} = H$.
\item[\num{4}]
$x \in \rfs{bd}(X)$,
$\trpl{\varphi}{A}{r}$ is a boundary chart element for~$x$,
$\gamma \in \itGamma$
and $\varphi$ is $\gamma$-bicontinuous.
\item[\num{5}]
$y \in \rfs{bd}(Y)$,
$\trpl{\psi}{B}{s}$ is a boundary chart element for~$y$,
$\sigma \in \itSigma$
and $\psi$ is $\sigma$-bicontinuous.
\item[\num{6}]
$\tau^{\srfs{cl}}(x) = y$
and
$\tau(\varphi(\rfs{BCD}^{E}(A,r))) \subseteq \psi(\rfs{BCD}^{F}(B,s))$.
\item[\num{7}] Set
$L = \rfs{bd}(A)$,
$\swhatX = \rfs{BCD}^{E}(A,r) \cup (L \cap B(0,r))$,
$\hattau = \psi\inverse \scirc \tau^{\srfs{cl}} \scirc \varphi$,
$\whatY = \hattau(\swhatX)$, $\wtildeY = \hattau(\rfs{BCD}^{E}(A,r))$
and $\calP = \setm{(v + L) \cap \swhatX}{v \in \swhatX}$.
\vspace{-05.7pt}
\end{itemize}
Assume further that
\begin{itemize}
\addtolength{\parskip}{-11pt}
\addtolength{\itemsep}{06pt}
\item[\num{8}]
$\itOmega$ is $\itSigma$-star-closed.
\vspace{-05.7pt}
\end{itemize}

\num{a} Let $M = \rfs{bd}(B)$.
Then there is $W \in \rfs{Nbr}^F(0)$ and $\omega \in \itOmega$
such that for every
$x \in (\whatY - M) \cap W$ and $y \in \whatY \cap M \cap W$,
$\norm{\hattau\inverse(x) - \hattau\inverse(y)} \leq
\omega(\norm{x - y})$.

\num{b} $\hattau\inverse$ is $\itOmega$-continuous at $\hattau(0)$.
\end{lemma}

\noindent
{\bf Proof } (a)  Let $\alpha \in \itSigma$ be such that
$\itOmega$ is $\alpha$-star-closed.
It is easy to see that $\alpha$ may be chosen to be in $\rfs{MBC}$.
Note that $\wtildeY = \whatY - M$.
Let $\hateta = \hattau\inverse$.
By Lemma \ref{ams-bddly-lip-bldr-l5.13},
there are $\rho \in \itOmega$ and $W_1 \in \rfs{Nbr}^F(0^F)$
such that for every $u,v \in W_1 \cap \wtildeY$:
if $\pair{u}{v}$ is $2 \alpha$-submerged in $\wtildeY$,
then $\norm{\hateta(u) - \hateta(v)} \leq \rho(\norm{u - v})$.
Let $\nu \in \itOmega$ and $a > 0$ be such that
$\rho \rsstar \alpha \nrestriction [0,a] \leq \nu  \nrestriction [0,a]$.

Let $\calP = \setm{(v + L) \cap \swhatX}{v \in \swhatX}$.
By Lemma \ref{ams-bddly-lip-bldr-l5.9},
$\hattau$ is inversely $\pair{\itDelta}{\calP}$-continuous at $0^E$.
Note that $\hattau(L \cap \swhatX) = M \cap \whatY$,
that is, $M \cap \whatY \in \hattau(\calP)$.
So there are $\sigma \in \itSigma$ and $W_2 \in \rfs{Nbr}^F(0^F)$
such that for every $u,v \in W_2 \cap M \cap \whatY$,
$\norm{\hateta(u) - \hateta(v)} \leq \sigma(\norm{u - v})$.
Choose $s_0 \in (0,\dgfrac{a}{2})$ such that
$\overB(0^F,6 s_0) \cap \rfs{BCD}^F(B,s) \subseteq
\whatY \cap W_1 \cap W_2$
and let  $W = B(0^F,s_0)$.

Let $x \in (\whatY - M) \cap W$ and
and $y \in \whatY \cap M \cap W$.
Let $\setm{x_i}{i \in \bbN}$ be the sequence assured by 
Proposition~\ref{ams-bddly-lip-bldr-p5.19} and $\barx = \lim_i x_i$.
Note that by Clause~(iii) of \ref{ams-bddly-lip-bldr-p5.19},
$\sum_{i \in \bbN} \norm{x_i - x_{i + 1}} \leq \norm{x - y} < 2 s_0$.
So
$\norm{x_n} \leq
\norm{x} + \sum_{i = 0}^{n - 1} \norm{x_i - x_{i + 1}} < 3 s_0$
for every $n \in \bbN$.
Similarly, $\norm{\barx} < 3 s_0$.
Hence
$\setm{x_i}{i \in \bbN} \subseteq W_1 \subseteq \rfs{Dom}(\hateta)$
and $\barx \in W_2 \subseteq \rfs{Dom}(\hateta)$.
We conclude that
$$
\norm{\hateta(x) - \hateta(y)} \leq
\hbox{$\sum\limits_{i = 0}^{\infty}$}
\norm{\hateta(x_i) - \hateta(x_{i + 1})} +
\norm{\hateta(\barx) - \hateta(y)} \eqdf \bfiA.
$$
Since $\barx,y \in W_2 \cap M \cap \whatY$,
we have that
$\norm{\hateta(\barx) - \hateta(y)} \leq  \sigma(\norm{\barx - y})$.

By Clause~(ii) of \ref{ams-bddly-lip-bldr-p5.19},
$\pair{x_i}{x_{i + 1}}$ is $2 \alpha$-submerged in $F - M$.
Using the facts that
$x_i \in B(0,3s_0)$ and that
$B(0^F,6 s_0) \cap \rfs{BCD}^F(B,s) \subseteq \whatY$,
it is easily seen that $\delta^{\wtildeY}(x_i) = \delta^{F -M}(x_i)$
for every $i \in \bbN$.
So $\pair{x_i}{x_{i + 1}}$ is $2 \alpha$-submerged in $\wtildeY$.
This, together with the fact that $x_i,x_{i + 1} \in W_1$, implies
that
$\norm{\hateta(x_i) - \hateta(x_{i + 1})} \leq
\rho(\norm{x_i - x_{i + 1}})$.
Hence
$$
\bfiA \leq \hbox{$\sum\limits_{i = 0}^{\infty}$}
\rho(\norm{x_i - x_{i + 1}}) + \sigma(\norm{\barx - y}) \eqdf \bfiB.
$$
By the increasingness of $q_{\alpha,i}$ and Clause (iii) in
Proposition \ref{ams-bddly-lip-bldr-p5.19},
$$
\hbox{$\sum\limits_{i = 0}^{\infty}$} \rho(\norm{x_i - x_{i + 1}}) \leq
\hbox{$\sum\limits_{i = 0}^{\infty}$} \rho(q_{\alpha,i}(\norm{x - y})) =
\rho \rsstar \alpha(\norm{x - y}).
$$
Clause (v) in \ref{ams-bddly-lip-bldr-p5.19} implies that
$\sigma(\norm{\barx - y}) \leq \sigma(2 \norm{x - y})$.
Hence
$$
\bfiB \leq \rho \rsstar \alpha(\norm{x - y}) +
\sigma \scirc (2 \mcdot \rfs{Id})(\norm{x - y}).
$$
Recall that $\nu \in \itOmega$,
$\rho \rsstar \alpha \nrestriction [0,a] \leq \nu  \nrestriction [0,a]$
and $s_0 < \dgfrac{a}{2}$.
Let $\omega = \nu + \sigma \scirc (2 \mcdot \rfs{Id})$.
It follows from the above that
$\omega \in \itOmega$
and $\norm{\hateta(x) - \hateta(y)} \leq \omega(\norm{x - y})$.
This proves Part (a).\smallskip

(b) We use the notations of Part (a).
Let $x,y \in W \cap \whatY$.
If $x,y \in M$,
then $\norm{\hateta(x) - \hateta(y)} \leq \sigma(\norm{x - y})$.
If $x \not\in M$ and $y \in M$ or vice versa,
then $\norm{\hateta(x) - \hateta(y)} \leq \omega(\norm{x - y})$.
Suppose that $x,y \not\in M$ and write $\beta  = 2 \alpha$.
If $\pair{x}{y}$ is $\beta$-submerged in $\wtildeY$
or $\pair{y}{x}$ is $\beta$-submerged in $\wtildeY$,
then $\norm{\hateta(x) - \hateta(y)} \leq \rho(\norm{x - y})$.

Suppose that neither $\pair{x}{y}$ nor $\pair{y}{x}$ are
$\beta$-submerged in $\wtildeY$.
Since
$x,y \in B(0,s_0)$ and
$B(0,6s_0) \cap \rfs{BCD}^F(B,s) \subseteq \whatY$, \,
$\delta^{F - M}(x) = \delta^{\wtildeY}(x)$ and
$\delta^{F - M}(y) = \delta^{\wtildeY}(y)$.
So by the non-submergedness of $\pair{x}{y}$ and $\pair{y}{x}$,
$\delta^{F - M}(x), \delta^{F - M}(y) <
\norm{x - y} + \beta\inverse(\norm{x - y})$.
Since $\beta \in \rfs{MBC}$, $\beta\inverse(t) \leq t$ for every $t$.
So
$\delta^{F - M}(x), \delta^{F - M}(y) < 2 \norm{x - y}$.

Let $\barx,\bary \in M$ be such that
$\norm{x - \barx} < 2 \delta^{F - M}(x)$
and $\norm{y - \bary} < 2 \delta^{F - M}(y)$.
Clearly, $\norm{\barx} < 3 \norm{x} < 3s_0$.
Hence $\barx \in W_2 \cap M \cap \whatY$.
Similarly, $\bary \in W_2 \cap M \cap \whatY$.
We also have that
$$\norm{\barx - \bary} \leq
\norm{\barx - x} + \norm{x - y} + \norm{y - \bary} \leq
2 \delta^{F - M}(x) + \norm{x - y} + 2 \delta^{F - M}(y) \leq
9 \norm{x - y}
$$
and $\norm{x - \barx},\norm{y - \bary} < 4 \norm{x - y}$.
The final estimate is
\vspace{1.5mm}
\newline
\rule{7pt}{0pt}
\renewcommand{\arraystretch}{1.5}
\addtolength{\arraycolsep}{-3pt}
$
\begin{array}{ll}
&
\norm{\hateta(x) - \hateta(y)} \leq
\norm{\hateta(x) - \hateta(\barx)} +
\norm{\hateta(\barx) - \hateta(\bary)} +
\norm{\hateta(\bary) - \hateta(y)}
\\
\leq
\rule{5pt}{0pt}
&
\omega(\norm{x - \barx}) + \sigma(\norm{\barx - \bary}) +
\omega(\norm{\bary - y})
\\
\leq
\rule{5pt}{0pt}
&
\omega(4 \norm{x - y}) + \sigma(9 \norm{x - y}) +
\omega(4 \norm{x - y}) \leq
8 \omega(\norm{x - y}) + 9 \sigma(\norm{x - y}).
\vspace{1.7mm}
\end{array}
$
\renewcommand{\arraystretch}{1.0}
\addtolength{\arraycolsep}{3pt}
\newline
Clearly, $\gamma \eqdf 8 \omega + 9 \sigma \in \itOmega$.
Obviously, $\sigma,\omega,\rho \leq \gamma$.
We have thus shown that for every $x,y \in W \cap \whatY$,
$\norm{\hateta(x) - \hateta(y)} \leq \gamma(\norm{x - y})$.
So $\hateta$ is $\itOmega$-continuous at $0^F$.
\smallskip
\hfill\myqed

We make a last trivial observation before proving the main theorem.

\begin{prop}\label{ams-bddly-lip-bldr-p5.21}
\num{a} Let $\itGamma$ be a modulus of continuity and
$\alpha \in \rfs{MBC} - \itGamma$.
Let $X$ be an open subset of a normed space $E$
and $x \in \rfs{bd}(X)$. Then there is $g \in H(E)\sprt{X}$
such that $g$ is $9 \mcdot \alpha \scirc \alpha$-bicontinuous and
$g$ is not $\itGamma$-bicontinuous at $x$.

\num{b}
Let $\itGamma,\itDelta$ be moduli of continuity,
$E,F$ be normed spaces,
$X \subsetneqq E$ be an open $\itGamma$-LIN-bordered set,
$Y \subseteq F$ be an open $\itDelta$-LIN-bordered set,
$G \leq \rfs{EXT}(X)$ and $H \leq \rfs{EXT}(Y)$\break
be respectively
$\itGamma$-appropriate and $\itDelta$-appropriate,
$\tau \in (H^{\srfs{BDR.LC}}_{\itDelta})^{\pm}(X,Y)$
and $G^{\tau} = H$.
\underline{Then} $\itGamma = \itDelta$.
\end{prop}

\noindent
{\bf Proof } (a)
For $r > 0$ define $\iso{g_r}{E}{E}$ as follows:
$g_r(0) = 0$,
$g_r(z) =
\frac{r}{\alpha(r)} \mcdot \alpha(\norm{z}) \mcdot \frac{z}{\norm{z}}$
if $\norm{z} \in (0,r)$, and $g_r(z) = z$ if $\norm{z} \geq r$.
Obviously, $\rfs{supp}(g_r) = B(0,r)$,
and it is left to the reader to check that $g_r$ is
$\frac{3r}{\alpha(r)} \mcdot \alpha$-bicontinuous,
and that if $\gamma \in \rfs{MC}$ is such that $g_r$ is
$\gamma$-bicontinuous, then
$\gamma \nrestriction [0,r] \geq
\frac{r}{\alpha(r)} \mcdot \alpha \nrestriction [0,r]$.
For $y \in E$ define $g_{y,r} = g_r^{\srfs{tr}_{\kern-1pty}}$.
Let $\setm{B(x_i,r_i)}{i \in \bbN}$ be a sequence
of pairwise disjoint balls
such that for every $i$, $B(x_i,r_i) \subseteq X$
and $\lim_i x_i = x$,
and let $g = \bcirc_i g_{x_i,r_i} \nrestriction X$.
Then $g$ is as required.

(b) First we show that $\itDelta \subseteq \itGamma$.
Suppose otherwise. 
Let $x \in \rfs{bd}(X)$
and $y = \tau^{\srfs{cl}}(x)$. So $y \in \rfs{bd}(Y)$.
There are $W \in \rfs{Nbr}(y)$
and $\beta \in \itDelta$ such that
$\tau\inverse \nrestriction (W \cap Y)$ is $\beta$-bicontinuous.
Let $V \in \rfs{Nbr}(y)$ such that $V \subseteq W$ and
$H^{\srfs{CMP.LC}}_{\itDelta}(Y)\sprt{V \cap Y} \subseteq H$.
Choose $\alpha \in \itDelta \cap \rfs{MBC} - \itGamma$
and define $\baralpha = 9 \mcdot \alpha \scirc \alpha$
and $\delta = \beta \scirc \baralpha \scirc \beta$.
Let $U = \tau\inverse(V \cap Y)$.
Hence $x \in \rfs{bd}(U)$.

Let $X'$ be an open subset of $U \cap X$ such that
$\rfs{cl}(X') \cap \rfs{bd}(X) = \sngltn{x}$.
By (a), there is $g\fprime \in H(E)\sprt{X'}$ such that
$g\fprime$ is $\baralpha$-bicontinuous,
and $g\fprime$ is not $\itGamma$-bicontinuous at $x$.
Let $g = g\fprime \nrestriction X$ and $h = g^{\tau}$.
Since $g$ is $E$-biextendible and $\tau$ is $(E,F)$-biextendible,
$h$ is $F$-biextendible.
From the fact that $\tau \nrestriction (U \cap X)$ is 
$\beta$-bicontinuous, it follows that
$h \nrestriction (V \cap Y)$ is $\delta$-bicontinuous.
We wish to conclude that $h$ is $\delta$-bicontinuous.
Indeed, this follows from the facts:
$\rfs{cl}^F(\rfs{supp}(h)) \subseteq (V \cap Y) \cup \sngltn{y}$
and $y \in \rfs{cl}(V \cap Y)$.
(The same argument appears in the proof \ref{ams-bddly-lip-bldr-l5.13},
where it is proved that $h \in H^{\srfs{CMP.LC}}_{\itSigma}(Y)$).
Obviously, $\delta \in \itDelta$,
so $h \in H^{\srfs{CMP.LC}}_{\itDelta}(Y)\sprt{V \cap Y} \subseteq H$.
Recall that $G^{\tau} = H$, hence $g = h^{\tau\inverse} \in G$.
But $g$ is not $\itGamma$-bicontinuous at $x$.
This contradicts the fact that $G$ is $\itGamma$-appropriate.
Hence $\itDelta \subseteq \itGamma$.

It follows that
$\tau \in (H^{\srfs{BDR.LC}}_{\itGamma})^{\pm}(X,Y)$
and hence $\tau\inverse \in (H^{\srfs{BDR.LC}}_{\itGamma})^{\pm}(Y,X)$.
We now repeat the above argument for $\tau\inverse$.
So the roles of $\itGamma$ and $\itDelta$ are interchanged,
and we conclude that $\itGamma \subseteq \itDelta$.
\hfill\myqed

\subsection{Final results.}\label{ss12.3}
\begin{theorem}\label{ams-bddly-lip-bldr-t5.22}
{\rm Main Theorem of Chapter \ref{s12}.}
Assume that
\begin{itemize}
\addtolength{\parskip}{-11pt}
\addtolength{\itemsep}{06pt}
\item[\num{1}]
$\itGamma,\itDelta$ are countably generated moduli of continuity,
$\itGamma \subseteq \itDelta$ and $\itDelta$ is
$\itGamma$-star-closed.
(Or assume the special cases:
\num{i} $\itGamma$ is principal and $\itDelta = \itGamma$,
or
\num{ii} $\itGamma = \itGamma^{\srfs{LIP}}$
and $\itDelta = \itGamma^{\srfs{HLD}}$).
\item[\num{2}]
$X \subsetneqq E$ and $Y \subseteq F$ are open subsets of
the normed spaces $E$ and $F$,
$X$ is $\itGamma$-LIN-bordered,
and $Y$ is $\itDelta$-LIN-bordered.
\item[\num{3}]
$G \leq \rfs{EXT}(X)$ is $\itGamma$-appropriate,
and $H \leq \rfs{EXT}(Y)$ is $\itDelta$-appropriate.
\item[\num{4}]
$\tau \in \rfs{EXT}^{\pm}(X,Y)$ and $G^{\tau} = H$.
\vspace{-05.7pt}
\end{itemize}
\underline{Then} $\itGamma = \itDelta$ and
$\tau \in (H^{\srfs{BDR.LC}}_{\itGamma})^{\pm}(X,Y)$.
\end{theorem}

\noindent
{\bf Proof }
That (i) is a special case of (1)
follows from Proposition~\ref{ams-bddly-lip-bldr-p5.15}(a) and (b),
and that (ii) is a special case,
follows from Proposition~\ref{ams-bddly-lip-bldr-p5.17}.

Since $\itGamma \subseteq \itDelta$, the modulus of continuity
$\itOmega$ which is generated by $\itGamma \cup \itDelta$ is
$\itDelta$, and since $\itDelta$ is $\itGamma$-star-closed and
$\itGamma \subseteq \itDelta$,
we have that $\itDelta$ is $\itDelta$-star-closed.
So $\itOmega$ is $\itDelta$-star-closed.
Let $x \in \rfs{bd}(X)$.
There are a boundary chart element for $x$,
$\trpl{\varphi}{A}{r}$ and $\gamma \in \itGamma$
such that $\varphi$ be $\gamma$-bicontinuous.
Let $y = \tau^{\srfs{cl}}(x)$.
Choose a boundary chart element for $y$,
$\trpl{\psi}{B}{s}$ and $\sigma \in \itSigma$
such that $\psi$ is $\sigma$-bicontinuous.
Also assume that
$\tau(\varphi(\rfs{BCD}^{E}(A,r))) \subseteq \psi(\rfs{BCD}^{F}(B,s))$.
Set $L = \rfs{bd}(A)$,
$\whatX = \rfs{BCD}^{E}(A,r) \cup (L \cap B(0,r))$,
$\hattau = \psi\inverse \scirc \tau^{\srfs{cl}} \scirc \varphi$
and $\whatY = \hattau(\swhatX)$.

By Theorem \ref{ams-bddly-lip-bldr-l5.20}(b), $\hattau\inverse$ is
$\itOmega$-continuous at~$0^F$.
That is, $\hattau\inverse$ is $\itDelta$-continuous at~$0^F$.
Since $\varphi,\psi$ are $\itDelta$-bicontinuous at $0^E$ and $0^F$
respectively, $\varphi \scirc \hattau\inverse \scirc \psi\inverse$
is $\itDelta$-continuous\break
at $y$.
Note that there is $V \in \rfs{Nbr}^F(y)$ such that
$\rfs{Dom}(\varphi \scirc \hattau\inverse \scirc \psi\inverse) \supseteq
V \cap Y$.
Also,
$\varphi \scirc \hattau\inverse \scirc \psi\inverse
\nrestriction (V \cap Y) =
\tau\inverse \nrestriction (V \cap Y)$.
Hence $\tau\inverse$ is $\itDelta$-continuous at $y$.
Since it is also given that $\tau \in \rfs{EXT}^{\pm}(X,Y)$,
it follows that
$\tau\inverse \in H^{\srfs{BDR.LC}}_{\itDelta}(Y,X)$.

We now reverse the roles of $X$ and $Y$. Let $\eta = \tau\inverse$.
So $\iso{\eta}{Y}{X}$, $H^{\eta} = G$ and the modulus of continuity
$\itOmega$ generated by $\itDelta \cup \itGamma$ is again $\itDelta$.
So $\itOmega$ is $\itGamma$-star-closed.

Let $y \in \rfs{bd}(Y)$ and $x = \eta(y)$.
We choose $\psi$ and $\varphi$ and define $\hateta$ in the same way
that $\varphi$, $\psi$ and $\hattau$ were defined in the preceding
argument.
We thus conclude that $\hateta\inverse$
is $\itOmega$-continuous at $x$.
That is, $\hateta\inverse$ is $\itDelta$-continuous at $x$.
There is $U \in \rfs{Nbr}^E(x)$ such that
$\psi \scirc \hateta \scirc \varphi\inverse \nrestriction (U \cap X) =
\tau \nrestriction (U \cap X)$.
Hence $\tau$ is $\itDelta$-continuous at $x$.
We also need to know that $\tau \in \rfs{EXT}^{\pm}(X,Y)$,
and this is indeed given.
Hence $\tau \in H^{\srfs{BDR.LC}}_{\itDelta}(X,Y)$.
We proved that $\tau \in (H^{\srfs{BDR.LC}}_{\itDelta})^{\pm}(X,Y)$.
By Proposition \ref{ams-bddly-lip-bldr-p5.21}(b),
$\itGamma = \itDelta$.
\smallskip
\rule{0pt}{1pt}\hfill\myqed

\noindent
{\bf Proof of Theorem \ref{ams-bddly-lip-bldr-t1.7} }
If $X = E$ then $Y = F$ and hence
$H^{\srfs{CMP.LC}}_{\itGamma}(X) = H^{\srfs{LC}}_{\itGamma}(X)$,
and the same holds for $Y$.
So in this case the claim of~\ref{ams-bddly-lip-bldr-t1.7}
is implied by Theorem~\ref{metr-bldr-t3.27}.

Assume that $X \neq E$.
We apply Theorem \ref{ams-bddly-lip-bldr-t5.22} to the special case that
$\itGamma = \itDelta$ and $\itGamma$ is principal,
and take $G,H$ to be 
$H^{\srfs{CMP.LC}}_{\itGamma}(X)$ and $H^{\srfs{CMP.LC}}_{\itDelta}(Y)$
respectively.
So
$\tau \in (H^{\srfs{BDR.LC}}_{\itGamma})^{\pm}(X,Y)$.
By Theorem \ref{metr-bldr-t3.27},
$\tau$ is locally $\itGamma$-bicontinuous.
Hence
$\tau \in (H^{\srfs{CMP.LC}}_{\itGamma})^{\pm}(X,Y)$.
\rule{0pt}{0pt}\hfill\myqed

\kern1.3mm

\noindent
{\bf The final recostruction theorems of
Chapters \ref{s8}\,-\,\ref{s12}}

\kern1.5mm

Combining the results of the previous sections in different ways,
one obtains various reconstruction theorems.
Parts (a) and (b) of the following theorem are such corollaries.
Part (a) is a restatement of
Theorem~\ref{ams-bddly-lip-bldr-t1.4}(a).
Indeed, the special case of Part (a),
in which $\itGamma = \itGamma_{\srfs{LIP}}$
motivated the whole work presented in Chapters \ref{s8}\,-\,\ref{s12}.

The reconstruction theorem for the group
$H^{\srfs{BDR.LC}}_{\itGamma}(X)$ which appears in Part (b),
is a byproduct of the proof of the main result. We thought it was worth
mentioning.

In Part (c) we tried to capture the essence of the argument.
Part (c) can be further strengthened.
But it seems to be a natural stopping point.

\begin{theorem}\label{ams-bddly-lip-bldr-t6.1}
Let $\itGamma,\itDelta$ be moduli of continuity,
$E$ and $F$ be normed spaces and
$X \subseteq E$, $Y \subseteq F$ be open.
Suppose that $X$ is locally $\itGamma$-LIN-bordered,
and $Y$ is locally $\itDelta$-LIN-bordered.

\num{a} Suppose that $\itGamma$ is principal.
If
$\iso{\varphi}{H^{\srfs{CMP.LC}}_{\itGamma}(X)}
{H^{\srfs{CMP.LC}}_{\itDelta}(Y)}$.
Then $\itGamma = \itDelta$
and there is $\tau \in (H^{\srfs{CMP.LC}}_{\itGamma})^{\pm}(X,Y)$
such that $\varphi(g) = g^{\tau}$
for every $g \in H^{\srfs{CMP.LC}}_{\itGamma}(X)$.

\num{b} Suppose that $\itGamma$ is principal.
If
$\iso{\varphi}{H^{\srfs{BDR.LC}}_{\itGamma}(X)}
{H^{\srfs{BDR.LC}}_{\itGamma}(Y)}$.
Then there is\break
$\tau \in (H^{\srfs{BDR.LC}}_{\itGamma})^{\pm}(X,Y)$
such that $\varphi(g) = g^{\tau}$
for every $g \in H^{\srfs{BDR.LC}}_{\itGamma}(X)$.

\num{c} Suppose that $\itGamma$ and $\itDelta$ are countably
generated, $\itGamma \subseteq \itDelta$ and
$\itDelta$ is $\itGamma$-star-closed.
Let $G \leq \rfs{EXT}(X)$ be $\itGamma$-approptiate and
$H \leq \rfs{EXT}(Y)$ be $\itDelta$-approptiate.
Assume further that $\rfs{LIP}^{\srfs{LC}}(X) \leq G$
and $\rfs{LIP}^{\srfs{LC}}(Y) \leq H$,
and suppose that $\iso{\varphi}{G}{H}$.
Then $\itGamma = \itDelta$,
and there is $\tau \in (H^{\srfs{BDR.LC}}_{\itGamma})^{\pm}(X,Y)$
such that $\varphi(g) = g^{\tau}$
for every $g \in G$.
\end{theorem}

\noindent
{\bf Proof } 
(a) By Theorem~\ref{t2.4}(b), there is $\tau \in H(X,Y)$
such that $\tau$ induces $\varphi$.
By Theorem~\ref{metr-bldr-t3.27},
$\itGamma = \itDelta$
and $\tau \in (H^{\srfs{LC}}_{\itGamma})^{\pm}(X,Y)$.
By Theorem~\ref{t-bddlip-1.8}(a),
$\tau \in \rfs{EXT}^{\pm}(X,Y)$.
By Theorem~\ref{ams-bddly-lip-bldr-t1.7},
$\tau \in (H^{\srfs{CMP.LC}}_{\itGamma})^{\pm}(X,Y)$.

(b) The proof is similar to the proof of Part (a).
However, we use Theorem~\ref{t-bddlip-1.8}(b)
and not \ref{t-bddlip-1.8}(a).

(c)  By Theorem~\ref{t2.4}(b), there is
$\tau \in H(X,Y)$ such that $\tau$ induces $\varphi$.
By Theorem~\ref{t-bddlip-1.8}(b),
$\tau \in \rfs{EXT}^{\pm}(X,Y)$.
By Theorem~\ref{ams-bddly-lip-bldr-t5.22},
$\itGamma = \itDelta$ and
$\tau \in (H^{\srfs{BDR.LC}}_{\itGamma})^{\pm}(X,Y)$.
\smallskip\hfill\myqed

\noindent
{\bf Proof of Theorem \ref{ams-bddly-lip-bldr-t1.4}(a) }
Theorem \ref{ams-bddly-lip-bldr-t1.4}(a) is restated as Part (a) of
\ref{ams-bddly-lip-bldr-t6.1} above.
\rule{0pt}{1pt}\hfill\myqed

\bigskip
\bigskip
\bigskip
\newpage

\newpage
\baselineskip 16.5pt
\begin{small}

\kern 1mm

\noindent
\noindent
{\Large\bf Index of symbols}

\kern1mm

\noindent
{\bf by order of appearance}

\kern1mm

\noindent
\indexentry{$\iso{f}{X}{Y}$.    This means $f$ is a homeomorphism between    $X$ and $Y$}{8}
\indexentry{$\iso{\varphi}{G}{H}$. This means    $\varphi$ is an isomorphism between $G$ and $H$}{8}
\indexentry{$\pair{a}{b}$. Notation of an ordered pair}{8}
\indexentry{$G \leq H$. $G$ is a subgroup of $H$}{12}
\indexentry{$f \leq g$. This means for every $t$,           $f(t) \leq g(t)$}{14}
\indexentry{$\alpha \preceq \beta\ \ \equiv\ \ $   for some $t > 0$, \, $\alpha \nrestriction [0,t] \leq   \beta \nrestriction [0,t]$}{14}
\indexentry{$g^{\sscirc n}$. Notation for           $g \scirc \ldots \scirc g$ \ $n$ times}{15}
\indexentry{$H\inverse = \setm{h\inverse}{h \in H}$}{17}
\indexentry{$\iso{\eta}{\rfs{MR}(X,G)}{\rfs{MR}(Y,H)}$}{34}
\indexentry{$G(x) = \setm{g(x)}{g\in G}$}{34}
\indexentry{$\sprtd{G}{B}$.           If $G \subseteq \setm{g}{\fnn{g}{A}{A}}$,           then $\sprtd{G}{B} \eqdf           \setm{g \in G}{g \nrestriction (A - B) = \rfs{Id}}$}{34}
\indexentry{$U \cong V$. This means $(\exists g \in G)(g(U) = V)$}{41}
\indexentry{$U \prec V$. This means $U$ is strongly small in $V$}{41}
\indexentry{$U \spr V$.           This means $U$ is strongly separated from $V$}{41}
\indexentry{$[x,y]$. The line segment with endpoints $x$ and $y$}{42}
\indexentry{$\rule{1.0pt}{0pt}\overE$.           The completion of a normed vector space $E$}{46}
\indexentry{$\bcirc \kern-1ptF =           \bigcup \setm{f \nrestriction \rfs{supp}(f)}{f \in F} \cup           \rfs{Id} \nrestriction   (Z - \bigcup \setm{\rfs{supp}(f)}{f \in F})$}{47}
\indexentry{$\bcirc_{n \in \sboldbbN} h_n$}{47}
\indexentry{$A_{\lambda}(x) = \setm{h_{\lambda}(x)}{h \in A}$}{64}
\indexentry{$\lambda \vdrest X$. The restriction of a partial group action $\lambda$ to an open set $X$}{65}
\indexentry{$\lambda \vdrest H_1$.           The restriction of a partial action $\lambda$ of $H$           to a subgroup $H_1$ of $H$}{66}
\indexentry{$a \approx^{\gamma} b$. This means that           $a \leq \gamma(b)$ and $b \leq \gamma(a)$}{70}
\indexentry{$f \sim^{\rho} g$. This means that           $\rfs{Dom}(f) \cup \rfs{Rng}(f) \subseteq \rfs{Dom}(\rho)$           and $g = \rho \scirc f \scirc \rho\inverse$}{70}
\indexentry{$U^{[n;W_1,W_2]}$}{88}
\indexentry{$g^{\sscirc n}$. This means $n$ times composition           of $g$}{95}
\indexentry{ {\thickmuskip=3mu \medmuskip=2mu \thinmuskip=1mu $g^{\srfs{cl}}_{M,N} = \setm{\pair{x}{y}} {x \in \rfs{cl}^M(A), \, y \in N \mbox{ and } g \cup \sngltn{\pair{x}{y}} \mbox{ is a continuous function}}$}}{111}
\indexentry{$g^{\srfs{cl}}$.           Abbreviation of $g^{\srfs{cl}}_{M,N}$}{111}
\indexentry{$g^{\srfs{cl}}_M$.           Abbreviation of $g^{\srfs{cl}}_{M,M}$}{111}
\indexentry{$H^{\srfs{cl}}  = \setm{h^{\srfs{cl}}}{h \in H}$}{111}
\indexentry{$A^{\geq n} = \setm{m \in A}{m \geq n}$}{141}
\indexentry{$A^{> n} = \setm{m \in A}{m > n$}}{141}
\indexentry{$A^{\leq n} = \setm{m \in A}{m \leq n$}}{141}
\indexentry{$A^{< n} = \setm{m \in A}{m < n$}}{141}
\indexentry{           {\thickmuskip=2mu \medmuskip=1mu \thinmuskip=1mu           $x \simeq^X y$. This means that $x$ and $y$ lie in the           same connected component of $X$}}{185}
\indexentry{$\vecx \simeq^X \vecy$.           This means that for every $n$, $x_n \simeq^X y_n$}{185}
\indexentry{$f^{\srfs{eni}} =           f^{\srfs{cl}} \nrestriction \rfs{ENI}(\rfs{cl}(X))$}{197}
\indexentry{$a \approx^K b$. This means   $\frac{1}{K} a \leq b \leq Ka$}{233}
\indexentry{$\norm{\ }^1 \approx^K \norm{\ }^2$.   This means for every $u \in E$,   $\norm{u}^1 \approx^K \norm{u}^2$}{233}
\indexentry{$E = L \oplus^{\srfs{alg}} S$.   The algebraic direct sum}{233}
\indexentry{$(x)_{L,S}$. The $L$-component of $x$ in   $L \oplus S$}{233}
\indexentry{$(x)_L$. Abbreviation of $(x)_{L,S}$}{233}
\indexentry{$\norm{u}^{L,S} =   \norm{(u)_S} + \norm{(u)_L}$}{233}
\indexentry{$H \perp^M F$.   This means for every $u \in H$,           $d(u,F) \geq \frac{1}{M} \norm{u}$}{233}
\indexentry{$s \approx^{\alpha} t$.           This means $t \leq \alpha(s)$ and $s \leq \alpha(t)$}{241}
\indexentry{ $u \perp F$.           This means: $F$ is linear subspace of a normed space $E$,           $u \in E$ and $\norm{u} = d(u,F)$}{273}
\indexentry{ $K \perp F$.           This means: $K,F$ are linear subspaces of $E$           and for every $u \in K$, $u \perp F$}{273}
\indexentry{$\vecx^{\srsp{\sigma}}$. A sequence whose domain is           $\sigma \subseteq \bbN$}{281}
\indexentry{$\sigma^{\geq n} =                 \setm{k \in \sigma}{k \geq n}$}{281}
\indexentry{ $\vecx^{\kern1pt\geq n} = \vecx \nrestriction \rfs{Dom}(\vecx)^{\geq n}$}{281}
\indexentry{$x \approx^{(\alpha,b)}_{(X,E)} y$. This means           $d(x,b) \approx^{\alpha} d(y,b)$ and           $\delta^X(x) \approx^{\alpha} \delta^X(y)$}{281}
\indexentry{$\vecx^{\srsp{\sigma}} \approx^{(\alpha,b)}_{(X,E)}           \vecy^{\srsp{\sigma}}$. This means: for every $n \in \sigma$, $x_n \approx^{(\alpha,b)}_{(X,E)} y_n$}{281}
\indexentry{$\approx^{(\alpha,b)}$.           Abbreviation of \,$\approx^{(\alpha,b)}_{(X,E)}$}{281}
\indexentry{$\vecx^{\srsp{\rho}} \neweq^A \vecy^{\srsp{\rho}}$}{281}
\indexentry{$\vecx^{\srsp{\rho}} \neweq^{\alpha}           \vecy^{\srsp{\rho}}$}{281}
\indexentry{$\calP \kern-1pt\restriction\kern-1pt T =   \setm{P \cap T}{P \in \calP}$}{296}
\indexentry{$a \sim^{\calP} b$.   This means: there is $P \in \calP$ such that $a,b \in P$}{296}
\indexentry{ $\beta \rsstar \alpha(t) = \sum_{n = 0}^{\infty}    \beta(q_{\alpha,n}(t))$}{310}
\end{small}

\baselineskip 16.5pt
\begin{small}
\noindent
{\Large\bf Index of notations}

\kern 1mm

\noindent
{\bf by alphabetic order}

\kern 1mm

\baselineskip 16.5pt
\noindent
\indexentry{$\bbA(E)$. The group of affine automorphisms of $E$}{65}
\indexentry{$\bbA(E;F) = \setm{A \in \bbA(E)}{A(F) = F}$}{66}
\indexentry{$\rfs{ABUC}(X,Y) =           \setm{h \in H(X,Y)}{\mbox{for every bounded set }           A \subseteq X,\ h \nrestriction A \mbox{ is UC}}$}{134}
\indexentry{$\rfs{acc}^X(U)$ the set of accumulation points           of $U$ in $X$}{33}
\indexentry{$B^Z(x,r) = \setm{y \in Z}{\,d(x,y) < r}$}{41}
\indexentry{$\rule{0pt}{0pt}\overB^E(x,r) =           \setm{y \in E}{\,d(x,y) \leq r}$}{41}
\indexentry{$B^Z(A,r) = \bigcup_{x \in A} B^Z(x,r)$}{41}
\indexentry{\rule{0pt}{0pt}$B(x,r)$.           An abbrviation of $B^X(x,r)$}{42}
\indexentry{$B(x;r,s) = \setm{y \in X}{r < d(x,y) < s}$}{233}
\indexentry{$\rfs{BCD}^E(A,r) = B^E(0,r) - A$.           A boundary chart domain based on $E$ and $A$ with radius $r$}{226}
\indexentry{$\rfs{bd}^X(U)$ boundary of $U$ in $X$}{33}
\indexentry{$\rfs{BUC}(X,Y) = \setm{g \in H(X,Y)}    {g \mbox{ is boundedness preserving and } g \nrestriction A    \mbox{ is UC for every}    \newline\phantom{\hbox{$\rfs{BUC}(X,Y) =\ $}}\mbox{bounded set }    A \subseteq X}$}{18}
\indexentry{$\rfs{cl}_{\spreceq}(\itGamma) =    \setm{\alpha \in \rfs{MC}}{\mbox{ for some } \gamma \in \itGamma,\ \,    \alpha \preceq \gamma}$.}{15}
\indexentry{$\rfs{cl}^X(U)$ closure of $U$ in $X$}{33}
\indexentry{$\rfs{Cmp}(X)$. The set of connected components of $X$}{174}
\indexentry{$\rfs{CMP.LUC}(X)$. The group of biextendible           homeomorphisms of $X$ which are bi-uniformly   \newline\indent   continuous at every $x \in \rfs{cl}(X)$}{10}
\indexentry{$\rfs{co-dim}^E(L)$. Co-dimension of $L$ in $E$.   Abbreviation $\rfs{co-dim}(L)$}{233}
\indexentry{$\delta^X(A) = d(A,E - X)$}{123}
\indexentry{$\delta^X(x) = \delta^X(\sngltn{x})$}{123}
\indexentry{$\delta(A)$. Abbreviation of $\delta^X(A)$}{123}
\indexentry{$\delta(x)$. Abbreviation of $\delta^X(x)$}{123}
\indexentry{$\delta_1^X(x)$}{141}
\indexentry{$\Delta(A)$. $\sup_{a \in A} d(a,E - X)$}{208}
\indexentry{$\delta^{X,E}(x) = d(x,E - X)$. Abbreviations:           $\delta^X(x)$, $\delta(x)$}{281}
\indexentry{$\rfs{diam}(A) = \sup_{x,y \in A} d(x,y)$}{53}
\indexentry{$e_H$. The unit of a group $H$}{63}
\indexentry{$\rfs{ENI}(X,\itPhi,G) =    \setm{g(x)}{x \in \rfs{NI}(X,\itPhi) \mbox{ and } g \in G}$.    Extended normed interior of $\trpl{X}{\itPhi}{G}$}{59}
\indexentry{$\rfs{ENI}(X,\itPhi) = \rfs{ENI}(X,\itPhi,H(X))$}{59}
\indexentry{$\rfs{ENI}(X) =    \setm{h(x)}{x \in \rfs{int}^E(X),\ h \in H(X)}$}{59}
\indexentry{$\eta_{(\rho,a)}(s,t)$}{242}
\indexentry{$\rfs{EXT}(X) =           \setm{g \in H(X)}   {g \mbox { and } g\inverse \mbox{ are extendible}}$}{20}
\indexentry{$\rfs{EXT}^{M,N}(X,Y) =    \setm{h \in H(X,Y)}{\rfs{Dom}(h^{\srfs{cl}}_{M,N}) = \rfs{cl}^M(X)}$}{111}
\indexentry{$\rfs{EXT}(X,Y)$.           Abbreviation of $\rfs{EXT}^{M,N}(X,Y)$}{111}
\indexentry{$\rfs{EXT}^M(X)$.           Abbreviation of $(\rfs{EXT}^{M,M})^{\pm}(X,X)$}{111}
\indexentry{$\rfs{FD}(X)$}{60}
\indexentry{$\rfs{FD.LIP}(X)$}{60}
\indexentry{$\rfs{Fld}(\lambda) = \rfs{Dom}(e_{\lambda})$}{64}
\indexentry{$\itGamma_{\alpha} =           \rfs{cl}_{\spreceq}(\setm{\alpha^{\sscirc n}}{n \in \bbN})$}{310}
\indexentry{$\itGamma^{\srfs{HLD}}_r =    \setm{\alpha \in \rfs{MC}}{\mbox{ for some $K > 0$, }    \alpha \preceq K x^r}$}{15}
\indexentry{$\itGamma^{\srfs{HLD}} =           \bigcup \setm{\itGamma^{\srfs{HLD}}_r}{r \in (0,1]}$.           The H\"{o}lder modulus}{15}
\indexentry{$\itGamma^{\srfs{LIP}} =    \setm{\alpha \in \rfs{MC}}{\mbox{ for some $K > 0$, }    \alpha \preceq Kx}$}{15}
\indexentry{$H(X)$. The group of all auto-homeomorphisms    of $X$}{8}
\indexentry{$H_{\itGamma}(X)$. The group of all    $\itGamma$-bicontinuous auto-hoeomorphisms of $X$}{17}
\indexentry{$H(X;F) = \setm{h \in H(X)}{h(X \cap F) = X \cap F}$}{66}
\indexentry{$H_{\itGamma}(X) =           \setm{h \in H(X)}{\mbox{there is } \gamma \in \itGamma           \mbox{ such that } h \mbox{ is } \gamma\mbox{-bicontinuous}}$}{81}
\indexentry{$H_{\itGamma}(X;F) =           \setm{h \in H_{\itGamma}(X)}{h(F \cap X) = F \cap X}$}{81}
\indexentry{$H_{\itGamma}(X;S,F) =           H_{\itGamma}(X,S) \cap H_{\itGamma}(X;F)$}{81}
\indexentry{$H_{\itGamma}(X;\calS,\calF)$.           The subgroup of $H(X)$ generated by           {\thickmuskip=2mu \medmuskip=1mu \thinmuskip=1mu           $\bigcup \setm{H_{\itGamma}(X;S,F_S)}{S \in \calS}$}}{81}
\indexentry{$H(X,Y) = \setm{h}{\iso{h}{X}{Y}}$}{109}
\indexentry{$H(X;D) = \setm{h \in H(X)}{h(D) = D}$}{109}
\indexentry{$H_{\itGamma}(X,S) = H_{\itGamma}(X)\sprt{S}$}{81}
\indexentry{$H_{\itGamma}^{\srfs{BD}}(X)$}{27}
\indexentry{$H^{\srfs{BDR.LC}}_{\itGamma}(X) =           \setm{g \in \rfs{EXT}(X)} {\mbox{for every }   x \in \rfs{bd}(X),\ g \mbox{ is }   \itGamma\mbox{-bicontinuous at } x}$}{229}
\indexentry{ $H^{\srfs{CMP.LC}}_{\itDelta,\itGamma}(X) =           H^{\srfs{LC}}_{\itDelta}(X) \cap           H^{\srfs{BDR.LC}}_{\itGamma}(X)$}{229}
\indexentry{$H_{\itGamma}^{\srfs{BPD}}(X,Y) =           \setm{f \in \rfs{BPD.P}(X,Y)}{\kern1pt   \mbox{for every BPD set } A \subseteq X,           f \nrestriction A \mbox{ is } \itGamma\mbox{-continuous}}$}{140}
\indexentry{$H^{\srfs{CMP.LC}}_{\itGamma}(X) =           \setm{g \in \rfs{EXT}(X)}{(\forall x \in \rfs{cl}^E(X))           (\exists U \in \rfs{Nbr}^E(X))(\exists \alpha \in \itGamma)           (g \nrestriction (U \cap X)\\\rule{27.3mm}{0pt}           \mbox{ is } \alpha\mbox{-bicontinuous)}}$}{227}
\indexentry{$H^{\srfs{CMP.LC}}_{\itGamma}(X,Y) =           \setm{g \in \rfs{EXT}(X,Y)}{g^{\srfs{cl}}           \mbox{ is locally $\itGamma$-continuous}}$}{23}
\indexentry{$H_{\itGamma}^{\srfs{LC}}(X)$.   The group of locally $\itGamma$-bicontinuous   auto-homeomorphisms of $X$}{16}
\indexentry{$H_{\itGamma}^{\srfs{LC}}(X)$.           The group of locally $\itGamma$-bicontinuous homeomorphisms   of $X$}{81}
\indexentry{$H_{\itGamma}^{\srfs{LC}}(X,S) =           H_{\itGamma}^{\srfs{LC}}(X)\sprt{S}$}{81}
\indexentry{$H_{\itGamma}^{\srfs{LC}}(X;F) =           \setm{h \in H_{\itGamma}^{\srfs{LC}}(X)}{h(F \cap X) =           F \cap X}$}{81}
\indexentry{$H_{\itGamma}^{\srfs{LC}}(X;S,F) =           H_{\itGamma}^{\srfs{LC}}(X,S) \cap   H_{\itGamma}^{\srfs{LC}}(X;F)$}{81}
\indexentry{ {\thickmuskip=2mu \medmuskip=1mu \thinmuskip=1mu           $H_{\itGamma}^{\srfs{LC}}(X;\calS,\calF)$.} The subgroup $H(X)$ generated by {\thickmuskip=2mu \medmuskip=1mu \thinmuskip=1mu $\bigcup \setm{H_{\itGamma}^{\srfs{LC}}(X;S,F_S)}{S \in \calS}$}}{81}
\indexentry{$H_{\itGamma}^{\srfs{LC}}(X,\itPhi) =    \setm{h \in H(X)}    {\forall x(\exists \varphi,\psi \in \itPhi)    (x \in \rfs{int}(\rfs{Rng}(\varphi)),\    h(x) \in \rfs{int}(\rfs{Rng}(\psi))    \mbox{ and }    \newline\rule{25mm}{0pt}\psi\inverse \scirc h \scirc \varphi    \mbox{ is $\itGamma$-bicontinuous at $\varphi\inverse(x)$}}$}{106}
\indexentry{$H^{\srfs{LC}}_{\itGamma}(X,\itPhi,\calS)$}{106}
\indexentry{$H_{\itGamma}^{\srfs{NBPD}}(X,Y) = \setm{h \in BPD.P(X,Y)}{h \mbox{ is nearly $\itGamma$-continuous} \mbox{ on BPD sets}}$}{140}
\indexentry{$H_{\itGamma}^{\srfs{PW}}(X)$}{27}
\indexentry{$H_{\itGamma}^{\srfs{RG}}(X)$}{27}
\indexentry{$H_{\itGamma}^{\srfs{WBPD}}(X,Y) =           \setm{h \in BPD.P(X,Y)}{h   \mbox{ is weakly $\itGamma$-continuous} \mbox{ on BPD sets}}$}{141}
\indexentry{$\rfs{int}^X(U)$ interior of $U$ in $X$}{33}
\indexentry{           $\overline{\rfs{int}}^E(X) =           \bigcup \setm{B^{\oversE}(x,r)}{x \in X \mbox{ and }           B^E(x,r) \subseteq X}$}{53}
\indexentry{$\rfs{IXT}^E(X)$. The group of bi-externally-extendible           auto-homeomorphisms of $X$}{53}
\indexentry{$\itbfK_{\calM,\calP}$. The category $\pair{\calM}{\setm{g}{\iso{g}{X}{Y},\ X,Y \in \calM \mbox{ and } g,g\inverse \mbox{ have property } \calP}}$ }{11}
\indexentry{$\itbfK_{\itGamma}$}{16}
\indexentry{$\Kappa^X(x,A) =           \setm{\kappa}{(\forall U \in \rfs{Nbr}(x))           (\exists B \subseteq A \cap U)           (\abs{B} = \kappa \mbox{ and } B \mbox{ is spaced}}$}{70}
\indexentry{$\kappa^X(x,A) = \sup(\Kappa^X(x,A))$}{70}
\indexentry{$\kappa(X) = \min_{x \in X} \kappa^X(x,X)$}{70}
\indexentry{$K_{\srfs{arc}}(\ell,t)$}{145}
\indexentry{$K_{\srfs{B}} =           \setm{\pair{X}{G}}{X \mbox{ is an open subset of a           Banach space and } \rfs{LIP}(X) \leq G \leq H(X)}$}{38}
\indexentry{$K_{\srfs{BM}}$}{59}
\indexentry{$K_{\srfs{BNM}}$}{59}
\indexentry{$K_{\srfs{BNO}}$}{38}
\indexentry{$K_{\srfs{BO}}$}{38}
\indexentry{$K_{\srfs{LCM}} =           \setm{\pair{X}{G}} { \mbox{ is Haussdorf, perfect           locally compact and for every}\newline\indent           \rule{25pt}{0pt}\mbox{ open } V \subseteq X \mbox{ and }           x \in V,\  G\sprt{V}(x) \mbox{ is somewhere dense}}$}{34}
\indexentry{$K_{\srfs{N}} =           \setm{\pair{X}{G}}{X \mbox{ is an open subset of a           normed space and }           \rfs{LIP}^{\srfs{LC}}(X) \leq G \leq H(X)}$}{38}
\indexentry{$K_{\srfs{NFCB}}$}{56}
\indexentry{$K_{\srfs{NL}}$}{53}
\indexentry{$K_{\srfs{NM}}$}{59}
\indexentry{$K_{\srfs{NO}}$}{38}
\indexentry{$K^{\calN\calO}_{\srfs{NMX}} =    \setm{\pair{X}{Z} \in  K^{\calN\calO}_{\srfs{NRM}}}    {X \mbox{ is BR.LC.AC and JN.AC with repect to } Z}$}{202}
\indexentry{    $K^{\calN\calO}_{\srfs{NRM}} =    \setm{\pair{X}{Z}}{X \in K^{\calO}_{\srfs{NRM}}    \mbox{ and } X \subseteq Z \subseteq \rfs{cl}(X)}$}{202}
\indexentry{$K_{\srfs{BCX}}^{\calO}$}{171}
\indexentry{           $K^{\calO}_{\srfs{BLPM}} =   \setm{Y}{Y   \mbox{ is an open subset of a Banach Lipschitz manifold}}$}{220}
\indexentry{$K^{\calO}_{\srfs{BNC}}$.           Class of all spaces which are an open subset of a   Banach space}{122}
\indexentry{$K_{\srfs{BX}}^{\calO}$}{185}
\indexentry{$K_{\srfs{IMX}}^{\calO}$.           Class of open finite-dimensional BR.IS.MV open sets}{192}
\indexentry{$K^{\calO}_{\srfs{NFCB}}$.           Class of open subsets of first category or complete   normed spaces}{122}
\indexentry{$K^{\calO}_{\srfs{NLPM}} =           \setm{Y}{Y \mbox{ is an open subset of a normed Lipschitz           manifold}}$}{204}
\indexentry{$K_{\srfs{NMX}}^{\calO}$}{171}
\indexentry{$K^{\calO}_{\srfs{NRM}}$.           Class of all spaces which are an open subset of   a normed space}{122}
\indexentry{$K_{\srfs{WFD.BNO}}$}{61}
\indexentry{$\bbL(E)$. The group of bounded linear automorphisms           of $E$}{65}
\indexentry{$\bbL(E,x) =           (\bbL(E))^{\srfs{tr}_x^E}$}{65}
\indexentry{$\bbL(E;F) = \setm{T \in \bbL(E)}{T(F) = F}$}{66}
\indexentry{$\lambda_{\bbA}^{E;F} =           \lambda_{\bbA}^E \vdrest \bbA(E;F)$}{66}
\indexentry{$\limti{i} A_i = x$. The limit of a sequence of sets}{124}
\indexentry{$\rfs{LIP}(X)$. The group of bilipschitz    auto-homeomorphisms of a metric space $X$}{12}
\indexentry{$\rfs{LIP}(X,S)$. For $S \subseteq X$,           $\rfs{LIP}(X,S) =           \setm{h \in \rfs{LIP}(X)}   {h \nrestriction (X - S) = \rfs{Id}}$}{12}
\indexentry{$\rfs{LIP}(X;F)$.           For a normed space $E$,   $X \subseteq E$ and a dense linear subspace           \newline\indent           $F$ of $E$, \           $\rfs{LIP}(X;F) =   \setm{h \in \rfs{LIP}(X)}{h(X \cap F) = X \cap F}$}{12}
\indexentry{$\rfs{LIP}(X;S,F) =   \rfs{LIP}(X;F) \cap \rfs{LIP}(X,S)$}{12}
\indexentry{$\rfs{LIP}(X;\calS,\calF)$.   Subgroup of $H(X)$ generated by           $\bigcup \setm{\rfs{LIP}(X;S,F_S)}{S \in \calS}$}{38}
\indexentry{$\rfs{LIP}(X,\calS)$. The subgroup of $H(X)$ generated by $\bigcup \setm{\rfs{LIP}(X,S)}{S \in \calS}$}{38}
\indexentry{$\rfs{LIP}(X;\itPhi,\calF)$}{58}
\indexentry{$\rfs{LIP}(X;\itPhi)$}{58}
\indexentry{$\rfs{LIP}(X,\itPhi) = \setm{h \in H(X)}    {\exists K \forall x(\exists \varphi,\psi \in \itPhi)    (x \in \rfs{int}(\rfs{Rng}(\varphi)),\    h(x) \in \rfs{int}(\rfs{Rng}(\psi))    \mbox{ and }    \newline\rule{24mm}{0pt}\psi\inverse \scirc h \scirc \varphi    \mbox{ is $K$-bilipschitz}}$}{105}
\indexentry{$\rfs{LIP}(X,\itPhi,\calS)$}{106}
\indexentry{$\rfs{LIP}_{00}(X) =    \setm{f \in \rfs{LIP}(X)}{\rfs{supp}(f)    \mbox{ is a BPD set}}$}{141}
\indexentry{$\rfs{LIP}^{\srfs{LC}}(X)$.           The group of locally bilipschitz auto-homeomorphisms of $X$}{12}
\indexentry{$\rfs{LIP}^{\srfs{LC}}(X,S)$. For $S \subseteq X$,           $\rfs{LIP}(X,S)^{\srfs{LC}} =           \setm{h \in \rfs{LIP}^{\srfs{LC}}(X)}{h \nrestriction (X - S)   = \rfs{Id}}$}{12}
\indexentry{$\rfs{LIP}^{\srfs{LC}}(X;F)$.           For a normed space $E$, $X \subseteq E$ and a dense linear   subspace           \newline\indent           $F$ of $E$, \           $\rfs{LIP}^{\srfs{LC}}(X;F) =           \setm{h \in \rfs{LIP}^{\srfs{LC}}(X)}{h(X \cap F) =   X \cap F}$}{12}
\indexentry{$\rfs{LIP}^{\srfs{LC}}(X;S,F) =           \rfs{LIP}(X;F) \cap \rfs{LIP}(X,S)$}{12}
\indexentry{           {\thickmuskip=2mu \medmuskip=1mu \thinmuskip=1mu   $\rfs{LIP}^{\srfs{LC}}(X;\calS,\calF)$.   Subgroup of $H(X)$ generated by           $\bigcup \setm{\rfs{LIP}^{\srfs{LC}}(X;S,F_S)}{S \in \calS}$}}{38}
\indexentry{$\rfs{LIP}^{\srfs{LC}}(X,\calS)$. The subgroup of $H(X)$ generated by $\bigcup \setm{\rfs{LIP}^{\srfs{LC}}(X,S)}{S \in \calS}$}{38}
\indexentry{$\rfs{LIP}^{\srfs{LC}}(X;\itPhi,\calF)$}{58}
\indexentry{$\rfs{LIP}^{\srfs{LC}}(X;\itPhi)$}{58}
\indexentry{$\lambda_{\bbL}^{E;F} =           \lambda_{\bbL}^E \vdrest \bbL(E;F)$}{66}
\indexentry{$\lambda_{\bbL}^{E,x;F} =           \lambda_{\bbL}^{E,x} \vdrest \bbL(E,x;F)$}{66}
\indexentry{$\rfs{lngth}(L)$. Length of an arc}{124}
\indexentry{$\lambda_{\bbT}^E$,           $\lambda_{\bbL}^E$, $\lambda_{\bbL}^{E,x}$,   $\lambda_{\bbA}^E$. Actions of $\bbT(E)$,   $\bbL(E)$, $\lambda_{\bbL}^{E,x}$ and $\bbA(E)$ on $E$}{66}
\indexentry{$\lambda_{\bbT}^{E;F} =           \lambda_{\bbT}^E \vdrest \bbT(E;F)$}{66}
\indexentry{$\rfs{LUC}(X,Y) =    \setm{h \in H(X,Y)}{h \mbox{ is locally UC}}$}{110}
\indexentry{$\rfs{LUC}^{\wpm}(X,Y) =    \setm{h \in H(X,Y)}{h \mbox{ is locally bi-UC}}$}{110}
\indexentry{$\rfs{LUC}(X) = \rfs{LUC}^{\wpm}(X,X)$}{110}
\indexentry{$\rfs{LUC}_{01}(X) =           \setm{h \in \rfs{LUC}(X)}{(\exists U)( U \mbox{ is $E$-open, }           U \supseteq \rfs{bd}(X) \mbox{ and }           \rfs{supp}(h) \subseteq X - U)}$}{172}
\indexentry{$M(X,G) =           \semisixtpl{X}{\tau^X}{G}{\in}{\scirc}{\rfs{Ap}}$}{291}
\indexentry{$M^{\srfs{aoc}}(n)$}{233}
\indexentry{$M^{\srfs{arc}}(t)$}{235}
\indexentry{$\rfs{MBC} =           \setm{\alpha \in \rfs{MC}}{\rfs{Id}_{[0,\infty)} \leq   \alpha}$}{63}
\indexentry{$M^{\srfs{bnd}}(K)$}{239}
\indexentry{$\rfs{MC} =    \setm{h \in H([0,\infty))}{h \mbox{ is concave}}$}{14}
\indexentry{$M^{\srfs{cmp}}$}{240}
\indexentry{$M^{\srfs{fdn}} = M^{\srfs{fdn}}(2)$}{235}
\indexentry{$M^{\srfs{fdn}}(n)$}{235}
\indexentry{$M^{\srfs{hlb}} = M^{\srfs{hlb}}(2)$}{233}
\indexentry{$M^{\srfs{hlb}}(n)$}{233}
\indexentry{$M^{\srfs{lift}}$}{276}
\indexentry{$M^{\srfs{ort}} = M^{\srfs{ort}}(2)$}{235}
\indexentry{$M^{\srfs{ort}}(n)$}{235}
\indexentry{$M^{\srfs{prj}}(n)$}{235}
\indexentry{$\rfs{MR}(X,G) =           \sixtpl{\rfs{Ro}(X)}{G}{+}{\cdot}{-}{\rfs{Ap}}$}{34}
\indexentry{$M^{\srfs{rot}}$}{236}
\indexentry{$M^{\srfs{rtn}}$}{304}
\indexentry{$M^{\srfs{seg}}$}{235}
\indexentry{$M^{\srfs{thn}} = M^{\srfs{thn}}(2)$}{233}
\indexentry{$M^{\srfs{thn}}(n)$}{233}
\indexentry{$\bbN^+ = \setm{n \in \bbN}{n > 0}$}{211}
\indexentry{$\rfs{Nbr}^X(x) = \setm{U}{x \in U \subseteq X           \mbox{ and $U$ is open}}$}{63}
\indexentry{$\rfs{NI}(X,\itPhi) =      \bigcup \setm{\varphi(B^{E_{\varphi}}(x_{\varphi},r_{\varphi}))}      {\varphi \in \itPhi}$. The normed interior of      $\pair{X}{\itPhi}$}{59}
\indexentry{$\rfs{opcl}(U) =          \rfs{int}^{\srfs{cl}(X)}(\rfs{cl}^{\srfs{cl}(X)}(U))$}{294}
\indexentry{$\calP(X,Y) = \setm{h}{\iso{h}{X}{Y}           \mbox{ and $h$ has property }\calP}$}{17}
\indexentry{$\calP^{\pm}(X,Y) = \calP(X,Y) \cap           (\calP(Y,X))\inverse$}{17}
\indexentry{$\calP(X) = \calP^{\pm}(X,X)$}{17}
\indexentry{$p_{\alpha,n}(t) =           ((\rfs{Id} + \alpha)\inverse)^{\sscirc n}(t)$}{310}
\indexentry{$\rule{1pt}{0pt}\rfs{PNT.UC}(X,x) = \setm{h \in H(X)}           {h(x) = x \mbox{ and } h \mbox{ is bi-UC at } x}$}{109}
\indexentry{$\rule{1.9pt}{0pt}q_{\alpha,n}(t) =           p_{\alpha,n}(t) - p_{\alpha,n + 1}(t)$}{310}
\indexentry{$R(u,v,g;\alpha,a,b,F)$}{243}
\indexentry{$R(u,v,g;M,a,b,F)$}{243}
\indexentry{$\rfs{Rad}^E_{\eta,z} =           z + \eta(\norm{x - z}) \frac{x - z}{\norm{x - z}}$.           The radial homeomorphism based on $\eta,z$}{81}
\indexentry{$\rfs{Rad}^E_{\eta} = \rfs{Rad}^E_{\eta,0^E}$}{81}
\indexentry{$\rfs{Ro}(X)$. The set of regular open subsets of $X$}{33}
\indexentry{$\rfs{Rot}^H_{\theta}$.   In a $2$-dimensional Hilbert space $H$,   rotation by the angle $\theta$}{235}
\indexentry{$\rfs{Rot}^{F,H}_{\theta}$. For a   $2$-dimensional Hilbert space $H$ and a normed space $F$,   the operator on\\   \indent$H \oplus F$ which is   $\rfs{Rot}^H_{\theta}$ on $H$ and $\rfs{Id}$ on $F$}{235}
\indexentry{$S^Z(x,r) = \setm{y \in Z}{\,d(x,y) = r}$}{41}
\indexentry{$S_{\calP} = \bigcup \calP$}{296}
\indexentry{$\rfs{supp}(h) = \setm{y \in Y}{h(y) \neq y}$}{42}
\indexentry{$\bbT(E) = \setm{\rfs{tr}_v^E}{v \in E}$.           The group of translations of $E$}{65}
\indexentry{$\bbT(E;F) = \setm{\rfs{tr}^E_v}{v \in F}$}{66}
\indexentry{$\rfs{tr}_v^E$. Translation by $v$. For $v,x \in E$,           $\rfs{tr}_v^E(x) = v + x$}{42}
\indexentry{$\rfs{tr}_v$. Abbreviation of $\rfs{tr}_v^E$}{42}
\indexentry{$\rfs{UC}(X,Y) = \setm{h \in H(X)}{h \mbox{ is UC}}$}{109}
\indexentry{$\rfs{UC}^{\wpm}(X,Y) =           \setm{h \in H(X)}{h \mbox{ is bi-UC}}$}{109}
\indexentry{$\rfs{UC}(X) = \rfs{UC}^{\wpm}(X,X)$}{109}
\indexentry{$\rfs{UC}(X;F) =           \setm{h \in \rfs{UC}(X)}{h(X \cap F) = X \cap F}$}{109}
\indexentry{$\rfs{UC}(X;S,F) =           \rfs{UC}(X)\sprt{S} \cap \rfs{UC}(X;F)$}{109}
\indexentry{$\rfs{UC}(X;S,F,x) =           \setm{h \in \rfs{UC}(X;S,F)}{h(x) = x}$}{109}
\indexentry{$\rfs{UC}(X,\calS)$.    The subgroup of $H(X)$ generated by    $\bigcup \setm{\rfs{UC}(X)\sprt{S}}{S \in \calS}$}{110}
\indexentry{$\rfs{UC}(X;\calS,\calF)$.    {\thickmuskip=2mu \medmuskip=1mu \thinmuskip=1mu    The subgroup of $H(X)$ generated by    $\bigcup \setm{\rfs{UC}(X;S,F_S)}{S \in \calS}$}}{110}
\indexentry{$\rfs{UC}_0(X) =           \setm{f \in \rfs{UC}(X) \cap \rfs{EXT}(X)}{ f^{\srfs{cl}} \nrestriction \rfs{bd}(X) = \rfs{Id}}$}{124}
\indexentry{$\rfs{UC}_{00}(X) =    \setm{f \in \rfs{UC}(X)}{\rfs{supp}(f)    \mbox{ is a BPD set}}$}{141}
\indexentry{$\rfs{UC}_{\rme}(X) =           \setm{h \in \rfs{UC}(X)}{h \mbox{ is strongly extendible}}$}{208}
\indexentry{$\rfs{UC}_0^{\srfs{eni}}(X) =           \setm{f^{\srfs{eni}}}{f \in  \rfs{UC}_0(X)}$}{197}
\indexentry{$\rfs{WFD}(X)$}{61}
\indexentry{$\rfs{WFD.LIP}^{\srfs{LC}}(X;\calS,\calF)$}{61}
\indexentry{$\rfs{WFD.LIP}(X)$}{61}
\indexentry{$\rfs{WFD.LIP}(X;\calS,\calF)$}{61}
\end{small}

\baselineskip 16.5pt
\begin{small}
\noindent
{\Large\bf Index of definitions}

\kern 1mm

\noindent
{\bf by alphabetic order}

\kern 1mm

\noindent
\indexentry{abiding sequence. $\alpha$-abiding sequence}{282}
\indexentry{affine-like partial action at $x$}{88}
\indexentry{affine-like partial action}{88}
\indexentry{almost $\alpha$-continuous at $x$}{95}
\indexentry{almost $\alpha$-continuous}{95}
\indexentry{almost $\beta$-continuous for $\alpha$-submerged pairs.           Abbreviation: $\scolonpair{\beta}{\alpha}$-almost-continuous}{306}
\indexentry{almost $\itGamma$-continuous at $x$}{95}
\indexentry{almost linear boundary chart domain}{231}
\indexentry{almost orthogonal complement}{233}
\indexentry{appropriate. A $\itGamma$-appropriate group}{229}
\indexentry{Banach manifold}{58}
\indexentry{BD.AC. Abbreviation of boundedly arcwise connected}{169}
\indexentry{BD.CW.AC. Abbreviation of           boundedly component-wise arcwise connected}{185}
\indexentry{BDD.P function. A function which takes bounded   sets to bounded sets}{124}
\indexentry{BDR.UC function}{124}
\indexentry{bi-UC at $x$. Abbreviation of bi-uniformly-continuous at $x$}{109}
\indexentry{bi-UC. Abbreviation of bi-uniformly-continuous}{108}
\indexentry{bi-uniformly-continuous at $x$}{109}
\indexentry{bi-uniformly-continuous}{108}
\indexentry{bicontinuous. $\alpha$-bicontinuous at $x \in \rfs{cl}(X)$}{281}
\indexentry{bicontinuous. $\alpha$-bicontinuous at $x$}{63}
\indexentry{bicontinuous. $\alpha$-bicontinuous homeomorphism}{63}
\indexentry{bicontinuous. $\itGamma$-bicontinuous at $x \in \rfs{cl}(X)$}{281}
\indexentry{bicontinuous. $\itGamma$-bicontinuous at $x$}{63}
\indexentry{bicontinuous. $\itGamma$-bicontinuous.    $h$ is $\itGamma$-bicontinuous,    if $(\exists \gamma \in \itGamma)(    h,\ h\inverse \mbox{ are } \gamma\mbox{-continuous})$}{17}
\indexentry{bicontinuous. $\pair{K}{\calP}$-bicontinuous}{296}
\indexentry{bilipschitz homeomorphism between    locally Lipschitz normed manifolds}{105}
\indexentry{bilipschitz homeomorphism}{12}
\indexentry{BI.UC function}{125}
\indexentry{BNO-system}{38}
\indexentry{boundary chart domain based on $E$ and $A$           with radius $r$. A set of the form $\rfs{BCD}^E(A,r)$}{226}
\indexentry{boundary chart element}{226}
\indexentry{boundary type. A group of boundary type $\itGamma$}{229}
\indexentry{bounded positive distance UC function}{124}
\indexentry{boundedly $\itGamma$-continuous}{27}
\indexentry{boundedly UC function. A function which is uniformly   continuous on every bounded set}{124}
\indexentry{boundedly arcwise connected. Abbreviated by BD.AC}{169}
\indexentry{boundedly component-wise arcwise connected, $X$ is}{185}
\indexentry{boundedly uniformly\,-\,in\,-\,diameter arcwise\,-\,connected}{135}
\indexentry{boundedness preserving function.   A function which takes bounded sets to bounded sets}{124}
\indexentry{BPD sequence. A sequence $\vecx$ such that    $\rfs{Rng}(\vecx)$ is a BPD set}{124}
\indexentry{BPD set. A bounded subset of $X$ whose distance    from the boundary of $X$ is positive}{124}
\indexentry{BPD-arcwise-connected}{141}
\indexentry{BPD.AC. Abbreviation of BPD-arcwise-connected}{141}
\indexentry{BPD.P function. A function which takes BPD sets to   BPD sets}{124}
\indexentry{BPD.UC function.   A function which is uniformly continuous on every BPD set}{124}
\indexentry{BR.CW.LC.AC. Abbreviation of $X$ is           component-wise locally arcwise connected at           \newline\indent           its boundary}{185}
\indexentry{BR.IS.MV. Abbreviation of    isotopically movable at the boundary}{191}
\indexentry{BR.LC.AC. Abbreviation of           locally arcwise connected at the boundary}{157}
\indexentry{BR.LUC function}{124}
\indexentry{BUC function. A function which is uniformly   continuous on every bounded set}{124}
\indexentry{BUD.AC.    Abbreviation of boundedly uniformly\,-\,in\,-\,diameter    arcwise\,-\,connected}{135}
\indexentry{closed half space. A set of the form           $\setm{x \in E}{\varphi(x) \geq 0}$,   where $\varphi \in E^*$}{226}
\indexentry{closed half subspace of a normed space}{231}
\indexentry{closed under $E$-discrete composition}{229}
\indexentry{CMP.LUC function.   An extendible function which is   UC at every $x \in \rfs{cl}(X)$}{124}
\indexentry{co-dimension $1$ at $x$.           $\rfs{bd}(X)$ has co-dimension $1$ at $x$}{230}
\indexentry{compatible. $\lambda$ is compatible with $G$ at $x$}{84}
\indexentry{compatible. $\lambda$ is compatible with $G$}{84}
\indexentry{complete cover. $\calU$ is a complete cover of $X$,           if $\bigcup \setm{\overfs{int}(U)}{U \in \calU} =   \overfs{int}(X)$}{53}
\indexentry{completely LUC function}{124}
\indexentry{completely discrete family of sets. A set $\calA$ of    pairwise disjoint sets such that\newline\indent    $\forall B((\forall A \in \calA)(\abs{B \cap A} \leq 1)    \rightarrow \rfs{acc}(B) = \emptyset)$}{169}
\indexentry{completely discrete sequence}{169}
\indexentry{completely discrete set}{169}
\indexentry{completely discrete track system}{208}
\indexentry{completely locally $\itGamma$-bicontinuous}{227}
\indexentry{completely locally $\itGamma$-continuous}{227}
\indexentry{component-wise locally arcwise connected at $x$}{185}
\indexentry{component-wise locally arcwise connected at the boundary}{185}
\indexentry{component-wise wide}{185}
\indexentry{continuous.           $\alpha$-continuous. $f$ is $\alpha$-continuous,           if for every $x,y$, \, $d(f(x),f(y)) \leq d(x,y)$}{14}
\indexentry{continuous.           $\pair{\alpha}{\calP}$-continuous at $x$}{296}
\indexentry{continuous. $(r,\alpha)$-continuous}{108}
\indexentry{continuous. $\alpha$-continuous at $x \in \rfs{cl}(X)$}{281}
\indexentry{continuous. $\alpha$-continuous at $x$. There is $U \in \rfs{Nbr}(x)$ such that $f \nrestriction U$ is $\alpha$-continuous}{63}
\indexentry{continuous. $\beta$-continuous for $\alpha$-submerged pairs.           Abbreviation: $\scolonpair{\beta}{\alpha}$-continuous}{306}
\indexentry{continuous. $\itDelta$-continuous for           $\itGamma$-submerged pairs.        Abbreviation: $\scolonpair{\itDelta}{\itGamma}$-continuous}{306}
\indexentry{continuous. $\itGamma$-continuous at $x$. There is           $\alpha \in \itGamma$ such that   $f$ is $\alpha$-continuous at $x$}{63}
\indexentry{continuous. $\pair{\alpha}{\calP}$-continuous}{296}
\indexentry{continuous. $\pair{\itGamma}{\calP}$-continuous at $x$}{296}
\indexentry{continuous. $\rho$ is $(n,\alpha)$-continuous}{241}
\indexentry{countably generated. $\itGamma$ is countably generated,           if for some           countable $\itGamma_0 \subseteq \itGamma$,   \newline\indent\rule{0pt}{1pt}\kern00mm   $\itGamma \subseteq   \setm{\alpha \in \rfs{MC}}{(\exists \gamma \in \itGamma_0)   (\alpha \preceq \gamma)}$}{15}
\indexentry{CP1 space}{93}
\indexentry{CP1. $X$ is CP1 at $x$}{93}
\indexentry{decayable action.           $\lambda$ is an $(a,\alpha,G)$-decayable action at $x$}{64}
\indexentry{decayable action. $(\alpha,G)$-decayable at $x$. This means $(\dghalf,\alpha,G)$-decayable at $x$}{64}
\indexentry{decayable action. $\alpha$-decayable at $x$. This means $(\dghalf,\alpha,H(X))$-decayable at $x$}{64}
\indexentry{decayable action. $\lambda$ is           an $(a,\itGamma,G)$-decayable action}{64}
\indexentry{decayable action. $\lambda$ is an           $(a,\alpha,G)$-decayable action in $A$}{64}
\indexentry{decayable action. $\lambda$ is an           $(a,\alpha,G)$-decayable action}{64}
\indexentry{determined class. $\calP$\,-\,determined class of           topological spaces}{17}
\indexentry{determining category}{8}
\indexentry{dimension $1$ at $x$.           $\rfs{bd}(X)$ is $1$-dimensional at $x$}{230}
\indexentry{discrete path property for BPD sets}{141}
\indexentry{discrete path property for large distances}{109}
\indexentry{discrete subset. $E$-discrete subset of $\rfs{EXT}^E(X)$}{229}
\indexentry{distinguishable categories}{8}
\indexentry{double boundary point}{185}
\indexentry{DPT. A metric space $X$ is DPT at $x \in X$}{93}
\indexentry{DPT. A metric space is DPT}{93}
\indexentry{e-track system}{208}
\indexentry{e-track. $\pair{\alpha}{\eta}$-e-track}{208}
\indexentry{evasive sequence. $\itGamma$-\kern1ptevasive sequence}{282}
\indexentry{extendible function. A function from $X$ to $Y$   that can be extended to a continuous\\   \indent function from   $\rfs{cl}(X)$ to $\rfs{cl}(Y)$}{124}
\indexentry{extendible homeomorphism}{10}
\indexentry{faithful class of space-group pairs}{8}
\indexentry{faithful class of topological spaces}{8}
\indexentry{fillable. $G$-fillable}{155}
\indexentry{filling. $G$-filling}{155}
\indexentry{finite-dimensional difference homeomorphism}{60}
\indexentry{generated. $\itGamma$ is $(\leq\kern-3pt \kappa)$-generated.           This means $\exists \itGamma_0(\abs{\itGamma_0} \leq \kappa           \mbox{ and } \itGamma = \rfs{cl}_{\preceq}(\itGamma_0))$}{70}
\indexentry{generates. $\itGamma_0$ generates $\itGamma$.           This means           $\itGamma = \rfs{cl}_{\preceq}(\itGamma_0)$}{70}
\indexentry{good semicover. $V$-good semicover}{46}
\indexentry{infinitely-closed. $\alpha$-infinitely-closed at $x$''}{70}
\indexentry{internally extendible in $E$.           A homeomrphism of $X \subseteq E$ which extends to a           \newline\indent\rule{0pt}{1pt}\kern00mm           continuous function on $\overline{\rfs{int}}^E(X)$}{53}
\indexentry{inversely $\pair{K}{\calP}$-continuous}{296}
\indexentry{isotopically movable at the boundary}{191}
\indexentry{isotopically movable with respect to $X$}{191}
\indexentry{JN.AC. Abbreviation of jointly arcwise connected}{169}
\indexentry{JN.ETC}{208}
\indexentry{JN.TC}{208}
\indexentry{joining system}{169}
\indexentry{jointly arcwise connected  }{169}
\indexentry{legal parametrization}{208}
\indexentry{limit-point. $\lambda$-limit-point}{64}
\indexentry{LIN-bordered. $\alpha$-LIN-bordered at $x$}{226}
\indexentry{linear boundary chart domain. A set of the form           $\rfs{BCD}^E(A,r)$, where $A$ is a   closed\\\indent subspace of $E$ different from $\sngltn{0}$   or a closed half space of $E$}{226}
\indexentry{Lipschitz function between    locally Lipschitz normed manifolds}{105}
\indexentry{Lipschitz homeomorphism}{12}
\indexentry{locally $\itGamma$-bicontinuous    with respect to $\itPhi$ and $\itPsi$}{106}
\indexentry{locally $\itGamma$-bicontinuous}{15}
\indexentry{locally $\itGamma$-continuous    with respect to $\itPhi$ and $\itPsi$}{106}
\indexentry{locally $\itGamma$-continuous}{15}
\indexentry{locally $\pair{\alpha}{\calP}$-continuous}{296}
\indexentry{locally almost $\itGamma$-continuous}{95}
\indexentry{locally arcwise connected at a boundary point}{157}
\indexentry{locally arcwise connected at the boundary. Abbreviated by           BR.LC.AC}{157}
\indexentry{locally bi-UC.           Abbreviation of locally bi-uniformly-continuous}{109}
\indexentry{locally bi-uniformly-continuous}{109}
\indexentry{locally bilipschitz homeomorphism}{12}
\indexentry{locally           Lipschitz homeomorphism}{12}
\indexentry{locally           Lipschitz normed manifold}{105}
\indexentry{locally movable at the multiple boundary}{185}
\indexentry{locally moving subgroup of $H(X)$}{33}
\indexentry{locally UC.           Abbreviation of locally uniformly continuous}{109}
\indexentry{locally uniformly continuous}{109}
\indexentry{locally-LIN-bordered. Locally $\itGamma$-LIN-bordered}{226}
\indexentry{LUC on $\rfs{bd}(X)$ function}{124}
\indexentry{manageable ball $B$ based on $S$}{42}
\indexentry{manageable ball (with respect to a BNO-system)}{42}
\indexentry{metrically dense subset}{88}
\indexentry{modulus of continuity}{16}
\indexentry{multiple boundary point}{185}
\indexentry{nearly $\itGamma$-continuous on BPD sets}{140}
\indexentry{nearly open set.    $Z$ is nearly open, if $Z \subseteq \rfs{cl}(\rfs{int}(Z))$}{202}
\indexentry{normed Lipschitz manifold}{204}
\indexentry{normed manifold}{58}
\indexentry{on different sides.           $u,v$ are on different sides of $\rfs{bd}(X)$           with respect to $\trpl{\psi}{A}{r}$}{230}
\indexentry{on the same side. $u,v$ are on the same side           of $\rfs{bd}(X)$ with respect to $\trpl{\psi}{A}{r}$}{230}
\indexentry{open sum partition with respect to $X$}{296}
\indexentry{order preserving at $x$}{282}
\indexentry{order reversing at $x$}{282}
\indexentry{order-irreversible.           $\rfs{bd}(X)$ is $G$-order-irreversible at $x$}{282}
\indexentry{order-reversible. $\rfs{bd}(X)$ is $G$-order-reversible           at $x$}{282}
\indexentry{pairwise disjoint family. A set of pairwise disjoint sets}{169}
\indexentry{partial action of a topological group on a topological space}{64}
\indexentry{PD set. A subset of $X$ whose distance from    the boundary of $X$ is $> 0$}{124}
\indexentry{PD.P function. A function which takes PD sets   to PD sets}{124}
\indexentry{PD.UC function. A function which is uniformly   continuous on every PD set}{124}
\indexentry{piecewise linearly radial. A radial homeomorphism           $\rfs{Rad}^E_{\eta}$ in which $\eta$ is piecewise linear}{235}
\indexentry{point pre-representative}{291}
\indexentry{pointwise $\itGamma$-continuous}{27}
\indexentry{positive distance UC function}{124}
\indexentry{positively distanced set. A subset of $X$ whose distance from    the boundary of $X$ is $> 0$}{124}
\indexentry{prinicipal. $\itGamma$ is principal, if for some           $\alpha \in \itGamma$,           $\itGamma \subseteq           \rfs{cl}_{\spreceq}(\setm{\alpha^{\sscirc n}}{n \in \bbN})$}{15}
\indexentry{Property MV1}{136}
\indexentry{radial homeomorphism based on $\eta$. $\rfs{Rad}^E_{\eta}$.}{81}
\indexentry{radial homeomorphism. $\rfs{Rad}^E_{\eta,z}$.           The radial homeomorphism based on $\eta,z$}{81}
\indexentry{RBM. A regional Banach manifold}{58}
\indexentry{regional Banach manifold (RBM)}{58}
\indexentry{regional normed atlas for $X$}{58}
\indexentry{regionally $\itGamma$-continuous}{27}
\indexentry{regionally normed manifold (RNM)}{58}
\indexentry{regionally translation-like action}{86}
\indexentry{regionally translation-like at $x$}{85}
\indexentry{regular open. A set is regular open, if it is equal to           the interior of its closure}{33}
\indexentry{restricted topological category}{8}
\indexentry{RNM. A regionally normed manifold}{58}
\indexentry{side preserving at $x$}{282}
\indexentry{side reversing at $x$}{282}
\indexentry{simple boundary point}{157}
\indexentry{SLIN-bordered.           $\alpha$-simply-linearly-bordered at $x$           ($\alpha$-SLIN-bordered at $x$)}{230}
\indexentry{small semicover. $V$-small semicover}{46}
\indexentry{small set}{41}
\indexentry{somewhere dense set. A set whose closure contains a           nonempty open set}{34}
\indexentry{space-group pair. $\pair{X}{G}$ is a space-group pair,   if $X$ is a topological space and $G \leq H(X)$}{8}
\indexentry{spaced set of sets. $r$-spaced set of sets}{70}
\indexentry{spaced subset of $X$. $A \subseteq X$ is spaced if           $(\exists r > 0)(\forall x,y \in A)((x \neq y)           \rightarrow (d(x,y) \geq r))$}{70}
\indexentry{spaced track system}{208}
\indexentry{star-closed. $\itGamma$ is $\alpha$-star-closed}{310}
\indexentry{star-closed. $\itGamma$ is $\itDelta$-star-closed}{310}
\indexentry{strongly extendible}{208}
\indexentry{strongly separated. $U$ is strongly separated from $V$,    if $\exists W(U \prec W \mbox{ and } W \cap V = \emptyset$)}{41}
\indexentry{strongly small set}{41}
\indexentry{submerged. $\pair{x}{y}$ is $\alpha$-submerged in $X$.           This means           $\delta^X(x) \geq \norm{x - y} +   \alpha\inverse(\norm{x - y})$}{306}
\indexentry{subspace choice for $\pair{X}{\itPhi}$}{58}
\indexentry{subspace choice system}{38}
\indexentry{subspace choice}{38}
\indexentry{tight Hilbert complementation}{234}
\indexentry{tight Hilbert norm}{233}
\indexentry{topological local movement system}{33}
\indexentry{track system}{208}
\indexentry{track. $\pair{\alpha}{\eta}$-track}{208}
\indexentry{translation-like partial action at $x$}{77}
\indexentry{translation-like partial action}{77}
\indexentry{translation-like.    $\pair{H}{\lambda}$ is $\calP$-translation-like in $L$}{299}
\indexentry{translation-like. $\pair{H}{\lambda}$           is $\calP$-translation-like at $x$}{298}
\indexentry{two-sided. $X$ is two-sided at $x$}{230}
\indexentry{UC around $\rfs{bd}(X)$}{124}
\indexentry{UC at $x$.           Abbreviation of uniformly continuous at $x$}{109}
\indexentry{UC-constant.           $M$ is a 1UC-constant for $\pair{a}{b}$}{274}
\indexentry{UC-constant.           $M$ is a UC-constant for $\pair{a}{b}$}{242}
\indexentry{UC. Abbreviation of uniformly continuous}{108}
\indexentry{UD.AC. Abbreviation of uniformly\,-\,in\,-\,diameter           arcwise\,-\,connected}{122}
\indexentry{uniformly continuous at $x$}{109}
\indexentry{uniformly continuous for all distances}{108}
\indexentry{uniformly continuous}{108}
\indexentry{uniformly\,-\,in\,-\,diameter arcwise\,-\,connected}{122}
\indexentry{weakly $\itGamma$-bicontinuous function}{141}
\indexentry{weakly $\itGamma$-continuous function}{141}
\indexentry{weakly $\itGamma$-continuous on BPD sets}{141}
\indexentry{weakly ``finite-dimensional difference'' homeomorphism}{61}
\indexentry{wide set}{169}
\end{small}
\newpage

\rule{3cm}{0pt}
\vbox{
{\bf Addresses}

\kern 3.0mm

Matatyahu Rubin
\newline\indent
Department of Mathematics,
Ben Gurion University
\newline\indent
Beer Sheva 84105, Israel
\newline\indent
Email: matti@math.bgu.ac.il

\kern3mm

Yosef Yomdin
\newline\indent
Department of Mathematics,
Weizmann Institute
\newline\indent
Rechovot 76100, Israel
\newline\indent
Email: yomdin@wisdom.weizmann.ac.il
}

\end{document}